\documentclass[12pt]{amsart}

\usepackage[colorlinks=true,urlcolor=blue,bookmarks=true,bookmarksopen=true,citecolor=blue
]{hyperref}

\usepackage{graphicx}

\usepackage{epsfig}
\usepackage{amscd}
\usepackage[mathscr]{eucal}
\usepackage{amssymb}
\usepackage{amsxtra}
\usepackage{amsmath}
\usepackage[all]{xy}
\usepackage{enumerate}

\usepackage{color}
\theoremstyle{plain}

\newtheorem{thm}{Theorem}[subsection]
\newtheorem{propspec}[thm]{Proposition$^*$}
\newtheorem{cor}[thm]{Corollary}
\newtheorem{corspec}[thm]{Corollary$^*$}
\newtheorem{lem}[thm]{Lemma}
\newtheorem{prop}[thm]{Proposition}

\theoremstyle{definition}
\newtheorem{dfn}[thm]{Definition}

\newtheorem{assumption}[thm]{Assumption}

\theoremstyle{remark}
\newtheorem{rem}[thm]{Remark}
\newtheorem*{remnonum}{Remark}

\newtheorem*{remsnonum}{Remarks}
\newtheorem{ex}[thm]{Example}

\theoremstyle{plain}
\newtheorem{statement}[thm]{}

\newcommand{\pbaddress}{biran@math.tau.ac.il}
\newcommand{\ocaddress}{cornea@dms.umontreal.ca}

\newcommand{\Qed}{\hfill \qedsymbol \medskip}

\newcommand{\Id}{{{\mathchoice {\rm 1\mskip-4mu l} {\rm 1\mskip-4mu l}
      {\rm 1\mskip-4.5mu l} {\rm 1\mskip-5mu l}}}}
\newcommand{\hooklongrightarrow}{\lhook\joinrel\longrightarrow}

\newcommand{\R}{\mathbb{R}}
\newcommand{\Z}{\mathbb{Z}}
\newcommand{\Q}{\mathbb{Q}}
\newcommand{\N}{\mathbb{N}}
\newcommand{\C}{\mathbb{C}}
\newcommand{\la}{\lambda}
\newcommand{\La}{\Lambda}

\newcommand{\Crit}{\rm{Crit\/}}

\newcommand{\FF}{{ \,  I  \hspace{-.11cm}   { F}  }}

\newcommand{\ind}{{\textnormal{ind}}}
\newcommand{\mubar}{{\bar{\mu}}}

\begin{document}

\title{Quantum structures for Lagrangian submanifolds}
\date{\today}

\thanks{The first author was partially supported by the ISRAEL
SCIENCE FOUNDATION  (grant No. 1227/06 *); the second author was supported by an NSERC Discovery grant and a FQRNT Group Research grant}

\author{Paul Biran and Octav Cornea} \address{Paul Biran, School of
  Mathematical Sciences, Tel-Aviv University, Ra mat-Aviv, Tel-Aviv
  69978, Israel} \email{\pbaddress} \address{Octav Cornea, Department
  of Mathematics and Statistics University of Montreal C.P. 6128 Succ.
  Centre-Ville Montreal, QC H3C 3J7, Canada} \email{\ocaddress}

\bibliographystyle{plain}

%
%

\maketitle

%
%
\tableofcontents

%

\section{Introduction.}
It is well-known that Gromov-Witten invariants are not defined, in
general, in the Lagrangian (or relative case): bubbling of disks is a
co-dimension one phenomenon and thus counting $J$-holomorphic disks,
possibly with various incidence conditions, produces numbers that
strongly depend on the particular choices of the almost complex
structure $J$ and of the geometry of the incidences. However, this lack
of invariance of the direct counts combined with the rich
combinatorial properties of the moduli spaces of disks indicates that
invariance can still be achieved by defining an appropriate homology
theory.

\

The purpose of this paper is to discuss systematically such a homology
theory and the related algebraic structures in the case of monotone
Lagrangians with minimal Maslov class at least $2$ (which we will
shortly call the {\em monotone case} below).  We will also discuss its
relations with Floer homology as well as various computations,
examples and applications. It is already important to underline the
fact that the point of view of this paper is not that of intersection
theory and, thus, not that of Floer theory. In particular, this
homology theory (and most of the other structures involved) is
associated to a single Lagrangian submanifold, never vanishes and is
invariant with respect to ambient symplectic isotopy.  Therefore, this
is a very rich structure and most of our various applications reflect
its rigidity.

\

It should be mentioned at the outset that there are two systematic
models for dealing with the general, non-monotone context: the
$A_{\infty}$ approach of Fukaya-Oh-Ohta-Ono \cite{FO3} and the cluster
homology approach of Cornea-Lalonde \cite{Cor-La:Cluster-1} (the last
one being closer to the point of view we take here).  The monotone
case, while remaining reasonably rich, has the property that many of
the technical complications which are present in these two models
disappear.  Indeed, as we shall see in this case, transversality
issues can be dealt with by elementary means and the various
invariants defined are all based on counting $J$-curves (for generic
almost complex structures $J$) and not perturbed objects, a fact which
is of invaluable help in computations.  Hence, this is a case which is
worth exploring in detail not only because many of the most relevant
examples of Lagrangians fit in this context but also because it shows
clearly what type of results and applications can be expected in
general and computations are efficient.

\

There are many relations between this work and the extensive
literature of the subject and they will be discussed explicitly later
in the paper. We feel we should mention at this point that at the
center of the construction is a chain complex called here the ``pearl
complex'' initially described by Oh in \cite{Oh:relative}. Thus, this
complex was known before and it was sometimes used in the literature.
Our intention in preparing this paper has been to focus on the various
structures related to this complex - a good number of which are first
introduced here - and on their applications.  However, we soon
realized that no complete proofs are available in what concerns even
the most basic parts of the construction, for example, for $d^{2}=0$
in the ``pearl complex''.  Therefore we have decided to provide
essentially complete arguments here. The structure of the paper is as
follows. Section \ref{sec:algstr} contains the statement of the main
algebraic properties of the structures that we are interested in here.
The next three sections \S\ref{S:transversality}, \S\ref{S:gluing},
\S\ref{S:proofalg} are focused, respectively, on transversality,
gluing and, finally, the proof of the algebraic properties announced
in \S\ref{sec:algstr}.  Then, in \S\ref{sec:appli} we describe our
applications of this structure together with their proofs.  A
cautionary word to the reader: while the paper is written in the
logical order needed to reduce redundancies, to more rapidly grasp the
power and the motivation behind the algebraic structures described in
the paper it might be useful to read the statements of the
applications in \S\ref{sec:appli} immediately after \S\ref{sec:algstr}
and before the technical chapters in between.

\subsubsection*{Acknowledgements.}
The first author would like to thank Kenji Fukaya, Hiroshi Ohta, and
Kaoru Ono for valuable discussions on the gluing procedure for
holomorphic disks. He would also like to thank Martin Guest and Manabu
Akaho for interesting discussions and great hospitality at the Tokyo
Metropolitan University during the summer of 2006. Special thanks to Leonid
Polterovich and Joseph Bernstein for interesting comments and their
interest in this project from its early stages. 

While working on this project our two children, 
Zohar and Robert, have been born. We would like to dedicate 
this work to them and to their lovely mothers, 
Michal and Alina. 



%

\section{Algebraic structures.}\label{sec:algstr}

\subsection{The main algebraic statement.}
We shall only work in this paper with, connected, closed, monotone
Lagrangians $L\subset (M^{2n},\omega)$ where $(M,\omega)$ is a tame
monotone symplectic manifold.  This means that the two
homomorphisms:
$$\omega:\pi_{2}(M,L)\to \Z, \quad \mu:\pi_{2}(M,L)\to \R$$
given
respectively by integration of $\omega$ and by the Maslov index
satisfy $$\omega(A) > 0 \quad \textnormal{iff} \quad \mu(A)>0, \quad
\forall\; A \in \pi_2(M,L).$$
It is easy to see that this is
equivalent to the existence of a constant $\tau > 0$ such that
\begin{equation}\label{eq:monotonicity}
   \omega(A)=\tau \mu(A),\ \forall \ A \in
   \pi_{2}(M,L)~.~
\end{equation}
We shall refer to $\tau$ as the {\em monotonicity constant} of $L
\subset (M,\omega)$. Define the {\em minimal Maslov number} of $L$ to
be the integer $$N_L = \min \{ \mu(A)>0 \mid A \in \pi_2(M,L) \}.$$
Throughout this paper we shall assume that $L$ is monotone with $N_L
\geq 2$. Since the Maslov numbers come in multiples of $N_L$ we shall
use sometimes the following notation:
\begin{equation}
   \mubar = \tfrac{1}{N_L}\mu : \pi_2(M,L) \to \mathbb{Z}.
\end{equation}

Let us now introduce the type of coefficient rings we shall work with.
Let $\Lambda = \mathbb{Z}_2[t,t^{-1}]$ the ring of Laurent polynomials
in the varible $t$. The grading on $\Lambda$ is given by $\deg t =
-N_L$. We shall also work with the positive version of $\Lambda$,
namely $\Lambda^+ = \mathbb{Z}_2[t]$ with the same grading.  The ring
$\Lambda$ should be viewed as a simplified version of the Novikov ring
over $\pi_2(M,L)$, commonly used in Floer theory. Since we work in the
monotone case the simplified Novikov rings $\Lambda$, $\Lambda^{+}$
are enough for our purposes.

There is a natural decreasing filtration of both $\La^{+}$ and $\La$
by the degrees of $t$, i.e.
\begin{equation}\label{eq:filtration}
   \mathcal{F}^{k}\La = \{P\in \mathbb{Z}_2[t,t^{-1}] \mid
   P(t)=a_{k}t^{k}+a_{k+1}t^{k+1}+\ldots \}~.~
\end{equation} We shall call this
filtration the \emph{degree filtration}. It induces an obvious
filtration on any free module over this ring.

Let $f:L\to \R$ be a Morse function on $L$ and let $\rho$ be a
Riemannian metric on $L$ so that the pair $(f,\rho)$ is Morse-Smale.
Fix also a generic almost complex structure $J$ compatible with
$\omega$.  It is well known that, under the above assumption of
monotonicity, the Floer homology of the pair $(L,L)$ is well defined
(see~\cite{Oh:HF1}) and we denote it by $HF(L)$ (the construction will
be rapidly reviewed later in the paper).

In what follows we shall also use the following version of the quantum
homology of $M$. Put $Q^{+}H(M) = H_*(M;\mathbb{Z}_2) \otimes
\Lambda^{+}$, $QH(M) = H_*(M;\mathbb{Z}_2) \otimes \Lambda$ with the
grading induced from $H_*(M;\mathbb{Z}_2)$ and $\Lambda^{+}$
(respectively $\Lambda$). We endow $Q^{+}H(M)$ and $QH(M)$ with the
quantum cap product (see~\cite{McD-Sa:Jhol-2} for the definition).
There are a few slight differences in our convention in comparison to
the ones common in the literature. The first is that the degree of the
variable $t$ in the quantum homology of $M$ is usually minus the
minimal Chern number $-N_M$ of $(M,\omega)$.  In our setting we have
$\deg t = -N_L$.  Since we are in the monotone case we have $N_L |
N_M$, thus our $QH(M)$ is actually a kind of extension of the usual
quantum homology of $M$. The second difference is that we work here
with coefficients over $\mathbb{Z}_2$ rather than $\mathbb{Q}$ or
$\mathbb{Z}$ which are more common in quantum homology theory. This is
not essential and has to do with technical issues concerning the
definition of the Floer homology of $L$ (see
Remark~\ref{rem:extensions} below). Finally note that we work here
with quantum {\em homology} (not cohomology), hence the quantum
product $QH_k(M) \otimes QH_l(M) \to QH_{k+l-2n}(M)$ has degree $-2n$.
The unit is $[M] \in QH_{2n}(M)$, thus of degree $2n$.

\begin{thm}\label{thm:alg_main}
   For a generic choice of the triple $(f,\rho,J)$ there exists a chain
   complex $$\mathcal{C}^{+}(L;f,\rho,J)=(\Z_2 \langle \Crit(f) \rangle
   \otimes \La^{+}, d)$$
   with the following properties:
   \begin{itemize}
     \item[i.] The homology of this chain complex is independent of
      the choices of $J, f, \rho$. It will be denoted by
      $Q^{+}H_{\ast}(L)$.  There exists a canonical (degree
      preserving) augmentation $\epsilon_{L}: Q^{+}H_{\ast}(L)\to
      \La^{+}$ which is a $\Lambda^{+}$-module map.
     \item[ii.] The homology $Q^{+}H(L)$ has the structure of a
      two-sided algebra with unit over the quantum homology of $M$,
      $Q^{+}H(M)$. More specifically, for every $i$, $j$, $k$ there exist
      $\Lambda^{+}$-bilinear maps:
      \begin{align*}
         & Q^{+}H_i(L) \otimes Q^{+}H_j(L) \to Q^{+}H_{i+j-n}(L),
         \quad \alpha\otimes \beta \mapsto \alpha*\beta, \\
         & Q^{+}H_k (M) \otimes Q^{+}H_i(L) \to Q^{+}H_{k+i-2n}(L),
         \quad a\otimes \alpha \mapsto a*\alpha.
      \end{align*}
      The first map endows $Q^{+}H(L)$ with the structure of a ring
      with unit. This ring is in general not commutative. The second
      map endows $Q^{+}H(L)$ with the structure of a module over the
      quantum homology ring $Q^{+}H(M)$. Moreover, when viewing these
      two structures together, the ring $Q^{+}H(L)$ becomes a
      two-sided algebra over the ring $Q^{+}H(M)$. The unit $[M]$ of
      $Q^{+}H(M)$ has degree $2n=\dim M$ and the unit of $Q^{+}H(L)$
      has degree $n=\dim L$.
     \item[iii.]  There exists a map
      $$i_{L}:Q^{+}H_{\ast}(L)\to Q^{+}H_{\ast}(M)$$
      which is a
      $Q^{+}H_{\ast}(M)$-module morphism and which extends the inclusion
      in singular homology. This map is determined by the relation:
      \begin{equation}\label{eq:inclusion_mod}
         \langle h^{\ast}, i_{L}(x) \rangle =\epsilon_{L}(h\ast x)
      \end{equation} for $x\in Q^{+}H(L)$,
      $h\in H_{\ast}(M)$, with $(-)^{\ast}$
      Poincar\'e duality and $\langle - , - \rangle$ the Kronecker pairing.
     \item[iv.] The differential $d$ respects the degree filtration and all
      the structures above are compatible with the resulting spectral
      sequences.
     \item[v.]The homology of the complex:
      $$\mathcal{C}(L;f,\rho,J)=\mathcal{C}^{+}(L;f,\rho,J)
      \otimes_{\La^{+}}\La$$
      is denoted by $QH_{\ast}(L)$ and all the
      points above remain true if using $QH(-)$ instead of
      $Q^{+}H(-)$.  The map $Q^{+}H(L)\to QH(L)$ induced in homology by the change of coefficients above is canonical. Moreover, there is an isomorphism
      $$QH_{\ast}(L)\to HF_{\ast}(L)$$
      which is also canonical up to a shift in
      grading.
   \end{itemize}

\end{thm}

By a two-sided algebra $A$ over a ring $R$ we mean that $A$ is a
module over $R$ which has an internal product $A\times A\to A$ so that
for any $r\in R$ and $a,b\in A$ we have $r(ab)=(ra)b=a(rb)$. The last
equality is non-trivial, of course, only when the product in $A$ is
not commutative. A more natural description is the following.  If $A$
is a (left)-module over $R$, define a right-action of $R$ on $A$ by
$ar = ra$. Then the ``two-sidedness'' of $A$ over $R$ means that {\em
  both} actions give $A$ the structure of a module over $R$.

Before going on any further we would like to point out that, the
existence of a module structure asserted by Theorem~\ref{thm:alg_main}
has already some non-trivial consequences. For example, the fact that
$QH_*(L) \cong HF_*(L)$ is a module over $QH_*(M)$ implies that if $a
\in QH_k(M)$ is an invertible element, then the map $a*(-)$ gives rise
to {\em isomorphisms} $HF_i(L) \to HF_{i+k-2n}(L)$ for every $i \in
\mathbb{Z}$. This clearly follows from the general algebraic
definition of a ``module over a ring with unit''.

\

We shall call the complex $\mathcal{C}(L;f,\rho,J)$ (respectively,
$\mathcal{C}^{+}(L;f,\rho,J)$) the {\em (positive) pearl complex}
associated to $f,\rho,J$ and we shall call the resulting homology the
(positive) {\em quantum homology} of $L$. In the perspective
of~\cite{Cor-La:Cluster-1, Cor-La:Cluster-2} the complex
$\mathcal{C}(L;f,\rho,J)$ corresponds to the {\em linear cluster
  complex}.

\

Parts of Theorem~\ref{thm:alg_main} appear already in the literature
and have been verified up to various degrees of rigor.  The complex
$\mathcal{C}(L;f,\rho,J)$ has been first introduced by
Oh~\cite{Oh:relative} (see also Fukaya~\cite{Fu:Morse-homotopy}) and
is a particular case of the cluster complex as described in
Cornea-Lalonde~\cite{Cor-La:Cluster-1}.  The module structure over
$Q^{+}H(M)$ discussed at point ii.  is probably known by experts - at
least in the Floer homology setting - but has not been explicitly
described yet in the literature.  The product at ii. is a variant of
the usual pair of pants product - it might not be widely known in this
form.  The map $i_{L}$ at point iii. is the analogue of a map first
studied by Albers in~\cite{Alb:extrinisic} in the absence of bubbling.
The spectral sequence appearing at iv. is a variant of the spectral
sequence introduced by Oh \cite{Oh:spectral}. The compatibility of
this spectral sequence with the product at point ii. has been first
mentioned and used by Buhovsky~\cite{Bu:products} and independently by
Fukaya-Oh-Ohta-Ono~\cite{FO3}. The positive Novikov ring $\Lambda^{+}$
is commonly used in algebraic geometry as well as in the closed case
and has appeared in the Lagrangian setting in
Fukaya-Oh-Ohta-Ono~\cite{FO3}. The comparison map at v. is an
extension of the Piunikin-Salamon-Schwarz construction~\cite{PSS}, it
extends also the partial map constructed by Albers in~\cite{Alb:PSS}
and a more general such map was described in~\cite{Cor-La:Cluster-1}
in the ``cluster" context. We also remark that this comparison map
identifies all the algebraic structures described above with the
corresponding ones defined in terms of the Floer complex.

\begin{rem} \label{rem:extensions} a.  It is quite clear that, with rather
   obvious modifications, all the structure described in this
   statement should carry over to the case when $L$ is non-monotone
   but orientable and relative spin. The coefficients in that case
   have to be rational - obviously this requires that the various
   moduli spaces involved be oriented coherently. One option to pursue
   the construction in this case is to further replace the Novikov
   ring $\La$ (or the positive Novikov ring $\La^{+}$) with a cluster
   complex $\mathcal{C}l(L;J)$ of $L$ \cite{Cor-La:Cluster-1}. Using
   these ``cluster" coefficients means that the complex
   $\mathcal{C}(L;f,\rho,J)$ is replaced with the fine Floer complex of
   \cite{Cor-La:Cluster-1} and, with the exception that $QH_{\ast}(L)$
   is replaced in all places by the fine Floer homology of $L$, $\FF
   H_{\ast}(L)$, the statement of Theorem \ref{thm:alg_main} should
   remain true and even have a ``positive" version.

   b. Another interesting point that we want to emphasize here - and
   will be exemplified later in the paper - is that the structures
   discussed in the statement of the theorem lead to the definition of
   certain Gromov-Witten type invariants.  The procedure is as
   follows.  Suppose first that $k\in \Z_2$ (or $\in\Z$, if we assume
   orientations) is some numerical invariant defined out of the
   algebra structure of $Q^{+}H_{\ast}(L)$ (this means that this
   number is left invariant by isomorphisms of the structure).  Assume
   also that, under certain circumstances, for special choices of the
   function $f$ and the almost complex structure $J$, the chain
   complex $\mathcal{C}^{+}(L;f,\rho,J)$ has a trivial differential. This
   could happen, for example, if $f$ is a perfect Morse function and
   if $J$ is a special ``symmetric" structure or, for example, as we
   shall see further in this paper, if $L$ is a torus with non trivial
   Floer homology. In that case, the ``counting" leading to $k$, which
   is invariant, in general, only after passage to homology will be
   invariant already at the chain level simply because, for these
   special choices of $f,J$, the chain level is isomorphic with the
   homology one.  But this means that, with these special choices, the
   count giving $k$ is invariant and this is exactly what is needed to
   define Gromov-Witten type invariants.  It is then another matter to
   interpret these numbers geometrically in a meaningful way.
 \end{rem}

\subsection{Other algebraic structures}\label{subsec:other_str}

\subsubsection{Duality}\label{susubsubsec:duality}
The first point that we want to discuss here is a form of
\emph{duality} which extends Poincar\'e duality.  We first fix some
notation.  Suppose that $(\mathcal{C},\partial)$ is a chain complex
over $\La^{+}$. In particular, it is a free module over $\La^{+}$,
$\mathcal{C}=G\otimes \La^{+}$ with $G$ some $\mathbb{Z}_2$ vector
space. We let
$$\mathcal{C}^{\odot}=\hom_{\La^{+}}(\mathcal{C},\La^{+})$$
graded so
that the degree of a morphism $g:\mathcal{C}\to \La^{+}$ is $k$ if $g$
takes $\mathcal{C}_{l}$ to $\La^{+}_{l+k}$ for all $l$.

Let $ \mathcal{C}'=\hom_{\Z_2}(G,\Z_2)\otimes \Lambda^{+}$ be graded
such that if $x$ is a basis element of $G$, then its dual $x^{\ast}\in
\mathcal{C}'$ has degree $|x^{\ast}|=-|x|$. There is an obvious degree
preserving isomorphism $\psi:\mathcal{C}^{\odot}\to \mathcal{C}'$
defined by $\psi(f)=\sum_{i}f(g_{i})g_{i}^{\ast}$ where $(g_{i})$ is a
basis of $G$ and $(g_{i}^{\ast})$ is the dual basis.  We define the
differential of $\mathcal{C}^{\odot}$, $\partial ^{\ast}$, as the
adjoint of $\partial$:
$$
\langle \partial^{\ast} y^{\ast},x \rangle = \langle
y^{\ast},\partial x \rangle \ ,\ \forall x,y\in G~.~$$
Clearly,
$\mathcal{C}^{\odot}$ continues to be a chain complex (and not a
co-chain complex).

An additional algebraic notion will be useful: the co-chain
complex  $\mathcal{C}^{\ast}$ associated to $\mathcal{C}$. To define it
we first let $(\La^{+})^{\ast}$ be the ring $\La^{+}$ with the reverse grading: the degree
of each element in $(\La^{+})^{\ast}$ is the opposite of the degree of the same element in $\La^{+}$.
For the free chain complex $\mathcal{C}=G\otimes
\La^{+}$ as before, we define
$\mathcal{C}^{\ast}=\hom_{\mathbb{Z}_2}(G,\mathbb{Z}_2)\otimes
(\La^{+})^{\ast}$ where the grading of the dual $x^{\ast}$ of a basis
element $x\in G$ is $|x^{\ast}|=|x|$.  The differential in
$\mathcal{C}^{\ast}$ is given as usual as the adjoint of the
differential in $\mathcal{C}$.  The complex $\mathcal{C}^{\ast}$ is
obviously a co-chain complex. The difference between this complex and
$\mathcal{C}^{\odot}$ is just that the grading is reversed in the sense
that if an element $x\otimes \la$ has degree $k$ in one complex, then
it has degree $-k$ in the other. The co-homology of $\mathcal{C}$ is then
defined as $H^{k}(\mathcal{C})=H^{k}(\mathcal{C}^{\ast})$. Obviously,
there is a canonical isomorphism:  $H_{-k}(\mathcal{C}^{\odot})
\cong H^{k}(\mathcal{C}^{\ast})$.

A particular case of interest here is when $\mathcal{C}=\mathcal{C}(L; f,\rho,J)$.
In this case we denote:  $$Q^{+}H^{n-k}(L)=H^{k}(\mathcal{C}^{+}(L;f,\rho,J)^{\ast})~.~$$

Notice that the chain morphisms $\eta: \mathcal{C}\to
\mathcal{C}^{\odot}$ of degree $-n$ are in correspondence with the
chain morphisms of degree $-n$:
$$\tilde\eta:\mathcal{C}\otimes_{\Lambda^{+}} \mathcal{C}\to
\La^{+}~.~$$
via the formula $\tilde{\eta}(x\otimes y)=\eta(x)(y)$.
Here the ring $\Lambda^{+}$ on the right handside is considered as a chain
complex with trivial differential. Moreover, if $\eta$ induces an
isomorphism in homology, then the pairing induced in homology by
$\tilde\eta$ is non-degenerate.

Fix now $n\in \N^{\ast}$. For any chain complex $\mathcal{C}$ as
before we let $s^{n}\mathcal{C}$ be its $n$-fold suspension. This is
a chain complex which coincides with $\mathcal{C}$ but its graded so
that the degree of $x$ in $s^{n}\mathcal{C}$ is $n+$ the degree of
$x$ in $\mathcal{C}$.

A particular useful case where these notions appear
is in the following sequence of obvious isomorphisms:
 $H_{k}(s^{n}\mathcal{C}^{\odot})\cong H_{k-n}(\mathcal{C}^{\odot})\cong H^{n-k}(\mathcal{C}^{\ast})$.

\begin{cor}\label{cor:duality} Set $n = \dim L$. There exists
   a degree preserving morphism of chain complexes:
   $$\eta: \mathcal{C}^{+}(L;f,\rho,J)\to
   s^{n}(\mathcal{C}^{+}(L;f,\rho,J))^{\odot}$$
   which induces an
   isomorphism in homology. In particular, we have an isomorphism: $\eta: Q^{+}H_{k}(L)\to Q^{+}H^{n-k}(L)$. The corresponding (degree $-n$) bilinear
   map
   $$H(\tilde\eta):Q^{+}H(L)\otimes Q^{+}H(L)\to \La^{+}$$
   coincides
   with the product described in Theorem \ref{thm:alg_main}-ii
   composed with the augmentation $\epsilon_{L}$.
    The same result continues to hold with $\Lambda^{+}$,
   $\mathcal{C}^{+}$, $QH^{+}$ replaced by $\Lambda$, $\mathcal{C}$,
   $QH$ respectively.
\end{cor}

\begin{rem}\label{rem:duality}
   a. The relation of the Corollary above with Poincar\'e duality is
   as follows: in case $\mathcal{C}^{+}(-)$ in the statement is
   replaced with the Morse complex $C(f)$ of some Morse function
   $f:L\to \R$ we may define the morphism $\eta:C(f)\to
   s^{n}(C(f)^{\odot})$ as a composition of two morphisms with the
   first being the usual comparison morphism $C(f)\to C(-f)$ and the
   second $C(-f)\to s^{n}(C(f)^{\odot})$ given by $x\in\Crit(f)\to
   x^{\ast}\in\hom_{\Z_2}(C(f),\Z_2)$. We have the
 identifications
   $H_{k}(s^{n}(C(f)^{\odot}))=H_{k-n}(C(f)^{\odot})=H^{n-k}(C(f))$
and the morphism $\eta$ described above induces in
 homology the Poincar\'e duality map: $H_{k}(L)\to H^{n-k}(L)$.

   b. The last Corollary also obviously shows that $Q^{+}H(L)$
   together with the bilinear map $\epsilon_{L}\circ (-\ast-)$ is a
   Frobenius algebra, though not necessarily commutative.
\end{rem}

\subsubsection{Action of the symplectomorphism group.}
This property is very useful in computations when symmetry is present.

\begin{cor}\label{cor:symm} Let $\phi:L\to L$ be a diffeomorphism which
   is the restriction to $L$ of an ambient symplectic diffeomorphism
   $\bar{\phi}$ of $M$. Let $f,\rho,J$ be so that the pearl complex
   $\mathcal{C}^{+}(L;f,\rho,J)$ is defined. There exists a chain map:
   $$\tilde{\phi}:\mathcal{C}^{+}(L;f,\rho,J)\to
   \mathcal{C}^{+}(L;f,\rho,J)$$
   which respects the degree filtration,
   induces an isomorphism in homology, and so that the morphism
   $E^{2}(\tilde{\phi})$ induced by $\tilde{\phi}$ at the $E^{2}$
   level of the degree spectral sequence coincides with
   $H_{\ast}(\phi)\otimes {id}_{\Lambda^{+}}$.  The map $\bar{\phi}\to
   \tilde{\phi}$ induces a representation:
   $$\hbar:Symp (M,L)\to Aut (Q^{+}H_{\ast}(L))$$
   where
   $Aut(Q^{+}H_{\ast}(L))$ are the augmented ring automorphisms of
   $Q^{+}H_{\ast}(L)$ and $Symp(M,L)$ are the symplectomorphisms of
   $M$ which restrict to diffeomorphisms of $L$. The restriction of
   $\hbar$ to $Symp_{0}(M)\cap Symp(M,L)$ takes values in the
   automorphisms of $Q^{+}H(L)$ as an algebra over $Q^{+}H(M)$.

   The same result continues to hold with $\Lambda^{+}$,
   $\mathcal{C}^{+}$, $QH^{+}$ replaced by $\Lambda$, $\mathcal{C}$,
   $QH$ respectively.
\end{cor}

\subsubsection{Minimal pearl complexes.} \label{subsec:minimal}

It is easy to see that all the calculations with the structures described above are much more
efficient if the Lagrangian $L$ admits a perfect Morse function - that is a Morse function
$f:L\to \R$ so that the Morse differential vanishes. We now want to notice that there exists
an algebraic procedure which allows one to treat any general $L$ in the same way. Moreover,
we will see that this produces another a chain complex which is a quantum invariant of $L$
and contains all the quantum specific properties that we generally want to study (a similar construction
in the cluster set-up has been sketched in \cite{Cor-La:Cluster-1}).

Let $G$ be a finite dimensional
graded $\Z_2$-vector space and let $\mathcal{C}=(G\otimes \La^{+}, d)$ be a chain complex.
For an element $x\in G$ let
$d(x)=d_{0}(x)+d_{1}(x)t$ with $d_{0}(x)\in G$. In other words $d_{0}$ is obtained from $d(x)$ by treating $t$ as
a polynomial variable and putting $t=0$. Clearly $d_{0}:G\to G$, $d_{0}^{2}=0$.
Let $\mathcal{H}$ be the homology of the complex $(G,d_{0})$.
Similarly, for a chain morphism
$\xi$ we denote by $\xi_{0}$ the $d_{0}$-chain morphism obtained by making $t=0$.

\begin{prop} \label{prop:min_model} With the notation above there exists a chain complex
 $$\mathcal{C}_{min}=(\mathcal{H}\otimes \La^{+},\delta), \ {\rm with}\ \delta_{0}=0$$ and chain maps $\phi:\mathcal{C}\to \mathcal{C}_{min}$, $\psi:\mathcal{C}_{min}\to \mathcal{C}$ so that
 $\phi\circ\psi=id$ and
 $\phi$ and $\psi$ induce isomorphisms in $d$-homology and $\phi_{0}$ and $\psi_{0}$ induce an isomorphism in $d_{0}$-homology.  Moreover, the properties above characterize
 $\mathcal{C}_{min}$ up to isomorphism.\end{prop}

Here is an important consequence of this result:

 \begin{cor}\label{cor:min_pearls}
 There exists a complex $\mathcal{C}_{min}(L)=(H_{\ast}(L;\Z_{2})\otimes \La^{+},\delta)$, with $\delta_{0}=0$ and so that
 for any  $(L,f,\rho,J)$  such that $\mathcal{C}(L;f,\rho,J)$
 is defined there are chain morphisms
 $\phi:\mathcal{C}(L;f,\rho,J)\to \mathcal{C}_{min}(L)$ and $\psi : \mathcal{C}_{min}(L)\to \mathcal{C}(L;f,\rho,J)$ which both induce isomorphisms in quantum homology as well as in Morse homology and so that $\phi\circ\psi=id$. The complex $\mathcal{C}_{min}(L)$ with these properties is unique up to isomorphism.
 \end{cor}

 We call the complex provided by this corollary the  {\em minimal pearl complex}. This terminology is
 justified by the use of minimal models in rational homotopy where a somewhat similar notion is central. There is a slight
 abuse in this notation as, while any two complexes as provided by the corollary are isomorphic this isomorphism is not canonical. Obviously, in case a perfect Morse function exists on $L$ any pearl complex associated
to such a function is already minimal.

 \begin{rem}\label{rem:product_min} a. An important consequence of the existence of the chain morphisms $\phi$ and $\psi$ is that
 all the algebraic structures described before (product, module structure etc) can be transported and computed
 on the minimal complex. For example, the product is the composition:
 $$\mathcal{C}_{min}(L)\otimes \mathcal{C}_{min}(L)\stackrel{\psi\otimes\psi}{\longrightarrow}\mathcal{C}(L;f,\rho,J)\otimes \mathcal{C}(L;f,\rho,J)\stackrel{\ast}{\to}\mathcal{C}(L;f,\rho,J)\stackrel{\phi}{\to}\mathcal{C}_{min}(L)~.~$$
 It is easy to see that - in homology -  the resulting product has as unit the fundamental class $[L]\in H_{n}(L)$.

 b. A consequence of point a.  is that $HF(L)\cong QH(L)=0$ iff there is some $x\in H_{\ast}(L; \Z_{2})$
 so that $\delta x=[L]t^{k}$ in $\mathcal{C}_{min}(L)$. Indeed, suppose that $QH(L)=0$. Then, as for degree reasons $[L]$ is a cycle
 in $\mathcal{C}_{min}(L)$, we obtain that it has to be also a boundary. Conversely, if $[L]$ is
 a boundary (which means $\delta x= [L] t^{k}$ for some $x$ and $k$) we have for any other
 cycle $c\in \mathcal{C}_{min}(L)$: $[c]=[c]\ast [L]=[c\ast [L]]=[\delta (c\ast x)]=0$ where
 we have denoted by $-\ast -$ the product on $\mathcal{C}_{min}(L)$ as defined above (see also
 \S\ref{Sb:criteria-QH} for other criteria of similar nature).

 c. It is also useful to note that there is an isomorphism $Q^{+}H(L)\cong H(L;\Z_{2})\otimes \La^{+}$
 iff the differential $\delta$ in $\mathcal{C}_{min}(L)$ is identically zero.
 \end{rem}

\subsubsection{Large and small coefficient rings.} \label{subsubsec:coeff}
We have seen before that both the ring $\La^{+}$ and the ring $\La$ can be used
in our constructions. Indeed, the interest of $\La^{+}$ is mainly that the resulting
homology never vanishes while the ring $\La$ is needed for the comparison with
Floer homology.

\begin{ex}\label{ex:circle_la+}
Consider $S^{1}\subset \C$ the standard circle in the complex plane. Obviously, as
$S^{1}$ is displaceable we have that $QH(S^{1})=HF(S^{1})=0$. However, the positive quantum
homology $Q^{+}H(S^{1})$ verifies: $Q^{+}H_{\ast}(S^{1})=0$ for $\ast\not=1$ and
$Q^{+}H_{1}(S^{1})=\Z_2$. Indeed, we may take on $S^{1}$ a Morse function   with a
single minimum $P$ and a single maximum $Q$. The standard almost complex structure is regular and
the standard disk which fills the circle is of Maslov class two. This disks obviously goes through the
minimum and the maximum and this shows that $d P = Qt$ in the pearl complex. It is easy to see
that $d Q=0$ for degree reasons.
Therefore, $Q$ is a cycle but not a boundary and this implies the claim.

This example generalizes easily
to show that for any monotone Lagrangian $L$ we have $Q^{+}H(L)\not=0$. Indeed, as $L$ is assumed
connected we may work with a Morse function $f:L\to \R$ with a single maximum
which we will again denote by $Q$. In this case we again have $d Q=0$. Indeed, the Morse differential
of $Q$ is null because $Q$ is  the unique maximum of $f$. Moreover, if $d Q=R t^{k}+...$ we
 need $|R|-kN_{L}=Q-1$ which is not possible for $k\not=0$ because $N_{L}\geq 2$.
 Thus, the unique maximum of such a function represents a cycle in the pearl complex and, given that $t$ is not invertible in $\La^{+}$, it follows that the homology
 class represented by the maximum is non-trivial. \end{ex}

 In a rather obvious way these are the minimal rings that one
can use for these purposes. Indeed, let $\pi_{2}(M,L)^{+}$ be the semi-group
of all the elements $u$ so that $\omega(u)\geq 0$. Then $\La^{+}=\Z_{2}[\pi_{2}(M,L)^{+}/\sim]$
with $\sim$ the equivalence relation $u\sim v$ iff $\mu(u)=\mu(v)$ and similarly
$\La=\Z_{2}[\pi_{2}(M,L)/\sim]$. For certain other applications it can be useful to also use large rings
which distinguish explicitely elements in $\pi_{2}(M,L)$.  For this purpose we remark now
that all the arguments in the paper carry over when replacing $\La^{+}$ with  $\tilde{\La}^{+}=
\Z_{2}[\pi_{2}(M,L)^{+}]$ and $\La$ with $\tilde{\La}=\Z_{2}[\pi_{2}(M,L)]$.

Indeed, with a single exception to be discussed below, for all  the constructions
in the paper to hold the coefficient ring needs to satisfy just two conditions: it needs to behave additively with respect to gluing and bubbling and it needs to distinguish disks with different symplectic areas (due to our monotonicity assumption this is, of course, the same as distinguishing disks with different Maslov classes).
In particular, any ring $\mathcal{R}$ such that there is a ring morphism $r:\Z_{2}[\pi_{2}(M,L)^{+}]\to \mathcal{R}$ with $ r(u)\not = r(v)$ whenever $\omega(u)\not=\omega(v) $ will do.
The single exception is the comparison with  Floer homology - and, in particular, to show that  $QH(-)$ vanishes for a displaceable Lagrangian.
For these additional properties to hold,  the ring $\mathcal{R}$ has to be stable with respect to the invertion of the elements in $\pi_{2}(M,L)$. For example, in the case described above, $r$ needs also to extend to a ring morphism $\Z_{2}[\pi_{2}(M,L)]\to \mathcal{R}$.

\subsection{Action estimates.} \label{subsec:action}
All the elements of the moduli spaces involved in our various
algebraic structures admit meaningful energy notions. It is therefore
easy (and essentially standard) to deduce various action estimates
from the non-triviality of these structures.

We shall give here a single example of such an application which is an
extension of the action estimates that appeared in the work of Albers
\cite{Alb:extrinisic} and apply to a version of the quantum inclusion
map $i_{L}$.

Let $\widetilde{ Ham (M)}$ be the universal cover of the group of
Hamiltonian diffeomorphisms and fix $\phi\in \widetilde{Ham (M)}$.
Recall that for any $\alpha\in H_{\ast}(M;\Z_2)$ there are spectral
invariants $\sigma(\alpha;\phi)\in\R$ and $\sigma(\alpha^{\ast};\phi)$
where $\alpha^{\ast}\in H^{\ast}(M;\Z_2)$ is the Poincar\'{e} dual of
$\alpha$. We refer the reader to~\cite{Sc:action-spectrum,
  Oh:spec-inv-1, Oh:spec-inv-2, Oh:spec-inv-3, Oh:spec-inv-4,
  Oh:spec-inv-5, McD-Sa:Jhol-2} for the foundations of the theory of
spectral invariants. We shall also recall the basic definitions in
\S\ref{subsec:action_proof}. We now define the \textit{depth} of
$\phi$ on $L$ by
$$\mathrm{depth}_{L}(\phi)=\sup_{[H]=\phi} \int_{0}^{1}(\inf_{x\in
  L}H(x,t))dt~.~$$
Similarly, we let the \textit{height} of $\phi$ on
$L$ be defined by:
$$\mathrm{height}_{L}(\phi)=\inf_{[H]=\phi}\int_{0}^{1}(\sup_{x\in
  L}H(x,t))dt~.~$$

\begin{cor} \label{cor:action}
   Assume that $\alpha\in H_{\ast}(M;\Z_2)$, $x,y\in Q^{+}H(L)$ are so
   that $y\not =0$ and $\alpha*x=yt^{k}+\mathrm{higher\ order\
     terms}$.  Then we have the following inequalities:

   $$\sigma(\alpha;\phi)- \mathrm{depth}_{L}(\phi)\geq - k\tau\leq
   \mathrm{height}_{L}(\phi)-\sigma(\alpha^{\ast};\phi) $$
   where
   $\tau$ is the monotonicity constant.

   As before, the same result continues to hold with $QH^{+}(L)$ replaced
   by $QH(L)=HF(L)$.
\end{cor}



%

\section{Transversality} \label{S:transversality}

The purpose of this section is to deal with the main transversality
issues that appear in the definition of our algebraic structures. The
pearl moduli spaces - they are at the heart of the definition of the
pearl complex - are introduced here and we shall see that
transversality is not difficult to achieve for them using the
structural results of Lazzarini~\cite{Laz:discs, Laz:decomp} combined
with some combinatorial arguments. The main ideas and technical lemmas
of this section will then be used for these and various other moduli
spaces of similar type in~\S\ref{S:proofalg}.

\subsection{Transversality for strings of pearls}
\label{Sb:tr-pearl-complex}

Let $(M^{2n}, \omega)$ be a tame symplectic manifold and $L^n \subset
(M^{2n}, \omega)$ a closed Lagrangian submanifold. Assume that $L$ is
monotone with minimal Maslov number $N_L \geq 2$. Denote by
$\mathcal{J}(M,\omega)$ the space of almost complex structures on $M$
which are compatible with $\omega$. Given a homology class $F \in
H_2(M,L;\mathbb{Z})$, denote by $\mathcal{M}(F,J)$ the space of
$J$-holomorphic disks $u:(D, \partial D) \to (M,L)$ in the class $F$.
(Here and in what follows $D \subset \mathbb{C}$ stands for the closed
unit disk.)

\begin{dfn} \label{D:simple-distinct}
   \begin{enumerate}
     \item A $J$-holomorphic disk $u:(D,\partial D) \to (M,L)$ is
      called simple if there exists an open dense subset $S \subset D$
      such that for every $z \in S$ we have $u^{-1}(u(z))=\{z\}$ and
      $du_{z} \neq 0$. We denote by $\mathcal{M}^*(F,J) \subset
      \mathcal{M}(F,J)$ the space of simple $J$-holomorphic disks
      $u:(D,\partial D) \to (M,L)$ in the class $F$.
     \item Let $v_i:(D, \partial D) \to (M,L)$, $i=1, \ldots, k$ be a
      sequence of $J$-holomorphic disks. We say that $(v_1, \ldots,
      v_k)$ are {\em absolutely distinct} if for every $1\leq i \leq
      k$ we have $v_i(D) \not \subset \bigcup_{j \neq i} v_j(D)$.
   \end{enumerate}
\end{dfn}

Let $f:L \to \mathbb{R}$ be a Morse function and $\rho$ a Riemannian
metric on $L$. We denote by $\Phi_t:L \to L$, $t \in \mathbb{R}$, the
{\em negative} gradient flow of $(f,\rho)$ (i.e.  the flow of the
vector field $-\textnormal{grad}_{\rho}f$).  Given critical points
$x,y \in \textnormal{Crit}(f)$ denote by $W_x^u$, $W_y^s$ the unstable
and stable submanifolds of $x$ and $y$ with respect to the negative
gradient flow of $f$.

Consider the (non-proper) embedding
\begin{equation} \label{Eq:Q-f}
   (L \setminus \textnormal{Crit}(f)) \times \mathbb{R}_{>0}
   \hooklongrightarrow L \times L, \quad (x,t) \longmapsto
   (x,\Phi_t(x)).
\end{equation}
Denote the image of this embedding by $Q_{f,\rho}
\subset L\times L$.

Denote by $G=\textnormal{Aut}(D) \cong PSL(2,\mathbb{R})$ the group of
biholomorphisms of $D$. Given points $p_1,\ldots, p_m \in D$ we denote
by $G_{p_1, \ldots, p_m} \subset G$ the subgroup of all the
automorphisms $\sigma \in G$ that fix all the $p_i$'s.

Let $\mathbf{A}=(A_1, \ldots, A_l)$ be a sequence of {\em non-zero}
homology classes $A_1, \ldots, A_l \in H_2(M,L;\mathbb{Z})$, $l \geq
1$. Set $\mu(\mathbf{A}) = \sum_{i=1}^l \mu(A_i)$. Put:
\begin{equation} \label{Eq:MA}
   \mathcal{M}(\mathbf{A},J) = \mathcal{M}(A_1,J)/G_{-1,1} \times
   \cdots \times \mathcal{M}(A_l,J)/G_{-1,1}.
\end{equation}
Denote by $\mathcal{M}^{*,\textnormal{d}}(\mathbf{A},J)$ the subspace
of all $(u_1, \ldots, u_l) \in \mathcal{M}(\mathbf{A},J)$ which are
simple and absolutely distinct. Consider the following evaluation map:
\begin{equation} \label{Eq:ev-1}
   ev_{\mathbf{A}}: \mathcal{M}(\mathbf{A},J)
   \longrightarrow L^{\times 2l}, \quad
   ev_{\mathbf{A}}(u_1, \ldots, u_l) = \bigl(u_1(-1),u_1(1), \ldots,
   u_l(-1),u_l(1)\bigr).
\end{equation}
For every $x,y \in \textnormal{Crit}(f)$, put
$$\mathcal{P}(x,y,\mathbf{A}; J,f,\rho) =
ev_{\mathbf{A}}^{-1}\Bigl(W_x^u \times \bigl(Q_{f,\rho}\bigr)^{\times
  (l-1)} \times W_y^s \Bigr).$$

\begin{figure}[htbp]
   \begin{center}
      \epsfig{file=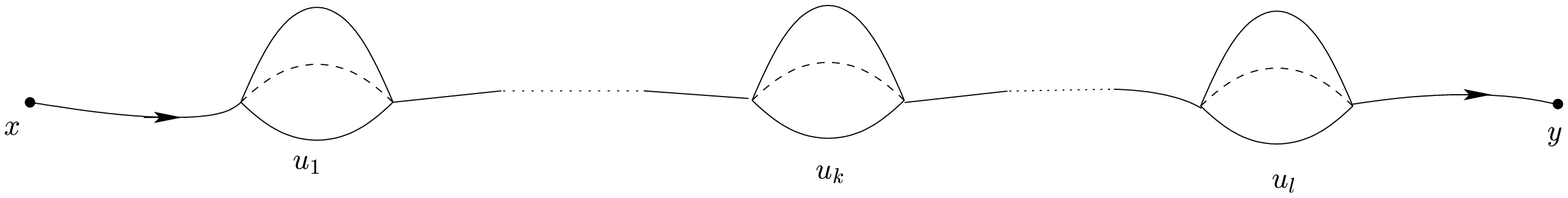, width=0.6\linewidth}
   \end{center}
   \caption{An element of $\mathcal{P}(x,y;\mathbf{A},J)$}
   \label{f:qm-pearls-1}
\end{figure}

We call $\mathcal{P}(- - -)$ the ``moduli space of pearls". They
consist of objects as in figure~\ref{f:qm-pearls-1}.

Finally, write $\mathcal{P}^{*,\textnormal{d}}(x,y,\mathbf{A};
J,f,\rho)= \mathcal{P}(x,y,\mathbf{A}; J,f,\rho) \cap
\mathcal{M}^{*,\textnormal{d}}(\mathbf{A},J)$, namely the subspace of
all $(u_1, \ldots, u_l) \in \mathcal{P}(x,y,\mathbf{A}; J,f,\rho)$
which are simple and absolutely distinct. Standard
arguments~\cite{McD-Sa:Jhol-2} show that:
\begin{statement} \label{St:tr-simple-dist-1}
   For every pair $(f,\rho)$ there exists a subset
   $\mathcal{J}_{\textnormal{reg}} \subset \mathcal{J}(M,\omega)$ of
   second category such that for every $J \in
   \mathcal{J}_{\textnormal{reg}}$, for every sequence of non-zero
   homology classes $\mathbf{A}$ and every $x,y \in
   \textnormal{Crit}(f)$, the restriction of $ev_{\mathbf{A}}$ to
   $\mathcal{M}^{*,\textnormal{d}}(\mathbf{A},J)$ is transverse to
   $W_x^u \times \bigl(Q_{f,\rho}\bigr)^{\times (l-1)} \times W_y^s$
   at every $\mathbf{u} \in
   \mathcal{M}^{*,\textnormal{d}}(\mathbf{A},J)$. In particular the
   space $\mathcal{P}^{*,\textnormal{d}}(x,y,\mathbf{A}; J,f,\rho)$ is
   either empty or a smooth manifold of dimension
   $\mu(\mathbf{A})+\ind_f(x) - \ind_f(y) - 1$.
\end{statement}

In this section we prove the following.
\begin{prop} \label{P:tr-1}
   Let $f:L \to \mathbb{R}$ be a Morse function and $\rho$ a
   Riemannian metric on $L$ such that the pair $(f, \rho)$ is
   Morse-Smale. Then there exists a subset of second category
   $\mathcal{J}_{\textnormal{reg}} \subset \mathcal{J}(M,\omega)$ with
   the following property. For every sequence of non-zero homology
   classes $\mathbf{A}=(A_1, \ldots, A_l)$ and every $x,y \in
   \textnormal{Crit}(f)$ with $\mu(\mathbf{A})+\ind_f(x) - \ind_f(y)-1
   \leq 1$:
   \begin{enumerate}
     \item $\mathcal{P}(x,y,\mathbf{A};
      J,f,\rho)=\mathcal{P}^{*,\textnormal{d}}(x,y,\mathbf{A};
      J,f,\rho)$. In other words all elements $(u_1, \ldots, u_l) \in
      \mathcal{P}(x,y,\mathbf{A}; J,f,\rho)$ are simple and absolutely
      distinct. Thus $\mathcal{P}(x,y,\mathbf{A}; J,f,\rho)$ is either
      empty or a smooth manifold of dimension
      $\mu(\mathbf{A})+\ind_f(x) - \ind_f(y)-1$. In particular, if
      $\mu(\mathbf{A})+\ind_f(x) - \ind_f(y)-1 <0$ we have
      $\mathcal{P}(x,y,\mathbf{A}; J,f,\rho)=\emptyset$.
     \item If $\mu(\mathbf{A})+\ind_f(x) - \ind_f(y)-1=0$ then
      $\mathcal{P}(x,y,\mathbf{A}; J,f,\rho)$ is a {\em compact}
      $0$-dimensional manifold hence consists of finite number of
      points.
   \end{enumerate}
\end{prop}
The proof of Proposition~\ref{P:tr-1} is given in
Section~\ref{Sb:prf-tr-1} below.

\begin{rem} \label{R:2nd-categ}
   Since a countable intersection of second category subset of
   $\mathcal{J}(M,\omega)$ is of second category too we shall denote
   various second category subsets of almost complex structure in this
   section by the same notation $\mathcal{J}_{\textnormal{reg}}$.
\end{rem}

Let $\mathbf{B}'=(B'_1, \ldots, B'_{l'})$, $\mathbf{B}''=(B''_1,
\ldots, B''_{l''})$ be two vectors of non-zero homology classes in
$H_2(M,L;\mathbb{Z})$. Put
\begin{equation*}
   \mathcal{M}(\mathbf{B}',\mathbf{B''},J) =
   \prod_{i=1}^{l'} \frac{\mathcal{M}(B'_i,J)}{G_{-1,1}} \times
   \prod_{j=1}^{l''}
   \frac{\mathcal{M}(B''_j,J)}{G_{-1,1}}.
\end{equation*}
Denote by $\mathcal{M}^{*,\textnormal{d}}(\mathbf{B}',\mathbf{B''},J)$
the subspace of all $(\mathbf{u}',\mathbf{u}'') \in
\mathcal{M}(\mathbf{B}',\mathbf{B''},J)$ for which the $J$-holomorphic
disks $(u'_1, \ldots, u'_{l'}, u''_1, \ldots, u''_{l''})$ are simple
and absolutely distinct.  Define an evaluation map
$ev_{\mathbf{B}',\mathbf{B}''}:
\mathcal{M}(\mathbf{B}',\mathbf{B''},J) \to L^{\times (2l'+2l'')}$ by:
\begin{equation*}
   ev_{\mathbf{B}',\mathbf{B}''}(u'_1, \ldots, u'_{l'},
   u''_1, \ldots, u''_{l''}) =
   \bigl(ev_{\mathbf{B}'}(u'_1, \ldots, u'_{l'}), \,
      ev_{\mathbf{B}''}(u''_1, \ldots, u''_{l''})\bigr),
\end{equation*}
where $ev_{\mathbf{B}'}$, $ev_{\mathbf{B}''}$ are defined as
in~\eqref{Eq:ev-1}. Put
$$\mathcal{P}(x,y, \mathbf{B}',\mathbf{B}''; J,f,\rho) =
ev_{\mathbf{B}', \mathbf{B}''}^{-1} \Bigl( W_x^u \times
\bigl(Q_{f,\rho}\bigr)^{\times (l'-1)} \times \textnormal{diag}(L)
\times \bigl(Q_{f,\rho}\bigr)^{\times (l''-1)} \times W_y^s \Bigr).$$
Finally, write $\mathcal{P}^{*,\textnormal{d}}(x,y,
\mathbf{B}',\mathbf{B}''; J,f,\rho) = \mathcal{P}(x,y,
\mathbf{B}',\mathbf{B}''; J,f,\rho) \cap
\mathcal{M}^{*,\textnormal{d}}(\mathbf{B}',\mathbf{B''},J)$.  Standrad
arguments~\cite{McD-Sa:Jhol-2} show that:
\begin{statement} \label{St:tr-simple-dist-2}
   For every pair $(f,\rho)$ there exists a subset
   $\mathcal{J}_{\textnormal{reg}} \subset \mathcal{J}(M,\omega)$ of
   second category such that for every $J \in
   \mathcal{J}_{\textnormal{reg}}$, for every two sequences of
   non-zero homology classes $\mathbf{B}'$, $\mathbf{B}''$ and every
   $x,y \in \textnormal{Crit}(f)$ the restriction of the map
   $ev_{\mathbf{B}',\mathbf{B}''}$ to
   $\mathcal{M}^{*,\textnormal{d}}(\mathbf{B}',\mathbf{B''},J)$ is
   transverse to $$W_x^u \times \bigl(Q_{f,\rho}\bigr)^{\times (l'-1)}
   \times \textnormal{diag}(L) \times \bigl(Q_{f,\rho}\bigr)^{\times
     (l''-1)} \times W_y^s.$$
   In particular, the space
   $\mathcal{P}^{*,\textnormal{d}}(x,y,\mathbf{B}',\mathbf{B}'';
   J,f,\rho)$ is either empty or a smooth manifold of dimension
   $\mu(\mathbf{B}')+\mu(\mathbf{B}'')+\ind_f(x) - \ind_f(y)-2$.
\end{statement}

With the above notation we have
\begin{prop} \label{P:tr-2}
   Let $f:L \to \mathbb{R}$ be a Morse function and $\rho$ a
   Riemannian metric on $L$ such that the pair $(f, \rho)$ is
   Morse-Smale. Then there exists a subset of second category
   $\mathcal{J}_{\textnormal{reg}} \subset \mathcal{J}(M,\omega)$ with
   the following property. For every two sequences of non-zero homology
   classes $\mathbf{B}'=(B'_1, \ldots, B'_{l'})$,
   $\mathbf{B}''=(B''_1, \ldots, B''_{l''})$ and every $x,y \in
   \textnormal{Crit}(f)$ with
   $\mu(\mathbf{B}')+\mu(\mathbf{B}'')+\ind_f(x) - \ind_f(y)-1 \leq 1$:
   \begin{enumerate}
     \item $\mathcal{P}(x,y,\mathbf{B}',\mathbf{B}''; J,f,\rho) =
      \mathcal{P}^{*,\textnormal{d}}(x,y,\mathbf{B}',\mathbf{B}'';
      J,f,\rho)$. In other words for every $$(u'_1, \ldots, u'_{l'},
      u''_1, \ldots, u''_{l''}) \in
      \mathcal{M}(x,y,\mathbf{B}',\mathbf{B}''; J,f,\rho)$$
      the disks
      $(u'_1, \ldots, u'_{l'}, u''_1, \ldots, u''_{l''})$ are simple
      and absolutely distinct.
     \item If $\mu(\mathbf{B}')+\mu(\mathbf{B}'')+\ind_f(x) -
      \ind_f(y)-1 \leq 0$ then
      $\mathcal{P}(x,y,\mathbf{B}',\mathbf{B}''; J,f,\rho)=\emptyset$.
     \item If $\mu(\mathbf{B}')+\mu(\mathbf{B}'')+\ind_f(x) -
      \ind_f(y)-1 = 1$ then $\mathcal{P}(x,y,\mathbf{B}',\mathbf{B}'';
      J,f,\rho)$ is either empty or a compact $0$-dimensional smooth
      manifold, hence consists of finite number of points.
   \end{enumerate}
\end{prop}
We shall not give a proof for Proposition~\ref{P:tr-2} since it can be
proved in a very similar way to Proposition~\ref{P:tr-1} which will be
proved below.

\subsection{Reduction to simple disks}
\label{Sb:reduction}
Let $F \in H_2(M,L; \mathbb{Z})$, $J \in \mathcal{J}(M,\omega)$.
Denote by $\mathcal{M}^*(F,J)$ the space of simple $J$-holomorphic
disks $u:(D,\partial D) \to (M,L)$ in the class $F$.  According to the
general theory of pseudo-holomorphic curves~\cite{McD-Sa:Jhol-2} for a
generic choice of $J \in \mathcal{J}(M,\omega)$ the space
$\mathcal{M}^*(F,J)$ is either empty or a smooth manifold of dimension
$\mu(F)+n$. This fails to be true for the space $\mathcal{M}(F,J)$
of all disks in the class $F$.  Therefore a crucial ingredient in the
proof of Proposition~\ref{P:tr-1} is a procedure which enables to
decompose a (general) $J$-holomorphic disk to simple ``pieces''. This
will make it possible to obtain transversality and to control
dimensions of moduli spaces of pseudo-holomorphic disks. There are two
(essentially equivalent) approaches to this decomposition, one due to
Kwon and Oh~\cite{Kw-Oh:discs} and the other to
Lazzarini~\cite{Laz:discs, Laz:decomp}. Below we shall follow
Lazzarini's approach.

Let $u:(D,\partial D) \to (M,L)$ be a non-constant $J$-holomorphic
disk. Put $\mathcal{C}(u)=u^{-1}(\{ du=0 \})$. Define a relation
$\mathcal{R}_u$ on pairs of points $z_1, z_2 \in \textnormal{Int\,}D
\setminus \mathcal{C}(u)$ in the following way:
\begin{equation*}
   z_1 \mathcal{R}_u z_2 \Longleftrightarrow
   \begin{cases}
      & \forall \textnormal{ neighbourhoods } V_1, V_2 \textnormal{ of
      } z_1,z_2, \\
      & \exists \textnormal{ neighbourhoods } U_1, U_2
      \textnormal{ such that:} \\
      & \textnormal{(i) } z_1 \in U_1 \subset V_1 ,
      z_2 \in U_2 \subset V_2. \\
      & \textnormal{(ii) } u(U_1)=u(U_2).
   \end{cases}
\end{equation*}
Denote by $\overline{\mathcal{R}}_u$ the closure of $\mathcal{R}_u$ in
$D \times D$. Note that $\overline{\mathcal{R}}_u$ is reflexive and
symmetric but it may fail to be transitive (see~\cite{Laz:decomp} for
more details on this). Define the {\em non-injectivity graph} of $u$
to be:
$$\mathcal{G}(u) = \{ z \in D \, | \, \exists z' \in \partial D
\textnormal{ such that } z \overline{\mathcal{R}}_u z' \}.$$
It is
proved in~\cite{Laz:decomp,Laz:discs} that $\mathcal{G}(u)$ is indeed
a graph and its complement $D \setminus \mathcal{G}(u)$ has finitely
many connected components. In what follows we shall use the following
theorem due to Lazzarini (See Proposition~4.1 in~\cite{Laz:decomp}.
See also~\cite{Laz:discs}).
\begin{thm}[Decomposition of disks]
   \label{T:decomposition}
   Let $u:(D,\partial D) \to (M,L)$ be a non-constant $J$-holomorphic
   disk. Then for every connected component $\mathfrak{D} \subset D
   \setminus \mathcal{G}(u)$ there exists a surjective map
   $\pi_{\overline{\mathfrak{D}}}: \overline{\mathfrak{D}} \to D$,
   holomorphic on $\mathfrak{D}$ and continuous on
   $\overline{\mathfrak{D}}$, and a simple $J$-holomorphic disk
   $v_{\mathfrak{D}}:(D,\partial D) \to (M,L)$ such that
   $u|_{\overline{\mathfrak{D}}} = \pi_{\overline{\mathfrak{D}}} \circ
   v_{\mathfrak{D}}$. The map
   $\pi_{\overline{\mathfrak{D}}}:\overline{\mathfrak{D}} \to D$ has a
   well defined degree $m_{\mathfrak{D}} \in \mathbb{N}$ and we have
   in $H_2(M,L;\mathbb{Z})$:
   $$[u]=\sum_{\mathfrak{D}} m_{\mathfrak{D}} [v_{\mathfrak{D}}],$$
   where the sum is taken over all connected components $\mathfrak{D}
   \subset D \setminus \mathcal{G}(u)$.
\end{thm}
\begin{remnonum}
   Some of the connected components $\mathfrak{D} \subset D \setminus
   \mathcal{G}(u)$ may not be disks. This happens if and only if
   $\mathcal{G}(u)$ is not connected. Nevertheless,
   $\overline{\mathfrak{D}} / \overline{\mathcal{R}}_{\mathfrak{D}}$
   is a disk, where $\overline{\mathcal{R}}_{\mathfrak{D}}$ is the
   relation defined similarly to $\overline{\mathcal{R}}_u$ but for
   pairs of points in $\overline{\mathfrak{D}}$.
\end{remnonum}

For the proof of Proposition~\ref{P:tr-1} we shall need the following
two lemmas. We state them here but defer their proofs to
Section~\ref{Sb:prf-L-non-simple}.
\begin{lem} \label{L:uv}
   Suppose $n=\dim L \geq 3$. Then there exists a second category
   subset $\mathcal{J}_{\textnormal{reg}} \subset
   \mathcal{J}(M,\omega)$ such that for every $J \in
   \mathcal{J}_{\textnormal{reg}}$ the following holds. If
   $u,v:(D,\partial D) \to (M,L)$ are simple $J$-holomorphic disks
   such that $u(D) \cap v(D)$ is an infinite set then:
   \begin{itemize}
     \item either $u(D) \subset v(D)$ and $u(\partial D) \subset
      v(\partial D)$; or
     \item $v(D) \subset u(D)$ and $v(\partial D) \subset u(\partial
      D)$.
   \end{itemize}
\end{lem}

\begin{lem} \label{L:non-simple}
   Suppose $n=\dim L \geq 3$. Then there exists a second category
   subset $\mathcal{J}_{\textnormal{reg}} \subset
   \mathcal{J}(M,\omega)$ such that for every $J \in
   \mathcal{J}_{\textnormal{reg}}$ the following holds. For every
   non-simple $J$-holomorphic disk $u:(D, \partial D) \to (M,L)$ with
   $u(-1) \neq u(1)$ there exists a $J$-holomorphic disk $u':(D,
   \partial D) \to (M,L)$ with the following properties:
   \begin{enumerate}
     \item $u'(-1)=u(-1)$, $u'(1)=u(1)$.
     \item $u'(D) = u(D)$ and $u'(\partial D) = u(\partial D)$.
     \item $u'$ is simple.
     \item $\omega([u']) < \omega([u])$. In particular, if $L$ is
      monotone we also have $\mu([u']) < \mu([u])$.
   \end{enumerate}
\end{lem}

\begin{rem} \label{R:not-monotone}
   None of the Lemmas~\ref{L:uv} and~\ref{L:non-simple} require $L$ to
   be monotone.
\end{rem}

\subsection{Proof of Proposition~\ref{P:tr-1}} \label{Sb:prf-tr-1}
We separate the proof of Proposition~\ref{P:tr-1} into two cases: $n =
\dim L \geq 3$ and $n = \dim L \leq 2$. We start with $n\geq 3$.
\begin{proof}[Proof of Proposition~\ref{P:tr-1} for $n \geq 3$]
   We start with statement~(1) of the Proposition. Let
   $\mathcal{J}_{\textnormal{reg}} \subset \mathcal{J}(M,\omega)$ be
   the intersection of the sets given by Lemma~\ref{L:non-simple} and
   by Statement~\ref{St:tr-simple-dist-1}. The proof is carried out by
   induction over the integer $\mu(\mathbf{A})/N_L$.

   Suppose $\mu(\mathbf{A})=N_L$. Let $\mathbf{u}=(u_1, \ldots, u_l)
   \in \mathcal{P}(x,y,\mathbf{A}; J,f,\rho)$. As $L$ is monotone we
   have $l=1$ (i.e. $\mathbf{u}$ consists of just one disk $u_1$).
   Since $\mu([u_1])=N_L$ it follows from
   Theorem~\ref{T:decomposition} that $u_1$ is simple, hence
   $\mathbf{u} \in \mathcal{P}^{*,\textnormal{d}}(x,y,\mathbf{A};
   J,f,\rho)$.

   Suppose that statement~(1) of our Proposition holds for every
   $\mathbf{A}$ with $\mu(\mathbf{A}) \leq kN_L$. Let
   $\mathbf{A}=(A_1, \ldots, A_l)$ be a sequence of non-zero homology
   classes with $\mu(\mathbf{A})=(k+1)N_L$, and such that
   $\mu(\mathbf{A})+\ind_f(x)-\ind_f(y)-1 \leq 1$. Let $\mathbf{u} \in
   \mathcal{P}(x,y,\mathbf{A}; J,f,\rho)$ and suppose by contradiction
   that $\mathbf{u} \notin
   \mathcal{P}^{*,\textnormal{d}}(x,y,\mathbf{A}; J,f,\rho)$.

   First note that for every $1 \leq i \leq l$ we have $u_i(-1) \neq
   u_i(1)$. Indeed if $u_j(-1)=u_j(1)$ for some $j$ then let
   $\mathbf{u}'$ be the sequence of disks obtained from $\mathbf{u}$
   by omitting $u_j$. Let $\mathbf{A}'$ be obtained from $\mathbf{A}$
   by omitting $A_j$. As $u_j(-1)=u_j(1)$ we have $\mathbf{u}' \in
   \mathcal{P}(x,y,\mathbf{A}'; J,f,\rho)$. But $\mu(\mathbf{A}') \leq
   \mu(\mathbf{A})-N_L$ hence by the induction hypothesis we have
   $\mathbf{u}' \in \mathcal{P}^{*,\textnormal{d}}(x,y,\mathbf{A}';
   J,f,\rho)$.  This leads to contradiction since $$\dim
   \mathcal{P}^{*,\textnormal{d}}(x,y,\mathbf{A}'; J,f,\rho) =
   \mu(\mathbf{A}) - \mu(A_j) + \ind_f(x) - \ind_f(y) -1 \leq 1 - N_L
   \leq -1.$$
   We assume from now on that $u_i(-1) \neq u_i(1)$ for
   every $i$.

   \smallskip \noindent \textbf{Case 1. There exists $1 \leq i_0 \leq
     l$ such that $u_{i_0}$ is not simple.}

   Apply Lemma~\ref{L:non-simple} with $u=u_{i_0}$ to obtain a
   (simple) disk $u'_{i_0}$ with $u'_{i_0}(-1)=u_{i_0}(-1)$,
   $u'_{i_0}(1)=u_{i_0}(1)$ and such that $\mu([u'_{i_0}]) <
   \mu([u_{i_0}])$. Let $\mathbf{u}'$ be the sequence of disks
   obtained from $\mathbf{u}$ by replacing $u_{i_0}$ by $u'_{i_0}$.
   Let $\mathbf{A}'$ be obtained from $\mathbf{A}$ by replacing
   $A_{i_0}$ by $[u'_{i_0}]$. Clearly $\mathbf{u}' \in
   \mathcal{P}(x,y,\mathbf{A}'; J,f,\rho)$. By the induction
   hypothesis we have $\mathbf{u}' \in
   \mathcal{P}^{*,\textnormal{d}}(x,y,\mathbf{A}'; J,f,\rho)$. But
   this leads to contradiction since
   \begin{equation} \label{Eq:dim-M*}
      \dim \mathcal{P}^{*,\textnormal{d}}(x,y,\mathbf{A}'; J,f,\rho) =
      \mu(\mathbf{A}') + \ind_f(x) - \ind_f(y) -1 \leq \mu(\mathbf{A})-N_L +
      \ind_f(x)-\ind_f(y) -1 \leq -1.
   \end{equation}

   \smallskip \noindent \textbf{Case 2. The disks $(u_1, \ldots, u_l)$
     are simple but not absolutely distinct.}

   In this case there exists $i_0$ such that $u_{i_0}(D) \subset
   \cup_{i\neq i_0} u_i(D)$. It follows that there exists $j_0$ such
   that $u_{i_0}(D) \cap u_{j_0}(D)$ is an infinite set. By
   Lemma~\ref{L:uv} we have:
   \begin{itemize}
     \item either $u_{i_0}(D) \subset u_{j_0}(D)$ and
      $u_{i_0}(\partial D) \subset u_{j_0}(\partial D)$; or
     \item $u_{j_0}(D) \subset u_{i_0}(D)$ and $u_{j_0}(\partial D)
      \subset u_{i_0}(\partial D)$.
   \end{itemize}
   Without loss of generality assume that the first possibility
   occurs.

   \smallskip \noindent \textbf{\em Subcase i. $i_0 < j_0$.}

   Denote by $\mathbf{u}'$ the sequence of disks obtained from
   $\mathbf{u}$ by omitting all the disks $u_{i_0}, \ldots,
   u_{j_0-1}$.  Denote by $\mathbf{A}'$ the corresponding vector of
   homology classes. There exists a point $p \in \partial D$ such that
   $u_{j_0}(p)= u_{i_0}(-1)$.

   In case $p \neq 1$ we can replace $u_{j_0}$ by $u_{j_0} \circ
   \sigma$, where $\sigma \in \textnormal{Aut}(D)$ is such that
   $\sigma(1)=1$ and $\sigma(-1)=p$.  Note that now $\mathbf{u}' \in
   \mathcal{P}(x,y,\mathbf{A}'; J,f,\rho)$.

   In case $p=1$ omit from $\mathbf{u}'$ also the disk $u_{j_0}$. If
   the resulting sequence of disks $\mathbf{u}'$ is empty we obtain a
   trajectory (of $-\textnormal{grad}_{\rho}f$) connecting $x$ to $y$.
   But this is impossible since $\ind_f(x) - \ind_f(y) \leq 1-
   \mu(\mathbf{A}) \leq -1$ and $(f,\rho)$ is Morse-Smale. Thus we may
   assume that $\mathbf{u}'$ is not empty and we have $\mathbf{u}' \in
   \mathcal{P}(x,y,\mathbf{A}'; J,f,\rho)$.

   Summing up, in both cases, $p=1$ and $p \neq 1$, we have
   $\mathbf{u}' \in \mathcal{P}(x,y,\mathbf{A}'; J,f,\rho)$ and
   $\mu(A') < \mu(A)$. The induction hypothesis implies that
   $\mathbf{u}' \in \mathcal{P}^{*,\textnormal{d}}(x,y,\mathbf{A}';
   J,f,\rho)$. We now obtain contradiction in the same way as in
   inequality~\eqref{Eq:dim-M*} above.

   \smallskip \noindent \textbf{\em Subcase ii. $i_0 > j_0$.}

   We argue similarly to Subcase~{\em i} only that now we omit from
   $\mathbf{u}$ the disks $u_{j_0+1}, \ldots, u_{i_0}$.

   This completes the proof of statement~(1) of
   Proposition~\ref{P:tr-1} in the case $n\geq 3$.

   The proof of statement~(2) of Proposition~\ref{P:tr-1} is based on
   similar arguments to the above. Note however that we shall need to
   reduce further the space $\mathcal{J}_{\textnormal{reg}}$ (e.g by
   intersecting it with the subset coming from
   statement~\ref{St:tr-simple-dist-2}).
\end{proof}

\begin{proof}[Proof of Proposition~\ref{P:tr-1} for $n \leq 2$]
   Again, we prove only statement~(1).

   Denote by $\mathcal{J}'$ be the set of all $J \in \mathcal{J}(M,
   \omega)$ for which the following holds: for every class $A \in
   H_2(M, L;\mathbb{Z})$ with $\mu(A)=2$ and every $x,y \in
   \textnormal{Crit}(f)$ the evaluation maps
   \begin{align} \label{Eq:ev-n=2} & ev'_A: \Bigl( \mathcal{M}^*(A,J)
      \times \textnormal{Int\,} D \Bigr) /G_1 \longrightarrow M \times
      L, \quad & &ev'_A(u,p) = (u(p),
      u(1)), \\
      & ev''_A: \Bigl( \mathcal{M}^*(A,J) \times \textnormal{Int\,} D
      \Bigr)/G_{-1} \longrightarrow M \times L, \quad & &ev''_A(u,p) =
      (u(-1), u(p)), \notag
   \end{align}
   are transverse to $W_x^u \times W_y^s$. Standard
   arguments~\cite{McD-Sa:Jhol-2} show that $\mathcal{J}' \subset
   \mathcal{J}(M, \omega)$ is of second category. We define the set
   $\mathcal{J}_{\textnormal{reg}} \subset \mathcal{J}(M, \omega)$ to
   be the intersection of $\mathcal{J}'$ with the sets given by
   Statements~\ref{St:tr-simple-dist-1},~\ref{St:tr-simple-dist-2}.

   Let $\mathbf{A}=(A_1, \ldots, A_l)$ be a sequence of non-zero
   classes, and $x, y \in \textnormal{Crit}(f)$ with $\mu(\mathbf{A})
   + \ind_f(x)-\ind_f(y)-1 \leq 1$. Since $n\leq 2$ we have
   $\mu(\mathbf{A}) \leq 4$.

   Suppose first that $\mu(\mathbf{A}) \leq 3$. Since $N_L \geq 2$,
   for every $(u_1, \ldots, u_l) \in \mathcal{P}(x,y,\mathbf{A};
   J,f,\rho)$ we must have $l=1$.  By Theorem~\ref{T:decomposition}
   $u_1$ is simple.  This completes the proof in the case
   $\mu(\mathbf{A}) \leq 3$.

   Suppose now that $\mu(\mathbf{A}) = 4$. Note that in this case we
   must have $n=2$.  Let $\mathbf{u} \in \mathcal{P}(x,y,\mathbf{A};
   J,f,\rho)$ and assume by contradiction that $\mathbf{u} \notin
   \mathcal{P}^{*,\textnormal{d}}(x,y,\mathbf{A}; J,f,\rho)$. By
   monotonicity of $L$ we either have $l=2$, $\mu(A_1)=\mu(A_2)=2$ or
   $l=1$, $\mu(A_1)=4$. Also, by a similar argument to the ones at the
   beginning of the proof for the case $n\geq 3$ we may assume that
   $u_i(-1) \neq u_i(1)$ for every $i$.

   \smallskip \noindent \textbf{The case $l=2$,
     $\mu(A_1)=\mu(A_2)=2$.}

   Since $\mu(u_1)=\mu(u_2)=2$, both $u_1$ and $u_2$ are simple.
   Thus $u_1, u_2$ are not absolutely distinct. Without loss of
   generality assume that $u_1(D) \subset u_2(D)$.

   Suppose first that $u_1(-1) \in u_2(\textnormal{Int\,} D)$.  Let $p
   \in \textnormal{Int\,}D$ such that $u_2(p)=u_1(-1)$.  Then $(u_2,p)
   \in (ev'_{A_2})^{-1}(W_x^u \times W_y^s)$, where $ev'_{A_2}$ is the
   evaluation map defined in~\eqref{Eq:ev-n=2}. Since $ev'_{A_2}$ is
   transverse to $W_x^u \times W_y^s$ a simple computation shows that
   $$\dim (ev'_{A_2})^{-1}(W_x^u \times W_y^s) =
   \mu(A_2)-n+\ind_f(x)-\ind_f(y) = \ind_f(x)-\ind_f(y) \leq 2-
   \mu(\mathbf{A}) = -2,$$
   a contradiction.

   Suppose now that $u_1(-1) \in u_2(\partial D)$. If $u_1(-1) \neq
   u_2(1)$ then after a suitable reparametrization of $u_2$ we may
   assume that $u_2 \in \mathcal{P}^{*,\textnormal{d}}(x,y,A_2;
   J,f,\rho)$ which is impossible since
   $$\dim \mathcal{P}^{*,\textnormal{d}}(x,y,A_2; J,f,\rho) = \mu(A_2)
   + \ind_f(x)-\ind_f(y) -1 = \mu(\mathbf{A}) + \ind_f(x)-\ind_f(y) -1
   -\mu(A_1) \leq -1.$$
   The remaining case to consider is $u_1(-1) =
   u_2(1)$. In this case we can omit both $u_1$ and $u_2$ and obtain a
   trajectory of $-\textnormal{grad}_{\rho}f$ going form $x$ to $y$.
   But this is impossible since $\ind_f(x) < \ind_f(y)$ and $(f,\rho)$
   is Morse-Smale.

   \smallskip \noindent \textbf{The case $l=1$, $\mu(A_1)=4$.}

   In this case $u_1$ is not simple. Let
   $\mathcal{G}=\mathcal{G}(u_1)$ be the non-injectivity graph of
   $u_1$. Since $\mu([u_1])=4$, $D \setminus \mathcal{G}$ may have
   at most two connected components.

   \smallskip \noindent \textbf{\em Subcase i. $D \setminus
     \mathcal{G}$ is connected.}

   By Theorem~\ref{T:decomposition}, $u_1$ factors through a simple
   $J$-holomorphic disk $v:(D, \partial D) \to (M,L)$ via a
   holomorphic map $\pi:D \to D$ of degree $\geq 2$. (In fact the
   degree is exactly $2$ here).  It follows that $\mu([v])=2$.
   Since $u_1(-1) \neq u_2(1)$ there exists two {\em distinct} points
   $p', p'' \in \partial D$ such that $v(p')=u_1(-1)$,
   $v(p'')=u_1(1)$. After a suitable reparametrization of $v$ we may
   assume that $p'=-1$, $p''=1$ and we have $v \in
   \mathcal{P}^{*,\textnormal{d}}(x,y,[v]; J,f,\rho)$. But this leads
   to contradiction since
   \begin{equation} \label{Eq:dim-M*-n=2}
      \dim \mathcal{P}^{*,\textnormal{d}}(x,y,[v]; J,f,\rho) =
      \mu([v]) + \ind_f(x)-\ind_f(y)-1 \leq \mu(A_1)-2 +
      \ind_f(x)-\ind_f(y) -1 \leq -1.
   \end{equation}

   It remains to deal with the case that $D \setminus \mathcal{G}$ has
   two connected components $\mathfrak{D}_1, \mathfrak{D}_2$. Let
   $\pi_i=\pi_{\overline{\mathfrak{D}}_i}$, $v_{\mathfrak{D}_i}$,
   $m_{\mathfrak{D}_j}$, $i=1,2$ be the maps and multiplicities given
   by Theorem~\ref{T:decomposition}. (In fact, since $\mu(A_1)=4$ we
   must have $m_1=m_2=1$ and $\mu([v_1])=\mu([v_2])=2$.)

   \smallskip \noindent \textbf{\em Subcase ii. $-1, 1 \in
     \overline{\mathfrak{D}}_1$ (see the left part of
     figure~\ref{f:subcases-ii-iii-n=2}).}

   There exists two {\em distinct} points $p', p'' \in \partial D$
   such that $v_1(p')=u_1(-1)$, $v_1(p'')=u_1(1)$.  After a suitable
   reparametrization of $v$ we may assume that $p'=-1$, $p''=1$ and we
   have $v_1 \in \mathcal{P}^{*,\textnormal{d}}(x,y,[v_1]; J,f,\rho)$.
   We now obtain contradiction by a dimension count similar
   to~\eqref{Eq:dim-M*-n=2}.

   \smallskip \noindent \textbf{\em Subcase iii. $-1 \in
     \overline{\mathfrak{D}}_1$, $1 \in \overline{\mathfrak{D}}_2$
     (see the right part of figure~\ref{f:subcases-ii-iii-n=2}).}

   Put $B_1=[v_1]$, $B_2=[v_2]$.  In case $v_1(D) \subset v_2(D)$ or
   $v_2(D) \subset v_1(D)$ we can argue in a similar way to {\em ``The
     case $l=2$, $\mu(A_1)=\mu(A_2)=2$''} above and arrive to
   contradiction.

   It remains to deal with the case that $v_1, v_2$ are absolutely
   distinct. Put $B_1=[v_1]$, $B_2=[v_2]$. After suitable
   reparametrizations of $v_1,v_2$ we may assume that
   $v_1(-1)=u_1(-1)$ and $v_2(1)=u_1(1)$. Since
   $\overline{\mathfrak{D}}_1 \cap \overline{\mathfrak{D}}_2$ must
   contain at a $1$-dimensional component there exists two arcs
   $\gamma_1, \gamma_2 \subset \partial D$ such that for every $p_1
   \in \gamma_1$, $p_2 \in \gamma_2$ we have $v_1(p_1)=v_2(p_2)$. It
   follows that $\{ (v_1, p_1, p_2, v_2)\}_{p_1\in \gamma_1, p_2 \in
     \gamma_2}$ lies in $\dim
   \mathcal{P}^{*,\textnormal{d}}(x,y,B_1,B_2; J,f,\rho)$ hence the
   latter space is at least $1$-dimensional. But this is impossible
   since
   $$\dim \mathcal{P}^{*,\textnormal{d}}(x,y,B_1,B_2; J,f,\rho) =
   \mu(B_1)+\mu(B_2)+\ind_f(x) - \ind_f(y)-2 \leq 0.$$

   This completes the proof of Statement~(1) of
   Proposition~\ref{P:tr-1} for the case $n\leq 2$.

   \begin{figure}[htbp]
      \begin{center}
         \psfig{file=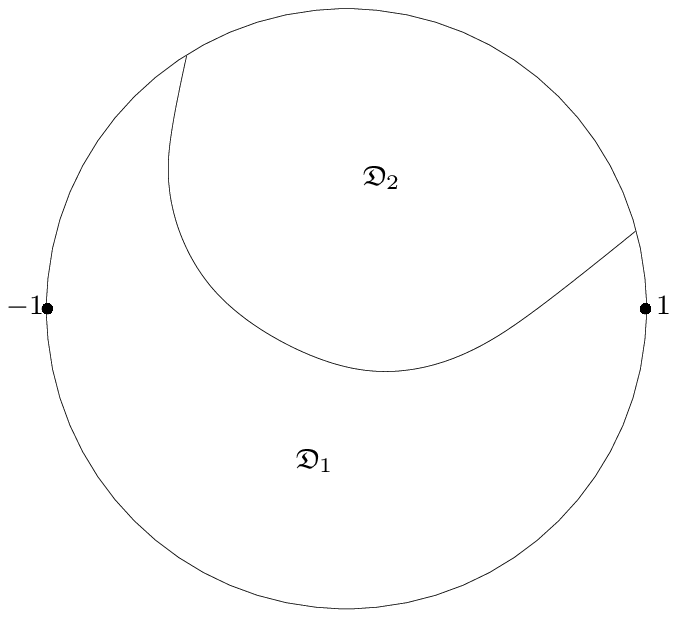} \qquad \psfig{file=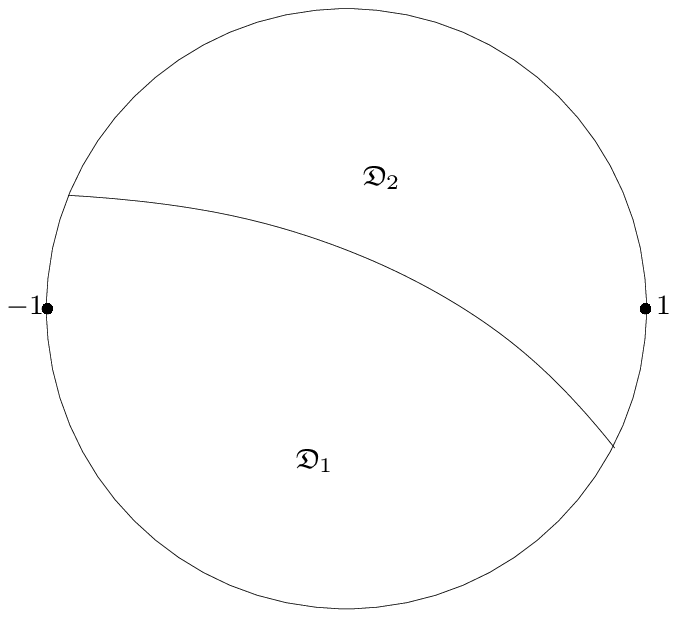}
      \end{center}
      \caption{Subcases ii and iii.}
      \label{f:subcases-ii-iii-n=2}
   \end{figure}

\end{proof}

\subsection{Proof of Lemmas~\ref{L:non-simple},~\ref{L:uv}}
\label{Sb:prf-L-non-simple}

\begin{proof}[Proof of Lemma~\ref{L:uv}]
   The Lemma is an immediate consequence of the following.
   \begin{prop} \label{P:n>=3}
      Suppose $n=\dim L \geq 3$. Then there exists a second category
      subset $\mathcal{J}_{\textnormal{reg}} \subset
      \mathcal{J}(M,\omega)$ such that for every $J \in
      \mathcal{J}_{\textnormal{reg}}$ the following holds:
      \begin{enumerate}
        \item For every simple $J$-holomorphic disk $u:(D, \partial D)
         \to (M,L)$ the set $u^{-1}(L) \cap \textnormal{Int\,} D$ is
         finite.
        \item For every two $J$-holomorphic disks $u, v:(D, \partial
         D) \to (M,L)$ which are simple and absolutely distinct the
         set $u(D) \cap v(D)$ is finite.
      \end{enumerate}
   \end{prop}

   Indeed, suppose that $u,v$ satisfy the assumptions of
   Lemma~\ref{L:uv}, i.e. $u,v$ are simple and $u(D) \cap v(D)$ is an
   infinite set. Statement~(2) of Proposition~\ref{P:n>=3} implies
   that $u,v$ are not absolutely distinct, namely $u(D) \subset v(D)$
   or $v(D) \subset u(D)$. Without loss of generality assume that
   $v(D) \subset u(D)$. It remains to show that $v(\partial D) \subset
   u(\partial D)$. To prove this, note that by statement~(1) of
   Proposition~\ref{P:n>=3} only finite number of points in $\partial
   D$ can be mapped by $v$ to $u(\textnormal{Int\,} D)$ for otherwise
   $u^{-1}(L) \cap \textnormal{Int\,} D$ would be an infinite set.
   Thus $v(\partial D \setminus \, \textnormal{finite set}) \subset
   u(\partial D)$. Since $v$ is continuous it easily follows that
   $v(\partial D) \subset u(\partial D)$.
\end{proof}

\begin{proof}[Proof of Proposition~\ref{P:n>=3}]
   Given two non-zero classes $F_1, F_2 \in H_2(M,L;\mathbb{Z})$, $J
   \in \mathcal{J}(M,\omega)$ and $\alpha, \beta \in \mathbb{Z}_{\geq
     0}$ denote by $\mathcal{M}^{*,\textnormal{d}}(F_1,F_2; J) \subset
   \mathcal{M}(F_1, J)\times \mathcal{M}(F_2,J)$ the subspace
   consisting of all pairs $(u,v)$ which are {\em simple and
     absolutely distinct}.

   Define $\mathcal{J}_{\textnormal{reg}}$ to be the subset of all $J
   \in \mathcal{J}(M,\omega)$ for which the following holds:
   \begin{itemize}
     \item For every non-zero homology class $F \in
      H_2(M,L;\mathbb{Z})$:
      \begin{itemize}
        \item The space $\mathcal{M}^*(F,J)$ is (either empty or) a
         smooth manifold of dimension $\mu(F)+n$.
        \item For every $\alpha \geq 1$ the evaluation map
         \begin{gather*}
            ev_{\alpha}: \mathcal{M}^*(F,J) \times (\textnormal{Int\,}
            D)^{\times \alpha} \longrightarrow M^{\times \alpha},\\
            ev_{\alpha}(u, p_1, \ldots, p_{\alpha}) = \bigl( u(p_1),
            \ldots, u(p_{\alpha}) \bigr)
         \end{gather*}
         is transverse to $L^{\times \alpha}$.
      \end{itemize}
     \item For every pair of non-zero homology classes $F_1, F_2 \in
      H_2(M,L;\mathbb{Z})$:
      \begin{itemize}
        \item The space $\mathcal{M}^{*,\textnormal{d}}(F_1,F_2; J)$
         is either empty or a smooth manifold of dimension $\leq
         2n+\mu(F_1)+\mu(F_2)$.
        \item For every $\alpha, \beta \in \mathbb{Z}_{\geq 0}$ the
         evaluation map
         $$ev_{\alpha,\beta}:
         \mathcal{M}^{*,\textnormal{d}}(F_1,F_2,\alpha,\beta; J)
         \times (\textnormal{Int\,}D)^{\times 2\alpha} \times
         (\partial D)^{\times 2\beta} \longrightarrow M^{\times
           2\alpha} \times L^{\times 2\beta},$$
         \begin{align*}
            & ev_{\alpha,\beta}(u, v, p_1, q_1, \ldots, p_{\alpha},
            q_{\alpha}, p'_1, q'_1, \ldots, p'_{\beta}, q'_{\beta}) \\
            & = \bigl(u(p_1),v(q_1), \ldots,
            u(p_{\alpha}),v(q_{\alpha}), u(p'_1),v(q'_1), \ldots,
            u(p'_{\beta}),v(q'_{\beta}) \bigr)
         \end{align*}
         is transverse to $\textnormal{diag}(M)^{\times \alpha} \times
         \textnormal{diag}(L)^{\times \beta}$.
      \end{itemize}
   \end{itemize}
   Standard arguments~\cite{McD-Sa:Jhol-2} show that the above subset
   $\mathcal{J}_{\textnormal{reg}} \subset \mathcal{J}(M,\omega)$ is
   indeed of second category.

   We now prove statement~(1) of Proposition~\ref{P:n>=3}. Let $J \in
   \mathcal{J}_{\textnormal{reg}}$ and let $u:(D,\partial D) \to
   (M,L)$ be a simple $J$-holomorphic curve.  Put $F=[u]$. By the
   transversality of the map $ev_{\alpha}$ we have $\dim
   ev_{\alpha}^{-1}(L^{\times \alpha}) = \mu(F)+n-\alpha(n-2)$.  As $n
   \geq 3$ it follows that for $\alpha \gg 1$, $\dim
   ev_{\alpha}^{-1}(L^{\times \alpha})<0$ hence $\dim
   ev_{\alpha}^{-1}(L^{\times \alpha})=\emptyset$. This proves
   statement~(1).

   We turn to the proof of statement~(2) of Proposition~\ref{P:n>=3}.
   Let $J \in \mathcal{J}_{\textnormal{reg}}$ and let $u,v:(D,\partial
   D) \to (M,L)$ be two simple $J$-holomorphic disks which are
   absolutely distinct. Put $F_1=[u]$, $F_2=[v]$. In view of
   statement~(1) of our proposition it is enough to prove that each of
   the following sets
   $$\bigl\{ (z_1,z_2) \in \textnormal{Int\,}D \times
   \textnormal{Int\,}D \,|\, u(z_1)=v(z_2) \bigr\}, \quad \bigl\{
   (z_1,z_2) \in \partial D \times \partial D \,|\, u(z_1)=v(z_2)
   \bigr\}$$
   is finite. By the transversality of the map
   $ev_{\alpha,\beta}$ we have $$\dim ev_{\alpha,\beta}^{-1}\bigl(
   \textnormal{diag}(M)^{\times \alpha} \times
   \textnormal{diag}(L)^{\times \beta} \bigr) = \mu(F_1)+\mu(F_2)
   - \beta (n-2) -2\alpha(n-2)+2n.$$
   As $n\geq 3$ it following that if
   $\alpha \gg 1$ or $\beta \gg 1$ then $$\dim
   ev_{\alpha,\beta}^{-1}\bigl( \textnormal{diag}(M)^{\times \alpha}
   \times \textnormal{diag}(L)^{\times \beta} \bigr)=\emptyset.$$
   This
   proves statement~(2).

   The proof of Proposition~\ref{P:n>=3} (hence of Lemma~\ref{L:uv}
   too) is complete.
\end{proof}

\begin{proof}[Proof of Lemma~\ref{L:non-simple}]
   Take $\mathcal{J}_{\textnormal{reg}} \subset \mathcal{J}(M,\omega)$
   to be the subset defined by Proposition~\ref{P:n>=3} and
   Lemma~\ref{L:uv}. Let $u:(D, \partial D) \to (M, L)$ be a
   non-simple $J$-holomorphic disk.

   Put $\mathcal{G} = \mathcal{G}(u)$. Let $\mathfrak{D}_1, \ldots,
   \mathfrak{D}_r \subset D \setminus \mathcal{G}$ be the connected
   components of the complement of $\mathcal{G}$. Let
   $\pi_{\overline{\mathfrak{D}}_j}:\overline{\mathfrak{D}}_j \to D$,
   $v_{\mathfrak{D}_j}:(D,\partial D) \to (M,L)$, $m_{\mathfrak{D}_j}
   \in \mathbb{N}$, $1\leq j \leq r$, be the maps and multiplicities
   given by Theorem~\ref{T:decomposition}. For simplicity we shall
   denote them by $\pi_j$, $v_j$, $m_j$, respectively.

   \smallskip \noindent \textbf{Case 1: $D \setminus \mathcal{G}$ has
     only one connected component (i.e. $r=1$).}

   Since $u$ is not simple, Theorem~\ref{T:decomposition} implies that
   $m_1 \geq 2$.  Put $u'=v_1$. Clearly $\omega([u']) < \omega([u])$.
   We also have $u'(\pi_1(-1))=u(-1)$ and $u'(\pi_1(1))=u(1)$.
   Finally, note that $\pi_1(-1) \neq \pi_1(1)$ (since $u(-1) \neq
   u(1)$) hence after a reparametrization $u'$ by an element of
   $\textnormal{Aut}(D)$ we may assume that $u'(-1)=u(-1)$,
   $u'(1)=u(1)$. The other properties claimed by the Lemma are
   obvious.

   \smallskip \noindent \textbf{Case 2: $D \setminus \mathcal{G}$ has
     more than one connected component (i.e. $r\geq 2$).}

   Define an abstract graph $\Gamma$ as follows (see
   figure~\ref{f:domains-1}). For each domain $\mathfrak{D}_i$ we
   assign a vertex $i \in \{1, \ldots, r\}$. We assign an edge between
   vertex $i'$ and vertex $i''$ if $\overline{\mathfrak{D}}_{i'} \cap
   \overline{\mathfrak{D}}_{i''}$ contains a $1$-dimensional component
   (in other words if the two domains have a $1$-dimensional common
   border). Note that $\Gamma$ is connected.

   \begin{figure}[htbp]
      \begin{center}
         \psfig{file=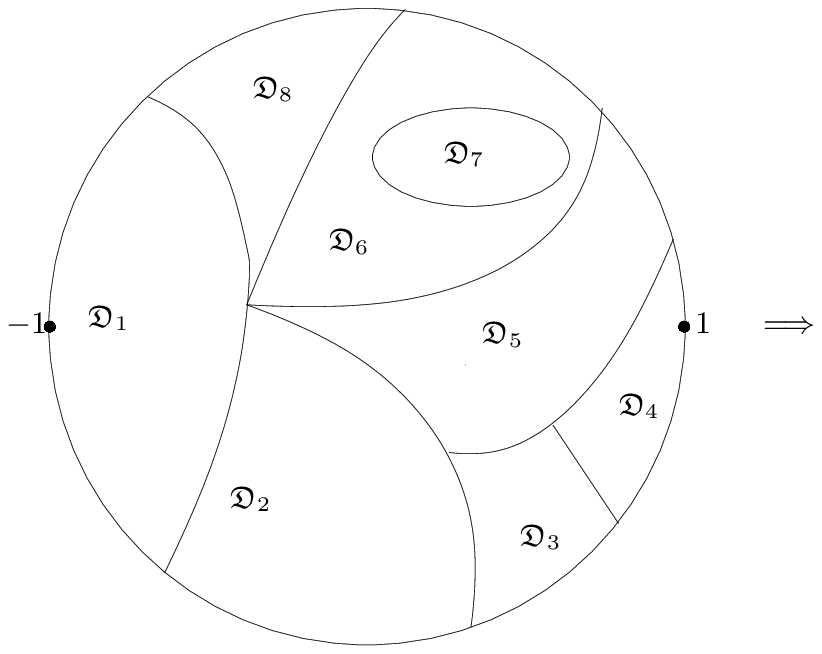} \psfig{file=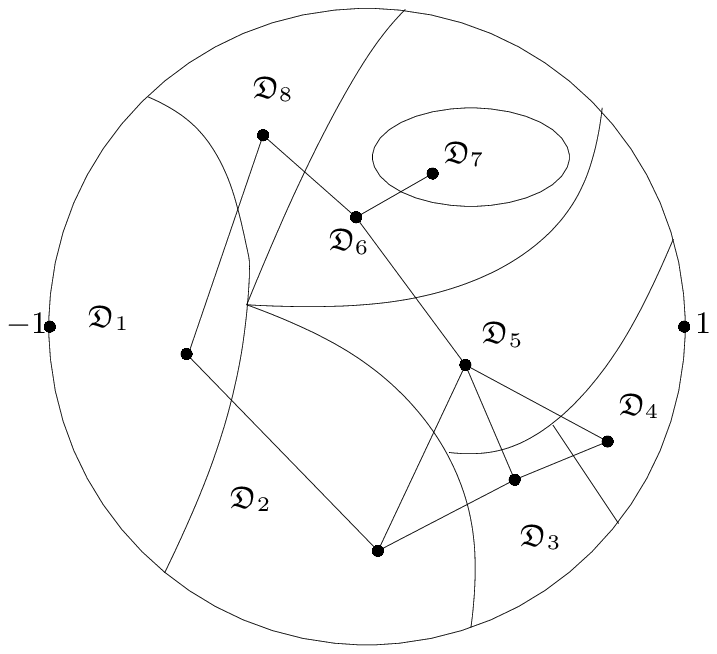}
      \end{center}
      \caption{The graph $\Gamma$.}
      \label{f:domains-1}
   \end{figure}

   Choose a path in $\Gamma$ that passes through each vertex of
   $\Gamma$ at least once. Denote the vertices of this path (in the
   order they appear in the path) by $t_1, \ldots, t_{\nu}$, $\nu \geq
   2$. We shall now construct a sequence of indices $1 \leq k_1 \leq
   \ldots \leq k_{\nu} \leq \nu$ such that for every $1 \leq q \leq
   \nu$:
   \begin{enumerate}[(I)]
     \item $v_{t_j}(D) \subset v_{t_{k_q}}(D)$ for every $1 \leq j
      \leq q$. \label{En:ind1}
     \item $v_{t_j}(\partial D) \subset v_{t_{k_q}}(\partial D)$ for
      every $1 \leq j \leq q$. \label{En:ind2}
   \end{enumerate}
   We construct the sequence $1 \leq k_1 \leq \ldots \leq k_{\nu} \leq
   \nu$ by induction as follows. Put $k_1=1$. By construction
   $v_{t_1}(D) \cap v_{t_2}(D)$ is an infinite set, hence
   Lemma~\ref{L:uv} implies that:
   \begin{itemize}
     \item either $v_{t_2}(D) \subset v_{t_1}(D)$ and
      $v_{t_2}(\partial D) \subset v_{t_1}(\partial D)$; or
     \item $v_{t_1}(D) \subset v_{t_2}(D)$ and $v_{t_1}(\partial D)
      \subset v_{t_2}(\partial D)$.
   \end{itemize}
   In the first case define $k_2=k_1=1$ and in the second case
   $k_2=2$.  Suppose that we have already constructed $k_1 \leq \ldots
   \leq k_q$ with properties~\ref{En:ind1},~\ref{En:ind2}. We define
   $k_{q+1}$ as follows. Since $v_{t_q}(D) \cap v_{t_{q+1}}(D)$ is
   infinite then $v_{t_{k_q}}(D) \cap v_{t_{q+1}}(D)$ is also an
   infinite set. By Lemma~\ref{L:uv} we have:
   \begin{itemize}
     \item either $v_{t_{q+1}}(D) \subset v_{t_{k_q}}(D)$ and
      $v_{t_{q+1}}(\partial D) \subset v_{t_{k_q}}(\partial D)$; or
     \item $v_{t_{k_q}}(D) \subset v_{t_{q+1}}(D)$ and
      $v_{t_{k_q}}(\partial D) \subset v_{t_{q+1}}(\partial D)$.
   \end{itemize}
   In the first case put $k_{q+1}=k_q$ and in the second case
   $k_{q+1}=q+1$. Clearly~\ref{En:ind1},~\ref{En:ind2} hold now with
   $q$ replaced by $q+1$. By induction we get the desired sequence $1
   \leq k_1 \leq \ldots \leq k_{\nu} \leq \nu$.

   Put $u'=v_{t_{k_{\nu}}}$. Properties~(2),(3) claimed by the lemma
   are obvious. Property~(4) follows from the fact that $r\geq 2$.
   Finally, since $u(-1) \neq u(1)$ we must have two {\em distinct}
   points $z_1, z_2 \in \partial D$ with $u'(z_1)=u(-1)$,
   $u'(z_2)=u(1)$.  Thus after a suitable reparametrization of $u'$ we
   may assume that $z_1=-1$ and $z_2=1$. This proves property~(1)
   claimed by the lemma.  The proof of Lemma~\ref{L:non-simple} is
   complete.
\end{proof}



%
%


\section{Gluing}
\label{S:gluing}

In essence, the gluing of $J$-holomorphic disks appears already in the
literature in the work of Fukaya-Oh-Ohta-Ono~\cite{FO3} (see
also~\cite{Aka:gluing}). However for the purposes of this paper we
need a small variation of the gluing theorem of~\cite{FO3}, and
moreover we also need the surjectivity of the gluing map which is not
explicitly discussed in~\cite{FO3}.  Therefore, for the sake of
completeness we felt it useful to include a detailed argument for
gluing in which we closely follow the original proof of
Fukaya-Oh-Ohta-Ono~\cite{FO3} as well as a proof for the surjectivity
of the gluing map.  We also discuss here the gluing for the pearls
introduced in \S\ref{S:transversality} and for some other of the
elements of the moduli spaces which will be used in
\S\ref{S:proofalg}. Other variants of these gluing statements will be used
sometimes in the paper - in particular, we focus here on the case of a fixed almost complex
structure but there is  sometimes a need to allow this structure to vary inside a family.
All of them are obtained by rather direct modifications of the gluing arguments presented here.

\subsection{Main statements} Let $(M^{2n}, \omega)$ be a tame
symplectic manifold endowed with an $\omega$-compatible almost complex
structure $J$. Let $L\subset M$ be a closed Lagrangian submanifold.
Let $u_1, u_0: (D, \partial D) \to (M,L)$ be two $J$-holomorphic
disks. Put $A_i=[u_i] \in H_2(M,L)$. Denote by $\mathcal{M}(A,J)$ the
space of $J$-holomorphic disks with boundary on L, in the class $A$.

Let $W$ be a manifold and $\mathbf{h}:W \to L \times M \times M \times
L$ a smooth map. We shall denote the components of $\mathbf{h}$ by
$h_{-}, h_1, h_0, h_{+}$ so that $\mathbf{h}(q) = (h_{-}(q), h_1(q),
h_0(q), h_{+}(q)) \in L \times M \times M \times L$ for every $q \in
W$. Fix two points lying on the {\em real} part of the disk $z_1,z_0
\in (\textnormal{Int\,}D) \cap \mathbb{R}$ and a point $q_* \in W$. In
what follows we shall put the following assumption on $u_1,u_0$ and
$J$:
\begin{assumption} \label{A:T1}
   \begin{enumerate}
     \item $u_1(1)=u_0(-1)$.
     \item $\mathbf{h}(q_*) = \bigl( u_1(-1), u_1(z_1), u_0(z_0),
      u_0(1) \bigr).$
     \item $J$ is regular for both $u_1$ and $u_0$ in the sense that
      the linearizations $D_{u_1}, D_{u_0}$ of the
      $\overline{\partial}$ operator at $u_1,u_0$ are surjective.
      \label{I:T1-reg-01}
     \item Let $ev: \mathcal{M}(A_1,J) \times \mathcal{M}(A_0,J) \to L
      \times M \times L \times L \times M \times L$ be the evaluation
      map
      $$ev(v_1,v_0)=\big( v_1(-1), v_1(z_1), v_1(1), v_0(-1),
      v_0(z_0), v_0(1) \big).$$
      Define a map $\mathbf{h}_{\Delta_L} :
      W \times L \to L \times M \times L \times L \times M \times L$
      by
      $$\mathbf{h}_{\Delta_L}(q,p)= \bigr( h_{-}(q), h_1(q), p, p,
      h_0(q), h_{+}(q) \bigr).$$
      Put $p_*=u_1(1)=u_0(-1)$. Then we
      assume that $ev$ and $\mathbf{h}_{\Delta_L}$ are mutually
      transverse at the points $(u_1, u_0) \in \mathcal{M}(A_1, J)
      \times \mathcal{M}(A_0, J)$ and $(q_*, p_*) \in W \times L$.
   \end{enumerate}
\end{assumption}

Put $A=A_1+A_0$. Consider the space of all $(u, r, q) \in
\mathcal{M}(A,J) \times (0,1) \times W$ such that
\begin{equation} \label{Eq:MB}
   \bigl( u(-1), u(-r), u(r), u(1) \bigr) = \mathbf{h}(q).
\end{equation}
We denote the space of $(u, r, q)$'s described in~\eqref{Eq:MB} by
$\mathcal{M}(A,J; \mathcal{C}(\mathbf{h}))$ (Here
$\mathcal{C}(\mathbf{h})$ stands for the configuration described by
conditions~\eqref{Eq:MB}.)

\begin{thm}[Gluing] \label{T:gluing}
   Under Assumption~\ref{A:T1} there exists a path $\{ (v_s, a(s),
   q_s) \}_{0 < s} \subset \mathcal{M}(A,J;\mathcal{C}(\mathbf{h}))$
   with the following properties:
   \begin{enumerate}
     \item $q_s \xrightarrow[s \to \infty]{} q_*$ and $a(s)
      \xrightarrow[s \to \infty]{} 1$.
     \item $v_s$ converges with the marked points $(-1,-a(s),a(s),1)$,
      as $s \to \infty$, to $(u_1,u_0)$ with the marked points
      $(-1,z_1),\,(z_0,1)$ in the Gromov topology.
   \end{enumerate}
   In particular, the point $\bigl( (u_1, z_1),(u_0,z_0), q_* \bigr)$
   lies in the boundary of the closure
   $\overline{\mathcal{M}}(A,J;\mathcal{C}(\mathbf{h}))$ of the space
   $\mathcal{M}(A,J;\mathcal{C}(\mathbf{h}))$ in the sense of (1), (2)
   above. Furthermore, if $\mu(A)+\dim W-5n=0$ then the above path is
   unique in the following sense. There exists a neighbourhood
   $\mathcal{U}$ of the point $\bigl((u_1,z_1),(u_0,z_0), q_*\bigr)$
   in $\overline{\mathcal{M}}(A,J;\mathcal{C}(\mathbf{h}))$ such that
   $\mathcal{U} \setminus \bigl\{ \bigl((u_1,z_1),(u_0,z_0),
   q_*\bigr)\bigr\}$ coincides with the path $\{ (v_s, a(s), q_s)\}$
   for $s \gg 0$. In other words, every path $\{ (w_s, a'(s), q'_s) \}
   \subset \mathcal{M}(A,J;\mathcal{C}(\mathbf{h}))$ with $q'_s
   \xrightarrow[s \to \infty]{} q_*$ and such that $w_s \xrightarrow[s
   \to \infty]{} (u_1, u_0)$ with marked points as in (2) is obtained
   from $\{ (v_s, a(s), q_s) \}$ by reparametrization in $s$, for $s
   \gg 0$.
\end{thm}

The proof of Theorem~\ref{T:gluing} will occupy
Section~\ref{Sb:overivew}-~\ref{S:auxiliary} below.

\begin{remsnonum}
   \begin{enumerate}
     \item The uniqueness statement seems to hold without the
      assumption $\mu(A)+ \dim W -5n=0$, however the proof is much
      more complicated in that case. Anyway, we shall not need this
      more general statement.
     \item The requirement that the points $z_1, z_0$ lie on the real
      axis of $D$ is not crucial. Indeed a similar theorem holds for
      any choice of $z_1, z_0 \in \textnormal{Int\,}D$ but the marked
      points $-r, r$ used to define
      $\mathcal{M}(A,J;\mathcal{C}(\mathbf{h}))$ must be changed
      accordingly.
   \end{enumerate}
\end{remsnonum}

\begin{figure}[htbp]
   \begin{center}
      \psfig{file=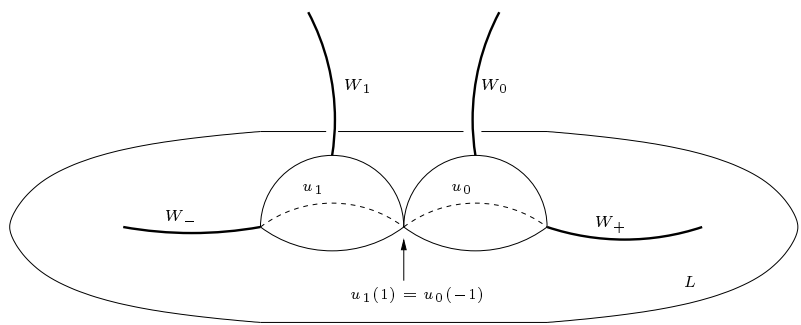}
   \end{center}
   \caption{Before gluing}
   \label{f:disks-1}
\end{figure}

\begin{figure}[htbp]
   \begin{center}
      \psfig{file=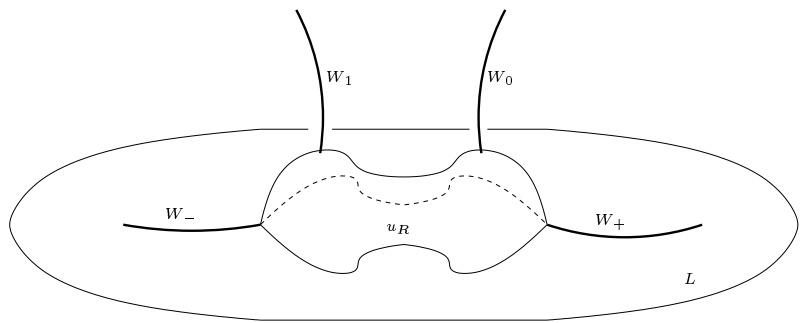}
   \end{center}
   \caption{After gluing}
   \label{f:disks-2}
\end{figure}

\subsection{Examples} \label{S:gluing-examples}
Before going to the proof of Theorem~\ref{T:gluing} we present several
typical applications of Theorem~\ref{T:gluing}. These will occupy
Section~\ref{Sbsb:glue-ex-1}-~\ref{Sbsb:glue-ex-2II} below.

In most of the situations the manifold $W$ will be a taken to be a
product $W=W_{-} \times W_1 \times W_0 \times W_{+}$ and
$\mathbf{h}=h_{-} \times h_1 \times h_0 \times h_{+}$ with
$h_{\pm}:W_{\pm} \to L$, $h_i:W_i \to M$, $i=1,0$. A realistic
illustration of the gluing process is given in
figures~\ref{f:disks-1},~\ref{f:disks-2}. In these figures $W_{\pm}$
are taken to be submanifold of $L$, $W_i$ submanifolds of $M$ and the
maps $h_{\pm}$, $h_i$ are the inclusions.

\subsubsection{Gluing two trajectories of pearls}
\label{Sbsb:glue-ex-1}

Let $f: L \to \mathbb{R}$ be a Morse function and $\rho$ a Riemannian
metric on $L$. Assume that $(f, \rho)$ is Morse-Smale. Given two
vectors of non-zero classes in $H_2(M,L;\mathbb{Z})$,
$\mathbf{B}'=(B'_1, \ldots, B'_{l'})$, $\mathbf{B}''=(B''_1, \ldots,
B''_{l''})$ we denote $\mathbf{B}' \# \mathbf{B}'' = (B'_1, \ldots,
B'_{l'-1}, B'_{l'}+B''_1, B''_2, \ldots, B''_{l''})$. For $x,y \in
\textnormal{Crit}(f)$ and a vector $\mathbf{C}$ of non-zero classes we
use the notation $\mathcal{P}(x,y,\mathbf{C}; J,f,\rho)$,
$\mathcal{P}(x,y, \mathbf{B}', \mathbf{B}''; J, f, \rho)$ introduced
in Section~\ref{Sb:tr-pearl-complex}.

The following is a corollary of Theorem~\ref{T:gluing}.
\begin{cor} \label{C:glue-pearls-1}
   Let $(f, \rho)$ be as above. There exists a second category subset
   $\mathcal{J}_{\textnormal{reg}} \subset \mathcal{J}(M,\omega)$ such
   that for every $J \in \mathcal{J}_{\textnormal{reg}}$, every pair
   of vectors of non-zero classes $\mathbf{B}', \mathbf{B}''$ and
   every $x,y\in \textnormal{Crit}(f)$ with
   $\ind_f(x)-\ind_f(y)+\mu(\mathbf{B}')+\mu(\mathbf{B}'')-1=1$ the
   following holds. For every $(\mathbf{v}', \mathbf{v}'') \in
   \mathcal{P}(x,y, \mathbf{B}',\mathbf{B}'';J,f,\rho)$ there exists a
   path $\{ \mathbf{u}_s \} \subset \mathcal{P} (x,y,\mathbf{B}' \#
   \mathbf{B}''; J, f, \rho)$ which converges in the Gromov topology
   as $s \to \infty$ to $(\mathbf{v}',\mathbf{v}'')$.  Moreover, the
   end of the $1$-dimensional manifold $\mathcal{P}(x,y,\mathbf{B}' \#
   \mathbf{B}'';J,f,\rho)$ parametrized by $\{ \mathbf{u}_s\}$ is
   unique in the sense that every other path $\{ \mathbf{w}_s\}$ in
   $\mathcal{P}(x,y,\mathbf{B}' \# \mathbf{B}'';J,f,\rho)$ that
   converges to $(\mathbf{v}',\mathbf{v}'')$ as $s \to \infty$, lies
   in the same end.
\end{cor}
\begin{proof}
   Fix an element $(\mathbf{v}', \mathbf{v}'') = (v'_1, \ldots,
   v'_{l'}, v''_1, \ldots, v''_{l''}) \in
   \mathcal{P}(x,y,\mathbf{B}',\mathbf{B}'';J,f,\rho)$. Recall from
   Proposition~\ref{P:tr-2} that by taking $J$ to be generic we may
   assume that $\mathcal{P}(x,y,\mathbf{B}',\mathbf{B}'';J,f,\rho)$ is
   a finite set and that the disks $(v'_1, \ldots, v'_{l'}, v''_1,
   \ldots, v''_{l''})$ are simple and absolutely distinct.

   We now define the manifold $W$ and the map $\mathbf{h}$ used for
   applying Theorem~\ref{T:gluing}. The manifold $W$ will be a product
   $W=W_{-} \times W_1 \times W_0 \times W_{+}$ defined as follows.
   If $l'=1$ put $W_{-}=W_x^u$. If $l'\geq 2$ define first
   $\widehat{W}_{-}$ to be the space of all $(w'_1, \ldots,
   w'_{l'-1},p')$ such that:
   \begin{itemize}
     \item $w'_i \in \mathcal{M}(B'_i,J)$ for every $1 \leq i \leq
      l'-1$.
     \item $w'_1(-1) \in W_x^u$.
     \item $(w'_i(1), w'_{i+1}(-1)) \in Q_{(f,\rho)}$ for every $1
      \leq i \leq l'-2$. (See formula~\eqref{Eq:Q-f} in
      Section~\ref{Sb:tr-pearl-complex} for the definition of
      $Q_{(f,\rho)}$.)
     \item $(w'_{l'-1}(1),p') \in Q_{(f,\rho)}$.
     \item $(w'_1, \ldots, w'_{l'-1})$ are simple and absolutely
      distinct.
   \end{itemize}
   Finally put $W_{-}=\widehat{W}_{-}/G_{-1,1}^{\times (l'-1)}$.
   Define $h_{-}:W_{-} \to L$ as follows. If $l'=1$, let $h_{-}$ be
   the inclusion. If $l'' \geq 2$ define $h_{-}(w'_1, \ldots,
   w'_{l'-1},p')=p'$.

   We define $W_{+}$ and $h_{+}:W_{+} \to L$ in an analogous way. We
   write elements of $W_{+}$ as $(p'', w''_2, \ldots, w''_{l''})$.
   Standard transversality results imply that for generic $J$, the
   spaces $W_{-}, W_{+}$ are smooth manifold of dimensions $\ind_f(x) +
   \sum_{i=1}^{l'-1} \mu(B'_i)$ and $n-\ind_f(y) + \sum_{i=2}^{l''}
   \mu(B''_i)$ respectively.

   Put $u_1=v'_{l'}$, $u_0=v''_1$, $p_*=u_1(1)=u_0(-1) \in L$ and set:
   \begin{itemize}
     \item $q^*_{-} = (v'_1, \ldots, v'_{l'-1},p')$, where
      $p'=u_1(-1)$,
     \item $q^*_{+}=(p'', v''_2, \ldots, v''_{l''})$, where
      $p''=u_0(1)$.
   \end{itemize}
   Note that $q^*_{-} \in W_{-}$, $q^*_{+} \in W_{+}$. Set $A_1 =
   B'_{l'}$ and $A_0 = B''_1$.  Consider the following maps:
   \begin{equation*}
      \begin{aligned}
         & h': W_{-} \times W_{+} \times L \to L \times L \times L
         \times L, \quad
         h'(q_{-},q_{+},p) = (h_{-}(q_{-}), p, p, h_{+}(q_{+})), \\
         & ev_{1,0}: \mathcal{M}^*(A_1,J) \times \mathcal{M}^*(A_0,J)
         \to L
         \times L \times L \times L, \\
         & ev_{1,0}(v_1, v_0) = (v_1(-1), v_1(1), v_0(-1), v_0(1)).
      \end{aligned}
   \end{equation*}
   Here $\mathcal{M}^*(A_i,J)$, $i=1,0$, stands for the space of {\em
     simple} $J$-holomorphic disks in the class $A_i$. Again, by
   taking $J$ to be generic we may assume that $\mathcal{M}^*(A_i,J)$
   are smooth manifolds and that the maps $h'$ and $ev_{1,0}$ are
   mutually transverse at the points $(q^*_{-}, q^*_{+},p_*) \in W_{-}
   \times W_{+} \times L$ and $(u_1, u_0) \in \mathcal{M}^*(A_1, J)
   \times \mathcal{M}^*(A_0, J)$. (For this to hold it is crucial to
   know that the disks $(v'_1, \ldots, v'_{l'}, v''_1, \ldots,
   v''_{l''})$ are simple and absolutely distinct.)

   We turn to the manifolds $W_1, W_0$. Let $z_1, z_0 \in
   (\textnormal{Int\,}D) \cap \mathbb{R}$ be two points for which
   ${du_1}_{(z_1)}, {du_0}_{(z_0)} \neq 0$. (Note that since $u_1,
   u_0$ are $J$-holomorphic and not constant, such two points $z_1,
   z_0$ do exist.) Fix two $(2n-1)$-dimensional manifolds $W_1, W_0
   \subset M$ with $u_1(z_1) \in W_1$, $u_0(z_0) \in W_0$ and such
   that $u_1, u_0:D \to M$ are transverse to $W_1, W_0$ at the points
   $z_1, z_0$ respectively.  Define $h_i:W_i \to M$ to be the
   inclusions.

   Define $W=W_{-} \times W_1 \times W_0 \times W_{+}$ and
   $\mathbf{h}:W \to L \times M \times M \times L$ to be
   $\mathbf{h}=(h_{-}, h_1, h_0, h_{+})$. Put $p_*=u_1(1)=u_0(-1) \in
   M$ and $q_*=(q^*_{-}, q^*_1, q^*_0, q^*_{+}) \in W$ where $q^*_{-},
   q^*_{+}$ are defined above and $q^*_1 = u_1(z_1)$,
   $q^*_0=u_0(z_0)$. A simple computation shows that the maps
   $\mathbf{h}_{\Delta_L}$ and $ev$ (defined in Assumption~\ref{A:T1})
   are mutually transverse at the points $(q_*, p_*) \in W \times L$
   and $(u_1, u_0) \in \mathcal{M}(A_1, J) \times \mathcal{M}(A_0,J)$.
   Clearly the rest of the assumptions in~\ref{A:T1} are also
   satisfied. Note that in the above construction once $J$ is fixed
   the manifolds $W_{-}$, $W_{+}$ are determined in a canonical way.
   The manifolds $W_1$, $W_0$ on the other hand are chosen a
   posteriori and depend on the element $(\mathbf{v}', \mathbf{v}'')$.

   We now apply Theorem~\ref{T:gluing}. We obtain from this Theorem a
   family of $J$-holomorphic disks $v_s \in \mathcal{M}(A_1+A_0,J)$
   together with marked points $-a(s), a(s) \in (\textnormal{Int\,}D)
   \cap \mathbb{R}$ and points $q_s = (q_{-,s}, q_{1,s}, q_{0,s},
   q_{+,s}) \in W_{-} \times W_1 \times W_0 \times W_{+}$, such that:
   \begin{itemize}
     \item $q_{-,s} = (u'_{1,s}, \ldots, u'_{l'-1, s}, p'_s)
      \xrightarrow[s \to \infty]{} q^*_{-} = (v'_1, \ldots,
      v'_{l'-1},p_*)$.
     \item $q_{+,s} = (p''_s, u''_{2,s}, \ldots, u'_{l'', s})
      \xrightarrow[s \to \infty]{} q^*_{+} = (p_*, v''_2, \ldots,
      v''_{l''})$.
     \item $v_s(-1)=p'_s$, $v_s(1)=p''_s$.
     \item $v_s(-a(s))=q_{1,s} \xrightarrow[s \to \infty]{} u_1(z_1)$,
      $v_s(a(s))=q_{0,s} \xrightarrow[s \to 1\infty]{} u_0(z_0)$.
     \item $v_s$ converges with the marked points $(-1, -a(s), a(s),
      1)$ to $(u_1, u_0)$ with the marked points $(-1, z_1)$, $(z_0,
      1)$.
   \end{itemize}
   Put $\mathbf{u}_s = (u'_{1,s}, \ldots, u'_{l'-1, s}, v_s,
   u''_{2,s}, \ldots, u'_{l'', s})$. Clearly $\mathbf{u}_s \in
   \mathcal{P}(x,y,\mathbf{B'}\# \mathbf{B}''; J, f, \rho)$ and
   $\mathbf{u}_s$ converges as $s \to \infty$ to $(\mathbf{v}',
   \mathbf{v}'')$.

   We turn to the uniqueness statement. Let $\{ \mathbf{w}_s \}$ be a
   path in $\mathcal{P}(x,y,\mathbf{B'}\# \mathbf{B}''; J, f, \rho)$
   such that $\mathbf{w}_s \xrightarrow[s \to \infty]{} (\mathbf{v}',
   \mathbf{v}'')$. Write $\mathbf{w}_s = (w_{1,s}, \ldots, w_{l'-1,s},
   w_{l',s}, \ldots, w_{l'+l''-1, s})$. Then we have
   $w_{l',s}\xrightarrow[s \to \infty]{} (u_1, u_0)$ in the Gromov
   topology. By applying a suitable family of holomorphic
   reparametrization to $w_{l',s}$ we may assume that the maps
   $w_{l',s}$ uniformly converge in the $C^{\infty}$-topology, as $s
   \to \infty$, to $u_1$ on compact subsets of $D \setminus \{ 1 \}$.
   Similarly, after (other) reparametrizations, $w_{l',s}$ uniformly
   converges in the $C^{\infty}$-topology, as $s \to \infty$, to $u_0$
   on compact subsets of $D \setminus \{ -1 \}$. Due to the
   transversality between $u_i$ and $W_i$ at $z_i \in D$, $i=1,0$, it
   follows that there exists points $b_1(s), b_0(s) \in
   (\textnormal{Int\,}D) \cap \mathbb{R}$ with $b_1(s) \xrightarrow[s
   \to \infty]{} -1$, $b_0(s) \xrightarrow[s \to \infty]{} 1$ such
   that $w_{l',s}(b_i(s)) \in W_i$ and $w_{l',s}(b_i(s))
   \xrightarrow[s \to \infty]{} u_i(z_i)$, $i=1,0$.

   After further reparametrizations we may assume that
   $b_1(s)=-b_0(s)$. As before we construct elements $q'_s \in W$
   using $(w_{1,s}, \ldots, w_{l'-1,s})$, $(w_{l'+1,s}, \ldots,
   w_{l'+l''-1,s})$ and $w_{l',s}(\pm b_0(s))$ such that $(w_{l',s},
   b_0(s), q'_s) \in \mathcal{M}(A,J;\mathcal{C}(\mathbf{h}))$, $q'_s
   \xrightarrow[s \to \infty]{} q_*$ and such that $w_{l',s}$
   converges with the marked points $(-1, -b_0(s), b_1(s), 1)$ as $s
   \to \infty$ to $(u_1, u_0)$ with the marked points $(-1, z_1)$,
   $(z_0, 1)$. Noting that
   \begin{align*}
      & \mu(A_1)+\mu(A_0) + \dim W -5n = \\
      & \mu(B'_{l'})+ \mu(B''_1) + \Bigl( \ind_f(x) + \sum_{i=1}^{l'-1}
      \mu(B'_i)\Bigr) + 2(2n-1) + \Bigl( n-\ind_f(y) + \sum_{i=2}^{l''}
      \mu(B''_i)\Bigr)
      -5n = \\
      & \ind_f(x)-\ind_f(y)+\mu(\mathbf{B}')+\mu(\mathbf{B}'') - 2 = 0,
   \end{align*}
   it follows from the uniqueness statement of Theorem~\ref{T:gluing}
   that for $s \gg 0$, $(w_{l',s}, b_0(s), q'_s)$ coincides with
   $(v_s, a(s), q_s)$ up to reparametrizations in $s$. It follows that
   the path $\mathbf{w}_s$ and $\mathbf{u}_s$ are the same up to
   reparametrization for $s \gg 0$. This completes the proof of
   Corollary~\ref{C:glue-pearls-1}.

   \begin{remnonum}
      The manifolds $W_1, W_0$ above were important only for the
      uniqueness statement. We had to choose them to be
      $(2n-1)$-dimensional in order to reduce the dimension of the
      space $\mathcal{M}(A,J; \mathcal{C}(\mathbf{h}))$ to be $1$,
      i.e. to assure that $\mu(A_1)+\mu(A_0) + \dim W -5n=0$,
      which is the assumption we need for the uniqueness statement in
      Theorem~\ref{T:gluing}.
   \end{remnonum}

\end{proof}

\subsubsection{Gluing a trajectory of pearls to a trajectory
  with external constrains I} \label{Sbsb:glue-ex-2I}

Here we show how to apply Theorem~\ref{T:gluing} in order to glue a
trajectory of pearls to a trajectory from the space
$\mathcal{P}_{III_i}(a,x,y;\mathbf{B}', \mathbf{B}'',J)$, $i=1,2$,
defined in Section~\ref{Sb:qm-spaces} in the context of the quantum
module structure.

Let $h:M \to \mathbb{R}$, $f:L \to \mathbb{R}$ be Morse functions and
$\rho_M$, $\rho_L$ Riemannian metrics on $M$, $L$. Assume that
$(f,\rho_L)$ and $(h, \rho_M)$ satisfy Assumption~\ref{A:MS}.

\begin{cor} \label{C:glue-pearls-2}
   Let $(f,\rho_L)$, $(h,\rho_M)$ be as above.  There exists a second
   category subset $\mathcal{J}_{\textnormal{reg}} \subset
   \mathcal{J}(M,\omega)$ such that for every $J \in
   \mathcal{J}_{\textnormal{reg}}$, every pair of vectors of non-zero
   classes $\mathbf{B}', \mathbf{B}''$ and every $x,y\in
   \textnormal{Crit}(f)$, $a \in \textnormal{Crit}(h)$ with $\ind_h(a)
   + \ind_f(x)-\ind_f(y)+\mu(\mathbf{B}')+\mu(\mathbf{B}'')-2n=1$ the
   following holds. For every $(\mathbf{v}', \mathbf{v}'') \in
   \mathcal{P}_{III_i}(a, x,y, \mathbf{B}',\mathbf{B}'';J)$ there
   exists a path $\{ \mathbf{u}_s \} \subset \mathcal{P}_I(a,
   x,y,\mathbf{B}' \# \mathbf{B}''; J)$ which converges in the Gromov
   topology as $s \to \infty$ to $(\mathbf{v}',\mathbf{v}'')$.
   Moreover, the end of the $1$-dimensional manifold
   $\mathcal{P}_I(a,x,y,\mathbf{B}' \# \mathbf{B}'';J)$ parametrized
   by $\{ \mathbf{u}_s\}$ is unique in the sense that every other path
   $\{ \mathbf{w}_s\}$ in $\mathcal{P}_I(a,x,y,\mathbf{B}' \#
   \mathbf{B}'';J)$ that converges to $(\mathbf{v}',\mathbf{v}'')$ as
   $s \to \infty$, lies in the same end.
\end{cor}
\begin{proof}
   We outline the proof for the space
   $\mathcal{P}_{III_1}(a,x,y;\mathbf{B}', \mathbf{B}'',J)$.  Suppose
   that $\mathbf{B}'=(B'_1, \ldots, B'_{l'})$, $\mathbf{B}''=(B''_1,
   \ldots, B''_{l''})$.  Recall that
   $\mathcal{P}_{III_1}(a,x,y;\mathbf{B}', \mathbf{B}'',J)$ is
   disjoint union of the spaces
   $\mathcal{P}_{III_1}(a,x,y;(\mathbf{B}',k'),\mathbf{B}'',J)$ where
   $k'$ goes from $1$ to $l'$.

   In case $(\mathbf{v}', \mathbf{v}'') \in
   \mathcal{P}_{III_1}(a,x,y;(\mathbf{B}',k'),\mathbf{B}'',J)$ with
   $1\leq k' \leq l'-1$ the proof is very similar to the proof of
   Corollary~\ref{C:glue-pearls-1}. The only difference is that the
   manifold $W_{-}$ is now defined by elements $(w'_1, \ldots,
   w'_{l'-1},p')$ with the additional condition that $w'_{k'}(0) \in
   W_a^u$. The rest of the proof goes in the same way as for
   Corollary~\ref{C:glue-pearls-1}.

   In case $k'=l'$, i.e. $(\mathbf{v}', \mathbf{v}'') \in
   \mathcal{P}_{III_1}(a,x,y;(\mathbf{B}',l'),\mathbf{B}'',J)$ the
   proof is even easier. We define the manifold $W_1$ to be $W_a^u$
   and take $z_1 = 0 \in \textnormal{Int\,} D$. The manifolds
   $W_{\pm}$, $W_0$ are defined in the same way as in the proof of
   Corollary~\ref{C:glue-pearls-1}.
\end{proof}

\subsubsection{Gluing trajectories with external constrains and
  Hamiltonian perturbations} \label{Sbsb:glue-ex-H}

Here we show how to perform gluing on elements from the space
$\mathcal{P}_{III}(a,x,y; \mathbf{B}', \mathbf{B}'',J, H)$ introduced
in Section~\ref{Sb:perturb-P}, in order to obtain elements from the
space $\mathcal{P}_I(a,x,y; \mathbf{B}'\# \mathbf{B}'',J,H)$.

Let $h:M \to \mathbb{R}$, $f:L \to \mathbb{R}$ be Morse functions and
$\rho_M$, $\rho_L$ Riemannian metrics on $M$, $L$. Assume that
$(f,\rho_L)$ and $(h, \rho_M)$ satisfy Assumption~\ref{A:MS}.
\begin{cor}\label{C:glue-pearls-H}
   Let $(f,\rho_L)$, $(h,\rho_M)$ be as above.  There exists a second
   category subset $\mathcal{J}_{\textnormal{reg}} \subset
   \mathcal{J}(M,\omega)$ such that for every $J \in
   \mathcal{J}_{\textnormal{reg}}$ there exists a second category
   subset $\mathcal{H}_{\textnormal{reg}}(J) \subset \mathcal{H}$ with
   the following properties. Let $H \in
   \mathcal{H}_{\textnormal{reg}}(J)$, $a \in \textnormal{Crit}(h)$,
   $x,y \in \textnormal{Crit}(f)$, $\mathbf{B}'=(B'_1, \ldots,
   B'_{l'})$, $\mathbf{B}''=(B''_1, \ldots, B''_{l''})$ be two vectors
   of non-zero classes with $\ind_h(a) + \ind_f(x) -\ind_f(y) +
   \mu(\mathbf{B}') + \mu(\mathbf{B}'') -2n = 1$. Let $(\mathbf{v}',
   \mathbf{v}'') \in \mathcal{P}_{III_i}(a,x,y; \mathbf{B}',
   \mathbf{B}'', J, H)$, $i=1,2$. Then there exists a path $\{
   \mathbf{u}_s \} \subset \mathcal{P}_I(a,x,y; \mathbf{B}' \#
   \mathbf{B}'',J, H)$ which converges in the Gromov topology, as $s
   \to \infty$, to $(\mathbf{v}', \mathbf{v}'')$. Moreover, the end of
   the $1$-dimensional manifold $\mathcal{P}_I(a,x,y; \mathbf{B}' \#
   \mathbf{B}'',J, H)$ parametrized by this path is unique in the
   sense that every other path with the same property for $s \to
   \infty$, lies in the same end.
\end{cor}

\begin{rem} \label{R:glue-II-H}
   A similar statement holds for the spaces $\mathcal{P}_{II_i}(a,x,y;
   \mathbf{A}, J, H)$ from Section~\ref{Sb:perturb-P}, whenever
   $\ind_h(a)+\ind_f(x)-\ind_f(y)+\mu(\mathbf{A})-2n=1$. The proof is
   almost the same as the one given below.
\end{rem}

\begin{proof}[Proof of Corollary~\ref{C:glue-pearls-H}]
   We prove the Corollary for the space $\mathcal{P}_{III_1}(a,x,y;
   \mathbf{B}', \mathbf{B}'', J, H)$. Recall that this space is a
   disjoint union of $\mathcal{P}_{III_1}(a,x,y; (\mathbf{B}',k'),
   \mathbf{B}'', J, H)$, where $k'$ goes from $1$ to $l'$. Assume that
   $(\mathbf{v}', \mathbf{v}'')$ lies in the $k'$'th component of this
   space for some $1 \leq k'\leq l'$.

   If $k'<l'$ the proof is essentially the same as for
   Corollary~\ref{C:glue-pearls-2}.

   Assume $k'=l'$. Write $\mathbf{v}'=(v'_1, \ldots, v'_{l'})$,
   $\mathbf{v}''=(v''_1, \ldots, v''_{l''})$. We have to glue the
   $(J,H)$-holomorphic disk $v'_{l'}$ to the (genuine) $J$-holomorphic
   disk $v''_1$, preserving the constrains imposed by
   $\mathcal{P}_I(a,x,y; (\mathbf{B}' \# \mathbf{B}'',k'),J, H)$.
   (Recall from Section~\ref{Sb:perturb-P} that the space
   $\mathcal{P}_I(a,x,y; \mathbf{B}' \# \mathbf{B}'',J, H)$ is a
   disjoint union of $\mathcal{P}_I(a,x,y; (\mathbf{B}' \#
   \mathbf{B}'',j),J, H)$ for $1 \leq j \leq l'+l''-1$.)

   We shall perform the gluing in $\widetilde{M} = D \times M$. Let
   $\widetilde{J}_H$ be the almost complex structure in $D \times M$
   associated to $J$ and $H$ (see Section~\ref{Sb:hampert}).  Put
   $\widetilde{L} = \partial{D} \times L$. Let $\widetilde{u}_1(z) =
   (z, v'_{l'}(z))$ be the graph of $v'_{l'}$ and let $u_0(z) = (1,
   v''_1(z))$ be a copy of $v''_1$ lying in the fibre over $1 \in D$.
   Clearly both $u_1, u_0$ are $\widetilde{J}_H$-holomorphic disks
   with boundary on $\widetilde{L}$ and $u_1(1) = u_0(-1)$. Moreover
   $A_1=[u_1] = B'_{l'}+[D] \in H_2(\widetilde{M}, \widetilde{L})$,
   $A_0=[u_0]=B''_1 \in H_2(\widetilde{M}, \widetilde{L})$.

   Again we take $W$ to be a product $W = W_{-} \times W_1 \times W_0
   \times W_{+}$. The manifolds $W_{\pm}$ and the maps $h_{\pm}$ are
   defined in a similar way as in the proof of
   Corollaries~\ref{C:glue-pearls-1},~\ref{C:glue-pearls-2} (only that
   now $h_{-}$ maps $W_{-}$ to $\{-1\} \times L$ and $h_{+}$ maps
   $W_{+}$ to $\{ 1 \} \times L$). As before, we have: $\dim W_{-} =
   \ind_f(x) + \sum_{i=1}^{l'-1} \mu(B'_i)$, $\dim W_{+} = n-\ind_f(y) +
   \sum_{i=2}^{l''} \mu(B''_i)$.

   Next, take $W_1$ to be $\{ 0 \} \times W_a^u \subset \widetilde{M}$
   and $W_0 = \widetilde{M}$. We take $h_i: W_i \to \widetilde{M}$ to
   be the inclusions. Note that
   \begin{equation} \label{Eq:dimW-H}
      \dim W = \ind_h(a)+\ind_f(x)-\ind_f(y) +
      \sum_{i=1}^{l'-1} \mu(B'_i)+ \sum_{i=2}^{l''} \mu(B''_i) +
      3n+2.
   \end{equation}

   As before, by taking $J$ and $H$ to be generic we may assume that
   the assumptions of Theorem~\ref{T:gluing} are satisfied for $u_1,
   u_0$. We obtain from this theorem a path $(\widehat{v}_s, a(s),
   q_s)$, where $\widehat{v}_s:(D, \partial{D}) \to (\widetilde{M},
   \widetilde{L})$ is in the class $A_1+A_0$ and satisfies:
   \begin{itemize}
     \item $(\widehat{v}_s(-1), \widehat{v}_s(-a(s)),
      \widehat{v}_s(a(s)), \widehat{v}_s(1)) = q_s \in W$.
     \item $\widehat{v}_s$ with the marked points $(-1, -a(s), a(s),
      1)$ converges, as $s \to \infty$, to $(u_1, u_0)$ with the
      marked points $(-1, 0)$, $(0, 1)$ in the Gromov topology.
   \end{itemize}
   In order to extract from $\widehat{v}_s$ a $(J,H)$-holomorphic disk
   write $\widehat{v}_s (z) = (\varphi_s(z), v_s(z))$, where
   $\varphi_s:(D,\partial{D}) \to (D, \partial{D})$ and $v_s:(D,
   \partial{D}) \to (M,L)$.  Note that $$\varphi_s(-1)=-1, \quad
   \varphi_s(1)=1, \quad \varphi_s(-a(s))=0.$$
   As $\widehat{v}_s$ is
   $\widetilde{J}_H$-holomorphic, $pr_D: (\widetilde{M},
   \widetilde{J}_H) \to (D,i)$ is holomorphic, and
   ${pr_{D}}_{_*}[\widehat{v}_s] = [D] \in H_2(D,\partial{D})$ it
   follows that $\varphi_s \in \textnormal{Aut}(D)$. Put $w_s=v_s
   \circ \varphi_s^{-1}$. Then $\widetilde{w}_s(z)=(z,w_s(z))$ is a
   $\widetilde{J}_H$-holomorphic section of $D \times M \to D$ hence
   by Proposition~\ref{P:JH} $w_s$ is $(J,H)$-holomorphic. Also note
   that
   $$w_s(\pm 1) = v_s(\pm 1), \quad w_s(0) = v_s(-a(s)) \in W_a^u.$$
   As in the proof of Corollary~\ref{C:glue-pearls-1} we extract from
   $q_s$ a path of chains of $J$-holomorphic disks $(u_{1,s}, \ldots,
   u_{k'-1,s}, u_{k'+1,s}, \ldots, u_{l'+l''-1,s})$ such that when we
   insert the disk $w_s$ to the $k'$'th entry we obtain
   $$\mathbf{u}_s = (u_{1,s}, \ldots, u_{k'-1,s}, w_s, u_{k'+1,s},
   \ldots, u_{l'+l''-1,s}) \in \mathcal{P}_I(a,x,y; (\mathbf{B}' \#
   \mathbf{B}'',k'), J, H),$$
   and $\mathbf{u}_s$ converges in the
   Gromov topology as $s \to \infty$ to $(\mathbf{v}', \mathbf{v}'')$.
   This concludes the proof of the existence statement.

   The uniqueness statement follows from Proposition~\ref{T:gluing} since
   by~\eqref{Eq:dimW-H} we have:
   $$\mu(A_1)+ \mu(A_0) + \dim W - 5 \dim \widetilde{M} = \ind_h(a) +
   \ind_f(x)-\ind_f(y) + \mu(\mathbf{B}')+\mu(\mathbf{B}'')-2n-1= 0.$$
\end{proof}

\subsubsection{Gluing a trajectory of pearls to a trajectory
  with external constrains II} \label{Sbsb:glue-ex-2II}

Consider the following situation encountered in the proof of
Proposition~\ref{prop:module_str} in Section~\ref{S:qm}. Let $f:L \to
\mathbb{R}$, $h', h'': M \to \mathbb{R}$ be Morse functions and
$\rho_L, \rho'_M, \rho''_M$ be Riemannian metrics on $L$, $M$. Let
$\mathbf{B}'=(B'_1, \ldots, B'_{l'})$, $\mathbf{B}''=(B''_1, \ldots,
B''_{l''})$ be two vectors of non-zero classes and $1 \leq k' \leq
l'$, $1 \leq k'' \leq l''$. Given $x,y \in \textnormal{Crit}(f)$, $a'
\in \textnormal{Crit}(h')$, $a'' \in \textnormal{Crit}(h'')$ define
$\widehat{\mathcal{P}}_{V}(a',a'', x, y, (\mathbf{B}',k'),
(\mathbf{B}'', k''); J)$ to be the space of all $(u'_1, \ldots,
u'_{l'}, u''_1, \ldots, u''_{l''})$ such that:
\begin{itemize}
  \item $u'_i \in \mathcal{M}(B'_i, J)$ for every $1\leq i\leq l'$,
   $u''_j \in \mathcal{M}(B''_j,J)$ for every $1\leq j\leq l''$.
  \item $u'_1(-1) \in W_x^u$, $u''_{l''}(1) \in W_y^s$.
  \item $(u'_i(1), u'_{i+1}(-1)) \in Q_{(f,\rho_L)}$ for every $1\leq
   i \leq l'-1$, $(u''_j(1), u''_{j+1}(-1)) \in Q_{(f,\rho_L)}$ for
   every $1\leq j \leq l''-1$.
  \item $u'_{l'}(1) = u''_1(-1)$.
  \item $u'_{k'}(0) \in W_{a'}^u$, $u''_{k''}(0) \in W_{a''}^u$.
\end{itemize}
Define now $$\mathcal{P}_{V}(a',a'', x, y, (\mathbf{B}',k'),
(\mathbf{B}'', k''); J) = \widehat{\mathcal{P}}_{V}(a',a'', x, y,
(\mathbf{B}',k'), (\mathbf{B}'', k''); J) / \mathbf{G}' \times
\mathbf{G}'',$$
where $\mathbf{G}'$ is $G_{-1,1}^{\times l'}$ with the
$k'$'th factor replaced by the trivial group, and $\mathbf{G}''$ is
defined in a similar way.

Similarly, for a vector of non-zero classes $\mathbf{C} = (C_1,
\ldots, C_l)$ and $1\leq p' < p'' \leq l$ consider the space of all
$(u_1, \ldots, u_l)$ such that:
\begin{enumerate}[(i)]
  \item $u_i \in \mathcal{M}(C_i,J)$ for every $1 \leq i \leq l$.
  \item $u_1(-1) \in W_x^u$, $u_l(1) \in W_y^s$.
  \item $(u_i(1), u_{i+1}(-1)) \in Q_{(f,\rho_L)}$ for every $1 \leq i
   \leq l-1$.
  \item $u_{p'}(0) \in W_{a'}^u$, $u_{p''}(0) \in W_{a''}^u$.
\end{enumerate}
The group $\mathbf{G}_{-1,1}^{\times (l-2)}$ acts on this space by
reparametrizations. Denote the quotient space by
$\mathcal{P}_{IV}(a',a'', x,y; (\mathbf{C},p',p'');J)$.

Finally, let $1 \leq p \leq l$. Consider the space of all $(u_1,
\ldots, u_l, r)$ such that:
\begin{itemize}
  \item $(u_1, \ldots, u_l)$ satisfy conditions (i)-(iii) above.
  \item $r \in (0,1)$.
  \item $u_p(-r) \in W_{a'}^u$, $u_p(r) \in W_{a''}^u$.
\end{itemize}
The group $\mathbf{G}_{-1,1}^{\times (l-1)}$ acts on this space by
reparametrizations. Denote the quotient space by
$\mathcal{P}_{IV}(a',a'', x,y; (\mathbf{C},p);J)$.

\begin{cor} \label{C:glue-pearls-3}
   Let $(f,\rho_L)$, $(h', \rho'_M)$, $(h'', \rho''_M)$ be as above
   where the triple of metrics $\rho_L, \rho'_M, \rho''_M$ is assumed
   to be generic.  There exists a second category subset
   $\mathcal{J}_{\textnormal{reg}} \subset \mathcal{J}(M,\omega)$ such
   that for every $J \in \mathcal{J}_{\textnormal{reg}}$, every pair
   of vectors of non-zero classes $\mathbf{B}', \mathbf{B}''$ and
   every $x,y\in \textnormal{Crit}(f)$, $a' \in \textnormal{Crit}(h')$
   $a'' \in \textnormal{Crit}(h'')$ with $\ind_{h'}(a') +
   \ind_{h''}(a'') +
   \ind_f(x)-\ind_f(y)+\mu(\mathbf{B}')+\mu(\mathbf{B}'')-4n=0$ the
   following holds. Let $1 \leq k' \leq l'$, $1 \leq k'' \leq l''$ and
   $(\mathbf{v}', \mathbf{v}'') \in \mathcal{P}_{V}(a',a'', x,y,
   (\mathbf{B}',k'),(\mathbf{B}'', k'');J)$.
   \begin{enumerate}
     \item If $(k',k'') \neq (l',1)$ then there exists a path $\{
      \mathbf{u}_s \} \subset \mathcal{P}_{IV}(a',a'',
      x,y,(\mathbf{B}' \# \mathbf{B}'', k', l'+k''-1); J)$ which
      converges in the Gromov topology as $s \to \infty$ to
      $(\mathbf{v}',\mathbf{v}'')$.
     \item If $k'=l'$, $k''=1$ then there exists a path $\{
      (\mathbf{u}_s, a(s)) \} \subset \mathcal{P}_{IV}(a',a'',
      x,y,(\mathbf{B}' \# \mathbf{B}'', k'); J)$ such that:
      \begin{enumerate}
        \item $\mathbf{u}_s$ converges in the Gromov topology as $s
         \to \infty$ to $(\mathbf{v}',\mathbf{v}'')$.
        \item $a(s) \xrightarrow[s \to \infty]{} 1$.
        \item The $k'$'th disk, $u_{k',s}$, in $\mathbf{u}_s$,
         converges with the marked points $(-1, -a(s), a(s), 1)$, as
         $s \to \infty$ to $(v'_{k'}, v''_1)$ with the marked points
         $(-1, 0)$, $(0,1)$, in the Gromov topology.
      \end{enumerate}
   \end{enumerate}
   Moreover, in both cases~(1) and~(2) above the end of the
   $1$-dimensional manifold $\mathcal{P}_{IV}$ parametrized by the
   path $\{ \mathbf{u}_s \}$ (resp. $\{ (\mathbf{u}_s, a(s)) \}$) is
   unique in the sense that every other path satisfying the same
   properties for $s \to \infty$, lies in the same end.
\end{cor}

\begin{proof}
   The case $(k',k'') \neq (l',1)$ is proved in a similar way to
   Corollaries~\ref{C:glue-pearls-1} and~\ref{C:glue-pearls-2}.

   As for the case $k'=l'$, $k''=1$, the manifolds $W_{\pm}$ are
   defined in a similar way as in the proofs of
   Corollaries~\ref{C:glue-pearls-1},~\ref{C:glue-pearls-2} and we
   take $W_1=W_{a'}^u$, $W_0=W_{a''}^u$.
\end{proof}

\subsection{Overview of the proof of Theorem~\ref{T:gluing}}
\label{Sb:overivew}
The proof of Theorem~\ref{T:gluing} is built from the following steps.
We first transform the domains of $u_1,u_0$ to the strip $S=\mathbb{R}
\times [0,1]$. The main reason for this is convenience, especially
when performing translations along the $\mathbb{R}$-axis. In
Section~\ref{S:analytic} we introduce the analytic setup for
performing the gluing. In particular we shall have to work with
weighted Sobolev spaces in order to make the linearization $D$ of the
non-linear $\overline{\partial}$ operator Fredholm.

The second step is usually called pregluing. Here we build an
approximate solutions $u_R$ of the $\overline{\partial}$ equation
depending on a parameter $R \gg 0$. The $u_R$'s are glued from $u_1,
u_0$ using partition of unity and coincide with suitable (larger and
larger) translates of $u_1, u_0$ near the ends of $S$. The $u_R$'s are
approximate solutions in the sense that $\overline{\partial}{u_R}$
becomes smaller and smaller as $R \to \infty$ in a suitable norm.  The
pregluing and the needed estimates of the $u_R$'s are done in
Section~\ref{S:pregluing}. In order for certain operators related to
$u_R$ to become uniformly bounded we shall have to deform our weighted
norms with $R$. These norms are introduced in
Section~\ref{S:pregluing}. For the reader convenience we included in
Section~\ref{S:auxiliary} a Sobolev-type inequality that will be used
frequently in the proof.

The third step is to construct a right inverse to the linearization
$D_{u_R}$ of the $\overline{\partial}$ operator at $u_R$. More
precisely we construct a family of operators $\{Q_R\}_{R \gg 0}$ such
that $D_{u_R} \circ Q_R = \Id$ for every $R \gg 1$ and such that the
$Q_R$'s are {\em uniformly bounded}. This is done in
Section~\ref{S:rinv}.

The next step is to use an implicit function theorem which will
correct the approximate solutions $u_R$ to genuine $J$-holomorphic
solutions $v_R$. This is carried out in Section~\ref{S:implicit}.  The
operators $Q_R$'s are needed for the implicit function theorem to
work.

In Section~\ref{S:implicit} we also prove that the $v_R$'s verify the
marked points conditions~\eqref{Eq:MB} and that $v_R$ converges as $R
\to \infty$ to $(u_1,u_0)$ with marked points in the Gromov topology.

The steps above prove the existence statement of
Theorem~\ref{T:gluing}. The final step is devoted to proving the
uniqueness statement of Theorem~\ref{T:gluing}. This is usually called
in the ``gluing literature'' {\em surjectivity of the gluing map}.  It
occupies the rest of Section~\ref{S:implicit}.

Throughout the proof of the existence part we essentially follow the
work of Fukaya-Oh-Ohta-Ono~\cite{FO3}, however we occasionally use
slightly different notation and normalizations.

\subsection{Analytic setting}  \label{S:analytic}

From now on we shall identify the disk $D$ with the compactified
strip $\widehat{S}=S \cup \{-\infty, \infty\}$ where $S=\mathbb{R}
\times [0,1] \approx \{z \in \mathbb{C}|\, 0\leq \textnormal{Im}z\leq
1\}$.  This is done via the biholomorphism
$$\lambda:\widehat{S} \to D, \quad \lambda(z) = \frac{e^{\pi
    z}-i}{e^{\pi z}+i}.$$
Note that $\lambda(\infty)=1,
\lambda(-\infty)=-1$.  Denote $S[a,b] = [a,b] \times [0,1] \subset S$.
Similarly we have $S[a,\infty)$ etc.  The reason for this change of
coordinates is that it is more handy for using translations. The price
we pay for this is that the Banach space norms the are well suited to
the strip $S$ are not conformally invariant. Moreover, our elliptic
boundary problem becomes an ``asymptotic boundary'' problem at the
ends of $S$. In order to control this asymptotic we shall need to
endow our Banach spaces with some weighted norms.

Denote by $g_{\omega,J}(\cdot,\cdot)=\omega(\cdot,J\cdot)$ be the
Riemannian metric associated to $(\omega,J)$. Choose a new metric
$g_{\omega,J,L}$ on $M$ for which $L$ is totally geodesic and which
coincides with $g_{\omega,J}$ outside a small neighbourhood of $L$.
{\em Henceforth, we shall use the metric $g_{\omega,J,L}$ as our main
  Riemannian metric rather than $g_{\omega,J}$.}  Thus all pointwise
norms, connections, and distances are to be understood with respect to
$g_{\omega,J,L}$, unless explicitly stated otherwise.

Fix $\delta>0$ small enough (see remark below) and $p>2$.
\begin{dfn} \label{D:banach-manifold}
   Denote by $W^{1,p;\delta}(M,L)$ the space of all pairs $(u,
   \underline{p})$ where:
   \begin{enumerate}
     \item $u:(S,\partial S) \to (M,L)$ is of class
      $W^{1,p}_{\textnormal{loc}}$.
     \item $\underline{p}=(p_{-\infty},p_{\infty}), p_{\pm \infty} \in
      L$.
     \item $\int_{S[0,\infty)} e^{\delta|\tau|}\bigl(
      \textnormal{dist}(u(\tau,t),p_{\infty})^p +
      |du_{(\tau,t)}|^p\bigr)d\tau dt < \infty$, \\
      $\int_{S(-\infty,0]} e^{\delta|\tau|}\bigl(
      \textnormal{dist}(u(\tau,t),p_{-\infty})^p +
      |du_{(\tau,t)}|^p\bigr)d\tau dt < \infty.$
   \end{enumerate}
\end{dfn}
\begin{rem} \label{R:analytic}
   \begin{enumerate}
     \item As $p > 2 = \dim S$ all $u \in W^{1,p;\delta}(M,L)$ are
      actually continuous.
     \item Standard arguments show that for every $(u,\underline{p})
      \in W^{1,p;\delta}(M,L)$ we have $$\lim_{\tau \to \pm \infty}
      u(\tau,t) = p_{\pm \infty}\quad \textnormal{uniformly in} \,\,
      t.$$
      Therefore the $\underline{p}$ in $(u,\underline{p})$ is
      superfluous as $\underline{p}$ can be recovered from $u$.
     \item A simple computation shows that there exists $\delta_0>0$
      such that for every smooth map $u_D:(D,\partial{D}) \to (M,L)$,
      the map $u=u_D \circ \lambda:(S,\partial{S}) \to (M,L)$ belongs
      to $W^{1,p;\delta}(M,L)$ for every $0<\delta<\delta_0$. This
      follows from the fact that there exists a constant $C$ such that
      $|\lambda'(R+it)| \leq C e^{-\pi |R|}$ for $|R| \gg 0$.
      \label{Rn:smooth-u}
     \item H\"{o}lder's inequality implies that each $u \in
      W^{1,p;\delta}(M,L)$ has finite energy $E(u)= \int_S u^*\omega =
      \tfrac{1}{2}\int_{S} |du|_{g_{\omega,J}}^2 < \infty$.
      \label{Rn:energy}
   \end{enumerate}
\end{rem}

We shall now endow $W^{1,p;\delta}(M,L)$ with a structure of a Banach
manifold. For this end fix $r_M>0$ small enough so that around every
$y \in L$ there exists a geodesic ball $B_{y}(r_M) \subset M$ of
radius $r_M$. Define a map $P:T(L) \to \textnormal{Vector fields}(M)$
in the following way: given $y \in L$, $v\in T_{y}(L)$, define the
vector field $P(v)$ by:
\begin{equation} \label{Eq:pal}
   P(v)(x)=
   \begin{cases}
      \chi(\textnormal{dist}(y,x))\textnormal{Pal}_x(v) &
      x \in B_{y}(r_M) \\
      0 & x \notin B_{y}(r_M)
   \end{cases}
\end{equation}
Here $\chi:[0,\infty) \to [0,1]$ is a cutoff function which equals $1$
on $[0, r_M/3]$ and $0$ on $[r_M/2, \infty)$.  $\textnormal{Pal}_x(v)$
stands for parallel transport of $v$ along the minimal geodesic that
connects $y$ to $x$.  Given a map $u:S \to M$ and $v \in T_{y}(L)$,
define $P_u(v):S \to T(M)$ by $P_u(v)(\tau,t) = P(v)(u(\tau,t))$.

Put $s_{\delta}(\tau) = e^{\delta|\tau|}$. Let $\xi \in
\Gamma(u^*T(M)$ be a section for which the limits $\xi_{\pm
  \infty}=\lim_{\tau \to \pm \infty} \xi(\tau,t)$ exist independently
of $t$. Define:
\begin{align*}
   \Vert \xi \Vert_{1,p;s_{\delta}}^p = & \int_{S(-\infty,0]}
   s_{\delta}(\tau) \Bigl( |\xi-P_u(\xi_{-\infty})|^p + |\nabla(\xi -
   P_u(\xi_{-\infty}))|^p \Bigr) d\tau dt \\
   & + \int_{S[0,\infty)} s_{\delta}(\tau) \Bigl(
   |\xi-P_u(\xi_{\infty})|^p + |\nabla(\xi - P_u(\xi_{\infty}))|^p
   \Bigr) d\tau dt + |\xi_{-\infty}|^p + |\xi_{\infty}|^p.
\end{align*}
Here and in what follows $\nabla$ stands for the Levi-Civita
connection of our metric $g_{\omega,J,L}$.

\begin{prop}[See~\cite{FO3}] \label{P:tangent-space}
   $W^{1,p;\delta}(M,L)$ is a Banach manifold with respect to the norm
   $\Vert \cdot \Vert_{1,p;s_{\delta}}$ on its tangent spaces.  Its
   tangent space $T_u^{1,p;\delta} = T_u(W^{1,p;\delta}(M,L))$ at $u$
   consists of all sections $\xi \in \Gamma(u^*T(M))$ such that:
   \begin{enumerate}
     \item $\xi \in W^{1,p}_{\textnormal{loc}}$.
     \item $\xi(x) \in T_{u(x)}(L)$ for every $x \in \partial S$.
     \item $\lim_{\tau \to \pm \infty} \xi(\tau,t)$ converges to a
      vector $\xi_{\pm \infty} \in T_{p_{\pm \infty}}(L)$
      independently of $t$, where $p_{\pm \infty} = \lim_{\tau \to \pm
        \infty} u(\tau,t)$.
     \item $\Vert \xi \Vert_{1,p;s_{\delta}} < \infty$.
   \end{enumerate}
\end{prop}
Again, note that since $p > 2 = \dim S$ all $\xi \in T^{1,p;\delta}_u$
are actually continuous.

Given $u:(S,\partial S) \to (M,L)$ we also define the Banach space
$\mathcal{E}_u^{0,p;\delta}$ consisting of all
$L^p_{\textnormal{loc}}$-section $\eta \in \Gamma(\Lambda^{0,1}(S)
\otimes u^*T(M))$ such that: $$\Vert \eta \Vert^p_{0,p;s_{\delta}} =
\int_S s_{\delta}(\tau) | \eta(\tau,t) |^p d\tau dt < \infty.$$

\begin{remnonum}
   If $0 < \delta_1 < \delta_2$ then $W^{1,p;\delta_1}(M,L) \supset
   W^{1,p;\delta_2}(M,L)$, and $\Vert \cdot \Vert_{1,p;s_{\delta_1}}
   \leq \Vert \cdot \Vert_{1,p;s_{\delta_2}}$, $\Vert \cdot
   \Vert_{0,p;s_{\delta_1}} \leq \Vert \cdot
   \Vert_{0,p;s_{\delta_2}}$.  In particular for every $u \in
   W^{1,p;\delta_2}(M,L)$ we have $T_u^{1,p;\delta_1} \supset
   T_u^{1,p;\delta_2}$, $\mathcal{E}_u^{0,p;\delta_1} \supset
   \mathcal{E}_u^{0,p;\delta_2}$.
\end{remnonum}

In what follows we shall need to decrease the size of $\delta$ several
times for various estimates to hold.  Nevertheless we shall do so only
a finite number of times and the final range of admissible $\delta$'s
will depend only on $(M,L,\omega,J,g_{\omega,J,L})$ and on our initial
$J$-holomorphic disks $u_1,u_0$.

\subsubsection{The linearization of the $\overline{\partial}$
  operator} \label{Sb:lin-dbar}

In order to linearize the $\overline{\partial}$ operator at an
arbitrary $u \in W^{1,p;\delta}(M,L)$ we follow~\cite{McD-Sa:Jhol-2}.
For this purpose we need to introduce a different connection
$\widetilde{\nabla}$ on $T(M)$ as follows. Let $\nabla'$ be the
Levi-Civita connection of the metric $g_{\omega,J}$ (this is the only
instance where we use the metric $g_{\omega,J}$ instead of
$g_{\omega,J,L}$).  Define $\widetilde{\nabla}_v X = \nabla'_v X -
\tfrac{1}{2}J(\nabla'_v J)X$.  Note that $\widetilde{\nabla}$
preserves $J$. Define
$$\mathcal{F}_u: T_u^{1,p;\delta} \to \mathcal{E}_u^{0,p;\delta},
\quad \mathcal{F}_u(\xi) = \Phi_u(\xi)^{-1}
\overline{\partial}_J(\exp_u \xi),$$
where $\Phi_u(\xi): u^*T(M) \to
\exp_u(\xi)^*T(M)$ is defined using parallel transport with respect to
$\widetilde{\nabla}$.  (See~\cite{McD-Sa:Jhol-2} for more details.)
However, in contrast to~\cite{McD-Sa:Jhol-2}, the $\exp$ is defined
here with respect to our main metric $g_{\omega,J,L}$. With the above
notation the linearization of $\overline{\partial}$ at $u$ is the
linearization of $\mathcal{F}_u$ at $\xi=0$, i.e. the operator
$$D_u := d \mathcal{F}_u (0): T_u^{1,p;\delta} \to
\mathcal{E}_u^{0,p;\delta}.$$
We have the following expression for
$D_u$ (see~\cite{McD-Sa:Jhol-2}):
$$D_u \xi = \frac{1}{2}(\nabla' \xi + J(u)\nabla' \xi \circ j) -
\frac{1}{2} J(u)(\nabla'_{\xi}J) \partial_Ju,$$
where $j$ is the
standard complex structure on $S$.  Note that when $u$ is
$J$-holomorphic $D_u$ does not depend on the choice of the
connections. Therefore assumption~\ref{A:T1}-(\ref{I:T1-reg-01}) makes
sense independently of any metric.

\subsubsection{Fredholm property}  \label{Sb:fredholm-1}
Without the weights in the norms (i.e. for $\delta=0$), the operator
$D_u$ is in general not Fredholm. The reason for this is that $D_u$
has degenerate asymptotic at the ends of $S$. The weighted norms
$\Vert \cdot \Vert_{1,p;s_{\delta}}$, $\Vert \cdot
\Vert_{0,p;s_{\delta}}$ are introduced mainly in order to rectify this
problem. Let us briefly explain how the weights are related to
Fredholmness. To distinguish between $\delta>0$ and $\delta=0$ denote
by $\mathcal{D}_u^0$, $\mathcal{D}_u^{\delta}$ the operator $D_u$ on
the spaces $T_u^{1,p;0}$, $T_u^{1,p;\delta}$ respectively.  (Here
$1,p;0$ means that we take $\delta=0$, i.e.  there is no weight.)  Let
$\mathcal{H}_u^{1,p;0} \subset T_u^{1,p;0}$,
$\mathcal{H}_u^{1,p;\delta} \subset T_u^{1,p;\delta}$ be the subspaces
consisting of all $\xi$ with $\xi(\pm \infty)=0$. Note that
$\mathcal{H}_u^{1,p;0}, \mathcal{H}_u^{1,p;\delta}$ have finite
codimension, hence Fredholmness is not affected. Define
\begin{align*}
   & \Theta_{\delta}: \mathcal{H}_u^{1,p;0} \to
   \mathcal{H}_u^{1,p;\delta},\quad
   \Theta_{\delta}(\xi)(\tau,t) = e^{-\delta|\tau|/p}\xi(\tau,t), \\
   & \rho_{\delta}:\mathcal{E}_u^{0,p;0} \to
   \mathcal{E}_u^{0,p;\delta}, \quad \rho_{\delta}(\eta)(\tau,t) =
   e^{-\delta|\tau|/p}\eta(\tau,t).
\end{align*}
These maps are bounded isomorphisms between the corresponding Banach
spaces. A simple computation shows that on $S[0,\infty)$ we have:
$$(\rho_{\delta}^{-1} \circ \mathcal{D}^{\delta}_u \circ
\Theta_{\delta}) \xi = \mathcal{D}^0_u \xi + \tfrac{\delta}{p} (d \tau
\otimes \xi - dt \otimes J\xi ).$$
An analogous formula holds for
$S(-\infty, 0]$. In other words, we get a perturbation of
$\mathcal{D}^0_u$ by a small $0$-order operator. This perturbation
term gives non-degenerate asymptotic at the ends of $S$ and so for
generic $\delta$ (and under certain assumptions on $u$) the operator
$(\rho_{\delta}^{-1} \circ \mathcal{D}^{\delta}_u \circ
\Theta_{\delta})$ becomes Fredholm.  Thus $\mathcal{D}^{\delta}_u$
will be Fredholm too.  Moreover for $\delta>0$ small enough the index
of this operator does not depend on $\delta$. We shall not use the
notation $\mathcal{D}^{\delta}_u$ anymore and continue to denote by
$D_u$ our operator defined on $T_u^{1,p;\delta}$ for $\delta>0$.

\subsubsection*{Conditions which assure Fredholmness}
Assume that $u \in W^{1,p;\delta}(M,L)$ has exponential decay at the
ends of $S$, (i.e.  $|du_{(\tau,t)}| \leq C e^{-c|\tau|}$ when $|\tau|
\gg 0$ for some constants $C,c>0$). Then there exists
$\delta_{\textnormal{Fred}}(u)>0$ such that for every $0< \delta \leq
\delta_{\textnormal{Fred}}(u)$, the operator $D_u$ is Fredholm.
Moreover for such $\delta$'s, $\textnormal{index}(D_u)=\dim L +
\mu([u])$ and the kernel and cokernel of $D_u$ are independent of
$\delta$.  The proofs of these statements follow in a straightforward
way from the theory developed in
e.g.~\cite{Fl:unregularized,Fl:relindex,Ro-Sa:spectral,Sc:thesis}.

An important situation in which we have exponential decay (hence also
Fredholmness for $0<\delta \ll 1$) is when $u$ is $J$-holomorphic near
the ends of $S$.  To see this recall that $u \circ \lambda^{-1}: (D
\setminus \{-1,1\}, \partial D \setminus \{-1,1\}) \to (M,L)$ has
finite energy (see Remark~\ref{R:analytic}-(\ref{Rn:energy})). Then by
the removal of boundary singularities theorem of Oh~\cite{Oh:rsing},
$u \circ \lambda^{-1}$ extends smoothly to $D$. It follows that $u$
has exponential decay at the ends of $S$. (See
Remark~\ref{R:analytic}-(\ref{Rn:smooth-u}).)

\subsubsection{A simplification} \label{Sb:gluing-simplification}
In order to simplify the notation we shall henceforth restrict
ourselves to the special case when $W=W_{-} \times W_1 \times W_0
\times W_{+}$ where $W_{\pm}$ are {\em sumbanifolds} of $L$, $W_1,
W_0$ are {\em submanifolds} of $M$, and furthermore the map
$\mathbf{h}$ is the inclusion $W_{-} \times W_1 \times W_0 \times
W_{+} \hooklongrightarrow L \times M \times M \times L$.  The proof of
the general case is not more complicated from the analytic point of
view, however the notation becomes heavier.

In view of the above simplification we denote from now on
$$k_{\pm} = \dim W_{\pm}, \qquad k_i = \dim W_i, \quad i=1,0.$$
Since
$\mathbf{h}$ is now the inclusion we can write elements of
$\mathcal{M}(A,J; \mathcal{C}(\mathbf{h}))$ as pairs $(u,r)$ (instead
of triples $(u,r,q)$) since $q=(u(-1), u(-r), u(r), u(1))$ is
determined by $(u,r)$.

Finally, note that the third point in Assumption~\ref{A:T1} can be now
simplified to {\sl ``$ev$ is transverse to $W_{-} \times W_1 \times
  \textnormal{diag}(L) \times W_0 \times W_{+}$ at the point $(u_1,
  u_0) \in \mathcal{M}(A_1, J) \times \mathcal{M}(A_0, J)$''.}

\subsubsection{Assumption~\ref{A:T1} in the new coordinates}
 \label{Sb:assumption-new} We return to our $J$-holomorphic
 disks $u_1, u_0$. Put $\tau_1=\textnormal{Re}\, \lambda^{-1}(z_1)$,
 $\tau_0=\textnormal{Re}\, \lambda^{-1}(z_0)$.  Note that
 $\textnormal{Im}\, \lambda^{-1}(z_i)=\tfrac{1}{2}$ and
 $\lambda^{-1}(\pm 1) = \pm \infty$. Replace $u_1,u_0:(D,\partial D)
 \to (M,L)$ by $u_1 \circ \lambda, u_0 \circ
 \lambda:(\widehat{S},\partial\widehat{S}) \to (M,L)$, respectively.
 We claim that for $\delta>0$ small enough Assumption~\ref{A:T1}
 implies the following one:
\begin{assumption} \label{A:T2}
   \begin{enumerate}
     \item $u_1(\infty)=u_0(-\infty)$.
     \item $u_1(-\infty) \in W_{-}, u_1(\tau_1, \tfrac{1}{2}) \in W_1,
      u_0(\tau_0, \tfrac{1}{2}) \in W_0, u_0(\infty) \in W_{+}$.
     \item $J$ is regular for both $u_1$ and $u_0$ in the sense that
      the linearizations $D_{u_1}, D_{u_0}$ of the
      $\overline{\partial}$ operator at $u_1,u_0$ are surjective.
      \label{A:T2-Du}
     \item The evaluation map $ev: \mathcal{M}(A_1,J) \times
      \mathcal{M}(A_0,J) \to L \times M \times L \times L \times M
      \times L$, $$ev(v_1,v_0)=\big( v_1(-\infty),
      v_1(\tau_1,\tfrac{1}{2}), v_1(\infty), v_0(-\infty), v_0(\tau_0,
      \tfrac{1}{2}), v_0(\infty) \big)$$
      is transverse to $W_{-}\times
      W_1 \times \textnormal{diag}(L) \times W_0 \times W_{+}$ at the
      point $(u_1, u_0)$. (See the remarks in
      Section~\ref{Sb:gluing-simplification}.)
   \end{enumerate}
\end{assumption}
That property~(\ref{A:T2-Du}) follows from the corresponding property
in Assumption~\ref{A:T1} is not completely obvious. This is because
the Banach spaces which are the domains and targets of $D_{u_1},
D_{u_0}$ change when we pass from the disk $D$ to the strip $S$.  The
other properties follow trivially from Assumption~\ref{A:T1}.
\begin{proof}[Proof that Assumption~\ref{A:T1} $\Longrightarrow$
   Assumption~\ref{A:T2}-(\ref{A:T2-Du})]

   Let $u_D:(D,\partial D) \to (M,L)$ be a $J$-holomorphic disk and
   denote $u=u_D \circ \lambda$. Let $T_{u_D}^{1,p}$ be the space of
   $W^{1,p}$-sections $\xi_D$ of the bundle $u_D^* T(M) \to D$ which
   satisfy $\xi_D(x) \in T(L)$ for every $x \in \partial D$. Put
   $$\mathcal{E}_{u_D}^{0,p} = L^p\bigl( \Lambda^{0,1}(D)\otimes u_D^*
   T(M) \bigr).$$
   We have to prove that if $D_{u_D}: T_{u_D}^{1,p} \to
   \mathcal{E}_{u_D}^{0,p}$ is surjective then for $\delta>0$ small
   enough the operator $D_{u}: T_{u}^{1,p;\delta} \to
   \mathcal{E}_u^{0,p;\delta}$ is surjective too.

   Choose $0< \delta < \pi p$ small enough so that $D_{u}$ is
   Fredholm. As $D_u$ has closed image and the subspace of smooth
   compactly supported sections $C_0^{\infty}\bigl(
   \Lambda^{0,1}(S)\otimes u^* T(M) \bigr) \subset
   \mathcal{E}_u^{0,p;\delta}$ is dense it is enough to show that this
   subspace lies in the image of $D_u$.

   Let $\eta \in \mathcal{E}_u^{0,p;\delta}$ be a smooth compactly
   supported section. Put $\eta_D = (\lambda^{-1})^* \eta \in
   \mathcal{E}_{u_D}^{0,p}$. By assumption~\ref{A:T1} there exists
   $\xi_D \in T_{u_D}^{1,p}$ such that $D_{u_D} \xi_D = \eta_D$.  Put
   $\xi = \xi_D \circ \lambda$. Then
   $$D_u \xi = D_u (\xi_D \circ \lambda) = \lambda^* D_{u_D} \xi_D =
   \lambda^* \eta_D = \eta.$$
   It remains to prove that $\xi \in
   T_u^{1,p;\delta}$, i.e. $\Vert \xi \Vert_{1,p;s_{\delta}} <
   \infty$.  To see this note that by elliptic regularity $\xi_D:D \to
   u_D^* T(M)$ is smooth. It follows that for $\tau \gg 0$ we have
   \begin{align} \label{Eq:xi-estim-1}
      & | \xi(\tau,t) - P_u \xi(\infty) | = | \xi_D \circ
      \lambda(\tau+it) -
      \textnormal{Pal}_{u_D \circ \lambda(\tau+it)} \xi_D(1)| \\
      & \leq C_1(\xi_D)|\lambda(\tau+it)-1| \leq C_1'(\xi_D) e^{-\pi
        \tau}, \notag
   \end{align}
   where the constants $C_1(\xi_D), C_1'(\xi_D)$ depend on $\xi_D$.
   Similarly, for $\tau \gg 0$ we have:
   \begin{equation} \label{Eq:xi-estim-2}
      | \nabla (\xi(\tau,t) - P_u \xi(\infty)) | \leq C_2(\xi_D) |
      \lambda(\tau+it)-1| \leq C_2'(\xi_D) e^{-\pi \tau},
   \end{equation}
   for some constants $C_2(\xi_D),C_2'(\xi_D)$ that depend on $\xi_D$.
   Analogous estimates to~\eqref{Eq:xi-estim-1}
   and~\eqref{Eq:xi-estim-2} hold for $\tau \ll 0$ too. As $\delta <
   \pi p$ it easily follows that $\Vert \xi \Vert_{1,p;s_{\delta}} <
   \infty$.
\end{proof}

\subsection{Pregluing}  \label{S:pregluing} Let
$u_1,u_0:(\widehat{S},\partial \widehat{S}) \to (M,L)$ be two
$J$-holomorphic disks satisfying Assumption~\ref{A:T2}.  Put
$p=u_1(\infty)=u_0(-\infty)$. At this point it will be convenient to
reparametrize $u_1, u_0$ in the following way. Pick
$\tau_{\textnormal{shift}}^1, \tau_{\textnormal{shift}}^0 > 0$ very
large so that $u_1(\tau+\tau_{shift}^1,t),
u_0(-\tau-\tau_{shift}^0,t)\in B_p(r_M/2)$ for every $\tau>0$.  Here
$r_M$ is the radius used to define the map $P$ in
formula~\eqref{Eq:pal} of Section~\ref{S:analytic}. We replace
$u_1(\tau, t), u_0(\tau,t)$ by
$u_1(\tau+\tau_{\textnormal{shift}}^1,t),
u_0(\tau-\tau_{\textnormal{shift}}^0,t)$ respectively. In order for
the transversality assumption in~\ref{A:T2} to continue to hold,
replace $(\tau_1,\tfrac{1}{2})$ by
$(\tau_1-\tau_{\textnormal{shift}}^1,\tfrac{1}{2})$ and
$(\tau_0,\tfrac{1}{2})$ by
$(\tau_0+\tau_{\textnormal{shift}}^0,\tfrac{1}{2})$. After these
replacements we may assume without loss of generality that
$u_1(\tau,t), u_0(-\tau,t) \in B_p(r_M/2)$ for every $\tau>0$ and that
the transversality assumption in~\ref{A:T2} holds with
$(\tau_1,\tfrac{1}{2})$ and $(\tau_0,\tfrac{1}{2})$. Moreover, by
choosing $\tau_{\textnormal{shift}}^1,\tau_{\textnormal{shift}}^0$
large enough and so that
$\tau_{\textnormal{shift}}^0-\tau_{\textnormal{shift}}^1 =
-(\tau_0+\tau_1)$ we may assume that may assume that now we have
$\tau_1 = -\tau_0$ and that $\tau_0>0$. We shall write from now on
$\tau_*=\tau_0=-\tau_1$.

Let $\zeta_1(\tau,t),\zeta_0(\tau,t) \in T_p(M)$ be such that
$u_1(\tau,t)=\exp_p(\zeta_1(\tau,t))$ for every $\tau\geq 0$ and
$u_0(\tau,t)=\exp_p(\zeta_0(\tau,t))$ for every $\tau \leq 0$.  Let
$\sigma^{+}_R, \sigma^{-}_R: \mathbb{R} \to [0,1]$, $R\in \mathbb{R}$,
be a family of smooth functions with the properties:
\begin{enumerate}
  \item $\sigma^{-}_R(\tau)=0$ for every $\tau\geq R+1$,
   $\sigma^{-}_R(\tau)=1$ for every $\tau \leq R-1$.
  \item ${\sigma^{-}_{R}}'(\tau) \leq 0$ for every $\tau$.
  \item $\sigma^{+}_R = 1-\sigma^{-}_R$.
\end{enumerate}
For every $R>0$ define $u_R:(\widehat{S},\partial \widehat{S}) \to
(M,L)$ by:
\begin{equation} \label{Eq:uR}
   u_R(\tau,t)=
   \begin{cases}
      u_1(\tau+5R,t) & \tau\leq -4R, \\
      \exp_p\Bigl(\sigma^{+}_{-R}(\tau)
      \zeta_0(\tau-5R,t)+\sigma^{-}_R(\tau) \zeta_1(\tau+5R,t)\Bigr) &
      |\tau|\leq 4R, \\
      u_0(\tau-5R,t) & \tau\geq 4R.
   \end{cases}
\end{equation}
Note that $u_R(\tau,t)$ is well defined for $R \gg 0$, and that
$u_R(\tau,t)=u_0(\tau-5R,t)$ for every $\tau \geq R+1$ while
$u_R(\tau,t)=u_1(\tau+5R,t)$ for every $\tau \leq -R-1$.  Also note
that $u_R$ maps $\partial S$ to $L$. This is because $L$ is totally
geodesic with respect to the metric we use hence $\zeta_i(\tau,0),
\zeta_i(\tau,1) \in T_p(L)$.

\subsubsection{The Fredholm property revisited}
\label{Sb:fredholm-2} There exists
$\delta_{\textnormal{Fred}}(u_1,u_0)>0$ and $R_0>0$ such that for
every $R_0 \leq R$, $0<\delta \leq
\delta_{\textnormal{Fred}}(u_1,u_0)$ the operator $D_{u_R}$ is
Fredholm. Moreover its index is $\dim L + \mu(A_1+A_0)$. This follows
from the arguments in~\cite{Ro-Sa:spectral, Fl:unregularized}. The
reason for this independence of $\delta$ on $R$ is roughly speaking as
follows. First note that each of the $u_R$'s has exponential decay at
the ends of $S$ since they are holomorphic near the ends of $S$
(see~\ref{Sb:fredholm-1}). For a given $R$, the set of $\delta$'s for
which $D_{u_R}$ is not Fredholm coincides with the spectrum of an
(unbounded) self-adjoint operator $A_R$ which can be derived from the
asymptotic behavior of $D_{u_R}$ as $|\tau| \to \infty$
(see~\cite{Ro-Sa:spectral, Fl:unregularized} for more details. See
also~\cite{Sa:HF-lectures,Sc:thesis} for expositions of the Fredholm
property in the case of cylindrical Floer trajectories). The point is
that, due to our definition of $u_R$, the asymptotic operators $A_R$
do not depend on $R$, hence their spectrum remains constant. In this
case this spectrum is a discrete set (which contains $0$ !)  hence
there exists $\delta_{\textnormal{Fred}}(u_1,u_0)>0$ such that $(0,
\delta_{\textnormal{Fred}}(u_1,u_0)]$ lies in its complement.

\subsubsection{Deformed weighted norms}  \label{Sb:weighted}
Given $R>1$, $\delta>0$ define $\alpha_{R,\delta}:S \to \mathbb{R}$ by
\begin{equation} \label{Eq:weight}
   \alpha_{R,\delta}(\tau,t)=\alpha_{R,\delta}(\tau)=
   \begin{cases}
      e^{\delta(|\tau|-5R)} & |\tau| \geq 5R, \\
      e^{\delta(5R-|\tau|)} & |\tau| \leq 5R.
   \end{cases}
\end{equation}
Using the weight $\alpha_{R,\delta}$ we define a new norm on
$T_{u_R}^{1,p;\delta}$ by:
\begin{align}
   \Vert \xi \Vert^p_{1,p;\alpha_{R,\delta}} = & \int_{S(-\infty,-5R]}
   \alpha_{R,\delta}\bigl(|\xi-P_{u_R}(\xi_{-\infty})|^p +
   |\nabla(\xi-P_{u_R}(\xi_{-\infty}))|^p \bigr) \\
   + & \int_{S[5R, \infty)}
   \alpha_{R,\delta}\bigl(|\xi-P_{u_R}(\xi_{\infty})|^p +
   |\nabla(\xi-P_{u_R}(\xi_{\infty}))|^p \bigr) \notag \\
   + & \int_{S[-5R,5R]}
   \alpha_{R,\delta}\bigl(|\xi-\textnormal{Pal}_{u_R}
   \xi(0,\tfrac{1}{2})|^p +
   |\nabla(\xi-\textnormal{Pal}_{u_R}\xi(0,\tfrac{1}{2}))|^p \bigr)
   \notag \\
   + & |\xi_{-\infty}|^p + |\xi_{\infty}|^p + |\xi(0,\tfrac{1}{2})|^p.
   \notag
\end{align}
We endow $\mathcal{E}_{u_R}^{0,p;\delta}$ with the weighted norm:
$$\Vert \eta \Vert_{0,p;\alpha_{R,\delta}} = \Bigl(\int_S
\alpha_{R,\delta} | \eta |^p \Bigr)^{1/p}.$$
\begin{remnonum}
   Using the inequalities from Section~\ref{S:auxiliary} it is easy to
   see that $\Vert \cdot \Vert_{1,p;\alpha_{R,\delta}}$ is equivalent
   to $\Vert \cdot \Vert_{1,p;s_{\delta}}$ however they are not
   uniformly equivalent as $R \to \infty$, i.e. it is impossible to
   find one constant $C$ such that $\Vert \cdot
   \Vert_{1,p;\alpha_{R,\delta}} \leq C \Vert \cdot
   \Vert_{1,p;s_{\delta}}$ for every $R$.
\end{remnonum}

\subsubsection{Comparison of norms}  \label{Sb:comparison}

\begin{prop} \label{P:comparison}
   Let $u \in W^{1,p;\delta}(M,L)$. There exists a constant
   $C_{\ref{P:comparison}}>0$ depending on $(u_1,u_0)$ and on $p$, but
   not on $R,\delta$, such that for every $\xi \in
   T_{u_R}^{1,p;\delta}$, $R,\delta>0$, we have:
   $$\Vert \xi \Vert_{L^{\infty}} \leq C_{\ref{P:comparison}} \Vert
   \xi \Vert_{1,p;s_{\delta}}, \quad \Vert \xi \Vert_{L^{\infty}} \leq
   C_{\ref{P:comparison}} \Vert \xi \Vert_{1,p;\alpha_{R,\delta}}.$$
\end{prop}
\begin{proof}
   Let $x \in S[5R, \infty)$. Since $\alpha_{R,\delta}\geq 1$, we get
   from Proposition~\ref{P:sblv}:
   \begin{align*}
      |\xi(x)| & \leq |\xi(x)-P_{u_R}(\xi_{\infty})(x)|+
      |\xi_{\infty}| \\
      & \leq
      C_{\ref{P:sblv}}\Bigl(\int_{S[5R,\infty)}\alpha_{R,\delta}
      \bigl(|\xi-P_{u_R}\xi_{\infty}|^p +
      |\nabla(\xi-P_{u_R}\xi_{\infty})|^p \Bigr)^{1/p}+|\xi_{\infty}| \\
      & \leq C \Vert \xi \Vert_{1,p;\alpha_{R,\delta}},
   \end{align*}
   where $C=2(C_{\ref{P:sblv}}+1)$.

   Let $x \in S[-5R,5R]$. Then by Proposition~\ref{P:sblv}:
   \begin{align*}
      |\xi(x)| & \leq
      |\xi(x)-\textnormal{Pal}_{u_R(x)}\xi(0,\tfrac{1}{2})|
      +|\xi(0,\tfrac{1}{2})| \\
      & \leq C_{\ref{P:sblv}}\Bigl(\int_{S[-5R,5R]}
      \bigl(|\xi-\textnormal{Pal}_{u_R} \xi(0,\tfrac{1}{2})|^p +
      |\nabla(\xi-\textnormal{Pal}_{u_R}\xi(0,\tfrac{1}{2}))|^p \bigr)
      \Bigr)^{1/p} + |\xi(0,\tfrac{1}{2})|, \\
      & \leq C \Vert \xi \Vert_{1,p;\alpha_{R,\delta}}.
      \end{align*}
      In a similar way one shows that for every $x \in
      S(-\infty,-5R]$, $|\xi(x)|\leq C \Vert\xi
      \Vert_{1,p;\alpha_{R,\delta}}$. This proves that $\Vert \xi
      \Vert_{L^{\infty}} \leq C \Vert \xi
      \Vert_{1,p;\alpha_{R,\delta}}$. The proof of the inequality for
      $\Vert \xi \Vert_{1,p;s_{\delta}}$ is similar.
\end{proof}

\subsubsection{Estimates on $u_R$}  \label{Sb:estimates-uR}
\begin{prop} \label{P:uR}
   There exists constants $\delta_0>0$, $R_0>1$, $C'_{\ref{P:uR}},
   c'_{\ref{P:uR}}, C''_{\ref{P:uR}}>0$ that depend only on
   $(u_1,u_0)$ such that for every $0<\delta\leq \delta_0$ and $R_0
   \leq R$ we have:
   $$\Vert \overline{\partial}_J u_R \Vert_{0,p;\alpha_{R,\delta}}
   \leq C'_{\ref{P:uR}}e^{-c'_{\ref{P:uR}} R}, \quad \Vert d u_R
   \Vert_{0,p;\alpha_{R,\delta}} \leq C''_{\ref{P:uR}}, \quad \Vert d
   u_R \Vert_{L^{\infty}} \leq C''_{\ref{P:uR}}.$$
\end{prop}

\begin{proof}[Outline of the proof]
   The proof is a straightforward computation combined with the
   following points:
   \begin{enumerate}
     \item We have an exponential decay estimate for $d u_i$, $i=1,0$,
      namely there exists constants $C,c>0$ such that
      $|d{u_i}_{(\tau,t)}| \leq C e^{-c|\tau|}$. Similar exponential
      decay estimates hold also for $\zeta_1,\zeta_0$.  It follows
      that $|d {u_0}_{(\tau-5R,t)}| \leq Ce^{-c|\tau-5R|}$ and $|d
      {u_1}_{(\tau+5R,t)}| \leq Ce^{-c|\tau+5R|}$.
     \item $\overline{\partial}_J u_R=0$ outside $S[-R-1,R+1]$, thus
      for estimating $\Vert \overline{\partial}_J u_R
      \Vert_{0,p;\alpha_{R,\delta}}$ it is enough to estimate $\Vert d
      u_R \Vert_{0,p;\alpha_{R,\delta}}$ on $S[-R-1,R+1]$.
     \item There exists $K>0$ such that $|d(\exp_p)_v|\leq K$ for
      every $v \in T_p(M)$ with $|v|\leq \epsilon$. Note that for $R$
      large enough the vectors appearing inside the $\exp_p$ in
      expression~\eqref{Eq:uR} have norm $\leq \epsilon$ due to
      exponential decay of the $\zeta_i$'s.
   \end{enumerate}
\end{proof}

\subsection{The main operators}  \label{S:main-op}

In view of the simplification in
Section~\ref{Sb:gluing-simplification}, fix four smooth maps $G_{\pm}:
L \to \mathbb{R}^{n-k_{\pm}}$, $G_i:M \to \mathbb{R}^{2n-k_i}$,
$i=1,0$ such that $G_{\pm}$ are submersions near $W_{\pm}, W_i$ and
such that near $W_{\pm},W_i$ we have $G_{\pm}^{-1}(0)=W_{\pm}$,
$G_i^{-1}(0)=W_i$. Put
\begin{equation} \label{Eq:T10}
   T_{u_1,u_0}^{1,p;\delta} = \{(\xi_1,\xi_0) \in
   T_{u_1}^{1,p;\delta} \oplus T_{u_0}^{1,p;\delta} \mid
   \xi_1(\infty)=\xi_0(-\infty) \}.
\end{equation}
Define the following operator:
\begin{align*}
   & D'_{u_1,u_0}:T_{u_1,u_0}^{1,p;\delta} \longrightarrow
   \mathcal{E}_{u_1}^{0,p;\delta} \oplus
   \mathcal{E}_{u_0}^{0,p;\delta} \oplus \mathbb{R}^{n-k_{-}} \oplus
   \mathbb{R}^{2n-k_1} \oplus
   \mathbb{R}^{2n-k_0} \oplus \mathbb{R}^{n-k_{+}}, \\
   & D'_{u_1,u_0}(\xi_1,\xi_0)=\bigl(D_{u_1}\xi_1, D_{u_0}\xi_0,
   dG_{-}\xi_1(-\infty), dG_1 \xi_1(-\tau_*,\tfrac{1}{2}), dG_0
   \xi_0(\tau_*,\tfrac{1}{2}), dG_{+}\xi_0(\infty) \bigr).
\end{align*}
For brevity denote by $E=\mathbb{R}^{n-k_{-}} \oplus
\mathbb{R}^{2n-k_1} \oplus \mathbb{R}^{2n-k_0} \oplus
\mathbb{R}^{n-k_{+}}$ the summand of the last four linear spaces in
the target of $D'_{u_1,u_0}$. We also write $D_{u_1,u_0}(\xi_1,\xi_0)$
for the first two coordinates of $D'_{u_1,u_0}(\xi_1,\xi_0)$ and $d
\mathcal{G}_{1,0}(\xi_1, \xi_0)$ for its last four coordinates.

\begin{prop} \label{P:surjective}
   Under Assumption~\ref{A:T2}, the operator $D'_{u_1,u_0}$ is
   surjective.
\end{prop}
\begin{proof}
   Let $\eta_i \in \mathcal{E}_{u_i}^{0,p;\delta}$, $i=1,0$ and
   $\overrightarrow{a}=(a_{-},a_1,a_0,a_{+}) \in E$.  Choose $w_{-}
   \in T_{u_1(-\infty)}(W_{-})$, $w_1 \in
   T_{u_1(-\tau_*,\tfrac{1}{2})}(W_1)$, $w_0 \in
   T_{u_0(\tau_*,\tfrac{1}{2})}(W_0)$, $w_{+} \in T_{u_0(\infty)}$
   such that $dG_{\pm}(w_{\pm})=a_{\pm}$, $dG_i(w_i)=a_i$, $i=1,0$.

   Since $D_{u_1}, D_{u_0}$ are surjective there exist $\xi'_1,
   \xi'_0$ such that $D_{u_i}\xi'_i=\eta_i$. By Assumption~\ref{A:T2}
   there exist $\zeta_i \in T_{u_i}(\mathcal{M}(A_i,J))$ such that
   \begin{align*}
      \Bigl(\zeta_1(-\infty), \zeta_1(-\tau_*,\tfrac{1}{2}),
      \zeta_1(\infty), \zeta_0(-\infty), \zeta_0(\tau_*,\tfrac{1}{2}),
      \zeta_0(\infty)\Bigr) \in &
      \bigl(w_-,w_1,\xi'_1(\infty),\xi'_0(-\infty),
      w_0,w_+\bigr) \\
      & + T_{\underline{x}}\bigl(W_- \times W_1 \times
      \textnormal{diag}(L) \times W_0 \times W_+\bigr),
   \end{align*}
   where $\underline{x}=(u_1(-\infty), u_1(-\tau_*,\tfrac{1}{2}),
   u_1(\infty), u_0(-\infty),u_0(\tau_*,\tfrac{1}{2}),u_0(\infty))$.

   Put $\xi_i=\xi'_i-\zeta_i$, $i=1,0$. Clearly $(\xi_1,\xi_0) \in
   T_{u_1,u_0}^{1,p;\delta}$ and $D'_{u_1,u_0}(\xi_1,\xi_0) = (\eta_1,
   \eta_0, \overrightarrow{a})$.
\end{proof}

Define a map
$$\mathcal{F}'_{u_R}: T_{u_R}^{1,p;\delta} \longrightarrow
\mathcal{E}_{u_R}^{0,p;\delta} \times \mathbb{R}^{n-k_-} \times
\mathbb{R}^{2n-k_1} \times \mathbb{R}^{2n-k_0} \times
\mathbb{R}^{n-k_+},$$
\begin{align*}
   \mathcal{F}'_{u_R}(\xi) = (\mathcal{F}_{u_R}(\xi),
   & G_-(\exp_{u_R} \xi(-\infty)), G_1(\exp_{u_R}\xi(-5R-\tau_*)), \\
   & G_0(\exp_{u_R}\xi(5R+\tau_*)), G_+(\exp_{u_R} \xi(\infty))).
\end{align*}
Let $D'_{u_R}=d \mathcal{F}'_{u_R}(0):T_{u_R}^{1,p;\delta} \to
\mathcal{E}_{u_R}^{0,p;\delta} \oplus E$ be its linearization at
$\xi=0$. We have $D'_{u_R} \xi = (D_{u_R}\xi, d \mathcal{G}_R \xi)$,
where $$d \mathcal{G}_R \xi = (d G_- \xi(-\infty), d G_1
\xi(-5R-\tau_*,\tfrac{1}{2}), d G_0 \xi(5R+\tau_*,\tfrac{1}{2}), d G_+
\xi(\infty)).$$
From now on we shall often endow the space
$\mathcal{E}_{u_R}^{0,p;\delta} \oplus E$ with the norm
\begin{equation} \label{Eq:norm'}
   \Vert(\eta, \overrightarrow{a}) \Vert'_{0,p;\alpha_{R,\delta}} =
   \Vert \eta \Vert_{0,p;\alpha_{R,\delta}} + |\overrightarrow{a}|.
\end{equation}
Clearly the operators $D'_{u_R}$ have the same Fredholm properties as
$D_{u_R}$ as described in Sections~\ref{Sb:fredholm-1}
and~\ref{Sb:fredholm-2}.

\subsection{A right inverse to $D_{u_R}$}  \label{S:rinv}

\begin{prop} \label{P:rinv}
   There exist $\delta_0>0$, $R_0>0$ such that for every $0<\delta\leq
   \delta_0$, there exists a family of operators $\bigl\{
   Q_R:\mathcal{E}_{u_R}^{0,p;\delta} \oplus E \to
   T_{u_R}^{1,p;\delta} \bigr\}_{R_0\leq R}$, with the following
   properties:
   \begin{enumerate}
     \item $D'_{u_R} \circ Q_R = \Id$.
     \item $Q_R$ is uniformly bounded in the
      $\alpha_{R,\delta}$-weighted norms, i.e.  there exists a
      constant $C_{\ref{P:rinv}}(\delta)>0$ that does not depend on
      $R$ such that for every $R_0 \leq R$, $(\eta,\overrightarrow{a})
      \in \mathcal{E}_{u_R}^{0,p;\delta}\oplus E$ we have $\Vert
      Q_R(\eta,\overrightarrow{a}) \Vert_{1,p;\alpha_{R,\delta}} \leq
      C_{\ref{P:rinv}} (\delta)(\Vert \eta
      \Vert_{0,p;\alpha_{R,\delta}}+|\overrightarrow{a}|)$.
     \item $Q_R$ depends smoothly on $R$. (See remark below.)
   \end{enumerate}
\end{prop}

\begin{remsnonum}
   \begin{enumerate}
     \item Uniform boundedness of the $Q_R$'s does not seem to hold
      for the $s_{\delta}$-weighted norms.
     \item In order to make sense of the smooth dependence of $Q_R$ on
      $R$ one has to identify the spaces $T^{1,p;\delta}_{u_R}$ for
      nearby $R's$. This can be done using parallel transport along
      short geodesics.
   \end{enumerate}
\end{remsnonum}

The rest of this section is devoted to the proof of
Proposition~\ref{P:rinv}. The reader who wishes to get a less
technical account of the gluing may skip to Section~\ref{S:implicit}
in which the proof of Theorem~\ref{T:gluing} is carried out.

\subsubsection{Some auxiliary operators}  \label{Sb:aux-op}
Put $$T_{u_1,u_0}^{1,p;\delta} = \{ (\xi_1,\xi_0)\in
T_{u_1}^{1,p;\delta} \oplus T_{u_0}^{1,p;\delta} \mid
\xi_1(\infty)=\xi_0(-\infty)\}.$$
Define a linear map
$I^R:T_{u_1,u_0}^{1,p;\delta} \to T_{u_R}^{1,p;\delta}$
\begin{align*}
   & I^R(\xi_1,\xi_0)(\tau,t) = \\
   & =
   \begin{cases}
      \xi_1(\tau+5R,t), & \tau\leq -4R, \\
      \textnormal{Pal}_{u_R(\tau,t)}\Bigl(v+\sigma^{+}_{-R}
      \bigl(\textnormal{Pal}_p(\xi_0(\tau-5R,t))-v\bigr) +
      \sigma^{-}_{R}\bigl(\textnormal{Pal}_p(\xi_1(\tau+5R,t))-v
      \bigr)
      \Bigr), & |\tau|\leq 4R, \\
      \xi_0(\tau-5R,t) & \tau\geq 4R,
   \end{cases}
\end{align*}
where:
\begin{enumerate}
  \item $p=u_1(\infty)=u_0(-\infty)$.
  \item $v=\xi_1(\infty)=\xi_0(-\infty)$.
  \item $\textnormal{Pal}_x(\eta)$ is defined for every $x \in
   B_p(r_M)$ and every $\eta \in T_y(M)$, $y \in B_p(r_M)$. It is the
   parallel transport of $\eta$ along the minimal geodesic connecting
   $y$ to $x$.
\end{enumerate}

Define also the following maps:
\begin{align*}
   & J^R_{0,*}: \mathcal{E}^{0,p;\delta}_{u_0} \to
   \mathcal{E}^{0,p;\delta}_{u_R}, & &J^R_{0,*} \eta(\tau,t)=
   \textnormal{Pal}_{u_R(\tau,t)}\bigl(\sigma^{+}_{-R}
   (\tau)\eta(\tau-5R,t)\bigr) \\
   &J^R_{1,*}: \mathcal{E}^{0,p;\delta}_{u_1} \to
   \mathcal{E}^{0,p;\delta}_{u_R}, & &J^R_{1,*} \eta(\tau,t)=
   \textnormal{Pal}_{u_R(\tau,t)}
   \bigl(\sigma^{-}_{R}(\tau)\eta(\tau+5R,t)\bigr)
\end{align*}
Note that $J^R_{0,*}\eta(\tau,t)=\eta(\tau-5R,t)$ for every $\tau\geq
R+1$ while $J^R_{0,*}\eta(\tau,t)=0$ for every $\tau\leq -R-1$.
Similarly define
\begin{align*}
   & J^R_{*,0}: \mathcal{E}^{0,p;\delta}_{u_R} \to
   \mathcal{E}^{0,p;\delta}_{u_0}, & &J^R_{*,0} \eta(\tau,t)=
   \textnormal{Pal}_{u_0(\tau,t)}\bigl(\sigma^{+}_{-R}(\tau+5R)
   \eta(\tau+5R,t)\bigr) \\
   &J^R_{*,1}: \mathcal{E}^{0,p;\delta}_{u_R} \to
   \mathcal{E}^{0,p;\delta}_{u_1}, & &J^R_{*,1} \eta(\tau,t)=
   \textnormal{Pal}_{u_1(\tau,t)}\bigl(\sigma^{-}_{R}(\tau-5R)
   \eta(\tau-5R,t)\bigr)
\end{align*}
Note that $J^R_{*,0}\eta(\tau,t)=\eta(\tau+5R,t)$ for every $\tau\geq
-4R+1$. The following identities are easily verified:
\begin{equation} \label{Eq:id-3}
   \begin{split}
      & J^R_{0,*} \circ J^R_{*,0}(\sigma^{+}_0(\tau)\eta)=
      \sigma^{+}_0 \eta, \\
      & J^R_{1,*} \circ J^R_{*,1}(\sigma^{-}_0(\tau)\eta)=
      \sigma^{-}_0 \eta, \\
      & J^R_{0,*} \circ J^R_{*,0}(\sigma^{+}_0(\tau)\eta) + J^R_{1,*}
      \circ J^R_{*,1}(\sigma^{-}_0(\tau)\eta) = \eta.
   \end{split}
\end{equation}

We now have the following estimate:
\begin{prop} \label{P:estim-DuI-JD}
   There exist constants $C_{\ref{P:estim-DuI-JD}}(R)>0$, with $\lim_{R\to
     \infty} C_{\ref{P:estim-DuI-JD}}(R)=0$ and such that:
   $$\Vert D_{u_R} \circ I^{R}(\xi_1,\xi_0)- J^R_{0,*}D_{u_0}\xi_0 -
   J^R_{1,*}D_{u_1}\xi_1 \Vert_{0,p;\alpha_{R,\delta}} \leq
   C_{\ref{P:estim-DuI-JD}}(R) (\Vert\xi \Vert_{1,p;s_{\delta}}+\Vert\xi_0
   \Vert_{1,p;s_{\delta}}).$$
\end{prop}
The proof is a straightforward computation, based on the definition of
the weighted norms $s_{\delta},\alpha_{R,\delta}$. See~\cite{FO3} for
more details.

\subsubsection{The approximate right inverse}
\label{Sb:approximate-rinv} Let
$Q_0':\mathcal{E}^{0,p;\delta}_{u_1} \oplus
\mathcal{E}^{0,p;\delta}_{u_0} \oplus E \to T^{1,p;\delta}_{u_1,u_0}$
be a bounded right inverse to $D'_{u_1,u_0}$ (see
Proposition~\ref{P:surjective}). Write $Q_0'$ as
\begin{equation} \label{Eq:Q0}
   Q_0'(\eta_1,\eta_0,\overrightarrow{a}) = Q_0(\eta_1,
   \eta_0)+A(\overrightarrow{a}).
\end{equation}

Define $\widetilde{Q}_R:\mathcal{E}^{0,p;\delta}_{u_R} \oplus E \to
T^{1,p;\delta}_{u_R}$ by
\begin{equation} \label{Eq:Qtilde}
   \widetilde{Q}_R(\eta,\overrightarrow{a})=I^R \Bigl(Q_0
   \bigl(J^R_{*,0} (\sigma^{+}_0 \eta), J^R_{*,1}(\sigma^{-}_0
   \eta) \bigr)+A(\overrightarrow{a})\Bigr).
\end{equation}

\begin{prop} \label{P:estim-Qtilde}
   There exist constants $R_0,\delta_0>0$ and constants
   $C_{\ref{P:estim-Qtilde}}(\delta)$ such that for every $R_0\leq R$,
   $0<\delta\leq \delta_0$, $\eta \in \mathcal{E}^{0,p;\delta}_{u_R}$,
   $\overrightarrow{a} \in E$ we have
   $$\Vert \widetilde{Q}_R (\eta, \overrightarrow{a})
   \Vert_{1,p;\alpha_{R,\delta}} \leq C_{\ref{P:estim-Qtilde}}(\delta)
   (\Vert \eta \Vert_{0,p;\alpha_{R,\delta}}+ |\overrightarrow{a}|).$$
\end{prop}

\begin{proof}[Proof of Proposition~\ref{P:estim-Qtilde}]
   A simple computation shows that:
   \begin{lem} \label{L:J-estim} There exists $C_{\ref{L:J-estim}}>0$
      such that for every $R\gg 0$:
      \begin{align*}
         \Vert J^R_{*,0}(\sigma^{+}_0\eta) \Vert_{0,p;s_{\delta}} \leq
         C_{\ref{L:J-estim}} \Vert \eta \Vert_{0,p;\alpha_{R,\delta}}, \\
         \Vert J^R_{*,1}(\sigma^{-}_0\eta) \Vert_{0,p;s_{\delta}} \leq
         C_{\ref{L:J-estim}} \Vert \eta \Vert_{0,p;\alpha_{R,\delta}}
      \end{align*}
   \end{lem}
   To prove Proposition~\ref{P:estim-Qtilde} it suffices to prove the
   following:
   \begin{prop} \label{P:estim-IR}
      There exist constants $R_0$, $\delta_0>0$ and
      $C_{\ref{P:estim-IR}}(\delta)$ such that for every $R_0 \leq R$,
      $0<\delta\leq\delta_0$, $(\xi_1,\xi_0)\in
      T^{1,p;\delta}_{u_1,u_0}$, we have
      $$\Vert I^R(\xi_1,\xi_0) \Vert_{1,p;\alpha_{R,\delta}} \leq
      C_{\ref{P:estim-IR}}(\delta)(\Vert \xi_1 \Vert_{1,p;s_{\delta}}
      + \Vert \xi_0 \Vert_{1,p;s_{\delta}}).$$
   \end{prop}
\end{proof}

\begin{proof}[Proof of Proposition~\ref{P:estim-IR}]
   Put $\xi_{\textnormal{new}}=I^R(\xi_1,\xi_0)$,
   $v=\xi_0(-\infty)=\xi_1(\infty)$.  It is easy to see that
   \begin{align*}
      & \int_{S[5R,\infty)} \alpha_{R,\delta}
      \Bigl(|\xi_{\textnormal{new}}-P_{u_R}\xi_0(\infty)|^p +
      |\nabla(\xi_{\textnormal{new}}-P_{u_R}\xi_0(\infty))|^p
      \Bigr) \\
      + & \int_{S(-\infty,5R]} \alpha_{R,\delta}
      \Bigl(|\xi_{\textnormal{new}}-P_{u_R}\xi_1(-\infty)|^p +
      |\nabla(\xi_{\textnormal{new}}-P_{u_R}\xi_1(-\infty))|^p
      \Bigr) \\
      \leq \; & \Vert \xi_0 \Vert_{1,p;s_{\delta}}^p+ \Vert\xi_1
      \Vert_{1,p;s_{\delta}}^p + |\xi_0(\infty)|^p +
      |\xi_1(-\infty)|^p.
   \end{align*}
   It remains to prove that
   \begin{align} \label{Eq:I}
      I:= & \int_{S[-5R,5R]}\alpha_{R,\delta}\Bigl(
      |\xi_{\textnormal{new}}-\textnormal{Pal}_{u_R}
      \xi_{\textnormal{new}}(0,\tfrac{1}{2})|^p +
      |\nabla(\xi_{\textnormal{new}}-\textnormal{Pal}_{u_R}
      \xi_{\textnormal{new}}(0,\tfrac{1}{2}))|^p \Bigr) \\
      & + |\xi_{\textnormal{new}}(0,\tfrac{1}{2})|^p \leq C_1(\Vert
      \xi_0 \Vert_{1,p;s_{\delta}}^p + \Vert \xi_1
      \Vert_{1,p;s_{\delta}}^p), \notag
   \end{align}
   for some constant $C_1$ that does not depend on $R$. Note that due
   to the reparametrization we have made in the beginning of
   Section~\ref{S:pregluing} we can write $\textnormal{Pal}_{u_R}$ in
   the two terms under the integral of~\eqref{Eq:I} instead of
   $P_{u_R}$. Now
   \begin{equation} \label{Eq:I-estim}
      I \leq C_2(I_1+I_2+|\xi_{\textnormal{new}}(0,\tfrac{1}{2})|^p),
   \end{equation}
   for some constant $C_2>0$, where
   \begin{align*}
      & I_1=\int_{S[-5R,5R]}\alpha_{R,\delta}
      \Bigl(|\xi_{\textnormal{new}}-\textnormal{Pal}_{u_R}v|^p+
      |\nabla(\xi_{\textnormal{new}}-\textnormal{Pal}_{u_R}v)|^p
      \Bigr), \\
      & I_2 = \int_{S[-5R,5R]}\alpha_{R,\delta}
      \Bigl(|\textnormal{Pal}_{u_R}
      \xi_{\textnormal{new}}(0,\tfrac{1}{2})-
      \textnormal{Pal}_{u_R}v|^p+|\nabla(\textnormal{Pal}_{u_R}
      \xi_{\textnormal{new}}(0,\tfrac{1}{2})-
      \textnormal{Pal}_{u_R}v)|^p \Bigr).
   \end{align*}

   We begin by estimating $I_2$.  Put $V(x)=\textnormal{Pal}_{u_R(x)}
   \xi_{\textnormal{new}}(0,\tfrac{1}{2})-\textnormal{Pal}_{u_R(x)}v$.
   Then, by the definition of $\xi_{\textnormal{new}}$ we have:
   \begin{align*}
      V(x)&=
      \textnormal{Pal}_{u_R(x)}\textnormal{Pal}_{u_R(0,\tfrac{1}{2})}
      \textnormal{Pal}_p \xi_0(-5R, \tfrac{1}{2}) +
      \textnormal{Pal}_{u_R(x)}\textnormal{Pal}_{u_R(0,\tfrac{1}{2})}
      \textnormal{Pal}_p \xi_1(5R,\tfrac{1}{2}) \\
      & -(\textnormal{Pal}_{u_R(x)}v +
      \textnormal{Pal}_{u_R(x)}\textnormal{Pal}_{u_R(0,\tfrac{1}{2})}v).
   \end{align*}

   Using the fact that for every $(\tau,t) \in S[-5R,5R]$ we have
   $\textnormal{dist}(u_R(\tau,t),p) \leq C_3e^{-c_3(5R-|\tau|)}$ for
   some constants $C_3,c_3>0$ (that do not depend on $R$) we obtain by
   standard arguments from ode's that:
   \begin{equation} \label{Eq:V}
         |V(\tau,t)| \leq |\xi_0(-5R,\tfrac{1}{2}) -
         \textnormal{Pal}_{u_0(-5R,\tfrac{1}{2})}v| +
         |\xi_1(5R,\tfrac{1}{2}) -
         \textnormal{Pal}_{u_1(5R,\tfrac{1}{2})}v|
          + C_4e^{-c_4(5R-|\tau|)}|v|,
   \end{equation}
   for some constants $C_4,c_4>0$.  Similar arguments show that:
   \begin{align} \label{Eq:NV}
      |\nabla V(\tau,t)| \leq & C_5 \bigl(|\xi_0(-5R,\tfrac{1}{2}) -
      \textnormal{Pal}_{u_0(-5R,\tfrac{1}{2})}v| +
      |\xi_1(5R,\tfrac{1}{2}) -
      \textnormal{Pal}_{u_1(5R,\tfrac{1}{2})}v| + |v| \bigr) \\
      + & C_5e^{-c_5(5R-|\tau|)}(|\xi_0(-5R,\tfrac{1}{2})| +
      |\xi_1(5R,\tfrac{1}{2})|+|v|) \notag
   \end{align}
   for some constants $C_5,c_5>0$. Note that in both~\eqref{Eq:V}
   and~\eqref{Eq:NV} the constants $C_4,c_4,C_5,c_5$ do not depend on
   $R$.

   We claim that there exists a constant $\delta_0>0$, independent of
   $R$, such that for every $0<\delta\leq \delta_0$ there exists
   $C_6=C_6(\delta)$ (independent of $R$) such that:
   \begin{equation} \label{Eq:I2}
      I_2 = \int_{S[-5R,5R]} \alpha_{R,\delta} (|V|^p+|\nabla V|^p) \leq
      C_6(\Vert \xi_0 \Vert_{1,p;s_{\delta}}^p+\Vert\xi_1
      \Vert_{1,p;s_{\delta}}^p+|v|^p).
   \end{equation}
   Indeed, by~\eqref{Eq:V},~\eqref{Eq:NV} there exist constants
   $C_7,c_7$ (independent of $R$ and of $\delta$) such that:
   \begin{align} \label{Eq:V11}
      & \int_{S[-5R,5R]} \alpha_{R,\delta} (|V|^p+|\nabla V|^p) \\
      \leq & C_7 \Bigl(|\xi_0(-5R,\tfrac{1}{2}) -
      \textnormal{Pal}_{u_0(-5R,\tfrac{1}{2})}v|^p +
      |\xi_1(5R,\tfrac{1}{2}) -
      \textnormal{Pal}_{u_1(5R,\tfrac{1}{2})}v|^p\Bigr)
      \int_{-5R}^{5R} \alpha_{R,\delta}(\tau) d\tau  \notag \\
      + & C_7 \bigl(|\xi_0(-5R,\tfrac{1}{2})|^p +
      |\xi_1(5R,\tfrac{1}{2})|^p+|v|^p\bigr) \int_{-5R}^{5R}
      e^{-c_7(5R-|\tau|)}\alpha_{R,\delta}(\tau)d\tau \notag \\
      \leq & \frac{2C_7}{\delta} e^{\delta 5R}
      \bigl(|\xi_0(-5R,\tfrac{1}{2}) -
      \textnormal{Pal}_{u_0(-5R,\tfrac{1}{2})}v|^p +
      |\xi_1(5R,\tfrac{1}{2}) -
      \textnormal{Pal}_{u_1(5R,\tfrac{1}{2})}v|^p\bigr) \notag \\
      + & \bigl(|\xi_0(-5R,\tfrac{1}{2})|^p +
      |\xi_1(5R,\tfrac{1}{2})|^p + |v|^p \bigr) C_7\int_{-5R}^{5R}
      e^{-c7(5R-|\tau|)+\delta(5R-|\tau|)} d \tau. \notag
   \end{align}
   Now $\int_{-5R}^{5R} e^{-c7(5R-|\tau|)+\delta(5R-|\tau|)} d \tau =
   \tfrac{2}{c_7-\delta}(1-e^{(\delta-c_7)5R})$, hence if we take
   $\delta_0 = c_7/2$ we obtain from~\eqref{Eq:V11} that for every
   $0<\delta<\delta_0$:
   \begin{align} \label{Eq:V12}
      & \int_{S[-5R,5R]} \alpha_{R,\delta} (|V|^p+|\nabla V|^p) \\
      \leq & \frac{2C_7}{\delta} e^{\delta 5R}
      \bigl(|\xi_0(-5R,\tfrac{1}{2}) -
      \textnormal{Pal}_{u_0(-5R,\tfrac{1}{2})}v|^p +
      |\xi_1(5R,\tfrac{1}{2}) -
      \textnormal{Pal}_{u_1(5R,\tfrac{1}{2})}v|^p\bigr) \notag \\
      + & \bigl(|\xi_0(-5R,\tfrac{1}{2})|^p +
      |\xi_1(5R,\tfrac{1}{2})|^p + |v|^p \bigr)\frac{4C_7}{c_7}.
      \notag
   \end{align}
   Applying Proposition~\ref{P:sblv} to~\eqref{Eq:V12} we obtain for
   every $0<\delta<\delta_0$:
   \begin{align} \label{Eq:V1}
      & \int_{S[-5R,5R]} \alpha_{R,\delta} (|V|^p+|\nabla V|^p) \\
      \leq & \frac{2C_7}{\delta} e^{\delta 5R} C_{\ref{P:sblv}}
      \int_{S[-5R-1,-5R]}\bigl(|\xi_0(\tau,t) -
      \textnormal{Pal}_{u_0(\tau,t)}v|^p + |\nabla(\xi_0 -
      \textnormal{Pal}_{u_0} v)|^p \bigr) e^{\delta(|\tau|-5R)}
      d\tau dt \notag \\
      + & \frac{2C_7}{\delta} e^{\delta 5R} C_{\ref{P:sblv}}
      \int_{S[5R,5R+1]}\bigl(|\xi_1(\tau,t) -
      \textnormal{Pal}_{u_1(\tau,t)}v|^p + |\nabla(\xi_1 -
      \textnormal{Pal}_{u_1} v)|^p \bigr) e^{\delta(|\tau|-5R)}
      d\tau dt \notag \\
      + & \bigl(|\xi_0(-5R,\tfrac{1}{2})|^p +
      |\xi_1(5R,\tfrac{1}{2})|^p + |v|^p \bigr) \frac{4C_7}{c_7}
      \notag \\
      =\, & \frac{2C_7}{\delta} C_{\ref{P:sblv}}
      \int_{S[-5R-1,-5R]}\bigl(|\xi_0(\tau,t) -
      \textnormal{Pal}_{u_0(\tau,t)}v|^p + |\nabla(\xi_0 -
      \textnormal{Pal}_{u_0} v)|^p \bigr) s_{\delta}(\tau)
      d\tau dt \notag \\
      + & \frac{2C_7}{\delta} C_{\ref{P:sblv}}
      \int_{S[5R,5R+1]}\bigl(|\xi_1(\tau,t) -
      \textnormal{Pal}_{u_1(\tau,t)}v|^p + |\nabla(\xi_1 -
      \textnormal{Pal}_{u_1} v)|^p \bigr) s_{\delta}(\tau)
      d\tau dt \notag \\
      + & \bigl(|\xi_0(-5R,\tfrac{1}{2})|^p +
      |\xi_1(5R,\tfrac{1}{2})|^p + |v|^p \bigr) \frac{4C_7}{c_7}.
      \notag
   \end{align}
   By Corollary~\ref{P:sblv}
   \begin{align*}
      & |\xi_0(-5R,\tfrac{1}{2})| \leq
      |\xi_0(-5R,\tfrac{1}{2})-P_{u_0(-5R,\tfrac{1}{2})}v|+|v| \leq
      C_{\ref{P:sblv}} \Vert \xi_0 \Vert_{1,p;s_{\delta}} + |v|, \\
      & |\xi_1(5R,\tfrac{1}{2})| \leq
      |\xi_1(5R,\tfrac{1}{2})-P_{u_1(5R,\tfrac{1}{2})}v|+|v| \leq
      C_{\ref{P:sblv}} \Vert \xi_1 \Vert_{1,p;s_{\delta}} + |v|
   \end{align*}
   Hence we obtain from~\eqref{Eq:V1}
   $$\int_{S[-5R,5R]} \alpha_{R,\delta} (|V|^p+|\nabla V|^p) \\
   \leq C_6 (\Vert \xi_0 \Vert_{1,p;s_{\delta}}^p+\Vert\xi_1
   \Vert_{1,p;s_{\delta}}^p+|v|^p),$$
   for some $C_6=C_6(\delta)$. This
   proves inequality~\eqref{Eq:I2} which estimates $I_2$.

   We turn to estimating $|\xi_{\textnormal{new}}(0,\tfrac{1}{2})|$:
   \begin{align} \label{Eq:xi05}
      |\xi_{\textnormal{new}}(0,\tfrac{1}{2})| & =
      |\textnormal{Pal}_{u_R(0,\tfrac{1}{2})} \bigl(\textnormal{Pal}_p
      \xi_0(-5R,\tfrac{1}{2})+\textnormal{Pal}_p
      \xi_1(5R,\tfrac{1}{2})\bigr) -
      \textnormal{Pal}_{u_R(0,\tfrac{1}{2})} v| \\
      & \leq |\xi_0(-5R,\tfrac{1}{2})|+|\xi_1(5R,\tfrac{1}{2})|+|v|
      \notag \\
      & \leq |\xi_0(-5R,\tfrac{1}{2}) -
      \textnormal{Pal}_{u_0(-5R,\tfrac{1}{2})}v|+|\xi_1(5R,\tfrac{1}{2})
      - \textnormal{Pal}_{u_1(5R,\tfrac{1}{2})}v|+3v \notag \\
      & \leq (C_{\ref{P:sblv}}+3)(\Vert \xi_0
      \Vert_{1,p;s_{\delta}}+\Vert\xi_1 \Vert_{1,p;s_{\delta}}+|v|).
      \notag
   \end{align}

   Finally we estimate $I_1$.
   $$\xi_{\textnormal{new}}(\tau,t) - \textnormal{Pal}_{u_R(\tau,t)}v
   =
   \textnormal{Pal}_{u_r(\tau,t)}\bigl(\sigma^{+}_{-R}(\textnormal{Pal}_p
   \xi_0(\tau-5R,t)-v)+\sigma^{-}_R(\textnormal{Pal}_p
   \xi_1(\tau+5R,t)-v)\bigr)$$

   Now $\textnormal{Pal}_{u_R(\tau,t)}(\textnormal{Pal}_p
   \xi_0(\tau-5R,t)-v) =
   \textnormal{Pal}_{u_R(\tau,t)}\textnormal{Pal}_p
   (\xi_0(\tau-5R,t)-\textnormal{Pal}_{u_0(\tau-5R,t)}v)$ and
   similarly for the term involving $\xi_1$, hence we have:
   $$|\xi_{\textnormal{new}}-\textnormal{Pal}_{u_R(\tau,t)}| \leq
   |\xi_0(\tau-5R,t) - \textnormal{Pal}_{u_0(\tau-5R,t)}v| +
   |\xi_1(\tau+5R,t) - \textnormal{Pal}_{u_1(\tau+5R,t)}v|.$$

   Using that fact that $\Vert d u_R \Vert_{L^{\infty}}$ is uniformly
   bounded (see Proposition~\ref{P:uR}), a straightforward computation
   shows that on $S[-5R,5R]$ we have:
   \begin{align*}
      & |\nabla(\xi_{\textnormal{new}} -
      \textnormal{Pal}_{u_R}v)(\tau,t)| \\
      & \leq C_8 \Bigl(
      |\xi_0(\tau-5R,t)-\textnormal{Pal}_{u_0(\tau-5R,t)}v| +
      |\xi_1(\tau+5R,t)-\textnormal{Pal}_{u_1(\tau+5R,t)}v| \\
      & + |\nabla(\xi_0(\tau-5R,t) -
      \textnormal{Pal}_{u_0(\tau-5R,t)}v)| + |\nabla(\xi_1(\tau+5R,t)
      - \textnormal{Pal}_{u_1(\tau+5R,t)}v)|\Bigr),
   \end{align*}
   where $C_8$ does not depend on $R$.  Therefore
   \begin{equation} \label{Eq:I1}
      \begin{split}
         I_1 = \int_{S[-5R,5R]} \alpha_{R,\delta}\bigl(&
         |\xi_{\textnormal{new}}-\textnormal{Pal}_{u_R}v|^p +
         |\nabla(\xi_{\textnormal{new}}-\textnormal{Pal}_{u_R}v)|^p
         \bigr) \\
         \leq C_9 \Bigl(\int_0^1 \int_{-5R}^{5R}
         \alpha_{R,\delta}(\tau) \bigl(& |\xi_0(\tau-5R,t) -
         \textnormal{Pal}_{u_0(\tau-5R,t)}v|^p \\
         & + |\nabla(\xi_0(\tau-5R,t) -
         \textnormal{Pal}_{u_0(\tau-5R,t)}v)|^p \bigr) d\tau dt \\
         + \int_0^1 \int_{-5R}^{5R} \alpha_{R,\delta}(\tau) \bigl( &
         |\xi_1(\tau+5R,t) -
         \textnormal{Pal}_{u_1(\tau+5R,t)}v|^p \\
         & + |\nabla(\xi_1(\tau+5R,t) -
         \textnormal{Pal}_{u_1(\tau+5R,t)}v)|^p \bigr) d \tau dt
         \Bigr) \\
         \leq C_{10}(\Vert \xi_0 \Vert_{1,p;s_{\delta}}^p + \Vert
         \xi_1 \Vert_{1,p;s_{\delta}}^p),
      \end{split}
   \end{equation}
   for some constants $C_9,C_{10}$. The last inequality here follows
   by comparing the weights $s_{\delta}$ and $\alpha_{R,\delta}$.

   Summing up, we obtain
   from~\eqref{Eq:I-estim},~\eqref{Eq:I2},~\eqref{Eq:xi05},~\eqref{Eq:I1}
   the desired inequality~\eqref{Eq:I}. This completes the proof of
   Proposition~\ref{P:estim-IR}, hence also of
   Proposition~\ref{P:estim-Qtilde}.
\end{proof}

\begin{proof}[Proof of Proposition~\ref{P:rinv}]
   Fix $\delta>0$ small enough (so that all the previous estimates
   hold for this $\delta$).  We first claim that $\Vert D'_{u_R} \circ
   \widetilde{Q}_R - \Id \Vert'_{0,p;\alpha_{R,\delta}} \to 0$ when $R
   \to \infty$.  To see this, we start with the following identities
   that easily follow from~\eqref{Eq:Q0}:
   \begin{align} \label{Eq:DuQ}
      D_{u_1,u_0} \circ Q_0 = \Id, \quad & D_{u_1,u_0} \circ A=0 \\
      d \mathcal{G}_{1,0} \circ Q_0=0, \quad & d \mathcal{G}_{1,0}
      \circ A = \Id. \notag
   \end{align}
   Let $\eta \in \mathcal{E}_{u_R}^{0,p;\delta}$, $\overrightarrow{a}
   \in E$. A simple computation gives:
   \begin{align} \label{Eq:DuRQR}
      D'_{u_R} \widetilde{Q}_R(\eta, \overrightarrow{a}) -
      (\eta,\overrightarrow{a}) = & \\
      \Bigl( & D_{u_R}\circ I^R \circ Q_0(J_{*,1}^{R}
      \sigma_0^{-}\eta, J_{*,0}^{R} \sigma_0^{+}\eta) - \eta +
      D_{u_R}\circ
      I^R \circ A(\overrightarrow{a}), \notag \\
      & d \mathcal{G}_R \circ I^R \circ Q_0(J_{*,1}^{R}
      \sigma_0^{-}\eta, J_{*,0}^{R} \sigma_0^{+}\eta) + d
      \mathcal{G}_R \circ I^R \circ
      A(\overrightarrow{a})-\overrightarrow{a}\Bigr). \notag
   \end{align}
   Using the definition of $I^R$ and~\eqref{Eq:DuQ} we have
   $$d \mathcal{G}_R \circ I^R \circ Q_0 = d \mathcal{G}_{1,0} \circ
   Q_0 = 0, \quad d \mathcal{G}_R \circ I^R \circ
   A(\overrightarrow{a})=d \mathcal{G}_{1,0} \circ
   A(\overrightarrow{a})=\overrightarrow{a}.$$
   Substituting this
   in~\eqref{Eq:DuRQR} we see that we have to estimate only the first
   component of~\eqref{Eq:DuRQR}. For the last term in the first
   component we have:
   \begin{align} \label{Eq:3trm-1}
      \Vert D_{u_R}\circ I^R \circ A(\overrightarrow{a})
      \Vert_{0,p;\alpha_{R,\delta}} \\
      \leq & \Vert D_{u_R}\circ I^R \circ A(\overrightarrow{a}) -
      (J_{1,*}^{R}D_{u_1} \oplus J_{0,*}^{R} D_{u_0}) \circ
      A(\overrightarrow{a}) \Vert_{0,p;\alpha_{R,\delta}} \notag \\
      & + \Vert (J_{1,*}^{R}D_{u_1} \oplus J_{0,*}^{R} D_{u_0}) \circ
      A(\overrightarrow{a}) \Vert_{0,p;\alpha_{R,\delta}}.  \notag
   \end{align}
   By~\eqref{Eq:DuQ} the last term in~\eqref{Eq:3trm-1} is $0$.
   Applying Proposition~\ref{P:estim-DuI-JD} to the first summand of
   inequality~\eqref{Eq:3trm-1} we get:
   \begin{equation} \label{Eq:3trm-2}
      \Vert D_{u_R}\circ I^R \circ A(\overrightarrow{a})
      \Vert_{0,p;\alpha_{R,\delta}} \leq C_{\ref{P:estim-DuI-JD}}(R)
      \Vert A \Vert |\overrightarrow{a}|
   \end{equation}
   for some constants $C_{\ref{P:estim-DuI-JD}}(R)$ that satisfy $\lim_{R
     \to \infty} C_{\ref{P:estim-DuI-JD}}(R)=0$.

   As for the first two terms in the first component
   of~\eqref{Eq:DuRQR} we have by Proposition~\ref{P:estim-DuI-JD}:
   \begin{align} \label{Eq:1trm}
      & \Vert D_{u_R} \circ \widetilde{Q}_R (\eta, \overrightarrow{0})
      - \eta \Vert_{0,p;\alpha_{R,\delta}} = \Vert D_{u_R} \circ
      I^{R}\circ Q_0(J_{*,1}^{R}\sigma_0^{-} \eta,
      J_{*,0}^{R}\sigma_0^{+}\eta)-\eta
      \Vert_{0,p;\alpha_{R,\delta}} \\
      \leq & \Vert(J_{1,*}^{R}\oplus J_{0,*}^{R})\circ D_{u_1,u_0}
      \circ Q_0(J_{*,1}^{R}\sigma_0^{-}\eta,
      J_{*,0}^{R}\sigma_0^{+}\eta)-\eta \Vert_{0,p;\alpha_{R,\delta}}
      \notag \\
      & + C_{\ref{P:estim-DuI-JD}}(R) \Vert Q_0(J_{*,1}^{R}
      \sigma_0^{-}\eta, J_{*,0}^{R}\sigma_0^{+}\eta)
      \Vert_{1,p;s_{\delta}} \notag \\
      \leq & \Vert J_{1,*}^{R} J_{*,1}^{R} \sigma_0^{-}\eta +
      J_{0,*}^{R} J_{*,0}^{R} \sigma_0^{+}\eta - \eta
      \Vert_{0,p;\alpha_{R,\delta}} \notag \\
      & + C_{\ref{P:estim-DuI-JD}}(R) \Vert Q_0 \Vert_{s_{\delta}}
      \bigl(\Vert J_{*,1}^{R} \sigma_0^{-} \eta \Vert_{0,p;s_{\delta}}
      + \Vert J_{*,0}^{R} \sigma_0^{+} \eta \Vert_{0,p;s_{\delta}}
      \bigr). \notag
     \end{align}
     By the 3'rd identity in~\eqref{Eq:id-3} the first term in the
     second to last line of~\eqref{Eq:1trm} is $0$, hence by
     Lemma~\ref{L:J-estim} we now get $$\Vert D_{u_R} \circ
     \widetilde{Q}_R \eta - \eta \Vert_{0,p;\alpha_{R,\delta}}\leq
     2C_{\ref{P:estim-DuI-JD}}(R) \Vert Q_0
     \Vert_{s_{\delta}}C_{\ref{L:J-estim}} \Vert\eta
     \Vert_{0,p;\alpha_{R,\delta}}.$$
     This together
     with~\eqref{Eq:3trm-2} proves that $\lim_{R\to \infty} \Vert
     D'_{u_R} \circ \widetilde{Q}_R - \Id
     \Vert'_{0,p;\alpha_{R,\delta}}=0$.

     Let $R_0>0$ be large enough so that $\Vert D'_{u_R}\circ
     \widetilde{Q}_R - \Id \Vert'_{0,p;\alpha_{R,\delta}} \leq
     \tfrac{1}{2}$ for every $R \geq R_0$. Define
     $$Q_R:\mathcal{E}_{u_R}^{0,p;\delta} \oplus E \to
     T_{u_R}^{0,p;\delta}, \qquad Q_R=\widetilde{Q}_R \circ
     \sum_{k=0}^{\infty}(\Id - D'_{u_R} \circ \widetilde{Q}_R)^k.$$
     Clearly $Q_R$ is a right inverse to $D'_{u_R}$ and $\Vert Q_R
     \Vert_{\alpha_{R,\delta}} \leq 2 \Vert \widetilde{Q}_R
     \Vert_{\alpha_{R,\delta}} \leq 2C_{\ref{P:estim-Qtilde}}(\delta)$
     for $R \gg 0$, i.e. $Q_R$ is uniformly bounded.
\end{proof}

\subsection{The implicit function theorem}
\label{S:implicit}
\begin{thm}[See Proposition~A.3.4 from~\cite{McD-Sa:Jhol-2}]
   \label{T:implicit}
   Let $X,Y$ be Banach spaces, $\mathcal{U} \subset X$ an open subset
   and $f:\mathcal{U} \to Y$ a $C^1$-map. Let $x_0\in \mathcal{U}$ be
   such that $df(x_0):X \to Y$ is surjective and has a bounded right
   inverse $Q:Y \to X$. Let $r, K>0$ be such that $B_{x_0}(r) \subset
   \mathcal{U}$, $\Vert Q \Vert \leq K$ and such that:
   \begin{enumerate}
     \item $\Vert x-x_0 \Vert < r \Longrightarrow \Vert df(x) -
      df(x_0) \Vert < \frac{1}{2K}$.
     \item $\Vert f(x_0) \Vert < \frac{r}{4K}$.
   \end{enumerate}
   Then there exists $x \in X$ with the following properties:
   \begin{enumerate}
     \item $f(x)=0$.
     \item $\Vert x-x_0 \Vert \leq r$.
   \end{enumerate}
   In fact we have $\Vert x-x_0 \Vert \leq 2K \Vert f(x_0) \Vert$.
\end{thm}

\subsubsection{Quadratic estimates}  \label{Sb:quadratic} In
the following proposition we endow the space
$\mathcal{E}_{u_R}^{0,p;\delta} \oplus E$ with the norm $\Vert(\eta,
\overrightarrow{a}) \Vert'_{0,p;\alpha_{R,\delta}} = \Vert \eta
\Vert_{0,p;\alpha_{R,\delta}} + |\overrightarrow{a}|$.

\begin{prop}[See Proposition~3.5.3 and Remark~3.5.5
   in~\cite{McD-Sa:Jhol-2}]
   \label{P:quadratic}
   For every $c_0>0$ there exist constants
   $C_{\ref{P:quadratic}}=C_{\ref{P:quadratic}}(c_0)>0$,
   $R_0=R_0(c_0)>1$, $\delta_0=\delta_0(c_0)>0$ such that for every
   $0<\delta\leq \delta_0$, $R\geq R_0$, $\xi \in
   T_{u_R}^{1,p;\delta}$ with $\Vert \xi \Vert_{L^{\infty}}\leq c_0$
   we have:
   \begin{equation} \label{Eq:quadratic-1}
      \Vert d \mathcal{F}'_{u_R}(\xi) \xi' - D'_{u_R} \xi'
      \Vert'_{0,p;\alpha_{R,\delta}} \leq C_{\ref{P:quadratic}} \Vert\xi
      \Vert_{1,p;\alpha_{R,\delta}} \Vert \xi' \Vert_{1,p;
        \alpha_{R,\delta}}, \quad \forall \xi' \in T_{u_R}^{1,p;\delta}.
   \end{equation}
   In other words, the norm of the operator $d
   \mathcal{F}'_{u_R}(\xi)-D'_{u_R}$ (with respect to $\Vert \cdot
   \Vert_{1,p; \alpha_{R,\delta}}$ and $\Vert \cdot
   \Vert'_{0,p;\alpha_{R,\delta}}$) satisfies $\Vert d
   \mathcal{F}'_{u_R}(\xi) - D'_{u_R} \Vert \leq C_{\ref{P:quadratic}}
   \Vert\xi \Vert_{1,p;\alpha_{R,\delta}}$ whenever $\Vert \xi
   \Vert_{L^{\infty}} \leq c_0$.

   Moreover, whenever $\xi_0, \xi \in T_{u_R}^{1,p;\delta}$ satisfy
   $\Vert \xi \Vert_{L^{\infty}}, \Vert \xi_0 \Vert_{L^{\infty}} \leq
   c_0$ we have for every $R \geq R_0$, $0<\delta \leq \delta_0$:
   \begin{equation} \label{Eq:quadratic-2}
      \Vert \mathcal{F}'_{u_R}(\xi_0+\xi) - \mathcal{F}'_{u_R}(\xi_0) -
      d \mathcal{F}'_{u_R}(\xi_0)\xi \Vert'_{0,p;\alpha_{R,\delta}}
      \leq C_{\ref{P:quadratic}} \Vert \xi \Vert_{L^{\infty}}
      \Vert \xi \Vert_{1,p;\alpha_{R,\delta}}.
   \end{equation}

\end{prop}

\begin{proof}[Outline of the proof]
   The proof of the estimate~\eqref{Eq:quadratic-1} is essentially
   given in~\cite{McD-Sa:Jhol-2} (Proposition~3.5.3). To adapt that
   proof to our case we use the following additional ingredients:
   \begin{enumerate}
     \item In~\cite{McD-Sa:Jhol-2} the pointwise norms $|\cdot|$ are
      taken with respect to the metric $g_{\omega,J}$ but we work with
      the metric $g_{\omega,J,L}$. However $L$ is compact and
      $g_{\omega,J,L}=g_{\omega,J}$ outside a small neighbourhood of
      $L$ hence the two pointwise norms are comparable.
     \item $\Vert du_R \Vert_{0,p;\alpha_{R,\delta}}$, $\Vert d u_R
      \Vert_{L^{\infty}}$ are uniformly bounded in $R$ by
      Proposition~\ref{P:uR}.
     \item Due to Proposition~\ref{P:comparison} we have $\Vert \cdot
      \Vert_{L^{\infty}} \leq C_{\ref{P:comparison}} \Vert\cdot
      \Vert_{1,p;\alpha_{R,\delta}}$.
     \item Standard arguments involving ode's
      show that for every $\xi_{\infty} \in T_{u_1(\infty)}(L)$ we
      have $|\nabla P_{u_R} \xi_{\infty}(\tau,t)| \leq C_1
      e^{-c_1(\tau-5R)} |\xi_{\infty}|$ for every $\tau\geq 5R$ for
      some constants $C_1,c_1$ that do not depend on $R$. An analogous
      estimate holds for $|\nabla P_{u_R} \xi_{-\infty}(\tau,t)|$.
   \end{enumerate}
   It follows that $|\nabla \xi (\tau,t)| \leq
   |\nabla(\xi-P_{u_R}\xi_{\infty}(\tau,t))| +
   C_1e^{-c_1(\tau-5R)}|\xi_{\infty}|$ for every $\tau\geq 5R$. By a
   simple computation it now follows that $\int_{S[5R, \infty)}
   \alpha_{R,\delta} |\nabla\xi |^p \leq C_2 \Vert \xi
   \Vert^p_{1,p;\alpha_{R,\delta}}$ for some constant $C_2$ that does
   not depend on $R$. A similar estimate holds for the integral over
   $S(-\infty, -5R]$.  The rest of the proof is as
   in~\cite{McD-Sa:Jhol-2} using the additional estimates 1-4.  The
   proof of estimate~\eqref{Eq:quadratic-2} is similar (see
   Remark~3.5.5 in~\cite{McD-Sa:Jhol-2}).
\end{proof}

\subsubsection{Proof of Theorem~\ref{T:gluing}}

\begin{proof}[Proof of the existence statement of Theorem~\ref{T:gluing}]
   Fix $\delta>0$ small enough so that all the previous estimates hold
   for this $\delta$. By Proposition~\ref{P:rinv} there exists $K>0$
   such that $\Vert Q_R \Vert_{\alpha_{R,\delta}} \leq K$ for $R \gg
   0$.  Applying Proposition~\ref{P:quadratic} with $c_0=K$ we obtain
   a constant $C(K)$ such that for every $R \gg 0$ we have
   \begin{equation} \label{Eq:imp-estim} \Vert \xi \Vert_{L^{\infty}}
      \leq C(K)
      \Longrightarrow \Vert d \mathcal{F}_{u_R}'(\xi)-D_{u_R}'
      \Vert_{\alpha_{R,\delta}} \leq C(K) \Vert \xi
      \Vert_{1,p;\alpha_{R,\delta}}.
   \end{equation}
   By Proposition~\ref{P:comparison}, there exists
   $C_{\ref{P:comparison}}>0$ such that $\Vert \xi \Vert_{L^{\infty}}
   \leq C_{\ref{P:comparison}} \Vert \xi
   \Vert_{1,p;\alpha_{R,\delta}}$ for every $R>0$, $\xi \in
   T_{u_R}^{1,p;\delta}$. Fix $0<r_0< \min \bigl\{\tfrac{1}{2KC(K)},
   \tfrac{K}{C_{\ref{P:comparison}}} \bigr\}$. If $\Vert \xi
   \Vert_{1,p;\alpha_{R,\delta}}<r_0$ then $\Vert \xi
   \Vert_{L^{\infty}} \leq C_{\ref{P:comparison}} r_0 < K$ hence
   $$\Vert d \mathcal{F}_{u_R}'(\xi)-D_{u_R}'
   \Vert_{\alpha_{R,\delta}} \leq C(K) \Vert \xi
   \Vert_{1,p;\alpha_{R,\delta}}<\frac{1}{2K}.$$
   Note that for $R \gg
   0$, $\mathcal{F}_{u_R}'(0)=\bigl(\overline{\partial}_J u_R,0
   \bigr)$, hence by Proposition~\ref{P:uR}
   \begin{equation} \label{Eq:FuR}
      \Vert \mathcal{F}_{u_R}'(0) \Vert'_{0,p;\alpha_{R,\delta}} =
      \Vert \overline{\partial}_J u_R \Vert_{0,p;\alpha_{R,\delta}} \leq
      C'_{\ref{P:uR}}e^{-c'_{\ref{P:uR}} R}
   \end{equation}
   for some constants $C'_{\ref{P:uR}}, c'_{\ref{P:uR}}>0$.  Thus by
   taking $R \gg 0$ we may assume that $\Vert \mathcal{F}_{u_R}'(0)
   \Vert'_{0,p;\alpha_{R,\delta}} \leq \tfrac{r_0}{4K}$. We now apply
   Theorem~\ref{T:implicit} with:
   $$\mathcal{U}=X=\bigl(T^{1,p;\delta}_{u_R},\Vert \cdot
   \Vert_{1,p;\alpha_{R,\delta}}\bigr), \quad
   Y=\bigl(\mathcal{E}^{0,p;\delta}_{u_R}\oplus E, \Vert \cdot
   \Vert'_{0,p;\alpha_{R,\delta}} \bigr), \quad x_0=0, \quad
   f=\mathcal{F}'_{u_R}.$$
   By Theorem~\ref{T:implicit} we obtain
   $\xi_R \in T_{u_R}^{1,p;\delta}$ with $\Vert \xi_R
   \Vert_{1,p;\alpha_{R,\delta}} \leq r_0$ such that
   $\mathcal{F}_{u_R}'(\xi_R)=(0,0)$, i.e.  $v_R=\exp_{u_R}(\xi_R)$
   satisfies:
   \begin{enumerate}
     \item $\overline{\partial}_J v_R=0$.
     \item $G_{\pm}(v_R(\pm \infty))=0$.
     \item $G_1(v_R(-5R-\tau_*,\tfrac{1}{2}))=0$,
      $G_0(v_R(5R+\tau_*,\tfrac{1}{2}))=0$.
   \end{enumerate}
   By elliptic regularity $v_R$ is actually smooth. Note that $v_R$
   has finite energy (see Remark~\ref{R:analytic}-(\ref{Rn:energy}))
   hence by the removal of boundary singularities theorem of
   Oh~\cite{Oh:rsing}, $v_R \circ \lambda^{-1}$ extends smoothly to
   $D$.

   Next, note that by Theorem~\ref{T:implicit} we actually have $\Vert
   \xi_R \Vert_{1,p;\alpha_{R,\delta}} \leq 2K \Vert
   \mathcal{F}'_{u_R}(0) \Vert'_{0,p;\alpha_{R,\delta}}$.
   From~\eqref{Eq:FuR} we obtain $\Vert \xi_R
   \Vert_{1,p;\alpha_{R,\delta}}\xrightarrow[R \to \infty]{} 0$. It
   follows from Proposition~\ref{P:comparison} that we also have
   $\Vert \xi_R \Vert_{L^{\infty}} \xrightarrow[R \to \infty]{} 0$.
   Clearly for $R \gg 0$, $[v_R]=[u_R]=A_1+A_0 \in H_2(M,L)$.

   We claim that for $R \gg 0$ we have $v_R(\pm \infty) \in W_{\pm}$,
   $v_R(-5R-\tau_*,\tfrac{1}{2}) \in W_1, v_R(5R+\tau_*,\tfrac{1}{2})
   \in W_0$.  Indeed, since $\Vert \xi_R \Vert_{L^{\infty}}
   \xrightarrow[R \to \infty]{} 0$ then for $R \gg 0$ the points
   $v_R(-\infty)$, $v_R(-5R-\tau_*,\tfrac{1}{2})$,
   $v_R(5R+\tau_*,\tfrac{1}{2})$, $v_R(\infty)$ lie in arbitrarily
   small neighbourhoods of $u_1(-\infty)$,
   $u_1(-\tau_*,\tfrac{1}{2})$, $u_0(\tau_*,\tfrac{1}{2})$,
   $u_0(\infty)$ respectively. By the definition of $G_{\pm}, G_1,
   G_0$ (see Section~\ref{S:main-op}) it follows that $v_R(\pm \infty)
   \in W_{\pm}$, $v_R(-5R-\tau_*,\tfrac{1}{2}) \in W_1,
   v_R(5R+\tau_*,\tfrac{1}{2}) \in W_0$.

   Put $s_R = \lambda^{-1}(5R+\tau_*, \tfrac{1}{2}) \in D$.  Then $s_R
   \xrightarrow[R \to \infty]{} 1$ and $(v_R \circ \lambda^{-1}, s_R)
   \in \mathcal{M}(A,J;\mathcal{C}(\mathbf{h}))$ for $R \gg 0$.

   To complete the proof of the existence part of the theorem it
   remains to show that $v_R$ together with the marked points
   $(-\infty, (-5R-\tau_*,\tfrac{1}{2}), (5R+\tau_*,\tfrac{1}{2}),
   \infty)$ converges to $(u_1,u_0)$ with the marked points $(-\infty,
   (-\tau_*,\tfrac{1}{2})),\, ((\tau_*,\tfrac{1}{2}), \infty)$ in the
   Gromov topology as $R \to \infty$. (See~\cite{Fr:msc} for more
   details on Gromov compactness for disks.)

   For this aim we have to find biholomorphisms $\rho^1_R,
   \rho^0_R:\widehat{S} \to \widehat{S}$ such that:
   \begin{enumerate}
     \item[($v$-1)] $v_R \circ \rho^1_R \xrightarrow[R \to \infty]{}
      u_1$, $v_R \circ \rho^0_R \xrightarrow[R \to \infty]{} u_0$
      uniformly on compact subsets of $\widehat{S} \setminus \{
      \infty\}$, of $\widehat{S} \setminus \{ -\infty\}$ respectively.
     \item [($v$-2)] $(\rho^1_R)^{-1} \circ \rho^0_R \xrightarrow[R
      \to \infty]{} \infty$, $(\rho^0_R)^{-1} \circ \rho^1_R
      \xrightarrow[R \to \infty]{} -\infty$ uniformly on compact
      subsets of $\widehat{S} \setminus \{-\infty \}$, of $\widehat{S}
      \setminus \{ \infty \}$ respectively.
     \item [($v$-3)] $(\rho^1_R)^{-1}(-\infty) \xrightarrow[R \to
      \infty]{} -\infty$, $(\rho^1_R)^{-1}(-5R-\tau_*,\tfrac{1}{2})
      \xrightarrow[R \to \infty]{} (-\tau_*,\tfrac{1}{2})$.
     \item [($v$-4)] $(\rho^0_R)^{-1}(\infty) \xrightarrow[R \to
      \infty]{} \infty$, $(\rho^0_R)^{-1}(5R+\tau_*,\tfrac{1}{2})
      \xrightarrow[R \to \infty]{} (\tau_*,\tfrac{1}{2})$.
   \end{enumerate}

   To prove this define $\rho^1_R(\tau,t)=(\tau-5R,t)$,
   $\rho^0_R(\tau,t)=(\tau+5R,t)$. Clearly $\rho^1_R,\rho^0_R$ satisfy
   properties ($v$-2)--($v$-4) above. Property ($v$-1) follows easily
   from the definition of $u_R$ (see~\eqref{Eq:uR} in
   Section~\ref{S:pregluing}) and the fact that
   $v_R=\exp_{u_R}(\xi_R)$ with $\Vert \xi_R \Vert_{L^{\infty}}
   \xrightarrow[R \to \infty]{} 0$. This concludes the proof of the
   existence statement in Theorem~\ref{T:gluing}.

\end{proof}

\begin{proof}[Proof of the uniqueness statement of
   Theorem~\ref{T:gluing}] We continue to use the notation introduced
   above in the proof of the existence part.

   Let $(w_n,\tau_n) \in \mathcal{M}(A,J;\mathcal{C}(\mathbf{h}))$ be
   a sequence that converges with the marked points $(-\infty,
   (-\tau_n,\tfrac{1}{2}), (\tau_n,\tfrac{1}{2}), \infty)$ to $(u_1,
   u_0)$ with the marked points $(-\infty, (-\tau_*,\tfrac{1}{2})),\,
   ((\tau_*,\tfrac{1}{2}), \infty)$ in the Gromov topology. We shall
   prove that for every large enough $n$ we have $(w_n, \tau_n) =
   (v_{R_n},5R_n+\tau_*)$ for some $R_n$'s with $R_n \xrightarrow[n
   \to \infty]{} \infty$.

   The implicit function Theorem~\ref{T:implicit} is not enough to
   prove this statement although there is an extension of
   Theorem~\ref{T:implicit} which states that the solution $x$ is
   unique among solutions that satisfy $x \in \textnormal{image\,}Q$
   and $x \in B_{x_0}(r)$. In our case this would require to prove
   that $w_n=\exp_{u_{R_n}}(\xi_n^{(w)})$ for some $\xi_n^{(w)} \in
   T_{u_{R_n}}^{1,p;\delta}$ with $\Vert \xi_n^{(w)}
   \Vert_{1,p;\alpha_{R_n,\delta}} \xrightarrow[n \to \infty]{} 0$.
   Instead, we shall prove a similar statement in the $C^0$-topology
   and use the quadratic estimates of Proposition~\ref{P:quadratic} to
   deduce the uniqueness.

   We first fix $\delta>0$ small enough so that all the previous
   estimates hold for this $\delta$ and such that all the operators
   $D'_{u_R}$, $R \gg 0$, are Fredholm for this $\delta$ (see
   Section~\ref{Sb:fredholm-2}.) Note that $u_R$, $v_R$, $w_n$, all
   extend smoothly to $\widehat{S} \approx D$ hence they belong to
   $W^{1,p;\delta}$ for {\em every} small enough $\delta>0$. (See
   Remark~\ref{R:analytic}-(\ref{Rn:smooth-u}).)

   By the definition of Gromov convergence with marked points (see
   e.g.~\cite{Fr:msc}) we have two sequences of biholomorphisms
   $\varphi^1_n, \varphi^0_n : \widehat{S} \to \widehat{S}$ such that:
   \begin{enumerate}
     \item[($w$-1)] $w_n \circ \varphi^1_n \xrightarrow[n \to
      \infty]{} u_1$, $w_n \circ \varphi^0_n \xrightarrow[n \to
      \infty]{} u_0$ uniformly on compact subsets of $\widehat{S}
      \setminus \{ \infty \}$, of $\widehat{S} \setminus \{ -\infty
      \}$ respectively.
     \item[($w$-2)] $(\varphi^1_n)^{-1}\circ \varphi^0_n
      \xrightarrow[n \to \infty]{} \infty$, $(\varphi^0_n)^{-1}\circ
      \varphi^1_n \xrightarrow[n \to \infty]{} -\infty$ uniformly on
      compact subsets of $\widehat{S} \setminus \{ -\infty \}$, of
      $\widehat{S} \setminus \{ \infty \}$ respectively.
     \item[($w$-3)] $(\varphi^1_n)^{-1}(-\infty) \xrightarrow[n \to
      \infty]{} -\infty$, $(\varphi^0_n)^{-1}(\infty) \xrightarrow[n
      \to \infty]{} \infty$.
     \item[($w$-4)] $(\varphi^1_n)^{-1}(-\tau_n,\tfrac{1}{2})
      \xrightarrow[n \to \infty]{} (-\tau_*,\tfrac{1}{2})$,
      $(\varphi^0_n)^{-1}(\tau_n,\tfrac{1}{2}) \xrightarrow[n \to
      \infty]{} (\tau_*,\tfrac{1}{2})$.
   \end{enumerate}
   Define $$R_n = \frac{\tau_n-\tau_*}{5}, \quad \textnormal{so that }
   \, 5R_n+\tau_*=\tau_n.$$
   We have:
   \begin{itemize}
     \item $(\varphi^1_n)^{-1} \circ \rho^1_{R_n}(-\infty)
      \xrightarrow[n \to \infty]{} -\infty$, $(\varphi^1_n)^{-1} \circ
      \rho^1_{R_n}(-\tau_*,\tfrac{1}{2}) \xrightarrow[n \to \infty]{}
      (-\tau_*,\tfrac{1}{2})$.
     \item $(\varphi^0_n)^{-1} \circ \rho^0_{R_n}(\infty)
      \xrightarrow[n \to \infty]{} \infty$, $(\varphi^1_n)^{-1} \circ
      \rho^0_{R_n}(\tau_*,\tfrac{1}{2}) \xrightarrow[n \to \infty]{}
      (\tau_*,\tfrac{1}{2})$.
   \end{itemize}
   A simple argument involving M\"{o}bius transformations implies that
   \begin{equation} \label{Eq:phi-rho}
      (\varphi^1_n)^{-1} \circ \rho^1_{R_n} \xrightarrow[n \to \infty]{}
      \Id, \quad (\varphi^0_n)^{-1} \circ \rho^0_{R_n} \xrightarrow[n \to
      \infty]{} \Id \quad \textnormal{uniformly.}
   \end{equation}
   Moreover, it follows from ($w$-2) that $R_n, \tau_n \xrightarrow[n
   \to \infty]{} \infty$.

   We now have the following
   \begin{lem} \label{L:duw-1}
      $\sup_{z \in \widehat{S}} \textnormal{dist} \bigl( w_n(z),
      u_{R_n}(z) \bigr) \xrightarrow[n\to \infty]{} 0$.
   \end{lem}
   We defer the proof of this lemma to Section~\ref{Sb:prf-lem-duw-1}
   and continue with the proof of our theorem. By Lemma~\ref{L:duw-1},
   for $n \gg 1$ there exist $\xi^{(w)}_n \in
   T^{1,p;\delta}_{u_{R_n}}$ with $\Vert \xi^{(w)}_n
   \Vert_{L^{\infty}} \xrightarrow[n\to \infty]{} 0$ such that $w_n =
   \exp_{u_{R_n}}(\xi^{(w)}_n)$. We shall prove that
   $\xi^{(w)}_n=\xi_{R_n}$ for $n \gg 1$, hence $w_n=v_{R_n}$.

   We start with the observation that for $n \gg 1$ the operator
   $Q_{R_n}:\mathcal{E}^{0,p;\delta}_{u_{R_n}} \oplus E \to
   T^{1,p;\delta}_{u_{R_n}}$ is bijective.  To see this, recall that
   $D'_{u_{R_n}}:T^{1,p;\delta}_{u_{R_n}} \to
   \mathcal{E}^{0,p;\delta}_{u_{R_n}} \oplus E$ is Fredholm and its
   index is $\mu(A_1+A_0) + k_- + k_+ + k_1 + k_0 - 5\dim L$ which by
   our assumptions is $0$. As $D'_{u_{R_n}}$ is surjective for $n \gg
   1$ it must be bijective, hence its inverse $Q_{u_{R_n}}$ is
   bijective too.


   To finish the proof, we use an argument due to McDuff and Salamon
   (see the proof of Corollary~3.5.6 in~\cite{McD-Sa:Jhol-2}). Put
   $\xi'_n = \xi^{(w)}_n - \xi_{R_n}$.  Since $Q_{u_{R_n}}$ is
   bijective we have $\xi'_n=Q_{u_{R_n}} \circ D'_{u_{R_n}} \xi'_n$
   and $\mathcal{F}'_{u_{R_n}}(\xi_{R_n})=0$,
   $\mathcal{F}'_{u_{R_n}}(\xi_{R_n}+\xi'_n)=0$.  By
   Proposition~\ref{P:rinv} there exists a constant
   $C_{\ref{P:rinv}}(\delta)$ such that:
   \begin{align} \label{Eq:xin-1}
      \Vert \xi'_n \Vert_{1,p;\alpha_{R_n,\delta}} \leq
      &C_{\ref{P:rinv}}(\delta)
      \Vert D'_{u_{R_n}} \xi'_n \Vert'_{0,p;\alpha_{R_n,\delta}} \\
      \leq &C_{\ref{P:rinv}}(\delta) \Vert
      \mathcal{F}'_{u_{R_n}}(\xi_{R_n}+\xi'_n) -
      \mathcal{F}'_{u_{R_n}}(\xi_{R_n}) -
      d\mathcal{F}'_{u_{R_n}}(\xi_{R_n})\xi'_n
      \Vert'_{0,p;\alpha_{R_n,\delta}} \notag \\
      & + C_{\ref{P:rinv}}(\delta) \Vert \bigl(
      d\mathcal{F}'_{u_{R_n}}(\xi_{R_n})-D'_{u_{R_n}} \bigr) \xi'_n
      \Vert'_{0,p;\alpha_{R_n,\delta}} \notag.
   \end{align}
   Applying Proposition~\ref{P:quadratic} to both terms of the very
   right-hand side of~\eqref{Eq:xin-1} we get
   \begin{align} \label{Eq:xin-2}
      \Vert \xi'_n \Vert_{1,p;\alpha_{R_n,\delta}} & \leq
      C_{\ref{P:rinv}}(\delta) C_{\ref{P:quadratic}} \bigl( \Vert
      \xi'_n \Vert_{L^{\infty}} + \Vert \xi_{R_n}
      \Vert_{1,p;\alpha_{R_n,\delta}} \bigr) \Vert
      \xi'_n \Vert_{1,p;\alpha_{R_n,\delta}} \\
      & \leq C_{\ref{P:rinv}}(\delta) C_{\ref{P:quadratic}} \bigl(
      \Vert \xi^{(w)}_n \Vert_{L^{\infty}} + \Vert \xi_{R_n}
      \Vert_{L^{\infty}} + \Vert \xi_{R_n}
      \Vert_{1,p;\alpha_{R_n,\delta}} \bigr) \Vert \xi'_n
      \Vert_{1,p;\alpha_{R_n,\delta}}. \notag
   \end{align}

   Recall that $\Vert \xi^{(w)}_n \Vert_{L^{\infty}}, \Vert \xi_{R_n}
   \Vert_{L^{\infty}}, \Vert \xi_{R_n} \Vert_{1,p;\alpha_{R_n,\delta}}
   \xrightarrow[n\to \infty]{} 0$. Therefore for $n \gg 1$,
   $$C_{\ref{P:rinv}}(\delta) C_{\ref{P:quadratic}} \bigl( \Vert
   \xi^{(w)}_n \Vert_{L^{\infty}} + \Vert \xi_{R_n} \Vert_{L^{\infty}}
   + \Vert \xi_{R_n} \Vert_{1,p;\alpha_{R_n,\delta}} \bigr) < 1.$$
   It
   immediately follows from~\eqref{Eq:xin-2} that $\Vert \xi'_n
   \Vert_{1,p;\alpha_{R_n,\delta}} = 0$ for $n \gg 1$. This completes
   the proof of the uniqueness statement of Theorem~\ref{T:gluing}
   (modulo the proof of Lemma~\ref{L:duw-1}).
\end{proof}

\subsubsection{Proof of Lemma~\ref{L:duw-1}}
\label{Sb:prf-lem-duw-1} We first claim that for every $a \in
\mathbb{R}$:
\begin{align}
   & \sup_{z \in \widehat{S}[-\infty,-5R_n+a]} \textnormal{dist}
   \bigl( w_n(z), u_{R_n}(z) \bigr) \xrightarrow[n\to \infty]{} 0,
   \label{Eq:wn-uR-1}
   \\
   & \sup_{z \in \widehat{S}[5R_n-a, \infty]} \textnormal{dist} \bigl(
   w_n(z), u_{R_n}(z) \bigr) \xrightarrow[n\to \infty]{} 0.
   \label{Eq:wn-uR-2}
\end{align}
\begin{proof}[Proof of~\eqref{Eq:wn-uR-1},~\eqref{Eq:wn-uR-2}]
   We start with the inequality:
   \begin{align} \label{Eq:wn-vn}
      & \sup_{z \in \rho^1_{R_n}(\widehat{S}[-\infty,a])}
      \textnormal{dist} \bigl( w_n(z), v_{R_n}(z) \bigr) \\
      \leq & \sup_{z \in \rho^1_{R_n}(\widehat{S}[-\infty,a])}
      \textnormal{dist} \bigl( w_n(z), u_1\circ (\varphi^1_n)^{-1}(z)
      \bigr) + \sup_{z \in \rho^1_{R_n}(\widehat{S}[-\infty,a])}
      \textnormal{dist} \bigl( v_{R_n}(z), u_1\circ
      (\rho^1_{R_n})^{-1}(z) \bigr) \notag \\
      + & \sup_{z \in \rho^1_{R_n}(\widehat{S}[-\infty,a])}
      \textnormal{dist} \bigl(u_1\circ (\varphi^1_n)^{-1}(z), u_1\circ
      (\rho^1_{R_n})^{-1}(z) \bigr). \notag
   \end{align}
   The last term in the right-hand side of this inequality equals
   $$\sup_{z \in \widehat{S}[-\infty,a]} \textnormal{dist}
   \bigl(u_1\circ (\varphi^1_n)^{-1}\circ \rho^1_{R_n}(z), u_1(z)
   \bigr),$$
   hence by~\eqref{Eq:phi-rho} it tends to $0$ as $n\to
   \infty$. Similarly by ($v$-1) the second term in the right-hand
   side of~\eqref{Eq:wn-vn} tends to $0$ too as $n \to \infty$.  As
   for the first term, due to~\eqref{Eq:phi-rho} we have $\rho^1_{R_n}
   (\widehat{S}[-\infty,a]) \subset
   \varphi^1_n(\widehat{S}[-\infty,a+1])$ for $n \gg 1$, therefore
   \begin{align*}
      & \sup_{z \in \rho^1_{R_n}(\widehat{S}[-\infty,a])}
      \textnormal{dist} \bigl( w_n(z), u_1\circ (\varphi^1_n)^{-1}(z)
      \bigr) \\
      & \leq \sup_{z \in \varphi^1_n(\widehat{S}[-\infty,a+1])}
      \textnormal{dist} \bigl( w_n(z), u_1\circ (\varphi^1_n)^{-1}(z)
      \bigr) = \sup_{z \in \widehat{S}[-\infty,a+1]} \textnormal{dist}
      \bigl( w_n\circ \varphi^1_n(z), u_1(z) \bigr) \xrightarrow[n\to
      \infty]{} 0.
   \end{align*}
   It follows that
   $$\sup_{z \in \widehat{S}[-\infty,-5R_n+a]} \textnormal{dist}
   \bigl( w_n(z), v_{R_n}(z) \bigr)= \sup_{z \in
     \rho^1_{R_n}(\widehat{S}[-\infty,a])} \textnormal{dist} \bigl(
   w_n(z), v_{R_n}(z) \bigr) \xrightarrow[n\to \infty]{} 0.$$
   But
   $\sup_{z \in \widehat{S}} \textnormal{dist} \bigl( v_{R_n}(z),
   u_{R_n}(z) \bigr) \xrightarrow[n\to \infty]{} 0$,
   hence~\eqref{Eq:wn-uR-1} follows.  The proof of~\eqref{Eq:wn-uR-2}
   is analogous. This concludes the proof
   of~\eqref{Eq:wn-uR-1},~\eqref{Eq:wn-uR-2}.
\end{proof}

In the rest of the proof we shall need to use some estimates
concerning energy of pseudo-holomorphic strips. Given a subset $T
\subset \widehat{S}$ and a pseudo-holomorphic map $u$ defined in a
neighbourhood of $T$ we denote by $E(u;T) = \int_{T} u^*\omega$ the
energy of $u$ along $T$. We first claim that:
\begin{align}
   & \lim_{d \to \infty} \lim_{n \to \infty} E\bigl(w_n;
   \widehat{S}[-5R_n+d, 5R_n-d]\bigr)=0, \label{Eq:limE-1} \\
   & \lim_{d \to \infty} \lim_{n \to \infty} E\bigl(v_{R_n};
   \widehat{S}[-5R_n+d, 5R_n-d]\bigr)=0. \label{Eq:limE-2}
\end{align}
\begin{proof}[Proof of~\eqref{Eq:limE-1},~\eqref{Eq:limE-2}]
   To prove~\eqref{Eq:limE-1}, write
   $E_{d,n}=E\bigl(w_n;\widehat{S}[-5R_n+d, 5R_n-d]\bigr)$.  Recall
   that $A_1,A_0 \in H_2(M,L)$ are the classes of the disks $u_1, u_0$
   respectively. We have:
   \begin{align} \label{Eq:Edn-1}
      E_{d,n} = \int_{A_1+A_0} \omega & -
      E(w_n;\widehat{S}[-\infty,-5R_n+d]) - E(w_n;\widehat{S}[5R_n-d,
      \infty]) \\
      = \int_{A_1+A_0}\omega & - E(w_n\circ
      \varphi^1_n;(\varphi^1_n)^{-1}
      \circ \rho^1_{R_n}(\widehat{S}[-\infty,d])) \notag \\
      & - E(w_n\circ \varphi^0_n;(\varphi^0_n)^{-1} \circ
      \rho^0_{R_n}(\widehat{S}[-d, \infty])) \notag
   \end{align}
   By~\eqref{Eq:phi-rho} we have for $n \gg 1$:
   \begin{align*}
      & (\varphi^1_n)^{-1} \circ \rho^1_{R_n}(\widehat{S}[-\infty,d])
      \supset \widehat{S}[-\infty,d-1], \\
      & (\varphi^0_n)^{-1} \circ \rho^0_{R_n}(\widehat{S}[-d,\infty])
      \supset \widehat{S}[-d+1,\infty].
   \end{align*}
   Putting this together with~\eqref{Eq:Edn-1} we obtain:
   \begin{align*}
      E_{d,n} \leq \int_{A_1}\omega - E(w_n\circ
      \varphi^1_n;\widehat{S}[-\infty,d-1]) + \int_{A_0}\omega -
      E(w_n\circ \varphi^0_n;\widehat{S}[-d+1,\infty])
   \end{align*}
   It follows from~($w$-1) that
   $$\lim_{n \to \infty} E_{d,n} \leq E(u_1;\widehat{S}[d-1,\infty]) +
   E(u_0;\widehat{S}[-\infty, -d+1]),$$
   hence $\lim_{d \to \infty}
   \lim_{n \to \infty} E_{d,n}=0$. This proves~\eqref{Eq:limE-1}. The
   proof of~\eqref{Eq:limE-2} is similar.
\end{proof}

We shall now need the following estimate relating the ``variation'' of
a pseudo-holomorphic strip to its energy:
\begin{lem}[See Lemma A.6 in~\cite{Fr:msc},
   Lemma~4.7.3 in~\cite{McD-Sa:Jhol-2}] \label{L:E-estimate}

   There exist positive constants $\delta_{\ref{L:E-estimate}} =
   \delta_{\ref{L:E-estimate}}(M,\omega,J)$, $C_{\ref{L:E-estimate}}$
   such that every $J$-holomorphic curve $u:S[r_1,r_2] \to M$ with $u
   \bigl( [r_1, r_2] \times \{0\} \cup [r_1, r_2] \times \{1\} \bigr)
   \subset L$ and with $E_{r_1,r_2}=E(u; S[r_1,r_2])\leq
   \delta_{\ref{L:E-estimate}}$ satisfies:
   $$\sup_{z,z' \in S[r_1+1, r_2-1]} \textnormal{dist}\bigl( u(z),
   u(z') \bigr) \leq C_{\ref{L:E-estimate}}\sqrt{E_{r_1,r_2}}.$$
\end{lem}

\begin{proof}[The final step in the proof of Lemma~\ref{L:duw-1}]
   In view of~\eqref{Eq:wn-uR-1},~\eqref{Eq:wn-uR-2}, in order to
   prove Lemma~\ref{L:duw-1} it is enough to prove the following: {\sl
     For every $\epsilon>0$ there exists $a_0=a_0(\epsilon)$ and
     $N_0=N_0(\epsilon)$ such that for every $n \geq N_0$ we have:}
   \begin{equation} \label{Eq:sup-eps}
      \sup_{z \in \widehat{S}[-5R_n+a_0,5R_n-a_0]}
      \textnormal{dist} \bigl(w_n(z), u_{R_n}(z) \bigr) <
      \epsilon.
   \end{equation}

   We prove~\eqref{Eq:sup-eps}. Let $\epsilon>0$.
   By~\eqref{Eq:limE-1},~\eqref{Eq:limE-2} there exists
   $a_0=a_0(\epsilon)$ and $N_0=N_0(\epsilon)$ such that for every
   $n\geq N_0$:
   \begin{align}
      & E\bigl(w_n; \widehat{S}[-5R_n+a_0-1, 5R_n-a_0+1]\bigr) \leq
      \min \Bigl\{ \delta_{\ref{L:E-estimate}}, \,
      \frac{\epsilon^2}{9C_{\ref{L:E-estimate}}^2} \Bigr\},
      \label{Eq:duw-3.1} \\
      & E\bigl(v_{R_n}; \widehat{S}[-5R_n+a_0-1, 5R_n-a_0+1]\bigr)
      \leq \min \Bigl\{ \delta_{\ref{L:E-estimate}}, \,
      \frac{\epsilon^2}{9C_{\ref{L:E-estimate}}^2} \Bigr\}.
      \label{Eq:duw-3.2}
   \end{align}

   We have for every $z \in \widehat{S}[-5R_n+a_0, 5R_n-a_0]$:
   \begin{align} \label{Eq:duw-4}
      \textnormal{dist} \bigl( w_n(z),u_{R_n}(z) \bigr) \leq &
      \textnormal{dist}
      \bigl(w_n(z),w_n(-5R_n+a_0,\tfrac{1}{2})\bigr) \\
      & + \textnormal{dist} \bigl(w_n(-5R_n+a_0,\tfrac{1}{2}),
      u_{R_n}(-5R_n+a_0,\tfrac{1}{2})\bigr)
      \notag \\
      & + \textnormal{dist} \bigl( u_{R_n}(-5R_n+a_0,\tfrac{1}{2}),
      u_{R_n}(z) \bigr). \notag
      \end{align}
      By~\eqref{Eq:duw-3.1} and Lemma~\ref{L:E-estimate}, the first
      term in the right-hand side of~\eqref{Eq:duw-4} satisfies:
      \begin{equation} \label{Eq:duw-5}
         \textnormal{dist} \bigl(w_n(z),w_n(-5R_n+a_0,\tfrac{1}{2})\bigr)
         < \frac{\epsilon}{3}, \quad \forall z \in
         \widehat{S}[-5R_n+a_0, 5R_n-a_0], \, \forall n \geq N_0.
      \end{equation}
      Increasing $N_0$ if necessary\footnote{This increase may depend
        on $a_0=a_0(\epsilon)$ which however has already been fixed.}
      we have from~\eqref{Eq:wn-uR-1} that the second term in the
      right-hand side of~\eqref{Eq:duw-4} satisfies:
      \begin{equation} \label{Eq:duw-6}
         \textnormal{dist}
         \bigl(w_n(-5R_n+a_0,\tfrac{1}{2}),u_{R_n}(-5R_n+a_0,\tfrac{1}{2})
         \bigr)
         < \frac{\epsilon}{3}, \quad \forall  n \geq N_0.
      \end{equation}
      As for the last term in the right-hand side of~\eqref{Eq:duw-4}
      we have:
      \begin{align*}
         & \textnormal{dist} \bigl(
         u_{R_n}(-5R_n+a_0,\tfrac{1}{2}), u_{R_n}(z) \bigr) \\
         & \leq 2 \sup_{z\in \widehat{S}} \textnormal{dist}
         \bigl(v_{R_n}(z), u_{R_n}(z) \bigr) + \textnormal{dist}
         \bigl(v_{R_n}(-5R_n+a_0,\tfrac{1}{2}), v_{R_n}(z) \bigr).
      \end{align*}
      Recall that $v_{R_n}=\exp_{u_{R_n}}(\xi_{R_n})$ with $\Vert
      \xi_{R_n} \Vert_{L^{\infty}} \xrightarrow[n\to \infty]{} 0$.
      Thus increasing $N_0$ once more if necessary, and
      by~\eqref{Eq:duw-3.2} and Lemma~\ref{L:E-estimate} we obtain
      \begin{equation} \label{Eq:duw-7}
         \textnormal{dist} \bigl(
         u_{R_n}(-5R_n+a_0,\tfrac{1}{2}), u_{R_n}(z) \bigr) <
         \frac{\epsilon}{3}, \quad \forall z \in
         \widehat{S}[-5R_n+a_0, 5R_n-a_0], \, \forall n \geq N_0.
      \end{equation}
      Substituting the
      estimates~\eqref{Eq:duw-5},~\eqref{Eq:duw-6},~\eqref{Eq:duw-7}
      in~\eqref{Eq:duw-4} proves~\eqref{Eq:sup-eps}, hence concludes
      the proof of Lemma~\ref{L:duw-1}.
   \end{proof}

\subsection{An auxiliary inequality}  \label{S:auxiliary}
Fix $p>2$. Let $E \to S$ be a vector bundle endowed with a Riemannian
metric with norm $|\cdot|$ and let $\nabla$ be a metric connection on
$E$.
\begin{prop} \label{P:sblv}
   There exists a constant $C_{\ref{P:sblv}}$ such that for every
   $W^{1,p}_{\textnormal{loc}}$-section $\eta:S \to E$ and every $x\in
   S$
   $$|\eta(x)| \leq C_{\ref{P:sblv}} \Bigl(\int_{K_x} |\eta|^p+|\nabla
   \eta|^p \Bigr)^{1/p},$$
   where $K_x$ can be taken to be either
   $S[x,x+1]$ or $S[x-1,x]$.
\end{prop}

\begin{proof}
   The proof for the case of a trivial vector bundle with trivial
   connection and standard metric is essentially contained in the
   proof of Lemma~B.1.16 in~\cite{McD-Sa:Jhol-2}. To adapt the proof
   to the case of a general connection use Remark~3.5.1
   in~\cite{McD-Sa:Jhol-2}.
\end{proof}





%

\section{Proof of the main algebraic statement.}\label{S:proofalg}

This section contains the proofs of the results stated in
\S\ref{sec:algstr}.  Each of the first eight sub-sections below is focused on a
part of the relations to be verified. Each such verification is based
on three steps. First, appropriate moduli spaces are defined. Then,
some regularity properties are established for these spaces along the
lines in \S\ref{S:transversality}. Finally, by making use of the
gluing results in \S\ref{S:gluing}, the desired relations are deduced
out of the description of the boundary of the compactification of the
moduli spaces.  As many types (more than a dozen) moduli spaces will
appear in these verifications it is useful to shortly summarize here
the basic idea in the construction of all of them. Very likely, this
construction will appear quite familiar to the reader in view of the
``bubble tree'' description of bubbling as it appears in
\cite{McD-Sa:Jhol-2} and, even more so, given the cluster moduli
spaces as introduced in \cite{Cor-La:Cluster-1}. It is obviously a
construction that naturally extends that of the pearl moduli spaces
already introduced in \S\ref{S:transversality}.

All the moduli spaces used below consist of configurations modeled on
planar trees with oriented edges with, roughly, the following
properties.  Each of the edges corresponds to a negative gradient
trajectory of a Morse function. This Morse function might be defined
on $L$ or, sometimes, might be a Morse function defined on the ambient
manifold $M$ and more than a single Morse function can be used in the
same tree. Each edge carries a marking indicating which function is
used along that edge.  The orientation of the edge corresponds to that
of the negative gradient flow.  The tree might have several entries
but has a single exit. The entries and the exit are the only vertices
of valence one and each of them corresponds to a critical point of the
function marking the edge ending there. Each internal vertex has a
marking by a certain Maslov class and corresponds to either a
$J$-holomorphic disk with boundary on $L$ or to a $J$ holomorphic
sphere in that class together with a number of marked points (situated
on the boundary or the interior) equal in number to the valence of the
vertex. These marked points correspond to the ends of the edges
reaching the corresponding vertex so that the boundary marked points
are ends of edges labeled by functions defined on $L$ and the interior
marked points are ends of edges labeled by functions defined on $M$.
Depending on the valence and on the particular moduli spaces in
question, specific incidence relations describe at which explicit
marked point on the disk arrives which edge.  Finally, constant $J$
spheres and disks are allowed as long as they are stable in the sense
that, the valence of the corresponding vertex is at least three if it
corresponds to a sphere and, if this valence is two, then the vertex
corresponds to a disk with one marked point in the interior and one on
the boundary.  All the vertices in this paper have valence at most
$4$.

\begin{figure}[htbp]
   \begin{center}
      \includegraphics[scale=.7]{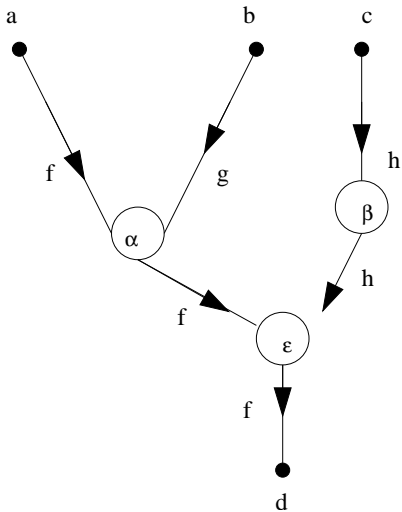}
   \end{center}
   \caption{The functions $f,g,h:L\to \R$,
     $a,d\in\Crit(f), b\in\Crit(g), c\in\Crit(h)$; the three disks are
     of Maslov indexes, respectively, $\alpha$, $\beta$ and $\epsilon$.}
   \label{f:tree}
\end{figure}

This general setup allows for a variety of moduli spaces but the basic
rules are simple: a differential or a morphism is defined by using
moduli spaces with one entry, an operation requires two entries,
associativity requires three, the structures involving the ambient
quantum homology require some internal marked points etc.

The proofs given below are very explicit in a number of cases, in the
sense that the relevant moduli spaces are described in all detail as
well as the related compactness arguments.  This happens for the pearl
moduli spaces in \S\ref{subsec:pearl} as well as for the module
structure in \S\ref{S:qm}. In these cases the ``tree'' language
summarized above is not explicitly used. However, this language does
appear in the description of many of the other moduli spaces. In those
cases, we only generally give the relevant special incidence relations
needed and indicate what variants of the standard arguments used in
\S\ref{subsec:pearl}, \S\ref{S:qm} are necessary in the proof.

\subsection{The pearl complex.}\label{subsec:pearl}

We recall that we have fixed an almost complex structure $J$ which is
compatible with $\omega$ in the sense that $\omega (X,JY)$ is a
Riemannian metric We also fix a Morse function $f:L\to \R$ together
with a metric $\rho$ on $L$ so that the pair $(f,\rho)$ is
Morse-Smale.  The induced negative gradient flow of $f$ is denoted by
$\gamma$. We start by reformulating (in a slightly more general
context) a definition from \S \ref{Sb:tr-pearl-complex}.

\begin{dfn}\label{def:moduli1}Given two points $a,b\in L$ we consider the
   moduli space $\mathcal{M}(a,b,\la;J,f,\rho)$ whose elements are
   objects $(t_{0},u_{0},t_{1},u_{1},t_{2},\ldots u_{k},t_{k+1}), k\in
   \N$ so that:
\begin{itemize}
  \item[i.] For $0\leq i\leq k$, $u_{i}$ is a non-constant
   $J$-holomorphic disk $u_{i}:(D,\partial D)\to (M,L)$.
  \item[ii.]  $\sum_{i}[u_{i}]=\la\in H_2(M,L;\mathbb{Z})$.
  \item[iii.] For $1\leq i\leq k$, $t_{i}\in (0,+\infty) $,
   $t_{0},t_{k+1}\in (0,+\infty]$, so that if we put $a_i =
   u_{i}(-1)$, $b_i = u_{i}(1)$, then $\gamma_{t_{i}}(b_{i-1})=a_{i}$
   for all $i\leq k$, moreover $\gamma_{t_{k+1}}(b_{k})=b$ and
   $\gamma_{-t_{0}}(a_{1})=a$ (of course, if $t_{0}$ is infinite this
   means that $a$ is a critical point of $f$ and similarly for
   $t_{k+1}$).
\end{itemize}
\end{dfn}
We then let $\mathcal{P}(x,y,\la;J,f,\rho)$ the moduli space obtained by
dividing $\mathcal{M}(x,y,\la;J,f,\rho)$ by the action of the obvious
reparametrization group. It is understood in the description above
that $k$ is not fixed. As a matter of convention we also want to
include in our moduli spaces objects for $k=-1$: these are just the
usual negative flow lines joining $x$ to $y$. We will mostly use these
moduli spaces when $x, y\in \Crit(f)$ and in this case this definition
is coherent with that in \S \ref{Sb:tr-pearl-complex}.  We will drop
the decorations $J,f,\rho$ when they are clear from the context.

We now define $$\mathcal{C}^{+}(L; f, \rho, J) = (\Z/2\langle \Crit(f)
\rangle \otimes \La^{+} ,d)$$
with the grading induced by
$|x|=ind_{f}(x)$ for all $x\in \Crit(f)$ and with the differential
defined by
$$dx=\sum_{|y|=|x|+\mu(\la)-1}\#_{\mathbb{Z}_2} \mathcal{P}(x,y,\la)y
t^{\mubar (\lambda)}~.~$$

The next proposition obviously implies the main part of point i.  of
Theorem \ref{thm:alg_main} - we will discuss the augmentation issue at
the end of the section.

\begin{prop}\label{prop:pearl_cplx} For a generic choice of $J$
   the map $d$ is well defined and is a differential. Two generic
   choices $(f,\rho,J)$ and $(f',\rho',J')$ produce chain complexes whose
   homologies are canonically isomorphic.
\end{prop}
The proof of this proposition occupies the rest of the section.  The
main ingredients are the transversality results in \S
\ref{S:transversality}, the gluing results in \S \ref{S:gluing} as
well as the Gromov compactness theorem for holomorphic disks as
described by Frauenfelder in \cite{Fr:msc}.

We will use further the notation of section \S\ref{S:transversality}.
In particular, we have:

\begin{equation}\label{eq:moduli_union}
   \mathcal{P}(x,y,\la; J,f,\rho)=\bigcup_{\substack{k,
       \mathbf{A}=(A_0, \ldots, A_k) \\
       \sum A_{i}=\lambda}} \mathcal{P}(x,y,\mathbf{A};J,f,\rho)~.~
\end{equation}

When $\delta=|x|-|y|+\mu(\la)-1\leq 1$ and for a generic $J$, we see
from Proposition \ref{P:tr-1} that the spaces
$\mathcal{P}(x,y,\mathbf{A};J,f,\rho)$ are manifolds of dimension
$\delta$ and, in case $\delta=0$, these manifolds are already compact.
By Gromov compactness there are only finitely many sequences
$\mathbf{A}=(A_{0},\ldots, A_{k})$ of non zero homology classes so
that $A_{i}$ contains a $J$-holomorphic disk $0\leq i\leq k$ and $\sum
A_{i}=\la$.  This means that the sum in the definition of the
differential $d$ is finite and so this sum is well defined.  For
further use we denote by
$\mathcal{P}(x,y,\emptyset)=\mathcal{P}(x,y,0; J, f,\rho)$ the moduli
space of negative flow lines of $f$ joining $x$ to $y$.  We will also
use the notaion $\mathcal{P}(x,y, (\mathbf{A}', \mathbf{A}''); J, f,
\rho)$ introduced in \S\ref{S:transversality}

\subsubsection{Verification of $d^{2}=0$} \label{subsubsec:verif_diff} Here is the main ingredient
in showing that $d$ is a differential:
\begin{lem}\label{lem:compact}
   Fix $\delta=1$ and $\mathbf{A}=(A_{0},\ldots,A_{k})$. For a generic
   $J$, the natural Gromov compactification
   $\overline{\mathcal{P}}(x,y,\mathbf{A})$ is a topological
   $1$-dimensional manifold whose boundary is described by:
   \begin{equation}\label{eq:bdry1}
      \partial \overline{\mathcal{P}}(x,y,\mathbf{A})=\bigcup_{z,
        \mathbf{A}=(\mathbf{A'},\mathbf{A''})}
      \mathcal{P}(x,z,\mathbf{A}')\times\mathcal{P}(z,y,\mathbf{A}'')\cup
      \mathcal{T}\cup \mathcal{T}'~.~
   \end{equation}
   Where the terms $\mathcal{T}$ and $\mathcal{T'}$ are given by:
   $$\mathcal{T}=\bigcup_{\mathbf{A}=(\mathbf{A'},\mathbf{A''}),\
     \mathbf{A}',\mathbf{A}''\not=\emptyset}
   \mathcal{P}(x,y,\mathbf{A'},\mathbf{A''})$$
   and
   $$\mathcal{T}'=\bigcup_{i\ ;\ A_{i}=A_{i}'+A_{i}'',\
     A_{i}',A_{i}''\not=0 } \mathcal{P}(x,y,(A_{0},\ldots
   A_{i-1},A_{i}'),(A_{i}'',A_{i+1},\ldots A_{k}))~.~$$
\end{lem}
We emphasize that in the first term of formula~\eqref{eq:bdry1}, in
the decomposition $\mathbf{A}=(\mathbf{A}',\mathbf{A''})$ it is
possible that $\mathbf{A}'=\emptyset$ or $\mathbf{A}''=\emptyset$.

We first remark that the lemma immediately implies $d^{2}=0$.  Indeed,
notice that, by the lemma, each term:
$$\mathcal{P}(x,y,\mathbf{A'},\mathbf{A''})$$
with
$\mathbf{A}'=(A_{0},\ldots, A_{i})$, $\mathbf{A''}= (A_{i+1},\ldots,
A_{k})$ appears as a boundary in precisely two moduli spaces:
$\overline{\mathcal{P}}(x,y,(\mathbf{A'},\mathbf{A}''))$ and
$\overline{\mathcal{P}}(x,y, (A_{0}\ldots,A_{i-1}, A_{i}+A_{i+1},
A_{i+2},\ldots,A_{k}))$.

This means that if we define:
$$\overline{\mathcal{P}}(x,y,\la)=\coprod_{\sum
  A_{i}=\la}\overline{\mathcal{P}}(x,y,\mathbf{A})~,~$$
then
$$0=\#_{\mathbb{Z}_2}(\partial\overline{\mathcal{P}}(x,y,\la)) =
\#_{\mathbb{Z}_2} (\bigcup_{z,\la'+\la''=\la}\mathcal{P}(x,z,\la')
\times \mathcal{P}(z,y,\la''))$$
which clearly means $d^{2}=0$.

{\em Proof of Lemma \ref{lem:compact}.} We start by making precise the
compactification $\overline{\mathcal{P}}(x,y,\mathbf{A})$. We denote
by $\mathcal{P}_L$ the space of continuous paths $\{\gamma: [0,b]\to
L: b\geq 0\}$. We consider the subspace $\mathcal{P}_{f,\rho}\subset
\mathcal{P}_{L}$ of those paths $\gamma:[0,b]\to L$ which are flow
lines of $-\nabla_{g} f$ reparametrized so that $f(\gamma (t))=
f(\gamma(0))-t$. We allow that one or both of the ends of such of a
flow line be a critical point of $f$. Exponential convergence close to
critical points insures that this choice of parametrization is
continuous when defined on the usual moduli space of negative flow
lines (parametrized by time). The space $\mathcal{P}_{f,\rho}$ has a
natural compactification $\overline{\mathcal{P}}_{f,\rho}$ formed by
adding all broken flow lines to $\mathcal{P}_{f,\rho}$.

We now notice that $\mathcal{P}(x,y,\mathbf{A})$ is included in:
$$\mathcal{L}=\mathcal{P}_{f,\rho}\times \mathcal{M}(A_0,J)/G_{-1,1}
\times \mathcal{P}_{f,\rho} \times \mathcal{M}(A_1,J)/G_{-1,1} \times
\cdots \mathcal{P}_{f,\rho} \times \mathcal{M}(A_l,J)/G_{-1,1} \times
\mathcal{P}_{f,\rho}~.~$$
There is a natural compactification
$\overline{\mathcal{L}}$ of $\mathcal{L}$ obtained by compactifying
each term $\mathcal{M}(A_i,J)$ in the sense of Gromov and by replacing
$\mathcal{P}_{f,\rho}$ by $\overline{\mathcal{P}}_{f,\rho}$. We take on
$\mathcal{P}(x,y,\mathbf{A})$ the induced topology and define its
compactification to be the closure of this space inside
$\overline{\mathcal{L}}$.

\emph{I. Structure of the compactification.}  The first part of the
proof is to show that, for $\delta=1$ and a generic $J$, if $u\in
\overline{\mathcal{P}}(x,y,\mathbf{A})\backslash
\mathcal{P}(x,y,\mathbf{A})$ then $u$ has the one of the forms given
in the right hand side of~\eqref{eq:bdry1}.  It is clear that
$u=(\gamma_{0},u_{0},\gamma_{1},\ldots u_{k},\gamma_{k+1})$ with
$u_{i}\in \overline{\mathcal{M}(A_{i},J)/G_{-1,1}}$,
$\gamma_{i}\in\overline{\mathcal{P}}_{f,\rho}$.  To proceed with the
proof it is useful to reinterpret in more combinatorial terms
Propositions \ref{P:tr-1} and \ref{P:tr-2}.  Let
$\mathcal{P}(x,y,\mathbf{A},s)\subset
\overline{\mathcal{P}}(x,y,\mathbf{A})$ be the moduli space whose
definition is identical with that of $\mathcal{P}(x,y,\mathbf{A})$
except that $s$ of the flow lines $\gamma_{i}$ are of length $0$.
Proposition \ref{P:tr-1} shows that, generically, if $\delta\leq 1$,
then $\dim (\mathcal{P}(x,y,\mathbf{A},0))=\delta$ and Proposition
\ref{P:tr-2} shows that $\dim
(\mathcal{P}(x,y,\mathbf{A},1))=\delta-1$ (both spaces being void if
the respective dimension is negative). The same argument also shows
that $\mathcal{P}(x,y,\mathbf{A},s)=\emptyset$ whenever $s\geq 2$.

Returning to our $u\in\overline{\mathcal{P}}(x,y,\mathbf{A})$, denote
by $\alpha_{i}$ the end of $\gamma_{i}$ and by $\beta_{i}$ the
beginning of $\gamma_{i+1}$. Clearly, $\alpha_{i},\beta_{i}\in
Image(u_{i})$. The structure of $u_{i}$ is described by a bubbling
tree (see \cite{Fr:msc} and \cite{McD-Sa:Jhol-2}) which carries two
marked points which are mapped to $\alpha_{i}$ and $\beta_{i}$. It is
important to notice that these two marked points are distinct but they
might both lie on a ``ghost" (or constant) disk.  We will denote these
two marked points by $-1$ and $+1$ respectively (by a slight abuse in
notation). In this tree let $u^{-1}_{i}$ and $u^{+1}_{i}$ be the two
curves carrying $-1$ and $+1$ respectively, and consider the path of
minimal length which joins the vertices $u^{-1}_{i}$ and $u^{+1}_{i}$.
In case $\alpha_{i}\not=\beta_{i}$ we denote by $v_{i}$ the nodal
curve formed by the union of all the vertices along this path. If
$\alpha_{i}=\beta_{i}$ we let $v_{i}$ be the constant map from the
disk equal to $\alpha_{i}$.  We notice that $v_{i}$ may be assumed not
to contain any sphere component.  Moreover $\omega(v_{i})\leq \omega
(u_{i})$ with equality if and only if the tree is linear with
$u^{-1}_{i}$ and $u^{+1}_{i}$ at its ends.  We notice that
$v=(\gamma_{0},v_{0},\ldots v_{k},\gamma_{k+1})
\in\overline{\mathcal{P}}(x,y,\mathbf{B})$ where
$\mathbf{B}=([v_{0}],\ldots, [v_{k}])$.  Let $l_{i}$ be the number of
components in $v_{i}$ less one. Let $b_{i}$ be the number of breaks in
$\gamma_{i}$ and let $s$ be the number of flow lines $\gamma_{j}$
which are of length $0$. Our statement follows if we show, $\forall i$
$\omega(v_{i})=\omega(u_{i})$ and $\sum l_{i}+\sum b_{i} +s=1$.

The key observation is that, because the components of the $v_{i}$'s
are disks, we may rewrite in a unique way: $$v\in
\mathcal{P}(x,z_{1},\mathbf{B}_{1},s_{1})\times\ldots
\times\mathcal{P}(z_{r},y,\mathbf{B}_{r+1},s_{r+1})$$
and, in this
case, $\sum s_{r}=\sum l_{i}+ s$; $\sum b_{i}=r$.

We also have $$2=\delta+1\geq
(|x|-|z_{1}|+\mu(B_{1}))+(|z_{1}|-|z_{2}|+\mu(B_{2}))+\ldots +
(|z_{r}|-|y|+\mu (B_{r+1}))$$
with equality if and only if $[v]=[u]$.
Now notice that it is not possible that $r>1$ because in that case at
least one of the parenthesis in the sum above would need to be $\leq
0$ which implies that the respective space is void. We are thus left
with two cases. The first is $r=1$ and $(|x|-|z_{1}|+\mu(B_{1}))=1$,
$(|z_{1}|-|y|+\mu(B_{2}))=1$.  But this implies that $s_{1}=0$ and
$s_{2}=0$ and $[v]=[u]$.  The last case is $r=0$. Then, as
$N_L=2$, we obtain $[v]=[u]$ and, as $u\not\in
\mathcal{P}(x,y,\mathbf{A})$ we get $s_{1}=1$.

\emph{II. Gluing.} The fact that $\partial
\overline{\mathcal{P}}(x,y,\mathbf{A})$ is contained in the union
appearing on the right side of~\eqref{eq:bdry1} follows from Gromov
compactness combined with standard arguments from Morse theory.

Next, we need to show that each element appearing on the right side
of~\eqref{eq:bdry1} is a boundary point of
$\overline{\mathcal{P}}(x,y,\mathbf{A})$ and moreover, that each such
element corresponds to a {\em unique} end of the $1$-dimensional manifold
$\mathcal{P}(x,y,\mathbf{A})$ (hence
$\overline{\mathcal{P}}(x,y,\mathbf{A})$ is indeed a compact manifold
with boundary).

By the results in \S\ref{S:transversality} we have that for a generic
choice of $J$ we have
$\mathcal{P}(x,y,\mathbf{A})=\mathcal{P}^{\ast,d}(x,y,\mathbf{A})$ for
all $x$, $y$, $\mathbf{A}$ so that $|x|-|y|+\mu(\mathbf{A})-1\leq 1$.
This implies by standard Morse theory that the elements in the first
term of the decomposition~\eqref{eq:bdry1} can be glued in the obvious
way. The uniqueness statement is also standard in this case.

The elements in $\mathcal{T}'$ are also seen as boundary points by
first applying the results in \S\ref{S:transversality} to show that
the transversality necessary for the gluing described in
\S\ref{S:gluing} is achieved generically and then gluing as there.
The uniqueness statement follows too from the results in
\S\ref{S:gluing}. (See Corollary~\ref{C:glue-pearls-1}.)

Finally, we are left to deal with the elements in $\mathcal{T}$.
Consider one such element
$u=(\gamma_{0},u_{0},\gamma_{1},\ldots,u_{i},u_{i+1},\ldots
u_{k},\gamma_{k+1}) \in \mathcal{P}(x,y,\mathbf{A}',\mathbf{A}'')$
with $\mathbf{A}'=(A_{0},\ldots A_{i})$ and
$\mathbf{A}''=(A_{i+1},\ldots, A_{k})$.  The transversality proved in
\S\ref{S:transversality} implies immediately the result in this case
if $p(u)=u_{i}(+1)=u_{i+1}(-1)\not\in\Crit(f)$. But, again for a
generic choice of $J$, it follows by inspecting the definition of
$\mathcal{P}(x,y,\mathbf{A}',\mathbf{A}'')$ that $p(u)$ avoids the
finite set $\Crit(f)$ whenever $\delta\leq 1$.  \qed

\subsubsection{Invariance}\label{subsubsec:inv}
The purpose of this subsection is to show that for two generic sets of
data $(J,f,\rho)$ and $(J',f',\rho')$ we can define a chain morphism
$$\phi: \mathcal{C}^{+}(L;f,\rho,J)\to \mathcal{C}^{+}(L;
f',\rho',J')$$
which induces an isomorphism in homology. This morphism
$\phi=\phi^{\overline{J},F,G}$ depends on auxiliary data: a generic
homotopy of almost complex structures $\overline{J}$ joining $J$ to
$J'$ and a Morse cobordism $(F,G)$ relating $(f,\rho)$ to $(f',\rho')$
(see \cite{Co-Ra:Morse-Novikov} for the formal definition of Morse
cobordisms).  Of course, all this construction is typical for
invariance proofs in Morse or Floer theories and, after defining the
relevant moduli spaces, we will only sketch the rest of the proof as
there are no new transversality issues compared with the
last section and the gluing arguments are analogous to the ones before.

We now fix a smooth family of almost complex structures
$\overline{J}=J_{t},\ t\in [0,1]$ so that $J_{0}=J$ and $J_{1}=J'$. We
also fix a Riemannian metric $G$ on $L\times [0,1]$ so that
$G|_{L\times\{0\}}=\rho$ and $G|_{L\times \{1\}}=\rho'$ as well as a Morse
function $H:L\times [0,1]\to \R$ so that $H|_{L\times\{0\}}(x)=f(x)+h$
(for some constant $h$) $H|_{L\times\{1\}}(x)=f'(x)$, $(H,G)$ is
Morse-Smale with
$\Crit_{i}(H)=\Crit_{i-1}(f)\times\{0\}\cup\Crit_{i}(f')\times\{1\}$
and, finally, $(\partial H/\partial t)(x,t)=0$ for $t=0,1$, $(\partial
H/\partial t)(x,t)< 0$ for $t\in (0,1)$.  It is easy to see how to
construct such $(H,\overline{J})$ (we refer to \cite{Co-Ra:Morse-Novikov} for an
explicit such construction).

Given two critical points $x\in \Crit(f)$ and $y'\in \Crit(f')$ we
define the moduli space $\mathcal{M}(x,y',\la; \overline{J},H,G)$ by
using a slightly modified version of Definition \ref{def:moduli1}. The
only changes are listed below:
 \begin{itemize}
   \item[i'.] The $u_{i}$'s are so that for each $u_{i}$ there is some
    $\tau_{i}\in [0,1]$ with $u_{i}:(D,\partial D)\to (M,L)\times
    \{\tau_{i}\}\subset (M,L)\times [0,1]$ and $u_{i}$ is
    $J_{\tau_{i}}$-holomorphic.
   \item[iii'.] The incidence relations at point iii. take place
    inside $L\times [0,1]$ with the flow of $-\nabla_{G}H$ on $L\times
    [0,1]$ taking the place of the flow $-\nabla_{\rho}(f)$.
\end{itemize}

We denote by $\mathcal{P}(x,y',\la; \overline{J},H,G)$ the resulting
moduli space after division by the reparametrization group. By
inspecting the construction and proofs in \S\ref{S:transversality} it
is easy to see that all the arguments carry over when replacing the
moduli spaces $\mathcal{M}(A,J)$ with moduli spaces
$\mathcal{M}(A,\overline{J})$ which are made of disks that are
$J_{\tau}$-holomorphic, $\tau\in [0,1]$ (and in the class $A\in \La$)
and, for each such disk, the relevant evaluation maps take values in
$L\times\{\tau\}$.  The key reason for this is that Lemmas \ref{L:uv}
and \ref{L:non-simple} still apply when $J$ is replaced with a generic one
parametric family of almost complex structures $\overline{J}$ and this
implies that the analogue of Proposition \ref{P:tr-1} remains true in
this setting (the same remains true in fact even if for higher parametric families).

The conclusion is that for a generic choice of $\overline{J}$, $H$ and
$G$, if $\tilde{\delta}=|x|-|y'|+\mu(\la)$ then, for
$\tilde{\delta}\leq 1$, the moduli space $\mathcal{P}(x,y',\la)$ is a
manifold of dimension $\tilde{\delta}$, is compact if $\tilde{\delta}=0$
and is void if $\tilde{\delta}\leq -1$. The compactification
$\overline{\mathcal{P}}(x,y',\la)$ of $\mathcal{P}(x,y',\la)$ can then
be defined as in the last section and, by the same method as there, we
see that when $\tilde{\delta}=1$, we have:
\begin{eqnarray}\label{eq:bdry2}
   0=\#_{\mathbb{Z}_2}
   (\partial\overline{\mathcal{P}}(x,y',\la)) =
   \#_{\mathbb{Z}_2}(\bigcup_{z\in\Crit(f),
     \la'+\la''=\la}
   \mathcal{P}(x,z,\la')\times \mathcal{P}(z,y',\la'')) + \nonumber\\ +\
   \#_{\mathbb{Z}_2}(\bigcup_{z'\in\Crit(f'),\la'+\la''=\la}
   \mathcal{P}(x,z',\la')\times\mathcal{P}(z',y',\la''))
\end{eqnarray}
 where in both unions we take into account the cases when $\la'$ or
 $\la''$ is null (in that case the relevant moduli spaces are just the
 usual moduli spaces of Morse trajectories).
 Of course, this identity also depends on a gluing argument in which $J$
 is no longer constant. However, this is a reasonably straightforward adaptation
 of the argument in \S\ref{S:gluing} and so, for brevity, we will not
 make it explicit here.

 For $x\in \Crit(f)$ we now put:
 $$\phi^{\overline{J},F,G}(x)=\sum_{y'\in\Crit(f'),\ \la\ ;\
   |x|-|y'|+\mu(\la)=0}\#_{\mathbb{Z}_2}(\mathcal{P}(x,y',\la;J,H,G))y
 t^{\mubar(\lambda)}~.~$$

 For generic choices of the data, formula~\eqref{eq:bdry2} implies
 that $\phi^{\overline{J},H,G}$ is a chain morphism.

 It is very easy to see that this chain morphism induces an
 isomorphism in homology. Indeed, both the pearl complex differential
 and the morphism $\phi$ clearly respect the degree filtration.
 Therefore, $\phi$ induces a morphism between the spectral sequences
 $E(\phi):E(J,f,\rho)\to(E'(J',f',\rho')$ associated to the two filtered
 complexes $\mathcal{C}^{+}(L;J,f,\rho)$ and
 $\mathcal{C}^{+}(L;J',f',\rho')$.  It is clear that the $E^{0}$ terms in
 these spectral sequences only involve moduli spaces of Morse
 trajectories so that $E^{1}(\phi)$ is already an isomorphism. But
 this means that $H_{\ast}(\phi)$ is an isomorphism.

The argument needed to prove that this isomorphism is canonical - while essentially standard -
 is more complicated.    First, we need to notice that if $f=f'$,
 $\overline{J}$ is the constant family $\overline{J}_{\tau}=J$, $G$ is also the constant metric $G_{\tau}=\rho$
  and $H$ is a trivial homotopy in the sense that $H(x,t)=f(x)+h(t)$ for some appropriate function $h:[0,1]\to \R$), then $\phi^{\overline{J}, H, G }=id$.  This is easily seen because on one hand these constant
  choices of homotopies are regular and, on the other hand, the only moduli spaces appearing in
  the definition of $\phi^{\overline{J}.H,G}$ which are $0$-dimensional are the flow lines in
  $L\times [0,1]$ which project to $L$ on constant paths (equal to critical points of $f$).

  We now fix a second set of homotopies $(H',G',\overline{J}')$
  relating $(f,\rho,J)$ and $(f',\rho',J')$ in the same way as above and we intend to show that the
  two morphisms $\phi^{\overline{J}, H, G }$ and $\phi^{\overline{J}', H', G' }$ are chain homotopic.
  The construction of this chain homotopy is perfectly similar to the construction of the chain
  morphism $\phi^{\overline{J}, H, G }$.  We start with a Morse homotopy $F:L\times [0,1]\times [0,1]\to \R$ relating $H$ to $H'$ so that $F(x,t,0)=H(x,t)+k$, $F(x,t,1)=H'(x,t)+k'$ (with $k$ and $k'$ appropriate
  constants) , $F(x,0,\tau)$ and $F(x,1,\tau)$ are trivial homotopies (of $f$ and of $f'$, respectively).
  The critical points of $F$ verify $\Crit_{k}(F)=\Crit_{k-2}(f)\times\{0,0\}\cup \Crit_{k-1}(f)\times\{1,0\}\cup\Crit_{k-1}(g)\times\{0,1\}\cup \Crit_{k}(g)\times\{1,1\}$ (see \cite{Co-Ra:Morse-Novikov}
  for the construction of such an $F$). We also consider homotopies
  $\tilde{J}_{t,\tau}$ between $\overline{J}$ and $\overline{J}'$ as well as $\tilde{G}$ between $G$ and $G'$. We now define pearl type moduli spaces in the same way as just above except that
  points i' and iii' are replaced by:
   \begin{itemize}
   \item[i''.] The $u_{i}$'s are so that for each $u_{i}$ there is some
    $(t_{i},\tau_{i})\in [0,1]\times [0,1]$ with $u_{i}:(D,\partial D)\to (M,L)\times
    \{(t_{i},\tau_{i})\}\subset (M,L)\times [0,1]\times [0,1]$ and $u_{i}$ is
    $\tilde{J}_{t_{i},\tau_{i}}$-holomorphic.
   \item[iii''.] The incidence relations at point iii. take place
    inside $L\times [0,1]\times [0,1]$ with the flow of $-\nabla_{\tilde{G}}F$ on $L\times
    [0,1]\times [0,1]$ taking the place of the flow $-\nabla_{\rho}(f)$.
\end{itemize}
 It is easy to see that a generic choice of $\tilde{G}$ and $\tilde{J}$ suffice for the transversality
 required to give the resulting moduli spaces $\mathcal{P}(x,y';\tilde{J},F,\tilde{G})$ a structure of
manifold of dimension $|k|-|q|+\mu(\la)-1$  where $x\in \Crit_{k}(F)$, $y'\in \Crit_{q}(F)$ and $\la$
is the total homotopy class in $\pi_{2}(M,L)$ of the configuration.
Let $\xi$ be the $\La^{+}$-module morphism obtained by counting
the elements in the $0$-dimensional such moduli spaces when $x\in \Crit_{k}(f)=\Crit_{k+2}(F)|_{L\times \{(0,0)\}}$ and $y\in\Crit_{k+1}(f')=\Crit_{k+1}(F)|_{L\times \{(1,1)\}}$. We want to remark that $\xi$ is precisely the wanted chain homotopy.

 \begin{equation}\label{equ:cube1} \begin{minipage}{5in}
\xy\xymatrix@C+15pt{
f\ar@{-}[rr]^{H}\ar@{-}[dd]\ar@{-}[dr]_{id} & & f'\ar@{..}[dd]\ar@{-}[dr]^{id} &
\\ & f \ar@{-}[rr]^{H'\ \ \ \ \ }\ar@{-}[dd] & & f'\ar@{-}[dd]\\
(0,0)\ar@{..}[rr]\ar@{-}[dr] &&(1,0)\ar@{..}[dr]& \\
& (0,1)\ar@{-}[rr] & & (1,1)}
\endxy

The verticals in this picture correspond to $L\times  $  the coordinate in the horizontal
square; $F$ is defined on the whole cube; trivial Morse homotopies are defined on the left
and right vertical faces; the Morse homotopy $H$ is defined on the back face and $H'$ is defined
on the front face.

\end{minipage}
\end{equation}

For this we need to analyze the boundary of the compactification of the $1$-dimensional moduli spaces
$\mathcal{P}(x,y';\tilde{J},F,\tilde{G})$  (clearly, we use here both Gromov compactness and  an appropriate variant of the gluing results in \S\ref{S:gluing})
and we see that the only boundary elements  which count are  of four types.  Indeed,
in one-dimensional such moduli spaces side bubbling is not possible and, by applying the exact
 same method as that used in \S\ref{subsubsec:verif_diff} to show $d^{2}=0$,  we see that, for generic choices
 of data, the only remaining terms correspond to a pearl configuration breaking at some critical point
 of $F$.   Thus there are four possibilities: a break in $z\in \Crit(f)\times \{(0,0)\}$ - this
 corresponds to an element in  $\mathcal{P}(x,z;J,f,\rho)$ followed by one counted in $\xi$;
 a break in $z'\in \Crit(f')\times\{(1,0)\}$- this corresponds to an element in $\mathcal{P}(x,z'; \overline{J},H, G)$ followed by one associated to the trivial homotopy $F|_{L\times\{1\}\times [0,1]}$;
 a break in $z''\in Crit(f)\times\{(0,1)\}$ - this corresponds to aan elelement associated to the
 trivial homotopy $F|_{L\times \{0\}\times [0,1]}$ followed by an element in $\mathcal{P}(z'',y';\overline{J}',H',G')$; finally, an element $z'''\in \Crit(f')\times\{(1,1)\}$ which corresponds to an element
 counted in $\xi$ followed by an element of $\mathcal{P}(z''',y'; J',f', \rho')$. Puting all of this together
 we obtain:
 $$d\xi+\xi d= I_{1}\circ \phi^{\overline{J},H,G}+\phi^{\overline{J}',H',G'}\circ I_{2}$$
 where both $I_{1}, I_{2}$ are induced by trivial homotopies. Given that trivial homotopies induce
 the identity at the level of the pearl complexes we obtain the claim.

\begin{rem}\label{rem:diff_m}
   a. In the proof of the invariance
   we may avoid entirely the use of spectral sequences by proving - again by constructing
   an appropriate chain homotopy - that given two homotopies $H$ from $f$ to $f'$ and $H'$ from
   $f'$ to $f''$ and assuming $H''=H\# H'$ is the concatenated homotopy from $f$ to $f''$, then
   $\phi^{H''}$ is chain homotopic to $\phi^{H'}\circ \phi^{H''}$.

   b. A useful feature of the pearl complex is the following. Let $f$
   be a Morse function on $L$ with a single maximum, $x_n$. Then any
   non-void moduli space $\mathcal{M}(x_n,x;\la)$,
   $|x_n|-|x|+\mu(\la)\leq 2$, is of dimension at least $1$ when
   $\la\not=0$. This means that in each of the complexes
   $\mathcal{C}^{+}(L;J,f)$, $\mathcal{C}(L;J,f)$ we have $d x_n=0$
   (similarly, it is easy to see that the minimum of $f$, if it is
   unique, can not be a boundary in this complex). This has the
   following interpretation: if $L$ as before has the property that
   there exists a point $p \in L$ and some $J$ so that no $J$-disk
   passes through $p$, then $QH_{\ast}(L)\not=0$.  Indeed, in this
   case, we may assume that $p$ is the maximum $x_n$ of some function
   $f$ and that $\mathcal{C}(L;f,J)$ is defined so that we know $d
   x_n=0$. If no $J$-disk goes through $x_n$, then $x_n$ can not be a
   boundary. Using point v. of Theorem \ref{thm:alg_main} this means
   in particular that if $HF_{\ast}(L)=0$ (for example if $L$ is
   displaceable by a Hamiltonian diffeomorphism), then through each
   point in $L$ passes a $J$-holomorphic disk.

   c. It is easy to see that $Q^{+}H_{\ast}(L)$ can never vanish.
   Indeed, assume that $f$ is a function with a single maximum $m$. In
   that case, as before, $dm=0$ in the (positive) pearl complex but
   $m$ can never be a boundary in $\mathcal{C}^{+}(L;f)$ (it obviously
   can be in $\mathcal{C}(L;f)$).
\end{rem}

\subsubsection{Augmentation.}

We now want to remark that there exists a chain morphism:
$$\epsilon_{L,f}:\mathcal{C}^{+}(L;f)\to \La^{+}$$
where the differential in the target is trivial.  Moreover, this
morphism will be easily seen to commute with the comparison maps
constructed above so that the induced map in homology is canonical and
will be denoted by $\epsilon_{L}$. The definition of $\epsilon_{L,f}$
is simple: it is a $\La^{+}$-module map so that $\epsilon_{L,f}(x)=0$
for all $x\in\Crit(f),\ |x|>0$ and $\epsilon_{L,f}(z)=1$ for all $z\in
\Crit(f), \ |z|=0$. The reason this definition produces a chain
morphism is similar to the point b. in Remark \ref{rem:diff_m}. Indeed
it is easy to see that for a minimum $x_0 \in \Crit_{0}(f)$ if a
moduli space $\mathcal{P}(x,x_0;\la)$ is non-void and
$x\not=x_0$,$\la\not=0$, then the dimension of this moduli space is
at least $1$. Given the form of our differential it follows that a
minimum can only appear in the Morse part of the differential of a
critical point. But, in that part minima always appear in pairs
(because, for a critical point of index $1$ the Morse differential is
always an even sum of minima). This implies that $\epsilon_{L}$ is a
chain map as claimed. To prove the invariance of $\epsilon_{L,f}$ it
is enough to notice that if $\phi:\mathcal{C}^{+}(L;f)\to
\mathcal{C}^{+}(L;f')$ is the morphism constructed in the last
section, then we have $\epsilon_{L,f'}\circ\phi=\epsilon_{L,f}$. This
happens because in a way similar as above we see that a ``comparison"
moduli space $\mathcal{P}(x,x_0';\la)$ with
$x_0'\in\Crit_{0}(f')$ $x\in \Crit(f)$ can be $0$-dimensional and
non-void only if $x\in\Crit_{0}(f)$ and $\la=0$.

\subsection{The quantum product.}\label{subsec:product}
We construct here an operation
$$\ast :
\mathcal{C}^{+}_{k}(L;f,J)\otimes\mathcal{C}^{+}_{l}(L;f',J)\to
\mathcal{C}^{+}_{k+l-n}(L;f'',J)$$
where $f'$ and $f''$ are generic
small deformations of the Morse function $f$. To simplify notation we
will assume that the critical points of both $f'$ and $f''$ coincide
with those of $f$. This morphism of chain complexes is defined by:
\begin{equation}\label{eq:prod}
   x\ast y=\sum_{y,\la}\#_{\mathbb{Z}_2}(\mathcal{P}(x,y,z;\la,
   J))z t^{\mubar(\lambda)},
\end{equation}
where the moduli space $\mathcal{P}(x,y,z;\la, J)$ is described as:
$$\begin{array}{l}
   \mathcal{P}(x,y,z;\la,J)=\{(l_{1},l_{2},l_{3}, u)\ :\\
   (l_{1},l_{2},l_{3})\in \mathcal{P}(x,a_{1},\la_{1}; f,J)\times
   \mathcal{P}(y,a_{3},\la_{3}; f',J)\times
   \mathcal{P}(a_{2},z,\la_{2};f'',J) \ , \\
   u:(D, \partial D)\to (M,L) \ , \bar{\partial}_{J} u=0, \
   u(e^{2ik\pi/3})=a_{k},\ k\in\{1,2,3\}\ , \\
   \la_{1}+\la_{2}+\la_{3}+[u]=\la\ \}~.~
\end{array}$$

Here, of course, the spaces $\mathcal{P}(a,b,\la;f,J)$ are pearl
moduli spaces which join the points $a$ and $b$ in $L$. The sum above
is taken only for those elements so that $
|x|+|y|+\mu(\la)-|z|-n=0~.~$ Notice that we do allow in this
description that $u$ be the constant map.

\begin{prop}\label{prop:quantum_prod}
   For generic choices of data, the operation defined in
   equation~\eqref{eq:prod} is well defined and a chain map. It
   induces in homology an associative product $$Q^{+}H_{k}(L)\otimes
   QH^{+}_{j}(L)\to Q^{+}H_{k+j-n}(L)$$ which is independent of the
   choices made in the construction.
\end{prop}

\begin{proof}
   We start the proof by describing the moduli spaces involved in a
   more precise way. We will use the notation of
   \S\ref{S:transversality}. Fix three sequences of nonvanishing
   homology classes: $\mathbf{A}=(A_{1},A_{2},\ldots A_{k})$,
   $\mathbf{A}'=(A'_{1},A'_{2},\ldots A'_{k'})$,
   $\mathbf{A}''=(A''_{1},A''_{2},\ldots A''_{k'})$ together with
   another homology class $U$ which can also be null.  We now define:
   \begin{equation}\label{eq:m-pr1}
      \begin{array}{l}
         \mathcal{P}(x,y,z,\mathbf{A},\mathbf{A'},\mathbf{A''},U;
         f,f',f'',J)=
         \{(u,u',u'',v)\in\mathcal{P}(x,a_{1},\mathbf{A})\times\\
         \times
         \mathcal{P}(y,a_{3},\mathbf{A}')
         \times\mathcal{P}(a_{2},z,\mathbf{A}'')
         \times\mathcal{M}(U, J) \ : \
         a_{k}=v(e^{2\pi ik/3}), \ k=1,2,3\}
      \end{array}
   \end{equation}
   It is clear that the moduli space $\mathcal{P}(x,y,z;\la,J)$ is the
   union of the one given above when $\sum A_{i}+\sum A'_{i}+\sum
   A''_{i}+U=\la$.  As the notation for all these moduli spaces will
   rapidly become hard to manipulate we will sometimes denote by
   $\mathcal{P}_{\eqref{eq:m-pr1}}$ a moduli space as defined in
   equation~\eqref{eq:m-pr1}. We will denote by
   $\mathcal{P}^{\ast,d}_{\eqref{eq:m-pr1}}$ the moduli space defined
   as in~\eqref{eq:m-pr1} but so that all the $J$-disks involved are
   simple and absolutely distinct. For an element $(u,u',u'',v)\in
   \mathcal{P}^{\ast,d}_{\eqref{eq:m-pr1}}$ we will call the disk $v$
   the core of the configuration.

   With these notation the first part of the proof is to show that:
   \begin{lem}\label{lem:transv-prod}
      For a generic set of almost complex structures $J$ we have
      $$\mathcal{P}_{\eqref{eq:m-pr1}} =
      \mathcal{P}^{\ast,d}_{\eqref{eq:m-pr1}}$$
      whenever
      $\delta'=|x|+|y|-|z|+\mu(\mathbf{A}+\mathbf{A}'+\mathbf{A}''+U)-n
      \leq 1$.
   \end{lem}
   If this is true it follows that $\mathcal{P}_{\eqref{eq:m-pr1}}$ is
   a manifold of dimension $\delta'$ when $\delta'=0,1$ and is void
   when $\delta'<0$.  Moreover, by Gromov compactness it also follows
   that for $\delta'=0$ this space is compact so that the sum
   in~\eqref{eq:prod} is well defined.

   \emph{Proof of Lemma \ref{lem:transv-prod}}. This is a
   straightforward adaptation of Proposition \ref{P:tr-1} and we will
   only indicate here the specific additional verifications that are
   needed in this case.  Clearly, an analogue of Proposition
   \ref{P:tr-2} is also needed - obviously the relevant moduli space
   in this case is defined in the same way as in
   $\mathcal{P}_{\eqref{eq:m-pr1}}$ except that precisely one of the
   edges (or flow lines of one of $f$, $f'$ or $f''$) is required to
   have length $0$. The dimension of these moduli spaces is controlled
   by $\delta'-1$. To show them, Lemma \ref{L:non-simple} can be
   expanded in the sense that, in the statement, the roots of order
   $2$ of unity, $+/-1$, may be replaced with the roots of order $3$
   of $1$ - the reason for that is, of course, that any three points
   on the boundary of a pseudo-holomorphic disk may be carried by
   reparametrization to the roots of order $3$. Compared with the
   proof of Lemma~\ref{L:non-simple} the only additional remark is
   that the cyclic order of the three points $z_{1}$, $z_{2}$, $z_{3}$
   on the disk $u'$ so that $u'(z_{k})=u(e^{2\pi i k})$ is the same as
   the cyclic order of the corresponding three roots of the unity on
   $u$. It is also needed to show that, generically, we may assume
   that the cores of the elements in our moduli spaces send the three
   roots of the unity to distinct points - this is again a simple
   exercise.  This argument covers the reduction to simple disks in
   the proof of Proposition~\ref{P:tr-1}. To pursue with the reduction
   to absolutely distinct disks it turns out that it is convenient to
   work with moduli spaces more general that those defined in
   equation~\eqref{eq:m-pr1}: they are again formed by $(u,u',u'',v)$
   except that the pearls $u$, $u'$, $u''$ are more general in the
   sense that along each of their edges we may use any one of the
   negative gradient flows $-\nabla f, -\nabla f',-\nabla f''$. More
   precisely, for the definition of the $(u,u',u'')$ in the incidence
   relation iii. in Definition \ref{def:moduli1} anyone of these three
   flows may be allowed with the assumption that, at the core, the two
   ``entry'' flow lines correspond to different flows.
   To see why these more general moduli spaces are needed recall the principle
   behind the reduction to absolutely distinct disks  in Case 2 in \S \ref{Sb:prf-tr-1}
   (for $n\geq 3$).  We assume that some disks
   are not absolutely distinct for some element of our moduli space and we show
   that, in this case, there exists a configuration of  Maslov index lower by at
   least two and with the same ends. To obtain this new configuration recall from
   Lemma \ref{L:uv} that,  if the disks are not absolutely distinct,
   then the image of one disk is included in that of some other and we replace the
   smaller disk with the bigger one. We then argue that this configuration belongs
   to a moduli space of negative virtual dimension and strictly lower Maslov class
   than the initial one and, by induction, such a configuration can not exist which leads to
   a contradiction.  Coming back to our moduli spaces $\mathcal{P}_{\eqref{eq:m-pr1}}$
   we see that if we assume that some disks in an element of this moduli space are not
   absolutely distinct, then the reduction described in Case 2 \S \ref{Sb:prf-tr-1}  leads to a
   configuration of lower Maslov class which does not necessarily lie in
    $\mathcal{P}_{\eqref{eq:m-pr1}}$ but is an element of the more general moduli spaces
    introduced above. With these
   more general moduli spaces the reduction to simple disks mentioned
   above still works without problems and, additionally, the
   combinatorial argument in Case 2 in the proof of Proposition
   \ref{P:tr-1} also adapts in an obvious way and this proves the
   statement when $\dim (L)=n\geq 3$.  In the case $n\leq 2$ it is
   important to note that the only Maslov indexes involved are still
   $2$ and $4$. This allows for the last part - valid for $n\leq 2$ -
   of the proof of Proposition \ref{P:tr-1} to adapt to this case.
   \qed

   The second step of the proof of the proposition is to show by using
   the previous lemma that the operation $\ast$ provides a chain map.
   To do this we need to consider the compactification
   $\overline{\mathcal{P}}_{\eqref{eq:m-pr1}}$ of our moduli spaces
   and we need to show the appropriate analogue of the
   Lemma~\ref{lem:compact}.  Here is the statement.
   \begin{lem}\label{lem:compact2}
      For $x,y,z\in \Crit (f)$, $\delta'=1$ and a generic almost
      complex structure $J$, the space
      $\overline{\mathcal{P}}_{\eqref{eq:m-pr1}}(x,y,z) =
      \overline{\mathcal{P}}(x,y,z,\mathbf{A},\mathbf{A'},
      \mathbf{A''},U; f,f',f'',J)$ is a compact, $1$-dimensional
      manifold whose boundary verifies:
      \begin{equation}\label{eq:prod-comp}
         \begin{array}{l}\partial
            \overline{\mathcal{P}}_{\eqref{eq:m-pr1}}(x,y,z) =
            \cup_{x'}\mathcal{P}(x,x')\times
            \mathcal{P}_{\eqref{eq:m-pr1}}(x',y,z)\bigcup \\
            \bigcup \ \cup_{y'}
            \mathcal{P}(y,y')\times
            \mathcal{P}_{\eqref{eq:m-pr1}}(x,y',z)\
            \bigcup\  \cup_{z'}\mathcal{P}_{\eqref{eq:m-pr1}}(x,y,z')
            \times
            \mathcal{P}(z',z)\cup
            \mathcal{R}\cup\mathcal{R}'\cup\mathcal{R}''
         \end{array}
      \end{equation}
      Where $\mathcal{R}$ assembles the terms $\xi=(u,u',u'',v)$ so
      that one edge in one of $u,u',u''$ is of length $0$;
      $\mathcal{R}'$ assembles the terms in which one of the disks in
      one of the classes $A_{i},A'_{j},A''_{k}$ splits in two (each
      piece being a non-constant disk carrying two incidence points);
      $\mathcal{R}''$ assembles the terms in which the core splits in
      two (with possibly one of the pieces being a stable ghost disk
      and each piece carrying two incidence points).  The unions in
      the first three terms are taken over all possible splittings
      such as to respect the homology classes
      $\mathbf{A},\mathbf{A'},\mathbf{A''},U$ (these classes have been
      omitted in the notation).
   \end{lem}

   The proof of this lemma is very similar to that of Lemma
   \ref{lem:compact} so that we will omit the details besides
   indicating the points where some differences occur. Of course, the
   condition $N_L \geq 2$ is crucial in insuring that if bubbling off
   of some disk occurs each piece will carry at least two incidence
   points. The main difference concerns the set $\mathcal{R}''$. This
   set consists of configurations associated to the bubbling off of
   the core. It is important to notice that as the core carries three
   marked points (which geometrically correspond to the three
   attachment points $a_{i}, i=1,2,3$) the bubbling off of a ``ghost"
   disk has to be allowed as long as this disk is stable.

   Given this result we proceed to show that $\ast$ induces a chain
   map.  For this it is enough to show that the number of elements in
   $\partial\overline{\mathcal{P}}(x,y,z;\la)$ is the same as the sum
   $S$ of the number of elements in the first three terms
   in~\eqref{eq:prod-comp} when $\mathbf{A},\mathbf{A}',\mathbf{A}'',
   U$ vary such that $$\sum A_{i}+\sum A_{j}'+\sum A_{k}''+U=\la~.~$$
   Indeed, this implies that $S=0$ which gives precisely the algebraic
   identity equivalent to $\ast$ being a chain map. Thus the proof is
   reduced to showing that each element coming from the terms
   $\mathcal{R}$, $\mathcal{R}'$, $\mathcal{R}''$ appears twice when
   $\mathbf{A},\mathbf{A}',\mathbf{A}'', U$ vary as above.  To check
   this, again, the only case which is different from the proof of
   Proposition \ref{prop:pearl_cplx} concerns the terms of type
   $\mathcal{R}''$.  But it is easy to see that the condition $N_L
   \geq 2$ insures that each element of type $\mathcal{R}''$ can also
   be viewed as an element in a set of type $\mathcal{R}$ obtained in
   a configuration with a trivial core when one of the flow lines
   attached to the core is of length zero. For this argument it is
   useful to notice that, generically, the elements of type
   $\mathcal{R}''$ have the property that if a ghost disk has bubbled
   off in the core then this ghost disk carries at most two of the
   marked points.

   \

   Therefore, $\ast$ is a chain map and thus induces an operation in
   homology. By using the same techniques as above combined with the
   invariance proof from \S\ref{subsubsec:inv} we obtain that at the
   homology level this product is independent of the choices made in
   its definition.

   \begin{lem} \label{L:QHL-unity}
      There exists a canonical element $w_L \in Q^{+}H_n(L)$ (resp.
      $QH_n(L)$) which is a unit with respect to the quantum cap
      product, i.e. $w_L * \alpha = \alpha$ for every $\alpha \in
      Q^{+}H(L)$ (resp. $\alpha \in QH(L)$).
   \end{lem}
   \begin{proof}
      It is easy to see that we may take $f''=f'$ in the definition of
      the product as before. Assume now that $f$ has a single maximum,
      $x_n$, and take generic $J$ so that the pearl complexes are
      defined. By Remark \ref{rem:diff_m}, $x_n$ is always a cycle and
      we will denote its homology class by $w_L$. Under these
      assumptions it is easy to see that, at the chain level, we have
      $x_n \ast y = y$ for all critical points $y \in \Crit(f')$ .
      This means that $w_L$ is the unit for our product $*$.

      By standard (essentially Morse-theoretic) arguments it follows
      that $w_L$ is indeed canonical in the sense that when we change
      our Morse function the identification morphism $\phi$ (described
      in \S\ref{subsubsec:inv}) preserves that homology class.
   \end{proof}
   \begin{rem}\label{rem:prod_topclass}
      a. In view of the proof of Lemma~\ref{L:QHL-unity} we will
      sometimes denote $w_L$ by abuse of notation also as $[L]$.

      b. Let $f$ be a Morse function with a single maximum $x_n$.
      Since $w_L=[x_n]$ is the unit it follows that $QH_{\ast}(L)=0$
      iff $x_n$ is a boundary in $\mathcal{C}(L;f,J)$.  Indeed,
      suppose that $x_n=d(\eta)$ and let $\alpha$ be a cycle. Then we
      have $a = x_n \ast a = (d\eta)\ast a=d(\eta\ast a)$ (see also Remark \ref{rem:product_min}
      for the same statement in the context of the minimal models).

      c. It is useful to note that the product described above is not
      commutative (even in the graded sense) in general. See
      Proposition~\ref{prop:clifford} in~\S\ref{Ss:clifford} for a
      concrete example.
   \end{rem}

   The only point left to conclude the proof of Proposition
   \ref{prop:quantum_prod} is that, in homology, this product is
   associative.

   \begin{lem}\label{lem:prodassoc}
      The product $$\ast : Q^{+}H_{k}(L)\otimes Q^{+}H_{i}(L)\to
      Q^{+}H_{i+k-n}(L)$$
      is associative.
   \end{lem}

   {\em Proof of Lemma \ref{lem:prodassoc}.} To prove this result some
   new moduli spaces need to be used. We remark that the pearl moduli
   spaces $\mathcal{P}(x,y,\la)$ may be viewed as modeled over linear
   trees (with oriented edges). Similarly, the moduli spaces
   $\mathcal{P}(x,y,z;\la)$ used to define the product are modeled
   over trees with two entries and one exit (and hence with a single
   vertex of valence three).

   The moduli spaces needed to prove the associativity are modeled
   over more general trees and we describe them rigorously now.

   We consider trees $\mathcal{T}$ with oriented edges embedded in
   $\R\times[0,1]\subset \R^{2}$ with three entries lying on the line
   $\R\times\{1\}$ and one exit on the line $\R\times\{0\}$ and so
   that the edges strictly decrease the $y$-coordinate.  The vertices
   of the tree - except for the entries and the exit - are labeled by
   elements of $H_2(M,L;\mathbb{Z})$ and the label of each such vertex
   will be called its {\em class}. The entries are labeled in order by
   the three critical points $x,y,z$ and the exit is labeled by $w$.
   Each edge is labeled by an element of the set $\{1,2,3\}$. Clearly,
   such a tree $\mathcal{T}$ has either two vertices of valence three
   or one vertex of valence four and each internal vertex has a single
   exit.  Fixing $x,y,z,w\in\Crit(f)$ and such a tree $\mathcal{T}$ we
   denote the associated moduli space by
   $\mathcal{P}_{\mathcal{T}}(x,y,z,w)$ (the rule here is that the
   last critical point is the exit; this is coherent with previous
   notation).  An element of this moduli space consists of a family of
   $J$-holomorphic disks one for each vertex of the tree, in the class
   of that vertex, together with a family of strictly positive real
   numbers, one for each edge in the tree, so that for each edge we
   have an incidence relation like the one in Definition
   \ref{def:moduli1} iii. This relates the vertices joined by the
   respective edge but instead of the flow $\gamma$ one might use any
   one of the negative gradient flows $\gamma_{1}=\gamma$,
   $\gamma_{2}$ induced by $-\nabla f'$ or $\gamma_{3}$ induced by
   $-\nabla f''$ so that the flow $\gamma_{i}$ is used precisely when
   the label of the edge is $i$.  The incidence points are as follows:
   for the vertices of valence two they are the points $-1,+1$ so that
   the entry corresponds to $-1$ and the exit to $+1$; for the
   vertices of valence three they are, in cyclic order, the roots of
   unity of order $3$ so that $e^{2i\pi/3}$ corresponds to the the
   entry at the left, $e^{4i\pi/3}$ corresponds to the exit, $1$
   corresponds to the entry on the right; the incidence points
   $z_{1},z_{2},z_{3},z_{4}$ on a disk of valence four are - in order
   - so that $z_{1},z_{2},z_{4}$ are the roots of unity of order $3$,
   ($z_{1}=e^{2i\pi/3}$), $z_{3}$ is a point strictly in between
   $z_{2}$ and $z_{4}$ and $z_{1},z_{3},z_{4}$ are, from left to
   right, the entries and $z_{2}$ is the exit.  All the $J$-disks are
   stable in the sense that if a disk is trivial then it carries at
   least three incidence (or marked) points and, moreover, the
   labeling of the edges arriving at any vertex are pairwise distinct.
   Finally, a last condition is necessary in defining
   $\mathcal{P}_{\mathcal{T}}(x,y,z,w)$:
   \begin{equation}\label{eq:rule}
      \begin{array}{l}
         {\rm the\ labeling\ of\ the\ edges\ arriving\ at
           \ each\ vertex\ respects}\\
         {\rm
           the\ planar\ order;\ the\ label\
           of\ the\ exiting\ edge\ equals\ the}\\
         {\rm smallest\ of\ the\ labels\ of\ the\
           arriving\ edges }; {\rm the\
           edge\ starting\ in\ } x\ {\rm has\ label\ } 1,\\
         {\rm the\ edge\ starting\ in}\
         y\ {\rm has\ label\ } 2\
         {\rm and\ the\ edge\ starting\ in\ } z\ {\rm has\ label\ } 3.
      \end{array}
   \end{equation}

   This last condition implies that the labeling of the edges of a
   tree $\mathcal{T}$ is completely determined by its topological
   type. It will be useful in the arguments below to also consider
   moduli spaces defined exactly as above but without imposing
   condition~\eqref{eq:rule}. These more general moduli spaces will be
   denoted by $\mathcal{G}_{\mathcal{T}}(x,y,z,w)$.

   \begin{rem}\label{rem:gen_moduli}Notice that these more general moduli
      spaces have already appeared - in the case of just two entries -
      in the proof of Lemma \ref{lem:transv-prod}. Clearly, the
      definition described above easily extends to trees with more
      entries.
   \end{rem}

   The way to proceed from this point is clear: we first define moduli
   spaces $\mathcal{P}_{\mathcal{T}}^{\ast,d}(x,y,z,w)$ which are as
   above except that all the $J$-disks involved are simple and
   absolutely distinct.  Similarly we define moduli spaces
   $\mathcal{G}_{\mathcal{T}}^{\ast,d}(x,y,z,w)$.  By arguments
   similar to those in \S\ref{S:transversality} and in the proof of
   Lemma \ref{lem:transv-prod} it is not difficult to see that for
   $n\geq 3$ and a generic choice of almost complex structure $J$
   $$\mathcal{G}_{\mathcal{T}}(x,y,z,w) =
   \mathcal{G}_{\mathcal{T}}^{\ast,d}(x,y,z,w)$$ whenever
   $\delta''=|x|+|y|+|z|+\mu(\mathcal{T})-|w|-2n+1\leq 1$ where
   $\mu(\mathcal{T})$ is the sum of the Maslov classes of all the
   labels of the internal vertices in $\mathcal{T}$. Moreover,
   $\mathcal{G}_{\mathcal{T}}^{\ast,d}(x,y,z,w)$ is a manifold of
   dimension $\delta''$. The same combinatorial arguments as in the
   proof of Proposition \ref{P:tr-1} case $n\leq 2$ are sufficient to
   also show the same statement for $n=2$ when $N_L \geq 3$. For $n=2$
   and $N_L =2$ a more involved argument is needed because we need to
   consider the case of $\mu(\mathcal{T})=6$. In all cases, we deduce
   that for a generic $J$ and $\delta''\leq 1$
   \begin{equation}\label{eq:moduli_assoc_tr}
      \mathcal{P}_{\mathcal{T}}(x,y,z,w) =
      \mathcal{P}_{\mathcal{T}}^{\ast,d}(x,y,z,w)
   \end{equation}
   and this moduli space is a manifold of dimension $\delta''$.

   It is then needed to consider the compactifications
   $\overline{\mathcal{P}}_{\mathcal{T}}(x,y,z,w)$ and establish a
   boundary formula as in Lemma \ref{lem:compact2}. We will not state
   this formula explicitly as it is very similar to the ones before
   but we will describe the boundary terms. They are defined as the
   usual elements in $\mathcal{P}_{\mathcal{T}}(x,y,z,w)$ except for
   precisely one modification which fits into one of the following
   categories:
   \begin{itemize}
     \item[$i_{\mathcal{T}}$.] One of the edges in $\mathcal{T}$
      corresponds to a flow line of $0$ length.
     \item[$ii_{\mathcal{T}}$.] One of the disks corresponding to a
      vertex in $T$ is replaced by a cusp curve with two components
      (due to bubbling off); the bubbling off of a ghost disk is
      possible if the vertex in question is of valence at least $3$.
      Each of the pieces appearing in a cusp curve carries at least
      two incidence points.
     \item[$iii_{\mathcal{T}}$.] One of the edges in $\mathcal{T}$
      corresponds to a flow line which is broken once.
   \end{itemize}

   Notice that, in this description, the coincidence of two marked
   points (which clearly may occur on a vertex of valence four)
   corresponds to the bubbling off of a ghost disk and, clearly, the
   reason why each component in a cusp curve as in $ii_{\mathcal{T}}$
   carries at least two incidence points is that $N_L \geq 2$ and
   $\delta''\leq 1$.

   The purpose of this construction is, of course, to define a chain
   homotopy:
   $$\xi : \mathcal{C}^{+}(L;f,J)\otimes\mathcal{C}^{+}(L;
   f',J)\otimes\mathcal{C}^{+}(L;f'',J)\to\mathcal{C}^{+}(L; f,J)$$
   so
   that $\xi : (( -\ \ast\ -)\ \ast\ -)\simeq (-\ \ast\ (-\ \ast\
   -))$.  The definition of $\xi$ is clear:
   $$\xi (x\otimes y\otimes z)=\sum_{\la,\ |\mathcal{T}|=\la,\ w}
   \#_{\mathbb{Z}_2} (\mathcal{P}_{\mathcal{T}}(x,y,z,w))we^{\la}$$
   where the sum is over all those terms so that
   $\delta''=|x|+|y|+|z|-|w|+1-2n+\mu(\la)=0$ and the sum of all the
   labels of the internal vertices of $\mathcal{T}$ is denoted by
   $|\mathcal{T}|$.  We now use the description of the boundary of
   $\overline{\mathcal{P}}_{\mathcal{T}}(x,y,z,w)$ given above to
   justify
   \begin{equation}\label{eq:chain_htpy}
      (d\xi+\xi d)(x\otimes y\otimes z)= (x\ast y)
      \ast z - x\ast (y\ast z)~.~
   \end{equation}

   We let $$\overline{\mathcal{P}}(x,y,z,w;\la)$$
   be the (disjoint)
   union of all the spaces
   $\overline{\mathcal{P}}_{\mathcal{T}}(x,y,z,w)$ where $\mathcal{T}$
   runs over all the topological types of trees so that
   $|\mathcal{T}|=\la$ . To prove~\eqref{eq:chain_htpy} we start by
   noticing that the coefficients of $w t^{\mubar(\lambda)}$ in this
   equation are in bijection with the union of the terms
   $iii_{\mathcal{T}}$.  Indeed, the terms on the left
   in~\eqref{eq:chain_htpy} correspond to those elements in which the
   break in the flow line disconnects the tree so that one of the
   remaining pieces is a linear tree (in the sense that it does not
   contain any vertex of valence three or more) and the terms on the
   right correspond precisely to those elements in which the tree is
   broken in two parts none of which is linear.  Thus, to finish the
   proof it is enough to notice that the terms of type
   $i_{\mathcal{T}}$ and $ii_{\mathcal{T}}$ all cancel out when
   $\mathcal{T}$ varies.

   For this we look at a configuration of type $i_{\mathcal{T}}$ and
   remark that if one of the edges of $\mathcal{T}$ becomes of zero
   length, then, generically, this configuration also appears as a
   configuration of type $ii_{\mathcal{T}'}$ corresponding to a tree
   $\mathcal{T}'$ so that in $\mathcal{T'}$ is obtained from
   $\mathcal{T}$ by replacing the edge in question together with the
   two bounding vertices by a single vertex. The fact that the
   labeling of the edges of each topological type of tree is unique is
   essential here as it implies that this tree $\mathcal{T}'$ is
   unique. Therefore, all elements of type $i_{\mathcal{T}}$ and
   $ii_{\mathcal{T}}$ cancel and this proves the lemma and concludes
   the proof of the proposition.

\end{proof}

\subsection{The quantum module structure} \label{S:qm}
In this section we define an operation
$$\ast : Q^{+}H_{k}(M)\otimes Q^{+}H_{s}(L)\to Q^{+}H_{k+s-2n}(L),
\quad a\otimes x \longmapsto a*x,$$
which will make $Q^+H_*(L)$ an
algebra over $Q^+H_*(M)$.

Let $h:M \to \mathbb{R}$, $f:L \to \mathbb{R}$ be Morse functions and
$\rho_M$, $\rho_L$ Riemannian metrics on $M$, $L$.  From now on we
will make the following assumptions on $(f,\rho_L)$ and $(h,\rho_M)$:
\begin{assumption} \label{A:MS}
   Each of the pairs $(f,\rho_L)$ and $(h,\rho_M)$ is Morse-Smale and
   $h$ has a single maximum.  Furthermore:
   \begin{enumerate}
     \item For every $a \in \textnormal{Crit}(h)$ the unstable
      submanifold $W_a^u$ is transverse to $L$.
     \item For every $a \in \textnormal{Crit}(h)$, $x,y \in
      \textnormal{Crit}(f)$, $W_a^u$ is transverse to $W_x^u$ and
      $W_y^s$.
   \end{enumerate}
\end{assumption}

\subsubsection{The external operation.} \label{Sb:external-op}
We will describe here the relevant moduli spaces very explicitly and
without using the ``tree" model which has been used at the end of the
previous section.

Let $J \in \mathcal{J}(M,\omega)$ and $A \in H_2(M,L;\mathbb{Z})$.
Define the following evaluation maps:
\begin{alignat*}{2}
   &ev_{-1,1}:\mathcal{M}(A,J) \to L \times L, & \quad
   &ev_{-1,0,1}:\mathcal{M}(A,J) \to L \times M \times L, \\
   &ev_{-1,1}(u) = (u(-1),u(1)), & \quad &ev_{-1,0,1}(u) =
   (u(-1),u(0),u(1)).
\end{alignat*}

Let $\mathbf{A}=(A_1, \ldots, A_l)$ be a vector of {\em non-zero}
classes in $H_2(M,L;\mathbb{Z})$. Put
\begin{align*}
   & \mathcal{M}(\mathbf{A},J) = \mathcal{M}(A_1,J) \times \cdots
   \times
   \mathcal{M}(A_l,J), \\
   & ev: \mathcal{M}(\mathbf{A},J) \to L^{\times 2l}, \quad ev(u_1,
   \ldots, u_l) = (ev_{-1,1}(u_1), \ldots, ev_{-1,1}(u_l)).
\end{align*}
Given $x,y \in \textnormal{Crit}(f)$ recall from previous sections the
notation for the pearl moduli-spaces:
\begin{equation*}
   \mathcal{M}(x,y;\mathbf{A},J) = ev^{-1}(W_x^u \times
   Q_{f,\rho_L}^{\times (l-1)} \times W_y^s), \quad
   \mathcal{P}(x,y;\mathbf{A},J) =
   \mathcal{M}(x,y;\mathbf{A},J)/\mathbf{G},
\end{equation*}
where $\mathbf{G} = G_{-1,1} \times \cdots \times G_{-1,1}$ (see
figure~\ref{f:qm-pearls-1}).

\begin{figure}[htbp]
   \begin{center}
      \epsfig{file=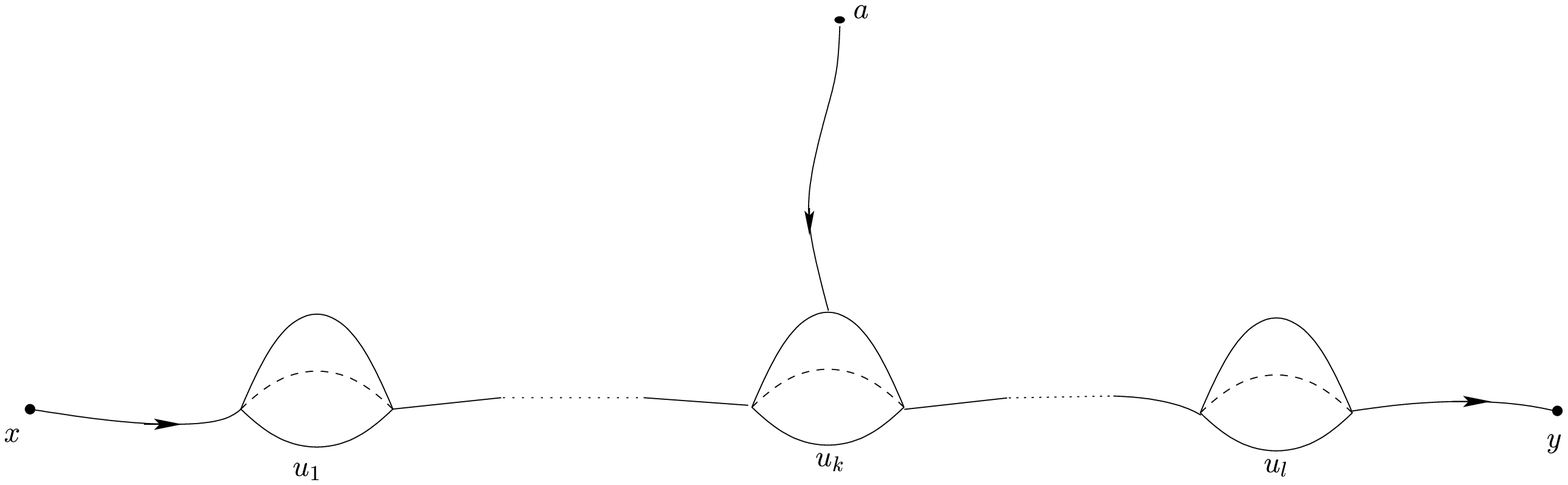, width=0.6\linewidth}
   \end{center}
   \caption{An element of $\mathcal{P}_I(a,x,y;(\mathbf{A},k),J)$}
   \label{f:qm-module-I}
\end{figure}

Now, given $1 \leq k \leq l$ define the evaluation map:
\begin{align*}
   & ev_{(k)}: \mathcal{M}(\mathbf{A},J) \to L^{\times (2k-1)}
   \times M \times L^{\times (2l-2k+1)}, \\
   & ev_{(k)}(u_1, \ldots, u_l) = (ev_{-1,1}(u_1), \dots,
   ev_{-1,1}(u_{k-1}), ev_{-1,0,1}(u_k), ev_{-1,1}(u_{k+1}), \dots,
   ev_{-1,1}(u_l)).
\end{align*}
Given $x,y \in \textnormal{Crit}(f)$ and $a \in \textnormal{Crit}(h)$
put
\begin{align} \label{Eq:qm-PI}
   & \mathcal{P}_I(a,x,y;(\mathbf{A},k),J) = ev_{(k)}^{-1}\bigl( W_x^u
   \times Q_{f,\rho_L}^{\times (k-1)} \times W_a^u \times
   Q_{f,\rho_L}^{\times (l-k)} \times W_y^s \bigr)/\mathbf{G}_k, \\
   & \mathcal{P}_I(a,x,y;\mathbf{A},J) = \bigcup_{k=1}^l
   \mathcal{P}_I(a,x,y;(\mathbf{A},k),J),
\end{align}
where $\mathbf{G}_k$ is taken here to be $(G_{-1,1})^{\times k}$ with
the $k$'th factor replaced by the trivial group.  See
figure~\ref{f:qm-module-I}. For each element
$u\in\mathcal{P}_{I}(a,x,y;(\mathbf{A},k),J)$ we will call the disk
$u_{k}$ the \emph{center} of $u$.

Denote by $\Phi_t:L \to L$, $t \in \mathbb{R}$, the {\em negative}
gradient flow of $(f,\rho_L)$ (i.e.  the flow of the vector field
$-\textnormal{grad}_{\rho_L}f$). Consider the (non-proper) embedding
$$(L \setminus \textnormal{Crit}(f)) \times \mathbb{R}_{>0} \times
\mathbb{R}_{>0} \hooklongrightarrow L \times L \times L, \quad (x,t,s)
\longmapsto (x,\Phi_t(x), \Phi_{t+s}(x)).$$
Denote the image of this
embedding by $Q'_{f,\rho_L} \subset L\times L \times L$. We will also
need the subset $Q_{f,\rho_L} \subset L \times L$ defined
by~\eqref{Eq:Q-f} in \S\ref{Sb:tr-pearl-complex}.

\begin{figure}[htbp]
   \begin{center}
      \epsfig{file=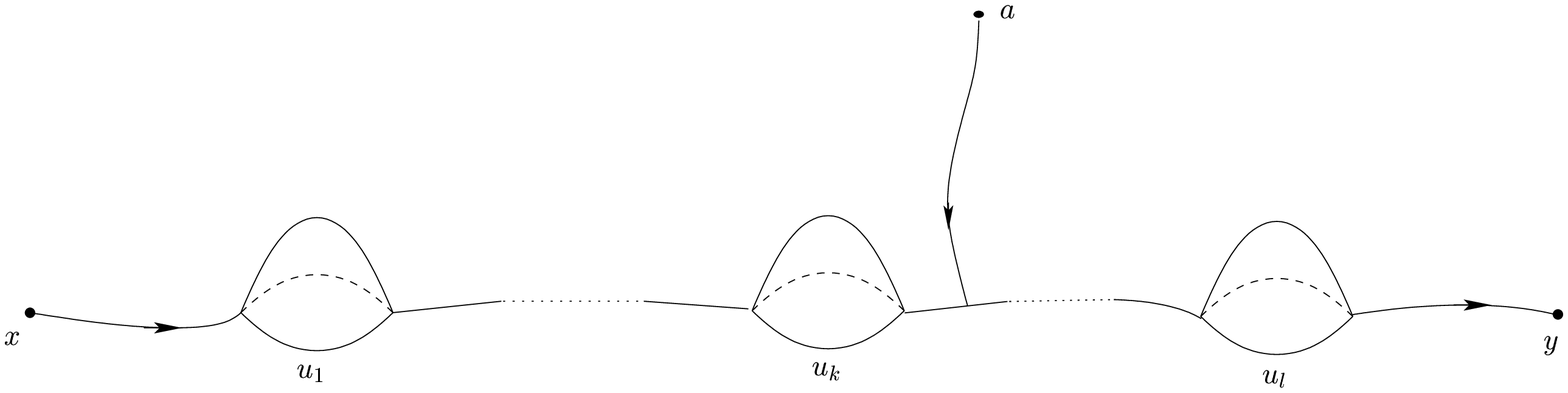, width=0.6\linewidth}
   \end{center}
   \caption{An element of $\mathcal{P}_{I'}(a,x,y;(\mathbf{A},k),J)$}
   \label{f:qm-module-Ip}
\end{figure}

Let $\mathbf{A}=(A_1, \ldots, A_l)$ be a vector of non-zero classes,
$0 \leq k \leq l$, $x, y \in \textnormal{Crit}(f)$ and $a \in
\textnormal{Crit}(h)$. Define the following space (see
figure~\ref{f:qm-module-Ip}):
\begin{align*}
   \mathcal{P}_{I'} &(a,x,y;(\mathbf{A},k),J) = \\
   &
   \begin{cases}
      \bigl\{ (\mathbf{u},p) \in \mathcal{P}(x,y;\mathbf{A},J) \times
      (W_a^u \cap L) \mid (u_k(1), p, u_{k+1}(-1)) \in
      Q'_{f,\rho_L} \bigr\} & \textnormal{if } 0< k <l, \\
      \bigl\{ (\mathbf{u},p) \in \mathcal{P}(x,y;\mathbf{A},J) \times
      (W_a^u \cap W_x^u) \mid (p,u_1(-1)) \in Q_{f,\rho_L} \bigr\} &
      \textnormal{if } k=0, \\
      \bigl\{ (\mathbf{u},p) \in \mathcal{P}(x,y;\mathbf{A},J) \times
      (W_a^u \cap W_y^s) \mid (u_l(1),p) \in Q_{f,\rho_L} \bigr\} &
      \textnormal{if } k=l.
   \end{cases} \\
   & \mathcal{P}_{I'}(a,x,y;\mathbf{A},J) = \bigcup_{k=0}^l
   \mathcal{P}_{I'}(a,x,y;(\mathbf{A},k),J).
\end{align*}
Note that the $p$ in $(\mathbf{u},p)$ cannot be neither $x$ nor $y$
due to the definition of $Q_{f,\rho_L}$.

\begin{rem}
   It is clear that we could as well define the space
   $\mathcal{P}_{I'}(--)$ as a subspace of $\mathcal{P}_{I}(--)$
   consisting of those configurations so that their center is the
   constant disk.
\end{rem}

\begin{prop} \label{P:qm-axy}
   Let $(f,\rho_L)$, $(h,\rho_M)$ be as above. There exists a second
   category subset $\mathcal{J}_{\textnormal{reg}} \subset
   \mathcal{J}(M, \omega)$ with the following properties. For every $J
   \in \mathcal{J}_{\textnormal{reg}}$, $\mathbf{A}$ as above, and
   every $x,y \in \textnormal{Crit}(f)$, $a \in \textnormal{Crit}(h)$
   with $|a| + |x|-|y| + \mu(\mathbf{A})-2n=0$ each of the spaces
   $\mathcal{P}_I(a,x,y;\mathbf{A},J)$,
   $\mathcal{P}_{I'}(a,x,y;\mathbf{A},J)$ is a compact $0$-dimensional
   manifold, hence a finite set. Furthermore, each element of these
   spaces consists of simple and absolutely distinct disks.
\end{prop}
The proof is given in Section~\ref{Sb:qm-prfs-props} below.

In case $|a| + |x|-|y| + \mu(\mathbf{A})-2n =0$ put:
\begin{align*}
   & n_I(a,x,y;\mathbf{A},J) = \#_{\mathbb{Z}_2}
   \mathcal{P}_I(a,x,y;\mathbf{A},J), \quad n_{I'}(a,x,y;\mathbf{A},J)
   = \#_{\mathbb{Z}_2}
   \mathcal{P}_{I'}(a,x,y;\mathbf{A},J), \\
   & n(a,x,y;\mathbf{A},J) = n_I(a,x,y;\mathbf{A},J) +
   n_{I'}(a,x,y;\mathbf{A},J).
\end{align*}
Finally, in case $|a| + |x|-|y| -2n =0$ put $n(a,x,y) =
\#_{\mathbb{Z}_2} (W_a^u \cap W_x^u \cap W_y^s)$.

Denote by $(C_*(h), \partial^h)$ the Morse complex of $(M; h,\rho_M)$
and by $(\mathcal{C}_*^{+}, d^f)$ the pearl complex associated to $(L;
f,\rho_L, J)$. Define an operation:
\begin{equation} \label{Eq:qm}
   \ast:C_{k}(h)\otimes \mathcal{C}^{+}_{q} \longrightarrow
   \mathcal{C}^{+}_{q+k-2n},
\end{equation}
by the formula
\begin{equation} \label{Eq:qm-*}
   a*x = \sum_{y} n(a,x,y)y +
   \sum_{y,\mathbf{A}}
   n(a,x,y;\mathbf{A},J) y t^{\mubar(\mathbf{A})},
\end{equation}
where the first sum is taken over all $y \in \textnormal{Crit}(f)$
with $|a| + |x|-|y| -2n =0$ and the second sum over all $y,
\mathbf{A}$ with $|a| + |x|-|y| + \mu(\mathbf{A})-2n =0$.  Note that
the first sum is the Morse theoretic interpretation of the classical
cap product operation of $H_*(M)$ on $H_*(L)$.

\begin{prop} \label{P:qm-chain}
   The homomorphism~\eqref{Eq:qm} is a chain map, namely
   $$d^f(a*x) = \partial^h(a) * x + a*d^f(x).$$
   In particular it
   induces a well defined operation:
   $$\ast : Q^{+}H_{k}(M)\otimes Q^{+}H_{s}(L)\to
   Q^{+}H_{k+s-2n}(L).$$
\end{prop}

The proof of Proposition~\ref{P:qm-chain} will occupy
Sections~\ref{Sb:qm-spaces} --~\ref{Sb:conclusion-prf-qm-ident-all-N}
below.

\subsubsection{Moduli spaces related to the external operation}
\label{Sb:qm-spaces}
In order to prove Proposition~\ref{P:qm-chain} we will need to
introduce several types of moduli spaces.

\subsubsection*{Type II}
\begin{figure}[htbp]
   \begin{center}
      \epsfig{file=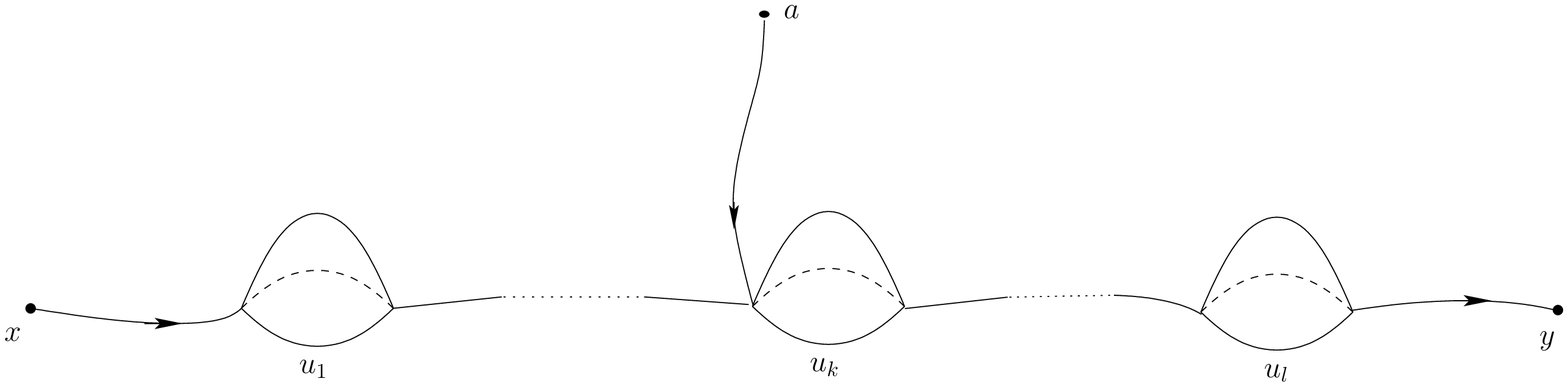, width=0.6\linewidth}
      \epsfig{file=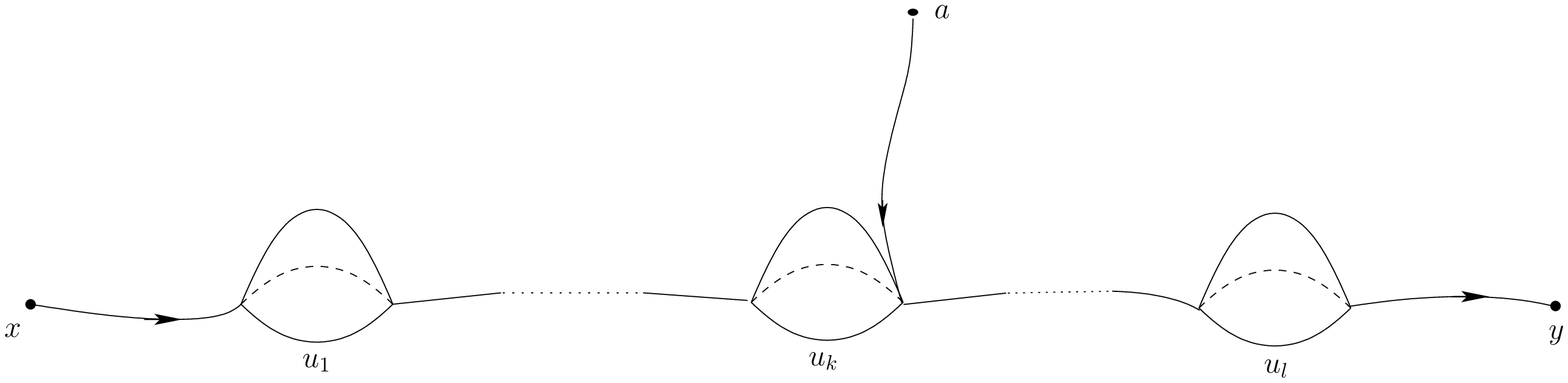, width=0.6\linewidth}
   \end{center}
   \caption{Elements of the spaces
     $\mathcal{P}_{II_1}(a,x,y;(\mathbf{A},k),J)$,
     $\mathcal{P}_{II_2}(a,x,y;(\mathbf{A},k),J)$.}
   \label{f:qm-module-II}
\end{figure}

Let $\mathbf{A}$ be a vector of non-zero classes and $1\leq k \leq l$.
Define the following spaces (see figure~\ref{f:qm-module-II}):
\begin{align*}
   & \mathcal{P}_{II_1}(a,x,y;(\mathbf{A},k),J) = \bigl\{ \mathbf{u}
   \in \mathcal{P}(x,y;\mathbf{A},J) \mid u_k(-1) \in W_a^u
   \bigr\}, \\
   & \mathcal{P}_{II_2}(a,x,y;(\mathbf{A},k),J) = \bigl\{ \mathbf{u}
   \in \mathcal{P}(x,y;\mathbf{A},J) \mid u_k(1) \in W_a^u
   \bigr\}, \\
   & \mathcal{P}_{II_1}(a,x,y;\mathbf{A},J) = \bigcup_{k=1}^l
   \mathcal{P}_{II_1}(a,x,y;(\mathbf{A},k),J), \quad
   \mathcal{P}_{II_2}(a,x,y;\mathbf{A},J) = \bigcup_{k=1}^l
   \mathcal{P}_{II_2}(a,x,y;(\mathbf{A},k),J).
\end{align*}

\subsubsection*{Type III}
\begin{figure}[htbp]
   \begin{center}
      \epsfig{file=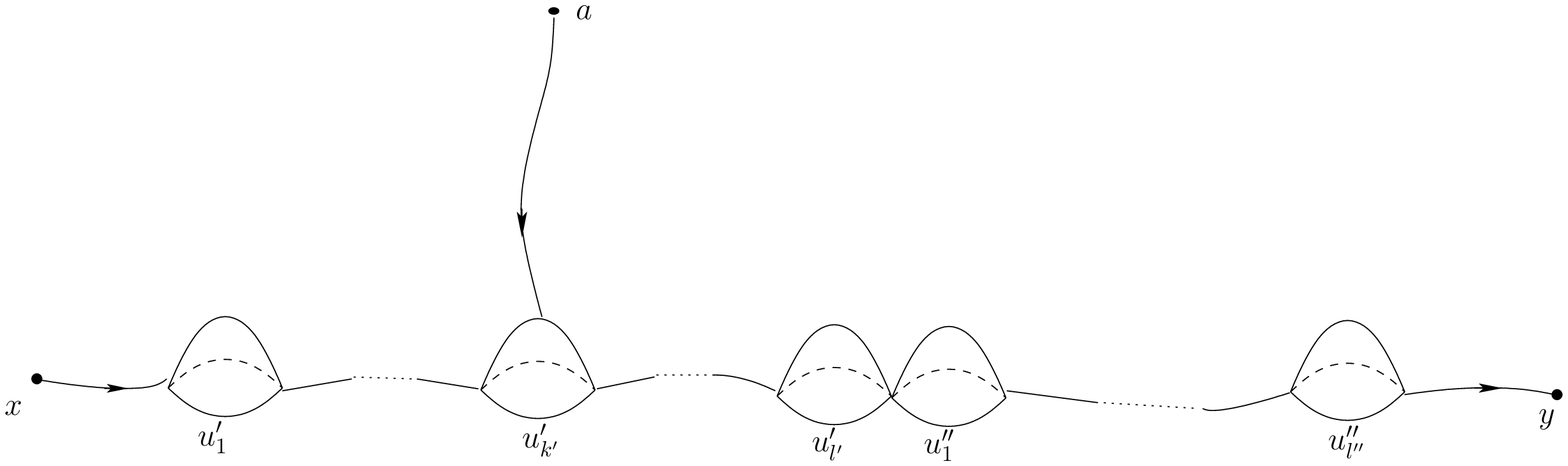, width=0.6\linewidth}
      \epsfig{file=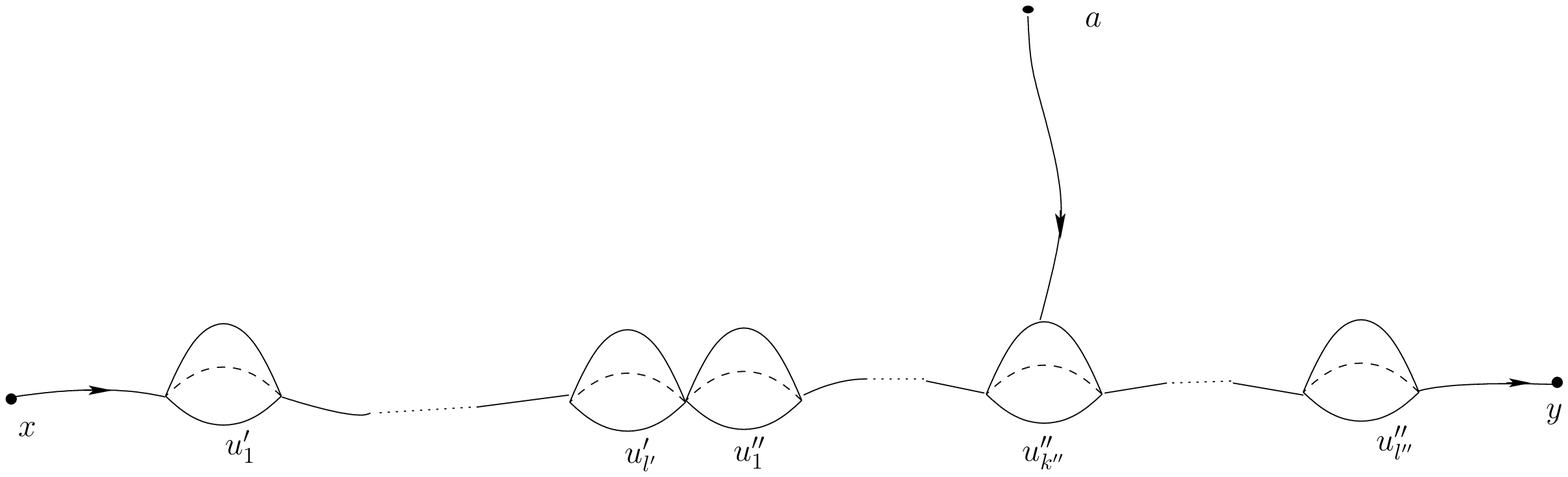, width=0.6\linewidth}
   \end{center}
   \caption{Elements of the spaces
     $\mathcal{P}_{III_1}(a,x,y;(\mathbf{B}',k'), \mathbf{B}'',J)$,
     $\mathcal{P}_{III_2}(a,x,y;\mathbf{B}',(\mathbf{B}'',k''),J)$.}
   \label{f:qm-module-III}
\end{figure}

Let $\mathbf{B}'=(B'_1, \ldots, B'_{l'})$, $\mathbf{B}''=(B''_1,
\ldots, B''_{l''})$ be two vectors of {\em non-zero} classes in
$H_2(M,L;\mathbb{Z})$. Let $1\leq k'\leq l'$ and $1\leq k''\leq l''$.
Put $\mathbf{A} = (\mathbf{B}',\mathbf{B}'') = (B'_1, \ldots, B'_{l'},
B''_1, \ldots, B''_{l''})$.  Define
\begin{equation} \label{Eq:P-III}
   \begin{aligned}
      & \mathcal{P}_{III_1}(a,x,y;(\mathbf{B}',k'), \mathbf{B}'',J)
      = \\
      & ev_{(k')}^{-1} \bigl( W_x^u \times Q_{f,\rho_L}^{\times
        (k'-1)} \times W_a^u \times Q_{f,\rho_L}^{\times (l'-k')}
      \times \textnormal{diag(L)} \times Q_{f,\rho_L}^{\times (l''-1)}
      \times
      W_y^s \bigr)/\mathbf{G}_{III_1}, \\
      & \mathcal{P}_{III_2}(a,x,y;\mathbf{B}',
      (\mathbf{B}'',k''),J) = \\
      & ev_{(l'+k'')}^{-1} \bigl( W_x^u \times Q_{f,\rho_L}^{\times
        (l'-1)} \times \textnormal{diag(L)} \times
      Q_{f,\rho_L}^{\times (k''-1)} \times W_a^u \times
      Q_{f,\rho_L}^{\times (l''-k'')} \times W_y^s
      \bigr)/\mathbf{G}_{III_2}.
   \end{aligned}
\end{equation}
Here $\mathbf{G}_{III_1}$ (respectively $\mathbf{G}_{III_2}$) is the
group $G_{-1,1}^{\times (l'+l'')}$ with the $k'$'th (respectively
$(l'+k'')$'th) component replaced by the trivial group. Finally put
\begin{align*}
   & \mathcal{P}_{III_1}(a,x,y;\mathbf{B}', \mathbf{B}'',J) =
   \bigcup_{k'=1}^{l'}\mathcal{P}_{III_1}(a,x,y;(\mathbf{B}',k'),
   \mathbf{B}'',J), \\
   & \mathcal{P}_{III_2}(a,x,y;\mathbf{B}', \mathbf{B}'',J) =
   \bigcup_{k''=1}^{l''}\mathcal{P}_{III_2}(a,x,y;\mathbf{B}',
   (\mathbf{B}'',k''),J).
\end{align*}
See figure~\ref{f:qm-module-III}.

\subsubsection*{Type III'}
Let $\mathbf{B}'=(B'_1, \ldots, B'_{l'})$, $\mathbf{B}''=(B''_1,
\ldots, B''_{l''})$ be two vectors of non-zero classes. Given $x, y
\in \textnormal{Crit}(f)$ put
\begin{equation*}
   \mathcal{P}(x,y; \mathbf{B}',\mathbf{B}'', J) =
   ev^{-1} \bigl( W_x^u \times
   Q_{f,\rho_L}^{\times (l'-1)} \times \textnormal{diag}(L) \times
   Q_{f,\rho_L}^{\times (l''-1)} \times W_y^s \bigr) / \mathbf{G},
\end{equation*}
where $ev:\mathcal{M}((\mathbf{B}',\mathbf{B}''),J) \to L^{\times
  (2l'+2l'')}$ is the evaluation map from Section~\ref{Sb:external-op}
and $\mathbf{G} = G_{-1,1}^{\times (l'+l'')}$.

Let $0 \leq k' < l'$.  Define the following space:
\begin{align*}
   \mathcal{P}_{III'_1} &(a,x,y;(\mathbf{B}',k'), \mathbf{B}'',J) = \\
   &
   \begin{cases}
      \bigl \{ (\mathbf{u}',\mathbf{u}'',p) \in \mathcal{P}(x,y;
      \mathbf{B}',\mathbf{B}'', J) \times (W_a^u \cap L) \mid
      (u'_{k'}(1),p,u'_{k'+1}) \in Q'_{f,\rho_L} \bigr\} &
      \textnormal{if } k'>0, \\
      \bigl\{ (\mathbf{u}', \mathbf{u}'',p) \in \mathcal{P}(x,y;
      \mathbf{B}',\mathbf{B}'', J) \times (W_a^u \cap W_x^u) \mid
      (p,u_1(-1)) \in Q_{f,\rho_L} \bigr\} & \textnormal{if } k'=0.
   \end{cases}
\end{align*}
Let $0<k''\leq l''$. Define:
\begin{align*}
   \mathcal{P}_{III'_2} &(a,x,y;\mathbf{B}',(\mathbf{B}'',k''),J) = \\
   &
   \begin{cases}
      \bigl \{ (\mathbf{u}',\mathbf{u}'',p) \in \mathcal{P}(x,y;
      \mathbf{B}',\mathbf{B}'', J) \times (W_a^u \cap L) \mid
      (u'_{k''}(1),p,u'_{k''+1}) \in Q'_{f,\rho_L} \bigr\} &
      \textnormal{if } k''<l'', \\
      \bigl\{ (\mathbf{u}', \mathbf{u}'',p) \in \mathcal{P}(x,y;
      \mathbf{B}',\mathbf{B}'', J) \times (W_a^u \cap W_y^s) \mid
      (u_{l''}(1), p) \in Q_{f,\rho_L} \bigr\} & \textnormal{if }
      k''=l''.
   \end{cases}
\end{align*}

\begin{figure}[htbp]
   \begin{center}
      \epsfig{file=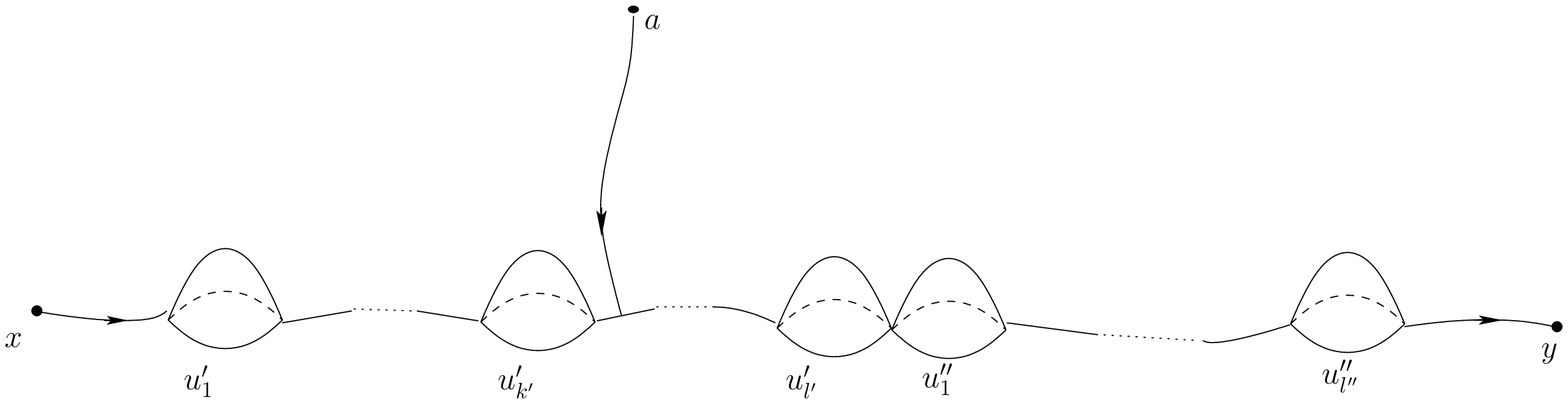, width=0.6\linewidth}
      \epsfig{file=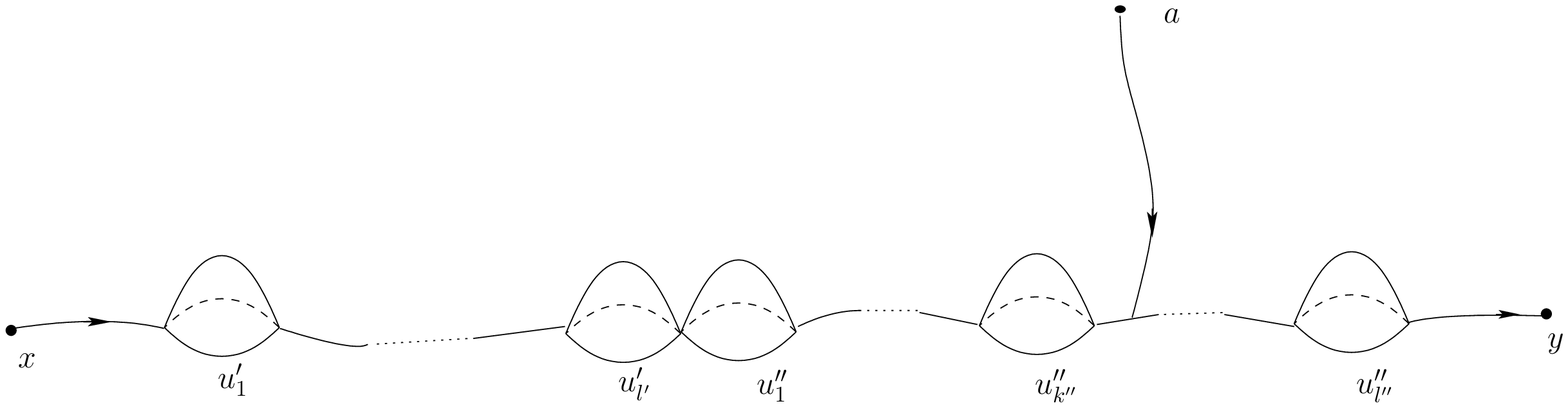, width=0.6\linewidth}
   \end{center}
   \caption{Elements of the spaces
     $\mathcal{P}_{III'_1}(a,x,y;(\mathbf{B}',k'), \mathbf{B}'',J)$,
     $\mathcal{P}_{III'_2}(a,x,y;\mathbf{B}',(\mathbf{B}'',k''),J)$.}
   \label{f:qm-module-III'}
\end{figure}

Finally put:
\begin{align*}
   & \mathcal{P}_{III'_1}(a,x,y;\mathbf{B}', \mathbf{B}'',J) =
   \bigcup_{k'= 0}^{l'-1}
   \mathcal{P}_{III'_1}(a,x,y;(\mathbf{B}',k'), \mathbf{B}'',J), \\
   & \mathcal{P}_{III'_2}(a,x,y;\mathbf{B}', \mathbf{B}'',J) =
   \bigcup_{k''= 1}^{l''}
   \mathcal{P}_{III'_2}(a,x,y;\mathbf{B}',(\mathbf{B}'',k''),J).
\end{align*}
See figure~\ref{f:qm-module-III'}.

\begin{prop} \label{P:qm-0-dim}
   There exists a second category subset
   $\mathcal{J}_{\textnormal{reg}} \subset \mathcal{J}(M,\omega)$ such
   that for every $J \in \mathcal{J}_{\textnormal{reg}}$ the following
   holds:
   \begin{enumerate}
     \item For every $a, x, y$ and every vector $\mathbf{A}$ of
      non-zero classes with $|a|+|x|-|y|+\mu(\mathbf{A})-2n=1$ each
      of the spaces $\mathcal{P}_{II_i}(a,x,y;\mathbf{A},J)$, $i=1,2$,
      is a compact $0$-dimensional manifold, hence a finite set.
      Moreover, every $(u_1, \ldots, u_l)$ in this space consists of
      simple and absolutely distinct disks.
     \item For every $a, x, y$ and every two vectors $\mathbf{B}'$,
      $\mathbf{B}''$ of non-zero classes with
      $|a|+|x|-|y|+\mu(\mathbf{B}')+\mu(\mathbf{B}'')-2n=1$ each
      of the spaces $\mathcal{P}_{III_i}(a,x,y;\mathbf{B}',
      \mathbf{B}'',J)$, $\mathcal{P}_{III'_i}(a,x,y;\mathbf{B}',
      \mathbf{B}'',J)$, $i=1,2$, is a compact $0$-dimensional
      manifold, hence a finite set. Moreover, every $(u'_1, \ldots,
      u'_{l'}, u''_1, \ldots, u''_{l''})$ in this space consists of
      simple and absolutely distinct disks.
   \end{enumerate}
\end{prop}

We put
\begin{align*}
   & n_{II_i}(a,x,y;\mathbf{A},J) = \#_{\mathbb{Z}_2}
   \mathcal{P}_{II_i}(a,x,y;\mathbf{A},J), \\
   & n_{III_i}(a,x,y;\mathbf{B}', \mathbf{B}'',J) = \#_{\mathbb{Z}_2}
   \mathcal{P}_{III_i}(a,x,y;\mathbf{B}', \mathbf{B}'',J), \\
   & n_{III_i}(a,x,y;\mathbf{B}', \mathbf{B}'',J) = \#_{\mathbb{Z}_2}
   \mathcal{P}_{III'_i}(a,x,y;\mathbf{B}', \mathbf{B}'',J),
\end{align*}
whenever $|a|+|x|-|y|+\mu(\mathbf{A})-2n=1$ or
$|a|+|x|-|y|+\mu(\mathbf{B}')+\mu(\mathbf{B}'')-2n=1$.

\subsubsection{Identities} \label{qm:identities}
Given $a, a' \in \textnormal{Crit}(h)$ with $|a'| = |a|-1$, set
$n(a,a')=\#_{\mathbb{Z}_2}(W_a^u \cap W_{a'}^s)/\mathbb{R}$, i.e. the
number modulo $2$ of (negative gradient) trajectories of $h$ going
from $a$ to $a'$. Similarly for $x,x' \in \textnormal{Crit}(f)$ we
have $n(x,x')$. In order to simplify the notation we will omit the
$J$'s from $n_I(a,x,y;\mathbf{A},J)$'s and from the $n_{II}$'s,
$n_{III}$'s etc. Given two vectors of non-zero classes
$\mathbf{B}'=(B'_1, \ldots, B'_{l'})$, $\mathbf{B}''=(B''_1, \ldots,
B''_{l''})$ write $$\mathbf{B}' \# \mathbf{B}'' = (B'_1, \ldots,
B'_{l'-1}, B'_{l'}+B''_1, B''_2, \ldots, B''_{l''}).$$

\begin{prop} \label{P:qm-idents}
   Let $(f,\rho_L)$, $(h,\rho_M)$ be as above.  There exists a second
   category subset $\mathcal{J}_{\textnormal{reg}} \subset
   \mathcal{J}(M,\omega)$ such that for every $J \in
   \mathcal{J}_{\textnormal{reg}}$, every $x,y \in
   \textnormal{Crit}(f)$, $a \in \textnormal{Crit}(h)$ and
   $\mathbf{A}$ with $|a|+|x|-|y| + \mu(\mathbf{A})-2n=1$ the
   following two identities hold:
   \begin{equation} \label{Eq:qm-ident-1}
      \begin{aligned}
         \sum_{|x'|=|x|-1} & n(x,x') n_I(a,x',y;\mathbf{A}) \, +
         \sum_{|y'|=|y|+1}
         n_I(a,x,y';\mathbf{A}) n(y',y) \, + \\
         \sum_{|a'|=|a|-1} & n(a,a') n_I(a',x,y;\mathbf{A}) \, +
         \sum_{\substack{(\mathbf{A'},\mathbf{A''})=\mathbf{A} \\
             |x'|=|x|+\mu(\mathbf{A}')-1}} n(x,x'; \mathbf{A}')
         n_I(a,x',y;\mathbf{A}'') \, + \\
         \sum_{\substack{(\mathbf{A'},\mathbf{A''})=\mathbf{A} \\
             |y'|=|y|-\mu(\mathbf{A}'')+1}} & n_I(a,x,y';\mathbf{A}')
         n(y',y; \mathbf{A}'') \, + n_{II_1}(a,x,y;\mathbf{A}) +
         n_{II_2}(a,x,y;\mathbf{A}) \,
         + \\
         \sum_{(\mathbf{B}',\mathbf{B''})=\mathbf{A}} \bigl( &
         n_{III_1}(a,x,y;\mathbf{B}',\mathbf{B}'') +
         n_{III_2}(a,x,y;\mathbf{B}',\mathbf{B}'')
         \bigr) \,+ \\
         \sum_{\mathbf{C}'\#\mathbf{C''}=\mathbf{A}} \bigl( &
         n_{III_1}(a,x,y;\mathbf{C}',\mathbf{C}'') \, +
         n_{III_2}(a,x,y;\mathbf{C}',\mathbf{C}'') \bigr) = 0.
      \end{aligned}
   \end{equation}
   \medskip
   \begin{equation} \label{Eq:qm-ident-2}
      \begin{aligned}
         \sum_{|x'|=|x|-1} & n(x,x') n_{I'}(a,x',y;\mathbf{A}) \, +
         \sum_{|y'|=|y|+1}
         n_{I'}(a,x,y';\mathbf{A}) n(y',y) \, + \\
         \sum_{|a'|=|a|-1} & n(a,a') n_{I'}(a',x,y;\mathbf{A}) \, +
         \sum_{\substack{(\mathbf{A'},\mathbf{A''})=\mathbf{A} \\
             |x'|=|x|+\mu(\mathbf{A}')-1}} n(x,x'; \mathbf{A}')
         n_{I'}(a,x',y;\mathbf{A}'') \, + \\
         \sum_{\substack{(\mathbf{A'},\mathbf{A''})=\mathbf{A} \\
             |y'|=|y|-\mu(\mathbf{A}'')+1}} &
         n_{I'}(a,x,y';\mathbf{A}') n(y',y; \mathbf{A}'') \, +
         n_{II_1}(a,x,y;\mathbf{A}) +
         n_{II_2}(a,x,y;\mathbf{A}) + \\
         \sum_{(\mathbf{B}',\mathbf{B''})=\mathbf{A}} \bigl( &
         n_{III'_1}(a,x,y;\mathbf{B}',\mathbf{B}'') +
         n_{III'_2}(a,x,y;\mathbf{B}',\mathbf{B}'')
         \bigr)\,+ \\
         \sum_{\mathbf{C}'\#\mathbf{C''}=\mathbf{A}} \bigl( &
         n_{III'_1}(a,x,y;\mathbf{C}',\mathbf{C}'') \, +
         n_{III'_2}(a,x,y;\mathbf{C}',\mathbf{C}'') \bigr) = 0.
      \end{aligned}
   \end{equation}
\end{prop}

\begin{proof}[Proof of Proposition~\ref{P:qm-chain}]
   Let $a \in \textnormal{Crit}(h)$, $x \in \textnormal{Crit}(f)$.  As
   we work with $\mathbb{Z}_2$-coefficients we have to show that
   $d^f(a*x) + \partial^h(a) * x + a*d^f(x) = 0$. Write
   $$d^f(a*x) + \partial^h(a) * x + a*d^f(x) = \sum_{|y|-\mu(\lambda)
     = |a|+|x|-2n-1} m_{y,\lambda} y t^{\mubar(\lambda)}, \quad
   m_{y,\lambda} \in \mathbb{Z}_2.$$
   Fix $y,\lambda$ with
   $|y|-\mu(\lambda) = |a|+|x|-2n-1$. We will show below that the
   coefficient $m_{y,\lambda}$ of $y t^{\mubar(\lambda)}$ in this sum
   is $0$.  Note that for $\lambda=0$ this follows from standard
   arguments from Morse theory (see e.g.~\cite{Sc:Equivalences}).
   Therefore we assume from now on that $\lambda \neq 0$.

   Take the sum of identities~\eqref{Eq:qm-ident-1} and
   ~\eqref{Eq:qm-ident-2}, then sum up the result over all possible
   vectors $\mathbf{A}=(A_1, \ldots, A_l)$ (of all possible lengths
   $l$) with $\sum A_i = \lambda$.

   First note that the summands $n_{II_1}(a,x,y;\mathbf{A}) +
   n_{II_2}(a,x,y;\mathbf{A})$ being present in both
   identities~\eqref{Eq:qm-ident-1} and ~\eqref{Eq:qm-ident-2} cancel
   out. Next, note that when summing over all possible $\mathbf{A}$'s
   the summands of the type
   $n_{III_1}(a,x,y;\mathbf{B}',\mathbf{B}'')$ and
   $n_{III_1}(a,x,y;\mathbf{C}',\mathbf{C}'')$ are in $1-1$
   correspondence: whenever the first one appears for $(\mathbf{B}',
   \mathbf{B}'')=\mathbf{A}$ the second one appears when summing over
   $\mathbf{\widetilde{A}} = \mathbf{B}' \# \mathbf{B}''$ and vice
   versa. The same holds for the summands of the type $n_{III_2}$,
   $n_{III'_1}$, $n_{III''_2}$. Thus after summing over all
   $\mathbf{A}$'s these summands cancel out and we obtain:
   \begin{equation} \label{Eq:qm-sum}
      \begin{aligned}
         \sum_{\mathbf{A}, \sum A_i = \lambda} \Bigl(\,
         \sum_{|x'|=|x|-1} & n(x,x') n(a,x',y;\mathbf{A}) \, +
         \sum_{|y'|=|y|+1}
         n(a,x,y';\mathbf{A}) n(y',y) \, + \\
         \sum_{|a'|=|a|-1} & n(a,a') n(a',x,y;\mathbf{A}) \, +
         \sum_{\substack{(\mathbf{A'},\mathbf{A''})=\mathbf{A} \\
             |x'|=|x|+\mu(\mathbf{A}')-1}} n(x,x'; \mathbf{A}')
         n(a,x',y;\mathbf{A}'') \, + \\
         \sum_{\substack{(\mathbf{A'},\mathbf{A''})=\mathbf{A} \\
             |y'|=|y|-\mu(\mathbf{A}'')+1}} & n(a,x,y';\mathbf{A}')
         n(y',y; \mathbf{A}'') \, \Bigr) = 0.
      \end{aligned}
   \end{equation}
   Note that the right-hand side of identity~\eqref{Eq:qm-sum} is
   exactly the coefficient $m_{y, \lambda}$ of the term $y
   t^{\mubar(\lambda)}$ in
   $$d^f(a*x) + \partial^h(a) * x + a*d^f(x).$$
\end{proof}

\subsubsection{Proof of Proposition~\ref{P:qm-idents}}
\label{Sb:prf-qm-idents-N3}
We will give two different proofs for identity~\eqref{Eq:qm-ident-1}.
The first proof is based on a natural compactification of the
$1$-dimensional manifolds $\mathcal{P}_I(a,x,y;\mathbf{A},J)$ when
$|a|+|x|-|y|+\mu(\mathbf{A})-2n=1$. This approach is quite
straightforward but has the drawback that it proves
identity~\eqref{Eq:qm-ident-1} only under the assumption that $N_L
\geq 3$.  The reason for this restriction comes from transversality
issues.  (As for identity~\eqref{Eq:qm-ident-2}, this proof works well
for every $N_L \geq 2$). In Sections~\ref{Sb:prf-qm-idents-2}
--~\ref{Sb:conclusion-prf-qm-ident-all-N} we will give an alternative
proof for identity~\eqref{Eq:qm-ident-1} that works for every $N_L
\geq 2$. This approach is based on {\em Hamiltonian perturbations}.
Although the second proof is more general we found it worth presenting
the first proof too since it is more geometric and is closer in spirit
to the theory of $J$-holomorphic disks.

We start with the first proof. An important ingredient in this proof
is the next proposition.
\begin{prop} \label{P:qm-1-dim}
   There exists a second category subset
   $\mathcal{J}_{\textnormal{reg}} \subset \mathcal{J}(M,\omega)$ such
   that for every $J \in \mathcal{J}_{\textnormal{reg}}$, $a \in
   \textnormal{Crit}(h)$, $x, y \in \textnormal{Crit}(f)$ and every
   vector $\mathbf{A}$ of non-zero classes with
   $|a|+|x|-|y|+\mu(\mathbf{A})-2n=1$ the following holds:
   \begin{enumerate}
     \item If $N_L \geq 3$ then the space $\mathcal{P}_I(a,x,y; J,
      \mathbf{A})$ is a smooth $1$-dimensional manifold. Moreover
      every $(u_1, \ldots, u_l)$ in this space consists of simple and
      absolutely distinct disks.
     \item If $N_L \geq 2$ then the space $\mathcal{P}_{I'}(a,x,y; J,
      \mathbf{A})$ is a smooth $1$-dimensional manifold. Moreover
      every $(u_1, \ldots, u_l)$ in this space consists of simple and
      absolutely distinct disks.
   \end{enumerate}
\end{prop}
We defer the proof of this proposition to
Section~\ref{Sb:qm-prfs-props}. Note that statement (1) of
Proposition~\ref{P:qm-1-dim} assumes $N_L \geq 3$.

\begin{proof}[Proof of Proposition~\ref{P:qm-idents}]

   We start with identity~\eqref{Eq:qm-ident-1}. For brevity we omit
   the $J$ from the notation of the moduli spaces $\mathcal{P}_I$ etc.

   By Proposition~\ref{P:qm-1-dim} the space $\mathcal{P}_I(a,x,y;
   \mathbf{A})$ is a $1$-dimensional manifold.  We claim that it
   admits a compactification, in the Gromov topology, into a {\em
     compact $1$-dimensional manifold with boundary}
   $\overline{\mathcal{P}_I(a,x,y; \mathbf{A})}$ whose boundary points
   consists of the following disjoint union:
   \begin{equation} \label{Eq:qm-moduli}
      \begin{aligned}
         & \partial \overline{\mathcal{P}_I(a,x,y; \mathbf{A})} =
         \overline{\mathcal{P}_I(a,x,y; \mathbf{A})} \setminus
         \mathcal{P}_I(a,x,y; \mathbf{A}) = \\
         & \quad \Bigl( \bigcup_{|x'|=|x|-1} \mathcal{P}(x',y)
         \times \mathcal{P}_I(a,x',y;\mathbf{A}) \Bigr) \,\coprod\,
         \Bigl( \bigcup_{|y'|=|y|+1}
         \mathcal{P}_I(a,x,y';\mathbf{A})
         \times \mathcal{P}(y',y)\Bigr) \, \coprod \\
         & \quad \Bigl( \bigcup_{|a'|=|a|-1} \mathcal{P}(a,a') \times
         \mathcal{P}_I(a',x,y;\mathbf{A}) \Bigr) \, \coprod \\
         & \quad \Bigl( \bigcup_{\substack{(\mathbf{A'},\mathbf{A''})=
             \mathbf{A} \\
             |x'|=|x|+\mu(\mathbf{A}')-1}} \mathcal{P}(x,x';
         \mathbf{A}') \times \mathcal{P}_I(a,x',y;\mathbf{A}'') \Bigr)
         \, \coprod \\
         & \quad \Bigl( \bigcup_{\substack{(\mathbf{A'},\mathbf{A''})=
             \mathbf{A} \\
             |y'|=|y|-\mu(\mathbf{A}'')+1}}
         \mathcal{P}_I(a,x,y';\mathbf{A}')\times \mathcal{P}(y',y;
         \mathbf{A}'') \Bigr) \, \coprod \\
         & \quad \mathcal{P}_{II_1}(a,x,y;\mathbf{A}) \, \coprod \,
         \mathcal{P}_{II_2}(a,x,y;\mathbf{A}) \, \coprod \\
         & \quad \bigcup_{(\mathbf{B}',\mathbf{B''})=\mathbf{A}}
         \Bigl( \mathcal{P}_{III_1}(a,x,y;\mathbf{B}',\mathbf{B}'') \,
         \coprod \,
         \mathcal{P}_{III_2}(a,x,y;\mathbf{B}',\mathbf{B}'') \Bigr) \,
         \coprod \\
         & \quad \bigcup_{\mathbf{C}' \# \mathbf{C''}=\mathbf{A}}
         \Bigl( \mathcal{P}_{III_1}(a,x,y;\mathbf{C}',\mathbf{C}'') \,
         \coprod \,
         \mathcal{P}_{III_2}(a,x,y;\mathbf{C}',\mathbf{C}'') \Bigr).
         \end{aligned}
   \end{equation}
   Note that this immediately proves identity~\eqref{Eq:qm-ident-1}
   since the number of points in $\partial
   \overline{\mathcal{P}_I(a,x,y; \mathbf{A})}$ is on the one hand
   exactly the left-hand side of the identity~\eqref{Eq:qm-ident-1}
   while on the other hand the number of boundary points of a compact
   $1$-dimensional manifold is always even.

   We turn to proving~\eqref{Eq:qm-moduli}. We claim that the boundary
   $\partial \overline{\mathcal{P}_I(a,x,y; \mathbf{A})}$ is included
   in the union that appears in~\eqref{Eq:qm-moduli}. To see this let
   $\mathbf{u}^{(n)} = (u_1^{(n)}, \ldots, u_l^{(n)})$ be a sequence
   of elements in $\mathcal{P}_I(a,x,y; \mathbf{A})$. By compactness
   there exists a subsequence, still denoted by $\mathbf{u}^{(n)}$
   that converges in the Gromov topology. Assume by contradiction that
   the limit of this subsequence is not an element of the spaces in
   the right-hand side of~\eqref{Eq:qm-moduli}. Without loss of
   generality we may assume that the whole sequence $\mathbf{u}^{(n)}$
   lies in the space $\mathcal{P}_I(a,x,y;(\mathbf{A},k),J)$ for the
   same $1 \leq k \leq l$. (See~\eqref{Eq:qm-PI} for the definition of
   this space). As $\mathbf{u}^{(n)}$ does not converge to an element
   of the right-hand side of~\eqref{Eq:qm-moduli} one of the following
   occurs for the limit of $\mathbf{u}^{(n)}$:
\begin{enumerate}
  \item More than one breaking occurs at the gradient trajectories
   involved in the definition of $\mathcal{P}_I(a,x,y;\mathbf{A})$.
  \item More than two gradient trajectories connecting two consecutive
   disks in $\mathbf{u}^{(n)}$ shrink to a point.
  \item Bubbling of a $J$-holomorphic sphere occurs at a point on one
   of the disks in $\mathbf{u}^{(n)}$.
  \item Bubbling of a holomorphic disk occurs at a point $p \neq -1,1$
   on the boundary of one of the disks $u_j^{(n)}$. See
   figures~\ref{f:qm-bubbles-1},~\ref{f:qm-bubbles-2}.
  \item Bubbling of a holomorphic disk occurs at a point $p=-1$ or
   $p=1$ but the marked point $p$ corresponds in the limit to the
   attaching point of the two bubble disks.
  \item $\mathbf{u}^{(n)}$ converges in the Gromov topology to a
   configuration of disks with at least {\em two} disks bubbling.
  \item A combination of the above.
\end{enumerate}
Possibilities (1),(2),(3) and (6) can be ruled out by a dimension
count. Indeed, using the techniques of section~\ref{S:transversality}
(See also Section~\ref{Sb:qm-prfs-props} below) it follows that for
generic $J$ each of the configurations in (1),(2),(3) and (6) must
have negative dimension, hence cannot occur.  We turn to possibility
(4). Assume first that $j \neq k$, i.e. the bubbling occurs in one of
the disks that is not connected to the gradient trajectory going from
$a \in \textnormal{Crit}(h)$ (See figure~\ref{f:qm-bubbles-1}.) It
follows that one of the components of the limit is an element of
$\mathcal{P}_I(a,x,y;\mathbf{A}')$ for some $\mathbf{A}'$ with
$\mu(\mathbf{A}') \leq \mu(\mathbf{A})-N_L$.  Again, the
techniques of section~\ref{S:transversality} and a dimension
computation show that this configuration has negative dimension, hence
impossible for generic $J$. Finally, assume that $j=k$. Recall that
$u_k^{(n)}(0) \in W_a^u$ for every $n$. Note that the points $-1, 0, 1
\in D$ lie on the same hyperbolic geodesic. Since the conformal
structure of the disks is preserved when passing to limits in the
Gromov topology, the marked points $-1,0,1$ on $u_k^{(n)}$ must
converge to the same bubble disk. (Thus the right-hand side of
figure~\ref{f:qm-bubbles-2} is impossible.)  It follows that one of
the components of the limit of $\mathbf{u}^{(n)}$ is an element of
$\mathcal{P}_I(a,x,y;\mathbf{A}')$ for some $\mathbf{A}'$ with
$\mu(\mathbf{A}') \leq \mu(\mathbf{A})-N_L$ and we arrive at
contradiction as before (See the left-hand side of
picture~\ref{f:qm-bubbles-2}.) This rules out possibility (4).
Possibility (5) is ruled out in a similar way.  Finally, possibility
(6) is ruled out by a combination of the above arguments.

\begin{figure}[htbp]
   \begin{center}
      \epsfig{file=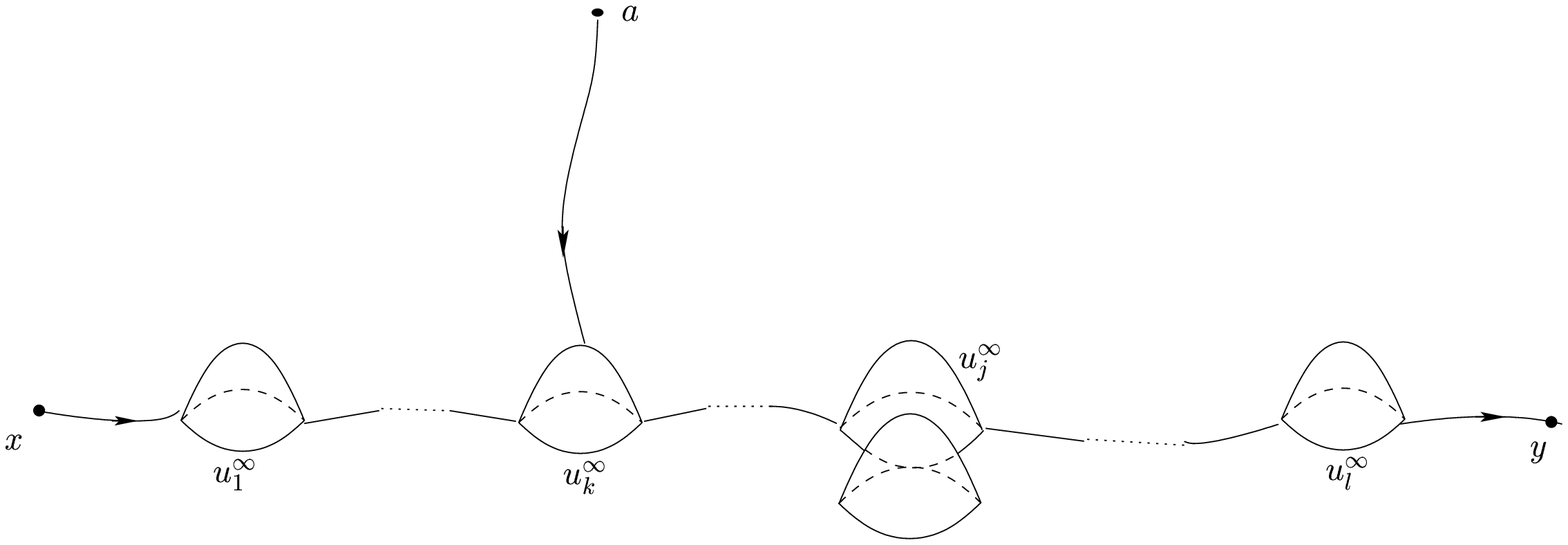, width=0.6\linewidth}
   \end{center}
   \caption{Limit $\mathbf{u}^{\infty}$ of $\mathbf{u}^{(n)}$ in
     possibility (4) with $j \neq k$.}
   \label{f:qm-bubbles-1}
\end{figure}

\begin{figure}[htbp]
   \begin{center}
      \epsfig{file=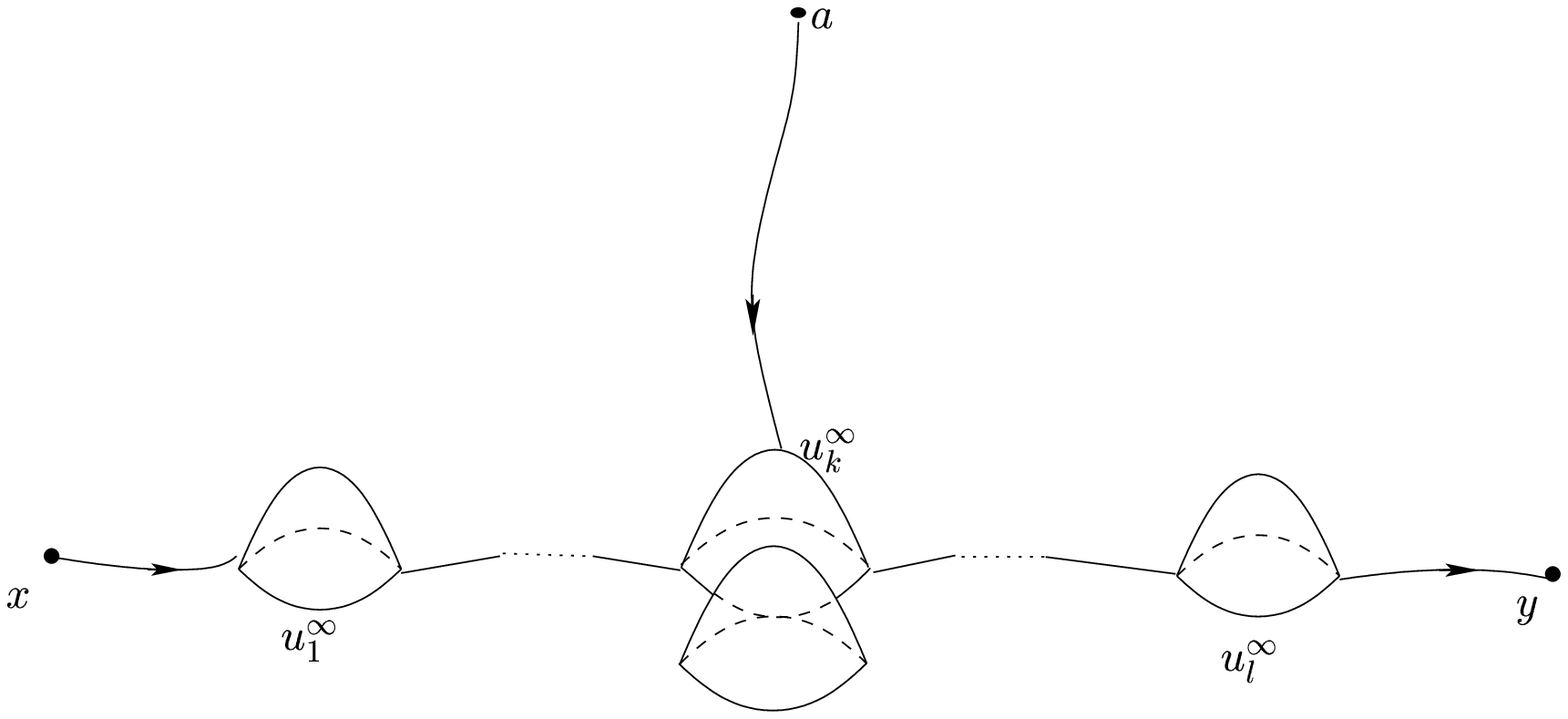, width=0.45 \linewidth} \quad
      \epsfig{file=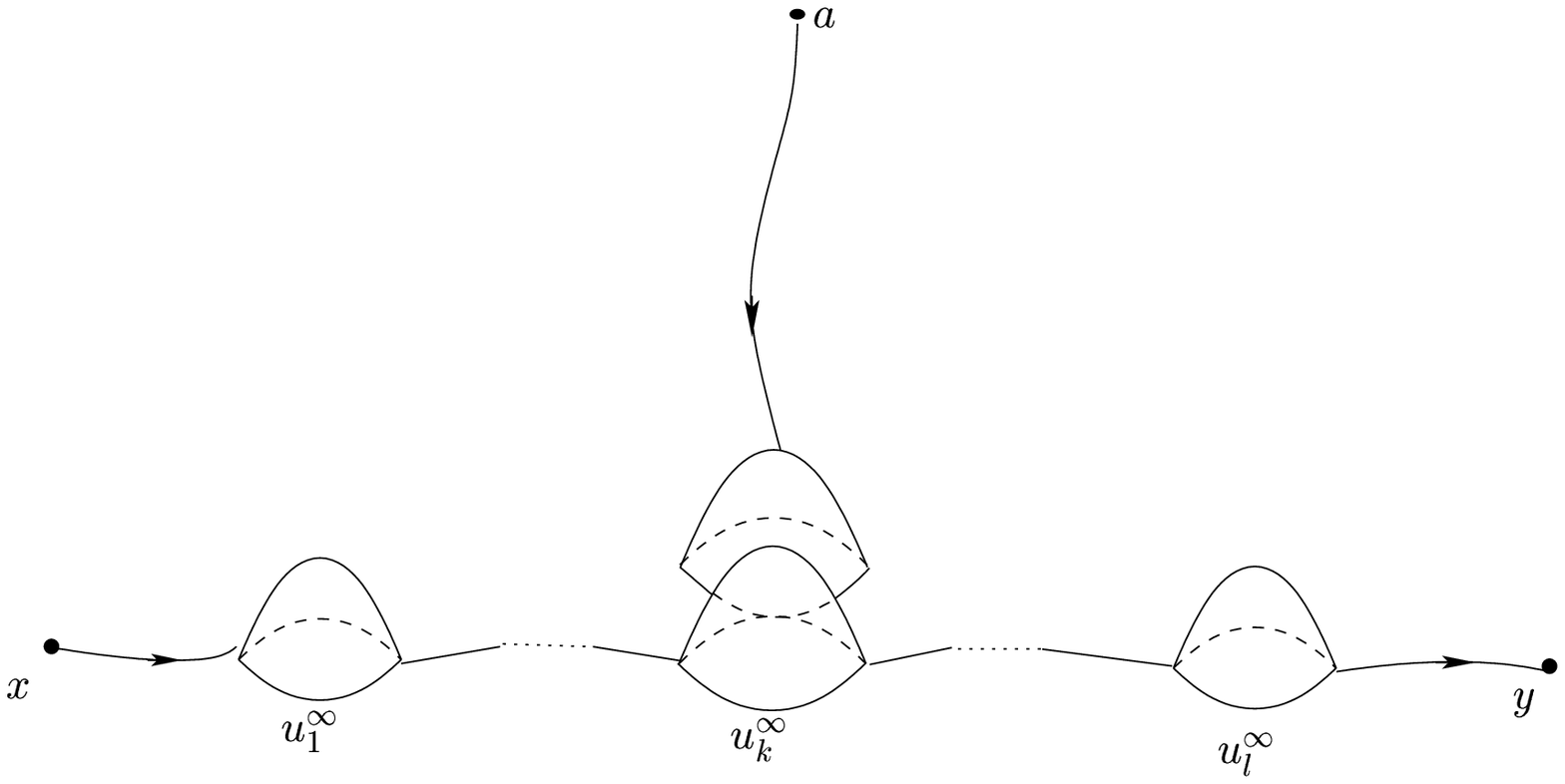, width=0.45\linewidth}
   \end{center}
   \caption{Possibility (4) for $j=k$.}
   \label{f:qm-bubbles-2}
\end{figure}

Finally we have to prove the following statement:
\begin{statement} \label{St:qm-boundary}
   Every point in the spaces from the right-hand side
   of~\eqref{Eq:qm-moduli} is indeed a limit of a sequence of elements
   from $\mathcal{P}_I(a,x,y;\mathbf{A})$. Moreover every such element
   corresponds to {\em exactly} one end of the $1$-dimensional
   manifold $\mathcal{P}_I(a,x,y;\mathbf{A})$.
\end{statement}

We start with the spaces appearing in the first 6 lines of the union
on the right-hand side of~\eqref{Eq:qm-moduli}.  In each of these 6
cases, Statement~\ref{St:qm-boundary} follows from general
transversality arguments and gluing theorems in the framework of Morse
theory.  (In each of these cases one has to add boundary components to
the submanifold $Q_{f,\rho_L}$ to include either a breaking of a
gradient trajectory or to add parts of $\textnormal{diag}(L \setminus
\textnormal{Crit})(f)$ to $Q_{f,\rho_L}$ and then make the suitable
evaluation maps transverse to this enlarged manifold with boundary).

Finally, it remains to prove Statement~\ref{St:qm-boundary} for
elements from $\mathcal{P}_{III_i}(a,x,y;\mathbf{C}',\mathbf{C}'')$,
$i=1,2$, where $\mathbf{C}' \# \mathbf{C''}=\mathbf{A}$. This follows
from the gluing procedure described in Section~\ref{Sbsb:glue-ex-2I}
(see Corollary~\ref{C:glue-pearls-2} there). Note that we need here
both the existence and uniqueness results of
Corollary~\ref{C:glue-pearls-2}.

We now turn to the proof of identity~\eqref{Eq:qm-ident-2}.  The proof
is very similar to the proof of the identity~\eqref{Eq:qm-ident-1},
but there are two slight changes.  First of all, we are using now
statement~(2) of Proposition~\ref{P:qm-1-dim} which assumes only $N_L
\geq 2$ and therefore the proof works well for every $N_L \geq 2$.
The second change is when $|a|=2n$. In this case the spaces
$\mathcal{P}_{I'}(a,x,y;(\mathbf{A},k))$, for $k=0$ and $k=l$ have
additional boundary points. Indeed, we may have a sequence
$(\mathbf{u}^{(n)},p_n) \in \mathcal{P}_{I'}(a,x,y; (\mathbf{A},0))$
with $p_n \to x$ as $n \to \infty$. Similarly we may have
$(\mathbf{u}^{(n)},p_n) \in \mathcal{P}_{I'}(a,x,y; (\mathbf{A},l))$
with $p_n \to y$ as $n \to \infty$. Therefore the boundary of
$\overline{\mathcal{P}_{I'}(a,x,y; \mathbf{A})}$ contains the
following additional points:
\begin{equation} \label{Eq:qm-additional}
   \Bigl( \mathcal{P}(x,y;\mathbf{A}) \times \{ x \} \Bigr)
   \, \coprod \, \Bigl( \mathcal{P}(x,y; \mathbf{A}) \times
   \{ y \} \Bigr).
\end{equation}
As before, standard transversality arguments show that all the points
in~\eqref{Eq:qm-additional} can indeed be realized as boundary points
of $\mathcal{P}_{I'}(a,x,y;\mathbf{A})$ and that each of them
corresponds to exactly one end of this $1$-dimensional manifold.

Note that since $|a|=2n$ we have $|x| - |y| + \mu(\mathbf{A})=1$
hence $\mathcal{P}(x,y; \mathbf{A})$ is a finite set. Finally note
that the number of points appearing in~\eqref{Eq:qm-additional} is
even (due to our assumption that $h$ has a single maximum) hence their
appearance will not change our identity.
\end{proof}

\subsubsection{Proof of Propositions~\ref{P:qm-axy},~\ref{P:qm-0-dim}
  and~\ref{P:qm-1-dim}} \label{Sb:qm-prfs-props}

\begin{proof}[Proof of Proposition~\ref{P:qm-axy}]
   Denote by $\mathcal{P}_I^{*,\textnormal{d}}(a,x,y;\mathbf{A},J)
   \subset \mathcal{P}_I(a,x,y;\mathbf{A},J)$ the subspace of
   sequences of disks that are simple and absolutely distinct.  By
   standard arguments, for generic $J$ this space is a smooth manifold
   of dimension $0$. We claim that
   $\mathcal{P}_I^{*,\textnormal{d}}(a,x,y;\mathbf{A},J) =
   \mathcal{P}_I(a,x,y;\mathbf{A},J)$.

   The main idea is that the dimension of a configuration of disks in
   $\mathcal{P}_I(a,x,y;\mathbf{A},J)$ that is not simple or not
   absolutely distinct is negative, hence such configurations cannot
   occur for generic $J$. Most of the arguments are very similar to
   those appearing in Section~\ref{S:transversality} therefore we will
   only give an outline of the proof explaining how to adjust the
   arguments from Section~\ref{S:transversality} to work in the
   present situation.

   As in Section~\ref{S:transversality} we separate the proof to the
   cases $n=\dim L \geq 3$ and $n=\dim L \leq 2$.  We start with
   $n\geq 3$. The proof in this case can be carried out in a similar
   way to the proof of Proposition~\ref{P:tr-1} from
   Section~\ref{S:transversality}, with Lemmas~\ref{L:uv}
   and~\ref{L:non-simple} being the main tools. The only adjustment
   needed is when applying Lemma~\ref{L:non-simple} to the $k$'th disk
   in $\mathbf{u} = (u_1, \ldots, u_l)$ where $\mathbf{u} \in
   \mathcal{P}_I(a,x,y;(\mathbf(A),k),J)$. Indeed, suppose (by
   contradiction) that $u_k$ is not simple. Then by
   Lemma~\ref{L:non-simple} we obtain a new disk $u'_k$ with the same
   image as $u$, with $u'_k(-1)=u_k(-1)$, $u'_k(1)=u_k(1)$ and
   $\mu[u'_k] \leq \mu([u_k])-N_L$, but we cannot assume that
   $u'_k(0) = u_k(0)$ anymore. Therefore the sequence of disks
   $\mathbf{u}'$ obtained from $\mathbf{u}$ by replacing $u_k$ by
   $u'_k$ is not necessarily an element of
   $\mathcal{P}_I(a,x,y;\mathbf{A}',J)$, where $\mathbf{A}'$ is the
   vector of classes of $\mathbf{u}'$. However,
   Lemma~\ref{L:non-simple} still implies that there exists a point $p
   \in \textnormal{Int\,}(D)$ such that $u'_k(p)=u_k(0) \in W_a^u$.
   Consider the space of disks $\mathbf{v}=(v_1, \ldots, v_l)$ in
   $\mathcal{P}(x,y; \mathbf{A}',J)$ with the additional condition
   that there exists a point $q \in \textnormal{Int\,}(D)$ with
   $v_k(q) \in W_a^u$. A simple computation shows that the dimension
   of this space is
   \begin{align} \label{Eq:dim-axy-1}
      & |a|+|x|-|y|-2n+\mu(\mathbf{A}') + 1 \leq \\
      & |a|+|x|-|y|-2n + \mu(\mathbf{A}) - N_L + 1 = -N_L+1 <0.
      \notag
   \end{align}
   On the other hand $\mathbf{u}'$ belongs to this space, a
   contradiction. It remains to show that the disks present in the
   configurations of $\mathcal{P}_{I}(a,x,y;\mathbf{A},J)$ are also
   absolutely distinct. Under the assumption $n\geq 3$ we may apply
   Lemma \ref{L:uv}. Assuming that two disks $u_{1}$ and $u_{2}$ in a
   configuration $u\in \mathcal{P}_{I}(a,x,y;\mathbf{A},J)$ are not
   absolutely distinct, this result allows us to eliminate the
   ``smallest" of the two disks - say $u_{1}$ - by $u_{2}$ thus
   reducing the total energy of the configuration. It is easy to see
   that this replacement can always be done in such a way that the
   resulting object $u'$ is still modeled on a tree.  However, it is
   also easy to notice that $u'$ might not be in any of the spaces
   $\mathcal{P}_{-}(a,x,y;\mathbf{A}',J)$ but only in a more general
   type of moduli space. These moduli spaces consist of configurations
   modeled on a tree with two entries and one exit - as before - so
   that one of the entering edges corresponds to a flow line of $h$,
   all the other edges are flow lines of $f$ and the vertex of valence
   three can be the end of the flow line of $h$ or, in contrast to the
   cases described before, it is also allowed to have its two entering
   edges correspond to flow lines of $f$. To achieve transversality,
   an important point is to notice that if the attaching points of the
   two entering edges happen to coincide in this last case, then the
   respective configuration can be rewritten as an element of the same
   type of moduli space only modeled on a tree in which one of these
   edges (the shortest one) is eliminated and a new vertex of valence
   three is introduced. This vertex corresponds to a ghost (or
   constant) disk to which one of the branches of the tree is attached
   by an edge of length $0$. Now we apply the proof before to first
   show that the disks present in these configurations are all simple
   for those moduli spaces whose virtual dimension is $0$ or $1$ and
   then we notice that the replacement argument before allows to show
   by recurrence on total energy that the disks are also absolutely
   distinct for a second category subset of $J$'s in the same range of
   dimensions. The condition $N_{L}\geq 3$ is used again in the
   ``replacement" argument as in this process the incidence condition
   for the flow line of $h$ with the center of the disk at its end
   might not be preserved.

   We now turn to the case $n=\dim L=2$.  Since $|a|+|x|-|y|-2n +
   \mu(\mathbf{A})=0$ we must have $\mu(\mathbf{A}) \leq 6$.
   Therefore the number of disks $l$ in an element of
   $\mathcal{P}_I(a,x,y; \mathbf{A},J)$ is $\leq 3$. Consider the case
   $l=1$ first. We have to show that for generic $J$ every
   $J$-holomorphic disk $u$ with $u(-1) \in W_x^u$, $u(1) \in W_y^s$,
   $u(0) \in W_a^u$ and with $|a|+|x|-|y|-2n+\mu(A)=0$ must be
   simple.

   Suppose by contradiction that $u$ is a $J$-holomorphic disk as
   above which is not simple. Applying Theorem~\ref{T:decomposition}
   we obtain domains $\mathfrak{D}_1, \ldots, \mathfrak{D}_r \subset
   D$ and $J$-holomorphic disks $v_1, \ldots, v_r$ through which $u$
   factors on each of the domains $\mathfrak{D}_i$. Denote by $m_i$
   the degree of the map $\pi_i:\overline{\mathfrak{D}}_i \to D$ for
   which ${u_i}|_{\overline{\mathfrak{D}}_i} = v_i \circ \pi_i$. As
   $N_L \geq 2$ and $\mu(A) \leq 6$ we have $r\leq 3$.

   \smallskip \noindent \textbf{Case 1. $l=1$, $r=3$.}  In this case
   we have three simple disks $v_1, v_2, v_3$ each of them with Maslov
   number $2$. In this case $\mu(A)=6$ and $m_1=m_2=m_3=1$. By
   assumption we have $|x|-|y| = -2-|a| \leq -2$. We may assume that
   $u(-1) \neq u(1)$ (for otherwise by omitting $u$ we obtain a
   negative gradient trajectory of $f$ going from $x$ to $y$ which is
   impossible).

   Assume first that the three disks $v_1, v_2, v_3$ are absolutely
   distinct. Note that there are several possibilities as to how the
   three points $-1, 0, 1$ are distributed among the domains
   $\mathfrak{D}_1, \mathfrak{D}_2, \mathfrak{D}_3$. In case the
   points $-1$, $1$ lie in the same domain, say $\mathfrak{D}_1$ we
   can omit the other two disks $v_2, v_3$ and after reparametrizing
   $v_1$ obtain an element of the space
   $\mathcal{P}^{*,\textnormal{d}}_I(x,y; B)$ where $\mu(B)=2$. Note
   that this is impossible since the dimension of this space is
   $|x|-|y|+\mu(B)-1 \leq -1$. Therefore we assume that $-1$, $1$
   lie in different domains, say $-1 \in \mathfrak{D}_1$, $1 \in
   \mathfrak{D}_2$. If the point $0$ lies in interior of
   $\mathfrak{D}_1$ then after reparametrizing $v_1$ we may assume
   that $v_1(-1) \in W_x^u$, $v_1(0) \in W_a^u$. On the other hand, a
   simple computation shows that the space of such disks (with
   $\mu=2$) has dimension $|x|+|a|-2$. But $|x|+|a|-2 = -4+|y| < 0$.
   A contradiction.  A similar argument shows that it is impossible
   for $0$ to lie on the boundary of $\mathfrak{D}_1$. The same
   arguments applied to $v_2$ show that it is impossible for $0$ to
   lie neither in $\textnormal{Int\,}\mathfrak{D}_2$ nor on the
   boundary of this domain.  We are left with the case $0 \in
   \textnormal{Int\,}\mathfrak{D}_3$. Note that the boundary of
   $\mathfrak{D}_3$ must have a non-trivial intersection with the
   boundary of at least one of the domains $\mathfrak{D}_1$,
   $\mathfrak{D}_2$. (In fact this intersection must be at least
   $1$-dimensional.) Without loss of generality assume that such an
   intersection occurs with the boundary of $\mathfrak{D}_1$.

   Reparametrizing $v_1, v_3$ we may assume that $v_1(-1) \in W_x^u$,
   $v_3(0) \in W_a^u$ and $v_1(1)=v_3(-1)$. Put $B'=[v_1]$,
   $B''=[v_3]$. It follows that the space
   \begin{align} \label{Eq:dim-v1-v3}
      \Bigl\{ & (w',w'') \in \mathcal{M}(B', J)/G_{-1,1} \times
      \mathcal{M}(B'', J) \, \Big| w'(-1) \in W_x^u, w''(0) \in W_a^u,
      \\
      & w'(1)=w''(-1), (w',w'') \textnormal{ are simple and absolutely
        distinct} \Bigr\} \notag
   \end{align}
   is not empty. On the other hand a simple computation shows that the
   dimension of this space is $|x|+|a|-1 = -3+|y| < 0$, a
   contradiction.
   \begin{rem} \label{R:v1-v3-1-dim}
      It is worth noting that the pair $(v_1, v_3)$ in fact gives rise
      to a $1$-parametric family of elements in the space
      in~\eqref{Eq:dim-v1-v3} since the boundary of at least one of
      the domains $\mathfrak{D}_1$, $\mathfrak{D}_2$ must intersect
      the boundary of $\mathfrak{D}_3$ along a $1$-dimensional piece.
      Therefore there is a $1$-parametric family of reparametrizations
      of $(v_1, v_3)$ which lies in the space in~\eqref{Eq:dim-v1-v3}.
      Here we have not used this fact as we got contradiction by
      showing that the dimension of this space is negative. However,
      this observation can used to rule out other configurations
      below.
   \end{rem}

   This concludes the proof of case 1 when $v_1, v_2, v_3$ are
   absolutely distinct. It remains to deal with the case when these
   three disks are not absolutely distinct. Again, there are several
   subcases to be considered:
   \begin{itemize}
     \item $v_1(D) \subset v_2(D) \cup v_3(D)$ but $v_2, v_3$ are
      absolutely distinct.
     \item $v_1(D), v_2(D) \subset v_3(D)$.
   \end{itemize}
   In each of these two subcases one has to deal with all possible
   distribution of the three points $-1, 0, 1$. Ruling out these
   possibilities is made by similar arguments to the above and we will
   not give the details here.

   \smallskip \noindent \textbf{Case 2. $l=1$, $r=2$.}  Here we have
   two simple disks $v_1, v_2$ and two multiplicities $m_1, m_2$. As
   $\mu(\mathbf{A})\leq 6$ we may assume that one of the following
   holds:
   \begin{itemize}
     \item $\mu([v_1])=2$, $m_1=1$, $\mu_([v_2])=2$, $m_2=1$;
      $\mu(\mathbf{A})=4$, $|a|+|x|-|y|=0$.
     \item $\mu([v_1])=2$, $m_1=1$, $\mu_([v_2])=2$, $m_2=2$;
      $\mu(\mathbf{A})=6$, $|a|+|x|-|y|=-2$.)
     \item $\mu([v_1])=2$, $m_1=1$, $\mu([v_2])=4$, $m_2=1$;
      $\mu(\mathbf{A})=6$, $|a|+|x|-|y|=-2$.)
   \end{itemize}
   Again, by arguments similar to those used above one can extract in
   each of these three cases a configuration which has negative
   dimension, hence deduce that it cannot occur for generic $J$.

   \smallskip \noindent \textbf{Case 3. $l=1$, $r=1$.}  In this case
   the disk $u$ is multiply covered, say with multiplicity $m \geq 2$.
   We denote by $v$ its reduction. By assumption
   $|a|+|x|-|y|=4-m\mu([v])\leq 0$. First note that $u(-1) \neq
   u(1)$ for otherwise we would get a negative gradient trajectory of
   $f$ going from $x$ to $y$ which is impossible since $|x| \leq |y|$.
   Let $p', p'' \in \partial{D}$ be points such that $v(p')=u(-1)$,
   $v(p'')=u(1)$. As $p' \neq p''$ we can reparametrize $v$ so that
   $v(-1) \in W_x^u$, $v(1) \in W_y^s$. Let $q \in D$ a point such
   that $v(q) \in W_a^u$ (note that we cannot assume anymore that
   $q=0$). Put $B=[v]$. Assume first that $q \in \textnormal{Int\,}D$.
   Then $(v,q)$ belongs to the space
   $$\{ (w,p) \in \bigr(\mathcal{M}^*(B,J) \times \textnormal{Int\,}D
   \bigr) /G_{-1,1} \mid w(-1) \in W_x^u, w(p) \in W_a^u, w(1) \in
   W_y^s \}.$$
   On the other hand a simple computation shows that the
   dimension of this space is: $$-3+|x|+|a|-|y| + \mu(B) = 1-
   (m-1)\mu(B) < 0,$$
   where the last inequality follows from $m\geq
   2$, $\mu(B) \geq 2$. A contradiction. This concludes the proof for
   the case $l=1$.

   The remaining cases are $l=2,3$. The proofs for these cases are
   again similar to the preceding ones hence we omit the details. Let
   us only list the variety of configurations needed to be ruled out.
   In the case $l=2$ we have two disks $u_1, u_2$ and there are four
   possibilities to be ruled out:
   \begin{enumerate}
     \item $u_1$ is simple with $\mu([u_1])=2$. $\mu([u_2])=4$ and
      $u_2$ is double covered.
     \item $u_1$ is simple with $\mu([u_1])=2$. $\mu([u_2])=4$ and
      $u_2$ can be decomposed into two simple disks $v'_2, v''_2$ each
      with $\mu=2$.
     \item $u_1, u_2$ are both simple with $\mu([u_1])=2$,
      $\mu([u_2])=4$, but $u_1(D) \subset u_2(D)$.
     \item $u_1, u_2$ are both simple with $\mu(u_1)=\mu(u_2)=2$
      but $u_1, u_2$ are not absolutely distinct, i.e. $u_1(D) \subset
      u_2(D)$ or $u_2(D) \subset u_1(D)$.
   \end{enumerate}

   Finally, in case $l=3$ we have three disks with
   $\mu([u_1])=\mu([u_2])=\mu([u_3])=2$ hence they are all simple.  In
   this case one has to rule out the possibility that they are not
   absolutely distinct. This concludes the proof for the case $n=\dim
   L=2$.

   Finally, in case $n=\dim L=1$ there is nothing to prove.  Indeed,
   the assumption $|a|+|x|-|y|-2n+\mu(\mathbf{A})=0$ implies that we
   have $\mu(\mathbf{A})\leq 3$ hence every element $\mathbf{u} \in
   \mathcal{P}_I(a,x,y;\mathbf{A},J)$ consists of exactly one simple
   disk.

   Next, we turn to compactness of the space
   $\mathcal{P}_I(a,x,y;\mathbf{A},J)$ (under the assumption
   $|a|+|x|-|y|-2n+\mu(\mathbf{A})=0$). The main argument is again a
   combination of Gromov compactness theorem with a dimension count.
   One lists all the possible limits, in the Gromov topology, of
   sequences $\mathbf{u}^{(n)} \in \mathcal{P}_I(a,x,y;\mathbf{A},J)$
   as in~\eqref{Eq:qm-moduli}.  (As in the proof of
   Proposition~\ref{P:qm-idents} one first has to show that none of
   the possibilities (1)-(7) listed there can appear.)  Then using the
   same methods as above one shows that under the assumption
   $|a|+|x|-|y|-2n+\mu(\mathbf{A})=0$ all the elements in the spaces
   appearing in the spaces in~\eqref{Eq:qm-moduli} are simple and
   absolutely distinct. Finally, a dimension computation shows that
   when $|a|+|x|-|y|-2n+\mu(\mathbf{A})=0$ all these spaces have
   negative dimension hence empty. This completes the proof of the
   statement on compactness, hence the proof
   Proposition~\ref{P:qm-axy} for the space
   $\mathcal{P}_I(a,x,y;\mathbf{A},J)$.

   The proof for the space $\mathcal{P}_{I'}(a,x,y;\mathbf{A},J)$ is
   similar (and in fact simpler since all disks in this space have
   only two marked points).
\end{proof}

The proof of Proposition~\ref{P:qm-0-dim} is similar to the proof of
Proposition~\ref{P:qm-axy}, hence we omit it.

\begin{proof}[Outline of proof of Proposition~\ref{P:qm-1-dim}]
   The proof goes along the same lines as the proof of
   Proposition~\ref{P:qm-axy} with one main difference for the space
   $\mathcal{P}_I(a,x,y; \mathbf{A},J)$. While in
   Proposition~\ref{P:qm-axy} we had the assumption
   $|a|+|x|-|y|+\mu(\mathbf{A})-2n=0$ now we have
   $|a|+|x|-|y|+\mu(\mathbf{A})-2n=1$. This is the reason that we
   have to assume that $N_L\geq 3$ rather than just $N_L\geq 2$. The
   point is that when $N_L=2$ non-simple disks might appear in
   $\mathcal{P}_I(a,x,y;\mathbf{A},J)$ but not for $N_L \geq 3$.  For
   example, the computation of dimension in~\eqref{Eq:dim-axy-1} gives
   us in the present situation:
   \begin{align*}
      & |a|+|x|-|y|-2n+\mu(\mathbf{A}') + 1 \leq \\
      & |a|+|x|-|y|-2n + \mu(\mathbf{A}) - N_L + 1 = -N_L+2.
   \end{align*}
   In order for this number to be negative we need $N_L \geq 3$.

   The rest of the proof continues in an analogous way to the proof of
   Proposition~\ref{P:qm-axy}.
\end{proof}

\subsubsection{A more general proof for identity~\eqref{Eq:qm-ident-1}
  of Proposition~\ref{P:qm-idents}} \label{Sb:prf-qm-idents-2} Below
is an alternative proof of identity~\ref{P:qm-idents} which works
under the more general assumption $N_L \geq 2$.

This approach consists of the following 3 steps.

{\bf Step 1.} We perturb the Cauchy-Riemann equation to a
non-homogeneous equation by adding a perturbation term $H$ (generated
by Hamiltonian functions). The perturbation procedure is applied to
elements of the spaces $\mathcal{P}_I, \mathcal{P}_{II_i},
\mathcal{P}_{III_i}$ in the following way. In each element of these
spaces, which is a chain of disks, we perturb only the single disk in
the chain which has 3 marked points. All the other disks in the chain
are left unperturbed hence remain $J$-holomorphic. The result of this
procedure is ``perturbations'' $\mathcal{P}_I(-,H),
\mathcal{P}_{II_i}(-,H), \mathcal{P}_{III_i}(-,H)$ of the original
spaces $\mathcal{P}_I, \mathcal{P}_{II_i}, \mathcal{P}_{III_i}$.

{\bf Step 2.} By counting the number of points in these spaces we
obtain ``perturbed'' versions $n_I(-,H), n_{II_i}(-,H),
n_{III_i}(-,H)$ of the numbers $n_I(-), n_{II_i}(-), n_{III_i}(-)$.
The advantage of the above perturbations is that now it is easy to
achieve transversality for the perturbed disk (without any simplicity
requirements). Therefore, by similar arguments to those from the older
proof of identity~\eqref{Eq:qm-ident-1} we conclude that
identity~\eqref{Eq:qm-ident-1} holds (for $N_L \geq 2$) with the $n_I,
n_{II_i}, n_{III_i}$'s replaced by their perturbed analogues.

{\bf Step 3.} Finally, we prove that for small enough perturbations
$H$ the numbers $n_I(-,H)$, $n_{II_i}(-,H)$, $n_{III_i}(-,H)$ coincide
with the numbers $n_I(-), n_{II_i}(-), n_{III_i}(-)$ hence
identity~\eqref{Eq:qm-ident-1} in fact holds in its original form.

\medskip

We now turn to the implementation of the proof.

\subsubsection{Hamiltonian perturbations} \label{Sb:hampert}
Here we briefly summarize the necessary ingredients from the theory of
Hamiltonian perturbations that will be needed later. The material of
this section is mostly based on~\cite{Ak-Sa:Loops, McD-Sa:Jhol-2}.  We
refer the reader to these texts for the foundations and more details.

Let $H \in \Omega^1(D, C_0^{\infty}(M))$ be a $1$-form on the disk $D$
with values in $C_0^{\infty}(M)$. If $(s,t)$ are coordinates on $D$ we
can write $H$ as $$H=F ds + G dt$$
for some compactly supported smooth
functions $F,G: D \times M \to \mathbb{R}$. Denote by $pr_M:D \times M
\to M$, $pr_{D}:D \times M \to D$ the projections. We write
$F_{s,t}(x) = F(s,t,x)$ and $G_{s,t}(x)=G(s,t,x)$. Define the
following $2$-form on $D \times M$:
$$\widetilde{\omega}_H = pr_M^* \omega - dH = pr_M^{*}\omega - d'F
\wedge ds - d'G \wedge dt + (\partial_t F - \partial_s G) ds \wedge
dt,$$
where $d'$ stands for exterior derivative in the $M$-direction.
Henceforth we will work with $H$'s that have the following additional
property:
\begin{equation} \label{Eq:Lag-H}
   \widetilde{\omega}_H|_{T(\partial D \times L)}=0.
\end{equation}
We denote by $\mathcal{H}$ the space of all $H$
satisfying~\eqref{Eq:Lag-H}.

Note that the form $\widetilde{\omega}_H$ may be degenerate however
for $\kappa>0$ large enough the form
\begin{equation} \label{Eq:om-hk}
   \widetilde{\omega}_{H,\kappa} = \widetilde{\omega}_H +
   pr_D^*(\kappa \, ds \wedge dt)
\end{equation}
is symplectic and $\partial D \times L$ is still Lagrangian with
respect to it. How large should $\kappa$ be taken for this purpose is
determined by the {\em curvature} function $R_H:D \times M \to
\mathbb{R}$:
$$R_H(s,t,x)=\partial_s G - \partial_t F + \{ F_{s,t}, G_{s,t} \},$$
where $\{ \cdot, \cdot \}$ is the Poisson bracket on $(M,\omega)$.  A
simple calculation (see~\cite{McD-Sa:Jhol-2, Ak-Sa:Loops}) shows that
if $\kappa > R_H(s,t,x)$ for every $(s,t,x) \in D \times M$ then
$\widetilde{\omega}_{H,\kappa}$ is symplectic.

Let $X_{F_{s,t}}$ be the Hamiltonian vector field corresponding to the
function $F_{s,t}$. Put $X_F(s,t,x)=X_{F_{s,t}}(x)$, $(s,t,x) \in D
\times M$. Similarly we have $G_{s,t}$, $X_{G_{s,t}}$ and
$X_G(s,t,x)$. Given an $\omega$-compatible almost complex structure
$J$ on $M$, consider the following elliptic boundary value problem:
\begin{equation} \label{Eq:delbar-pert}
   \begin{cases}
      u:(D, \partial D) \to (M,L) \\
      \partial_s u + J(u)\partial_t u = -X_F(s,t,u)-J(u)X_G(s,t,u)
   \end{cases}
\end{equation}
Solutions $u$ of this equation will be called {\em $(J,H)$-holomorphic
  disks}. We denote by $\mathcal{M}(A,J,H)$ the space of
$(J,H)$-holomorphic disks $u$ with $u_*([D])=A \in
H_2(M,L;\mathbb{Z})$.

As noted by Gromov~\cite{Gr:phol} solutions of~\eqref{Eq:delbar-pert}
are in 1-1 correspondence with $\widetilde{J}$-holomorphic disks in $D
\times M$ for some $\widetilde{J}$. To describe this correspondence
define the following endomorphism $\widetilde{J}_H:T(D \times M) \to
T(D \times M)$:
\begin{equation} \label{Eq:JH}
   \widetilde{J}_H =
   \left(
      \begin{matrix}
         0 & -1 & 0 \\
         1 &  0 & 0 \\
         JX_F - X_G & JX_G + X_F & J
      \end{matrix}
      \right).
\end{equation}

\begin{prop}[See~\cite{McD-Sa:Jhol-2}] \label{P:JH}
   \begin{enumerate}
     \item $\widetilde{J}_H$ is an almost complex structure on $D
      \times M$. Moreover, if $\kappa > R_H(s,t,x)$ for every $(s,t,x)
      \in D \times M$ then $\widetilde{J}_H$ is compatible with
      $\widetilde{\omega}_{H,\kappa}$.
     \item $u:(D, \partial D) \to (M,L)$ is a solution
      of~\eqref{Eq:delbar-pert} if and only if its graph
      $\widetilde{u}:(D,\partial D) \to (D \times M, \partial D \times
      L)$ defined by $\widetilde{u}(z)=(z,u(z))$ is
      $\widetilde{J}_H$-holomorphic.
     \item There exists a subset $\mathcal{H}_{\textnormal{reg}}
      \subset \mathcal{H}$ of second category such that for every $H
      \in \mathcal{H}$ and every $A \in H_2(M,L;\mathbb{Z})$ the space
      $\mathcal{M}(A,J,H)$ is either empty or a smooth manifold of
      dimension $n+\mu(A)$.
   \end{enumerate}
\end{prop}

\begin{rem}
   \begin{enumerate}
     \item In contrast to the case of genuine $J$-holomorphic disks,
      $(J,H)$-holomorphic disks in the class $0 \in
      H_2(M,L;\mathbb{Z})$ are not constant for general $H$. Thus we
      may have a (non-constant) $(J,H)$-holomorphic disk $u$ with
      $\omega([u])=0$.
     \item In contrast to the case of genuinely $J$-holomorphic disks
      we do not have to require the $(J,H)$-holomorphic disks to be
      simple in order for $\mathcal{M}(A,J,H)$ to be a smooth
      manifold.  See~\cite{McD-Sa:Jhol-2} for more details.
   \end{enumerate}
\end{rem}

Given $c>0$ put
\begin{equation}
   \mathcal{H}_c = \{ H \in \mathcal{H} \mid R_H(s,t,x) \leq c,
   \, \forall \, (s,t,x) \in D \times M \}.
\end{equation}

Define
$$A_L = \inf \{ \omega(A) \mid A \in \pi_2(M,L), \, \omega(A)>0 \}.$$
Note that when $L$ is monotone the $\inf$ is attained and $A_L > 0$.
\begin{prop} \label{P:om>=0}
   Suppose $L \subset (M, \omega)$ is monotone.  Let $c <
   \tfrac{A_L}{\pi}$ and $H \in \mathcal{H}_c$. Then for every
   $(J,H)$-holomorphic disk $u$ we have $\omega([u]) \geq 0$ and
   $\mu([u]) \geq 0$.
\end{prop}
\begin{proof}
   Pick $ c < \kappa < \tfrac{A_L}{\pi}$. Then the form
   $\widetilde{\omega}_{H,\kappa}$ is symplectic and $\widetilde{J}_H$
   is compatible with it. Let $u:(D,\partial D) \to (M,L)$ be a
   $(J,H)$-holomorphic disk and let $\widetilde{u} = (z, u(z))$ be its
   graph. Then $\widetilde{u}$ is a {\em non-constant}
   $\widetilde{J}_H$-holomorphic disk with boundary on $\partial D
   \times L$ hence
   $$0< \widetilde{\omega}_{H,\kappa}([\widetilde{u}]) = \omega([u])+
   \pi \kappa.$$
   It follows that $\omega([u]) > -\pi \kappa > -A_L$,
   hence $\omega([u])\geq 0$. The statement on $\mu$ follows from
   monotonicity.
\end{proof}

\subsubsection{Compactness} \label{Sb:compactness}
Let $u_{\nu}$ be a sequence in $\mathcal{M}(A,J,H)$. Then either there
exists a $C^{\infty}$-convergent subsequence or there exists a
subsequence $u_{\nu_k}$ that converges in the sense of Gromov to a
bubble tree consisting of maps $v_0, v_1, \ldots, v_l:(D,\partial D)
\to (M,L)$, $w_1, \ldots, w_r: \mathbb{C}P^1 \to M$, where $v_0$ is
$(J,H)$-holomorphic, $v_1, \ldots, v_l$, $w_1, \ldots, w_r$ are
$J$-holomorphic, at least one of $l,r$ is $\geq 1$ and $[v_0] + [v_1]
+ \ldots + [v_l] + j_*[w_1] + \ldots + j_*[w_r] = A \in
H_2(M,L;\mathbb{Z})$. Here $j_*: H_2(M;\mathbb{Z}) \to
H_2(M,L;\mathbb{Z})$ is the natural homomorphism. The root of the tree
can be thought of as $v_0$.

The relation to Gromov's compactness theorem for genuine
pseudo-holomorphic disks can be understood as follows. Let
$\widetilde{u}_{\nu}:(D,\partial D) \to (D\times M, \partial D \times
L)$ be the sequence of $\widetilde{J}_H$-holomorphic disks obtained as
graphs of the $u_{\nu}$'s. Note that $[\widetilde{u}_{\nu}]=[D]+A \in
H_2(D \times M, \partial D \times L;\mathbb{Z})$ for every $\nu$.  If
$u_{\nu}$ does not have a $C^{\infty}$-convergent subsequence then by
Gromov compactness theorem there exists a subsequence
$\widetilde{u}_{\nu_k}$ which converges in the Gromov topology to a
bubble tree consisting of $\widetilde{J}_H$-holomorphic disks
$\widehat{v}_0, \ldots, \widehat{v}_l$ and
$\widetilde{J}_H$-holomorphic spheres $\widehat{w}_1, \ldots,
\widehat{w}_r$ such that $[\widehat{v}_0]+[\widehat{v}_1] + \cdots +
[\widehat{v}_l] + j'_*[\widehat{w}_1] + \ldots + j'_*[\widehat{w}_r] =
[D]+A$, where $j':H_2(D \times M; \mathbb{Z}) \to H_2(D\times M,
\partial D \times L; \mathbb{Z})$ is the natural map. As $pr_D
:(D\times M, \widetilde{J}_H) \to (D, i)$ is holomorphic it follows
that {\em exactly} one of the $\widehat{v}_i$'s, say $\widehat{v}_0$,
projects surjectively to $D$ while each of $\widehat{v}_1, \ldots,
\widehat{v}_l, \widehat{w}_1, \ldots, \widehat{w}_r$ have constant
projection to $D$. It follows that all the disks and spheres but
$\widehat{v}_0$ are in fact $J$-holomorphic maps lying in fibres of
$pr_D$. Moreover it is easy to see that $\widehat{v}_0$ is, after a
suitable reparametrization, of the form $\widetilde{v}_0$ with $v_0$
being $(J,H)$-holomorphic.

\subsubsection{Perturbation of the spaces $\mathcal{P}_I, \mathcal{P}_{II},
  \mathcal{P}_{III}$} \label{Sb:perturb-P}

Let $J \in \mathcal{J}(M,\omega)$ and $H \in \mathcal{H}$. We will
perturb each of the spaces $\mathcal{P}_I$, $\mathcal{P}_{II_i}$,
$\mathcal{P}_{III_i}$, $i=1,2$ by requiring that in every chain of
disks forming an element of these spaces the (single) disk that has 3
marked points is $(J,H)$-holomorphic. All the other disks in each
chain will remain $J$-holomorphic.

We start with the spaces $\mathcal{P}_I$. Let $\mathbf{A}=(A_1,
\ldots, A_l)$ be a vector of non-zero classes and $1\leq k \leq l$.
Let $x, y \in \textnormal{Crit}(f)$, $a \in \textnormal{Crit}(h)$.
Define $\mathcal{P}_{I}(a, x, y; (\mathbf{A},k), J, H)$ to be the
space of all $(u_1, \ldots, u_l)$ such that:
\begin{enumerate}
  \item $u_i \in \mathcal{M}(A_i,J)/G_{-1,1}$ for every $i \neq k$.
  \item $u_k \in \mathcal{M}(A_k,J,H)$.
  \item $u_1(-1) \in W_x^u$, $u_l(1) \in W_y^s$.
  \item $(u_i(1),u_{i+1}(-1)) \in Q_{f,\rho_L}$ for every $1\leq i
   \leq l-1$, $u_k(0) \in W_a^u$.
\end{enumerate}
In other words elements $(u_1, \ldots, u_l)$ of $\mathcal{P}_I(a,x,y;
(\mathbf{A},k), J,H)$ are the same as those of $\mathcal{P}_I(a,x,y;
(\mathbf{A},k),J)$ only that now $u_k$ is $(J,H)$-holomorphic.  Put
$$\mathcal{P}_I(a,x,y; \mathbf{A},J,H)=\bigcup_{k=1}^l
\mathcal{P}_I(a,x,y;(\mathbf{A},k),J,H).$$

Next define $\mathcal{P}_{II_1}(a,x,y;(\mathbf{A},k),J,H)$ to be the
space of all $(u_1, \ldots, u_{k-1}, v, u_k, \ldots, u_l)$ such that:
\begin{enumerate}
  \item $u_i \in \mathcal{M}(A_i,J)/G_{-1,1}$ for every $1\leq i \leq
   l$.
  \item $v \in \mathcal{M}(0,J,H)$.
  \item $ev_{(k)}(u_1, \ldots, u_{k-1}, v, u_k, \ldots, u_l) \in W_x^u
   \times Q_{f,\rho_L}^{\times(k-1)} \times W_a^u \times
   \textnormal{diag(L)} \times Q_{f,\rho_L}^{\times(l-k)} \times
   W_s^y$, where $ev_{(k)}$ is the evaluation map defined in
   Section~\ref{Sb:external-op}.
\end{enumerate}
Note that there are $(l+1)$ disks in every element of this space and
the disk $v$ which has $3$ marked points is in the class $0$. The
relation to the unperturbed space
$\mathcal{P}_{II_1}(a,x,y;(\mathbf{A},k),J)$ is the following. In the
unperturbed space the disk $v$ corresponds to a constant disks, namely
a point, where this point lies on the boundary of the disk $u_k$ and
is connected by a gradient trajectory coming from $a$. (See
figure~\ref{f:qm-module-II}.)

Similarly define $\mathcal{P}_{II_2}(a,x,y;(\mathbf{A},k),J,H)$ to be
the space of all \linebreak $(u_1, \ldots, u_k, v, u_{k+1}, \ldots,
u_l)$ such that:
\begin{enumerate}
  \item $u_i \in \mathcal{M}(A_i, J)/G_{-1,1}$ for every $1\leq i \leq
   l$.
  \item $v \in \mathcal{M}(0,J,H)$.
  \item $ev_{(k+1)}(u_1, \ldots, u_k, v, u_{k+1}, \ldots, u_l) \in
   W_x^u \times Q_{f,\rho_L}^{\times k} \times \textnormal{diag}(L)
   \times W_a^u \times Q_{f,\rho_L}^{\times(l-k-1)} \times W_y^s$.
\end{enumerate}
Finally for $i=1,2$ put
$$\mathcal{P}_{II_i}(a,x,y;\mathbf{A},J,H)=\bigcup_{k=1}^l
\mathcal{P}_{II_i}(a,x,y;(\mathbf{A},k),J,H).$$

We now turn to the spaces $\mathcal{P}_{III}$.  Let
$\mathbf{B}'=(B'_1, \ldots, B'_l)$, $\mathbf{B}''=(B''_1, \ldots,
B''_{l''})$ two vectors of non-zero classes. Let $1\leq k'\leq l'$,
$1\leq k''\leq l''$. We define
$\mathcal{P}_{III_1}(a,x,y;(\mathbf{B}',k'),\mathbf{B}'',J,H)$ (resp.
$\mathcal{P}_{III_1}(a,x,y;(\mathbf{B}',k'),\mathbf{B}'',J,H)$) in the
same way as in~\eqref{Eq:P-III} only that now $u'_{k'}$ (resp.
$u''_{k''}$) is $(J,H)$-holomorphic. Put
\begin{align*}
   & \mathcal{P}_{III_1}(a,x,y;\mathbf{B}',\mathbf{B''},J,H) =
   \bigcup_{k'=1}^{l'}
   \mathcal{P}_{III_1}(a,x,y;(\mathbf{B}',k'),\mathbf{B}'',J,H) \\
   & \mathcal{P}_{III_2}(a,x,y;\mathbf{B}',\mathbf{B}'',J,H) =
   \bigcup_{k''=1}^{l''}
   \mathcal{P}_{III_2}(a,x,y;\mathbf{B}',(\mathbf{B}'',k''),J,H).
\end{align*}
Note that the space $\mathcal{P}_{II_1}(a,x,y;(\mathbf{A},k),J,H)$
could be viewed as a special case of the space
$\mathcal{P}_{III_1}(a,x,y;(\mathbf{B'},l'),\mathbf{B}'',J,H)$ if we
would have allowed the $l'$'th class $B'_{l'}$ in $\mathbf{B}'$ to be
$0$. Namely the space $\mathcal{P}_{II_1}(a,x,y;(\mathbf{A},k),J,H)$
is the same as
$\mathcal{P}_{III_1}(a,x,y;(\mathbf{B}',l'),\mathbf{B}'',J,H)$ with
$\mathbf{B}'=(A_1, \ldots, A_{k-1},0)$, $\mathbf{B}''=(A_k, \ldots,
A_l)$. A similar remark holds for the perturbations of
$\mathcal{P}_{II_2}$ and $\mathcal{P}_{III_2}$. Nevertheless it will
be more convenient for us to work with vectors $\mathbf{B}'$,
$\mathbf{B}''$ of {\em non-zero classes}, hence separate the
perturbations of $\mathcal{P}_{II_i}$ from those of
$\mathcal{P}_{III_i}$.

The following proposition is analogous to Propositions~\ref{P:qm-axy}
and~\ref{P:qm-0-dim}.
\begin{prop} \label{P:qm-axy-pert-1}
   Let $(f,\rho_L)$, $(h,\rho_M)$ be as in assumption~\ref{A:MS}. Then
   there exists a second category subset
   $\mathcal{J}_{\textnormal{reg}} \subset \mathcal{J}(M,\omega)$ such
   that for every $J \in \mathcal{J}_{\textnormal{reg}}$ there exists
   a second category subset $\mathcal{H}_{\textnormal{reg}}(J) \subset
   \mathcal{H}$ with the following properties. For every $x,y \in
   \textnormal{Crit}(f)$, $a \in \textnormal{Crit}(h)$ and every
   vectors of non-zero classes $\mathbf{A}$, $\mathbf{B}'$,
   $\mathbf{B}''$ the following holds:
   \begin{enumerate}
     \item If $|a|+|x|-|y|+\mu(\mathbf{A})-2n=0$ then
      $\mathcal{P}_{I}(a,x,y; \mathbf{A},J,H)$ is a compact
      $0$-dimensional manifold.
     \item If $|a|+|x|-|y|+\mu(\mathbf{A})-2n=1$ then each of
      $\mathcal{P}_{II_i}(a,x,y;\mathbf{A},J,H)$, $i=1,2$, is a
      compact $0$-dimensional manifold.
     \item if
      $|a|+|x|-|y|+\mu(\mathbf{B}')+\mu(\mathbf{B}'')-2n=1$ then
      each of $\mathcal{P}_{III_i}(a,x,y;\mathbf{B}',
      \mathbf{B}'',J,H)$, $i=1,2$, is a compact $0$-dimensional
      manifold.
   \end{enumerate}
   Moreover, in each of the above three cases the $J$-holomorphic part
   of every element of these spaces consist of simple and absolutely
   distinct disks.
\end{prop}
The proof is very similar to that of
Propositions~\ref{P:qm-axy},~\ref{P:qm-0-dim} hence we omit it.

In view of Proposition~\ref{P:qm-axy-pert-1}, for $J \in
\mathcal{J}_{\textnormal{reg}}$, $H \in
\mathcal{H}_{\textnormal{reg}}(J)$ define:
\begin{align*}
   & n_I(a,x,y;\mathbf{A},H) = \#_{\mathbb{Z}_2}
   \mathcal{P}_I(a,x,y;\mathbf{A},J,H), \quad \textnormal{when }
   |a|+|x|-|y|+\mu(\mathbf{A})-2n=0, \\
   & n_{II_i}(a,x,y;\mathbf{A},H) = \#_{\mathbb{Z}_2}
   \mathcal{P}_{II_i}(a,x,y;\mathbf{A},J,H), \quad \textnormal{when }
   |a|+|x|-|y|+\mu(\mathbf{A})-2n=1, \\
   & n_{III_i}(a,x,y;\mathbf{B}',\mathbf{B}'',H) =
   \#_{\mathbb{Z}_2}\mathcal{P}_{III_i}(a,x,y;
   \mathbf{B}',\mathbf{B}'',J,H), \\
   & \textnormal{when }
   |a|+|x|-|y|+\mu(\mathbf{B}')+\mu(\mathbf{B}'')-2n=1.
\end{align*}
(For simplicity we have suppressed here the $J$'s from the notation.)

We will also need the following analogue of statement~(1) of
Proposition~\ref{P:qm-1-dim}. Note however that now we do not require
$N_L$ to be $\geq 3$ anymore.
\begin{prop} \label{P:qm-axy-pert-2}
   Let $(f,\rho_L)$, $(h,\rho_M)$ be as in assumption~\ref{A:MS}. Then
   there exists a second category subset
   $\mathcal{J}_{\textnormal{reg}} \subset \mathcal{J}(M,\omega)$ such
   that for every $J \in \mathcal{J}_{\textnormal{reg}}$ the following
   holds. There exists a second category subset
   $\mathcal{H}_{\textnormal{reg}}(J) \subset \mathcal{H}$ such that
   for every $H \in \mathcal{H}_{\textnormal{reg}}(J)$, every vector
   of non-zero classes $\mathbf{A}$ and every $x,y \in
   \textnormal{Crit}(f)$, $a \in \textnormal{Crit}(h)$ with
   $|a|+|x|-|y|+\mu(\mathbf{A})-2n=1$ the space
   $\mathcal{P}_I(a,x,y; \mathbf{A},J,H)$ is a smooth $1$-dimensional
   manifold. Moreover, the $J$-holomorphic part of every element of
   this space consist of simple and absolutely distinct disks.
\end{prop}

\begin{proof}
   The proof is very similar to that of Propositions~\ref{P:qm-1-dim}
   and~\ref{P:qm-axy} however in contrast to
   Proposition~\ref{P:qm-1-dim} we do not need to assume here that
   $N_L \geq 3$ anymore (i.e. $N_L \geq 2$ is enough). The point is
   that in order to assure that the space
   $\mathcal{P}_I(a,x,y;(\mathbf{A},k),J,H)$ is a smooth manifold it
   is enough to require that the disks $(u_1, \ldots, u_{k-1},
   u_{k+1}, \ldots, u_l)$ are simple and absolutely distinct without
   any requirement on the (perturbed) disk $u_k$. Recall that it was
   due to {\em this} disk that we had to assume in
   Proposition~\ref{P:qm-1-dim} that $N_L \geq 3$ in order to assure
   simplicity. Finally note that absolutely distinctiveness of the
   $J$-holomorphic part of each element of $\mathcal{P}_I(a,x,y;
   \mathbf{A},J,H)$ does not require $N_L \geq 3$ either and can be
   proved in a similar way as in Proposition~\ref{P:qm-1-dim}.
\end{proof}

\subsubsection{Identities for the perturbed spaces}
\label{Sb:identities-pert}

\begin{prop} \label{P:qm-idents-pert}
   Let $(f,\rho_L)$, $(h,\rho_M)$ be as in assumption~\ref{A:MS}. Then
   there exists a second category subset
   $\mathcal{J}_{\textnormal{reg}} \subset \mathcal{J}(M,\omega)$ such
   that for every $J \in \mathcal{J}_{\textnormal{reg}}$ there exists
   a second category subset $\mathcal{H}_{\textnormal{reg}}(J) \subset
   \mathcal{H}$ with the following property. For every $H \in
   \mathcal{H}_{\textnormal{reg}}(J)$, every $x,y \in
   \textnormal{Crit}(f)$, $a \in \textnormal{Crit}(h)$ and
   $\mathbf{A}$ with $|a|+|x|-|y|+\mu(\mathbf{A})-2n=1$
   identity~\eqref{Eq:qm-ident-1} holds with $n_I(-)$ replaced by
   $n_I(-,H)$ $n_{II_i}(-)$ by $n_{II_i}(-, H)$ and $n_{III_i}(-)$ by
   $n_{III_i}(-,H)$.
\end{prop}

\begin{proof}
   The proof is very similar to that of Proposition~\ref{P:qm-idents}
   with the following two slight changes. Now we use the compactness
   statement from Section~\ref{Sb:compactness} rather than compactness
   for genuinely pseudo-holomorphic disks. The second difference is
   that the gluing procedure between a $J$-holomorphic disk $u'$ and a
   $(J,H)$-holomorphic disks $u''$ with $u'(1)=u''(-1)$ is now carried
   out in $D \times M$. Namely we consider $u'$ as a
   $\widetilde{J}_H$-holomorphic disk in the fibre $\{-1\} \times M$
   and glue it to the $\widetilde{J}_H$-holomorphic section
   $\widetilde{u}''$ corresponding to $u''$. Standard arguments imply
   that the result of the gluing is (after a suitable
   reparametrization) again a ($\widetilde{J}_H$-holomorphic) {\em
     section} hence descends to a $(J,H)$-holomorphic disk. The
   precise details of this type of gluing are given in
   Section~\ref{Sbsb:glue-ex-H}.

   The rest of the proof is a straightforward adaption of the proof of
   Proposition~\ref{P:qm-idents}
\end{proof}

\begin{prop} \label{P:n-pert=no-pert}
   Let $(f,\rho_L)$, $(h,\rho_M)$ be as in assumption~\ref{A:MS}.
   Then there exists a second category subset
   $\mathcal{J}_{\textnormal{reg}} \subset \mathcal{J}(M,\omega)$ and
   a neighbourhood (in the $C^{\infty}$-topology) $\mathcal{U}$ of $0
   \in \mathcal{H}$ such that for every $J \in
   \mathcal{J}_{\textnormal{reg}}$ there exists a second category
   subset $\mathcal{H}_{\textnormal{reg}}(J) \subset \mathcal{H}$ with
   the following properties. For every $H \in
   \mathcal{H}_{\textnormal{reg}}(J) \cap \mathcal{U}$, every $x,y \in
   \textnormal{Crit}(f)$, $a \in \textnormal{Crit}(h)$ and every
   vectors $\mathbf{A}$, $\mathbf{B}'$, $\mathbf{B}''$ of non-zero
   classes the following holds:
   \begin{enumerate}
     \item If $|a|+|x|-|y|+\mu(\mathbf{A})-2n=0$ then
      $n_I(a,x,y;\mathbf{A},H) = n_I(a,x,y;\mathbf{A})$.
     \item If $|a|+|x|-|y|+\mu(\mathbf{A})-2n=1$ then
      $n_{II_i}(a,x,y;\mathbf{A},H) = n_{II_i}(a,x,y;\mathbf{A})$,
      $i=1,2$.
     \item If $|a|+|x|-|y| +
      \mu(\mathbf{B}')+\mu(\mathbf{B}'')-2n=1$ \\
      then $n_{III_i}(a,x,y;\mathbf{B}',\mathbf{B}'',H) =
      n_{II_i}(a,x,y;\mathbf{B}',\mathbf{B}'')$, $i=1,2$.
   \end{enumerate}
\end{prop}

\begin{proof}
   The proof is based on a standard cobordism argument. We outline the
   main ideas. Suppose that we are under the assumption appearing in
   statement~(1) of the Proposition. We define the set
   $\mathcal{J}_{\textnormal{reg}}$ to be the intersection of all the
   $\mathcal{J}_{\textnormal{reg}}$'s from the previous propositions
   and require in addition that for every $J \in
   \mathcal{J}_{\textnormal{reg}}$ the following holds:
   \begin{enumerate}
     \item[(i)] For every two vectors $\mathbf{C}'$, $\mathbf{C}''$ of
      non-zero classes and every $x',y' \in \textnormal{Crit}(f)$,
      $a'\in \textnormal{Crit}(h)$ with
      $|a'|+|x'|-|y'|+\mu(\mathbf{C}') +
      \mu(\mathbf{C}'')-2n\leq 0$ the spaces
      $\mathcal{P}_{III_i}(a',x',y';\mathbf{C}',\mathbf{C}'',J)$,
      $i=1,2$, are {\em empty}.
     \item[(ii)] For every vector of non-zero classes $\mathbf{A}'$
      and every $x',y' \in \textnormal{Crit}(f)$, $a \in
      \textnormal{Crit}(h)$ with
      $|a'|+|x'|-|y'|+\mu(\mathbf{A}')-2n\leq 0$ the spaces
      $\mathcal{P}_{II_i}(a',x',y';\mathbf{A}',J)$, $i=1,2$, are {\em
        empty}.
   \end{enumerate}
   Standard arguments imply that the resulting set
   $\mathcal{J}_{\textnormal{reg}}$ is of second category. We now
   claim the following:
   \begin{statement} \label{St:pert-U}
      There exists a neighbourhood $\mathcal{U}$ of $0 \in
      \mathcal{H}$ such that for every $H \in \mathcal{U}$ and every
      $a,x,y$, $\mathbf{A}$ with $|a|+|x|-|y| + \mu(\mathbf{A})-2n
      \leq 0$ the following holds:
      \begin{enumerate}
        \item For every splitting
         $\mathbf{A}=(\mathbf{C}',\mathbf{C}'')$ of $\mathbf{A}$ into
         two vectors $\mathbf{C}'$, $\mathbf{C}''$ of non-zero classes
         we have
         $\mathcal{P}_{III_i}(a,x,y;\mathbf{C}',\mathbf{C}'',J,H) =
         \emptyset$ for $i=1,2$.
        \item $\mathcal{P}_{II_i}(a,x,y;\mathbf{A},J,H)=\emptyset$,
         for $i=1,2$.
      \end{enumerate}
   \end{statement}
   We prove~\ref{St:pert-U}. Indeed, if the contrary to
   statement~\ref{St:pert-U} happens then there exists a sequence
   $H_{\nu} \in \mathcal{H}$ with $H_{\nu} \longrightarrow 0$ in the
   $C^{\infty}$-topology and elements $\mathbf{u}_{\nu}$ in either
   $\mathcal{P}_{III_i}(a_{\nu}, x_{\nu},
   y_{\nu};\mathbf{C}'_{\nu},\mathbf{C}''_{\nu},J,H_{\nu})$ or
   $\mathcal{P}_{II_i}(a_{\nu},x_{\nu},y_{\nu};\mathbf{A}_{\nu},J,H_{\nu})$,
   where $x_{\nu}, y_{\nu}, a_{\nu}$, $\mathbf{A}_{\nu}$,
   $\mathbf{C}'_{\nu}$, $\mathbf{C}''_{\nu}$ satisfy $|a_{\nu}| +
   |x_{\nu}| - |y_{\nu}| + \mu(\mathbf{A}_{\nu})-2n \leq 0$ or
   $|a_{\nu}| + |x_{\nu}| - |y_{\nu}| +
   \mu(\mathbf{C}'_{\nu})+\mu(\mathbf{C}''_{\nu})-2n \leq 0$.  Viewing
   the disks in $\mathbf{u}_{\nu}$ as
   $\widetilde{J}_{H_{\nu}}$-holomorphic disks in $D \times M$ (one as
   a section and the others lying in the fibres of $D \times M \to D$)
   it is easy to see that their energy is uniformly bounded for
   symplectic forms of the type $\widetilde{\omega}_{H_{\nu},\kappa}$
   with fixed $\kappa$.  (See~\eqref{Eq:om-hk} and
   Proposition~\ref{P:om>=0}).  Therefore by Gromov compactness
   theorem, after passing to a subsequence, the sequence converges to
   a chain of bubble trees $\mathbf{u}_{\infty}$ all of whose elements
   are $(J,0)$-holomorphic (i.e. genuinely $J$-holomorphic) disks and
   spheres. As $J \in \mathcal{J}_{\textnormal{reg}}$, by a dimension
   count argument we can rule out all possible configurations for
   $\mathbf{u}_{\infty}$ except of maybe the cases
   $\mathbf{u}_{\infty} \in
   \mathcal{P}_{III_i}(a,x,y;\mathbf{C}',\mathbf{C}'',J)$ or
   $\mathbf{u}_{\infty} \in \mathcal{P}_{II_i}(a,x,y;\mathbf{A}',J)$
   for some $a'$, $x'$, $y'$, $\mathbf{A}$ or $\mathbf{C}',
   \mathbf{C}''$ that satisfy $|a'|+|x'|-|y'|+\mu(\mathbf{C}') +
   \mu(\mathbf{C}'')-2n\leq 0$ or
   $|a'|+|x'|-|y'|+\mu(\mathbf{A}')-2n\leq 0$. However these two last
   cases are again impossible by the definition of
   $\mathcal{J}_{\textnormal{reg}}$ in~(i) and~(ii) above.  This
   proves statement~\ref{St:pert-U}.

   We continue with the proof of Proposition~\ref{P:n-pert=no-pert}.
   Fix $J \in \mathcal{J}_{\textnormal{reg}}$ and let
   $\mathcal{H}_{\textnormal{reg}}(J)$ be the intersection of all the
   $\mathcal{H}_{\textnormal{reg}}(J)$'s from the previous
   Propositions. Let $\mathcal{U}$ be the neighbourhood of $0 \in
   \mathcal{H}$ defined by statement~\ref{St:pert-U}. Replace
   $\mathcal{U}$ if needed by its connected component containing $0$.
   Let $H \in \mathcal{H}_{\textnormal{reg}}(J) \cap \mathcal{U}$.
   Take a generic path $\{H_{\lambda}\}_{0 \leq \lambda \leq 1}$ in
   $\mathcal{U}$ with $H_0=H$ and $H_1=0$. Then
   $$\mathcal{P}_{I}(a,x,y;\mathbf{A},J,\{ H_{\lambda}\}) = \{
   (\lambda, \mathbf{u}) \mid 0 \leq \lambda \leq 1, \, \mathbf{u} \in
   \mathcal{P}_{I}(a,x,y;\mathbf{A},J,H_{\lambda})\}$$
   is a
   $1$-dimensional smooth manifold whose boundary is
   \begin{equation} \label{Eq:bndry-cobord-H}
      \partial
      \mathcal{P}_{I}(a,x,y;\mathbf{A},J,\{ H_{\lambda}\}) =
      \mathcal{P}_I(a,x,y;\mathbf{A},J,H_0) \coprod
      \mathcal{P}_I(a,x,y;\mathbf{A},J).
   \end{equation}
   It is important to note here that the space lying over $\lambda=1$,
   namely $\mathcal{P}_I(a,x,y;\mathbf{A},J)$ consists of chains of
   disks that are simple and absolutely distinct. This is assured by
   Proposition~\ref{P:qm-axy} (which holds for $N_L\geq 2$).

   We claim that the $1$-dimensional manifold
   $\mathcal{P}_{I}(a,x,y;\mathbf{A},J,\{ H_{\lambda}\})$ is compact.
   This follows form a straightforward dimension count. Indeed
   occurrence of bubbling of spheres reduces the dimension by at least
   $2$ hence cannot appear in a generic $1$-dimensional family.
   Similarly, bubbling of a $J$-holomorphic disk at a point $p \in
   \partial D$ with $p \neq -1,1$ also yields negative dimension since
   $N_L \geq 2$. It remains to rule out the following three cases:
   \begin{itemize}
     \item Shrinking to a point of a gradient trajectory connecting
      two consecutive disks in $\mathcal{P}_{I}(a,x,y;\mathbf{A},J,\{
      H_{\lambda}\})$.
     \item Bubbling of a $J$-holomorphic disk at $p = -1$ or $p=1$ in
      one of the disks participating in
      $\mathcal{P}_{I}(a,x,y;\mathbf{A},J,\{ H_{\lambda}\})$.
     \item Breaking of a gradient trajectories (of $f$ or $h$)
      involved in $\mathcal{P}_{I}(a,x,y;\mathbf{A},J,\{
      H_{\lambda}\})$.
   \end{itemize}
   The first possibility is ruled out by~(1) of
   statement~\ref{St:pert-U}. The second possibility is ruled out
   by~(1) and~(2) of statement~\ref{St:pert-U}. Note that we have to
   apply statement~\ref{St:pert-U} for a different $\mathbf{A}$ here.
   Namely we have to take the $\mathbf{C}'$, $\mathbf{C}''$ coming
   from the bubbling and change $\mathbf{A}$ to be
   $(\mathbf{C}',\mathbf{C}'')$. Also note that when bubbling of a
   disk occurs in the perturbed disk, the remaining perturbed
   component might have now $\mu=0$. This is why we needed also~(2)
   in statement~\ref{St:pert-U}. The third possibility can be ruled
   out in a similar way to the first one. This completes the proof
   that $\mathcal{P}_{I}(a,x,y;\mathbf{A},J,\{ H_{\lambda}\})$ is
   compact.

   To conclude the proof of Proposition~\ref{P:n-pert=no-pert} note
   that since $\mathcal{P}_{I}(a,x,y;\mathbf{A},J,\{ H_{\lambda}\})$
   is compact it follows from~\eqref{Eq:bndry-cobord-H} that
   $n_I(a,x,y;\mathbf{A},H) = n_I(a,x,y;\mathbf{A})$.

   The proofs of statements~(2) and~(3) are similar.
\end{proof}

\subsubsection{Conclusion of the Proof of identity~\eqref{Eq:qm-ident-1}
  for every $N_L \geq 2$} \label{Sb:conclusion-prf-qm-ident-all-N}
Identity~\eqref{Eq:qm-ident-1} follows now immediately by combining
Propositions~\ref{P:qm-idents-pert} and~\ref{P:n-pert=no-pert}.

\

By similar methods, it is not difficult to see that, in homology, the
operation defined by formula~\eqref{Eq:qm-*} is independent of the
choices made in its definition.

\

\subsubsection{Module structure.} The purpose of this section is to show
the following identities:
\begin{prop}\label{prop:module_str}
   For any $\alpha \in Q^{+}H_{\ast}(L)$, $a,b\in
   Q^{+}H_{\ast}(M)$ we have
   \begin{enumerate}
     \item $(a\ast b)\ast \alpha = a\ast (b\ast \alpha)$.
      \label{I:mod-assoc}
     \item $[M]*\alpha = \alpha$. (Recall that $[M] \in Q^+H_{2n}(M)$
      is the unit with respect to the quantum cap product).
      \label{I:mod-unit}
   \end{enumerate}
\end{prop}

We will present first the proof of point~(\ref{I:mod-assoc}) in a
particular case: when $N_{L}\geq 4$. The proof in this case is the
cleanest as it does not require the use of perturbations.  The case
when $N_{L}\leq 4$ will be treated later.  It is also useful to note
that if using the full Novikov ring instead of $Q^{+}H(-)$, the
formula in point~(\ref{I:mod-assoc}) of the Proposition follows (for
any $N_L \geq 2$) from the comparison with Floer homology described in
the next section.

\begin{proof} We start with the proof of point~(\ref{I:mod-assoc}).
   The proof is very similar to the one for the associativity of the
   quantum product as described in the proof of Lemma
   \S\ref{lem:prodassoc}.  One important point which is useful to note
   before going further is that it is essential for this formula that
   the product used on $Q^{+}H_{\ast}(M)$ is the quantum cap product
   (rather than the classical cap-product).

   We now proceed with the description of the necessary moduli spaces
   (this is, obviously, independent of our restriction on $N_{L}$).

   \

   It is again useful to use the language of trees as we already did
   in the section on the quantum product \S \ref{subsec:product}. We
   first notice that the moduli spaces in the first subsections of \S
   \ref{S:qm} can all be described in these terms. For example
   $\mathcal{P}_{I}(a,x,y;\mathbf{A},J)$ is modeled on a tree with two
   entries so that one of the entries is connected by an edge to the
   unique vertex of valence three which is labeled by a non-zero
   element of $H_2(M,L;\mathbb{Z})$. This edge corresponds to a flow
   line of the function $h$.  It is useful to consider the edges of
   the tree being labeled by either the function $f$ or the function
   $h$.  The same type of description applies to the other moduli
   spaces each of which is obtained by adding some auxiliary
   conditions: $\mathcal{P}_{I'}(a,x,y;\mathbf{A},J)$ has the same
   description as $\mathcal{P}_{I}(a,x,y;\mathbf{A},J)$ except that
   the vertex of valence three is a constant disk (in other words the
   labeling of this vertex is by the element $0\in
   H_2(M,L;\mathbb{Z}$).  Further, the two moduli spaces
   $\mathcal{P}_{II_{1}}(a,x,y;\mathbf{A},J)$ and
   $\mathcal{P}_{II_{2}}(a,x,y;\mathbf{A},J)$ are both obtained as
   before but with the additional conditions that one of the edges
   labeled by $f$ is of length $0$ and one of the two adjacent
   vertices of this edge is a constant disk of valence three.  For the
   moduli spaces $\mathcal{P}_{III_{1}}(a,x,y;\mathbf{A},J)$,
   $\mathcal{P}_{III_{2}}(a,x,y;\mathbf{A},J)$ the requirement is that
   one of the edges labeled by $f$ is of length $0$ and, finally, the
   spaces $\mathcal{P}_{III'_{1}}(a,x,y;\mathbf{A},J)$ and
   $\mathcal{P}_{III'_{2}}(a,x,y;\mathbf{A},J)$ are characterized by
   the fact that one edge labeled by $f$ is of length $0$ and the
   vertex of valence three is a constant disk.

   \

   We now start to discuss the specific moduli spaces required in the
   proof. These moduli spaces are similar to the one used to show the
   associativity of the quantum product in \S\ref{subsec:product}. As
   there, they are modeled on trees with three entries and one exit.
   We will also need to use two Morse-Smale functions on $M$, $h_{1}$
   and $h_{2}$, together with the Morse function $f$ on $L$. We will
   assume that $h_{1}$ and $h_{2}$ have the same critical points and
   are in generic position. We consider planar trees $\mathcal{T}$
   with three entries and one exit and which are labeled as follows:
   each edge is labeled by a number from $\{1,2,3\}$ and each vertex
   is labeled by an element of $H_2(M,L;\mathbb{Z})$.

   For each topological type of such a tree, $\mathcal{T}$, and for
   $a,b\in \Crit(h_{1})$, $x,y\in \Crit(f)$ we denote by
   $\mathcal{P}_{\mathcal{T}}(a,b,x,y;J)$ the moduli space of
   configurations which satisfy:
   \begin{itemize}
     \item[(1)] There are four vertices of valence $1$, the three
      entries and the single exit. The graph is planar with the three
      entry points on the line $\R\times\{1\}$ the $i$-th entry point
      (along this line) is the origin of an edge labeled by $i$, the
      edge arriving in the exit point is labeled by $3$.
     \item[(2)]Each vertex corresponds to a pseudo-holomorphic disk
      with boundary on $L$ or to a $J$-holomorphic sphere.
     \item[(3)] Each edge corresponds to a flow line (of non-zero length)
      so that if the label of that edge is $1$, the corresponding
      function is $h_{1}$, if the label is $2$ the function is $h_{2}$
      and if the label is $3$ then the function is $f$.
     \item[(4)] A vertex of valence $2$ corresponds to a
      pseudo-holomorphic disk of non-vanishing class.
     \item[(5)] For each vertex of valence $2$ the incidence relations
      are as in the definition of the moduli spaces used in the pearl
      complex, in particular, both the entry and exit edges are labeled
      by $3$.
     \item[(6)] Any vertex of trivial class is at least of valence $3$
     \item[(7)] If a vertex is of valence $3$ and is a pseudo-holomorphic
      disk, then one of the entering edges and the exit edge are labeled
      by $3$ and the incidence relations at this vertex are as in the
      definition of the moduli space $\mathcal{P}_{I}(----)$.
     \item[(8)] If a vertex of valence $3$ corresponds to a
      pseudo-holomorphic sphere, then the two entry edges are labeled
      by $1$ and $2$ and the exit edge by $1$ so that the incidence
      points are $1$ - which is the attaching point for the entering
      edge labeled by $1$, $e^{2\pi i/3}$ - which corresponds to the
      exit edge and $e^{4\pi i /3}$ - which is the attaching point of
      the edge labeled by $3$ (we view here $S^{2}$ as $\mathbb{C}
      \cup \{ \infty \}$).
     \item[(9)] If a vertex is of valence $4$, then it corresponds to
      a $J$-holomorphic disk with boundary on $L$ and the incidence
      relations are so that the point $-1$ is the attaching point of
      an entering flow line of $f$, the point $+1$ is the attaching
      point of the exiting edge which is again a flow line of $f$, the
      point $0$ is the attaching point of an entering flow line
      corresponding to $h_{2}$, there is a real point of coordinates
      $(x,0)$ with $0<x<1$ which is the attaching point of a flow line
      of $h_{1}$.
   \end{itemize}

   Notice that, with this definition, there can not be more than a
   single $J$-holomorphic sphere in such a configuration and,
   moreover, if such a sphere appears there are a $h_{1}$-flow line
   starting from $a$ and a $h_{2}$-flow line starting from $b$ which
   both arrive at the sphere as as well as an exiting $h_{1}$-flow
   line which starts from the sphere and ends at the center of a disk
   with boundary on $L$ (this disk might be constant).

   These moduli spaces will be used to define a chain homotopy $\xi$
   so that we have:
   \begin{equation}\label{eq:htpy}
      (d\xi+\xi d)(a\otimes b\otimes x)=(a\ast b)\ast x
      + a\ast (b\ast
      x)~.~
   \end{equation}
   For this we need to first study the regularity properties of our
   moduli spaces, $\mathcal{P}_{\mathcal{T}}(a,b,x,y)$. To do so we
   let $\mathcal{P}^{\ast,d}_{\mathcal{T}}(a,b,x,y)$ be the subset of
   $\mathcal{P}_{\mathcal{T}}(a,b,x,y)$ which is constituted by the
   configurations in which each disk and sphere is simple and they are
   all absolutely distinct. It is clear, by the standard
   transversality arguments used earlier in the paper that, for a
   generic almost complex structure $J$,
   $\mathcal{P}^{\ast,d}_{\mathcal{T}}(a,b,x,y)$ is a manifold of
   dimension $|a|+|b|+|x|-|y|-4n+1+\mu(\mathcal{T})$. We denote
   $\delta'''=\delta'''(a,b,x,y,\mathcal{T}) =
   |a|+|b|+|x|-|y|-4n+1+\mu(\mathcal{T})$ and we intend to show:

   \begin{lem}\label{lem:tpras}
      For $\delta'''\leq 1$ and if $N_{L}\geq 4$, and a generic choice
      of almost complex structure, we have:
      $$\mathcal{P}_{\mathcal{T}}(a,b,x,y) =
      \mathcal{P}^{\ast,d}_{\mathcal{T}}(a,b,x,y)~.~$$
   \end{lem}

   \begin{proof}[Proof of Lemma~\ref{lem:tpras}]
      The first step is to reduce to the case of simple disks and
      spheres. The argument is very similar to the proof of the
      associativity for the quantum product combined with some
      arguments already used to show that the external product is a
      chain map.  First we notice that as spheres are either simple or
      multiple covered it is easy to see that we may assume that the
      sphere appearing in a configuration $u\in
      \mathcal{P}_{\mathcal{T}}(a,b,x,y)$ for $\delta'''\leq 1$ is
      simple. So we are left with proving that all the disks can also
      be assumed to be simple and that all the objects involved are
      absolutely distinct. In the proof we also need a type of
      slightly more general moduli spaces than the ones described
      before.  They will be denoted by
      $\mathcal{P}'_{\mathcal{T}}(a,b,x,y)$ and they are characterized
      by the the fact that in condition (7) above we allow for the
      interior marked point (inside the disk) to be different from $0$
      and, similarly, in condition (9) the two marked points are again
      allowed to be arbitrary.  The virtual dimension of this moduli
      space is $\delta'''+2$. The proof proceeds by recurrence: we
      intend to show by induction over $\mu(\mathcal{T})$ that if
      $\delta'''=1$, then the claim is true. For this we will need to
      also show that if $\delta'''-2<0$, then
      $\mathcal{P}'_{\mathcal{T}}(a,b,x,y)$ can only contain
      configurations made out of simple, absolutely distinct curves
      and thus $\mathcal{P}'_{\mathcal{T}}(a,b,x,y)=\emptyset$ in this
      case.  As mentioned before, we may already assume that the
      spheres appearing in the configurations described above are
      simple.  To shorten notation we will denote in the paragraph
      below: $\mathcal{P}=\mathcal{P}_{\mathcal{T}}(a,b,x,y)$ and
      $\mathcal{P}'=\mathcal{P}'_{\mathcal{T}}(a,b,x,y)$. To prove our
      claim we will actually need to work with a type of even more
      general moduli spaces (this is similar to the first part of the
      proof of Proposition). We will denote these moduli spaces by $
      \mathcal{G}=\mathcal{G}_{\mathcal{T}}(a,b,x,y)$, $\mathcal{G}'=
      \mathcal{G}'_{\mathcal{T}}(a,b,x,y)$ they consist of
      configurations as those in $\mathcal{P}$ and, respectively, in
      $\mathcal{P}'$ except that: condition (5) is modified so as to
      allow vertices of valence two so that the entering edge is
      labeled by $1$ or $2$ and the exit one by $3$; condition (7) and
      (8) are modified so that vertices of valence three or four
      (which correspond to disks) are allowed to have more than one
      (possibly even three) entering edges labeled by $3$ (the
      incidence points on such a disk are assumed distinct, two of
      them being $+1$ and $-1$). We intend to show that, for a generic
      $J$, these moduli spaces are formed by configurations containing
      only simple absolutely distinct curves. This obviously implies
      the statement for $\mathcal{P}$ and $\mathcal{P}'$.

      We now consider a configuration $u$ that might appear either in
      $\mathcal{G}$ or in $\mathcal{G}'$, we assume that $\delta'''
      \leq 1$ if it is in $\mathcal{G}$ and that $\delta'''-2< 0$ in
      the second case. We also assume the statement proved for the
      moduli spaces with $\mu(\mathcal{T})< k$ so that we take here
      $\mu(\mathcal{T})=k$.  Suppose that one disk $u_{0}$ appearing
      in $u$ is not simple. There are two cases to consider. The first
      is when $n = \dim L \geq 3$ so that we may use
      Lemma~\ref{L:non-simple}.  Assuming now that $u$ is in
      $\mathcal{G}'$ we see that we may replace $u_{0}$ by another
      disk $u_{0}'$ which is simple so that
      $u_{0}(-1/+1)=u_{0}'(-1/+1)$, $u_{0}(D)= u_{0}'(D)$ and
      $\mu([u_{0}'])<\mu([u_{0}])$.  and so, by replacing $u_{0}$ by
      $u_{0}'$, we obtain a configuration of lower Maslov class than
      that of $u$ and which belongs to $\mathcal{G}'$ which
      contradicts the induction hypothesis. Suppose now that $u\in
      \mathcal{G}$.  It is only here that the hypothesis $N_{L}\geq 4$
      intervenes. Indeed, we suppose again that the disk $u_{0}$
      appearing in $u$ is not simple and we let $u_{0}'$ be the simple
      disk provided as above by the Lemma \ref{L:non-simple}.  We
      consider the configuration $u'$ obtained by replacing $u_{0}$ by
      $u_{0}'$.  We obviously have $\mu(u')\leq \mu(u)-N_{L} \leq
      \mu(u)-4$. The key remark is that, in general, $u'$ is not an
      element in $\mathcal{G}$ but is an element of
      $\mathcal{G}'_{\mathcal{T}'}$ (for a suitable tree
      $\mathcal{T}'$).  The virtual dimension of this last moduli
      space is
      $$|a|+|b|+|x|-|y|-4n+3+\mu(u')\leq |a|+|b|+|x|-|y| - 4n +3
      +\mu(u)-4=\delta'''(a,b,x,y,\mathcal{T})-2\leq -1.$$
      As
      $\mu(\mathcal{T}')< k$ we again see that this leads to a
      contradiction and so all disks may be assumed simple at least
      for $n\geq 3$. We pursue under this assumption to show that all
      the curves involved may also be assumed to be absolutely
      distinct.  Lemma \ref{L:uv} plays an important role here.
      Indeed, it implies that if the curves are not absolutely
      distinct, then there are two disks $u_{1}$ and $u_{2}$ so that
      $u_{1}(D)\subset u_{2}(D)$ and the same relation is valid for
      the respective boundaries.  We intend to use this fact to
      eliminate $u_{1}$ from the tree. It is not difficult to see that
      by appropriately replacing $u_{1}$ by $u_{2}$ in the tree and
      possibly eliminating the chain connecting $u_{1}$ to $u_{2}$ in
      $\mathcal{T}$ we obtain a configuration which still belongs to
      $\mathcal{G}'_{\mathcal{T}'}$ (for some other appropriate tree
      $\mathcal{T}'$) but the total Maslov index of this configuration
      is at least smaller by $4$ compared to the initial one so that,
      by induction, generically this is not possible.

      We are now left to treat the case $n=\dim L =2$. In this case
      the condition $N_{L}\geq 4$ implies that at most two disk can
      appear in the configurations that are relevant to us and the
      combinatorial arguments used before (in \S
      \ref{S:transversality}) for this type of configuration suffice
      to conclude.
   \end{proof}

   We will return to the case $N_{L}\leq 3$ later but will now pursue
   with the proof of the identity~\eqref{eq:htpy} for $N_{L}\geq 4$.

   We define $$\xi(a\otimes b\otimes x)=\sum_{ y,\mathcal{T}}
   \#_{\mathbb{Z}_2} (\mathcal{P}_{\mathcal{T}}(a,b,x,y))y
   t^{\mubar(\mathcal{T})},$$
   where we recall that
   $\mubar(\mathcal{T})$ is $\tfrac{1}{N_L}$ of the total Maslov class
   of the tree $\mathcal{T}$ and the sum is taken over all elements so
   that $|a|+|b|+|x|-|y|+\mu(\mathcal{T})-4n+1=0$.

   As always, we will deduce formula~\eqref{eq:htpy} from the study of
   the boundary of the Gromov compactification of the space
   $\mathcal{P}_{\mathcal{T}}(a,b,x,y)$. Let
   $\overline{\mathcal{P}}_{\mathcal{T}}(a,b,x,y)$ denote this
   compactification. We now describe the types of configurations to
   $u\in\overline{\mathcal{P}}_{\mathcal{T}}(a,b,x,y)\backslash
   \mathcal{P}_{\mathcal{T}}(a,b,x,y)$ when $|a|+|b|+|x|-|y| +
   \mu(\mathcal{T}) - 4n + 1 = 1$.  These configurations are described
   in the same way as in the points (1) - (9) before except that,
   additionally one of the following conditions is satisfied:
   \begin{itemize}
     \item[(1')] Precisely one edge in $\mathcal{T}$ can be
      represented by a flow line of $0$ length.
     \item[(2')] Precisely one vertex in the tree $\mathcal{T}$ is
      represented by a cusp curve with two components - one a disk an
      the other a sphere.  In this case, the sphere carries two marked
      points corresponding to the incidence points of a flow lines of
      $h_{1}$ and one of $h_{2}$ and the disk carries an internal
      marked point (which is its junction with the sphere) which is
      situated on the real line joining $-1$ to $+1$. It is possible
      for one of the disk or the sphere to be a constant one as long
      as it is stable.
     \item[(3')] Precisely one vertex in the tree is represented by a
      cusp curve with two components both of which are disks. In this
      case each of the two components carry two boundary marked points
      (one of which is their junction point) and each of them might
      also carry some internal marked points situated on the real line
      joining $-1$ to $+1$. It is possible for one of the disks to be
      a constant one as long as it is stable.
     \item[(4')] Precisely, one edge in the tree is represented by a
      broken flow line of either $h_{1}$, $h_{2}$ or $f$.
   \end{itemize}

   In view of the descriptions of the compactifications of the various
   moduli spaces described earlier in the paper, this statement is
   rather obvious. The only point worth explicit mention here is that
   if a sphere bubbles off, then for dimension reasons, it necessarily
   has to carry three marked points (one of them being a junction
   point with a disk). In view of this, as the internal marked points
   on each disk are along the line joining $-1$ to $+1$ the only such
   possible bubbling fits in case (2').

   The techniques described earlier in the paper (in particular, the
   gluing results) insure that each of the configurations described at
   (1'), (3'), (4') is in fact a boundary point of the $1$-dimensional
   manifold with boundary
   $\overline{\mathcal{P}}_{\mathcal{T}}(a,b,x,y)$. An additional
   gluing argument applies for the configurations at (2') which is
   concerned with the gluing of a sphere to a disk in an internal
   point of the disk. However, this type of gluing is already well
   understood as it coincides, essentially, with the closed case (the
   gluing of two $J$-holomorphic spheres). See~\cite{McD-Sa:Jhol-2}
   for example.

   To prove relation~\eqref{eq:htpy} we now consider the sum $S$ of
   all these boundary points when $\mathcal{T}$ takes all possible
   values (and $a,b,x,y$ are fixed). We first notice that each of the
   terms described by (1'), (2') and (3') appears twice in this sum.
   Indeed, each cusp configuration also appears from the ``collapsing"
   of an edge in a tree $\mathcal{T}'$ with one more vertex and one
   more edge than $\mathcal{T}$.  Therefore the sum $S$, which
   vanishes, equals the sum of the terms of type (4').  We now
   consider the terms appearing in $S$ and notice that they are of the
   following types:
   \begin{itemize}
     \item[(i)] the broken flow line is associated to $f$ and is
      broken below the vertices of valence strictly greater than two.
     \item[(ii)] the broken flow line is associated to $h_{1}$ and is
      broken above all vertices of valence two.
     \item[(iii)] the broken flow line is associated to $h_{2}$ and is
      broken above all the vertices of valence two.
     \item[(iv)] the broken flow line is associated to $h_{1}$ and is
      broken below a vertex of valence three (which, obviously,
      corresponds to a $J$-sphere).
     \item[(v)] the broken flow line is associated to $f$ and is
      broken above all the vertices of valence strictly greater than
      two.
     \item[(vi)] the broken flow line is associated to $f$ and it is
      broken above some vertex of valence three and below some other
      vertex of valence three.
   \end{itemize}

   We now observe that the terms counted in (i) are precisely those
   appearing in $d\xi(a\otimes b\otimes x)$; the terms in (ii)
   correspond to those in $\xi ((da)\otimes b\otimes x)$; the terms in
   (iii) correspond to those in $(\xi a\otimes (db)\otimes x)$; the
   terms in (iv) correspond to those in $(a\ast b)\ast x$; the terms
   in (v) correspond to those in $\xi (a\otimes b\otimes (dx))$;
   finally, the terms in (vi) correspond to those appearing in $a\ast
   (b\ast x)$.  This concludes our proof when $N_{L}\geq 4$.

   We now give the argument for the case $2\leq N_{L}\leq 3$ under the
   assumption that that $n = \dim L \geq 3$. We will not give an
   explicit proof of the formula~\eqref{eq:htpy} in the case
   $N_{L}\leq 3$ and $n=2$.  However, we mention that if coefficients
   are taken in $\Lambda$ instead of $\Lambda^{+}$ then the module
   axiom for our exterior operation follows from the comparison with
   Floer homology which is explained in~\S\ref{Sb:comparison-HF}.

   Returning to the case $N_{L}\leq 3$ and $n\geq 3$ the general
   argument follows the lines of the proof above except that the
   perturbation techniques described in \S\ref{Sb:hampert} are also
   required. Indeed, in the definitions of the various moduli spaces -
   which are perfectly similar to the ones given before - we need that
   the disks carrying one or two marked points in their interior be
   perturbed disks in the sense of that section. All the arguments
   work in the same way as in the case treated above (for $N_{L}\geq
   4$) and are in fact simpler because in all replacement arguments as
   well as in the reduction to simple disks, the disks carrying marked
   points in their interior do not intervene (because transversality
   is insured for these perturbed disks by a generic choice of
   perturbation).  There is only one point where this proof requires a
   new argument.  Indeed, we consider the sum $S'$ which is defined as
   $S$ above except that the disks of valence strictly greater than
   $2$ are $(J,H)$-disks (for a small, generic $H$).  This sum $S'$
   again vanishes. However, besides the terms of type (4') which
   appear in $S'$ there are some other terms which contribute to this
   sum. These are all the configurations obtained by the bubbling off
   of a (perturbed) disk of valence four giving as result a cusp curve
   made of a $(J,H)$-disk of valence three which is joined to a
   genuine pseudo-holomorphic disk of valence three together with all
   the configurations obtained when an edge joining two $(J,H)$-disks
   of valence three reduces to $0$.  One point is crucial here: in the
   first type of configuration described here - which is obtained by
   bubbling off - the $(J,H)$-disk appears before the $J$-disk in the
   tree (this happens because, the first internal incidence point is
   the center of the disk and the second incidence point is in the
   interval $(0,+1)\subset D$) !

   Denote by $S''$ the number of configurations of these two types.
   Clearly, if $H=0$ these two types of configurations coincide and,
   by the usual gluing results, $S''=0$.  In general, another argument
   is needed to show that $S''=0$. For this, one considers the moduli
   space $\mathcal{P}_{tan, H, H'}(a,b,x,y)$ consisting of
   configurations of the usual type but with two disks of valence
   three which share an incidence point and so that the first disk (in
   the order of the tree) is a $(J,H)$-disk and the second disk is a
   $(J,H')$-disk. The result follows if we can show that, when the
   virtual dimension of the space $\mathcal{P}_{tan,H,H'}(a,b,x,y)$ is
   equal to $0$, for all sufficiently small $H$ and $H'$ these moduli
   spaces are constituted by simple, absolutely distinct disks and the
   number of elements in $\mathcal{P}_{tan,H,H'}(a,b,x,y)$ equals that
   of $\mathcal{P}_{tan,0,0}(a,b,x,y)$. It is not difficult to check
   that this result follows by a usual cobordism argument similar to
   Proposition \ref{P:n-pert=no-pert} if we prove that
   $\mathcal{P}_{tan,0,0}(a,b,x,y)$ is made of simple, absolutely
   distinct disks for generic almost complex structures. In turn the
   key point for this is that, as we only need this result when the
   virtual dimension of the respective moduli space is $0$, we can
   reason as in Proposition \ref{P:qm-axy}. In other, words assuming
   that the disks are not all simple and absolutely distinct we may
   construct a new configuration of strictly lower Maslov index (which
   leads to a drop of virtual dimension of at least $2$ because
   $N_{L}\geq 2$) but the condition on one of the marked points lying
   in the center of one of the disks of valence three might be lost
   which leads to a potential increase in virtual dimension by $1$.
   As a whole, we still obtain a configuration belonging to a moduli
   space of virtual dimension at most $-1$ which leads to a
   contradiction because, inductively, this moduli space is assumed to
   be void.

   \

   It remains to prove point~(\ref{I:mod-unit}) of
   Proposition~\ref{prop:module_str}. Let $h:M\to \R$ be a Morse
   function with a single maximum $\Theta$. Let $f: L \to \mathbb{R}$
   be a Morse function and $x \in \textnormal{Crit}(f)$. Then we have
   (at the chain level): $\Theta \ast x=x\in \mathcal{C}^{+}(L;f,J)$.
   The reason for this is that, because $\Theta$ is a maximum, the
   moduli spaces $\mathcal{P}(\Theta,x,y;\la)$ which define the module
   operation as described in~\S\ref{Sb:external-op} can only be
   $0$-dimensional and non-void if $x=y$, $\la=0$ and, in fact, the
   ``pearls" appearing in these configurations are constant equal to
   $x$.

   Since $[M]=[\Theta]$ point~(\ref{I:mod-unit}) follows immediately.
\end{proof}

\subsubsection{Algebra structure.} \label{Sb:algebra-struct}
In this section we will show the following identity:

\begin{prop}\label{prop:alg}
   For any $x,y\in Q^{+}H_{\ast}(L)$, $a\in
   Q^{+}H(M)$ we have $$a\ast (x\ast y)=(a\ast x)\ast y~.~$$
\end{prop}

\begin{proof}
   The argument is very similar to the one used for Proposition
   \ref{prop:module_str}.  We will again make use of moduli spaces
   modeled on trees with three entries and one exit so that the first
   entry is a critical point of a Morse function $h:M\to \R$ but the
   second and third entries are both critical points of Morse
   functions, $f_{1}$ and $f_{2}$, on $L$. We will also assume that
   $h,f_{1},f_{2}$ are all in generic position and that $f_{1}$ and
   $f_{2}$ have the same critical points. The moduli spaces needed for
   our argument are similar to those used to prove the associativity
   of the quantum product, in particular, in the proof of Lemma
   \ref{lem:prodassoc}. For $a\in\Crit(h)$, $y,z,w\in \Crit(f_{1})$
   and $\mathcal{T}$ a planar tree as there we will denote by
   $\mathcal{P}_{\mathcal{T}}(a,y,z,w)$ the relevant moduli spaces.
   They are defined as the moduli spaces
   $\mathcal{P}_{\mathcal{T}}(-,-,-,-)$ in that Lemma except that the
   labeling of the edges is by the symbols $0$, $1$, $2$ and is so
   that there is a single edge labeled by $0$ which corresponds to a
   flow line of $-\nabla h$, the edges labeled by $1$ and $2$
   correspond, respectively, to flow lines of $-\nabla f_{1}$ and
   $-\nabla f_{2}$. The various incidence conditions are as follows.
   At the vertices of valence $2$ the incidence points are $+1$ and
   $-1$ and the entry edge is labeled by the same symbol as the
   exiting edge; at a vertex of valence three if all the incidence
   points are on the boundary, then they are the roots of order three
   of the unity (as in Lemma \ref{lem:prodassoc}), the entrance edges
   are labeled by $1$ and $2$ and the exit edge is labeled by $1$; a
   vertex of valence three might also have only two incidence points
   on the boundary and one in the interior - in this case the edge
   arriving in the interior is labeled by $0$ and the other entering
   edge is labeled by $1$ as well as the exiting edge, the incidence
   relations in this case are as in \S \ref{Sb:external-op}; finally a
   vertex of valence four has three marked points on the boundary so
   that two correspond to $-1$ and $+1$ and there is an entering edge
   labeled by $1$ arriving at $-1$ and the exiting edge is again
   labeled by $1$ and is attached at $+1$, there is an additional
   incidence point $\theta\in (0,\pi)\subset S^{1}$ where is attached
   an entering edge labeled by $2$, the fourth incidence point is
   $0\in D$ and is the arrival point of the edge labeled by $0$.

   The virtual dimension of these moduli spaces is
   $|a|+|y|+|z|-|w|+\mu(\mathcal{T})-3n+1$ and they are used to define
   a chain homotopy
   $$\xi':C_k(h) \otimes \mathcal{C}^{+}_q(f_{1}) \otimes
   \mathcal{C}^{+}(f_{2})_{p}\to \mathcal{C}^{+}_{k+q+p-3n+1}(f_{1})$$
   between $(- \ast -) \ast -$ and $-\ast (-\ast -)$.

   The argument goes again in two stages: first we need to show that
   the necessary

   transversality is satisfied for such moduli spaces of virtual
   dimension at most $1$ and, as a second step, an appropriate
   boundary formula needs to be justified.  This second step is
   essentially identical with the one used to prove
   formula~\eqref{eq:htpy} so we will leave it to the reader. The
   first step reduces to showing that
   $$
   \mathcal{P}_{\mathcal{T}}(a,y,z,w) =
   \mathcal{P}_{\mathcal{T}}^{\ast,d}(a,y,z,w)$$
   where, as always, the
   superscript $\ast,d$ indicates those configurations which consist
   of simple absolutely disjoint disks.  To show this equality the
   procedure is again as before in the paper. Without the use of any
   perturbations it is easy to treat the case $N_{L}\geq 3$ and $n\geq
   3$. In case $N_{L}=2$, $n\geq 3$ the use of perturbations is
   necessary. Finally, for $N_{L}=2$, $n=2$ the argument is much more
   combinatorial and we will not give it - rather we refer to the
   section concerning the comparison with Floer homology: the relation
   for $QH(-)$ (but not for $Q^{+}H(-)$) follows in general from
   there.
\end{proof}

 \subsubsection{Two-sided algebra structure.}
 The purpose in this paragraph is to prove:

\begin{prop}\label{prop:two-sided}
   For any $a\in H^{+}Q(M)$ and $x,y\in H^{+}Q(L)$ we have:
   $$a\ast (x\ast y) = x\ast (a\ast y)~.~$$
\end{prop}

\begin{proof}
   The proof is similar to the one in the last subsection.  We will
   use the same type of moduli spaces - now denoted by
   $\mathcal{P}_{\mathcal{T}}(y,a,z,w)$ - except for a couple of very
   simple modification.
\begin{itemize}
  \item[i.] We require the vertices of valence three with three
   incidence points on the boundary to correspond to edges which, in
   the cyclic order, verify: the first is labeled by $2$ and is an
   entry edge, the second is an exit edge and is labeled by $1$, the
   third in an entry edge and is labeled by $1$.
  \item[i.] For a vertex of valence four we require the incidence
   points to be $-1$, $+1$ and $0$ and to satisfy the same properties
   as in the proof of Proposition \ref{prop:alg} but the fourth
   incidence point which again corresponds to an entry edge labeled
   by $2$ is now a point in $(\pi,2\pi)\subset S^{1}$.
\end{itemize}

Using these moduli spaces, the usual arguments, applied as in the last
section produce a chain homotopy
$$\xi'' :C(h)_{k}\otimes
\mathcal{C}^{+}(f_{2})_{q}\otimes\mathcal{C}^{+}(f_{1})_{p}\to
\mathcal{C}^{+}(f_{1})_{k+q+p-3n+1}$$
between the two order three
products which appear in the statement.
\end{proof}

\subsection{Quantum inclusion.} \label{Sb:q-inclusion}
There is a canonical map:
\begin{equation} \label{eq:quantum_incl}
   i_{L}:Q^{+}H_{\ast}(L)\to
   Q^{+}H_{\ast}(M)
\end{equation}
which is defined at the level of chain complexes by:
\begin{equation} \label{Eq:iL}
   \begin{gathered}
      i_{L}:\mathcal{C}^{+}_{k}(L;f,J)\to (\mathbb{Z}_2 \langle \Crit
      (h) \rangle \otimes \La^{+})_{k}, \\
      i_{L}(x) = \sum_{\substack{a\in \textnormal{Crit}(h) \\
          \la,\mathcal{T}}} \#_{\mathbb{Z}_2}
      \mathcal{P}_{\mathcal{T}}(x,a)a t^{\mubar(\mathcal{T})},
   \end{gathered}
\end{equation}
where $\mathcal{P}_{\mathcal{T}}(x,a)$ is again a moduli space modeled
on a tree which will be linear here (that is a tree with one entry and
one exit) which we will describe in more detail below.  The sum is
taken over all such trees $\mathcal{T}$ with
$|x|-|a|+\mu(\mathcal{T})=0$.

The first part of point iii. of Theorem \ref{thm:alg_main} comes down
to:
\begin{prop}\label{prop:inclusion}
   The map in equation~\eqref{eq:quantum_incl} is well defined and it
   verifies:
   $$i_{L}(a\ \ast\ x) = a\ast i_{L}(x)\ , \forall \ a\in
   Q^{+}H_{\ast}(M)\ ,\ x\in Q^{+}H_{\ast}(L)~.~$$
\end{prop}

Here is now the more explicit description of the moduli space
$\mathcal{P}_{\mathcal{T}}(x,a)$.  It consists of configurations that
are similar to those lying in the moduli spaces of pearls - as in
Definition \ref{def:moduli1} - except that the condition imposed to
the last disk in the string, $u_{k}$, is modified so that
$\gamma'_{+\infty}(b_{k})=a$, where $b_{k}=u_{k}(0)$ and $\gamma'$ is
the negative gradient flow of $h$, and the disk $u_{k}$ may also be
constant.  In other words, the ``exit" incidence point on the last
disk is the interior point $0\in D$ (instead of $+1\in \partial D$)
and the exit edge is a negative gradient line of $h$ which ends in
$a\in \Crit(h)$. When $n\geq 3$, the methods described before
immediately show that with the previous definition $i_{L}$ is a chain
map. When $n=2$ a simple combinatorial argument suffices. It is also
easy to see that, in homology, this map does not depend of the choices
made in its construction.  It is useful to note that, as on the last
disk $u_{k}$ there are just two marked points - one interior and one
on the boundary, the use of perturbations is not necessary in these
arguments.

\

We now justify the equation in Proposition~\ref{prop:inclusion}. This
is based on the construction of an appropriate chain homotopy. In
turn, this depends on defining yet other moduli spaces
$\mathcal{P}_{\mathcal{T}}(a,y,b)$. Here $a, b\in\Crit(h)$ and $y\in
\Crit(f)$ and $\mathcal{T}$ is a tree with two entries and one exit.
As for the moduli spaces $\mathcal{P}_{-}(x,a)$ above, the key point
is that the last edge in the tree is again a flow line of the negative
gradient of $h$ which arrives in $b$. There is also an edge which
corresponds to a negative gradient flow line of $h$ which leaves from
$a$. All the vertices of valence at least two correspond to $J$-disks
with the possible exception of a single one which is of valence three
and may correspond to a $J$-sphere. We now describe the incidence
relations.  There are two types of vertices of valence two: $J$-disks
with the two marked points on their boundary - in this case the
incidence points are $-1$ which is the entry point of a flow line of
$-\nabla f$ and $+1$ which is the exiting point of a flow line of
$-\nabla f$; a $J$-disk with one marked point on the boundary $-1$
which is an entry point of a flow line of $-\nabla f$ and an interior
marked point $0$ which is the exiting point of a flow line of $-\nabla
h$. There also are three possibilities for the vertex of valence $3$.
In the first, this vertex corresponds to a $J$-disk and the incidence
points are $-1$, $0$, $+1$ so that $-1$, $0$ are entry points for flow
lines of respectively, $-\nabla f$, $-\nabla h$ and $+1$ is an exiting
point of a flow line of $-\nabla f$. In the second case, again the
vertex corresponds to a $J$-disk but this time the marked points are
$-1,p,0$ where $p\in (-1,0)$ and $-1$ is again an entry point for a
flow line of $-f$, $p$ is an entry point for a flow line of $-h$ and
$0$ is an exit point for a flow line of $-h$. Finally, in the third
case the vertex in question corresponds to a $J$-sphere the marked
points in that case are roots of order three of the unity so that the
first root is an entrance point of a flow line of $-h$ originating in
$a$, the second is an exiting flow line of $-h$ arriving in $b$ and
the third is an entering flow line of $-h$ (whose origin is,
necessarily, at the center of a $J$-disk).  Trivial disks and spheres
can appear in these configurations as long as they are stable. It is
easy to see that the virtual dimension of these moduli spaces is
$|x|+|a|-2n+\mu(\mathcal{T})-|b|+1$. As always we will only need to
use those moduli spaces of virtual dimension $0$ and $1$ and the chain
homotopy in question is defined by
$$\xi''(a\otimes x)=\sum \#_{\mathbb{Z}_2}
(\mathcal{P}_{\mathcal{T}}(a,x,b))b t^{\mubar(\mathcal{T})}$$
where the
sum is over all the trees $\mathcal{T}$ so that the virtual dimension
is $0$. Of course, this definition is only valid generically: in that
case, the methods described earlier in the paper can be easily applied
here to show that $\mathcal{P}_{\mathcal{T}}(a,x,b)$ is a
$0$-dimensional compact manifold. Finally, to show that:
$$(d \xi'' + \xi'' d)(a\otimes x)=i_{L}(a\ast x)-a\ast i_{L}(x)$$
the
regularity of the moduli spaces of dimension $1$ is needed together
with a boundary description. This follows again by the methods
described earlier in the paper. However, notice that the use of
perturbations is necessary in this case.

The last step to conclude point iii of Theorem \ref{thm:alg_main} is
to show relation~\eqref{eq:inclusion_mod}:
$$<h^{\ast},i_{L}(x)>=\epsilon_{L}(h\ast x)$$
for all $h\in
H_{\ast}(M)$ and $x\in Q^{+}H_{\ast}(L)$. First denote by $m$ the
minimum of $f$ (we assume it is unique to simplify the discussion).
Notice that there is a bijection
$$b: \mathcal{P}_{\mathcal{T}}(x,a)\to
\mathcal{P}_{\mathcal{T}'}(a,x,m)$$ where
$\mathcal{P}_{\mathcal{T}'}(a,x,m)$ are the moduli spaces of the type
used in the definition of the module operation in \S\ref{S:qm}
associated to the function $f$ and its critical points $x$ and $m$ and
{\em to the function} $-h$ together with its critical point $a$; the
tree $\mathcal{T}'$ is obtained by inverting the last edge in the tree
$\mathcal{T}$ and adding one edge going from the last disk in
$\mathcal{T}$ to the critical point $m$. This also shows how to define
the bijection $b$: for each configuration $v$ in
$\mathcal{P}_{\mathcal{T}}(x,a)$ we consider the last disk $u_{k}$ of
$v$ and the point $u_{k}(+1)$.  This point belongs to the unstable
manifold of $m$ (generically, as always) and so, by replacing the last
edge in $v$ by a negative gradient flow line of $-h$ which goes from
$a$ to $u_{k}(0)$ and adding to $v$ one edge given by a negative
gradient flow line of $f$ which joins $u_{k}(+1)$ to $m$ we obtain
$b(v)$.  It is then clear that this map $b$ so defined is bijective.
Combining this with the definition of $\epsilon_{L}$ the
relation~\eqref{eq:inclusion_mod} follows. \Qed

\subsection{Spectral sequences.} \label{Sb:spect-seq}
We now notice that all the structures defined above are compatible
with the degree filtration so that the point iv. of Theorem
\ref{thm:alg_main} is trivial.

\subsection{Comparison with Floer homology.} \label{Sb:comparison-HF}

We deal here with the last point of Theorem \ref{thm:alg_main}.
Recall from the statement of Theorem \ref{thm:alg_main} that we denote
by
$$\mathcal{C}(L;f,\rho,J)=\mathcal{C}^{+}(L;f,\rho,J)
\otimes_{\La^{+}}\La~.~$$

\subsubsection{Comparing complexes.} \label{sbsb:comparing}
The version of Floer homology which we need is defined with the help
of an auxiliary Hamiltonian $H:M\times [0,1]\to M$ and its
construction is standard (see~\cite{Oh:HF1}).

Put $H_t(x) = H(x,t)$. Denote
$$\mathcal{P}_{0}(L,L)=\{\gamma:[0,1]\to M : \gamma(0)\in L,\
\gamma(1)\in L\ \ \gamma\simeq\ast\}~.~$$

The generators of the Floer complex, $CF(L;H,J)$, are elements of
$\mathcal{P}_{0}(L,L)$ which are orbits of the Hamiltonian vector
field $X^H_t$ (defined by $\omega(X^H_t,Y)=-dH_t(Y)$). We denote this
set of orbits by $I_{H}$. There is a natural map $\pi_{2}(M,L)\to
\pi_{1}(\mathcal{P}_{0}(M,L),*)$ and extensions to
$\pi_{1}(\mathcal{P}_{0}(M,L),*)$ of each of the maps $\omega$ and
$\mu$ (it is easy to see that these extensions continue to be
proportional in this case). We need to fix a base point, $\bar{\eta}$,
in the space $\bar{\mathcal{P}}_{0}(L,L)$ which is the abelian cover
of $\mathcal{P}_{0}(L,L)$ which is associated to the kernel of $\mu$.
Once this choice is made the Floer complex is defined by
$$CF(L;H,J)=\mathbb{Z}_2 \langle I_{H} \rangle \otimes \La~.~$$

The differential is defined by first fixing one lift
$\bar{x}\in\bar{\mathcal{P}}_{0}(L,L)$ for each of the elements $x\in
I_{H}$ and there is a Maslov index $\mu(\bar{x}), 2\mu(\bar{x})\in \Z$
which is well defined. This determines a grading on $CF(L;H,J)$ by the
formula $|x\otimes t^r|= \mu(\bar{x})-r N_L$. The differential is
defined by
\begin{equation}\label{eq:differential}
   dx=\sum_{y,\la}\#_{\mathbb{Z}_2} \mathcal{M}(\bar{x},\bar{y};\la)
   y t^{\mubar(\lambda)},
\end{equation}
where $x,y\in I_{H}$, $\la \in \pi_{1}(\mathcal{P}_{0}(M,L),*) / \ker
\mu$ and $\mathcal{M}(\bar{x},\bar{y};\la)$ is described as follows.
First let $\mathcal{M}'(\bar{x},\bar{y};\la)$ be the moduli space of
paths $\bar{u}:\R\to \bar{\mathcal{P}}_{0}(L,L)$ so that
$\bar{u}(-\infty)=\bar{x}$, $\bar{u}(+\infty)=\la \bar{y}$ and the
projection $u$ of $\bar{u}$ to $\mathcal{P}_{0}(L,L)$ verifies
\begin{equation}\label{eq:floer}
   \partial u/\partial s + J\partial u/\partial t +\nabla H_t (u)=0,
\end{equation}
where $\nabla H_t$ is defined with respect to the metric $g_{J}$
associated to $(\omega, J)$. There is an obvious action of $\R$ on
this moduli space and we let
$$\mathcal{M}(\bar{x},\bar{y};\la)=
\mathcal{M}'(\bar{x},\bar{y};\la)/\R~.~$$
Generically, this is a
manifold of dimension
$$\mu(\bar{x})-\mu(\bar{y})+\mu(\la)-1$$
and the sum in~\eqref{eq:differential} is taken over all of those
$\bar{y},\la$ so that the respective moduli space is $0$-dimensional.
The homology of this complex is the Floer homology of $L$, $HF(L)$,
and does not depend on $H$ and $J$. It depends on the choice of the
base point $\bar{\eta}$ up to translation.

\begin{rem}\label{rem:coefficients_Floer}
   We leave it to the reader to verify that using the positive Novikov
   ring in this case is not possible.
\end{rem}

The point ii. of Theorem \ref{thm:alg_main} is a consequence of the
following:

\begin{prop}\label{prop:comparison} For generic $(f,\rho,H,J)$
   there are chain morphisms $$\psi:\mathcal{C}(L;f,\rho,J)\to CF(L;
   H,J)$$
   and $$\phi:CF(L; H,J)\to \mathcal{C}(L;f,\rho,J)$$
   which
   induce canonical isomorphisms in homology. These induced maps are
   inverse one to the other.
\end{prop}

\begin{rem} These morphisms are constructed by the Piunikin-Salamon-Schwarz
   method and, indeed, they are the exact Lagrangian counterpart of
   the PSS morphisms.  Such morphisms have been discussed in the
   Lagrangian setting - when bubbling is avoided - in
   \cite{Bar-Cor:NATO, Alb:PSS, Kat-Mil:PSS} and, in the general
   cluster setting, in \cite{Cor-La:Cluster-1}.
\end{rem}
\begin{proof}[Proof of Proposition~\ref{prop:comparison}]
   Notice that, given $\gamma\in I_{H}$ we may view each element
   $\bar{\gamma}\in \bar{\mathcal{P}}_{0}(L,L)$ which covers $\gamma$
   as a pair $(\gamma, u)$ where $u$ is a ``half-disk" capping
   $\gamma$. For our fixed lifts $\bar{\gamma}$ of the orbits
   $\gamma\in I_{H}$ we denote $\bar{\gamma}=(\gamma,u_{\gamma})$.

   To define the morphisms $\psi$ and $\phi$ some new moduli spaces
   are necessary. Given $x\in \Crit(f)$ and $\gamma\in I_{H}$ we will
   define next the moduli spaces $\mathcal{P}_{\mathcal{T}}(x,\gamma)$
   and $\mathcal{P}_{\mathcal{T}}(\gamma,x)$.  In both cases
   $\mathcal{T}$ is a linear tree as those used in the definition of
   the pearl moduli spaces. Compared to the definition of the pearl
   moduli spaces - Definition \ref{def:moduli1} - the configurations
   assembled in $\mathcal{P}_{\mathcal{T}}(x,\gamma)$ have the
   property that the last vertex in the chain does not correspond to a
   $J$-disk but rather to an element $u_{k}:\R\times [0,1]\to M$ so
   that we have:
   $$u_{k}(\mathbb{R} \times\{0,1\})\subset L\ , \
   \partial_{s}(u_{k})+ J(u_{k})\partial_{t}(u_{k})+\beta(s)\nabla
   H_t(u)=0\ , u_{k}(+\infty)=\gamma$$
   and
   $r_{t_{k}}(u_{k-1}(+1))=u_{k}(-\infty)$ where $r_{t}$ is the
   negative gradient flow of $f$ and $\beta$ is a smooth cut-off
   function which is increasing and vanishes for $s\leq 1/2$ and
   equals $1$ for $s\geq 1$. The virtual dimension of this moduli
   space is $|x|-|\bar{\gamma}|+\mu(\mathcal{T})-1$ where
   $|\bar{\gamma}|=\mu(\bar{\gamma})$ and $\mu(\mathcal{T})$ is by
   definition the sum of the Maslov indices of the disks corresponding
   to the vertices appearing in the tree $\mathcal{T}$ to which we add
   the Maslov class of the disk obtained by gluing $u_{k}$ to
   $u_{\gamma}$ along $\gamma$.  Similarly, the moduli space
   $\mathcal{P}_{\mathcal{T}}(\gamma,x)$ has an analogue definition
   except that the end conditions are reversed. More precisely, the
   first disk $u_{0}$ is the solution of an equation:
   $$u_{0}(R\times\{0,1\})\subset L\ , \ \partial_{s}(u_{0})+
   J(u_{0})\partial_{t}(u_{0})+\xi(s)\nabla H(u,t)=0\ ,
   u_{0}(-\infty)=\gamma$$, $r_{-t_{1}}(u_{1}(-1))=u_{0}(+\infty)$
   where $\xi$ is a smooth curt-off function which is decreasing and
   vanishes for $s\geq 1/2$ and equals $1$ for $s\leq 0$. The same
   methods as those used earlier in the paper show that for a generic
   $J$ we have an equality
   $$\mathcal{P}_{\mathcal{T}}(x,\gamma) =
   \mathcal{P}_{\mathcal{T}}^{\ast,d}(x,\gamma)$$
   (and similarly for
   the moduli spaces $\mathcal{P}_{\mathcal{T}}(\gamma,x)$) where the
   moduli space on the left consists of configurations containing only
   simple, absolutely distinct disks. For $H$ generic the appropriate
   transversality of the evaluation maps can be achieved when the
   virtual dimension is at most $1$.  The definition of the morphism
   $\psi$ is now as follows:
   $$\psi
   (x)=\sum_{\gamma,\mathcal{T}}\#(\mathcal{P}_{\mathcal{T}}(x,\gamma))
   \gamma t^{\mubar(\mathcal{T})}$$
   where the sum is taken over all
   trees so that the dimension of the respective moduli spaces is $0$.
   Similarly,
   $$\phi(\gamma) =
   \sum_{x,\mathcal{T}}\#(\mathcal{P}_{\mathcal{T}}(\gamma ,x))x
   t^{\mubar(\mathcal{T})}$$
   where again we only take into account
   moduli spaces of dimension $0$. By analyzing the boundary of the
   compactifications of the moduli spaces
   $\mathcal{P}_{\mathcal{T}}(x,\gamma)$ and
   $\mathcal{P}_{\mathcal{T}}(\gamma,x)$ it is easy to show that both
   $\phi$ and $\psi$ are chain morphisms. We then need to show that
   the compositions $\phi\circ \psi$ and $\psi\circ \phi$ are both
   chain homotopic with the respective identities. As in the
   non-bubbling case, this proof is based on a gluing argument - which
   allows to view the product
   $\mathcal{P}_{\mathcal{T}}(x,\gamma)\times
   \mathcal{P}_{\mathcal{T'}}(\gamma,y)$ as part of the boundary of a
   moduli space $\mathcal{P}_{\mathcal{T\# T'}}(x,y,k,H)$ where this
   last moduli space is modeled on the tree $\mathcal{T\# T'}$ which
   is obtained by gluing $\mathcal{T}'$ at the end of $\mathcal{T}$
   and it consists of pearl-like objects except that the $k$-th disk
   verifies a perturbed Floer type equation of the form
   $\partial_{s}(u_{0})+ J(u_{0})\partial_{t}(u_{0})+\nu_{R}(s)\nabla
   H(u,t)=0$ where $\nu_{R}$ is a family of smooth functions so that
   $\nu_{R}(s)=0$ for $|s|>R$, is increasing for $s<0$ and decreasing
   for $s>0$ and, for $R$ sufficiently big, it is equal to $1$ for
   $|s|< R-1$.

\end{proof}

\subsubsection{Module action and internal product on $HF(L)$.}
In this subsection we want to notice that there exists a natural
action of $QH(M)$ on $HF(L)$ which is identified via the PSS maps with
the action discussed in \S\ref{S:qm}.  The definition of this module
structure is completely similar to the $\cap$- action of singular
homology on Hamiltonian Floer homology as it is described for example
in~\cite{Fl:fp-hol-spheres} or in~\cite{Sc:action-spectrum}.
Similarly, we also have the ``half"-pair of pants product on $HF(L)$.

\

Given a Morse function $h:M\to \R$ together with a Riemannian metric
$\rho_M$ so that the pair $(h,\rho_M)$ is Morse-Smale as in
\S\ref{S:qm} we define an operation:
\begin{equation}\label{eq:operation}
   \ast : C_{k}(h)\otimes CF_{l}(L;H,J)\to CF_{l+k-2n}(L;H,J)
\end{equation}as follows
$$a\ast x=\sum_{y,\la}\#\mathcal{M}(a;\bar{x},\bar{y};\la)y
t^{\mubar(\la)}~.~$$
Here
\begin{equation}
   \label{eq:moduli_cap}\mathcal{M}(a;\bar{x},\bar{y};\la)=\{(u,p)\in
   \mathcal{M}'(\bar{x},\bar{y};\la)\times W_a^u  \ : \
   u(0,1/2)=p\}
\end{equation}
and $W_a^u$ is the unstable submanifold of the critical point $a$ for
the flow of $-\nabla h$. As always, the sum above is understood to be
taken only when the moduli spaces in question are finite.

\

The Floer intersection product is an associative operation with unit
$$HF_{p}(L)\otimes HF_{q}(L)\to HF_{p+q-n}(L)~.~$$

We first recall that this product is defined by a chain map:
$$\ast : CF_{k}(L;H, J)\otimes CF_{l}(L;H',J)\to CF_{k+l-n}(L;H'',J)$$
where $H'$ and $H''$ are small deformations of $H$ and
$$x\ast y=\sum_{z,\la} \#(\mathcal{M}^{x,y}_{z}(\la))z
t^{\mubar(\la)},$$
where the moduli spaces $\mathcal{M}^{x,y}_{z}(\la)$
consist of semi-pants (or half pants) with their boundaries on $L$ and
which are otherwise similar to the usual pair of pants used to define
the standard product in Hamiltonian Floer homology.

Again, both structures are well defined for generic choices of data
and independent of these choices and the PSS maps identify them in
homology with the ones described in \S \ref{S:qm} and \S
\ref{subsec:product}. We leave this verification to the reader. It is
also easy to verify that analogue statements are valid for the map
$i_{L}$ and the spectral sequence compatibility.

\subsection{Duality} \label{subsec:proof_duality}
The purpose of this section is to prove the Corollary
\ref{cor:duality}. We start by recalling the conventions and the
notation in \S\ref{susubsubsec:duality}. Fix generic $f,\rho,J$ so
that the pearl complex $\mathcal{C}^{+}(L;f,\rho,J)$ is defined.  Let
$\mathcal{C}^{\odot}=\hom_{\La}(\mathcal{C}(L;f,\rho,J),\La)$ endowed
with the differential which is the adjoint of the differential of
$\mathcal{C}^{+}(L;f,\rho,J)$ and with the grading $|x^{\ast}|=-|x|$.
Recall also that $s^{k}\mathcal{C}$ indicates the $k$-th suspension of
the complex $\mathcal{C}$.

\emph{Proof of Corollary \ref{cor:duality}}.  The first step is to
notice that the pearl complex $\mathcal{C}^{+}(L;-f,\rho,J)$ is also
defined and there is a basis preserving isomorphism $i$ between
$\mathcal{C}^{+}(L;-f,\rho,J)$ and $s^{n}\mathcal{C}^{\odot}$ defined in
the same way as in the Morse case (see Remark \ref{rem:duality} a): it sends
each generator represented by a critical point $x$ of $f$, $ind_{f}(x)=k$
to a generator of degree $n-k$ represented by the same critical point $x$
only viewed as critical point of $-f$ .

At the same time, there is a comparison chain morphism
$$\phi:\mathcal{C}^{+}(L;f,\rho,J)\to \mathcal{C}^{+}(L;-f,\rho,J)$$
as in
\S\ref{subsubsec:inv} which induces a (canonical) isomorphism in
homology. Therefore we get a morphism
$$\eta: \mathcal{C}^{+}(L;f,\rho,J)\to s^{n}\mathcal{C}^{\odot}\ , \ \eta=i\circ\phi$$
which
induces an isomorphism in homology and this concludes the first part
of the Corollary.  For the second part we need to show that the
bilinear map $\tilde\eta$ associated to $\eta$ coincides with the
reduced quantum product $Q^{+}H(L)\otimes Q^{+}H(L)
\stackrel{*}{\longrightarrow}
Q^{+}H(L)\stackrel{\epsilon_{L}}{\longrightarrow}\La$.

The first step for this is to notice that it is enough to work with a
comparison morphism $\phi : \mathcal{C}^{+}(L;f',J)\to
\mathcal{C}^{+}(L;-f,J)$ where $f'$ and $f$ are in generic position. In fact, it is clear that
we may even replace $\phi$ with any other chain morphism which is chain homotopic to it.
There is a specific such chain morphism which will be useful in the proof.
Its definition is quite general so we now take $h$ another generic Morse function
and we describe this new
 comparison morphism $\phi':\mathcal{C}^{+}(L;f'J)\to
\mathcal{C}^{+}(L;h,J)$. It will be used in our proof for $h=-f$.
The construction of $\phi'$ is based on counting the elements of
certain moduli spaces $\mathcal{P}^{!}(f',h; x,y;J,\la)$ which are modeled on
linear trees, as the pearl moduli spaces, except that there is an
additional marked point on the tree which  is placed in the interior of an edge
in the tree. The key property of these moduli spaces is that all
the edges (or segments of edges) above this marked point
correspond to negative flow lines of $f'$ and all the edges (or
segments) below this marked point - in the tree - correspond to negative flow
lines of $h$.  The virtual dimension of these moduli spaces is $|x|-|y|+\mu(\la)$
where $\la$ is the total homotopy class of the configuration.
\begin{figure}[htbp]
   \begin{center}
      \includegraphics[scale=.25]{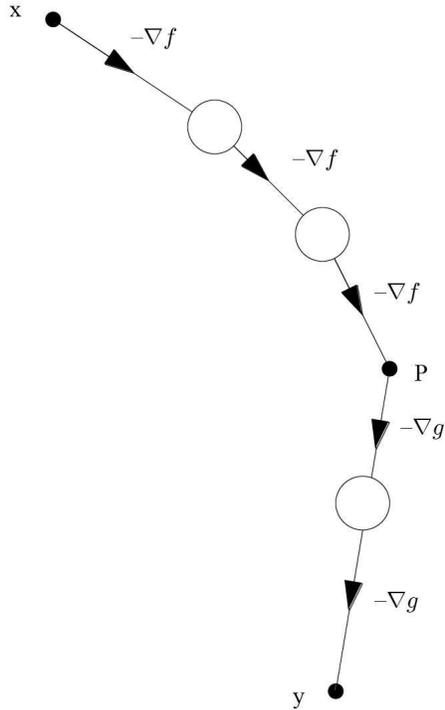}
   \end{center}
   \caption{ $x\in \Crit(f')$, $y\in \Crit(h), \ P$ is the new marked point.}
   \label{f:comp2}
\end{figure}
The same methods used earlier in the paper show that, with generic choices of defining data,
 counting the elements in the $0$-dimensional such moduli spaces  does indeed define a chain morphism.
The proof now consists of two steps. The first, is to remark that
$\phi$ and $\phi'$ are chain homotopic. We will postpone this argument
and proceed to describe the last step in the proof. For this we first fix a third Morse
function $f'':L\to \R$ which has a single
minimum $m$ and is in generic position with respect to $f'$ and $-f$.
We now notice, by reviewing the definition of the moduli spaces
$\mathcal{P}(x,y,m; f',f,f'')$ from \S\ref{subsec:product} that they
coincide in dimension $0$ with the moduli spaces
$\mathcal{P}^{!}(f', -f; x,y;J)$.  Indeed, if $\mathcal{P}(x,y,m;f',f,f'')$
is of dimension $0$, then the elements of this moduli space have the property
that their unique vertex of valence three is a constant disk and, moreover,
the chain of pearls associated to $f''$ and arriving in $m$ is a
single flow line of $-\nabla f''$ joining the vertex of valence three to $m$.
But  as $\mathcal{P}(x,y,m;f',f,f'')$ are precisely the moduli spaces which compute
$\epsilon_{L}\circ (-\ast-)$ we obtain $\epsilon_{L}(x\ast y)=<\phi'(x),y>=(i\circ\phi'(x))(y)$
which implies our claim.

To conclude our proof, we need to justify that
$\phi'$ and $\phi$ are chain homotopic. Recall that $\phi$ is defined
by making use of a Morse homotopy between $f'$ and $h$. Let $H:L\times
[0,1]\to \R$ be such a Morse homotopy.  We intend to apply a method similar to the
proof of the fact that the comparison morphism is canonical in homoloy as in \S\ref{subsubsec:inv}.
However, we apply that method in a more general situation: the place of the second Morse homotopy
$H'$ will be taken by the discontinous function $H'(x,t)=f'(x)$ for $t\in [0,1/2)$ and $H'(x,t)=h(x)$ for
$t\in [1/2,1]$. The reason that this method still works is that, despite the discontinuity of $H'$,
the vector field $\nabla H'$ is still well defined on $L\times [0,1]$ and
its negative flow lines are still well defined and continuous - they follow the negative
gradient of $f'$ for $t\in [0,1/2)$ and then follow along the negative gradient
of $h$ for $t\in [1/2, 1]$.    The only additional ingredient with respect to the method
described in \S\ref{subsubsec:inv} is a gluing statement. In essence, this is already seen in the
case of the purely Morse theoretic version of our statement (that is if no pseudo-holomorphic disks
are present): each  trajectory of $-\nabla H'$  has to be shown to be precisely the end of a one parametric
family of flow lines of the Morse homotopy relating $H$ to $H'$. We leave this last step as exercise for the
reader.
\qed

\subsection{Action of the symplectomorphism group.} \label{Sb:act-symp}

The purpose of this section is to prove Corollary \ref{cor:symm}.

\emph{Proof of Corollary \ref{cor:symm}.}  We assume that $\phi:L\to
L$ is a diffeomorphism which is the restriction to $L$ of the
symplectomorphism $\bar{\phi}$ and $f,\rho,J$ are such that the chain
complex $\mathcal{C}^{+}(L;f,\rho,J)$ is defined. Let
$f^{\phi}=f\circ\phi^{-1}$.  There exists a basis preserving
isomorphism
$$h^{\phi}:\mathcal{C}^{+}(L;f,\rho,J)\to
\mathcal{C}^{+}(L;f^{\phi},\rho^{\ast},J^{\ast})$$
induced by $x \to
\phi(x)$ for all $x\in \Crit(f)$ where $\rho^{\ast},J^{\ast}$ are
obtained by the push-forward of $\rho,J$ by means of $\phi$ and the
symplectomorphism $\bar{\phi}$. The isomorphism $h^{\phi}$ acts in
fact as an identification of the two complexes.

Finally, there is also the standard comparison chain morphism
$$c: \mathcal{C}^{+}(L;f^{\phi},\rho^{\ast}, J^{\ast})\to
\mathcal{C}^{+}(L;f,\rho,J)~.~$$

We now consider de composition $k=c\circ h^{\phi}$.  It is clear that
this map induces an isomorphism in homology and that it preserves the
augmented ring structure (as each of its factors does so). We now
inspect the Morse theoretic analogue of these morphisms - in the sense
that we consider instead of the complexes $\mathcal{C}(L, f,- )$ the
respective Morse complexes $C(f,-)$. It is easy to see that, by
possibly redescribing $H_{\ast}(c)$ as a morphism induced by a chain
morphism given in the same way as $\phi'$ in the proof of Corollary
\ref{cor:duality}, the Morse theoretic version of $k$ induces in Morse
homology precisely $H_{\ast}(\phi)$. But this means that at the
$E^{2}$ stage of the degree filtration the morphism induced by $k$ has
the form $H_{\ast}(\phi)\otimes id_{\Lambda ^{+}}$.

We now denote $k=\hbar(\bar{\phi})$ and we need to verify that for any
two elements $\bar{\phi}, \bar{\psi}\in Symp(M)$ we have
$\hbar(\bar{\phi}\circ\bar{\psi})=\hbar(\phi)\circ\hbar(\psi)$.  It is
easy to see that this is implied by the commutativity up to homotopy
of the following diagram:
$$\xy\xymatrix@+10pt{ \mathcal{C}^{+}(L;f')
  \ar[r]^{h^{\phi}}\ar[d]_{c} &
  \mathcal{C}^{+}(L;(f')^{\phi})\ar[d]^{c'}\\
  \mathcal{C}^{+}(L;f)\ar[r]_{h^{\phi}}& \mathcal{C}^{+}(L,f^{\phi}) }
\endxy $$
for any two Morse function $f$ and $f'$ so that the
respective complexes are defined. To see this, first we use some
homotopy $H$, joining $f$ to $f'$, to provide the comparison morphism
$c$ and we then use the homotopy $H\circ \phi^{-1}$ to define $c'$.

Finally, recall that the module structure of $Q^{+}H_{\ast}(L)$ over
$Q^{+}H(M)$ is defined by using an additional Morse function $F:M\to
\R$. If we put $F^{\phi}=F\circ \phi^{-1}$ we see easily that the
external operations defined by using $f,F,\rho,J$ and
$f^{\phi},F^{\phi},\rho^{\ast},J^{\ast}$ are identified one to the
other via the application $h^{\phi}$ (extended in the obvious way to
the critical points of $F$). The usual comparison maps then are used,
as before, to compare (by using appropriate homotopies)
$f^{\phi},F^{\phi},\rho^{\ast},J^{\ast}$ to $f,F,\rho,J$.  At the
level of the Morse quantum homology on $M$ the composition of these
two maps induces $H_{\ast}(\bar{\phi})\otimes id_{\Lambda^{+}}$.
Therefore, if $\bar{\phi}\in Symp_{0}(M)$, it follows that this last
map is the identity and proves the claim.  \qed

\begin{rem}
   It results from the proof above that for $\hbar(\phi)$ to be an
   algebra automorphism it is sufficient that $\bar{\phi}$ induce the
   identity at the level of the singular homology of $M$, e.g. $\phi$
   is homotopic to the identity.
\end{rem}

\subsection{Minimality for the pearl complex.}\label{subsec:proof_minimal}
This subsection is purely algebraic and its purpose is to show the statements in
\S\ref{subsec:minimal}.

\label{subsec:minimality_proof}
{\em Proof of Proposition \ref{prop:min_model}.}  We start by choosing generators for
  the complex $(G,d_{0})$ as follows: $G=\Z_2 <x_{i}: i\in I >\oplus \Z_2<y_{j}:j\in J>\oplus \Z_2<y'_{j}: j\in J>$ so that $d_{0}x_{i}=0$, $d_{0}(y_{j})=0$, $d_{0}y'_{j}=y_{j}$, $\forall j\in J$.
  The index  families $I$ and $J$ are finite. Clearly, $\mathcal{H}\cong \Z_2 <x_{i}>$ and we will identify
  further these two vector spaces and denote $\mathcal{C}_{min}=\Z_2 <\tilde{x}_{i}>\otimes \La^{+}$
  where $\tilde{x}_{i}, i\in I$ are of the same degree as the $x_{i}$'s (
  obviously, the differential on $\mathcal{C}_{min}$ remains to be defined).
  We  will construct $\phi$ and $\psi$ and $\delta$ so that
  $\phi_{0}(x_{i})=\tilde{x}_{i}$, $\phi_{0}(y_{j})=\phi_{0}(y'_{j})=0$ and
  $\psi_{0}(\tilde{x}_{i})=x_{i}$.
  The construction is by induction. We fix the following notation: $\mathcal{C}^{k}=
  \Z_2< x_{i},\ y'_{j},\ y_{j}\ : \ \ |x_{i}|\geq k,  |y'_{j}|\geq k\  >\otimes \La^{+}$.
  Similarly, we put $\mathcal{C}^{k}_{min}=\Z/2<\tilde{x}_{i} \ : \ |x_{i}|\geq k \ >\otimes \La^{+}$.

  Notice that
  there are some generators in $\mathcal{C}^{k}$ which are of degree $k-1$, namely the $y_{j}$'s
  of that degree. With this notation we also see that $\mathcal{C}^{k}$ is a sub-chain complex of $\mathcal{C}$.  To simplify notation we will identify the generators of these complexes by their
  type - $x$, $y'$, $y$, $\tilde{x}$ and their degree. Assume that $n$ is the maximal degree of the generators
  in $G$.  For the generators of $\mathcal{C}^{n}$ we let $\phi$ be equal to $\phi_{0}$, we put
  $\delta=0$ on $\mathcal{C}_{min}^{n}$ and we also let $\psi=\psi_{0}$
  on $\mathcal{C}_{min}^{n}$. To see that $\psi:\mathcal{C}^{n}\to \mathcal{C}_{min}^{n}$ is
  a chain morphism with these definitions it suffices to remark that if a generator of type $y$ has
  degree $n-1$, then $y=d_{0}y'=dy'$ and so $dy=0$.
  We now assume $\phi,\delta,\psi$ defined on $\mathcal{C}^{n-s+1}$, $\mathcal{C}_{min}^{n-s+1}$
  so that $\phi,\psi$ are chain morphisms, they induce isomorphisms in homology and $\phi\circ \psi=id$.

  We now indend to extend these maps to $\mathcal{C}^{n-s}$, $\mathcal{C}_{min}^{n-s}$.
  We first define $\phi$ on the generators of  type $x$ and $y'$ which are of degree $n-s$:
  $\phi(x)=\tilde{x}$, $\phi(y')=0$. We  let $\delta(\tilde{x})=\phi^{n-s+1}(dx)$ (
  when needed, we use the superscript $(-)^{n-s+1}$ to indicate the maps previously
  constructed by induction). Here it is important
  to note that, as $d_{0}x=0$, we have that $dx\in\mathcal{C}^{n-s+1}$. We consider now the
  generators of type $y\in \mathcal{C}^{n-s}$ which are of degree $n-s-1$ and we put
   $\phi(y)=\phi^{n-s+1}(dy'-y)$. This makes sense because $dy'-y\in \mathcal{C}^{n-s+1}$.
   If we write $dy'=y+y''$ we see  $\phi (dy')= \phi(y) + \phi(y'')=0=\delta(\phi(y'))$ so that, to make sure
   that $\phi^{n-s}$ is a chain morphism with these definitions, it remains to check
   that $\delta\phi(y)=\phi(d y)$ for all  generators of type $y$ and of degree $n-s-1$.
   But $\delta \phi(y)=\delta \phi^{n-s+1}(y'')$ and as $\phi^{n-s+1}$ is a chain morphism,
   we have $\delta \phi^{n-s+1}(y'')=\phi^{n-s+1} d(y'')$ which implies our identity because
   $d y''+dy=d^{2}y'=0$.  It is clear that $\phi$ so defined induces an isomorphism
    on the $d$-homology of $\mathcal{C}^{n-s}$ because the kernel of $\phi$ is generated
   by couples $(y',dy')$ so that it is acyclic. To conclude our induction step it remains to construct
   the map $\psi$ on the generators $\tilde{x}$ of degree $n-s$.   We now consider the
   difference $dx-\psi^{n-s+1}(\delta \tilde{x})$ and we want to show that
   there exists $\tau\in \mathcal{C}^{n-s+1}$ so that $d\tau=dx-\psi^{n-s+1}(\delta \tilde{x})$
   and $\tau\in \ker (\psi^{n-s+1})$. Assuming the existence of this $\tau$ we will put
   $\psi(\tilde{x})=x-\tau$ and we see that $\psi$ is a chain map and $\phi\circ\psi=id$.
   To see that such a $\tau$ exists remark that first $w=dx-\psi^{n-s+1}(\delta \tilde{x})\in \mathcal{C}^{n-s+1}$ and $dw=d(\psi^{n-s+1}(\delta \tilde{x}))=\psi^{n-s+1}(\delta\circ \delta \tilde{x})=0$ (because
   $\psi^{n-s+1}$ is a chain map).
   Moreover, $\phi(w)=\phi^{n-s+1}(dx)-\delta\tilde{x}=0$ because $\phi^{n-s+1}\circ\psi^{n-s+1}=id$.
   Therefore $w$ is a cycle belonging to $\ker (\phi^{n-s+1})$. But $\phi^{n-s+1}$ is a chain morphism
   which induces an isomorphism in homology and which is surjective. Therefore $H_{\ast}(\ker(\phi^{n-s+1}))=0$. Thus there exists $\tau\in\ker(\phi^{n-s+1})$ so that $d\tau=w$ and this concludes the induction
   step. With this construction it is clear that $\phi_{0}$ induces an isomorphism in $d_{0}$-homology
   and as $\phi_{0}\circ\psi_{0}=id$ we deduce that so does $\psi_{0}$.

   This construction concludes the first part of the statement and to finish the proof of the proposition
   we only need to prove the uniqueness result. The following lemma is useful.

 \begin{lem}\label{lem:isomorphisms_min} Let $G$, $G'$ be finite dimensional, graded
 $\Z_2$-vector spaces. A morphism $\xi:G\otimes \La^{+}\to G'\otimes \La^{+}$ is an isomorphism iff $\xi_{0}$ is an isomorphism.
 \end{lem}
 \begin{proof} Indeed, it is immediate
 to see that if  $\xi$ is an isomorphism, then $\xi_{0}$ is one: the surjectivity of $\xi$ implies that
 of $\xi_{0}$ and a dimension count concludes this direction. Conversely, if $\xi_{0}$ is surjective
 a simple induction argument shows that $\xi$ is surjective too.  Assume that the maximal degree in
 $G$ is $k$. Obviously $\xi=\xi_{0}|_{G_{k}\otimes \La^{+}}$. We now suppose that $\xi$ is
 surjective when restricted to  $G_{k-s+1}\otimes \La^{+}\to G'_{k-s+1}\otimes \La^{+}$. Take
 $g'\in G'_{k-s}$. Then $g'=\xi_{0}(g)$ for some $g\in G_{k -s}$ so that we may write
 $\xi(g)=g'+g''t$ for some $g''\in G'_{k-s+1}\otimes \La^{+}$. But, by the induction hypothesis,
  $g''\in Im (\xi|_{G_{k-s+1}\otimes \La^{+}})$ so that $g'=\xi(g)-\xi(g''')t$ with $\xi(g''')=g''$,
  $g'''\in G_{k-s+1}\otimes \La^{+}$. A dimension count again shows the injectivity of $\xi$.
  \end{proof}

  To end the proof of the proposition, suppose $\phi': \mathcal{C}\to \mathcal{C}'$
  and $\psi':\mathcal{C}'\to \mathcal{C}$ are chain morphisms so that $\phi'\circ\psi'=id$ with
  $\mathcal{C}'=(H\otimes \La^{+},\delta')$, $\delta'_{0}=0$ and $H$ some graded, $\Z_{2}$-vector space and $\phi'$, $\psi'$, $\phi'_{0}$, $\psi'_{0}$ induce isomorphisms in the (respective) homology.
 We want to show that there exists a chain map $c:\mathcal{C}_{min}\to \mathcal{C}'$ so
 that $c$ is an isomorphism.  This is quite easy: we define $c(u)=\phi'\circ\psi(u)$, $\forall u\in \mathcal{C}_{min}$.  Now $H_{\ast}(\phi_{0})$ and $H_{\ast}(\phi'_{0})$, $H_{\ast}(\psi_{0})$,
 $H_{\ast}(\psi'_{0})$ are all isomorphisms (in $d_{0}$-homology). So $H(c_{0})$ is
 an isomorphism but as $\delta_{0}=0=\delta'_{0}$ we deduce that $c_{0}$ is an isomorphism.
 \qed

\

 {\em Proof of Corollary \ref{cor:min_pearls}.}
 Suppose that $\mathcal{C}(L;f,\rho,J)$ is defined and apply the Proposition \ref{prop:min_model}
 to it.  Denote by $(\mathcal{C}_{min},\phi,\psi)$ the resulting minimal complex and
 the chain morphisms as in the statement of \ref{prop:min_model}. The only part of the statement which
 remains to be shown is that given a different set of data $(f',\rho'.J')$ so that $\mathcal{C}(L;f'\rho',J')$
 is defined, there are appropriate morphisms $\phi',\psi'$ as in the statement.
 There are comparison morphisms: $h:\mathcal{C}(L;f'\rho',J')\to
\mathcal{C}(L;f\rho,J)$ as well as $h':\mathcal{C}(L;f\rho,J)\to \mathcal{C}(L;f'\rho',J')$
so that, by construction, both $h$ and $h'$ are inverse in homology and  both
induce an isomorphism in Morse homology and again these two isomorphisms are inverse
(see \S \ref{subsubsec:inv}). Define $\phi':\mathcal{C}(L;f',\rho',J')\to \mathcal{C}_{min}$,
$\psi'':\mathcal{C}_{min}\to \mathcal{C}(L;f',\rho',J')$ by  $\phi'=\phi\circ h$ and
$\psi''=h'\circ \psi$. It is clear that $\phi'$, $\psi''$, $\phi'_{0}$ and $\psi''_{0}$ induce
isomorphism in homology. Moreover, as $h_{0}$ and $h'_{0}$ are inverse in homology and
$\delta_{0}=0$ in $\mathcal{C}_{min}$ it follows that $\phi'_{0}\circ \psi''_{0}=id$.
This means by the Lemma \ref{lem:isomorphisms_min} that $v=\phi'\circ\psi''$ is a chain isomorphism
so that $v_{0}$ is the identity. We now put $\psi'=\psi''\circ v^{-1}$ and this verifies all the
needed properties.  The uniqueness of $\mathcal{C}_{min}(L)$ now follows from the uniqueness
part in the Proposition \ref{prop:min_model}.
  \qed

\subsection{Proof of the action estimates.}
\label{subsec:action_proof}
We first recall the definition of the two spectral invariants
involved.  Fix a generic pair $(H,J)$ consisting of a $1$-periodic
Hamiltonian $H:M\times S^{1}\to \R$ and an almost complex structure
$J$ so that the Floer complex $CF_{\ast}(H,J)$ is well defined. We
will assume the coefficients of this complex to be in the usual
Novikov ring $\Lambda$. We recall that the generators of
$CF_{\ast}(H,J)$ as a module over $\Lambda$ are pairs formed by
contractible orbits of $X^H$ together with fixed cappings.  Fix also a
Morse function $f:L\to \R$ as well as a Riemannian metric $g$ on $L$
so that the pearl complex $\mathcal{C}^{+}(L;f,\rho,J)$ is well
defined.

We first need to provide a description of our module external
operation which involves the two complexes above. This is based on
moduli spaces $\mathcal{P}(\gamma, x,y; \lambda)$ similar to the ones
used in \S\ref{S:qm} except that the vertex of valence three in the
string of pearls is now replaced by a half-tube with boundary on $L$
and with the $-\infty$ end on $\gamma$.  The total homotopy class
$\lambda$ of the configuration obtained in this way is computed by
using the capping associated to $\gamma$ to close the semi-tube to a
disk and adding up the homotopy class of this disk to the homotopy
classes of the other disks in the string of pearls.  More explicitly,
a half tube as before is a solution
$$u:(-\infty,0]\times S^{1}\to M$$
of Floer's equation $\partial
u/\partial s +J\partial u/\partial t+\nabla H(u,t)=0$ with the
boundary conditions $$u(\{0\}\times S^{1})\subset L \ \
\lim_{t\to-\infty}u(s,t)=\gamma(t)~.~$$
The incidence points on the
``exceptional" vertex which corresponds to $u$ are so that the point
$u(0,1)$ is an exit point for a flow line and $u(0,-1)$ is the entry
point. Both compactification and bubbling analysis for these moduli
spaces are similar to what has been discussed before to which is added
the study of transversality and bubbling for the spaces of half-tubes
as described by Albers in \cite{Alb:extrinisic} and, as described in
\cite{Alb:extrinisic}, an additional assumption is needed for this
part: $H$ is assumed to be such that no periodic orbit of $X^H$ is
completely included in $L$.

Counting elements in these moduli spaces defines an operation:
$$CF(H,J)\otimes \mathcal{C}(L;f,\rho,J)\to \mathcal{C}(L;f,\rho,J)$$
and, by using the Hamiltonian version of the Piunikin-Salamon-Schwarz
homomorphism, it is easy to see that, in homology, this operation is
canonically identified with the module action as described in
\S\ref{S:qm}.

\

The Floer complex $CF_{\ast}(H,J)$ is filtered by the values of the
action functional
$$\mathcal{A}_{H}(\overline{x})=\int_{0}^{1}H(x(t),t)dt -\int_{D}
\hat{x}^{\ast}\omega$$
where $\overline{x}=(x,\hat{x})$ with $x$ a
$C^{\infty}$ loop in $M$ and $\hat{x}$ a cap of this loop.  This
action is compatible with the Novikov ring in the sense that:
$\mathcal{A}_{H}(\gamma \otimes t^{k})=\mathcal{A}_{H}(\gamma)-k\tau$
(where $\tau$ is the monotonicity constant). The filtration of order
$\nu\in\R$ of the the Floer complex, $CF^{\nu}$, is the graded
$\mathbb{Z}_2$-vector space generated by all the elements
$\gamma\otimes\lambda$ of action at most $\nu$. This is a sub-complex
because the differential decreases action.

We now fix $\alpha\in H_{\ast}(M;\mathbb{Z}_2)$ and define
$\sigma(\alpha,H)$ by:
\begin{equation}\label{eq:spectral_hlgy}
   \sigma(\alpha,H) =
   \inf\{\nu : PSS(\alpha)\in
   \mathrm{Image}(\ H(CF^{\nu})\to HF(H,J)\ )\}.
\end{equation}
Here $PSS:H_{\ast}(M;\mathbb{Z}_2)\to HF(H,J)$ is the
Piunikin-Salamon-Schwarz homomorphism.  We also need a similar
definition for a cohomology class $\beta\in H^{\ast}(M;\mathbb{Z}_2)$.
A little more notation is needed for this.  We recall that we work
here over the Novikov ring $\Lambda$. We also recall the algebraic notation from
\S \ref{susubsubsec:duality}: $\Lambda^{\ast}$ is
be the ring $\Lambda$ with reverse grading (in short the element $t$
has now degree $N_{L}$) and given a free chain complex $\mathcal{C}=G\otimes
\La$ with $G$ a $\mathbb{Z}_2$ vector space,
$\mathcal{C}^{\ast}=\hom_{\mathbb{Z}_2}(G,\mathbb{Z}_2)\otimes
\La^{\ast}$ where the grading of the dual $x^{\ast}$ of a basis
element $x\in G$ is $|x^{\ast}|=|x|$ with the differential given
 as the adjoint of the differential in $\mathcal{C}$.

From the Floer complex $CF(H,J)$ we define the associated Floer
co-homology by $H^{k}F(H,J)=H^{k}(CF(H,J)^{\ast})$ and
similarly for all the various subcomplexes involved. In particular, we
have $H^{k}(CF^{\nu})=H^{k}((CF^{\nu})^{\ast})$ as well as
morphisms $$p_{\nu}:CF(H,J)^{\ast}\to (CF^{\nu})^{\ast}$$
which are
induced by the inclusions $CF^{\nu}\hookrightarrow CF(H,J)$. Moreover,
the inverse PSS map induces also a comparison morphism
$H^{\ast}(M;\mathbb{Z}_2)\to H^{\ast}F(H,J)$ which we will denote by
$PSS'$.  In view of this we may now define:
\begin{equation}\label{eq:spectral_chlgy}
   \sigma(\beta,H)=\sup\{\nu : \  PSS'(\beta)\in
   \ker( H^{\ast}F(H,J)\to H^{\ast}(CF^{\nu}))\}.
\end{equation}

Assuming that $H$ is normalized, it is well known that
$\sigma(\alpha,H)$, the spectral invariant of $\alpha$, only depends
on the class $[\phi^H]\in \widetilde{Ham(M)}$. The same holds for the
spectral invariant of the co-homology class $\beta$.

\begin{rem}\label{rem:co-ho-spec}
   The spectral invariant of a co-homology class $\beta$ satisfies te
   following property which will be useful in the following.  For any
   $\epsilon > 0$ there exists a co-chain $c=\sum_{i}
   \gamma_{i}^{\ast}$ so that $[c]=\beta$, $\gamma_{i}^{\ast}$ are
   dual to orbits $\gamma_{i}$ (the $\Lambda$-coefficients are
   integrated in $\gamma_{i}$ by a possible change of capping - in
   other words we view $\mathcal{C}$ as a $\mathbb{Z}_2$ vector space)
   and $\mathcal{A}_{H'}(\gamma_{i})\geq \sigma(\beta,H)-\epsilon$ for
   all $i$. For this, first notice that the map
   $p_{\nu}:CF(H,J)^{\ast}\to (CF^{\nu})^{\ast}$ is surjective and
   that its kernel is generated by those $\gamma^{\ast}$ with
   $\mathcal{A}_{H}(\gamma)>\nu$.  Fix $\nu=\sigma(\beta,H)-\epsilon$
   and let $c'=\sum \gamma^{\ast}_{i}$ be so that $[c']=\beta$.  Let
   $p_{\nu}(c')=\delta$. Then as $\nu< \sigma(\beta,H)$ we deduce
   $\delta=\partial^{\ast}h$ and as $p_{\nu}$ is surjective, we may
   view $h$ as an element of $CF(H,J)^{\ast}$ so that
   $p_{\nu}(c'-\partial^{\ast}h)=0$. But $\beta=[c'-\partial^{\ast}h]$
   so that our claim follows.
\end{rem}

\begin{rem}
   To avoid possible confusion with various other conventions used in
   the literature, notice that, with the definitions above, we have in
   general $\sigma(\alpha^{\ast},H)\not=\sigma(\alpha, H)$ where
   $\alpha^{\ast}$ is the co-homology class Poincar\'e dual to
   $\alpha$.
\end{rem}

\subsubsection{Proof of Corollary \ref{cor:action}.}
We recall that $\alpha\in H_{\ast}(M;\mathbb{Z}_2)$ is fixed as well
as $x,y\in Q^{+}H(L)$ so that $y\not=0$ and $\alpha *
x=yt^{k}+\mathrm{higher\ order\ terms}$.  We also fix
$\phi\in\widetilde{Ham(M)}$. We first intend to show that:
$\sigma(\alpha,\phi)-\mathrm{depth}_{L}(\phi)+k\tau\geq 0$. By
inspecting the definition of \textit{depth} in \S \ref{subsec:action}
we see that this reduces to showing that for every normalized
Hamiltonian $H$ with $[H]=\phi$, we have
$$\sigma(\alpha, H)-\int_{0}^{1}H(\gamma(t),t)dt+k\tau\geq0$$
for some
loop $\gamma:S^{1}\to L$. By a small perturbation of $H$ we may assume
that no closed orbit of $H$ is contained in $L$.

Now assume that $\eta=\sum \gamma_{i}\otimes \lambda_{i}$ is a cycle
in $CF(H,J)$ so that $[\eta]=PSS(\alpha)$ where $\gamma_{i}$ are
generators of $CF(H,J)$ and $\lambda_{i}\in\Lambda$. The relation
$\alpha\ast x=yt^{k}+...$ implies that there exits some periodic orbit
$\gamma_{i}$, say $\gamma_{1}$, and critical points $x_{1}, y_{1}\in
\Crit(f)$, $|x_{1}|=|x|$, $|y_{1}|=|y|$ so that the moduli space
$\mathcal{P}(\gamma_{1}',x_{1},y_{1}; t^{k})\not =\emptyset$ where
$\gamma'_{1}=\gamma_{1}\# \lambda_{1}$ is the same orbit of $X^H$ as
$\gamma_{1}$ but with the capping changed by $\lambda_{1}$.  We now
consider an element $v\in\mathcal{P}(\gamma_{1}',x_{1},y_{1}; t^{k})$
and we focus on the corresponding half-tube $u$ (which is part of
$v$). The usual energy estimate for this half-tube gives:

\begin{equation}\label{eq:energy_est}
   0\leq\int_{-\infty}^{0}\int_{0}^{1}||\partial u/\partial s||^{2}dtds=
   \int_{[-\infty,0]\times S^{1}}u^{\ast}\omega +
   \int_{S^{1}}H(\gamma_{1}(t),t)-\int_{S^{1}}H(u(0,t),t)~.~
\end{equation}
We now want to remark that:
$$\mathcal{A}_{H}(\gamma_{1}')+k\tau\geq \int_{[-\infty,0]\times
  S^{1}}u^{\ast}\omega+\int_{S^{1}}H(\gamma_{1}(t),t)~.~$$
Indeed,
this is obvious in view of the definition of the action and given that
$k\tau$ equals the symplectic area of all the disks in $v$ $+$ the
area of the tube $u$ $+$ the area of the cap corresponding to
$\gamma_{1}'$ (we have equality here iff no $J$-disks appear in $v$).
Given any $\epsilon>0$, in view of the definition of $\sigma(\alpha,
H)$, it follows that we may find in $CF^{\sigma(\alpha,H)+\epsilon}$ a
cycle $\eta$ with $[\eta]=PSS(\alpha)\in HF(H,J)$. Applying the
discussion above to this $\eta$ means
$\mathcal{A}_{H}(\gamma_{1}')\leq \sigma(\alpha,H)+\epsilon$ and this
implies the claimed inequality.

We now want to show the second inequality in Corollary
\ref{cor:action}:
$\mathrm{height}_{L}(\phi)-\sigma(\alpha^{\ast},\phi)+k\tau\geq 0$.
Again by taking a look at the definition of
$\mathrm{height}_{L}(\phi)$ we see that it is enough to show that for
some normalized Hamiltonian $H'$ we have
$\int_{0}^{1}H'(\gamma(t),t)-\sigma(\alpha^{\ast},H')+k\tau\geq 0$ for
some loop $\gamma$ in $L$.  We now return to the choices of $H,v,u,
\gamma'_{1}, x_{1},y_{1}$ used when establishing
formula~\eqref{eq:energy_est}. We now define $H'(x,t)=-H(x,t)$. We
notice that the the periodic orbits of $H'$ are related to those of
$H$ by the formula $\gamma(t)\to \gamma(1-t)$ and, moreover, the
complex $CF(H',J)$ is in fact identified with $s^{n}(CF(H,J))^{\odot}$
(to recall the algebraic notation $(-)^{\odot}$ etc see
\S\ref{susubsubsec:duality} and \S\ref{subsec:proof_duality}). The
equation~\eqref{eq:energy_est} can be interpreted by looking at $u$ as
a Floer half tube for $H'$ parametrized by $[0,+\infty)\times S^{1}$
with the $0$ end on $L$ and the $+\infty$ end on $\tilde{\gamma}'_{1}$
(where for a loop $\gamma$, $\tilde{\gamma}$ is the loop
$\gamma(1-t)$). We obtain:
\begin{equation}\label{eq:energy_height}
   0\leq \int u^{\ast}\omega -\int_{S^{1}}H'(\tilde{\gamma}_{1}(t),t)dt
   +\int_{S^{1}}H'(u(0,-t),t)dt~.~
\end{equation}
So that is suffices to show $k\tau-\sigma(\alpha^{\ast},H')\geq \int
u^{\ast}\omega -\int_{S^{1}}H'(\tilde{\gamma}_{1}(t),t)dt$.  By
noticing that the capping corresponding to $\tilde{\gamma}'_{1}$ is
changed with respect to that of $\gamma'_{1}$ by $(s,t)\to (s,1-t)$,
it is now easy to see that:

$$\mathcal{A}_{H'}(\tilde{\gamma}'_{1})
\leq\int_{S^{1}}H'(\tilde{\gamma}_{1}(t),t)dt -\int u^{\ast}\omega +
k\tau~.~$$ And the proof of the desired inequality reduces to showing
that, for any $\epsilon>0$ there exists $\gamma_{1}$ as before so that
$\sigma(\alpha^{\ast},H')-\epsilon\leq
\mathcal{A}_{H'}(\tilde{\gamma}'_{1})$.  But this follows immediately
from Remark \ref{rem:co-ho-spec}.
\qed

\subsection{Replacing $\mathcal{C}^{+}$, $Q^+H$ by
  $\mathcal{C}$, $QH$} \label{Sb:Q-Q+}

At this stage we indicate that everything proved in
\S\ref{subsec:pearl} --~\ref{Sb:spect-seq} as well as
in~\S\ref{subsec:proof_duality},~\ref{Sb:act-symp}, ~\ref{subsec:minimality_proof} continues to hold
(with the same proofs) if we replace $\mathcal{C}^{+}$ and $Q^+H$ by
$\mathcal{C}$ and $QH$ respectively everywhere. The same is true
for~\S\ref{subsec:action_proof} if we replace $\mathcal{C}^+(L)$ by
$\mathcal{C}(L)$ but in that section we still have to work with the
full version of $HF(M)$.

Note however, that in contrast to the above, in
\S\ref{Sb:comparison-HF} it is essential to work with $\mathcal{C}$
and $QH$ rather than their positive versions.



%

\section{Applications and examples.}\label{sec:appli}

Our applications are grouped in three categories. Ths first has to do
with the algebraic constraints coming from the interplay between
singular and quantum structures.  There are many such examples in the
paper. In \S\ref{subsubsec:tori} we see, for example, that when the
singular homology of $L$ is generated as an algebra (with the
intersection product) by classes of sufficiently high degree
(depending on the minimal Maslov number), then, either, the quantum
homology of $L$ is additively just singular homology with Novikov
coefficients or it vanishes. Other examples come from the interplay of
the topology of the ambient manifold and that of the Lagrangian.
Thus, in \S\ref{S:cpn}, we see that if the ambient manifold is $\C
P^{n}$, then the resulting restrictions on the Floer homology of the
Lagrangian are stringent.  As a concrete example we discuss the
Clifford torus in detail in \S\ref{Ss:clifford}.  Additional
assumptions on the Lagrangian - for example, $2H_{1}(L;\mathbb{Z})=0$
as in \S\ref{Sb:rpn-cpn} - lead to more homological rigidity as such
Lagrangians are seen to have a homology algebra very much like that of
$\R P^{n}$. We also discuss Lagrangian submanifolds of the quadric in
\S\ref{Sb:quadric}, as well as in complete intersections in
\S\ref{Sb:cinter}. In \S\ref{Sb:alg-qh} we also describe some examples
that go in the opposite direction: the existence of certain Lagrangian
submanifolds is seen to imply restrictions on the quantum homology of
the ambient manifold.

The second type of application is related to ways of measuring the
size of Lagrangians as well as that of the space surrounding them
inside the ambient manifold. Indeed, as introduced in
\cite{Bar-Cor:NATO} and in \cite{Cor-La:Cluster-2} there is a natural
notion of Gromov width of a Lagrangian and obviously one can also
consider the width of the complement of that Lagrangian. Moreover, one
can also define packing numbers for the Lagrangian, its complement as
well as mixed numbers.  Our techniques allow us to give estimates - in
\S\ref{Sb:pack} - for many of these numbers in all the cases mentioned
above.

The third type of application - presented in \S\ref{Sb:clif-enum} - is
concerned with the fact that the properties of our machinery can be
used to define certain numerical invariants roughly of Gromov-Witten
type which are associated to configurations different from the usual
ones. We discuss this construction in a very explicit way for two
dimensional monotone tori - this turns out to be a surprisingly rich
case. The invariants in question are associated to triangles lying on
the tori and are expressed as polynomials involving numbers of
$J$-holomorphic disks passing through the vertexes and/or the edges of
the triangle.

\

There is an underlying unifying idea for all of these diverse
applications.  Lagrangian submanifolds exhibit considerable rigidity:
topological (or algebraic) for our first class of applications,
geometric for the second and arithmetic for the third.

\subsection{Full Floer homology}\label{subsubsec:tori}

\begin{prop}\label{cor:floer=0-or-all}
   Let $L^n \subset (M^{2n}, \omega)$ be a monotone Lagrangian with
   $N_L \geq 2$.  Assume that its singular homology
   $H_*(M;\mathbb{Z}_2)$ is generated as an algebra by $H_{\geq
     n-k}(L;\mathbb{Z}_2)$. Suppose further that $N_L > k$.
   Then:
   \begin{enumerate}
     \item either $QH_*(L)=0$; or \label{I:QH-0-or-all-1}
     \item there exist isomorphisms of graded vector spaces $QH_*(L)
      \cong (H(L;\mathbb{Z}_2) \otimes \Lambda)_*$ and $Q^{+}H_*(L)
      \cong (H(L;\mathbb{Z}_2) \otimes \Lambda^{+})_*$. These
      isomorphisms are in general not canonical. Moreover, these
      isomorphisms, in general, do not respect the ring structures.
      \label{I:QH-0-or-all-2}
   \end{enumerate}
   If $N_L > k+1$ only the second alternative occurs.  Furthermore,
   when the second alternative occurs (whether $N_L > k+1$ or
   $N_L=k+1$) there exist canonical injections of $H_{\geq
     n-k}(L;\mathbb{Z}_2)$ into $QH_{\geq n-k}(L)$ and into
   $Q^{+}H_{\geq n-k}(L)$ which generate these algebras over $\Lambda$
   and over $\Lambda^+$ with respect to the quantum product.
\end{prop}

\begin{ex}
   Tori with minimal Maslov class at least $2$ furnish a nice example.
   Another immediate example is $\R P^{n}\subset \C P^{n}$. Other
   examples will appear later in this section.
\end{ex}

\begin{proof}[Proof of Proposition~\ref{cor:floer=0-or-all}] We will provide
two proofs for this proposition. The first one is based on a spectral sequence argument while
the second uses the minimal model machinery described in \S\ref{subsec:minimal}. We include
both arguments precisely to illustrate the role of these minimal models.

\

   {\bf A}. We will use here an argument involving spectral sequences.  Choose
   a generic $J \in \mathcal{J}(M,\omega)$. Let $f: L \to \mathbb{R}$
   be a Morse function with exactly one maximum $x_n$ and fix a
   generic Riemannian metric on $L$. Denote by $(CM_*(f),
   \partial_0)$, $(\mathcal{C}_*(f,J), d)$ the Morse and pearl
   complexes associated to $f$, $J$ and the Riemannian metric.

   Recall that $\mathcal{C}(f,J)$ is filtered by the degree filtration
   $\mathcal{F}^p \mathcal{C}$ (which is a decreasing filtration). It
   will be more convenient to work here with an increasing version of
   the same filtration. Put
   $$\mathcal{F}_p \mathcal{C}_i(f,J) = \bigoplus_{j \geq -p}
   \mathcal{C}_{i+jN_L}(f,J) t^j.$$
   Clearly this is a bounded
   increasing filtration. It gives rise to a spectral sequence $\{
   E^{r}_{p,q}, d_r\}_{r \geq 0}$ which converges to $QH_*(L) =
   H_*(\mathcal{C}(f,J),d)$. A simple computation shows that:
   \begin{enumerate}
     \item $E^0_{p,q} = CM_{p+q-pN_L} t^{-p}$, $d_0 = \partial_0$.
     \item $E^1_{p,q} = H_{p+q-pN_L}(L; \mathbb{Z}_2) t^{-p}$.
     \item The sequence collapses after a finite number of steps. In
      fact this number of steps is $\leq [\frac{n+1}{N_L}]+1$.
   \end{enumerate}
   Let $f': L \to \mathbb{R}$ be a small perturbation of $f$. Note
   that the filtration $\mathcal{F}_p$ is compatible with the quantum
   product $*$ in the sense that
   $$* : \mathcal{F}_{p_1} \mathcal{C}_{i_1}(f,J) \otimes
   \mathcal{F}_{p_2} \mathcal{C}_{i_2}(f',J) \longrightarrow
   \mathcal{F}_{p_1+p_2} \mathcal{C}_{i_1+i_2 - n}(f,J).$$ By taking
   $f'$ close enough to $f$ we may assume that the canonical
   quasi-isomorphism between $\mathcal{C}_*(f,J)$ and
   $\mathcal{C}_*(f',J)$ is in fact a base preserving isomorphism (see
   \S\ref{subsubsec:inv}). Thus we may identify the spectral sequences
   $\{ E^r_{p,q}(f,J), d_r\}$ and $\{ E^r_{p,q}(f',J), d_r\}$ and
   denote both of them by $\{ E^r_{p,q}, d_r\}$. It follows that this
   spectral sequence is multiplicative.

   Denote by $*_r$ the product induced by $*$ on $\{E^r_{*,*} \}$ and
   by $*_{\infty}$ the product induced by $*$ on $\{E^{\infty}_{*,*}
   \}$. Note that although there exists an isomorphism
   $H_l(\mathcal{C}, d) \cong \oplus_{p+q=l} E^{\infty}_{p,q}$ for
   every $l \in \mathbb{Z}$, the products $*$ on $H_*(\mathcal{C}, d)$
   and $*_{\infty}$ on $\oplus_{p+q=l} E^{\infty}_{p,q}$ are in
   general not isomorphic.

   A direct computation shows that the product $*$ induces on the
   $E^1$ level the classical cap product, namely the product $$*_1:
   E^1_{p_1, q_1} \otimes E^1_{p_2, q_2} \to E^1_{p_1+p_2,
     q_1+q_2-n}$$
   is the cap product $$\cap : H_{p_1+q_1 - p_1
     N_L}(L;\mathbb{Z}_2) t^{-p_1} \otimes H_{p_2+q_2 - p_2
     N_L}(L;\mathbb{Z}_2) t^{-p_2} \to H_{p_1+p_2+q_1+q_2-(p_1+p_2)N_L
     - n}(L;\mathbb{Z}_2)t^{-p_1-p_2}$$
   and $d_1$ satisfies Leibniz
   rule with respect to $\cap$.  It follows that
   $\{E^1_{p,q}\}_{p+q-pN_L \geq n-k}$ generate with respect to $*_1 =
   \cap$ the whole $\{E^1_{*,*}\}$.

   Assume now that $N_L > k+1$. For degree reasons $d_1$ vanishes on
   $E^1_{p,q}$ whenever $p+q-pN_L \geq n-k$. Since $d_1$ satisfies
   Leibniz rule with respect to $*_1$ it follows that $d_1=0$
   everywhere. Therefore $E^2_{*,*}=E^1_{*,*}$ and $*_2 = *_1=\cap$.
   The same argument applied to $d_2$ shows that $d_2=0$, hence
   $E^3_{*,*}=E^2_{*,*}=E^1_{*,*}$. Proceeding by induction we obtain
   that $d_r=0$ for every $r \geq 1$ hence $E^{\infty}_{*,*} = \cdots
   = E^1_{*,*}$ and $*_{\infty}= \cdots = *_1 = \cap$.

   It follows that there exists an isomorphism $$H_*(\mathcal{C},d)
   \cong \bigoplus_{p+q=*} E^1_{p,q} = \bigoplus_{p \in
     \mathbb{Z}}H_{*-p N_L}(L;\mathbb{Z}_2)t^{-p} = \bigl(
   H(L;\mathbb{Z}_2) \otimes \Lambda \bigr)_*.$$

   We now turn to the case $N_L = k+1$. First note that exactly as in
   the in the case $N_L > k+1$ we have that $d_1=0$ on all $E^1_{p,q}$
   with $p+q-pN_L \geq n-k+1$. Let us examine now the behavior of
   $d_1$ on $E^1_{p,q}$ for $p+q - pN_L=n-k$.  For such $p, q$ we have
   $$E^1_{p,q} = H_{n-k}(L;\mathbb{Z}_2) t^{-p}, \quad
   E^1_{p-1,q}=H_n(L;\mathbb{Z}_2) t^{-p+1}.$$
   Moreover, it is easy to
   see that $$d_1: H_{n-k}(L;\mathbb{Z}_2) t^{-p} \to
   H_n(L;\mathbb{Z}_2)t^{-p+1}$$
   takes the form $d_1 = \delta_1 t$
   where $\delta_1: H_{n-k}(L;\mathbb{Z}_2) \to H_n(L;\mathbb{Z}_2)$
   does not depend on $p$. There are now two cases to consider:

   \smallskip \noindent \textbf{Case I. $\delta_1 \neq 0$.} Since
   $H_n(L;\mathbb{Z}_2) = \mathbb{Z}_2 [L]$ is $1$-dimensional (and we
   work here over $\mathbb{Z}_2$ which is a {\em field}) it follows
   that $[L]$ is in the image of $\delta_1$. Therefore $[L] \in
   E^1_{0,n} = H_n(L;\mathbb{Z}_2)$ is the image under $d_1$ of some
   element in $E^1_{1,n} = H_{n-k}(L;\mathbb{Z}_2)t^{-1}$.  It follows
   that the homology class of $[L]$ in $E^2_{0,n}$ is zero. But $[L]$
   is the unit of $E^1_{*,*}$ with respect to $*_1=\cap$ and so its
   homology class in $E^2_{0,n}$ is the unit of $E^2_{*,*}$ with
   respect to $*_2$. As this class is zero we have $E^2_{*,*}=0$. It
   follows that $H_*(\mathcal{C},d)=QH_*(L)=0$.

   \smallskip \noindent \textbf{Case II. $\delta_1 = 0$.} In this case
   $d_1=0$ on $E^1_{p,q}$ for $p+q-pN_L=n-k$. The proof now continues
   exactly as in the case $N_L > k+1$ discussed above, showing that
   the spectral sequence degenerates at the $E^1$ level. It follows
   that $QH_*(L) = H_*(\mathcal{C}, d) \cong \bigl( H(L;\mathbb{Z}_2)
   \otimes \Lambda)_*$. This concludes the proof of the two
   alternatives for $QH(L)$.

   We now prove that $H_{\geq n-k}(L;\mathbb{Z}_2)$ canonically
   injects into $QH_{\geq n-k}(L)$ under the assumption that either
   $N_L > k+1$, or $N_L = k+1$ and $QH(L) \neq 0$. To see this note
   first that we have a {\em canonical} homomorphism $\sigma:
   H_l(L;\mathbb{Z}_2) \to QH_l(L)$ induced by the inclusion $CM_l(f)
   \subset \mathcal{C}_l(f,J)$ for $l \geq n-k$. Indeed, assume first
   that $N_L > k+1$ and let $x \in CM_l(f)$, $l \geq n-k$, be a
   $\partial_0$-cycle. For degree reasons $dx = \partial_0x = 0$.
   Similarly, if $x$ is a $\partial_0$-boundary, say $x=\partial_0 y$,
   $y \in CM_{l+1}(f)$, then again by degree reasons $x=dy$. This
   shows that the inclusion $CM_l(f) \subset \mathcal{C}_l(f,J)$
   induces a homomorphism $\sigma: H_l(L;\mathbb{Z}_2) \to
   H_l(\mathcal{C},d) = QH_l(L)$. Suppose now that $N_L = k+1$ and
   $QH(L) \neq 0$. For degree reasons $d=\partial_0$ on $CM_l(f)$ for
   every $l \geq n-k+1$. As for $CM_{n-k}(f)$ we can write $d =
   \partial_0 + \partial_1 t$, where $\partial_1: CM_{n-k}(f) \to
   CM_n(f)$. It follows that $\partial_1$ vanishes on all
   $\partial_0$-cycles, for otherwise the maximum $x_n \in CM_n(f)$
   would be a $d$-boundary implying that $QH(L)=0$. It follows that
   every $\partial_0$-cycle $x \in CM_{n-k}(f)$ is also a $d$-cycle.
   For degree reasons every $\partial_0$-boundary in $CM_{n-k}(f)$ is
   also a $d$-boundary. Thus in this case too we have the homomorphism
   $\sigma$ induced by the inclusion $CM_l(f) \subset
   \mathcal{C}_l(f,J)$, $l\geq n-k$.

   That $\sigma$ is canonical follows from the definition of the
   canonical identifications on $QH$ described in
   \S\ref{subsubsec:inv}. The point is that, for degree reasons, the
   chain morphism $\phi : \mathcal{C}_l(f,J)\to \mathcal{C}_l(f',J')$
   defined in \S\ref{subsubsec:inv} coincides with the analogous chain
   morphism in Morse theory $CM_l(f)\to CM_l(f')$ for for $l \geq n-k$
   (because $\phi$ involves only non-negative powers of $t$). This
   completes the proof that the homomorphism $\sigma$ is well defined
   and canonical.

   Next, we prove that $\sigma: H_l(L;\mathbb{Z}_2) \to QH_l(L)$ is
   injective for $l \geq n-k$. Again, we assume that either $N_L >
   k+1$, or $N_L=k+1$ and $QH(L) \neq 0$. To prove this first note
   that for degree reasons $CM_l(f) = \mathcal{F}_0
   \mathcal{C}_l(f,J)$ and in fact $d$ coincides with $\partial_0$
   here. Therefore this equality induces an isomorphism $\sigma':
   H_l(L;\mathbb{Z}_2) \to H_l(\mathcal{F}_0 \mathcal{C},d)$. Denote
   by $\mathcal{F}_p H_l(\mathcal{C}, d)$ the $p$-level of the
   associated filtration on the homology $H_l(\mathcal{C}, d)$, i.e.
   the image of the map $\iota : H_l(\mathcal{F}_p \mathcal{C}, d) \to
   H_l(\mathcal{C}, d)$ induced from the inclusion. For degree reasons
   we have $\mathcal{F}_{-1} H_l(\mathcal{C},d) = 0$, hence
   $E^{\infty}_{0,l} = \mathcal{F}_0 H_l(\mathcal{C},d)/
   \mathcal{F}_{-1} H_l(\mathcal{C},d) = \mathcal{F}_0
   H_l(\mathcal{C},d)$. On the other hand by what we have previously
   proved we know that $E^{\infty}_{0,l} = E^1_{0,l} =
   H_l(L;\mathbb{Z}_2)$. Putting all these together we obtain:
   \begin{equation}
      H_l(L;\mathbb{Z}_2)
      \xrightarrow[\cong]{\sigma'}H_l(\mathcal{F}_0 \mathcal{C},d)
      \xrightarrow[\textnormal{surjective}]{\iota} \mathcal{F}_0
      H_l(\mathcal{C}, d) = E^{\infty}_{0,l} = H_l(L;\mathbb{Z}_2).
   \end{equation}
   It follows that $\iota \circ \sigma'$ is surjective. By dimensions
   reasons $\iota \circ \sigma'$ is an isomorphism.  But the image of
   $\iota \circ \sigma'$ is the same as the image of $\sigma$, hence
   $\sigma$ is injective.

   It remains to show that the image of $\sigma$ in $QH_{\geq n-l}(L)$
   generates $QH_*(L)$ with respect to the quantum product.  Fix $i
   \in \mathbb{Z}$. Denote by $\mathcal{A}_* \subset
   H_*(\mathcal{C},d)$ the subalgebra generated (over $\Lambda$) by
   $\sigma \bigl(H_{\geq n-k}(L;\mathbb{Z}_2) \bigr)$ with respect to
   the quantum product $*$. We want to prove that $\mathcal{A}_* =
   H_*(\mathcal{C},d)$.

   Recall that $H_i(\mathcal{C},d)$ is filtered by the induced
   increasing filtration $\{ \mathcal{F}_pH_i(\mathcal{C},d)\}_{p \in
     \mathbb{Z}}$. This filtration is bounded so that $\mathcal{F}_p
   H_i(\mathcal{C},d)=0$ for every $p \leq p_0$ for some $p_0$, and
   $\mathcal{F}_p H_i(\mathcal{C},d)=H_i(\mathcal{C},d)$ for $p \gg
   0$. Therefore it is enough to prove that
   \begin{equation} \label{Eq:A=H}
      \mathcal{F}_p H_i(\mathcal{C},d) \subset \mathcal{A}_i, \quad
      \forall\; p.
   \end{equation}
   We will prove~\eqref{Eq:A=H} by induction on $p$.  Put
   \begin{equation} \label{Eq:GrQH}
      \mathcal{G}_i = \bigoplus_{p \geq p_0} \mathcal{F}_p
      H_i(\mathcal{C},d) / \mathcal{F}_{p-1}H_i(\mathcal{C},d) =
      \bigoplus_{p+q = i} E^{\infty}_{p,q} = \bigoplus_{p \geq p_0}
      H_{i-pN_L}(L;\mathbb{Z}_2) t^{-p}.
   \end{equation}
   By what we have proved the quantum product $*$ descends to a
   product $[*]$ on $\mathcal{G}_*$ which is identified with the
   classical cap product $\cap$ on the right-hand side
   of~\eqref{Eq:GrQH} (here we extend $\cap$ in an obvious way over
   $\Lambda$). By the assumption of the proposition
   $E^{\infty}_{0,\geq n-k}$ generates $\mathcal{G}_*$ (over
   $\Lambda$) with respect to $[*]=\cap$.

   Obviously~\eqref{Eq:A=H} holds for $p \leq p_0$ since
   $\mathcal{F}_p H_i(\mathcal{C} ,d)=0$ for such $p$'s. Assume the
   statement is true for every $p \leq p'$. Let $x \in
   \mathcal{F}_{p'+1} H_i(\mathcal{C},d)$. We want to prove that $x
   \in \mathcal{A}$.  The corresponding element in~\eqref{Eq:GrQH},
   $[x] = x \bigl( \bmod{ \mathcal{F}_{p'} H_i(\mathcal{C},d)}\bigr)$,
   can be identified with an element in
   $H_{i-{p'+1}N_L}(L;\mathbb{Z}_2)t^{-p'-1}$. By the assumption of
   the proposition $[x]$ can be expressed as a linear combination
   (over $\Lambda$) of $[*]$-products of elements in
   $E^{\infty}_{0,\geq n-k}=H_{\geq n-k}(L;\mathbb{Z}_2)$. Keeping in
   mind that $[*]$ is induced from $*$ this means that $$x = a +
   x_{(p')}$$
   for some elements $a \in \mathcal{A}$ and $x_{(p')} \in
   \mathcal{F}_{p'}H_i(\mathcal{C},d)$. By the induction hypothesis
   $x_{(p')} \in \mathcal{A}$ hence $x \in \mathcal{A}$ too. The
   desired statement follows now by induction.

   The analogous statements of the proposition for $Q^{+}H(L)$ are
   proved essentially in the same way. However, it is important to
   note that alternative~\ref{I:QH-0-or-all-1} holds only for $QH(L)$
   (in fact $Q^{+}H(L)$ cannot vanish). The reason for this lies in
   ``\textbf{Case I.}'' in the proof above where we had to use
   negative powers of $t$.

   \

   {\bf B.} Here is now the second argument based on the minimal model machinery
   from \S \ref{subsec:minimal}. Consider the pearl complex $\mathcal{C}(f,J)$ and recall
   from \S\ref{subsec:minimal} that there exists a chain complex $(\mathcal{C}_{min},\delta)$,
   unique up to isomorphism, and chain morphisms $\phi:\mathcal{C}(f,J)\to \mathcal{C}_{min}$, $\psi:\mathcal{C}_{min}\to \mathcal{C}(f,J)$ so that $\phi\circ\psi=id$, $\mathcal{C}_{min}=H_{\ast}(L;\Z_{2})\otimes \La^{+}$, $\delta_{0}=0$ (where $\delta_{0}$ is obtained from $\delta$ by putting $t=0$)
   and $\phi$ ($\phi_{0}$), $\psi$ ($\psi_{0}$) induce isomorphisms  in quantum (respectively, Morse) homology.
   By Remark \ref{rem:product_min} the quantum product in $\mathcal{C}(f,J)$ can be transported by the
   morphisms $\phi$ and $\psi$ to a product $\ast: \mathcal{C}_{min}\otimes\mathcal{C}_{min}\to \mathcal{C}_{min}$ which is a chain map and a quantum deformation of the singular intersection
   product (notice though that, as the maps $\phi$ and $\psi$ are not canonical, this product is not canonical either at the chain level).

   We will now show by induction that either $QH(L)=0$ or $\delta=0$.  Let $x\in H_{n-k+s}(L)$
  with $s\geq 0$. We indentify $H_{\ast}(L; \Z_{2})$ with the generators of $\mathcal{C}_{min}$ and
  we notice that $N_{L}>k$ implies for degree reasons that $\delta x=0$ when $|x|>n-k$ and
  $\delta x=\epsilon_{x} [L]t$ when $|x|=n-k$, $N_{L}=k+1$ with $\epsilon_{x}\in \{0,\}$.
  If for some such $x$ we have $\epsilon_{x}=1$, then, by Remark \ref{rem:product_min},
  we deduce $QH(L)=0$.  Thus, we now assume $QH(L)\not=0$ and we assume, by induction, that
  $\delta y=0$ for all $y\in H_{\ast}(L; \Z_{2})$ such that $|y|> n-k-s$, $s\geq 1$. Consider
  $x\in H_{n-k-s}(L; \Z_{2})$  so that $x=x_{1}\cdot x_{2}\cdot\ldots \cdot x_{r}$ with
  $x_{i}\in H_{\geq n-k}(L; \Z_{2})$. We then have $\delta (x_{i})=0$ and
  we write  $\delta(x_{1}\ast x_{2}\ast\ldots\ast x_{r})=
  \sum_{i} x_{1}\ast \ldots \delta(x_{i})\ast\ldots\ast x_{r}=0$. At the same time
  \begin{equation}\label{eq:prod_quantprod}
  x_{1}\ast x_{2}\ast\ldots\ast x_{r}= x+\sum_{j}z_{j}t^{j}
  \end{equation} with $z_{j}\in H_{> n-k-s}(L)$ so that, by the induction hypothesis, $\delta (z_{j})=0$.
  We conclude $\delta x=0$. But given our assumption on the structure of the singular homology
  of $L$, this implies that $\delta=0$. This means that $Q^{+}H(L)\cong H(L)\otimes\La^{+}$.
  The equation (\ref{eq:prod_quantprod}) immediately implies the rest of the statement at (2).
  \end{proof}

\begin{rem}
   a. It is important to notice that the dichotomy that we have proved
   when $N_L=k+1$ depends on the fact that we work over a field.  This
   appears in ``case \textbf{I}'' in the proof above.

   b. As shown by Cho~\cite{Cho:Clifford}, the Clifford torus has the
   property that for a certain choice of spin structure the associated
   Floer homology with rational coefficients vanishes.
   In~\cite{Cho-Oh:Floer-toric, Cho:products} there are some examples
   of non-monotone tori which satisfy similar dichotomy type
   properties (the argument used there is different from the one
   here).

   c. The proof above extends in obvious ways to other examples of
   Lagrangians with particular singular homology.

   d. Some partial results of the type above have also been obtained
   by Buhovsky in~\cite{Bu:products}.
\end{rem}

\subsubsection{Criteria for vanishing and non-vanishing of Floer homology}
\label{Sb:criteria-QH}

Here is a related but slightly different, and possibly more explicit
point of view on the same phenomenon from
Proposition~\ref{cor:floer=0-or-all} (see also Remark \ref{rem:product_min}).
  Let $L^n \subset (M^{2n},
\omega)$ be a monotone Lagrangian submanifold with $N_L \geq 2$.
Denote by $H_2^D \subset H_2(M,L;\mathbb{Z})$ the image of the
Hurewicz homomorphism $\pi_2(M,L) \to H_2(M,L;\mathbb{Z})$. Denote by
$\partial:H_2(M,L;\mathbb{Z}) \to H_1(L;\mathbb{Z})$ the boundary
homomorphism and by $\partial_{\mathbb{Z}_2}:H_2(M,L;\mathbb{Z}) \to
H_1(L;\mathbb{Z}_2)$ the composition of $\partial$ with the reduction
mod $2$, $H_1(L;\mathbb{Z}) \to H_1(L;\mathbb{Z}_2)$. Given $A \in
H_2^D$ and $J \in \mathcal{J}(M,\omega)$ consider the evaluation map
$$ev_{A,J}: (\mathcal{M}(A,J) \times \partial D)/G \longrightarrow L,
\quad ev_{A,J}(u,p) = u(p),$$
where $G = \textnormal{Aut}(D) \cong
PSL(2,\mathbb{R})$ is the group of biholomorphisms of the disk.

For every $J \in \mathcal{J}(M,\omega)$ let $\mathcal{E}_2(J)$ be the
set of all classes $A \in H_2^D$ with $\mu(A)=2$ for which there exist
$J$-holomorphic disks with boundaries on $L$ in the class $A$:
$$\mathcal{E}_2(J) = \{ A \in H_2^D \mid \mu(A)=2, \quad
\mathcal{M}(A,J) \neq \emptyset \}.$$
Define:
$$\mathcal{E}_2 = \bigcap_{J \in \mathcal{J}(M,\omega)}
\mathcal{E}_2(J).$$
Standard arguments show that:
\begin{enumerate}
  \item $\mathcal{E}_2(J)$ is a finite set for every $J$.
  \item There exists a second category subset
   $\mathcal{J}_{\textnormal{reg}} \subset \mathcal{J}(M,\omega)$ such
   that for every $J \in \mathcal{J}_{\textnormal{reg}}$,
   $\mathcal{E}_2(J) = \mathcal{E}_2$. In other words, for generic
   $J$, $\mathcal{E}_2(J)$ is independent of $J$.
  \item For every $J \in \mathcal{J}$ and every $A \in
   \mathcal{E}_2(J)$ the space $\mathcal{M}(A,J)$ is compact and all
   disks $u \in \mathcal{M}(A,J)$ are simple.
  \item For $J \in \mathcal{J}_{\textnormal{reg}}$ and $A \in
   \mathcal{E}_2$, the space $(\mathcal{M}(A,J) \times \partial D)/G$
   is a compact smooth manifold without boundary. Its dimension is $n
   = \dim L$. In particular, for generic $x \in L$, the number of
   $J$-holomorphic disks $u \in \mathcal{M}(A,J)$ with $u(\partial D)
   \ni x$ is finite.
  \item For every $A \in \mathcal{E}_2$ and $J_0, J_1 \in
   \mathcal{J}_{\textnormal{reg}}$ the manifolds $(\mathcal{M}(A,J_0)
   \times \partial D)/G$ and $(\mathcal{M}(A,J_1) \times \partial
   D)/G$ are cobordant via a compact cobordism. Moreover, the
   evaluation maps $ev_{A,J_0}$, $ev_{A,J_1}$ extend to this
   cobordism, hence $\deg_{\mathbb{Z}_2} ev_{A,J_0} =
   \deg_{\mathbb{Z}_2} ev_{A,J_1}$.  In other words
   $\deg_{\mathbb{Z}_2} ev_{A,J}$ depends only on $A \in
   \mathcal{E}_2$.
  \item In fact, the set $\mathcal{J}_{\textnormal{reg}}$ above can be
   taken to be the set of all $J \in \mathcal{J}(M,\omega)$ which are
   regular for all classes $A \in H_2^D$ in the sense that the
   linearization of the $\overline{\partial}_J$ operator is surjective
   at every $u \in \mathcal{M}(A,J)$.
\end{enumerate}

Let $J \in \mathcal{J}_{\textnormal{reg}}$ and let $x \in L$ be a
generic point. Define a one dimensional $\mathbb{Z}_2$-cycle
$\delta_x(J)$ to be the sum of the boundaries of all $J$-holomorphic
disks with $\mu=2$ whose boundaries pass through $x$. Of course, if a
disk meets $x$ along its boundary several times we take its boundary
in the sum with appropriate multiplicity. Thus the precise definition
is:
\begin{equation} \label{Eq:delta-x}
   \delta_x(J) = \sum_{A \in \mathcal{E}_2} \; \sum_{(u,p) \in
     ev_{A,J}^{-1}(x)} u(\partial D).
\end{equation}
By the preceding discussion the homology class $D_1 = [\delta_x(J)]
\in H_1(L;\mathbb{Z}_2)$ is independent of $J$ and $x$. In fact
\begin{equation} \label{Eq:D1}
   D_1 = \sum_{A \in \mathcal{E}_2} (\deg_{\mathbb{Z}_2} ev_{A,J})
   \partial_{\mathbb{Z}_2}A.
\end{equation}

In view of the proof of Proposition~\ref{cor:floer=0-or-all} the next
result follows easily.

\begin{prop} \label{P:criterion-QH=0-1}
   Let $L \subset (M,\omega)$ be a monotone Lagrangian submanifold
   with $N_L \geq 2$.
   \begin{enumerate}
     \item If $D_1 \neq 0$ then $QH_*(L)=0$. \label{I:criterion-QH-1}
     \item Suppose that $D_1 = 0$ and $H_*(L;\mathbb{Z}_2)$ is
      generated as an algebra by $H_{n-1}(L;\mathbb{Z}_2)$ with
      respect to the classical cap product.  Then $QH_*(L) \cong
      (H(L;\mathbb{Z}_2) \otimes \Lambda)_*$. This isomorphism is
      neither canonical nor multiplicative. However for $* \geq n-1$,
      there is a canonical injection $H_{\geq n-1}(L;\mathbb{Z}_2)
      \hookrightarrow QH_{\geq n-1}(L)$.
      \label{I:criterion-QH-2}
   \end{enumerate}
   In particular, if $H_*(L;\mathbb{Z}_2)$ is generated as an algebra
   by $H_{n-1}(L;\mathbb{Z}_2)$ then $QH_*(L)$ can be either
   $(H(L;\mathbb{Z}_2) \otimes \Lambda)_*$ or $0$ according to whether
   $D_1$ vanishes or not. When $D_1=0$ all the above continues to hold
   for $QH(L)$ replaced by $Q^{+}H(L)$ and $\Lambda$ replaced by
   $\Lambda^{+}$.
\end{prop}
\begin{rem} \label{R:nu-classes}
   When $D_1=0$ but $H_*(L;\mathbb{Z}_2)$ is not generated as an
   algebra by $H_{n-1}(L;\mathbb{Z}_2)$, the theorem does not say
   anything on $QH(L)$.  In this case it is possible to define in a
   similar way higher classes $D_j \in H_j(L;\mathbb{Z}_2)$, $1 \leq j
   \leq N_L-1$, which sometimes give more information.
\end{rem}

\begin{proof}[Proof of Proposition~\ref{P:criterion-QH=0-1}]
   Choose a generic $J \in \mathcal{J}(M,\omega)$. Let $f:L \to
   \mathbb{R}$ be a generic Morse function with precisely one local
   minimum $x \in L$ and fix a generic Riemannian metric on $L$.
   Denote by $(CM_*(f), \partial_0)$, $(\mathcal{C}_*(f,J),d)$ the
   Morse and pearl complexes associated to $f$, $J$ and the chosen
   Riemannian metric.

   For degree reasons the restriction of $d$ to $CM_{n-1}(f) \subset
   \mathcal{C}_{n-1}(f,J)$ is given by $d=\partial_0 + \partial_1 t$,
   where $\partial_1:CM_{n-1}(f) \to CM_n(f)=\mathbb{Z}_2 x$ counts
   pearly trajectories with holomorphic disks of Maslov index $2$ (of
   course, if $N_L > 2$ then $\partial_1=0$ and also $\delta_x(J)=0$,
   $D_1=0$). Since $x$ is a maximum of $f$ no $-\textnormal{grad}(f)$
   trajectories can enter $x$ (i.e. $W_x^s = \{x\}$). Therefore for
   every $y \in \textnormal{Crit}_{n-1}(f)$ we have
   \begin{equation} \label{Eq:del-1}
      \partial_1 y = \#_{\mathbb{Z}_2} \bigl( W_y^u \cap
      \delta_x(J)\bigr)x.
   \end{equation}

   We prove statement~\ref{I:criterion-QH-1}. Suppose that $D_1 \neq
   0$. By Poincar\'{e} duality there exists an $(n-1)$-dimensional
   cycle $C$ in $L$ such that
   $$\#_{\mathbb{Z}_2} C \cap \delta_x(J) \neq 0.$$
   Let $z \in
   CM_{n-1}(f)$ be a $\partial_0$-cycle representing $[C] \in
   H_{n-1}(L;\mathbb{Z}_2)$. Then $$d(z) = \partial_1(z)\otimes t =
   \#_{\mathbb{Z}_2} \bigl( W_z^u \cap \delta_x(J)\bigr)x \otimes t =
   \#_{\mathbb{Z}_2} \bigl(C \cap \delta_x(J)\bigr)x\otimes t = ax
   \otimes t$$
   for some non-zero scalar $a$. (Of course, $a \neq 0$ is
   the same as $a=1$ here, since we work over $\mathbb{Z}_2$. However
   we wrote $ax$ to emphasize that the argument works well over every
   field.)  It follows that $[x] = 0 \in QH_n(L)$. But $[x]$ is the
   unity of $QH_*(L)$, hence $QH_*(L)=0$.

   We prove statement~\ref{I:criterion-QH-2}.We will use here an
   argument involving spectral sequences, similar to the proof of
   Proposition~\ref{cor:floer=0-or-all}. Recall from the proof of
   Proposition~\ref{cor:floer=0-or-all} that the degree filtration
   gives rise to a spectral sequence $\{E^r_{p,q}, d_r\}_{p,q}$ that
   converges to $QH_*(L)$ and $E^1_{p,q} = H_{p+q-pN_L}(L;
   \mathbb{Z}_2) t^{-p}$. A simple computation shows that the
   differential $d_1$ is induced from the operator $\partial_1$
   mentioned earlier in the proof. Moreover, as explained in the proof
   of Proposition~\ref{cor:floer=0-or-all}, the quantum product $*$
   endows $\{E^r_{p,q}, d_r\}_{r \geq 1}$ with a multiplicative
   structure which coincides for $r=1$ with the classical cap product
   $\cap$. In particular $d_1$ satisfies Leibniz rule with respect to
   $\cap$, and the $d_r$, $r \geq 2$, satisfy too Leibniz rule with
   respect to the products induced on $E^r$. Since
   $H_*(L;\mathbb{Z}_2)$ is generated by $H_{n-1}(L;\mathbb{Z}_2)$ we
   conclude that $\{E_1^{p,q}\}_{p+q-pN_L=n-1}$ generate with respect
   to the cap product the whole of $E^1_{*,*}$. As $D_1=0$ we obtain
   from formula~\eqref{Eq:del-1} that $d_1$ vanishes on $E^1_{p,q}$
   whenever $p+q-pN_L=n-1$. It follows that $d_1=0$ on all
   $E_1^{p,q}$.

   The above implies that $E^2_{p,q}=E^1_{p,q}$ and the product
   induced on $E^2$ is still the cap product. In particular
   $E^2_{*,*}$ is generated with respect to $\cap$ by
   $\{E^2_{p,q}\}_{p+q-pN_L=n-1}$. For degree reasons $d_2$ vanishes
   on $\{E^2_{p,q}\}_{p+q-pN_L=n-1}$. As $d_2$ satisfies Leibniz rule
   too it follows that $d_2=0$ everywhere. Thus
   $E^3_{*,*}=E^2_{*,*}=E^1_{*,*}$. The same argument shows that
   $d_r=0$ for every $r\geq 2$ hence $E^r_{*,*} = E^1_{*,*}$. In other
   words the spectral sequence collapses at level $r=1$. Since this
   spectral sequence converges to $H_*(\mathcal{C}(f,J),d) = QH_*(L)$,
   we conclude that $QH_*(L) \cong (H(L;\mathbb{Z}_2) \otimes
   \Lambda)_*$.
\end{proof}

Let us turn to some examples. First of all, if $N_L \geq 3$ then
$\mathcal{E}_2 = \emptyset$ hence $D_1=0$. Therefore if
$H_*(L;\mathbb{Z}_2)$ is generated as an algebra by
$H_{n-1}(L;\mathbb{Z}_2)$ we must have $QH_*(L) \cong
(H(L;\mathbb{Z}_2) \otimes \Lambda)_*$. An example of such a
Lagrangian is $\mathbb{R}P^n \subset {\mathbb{C}}P^n$, $n \geq 2$.

\begin{ex} \label{Sb:tclif-QH-1}
   Let $\mathbb{T}_{\textnormal{clif}} = \{ [z_0: \cdots :z_n] \in
   {\mathbb{C}}P^n \mid |z_0|= \cdots =|z_n| \}$ be the
   $n$-dimensional Clifford torus.  This is a monotone Lagrangian
   torus with $N_L=2$.  The Floer homology of
   $\mathbb{T}_{\textnormal{clif}}$ was computed by
   Cho~\cite{Cho:Clifford} by a direct computation of the Floer
   complex.  Below we will review this computation from the
   perspective of Proposition~\ref{P:criterion-QH=0-1}.

   A simple computation shows that $H_2^D \cong \pi_2({\mathbb{C}}P^n,
   \mathbb{T}_{\textnormal{clif}}) \cong \mathbb{Z} A_0 \oplus \cdots
   \oplus \mathbb{Z} A_n$, where $A_i$, is represented by the map
   $v_i: (D, \partial D) \to ({\mathbb{C}}P^n,
   \mathbb{T}_{\textnormal{clif}})$ given by $v_i(z) = [1: \cdots: z:
   \cdots: 1]$ (here the $z$ stands in the $i$'th entry).  A
   straightforward computation shows that $\mu(A_i)=2$ for every $i$.

   Let $J_0$ be the standard complex structure of ${\mathbb{C}}P^n$.
   We will use the following facts proved by Cho~\cite{Cho:Clifford}:
   \begin{enumerate}
     \item $\mathcal{E}_2(J_0) =\{A_0, \ldots, A_n\}$.
     \item $J_0$ is regular for each of the classes $A_i$.
     \item $ev_{A_i,J_0}:(\mathcal{M}(A_i,J_0) \times \partial D)/G
      \to \mathbb{T}_{\textnormal{clif}}$ is a diffeomorphism, hence
      $\deg_{\mathbb{Z}_2} ev_{A_i,J_0}=1$. In fact, given a point
      $\xi=[\xi_0: \cdots: \xi_n] \in \mathbb{T}_{\textnormal{clif}}$,
      the unique disk (up to reparametrization) $u:(D, \partial D) \to
      ({\mathbb{C}}P^n, \mathbb{T}_{\textnormal{clif}})$ in the class
      $A_i$ with $u(\partial D) \ni \xi$ is given by $u(z) = [\xi_0:
      \cdots :\xi_{i-1} : \xi_i z : \xi_{i+1}: \cdots : \xi_n]$.
   \end{enumerate}
   It follows from the discussion in \S\ref{Sb:criteria-QH} that
   $\mathcal{E}_2 = \{A_0, \ldots, A_n\}$ and that for every $J \in
   \mathcal{J}_{\textnormal{reg}}$, $\deg_{\mathbb{Z}_2}
   ev_{A_i,J}=1$.  A simple computation shows that $\partial A_0 +
   \cdots + \partial A_n = 0 \in H_1(L;\mathbb{Z})$ hence we have:
   $$D_1 = \sum_{i=0}^n (\deg_{\mathbb{Z}_2} ev_{A_i,J})
   \partial_{\mathbb{Z}_2} A_i = 0.$$
   By
   Proposition~\ref{P:criterion-QH=0-1} we have
   \begin{equation} \label{Eq:QH-clif}
      QH_*(\mathbb{T}_{\textnormal{clif}}) \cong
      (H(\mathbb{T}_{\textnormal{clif}};\mathbb{Z}_2) \otimes \Lambda)_*,
      \quad
      Q^{+}H_*(\mathbb{T}_{\textnormal{clif}}) \cong
      (H(\mathbb{T}_{\textnormal{clif}};\mathbb{Z}_2)
      \otimes \Lambda^{+})_*,
   \end{equation}
   and there are canonical injections
   \begin{equation} \label{Eq:can-inj-clif}
      H_{n-1}(\mathbb{T}_{\textnormal{clif}};\mathbb{Z}_2)
      \hookrightarrow QH_{n-1}(\mathbb{T}_{\textnormal{clif}}), \quad
      H_n(\mathbb{T}_{\textnormal{clif}};\mathbb{Z}_2) \hookrightarrow
      QH_n(\mathbb{T}_{\textnormal{clif}})
   \end{equation}
   and similarly for $Q^{+}H(\mathbb{T}_{\textnormal{clif}})$.
\end{ex}

\subsection{Lagrangian submanifolds of ${\mathbb{C}}P^n$} \label{S:cpn}

Endow ${\mathbb{C}}P^n$ with the standard K\"{a}hler symplectic
structure $\omega_{\textnormal{FS}}$, normalized so that
$\int_{\mathbb{C}P^1} \omega_{\textnormal{FS}}=\pi$. Let $L \subset
{\mathbb{C}}P^n$ be a monotone Lagrangian submanifold with minimal
Maslov number $N_L \geq 2$. Below we will carry out computations
involving the quantum homology $QH({\mathbb{C}}P^n)$ and the Floer
homology $QH(L)$. We will work with the following (simplified) version
of the Novikov ring $\Lambda = \mathbb{Z}_2[t,t^{-1}]$ where $\deg t =
-N_L$.  Put $QH({\mathbb{C}}P^n) = H({\mathbb{C}}P^n;\mathbb{Z}_2)
\otimes \Lambda$ with the grading induced from both factors. Denote by
$h \in H_{2n-2}({\mathbb{C}}P^n;\mathbb{Z}_2)$ the class of the
hyperplane, and by $u \in H_{2n}({\mathbb{C}}P^n;\mathbb{Z}_2)$ the
fundamental class. It is well known that (see
e.g.~\cite{McD-Sa:Jhol-2}):
\begin{equation} \label{Eq:QHcpn}
   h^{*j} =
   \begin{cases}
      h^{\cap j}, & 0 \leq j \leq n \\
      u \otimes t^{\frac{2(n+1)}{N_L}}, & j=n+1
   \end{cases}
\end{equation}
Note that our choice of grading is somewhat different than the
convention usually taken in quantum homology theory. For example, if
$N_L = n+1$ then $\deg t = -(n+1)$ and we get from~\eqref{Eq:QHcpn}
that $h^{*(n+1)} = u \otimes t^2$ ({\em not} $u \otimes t$ !). Usually
in the theory of quantum homology the degree of $t$ is taken to be
$-2N_M$ where $N_M$ is the minimal Chern number of $(M,\omega)$.  Here
we have defined $\deg t = -N_L$ in order to keep compatibility with
the Novikov ring used for Floer homology. Note however that we have
$N_L | 2N_M$ thus our ring $\Lambda_*$ is obtained from the
``conventional'' Novikov ring by a variable change.

It follows from~\eqref{Eq:QHcpn} that $h$ is an invertible element.
Therefore by Theorem~\ref{thm:alg_main} we have:
\begin{cor} \label{T:Lag-CPn}
   Let $L \subset {\mathbb{C}}P^n$ be a monotone Lagrangian with $N_L
   \geq 2$. Then $QH_*(L)$ is $2$-periodic, i.e. $QH_i(L) \cong
   QH_{i-2}(L)$ for every $i \in \mathbb{Z}$. In fact the homomorphism
   $QH_i(L) \to QH_{i-2}(L)$ given by $\alpha \mapsto h*\alpha$ is an
   isomorphism for every $i \in \mathbb{Z}$.
\end{cor}
\begin{rem}
   \begin{enumerate}
     \item The first part of Theorem~\ref{T:Lag-CPn} was proved before
      by Seidel using the theory of graded Lagrangian
      submanifolds~\cite{Se:graded}.  The $2$-periodicity
      in~\cite{Se:graded} follows from the fact that ${\mathbb{C}}P^n$
      admits a Hamiltonian circle action which induces a shift by $2$
      on graded Lagrangian submanifolds. Note that this is compatible
      with our perspective since that $S^1$-action gives rise to an
      invertible element in $QH({\mathbb{C}}P^n)$ (the Seidel
      element~\cite{Se:pi1, McD-Sa:Jhol-2}) whose degree is exactly
      $2n$ minus the shift induced by the $S^1$-action. In our case
      the Seidel element turns out to be $h$.
     \item Let $\widetilde{\Lambda}=\mathbb{Z}_2[t]]$ be the ring of
      formal Laurent series with finitely many negative terms, i.e.
      elements of $\widetilde{\Lambda}$ are of the form
      $p(t)=\sum_{i=N}^{\infty} a_i t^i$, $N \in \mathbb{Z}$. Note
      that $\widetilde{\Lambda}$ is a field. If we define $QH$ with
      coefficients in $\widetilde{\Lambda}$ then
      $QH({\mathbb{C}}P^n;\widetilde{\Lambda})$ is isomorphic to the
      ring $\widetilde{\Lambda}[x]/\{x^{n+1}=t\}$. It is easy to see
      that this ring is in fact a field. Thus if we define
      $QH(L;\widetilde{\Lambda})$ to be $QH$ with coefficients in
      $\widetilde{\Lambda}$ we obtain for every $0 \neq \alpha \in
      QH(L;\widetilde{\Lambda})$ an {\em injective} homomorphism
      $QH({\mathbb{C}}P^n; \widetilde{\Lambda}) \hookrightarrow
      QH(L;\widetilde{\Lambda})$, defined by $a \mapsto a* \alpha$.
   \end{enumerate}
\end{rem}

\subsubsection{The Clifford torus}\label{Ss:clifford}

We now consider the $2$-dimensional Clifford torus
$\mathbb{T}^2_{\textnormal{clif}} \subset \C P^{2}$ and compute all
our structures in this case. We denote $\Lambda = \mathbb{Z}_2[t,
t^{-1}]$, $\Lambda^{+}= \mathbb{Z}_2 [t]$ where $deg(t)=-2$. We denote
by $h\in H_{2}(\C P^{2};\mathbb{Z}_2)$ the generator. Recall
from~\eqref{Eq:QH-clif} that $QH_*(\mathbb{T}^2_{\textnormal{clif}})
\cong (H(\mathbb{T}^2_{\textnormal{clif}};\mathbb{Z}_2) \otimes
\Lambda)_*$. In particular:
\begin{align}
   & QH_0(\mathbb{T}^2_{\textnormal{clif}}) \cong
   H_0(\mathbb{T}^2_{\textnormal{clif}};\mathbb{Z}_2) \oplus
   H_2(\mathbb{T}^2_{\textnormal{clif}};\mathbb{Z}_2) t,
   \label{Eq:QH-0-iso}\\
   & QH_1(\mathbb{T}_{\textnormal{clif}}) \cong
   H_1(\mathbb{T}_{\textnormal{clif}};\mathbb{Z}_2). \label{Eq:QH-1-iso}
\end{align}
Recall from~\eqref{Eq:can-inj-clif} that the isomorphism
in~\eqref{Eq:QH-1-iso} is canonical and the second summand
in~\eqref{Eq:QH-0-iso} is canonical too. (Note however that the first
summand in~\eqref{Eq:QH-0-iso} is {\em not} canonical. See \S
\ref{Sb:clif-enum} for more details on that).

\begin{prop}\label{prop:clifford}
   Let $w \in H_2(\mathbb{T}^2_{\textnormal{clif}};\mathbb{Z}_2)$ be
   the fundamental class. There are generators $a,b\in
   H_{1}(\mathbb{T}^2_{\textnormal{clif}}, \mathbb{Z}_2)$, and $m \in
   QH_0(\mathbb{T}^2_{\textnormal{clif}}) \cong
   (H_{\ast}(\mathbb{T}^2_{\textnormal{clif}};\mathbb{Z}_2)\otimes
   \Lambda)_{0}$ which together with $w$ generate
   $QH(\mathbb{T}^2_{\textnormal{clif}})$ as a $\Lambda$-module and
   verify the following relations:
   \begin{itemize}
     \item[i.] $a\ast b=m+ w t$, $b\ast a=m$, $a\ast a=b\ast b= w t$,
      $m*m = mt + wt^2$.
     \item[ii.] $h\ast a=at$, $h\ast b=bt$, $h\ast w=wt$, $h\ast m = mt$.
     \item[iii.] $i_{L}(m)=[pt]+ht+[\C P^{2}]t^{2}$,
      $i_{L}(a)=i_{L}(b)=i_{L}(w)=0$.
   \end{itemize}
   All the above continues to hold for the positive version of $QH$,
   namely with $\Lambda$ replaced by $\Lambda^{+}$ and
   $QH(\mathbb{T}^2_{\textnormal{clif}})$ replaced by
   $Q^{+}H(\mathbb{T}^2_{\textnormal{clif}})$.
\end{prop}
\begin{remsnonum}
   \begin{enumerate}
     \item As the formulae in i clearly show, the Lagrangian quantum
      product is {\em not commutative} (even when working over
      $\mathbb{Z}_2$).
     \item Point i of Proposition~\ref{prop:clifford} has been
      obtained before by Cho~\cite{Cho:products} by a different
      approach. From the perspective of that paper the Clifford torus
      is a special case of a torus which appears as a fibre of the
      moment map defined on a toric variety. See
      also~\cite{Cho-Oh:Floer-toric} for related results in this
      direction.
   \end{enumerate}
\end{remsnonum}

\begin{proof}[Proof of Proposition~\ref{prop:clifford}]
   We will use the following two geometric properties of the Clifford
   torus. The first is that, through each point of the Clifford torus,
   there are three different pseudo-holomorphic disks of Maslov index
   two.  They belong to three families that we denote by $\gamma_{1}$,
   $\gamma_{2}$ and $\gamma_{3}$.  Up to a possible change of basis,
   we may assume that the homotopy class of the boundaries of the
   elements in $\gamma_{1}$ is $a$, for $\gamma_{2}$ the same class is
   $b$ and for $\gamma_{3}$ this class is $-a-b$. See
   figure~\ref{f:clif-3-disks}.
   \begin{figure}[htbp]
      \begin{center}
         \epsfig{file=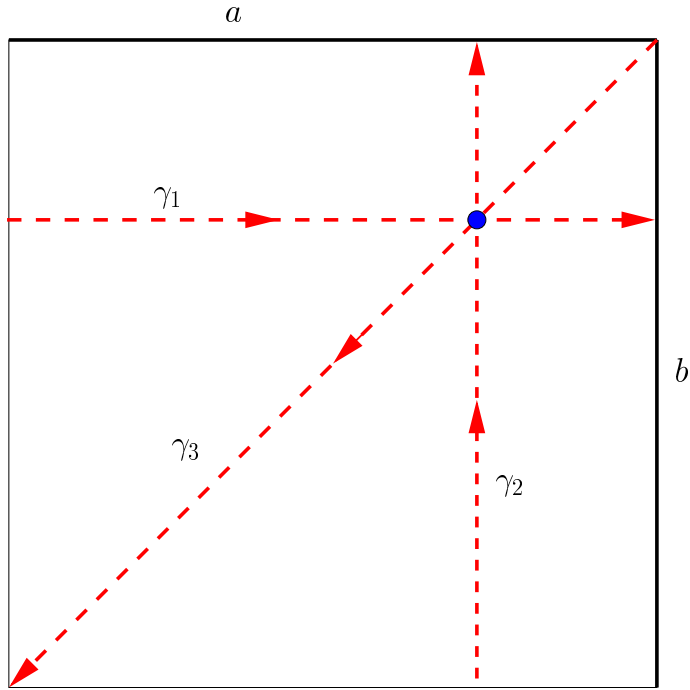, width=0.6\linewidth}
      \end{center}
      \caption{The boundaries of the 3 holomorphic disks with
        $\mu=2$ through every point on
        $\mathbb{T}^2_{\textnormal{clif}}$}
      \label{f:clif-3-disks}
   \end{figure}
   The second geometric fact is that there is a symplectomorphism
   homotopic to the identity, $\bar{\phi}: \C P^2\to \C P^2$, whose
   restriction to $\mathbb{T}^2_{\textnormal{clif}}$ is the
   permutation of the two factors in $\mathbb{T}^2_{\textnormal{clif}}
   \approx S^{1}\times S^{1}$.  We now consider a perfect Morse
   function $f:\mathbb{T}^2_{\textnormal{clif}}\to \R$ and, by a
   slight abuse in notation, we let its minimum be $m$, we let the
   maximum be $w$ and we let the two critical points of index $1$ be
   denoted by $a'$ and $b'$ so that the unstable manifold of $a'$ has
   the homotopy type $a\in
   H_{1}(\mathbb{T}^2_{\textnormal{clif}};\mathbb{Z}_2)$ and,
   similarly, the unstable manifold of the critical point $b'$ has
   homotopy type $b$.  We denote the disk in the family $\gamma_{i}$
   that passes through $w$ by $d_{i}$. See figure~\ref{f:clif-traj-1}.
   \begin{figure}[htbp]
      \begin{center}
         \epsfig{file=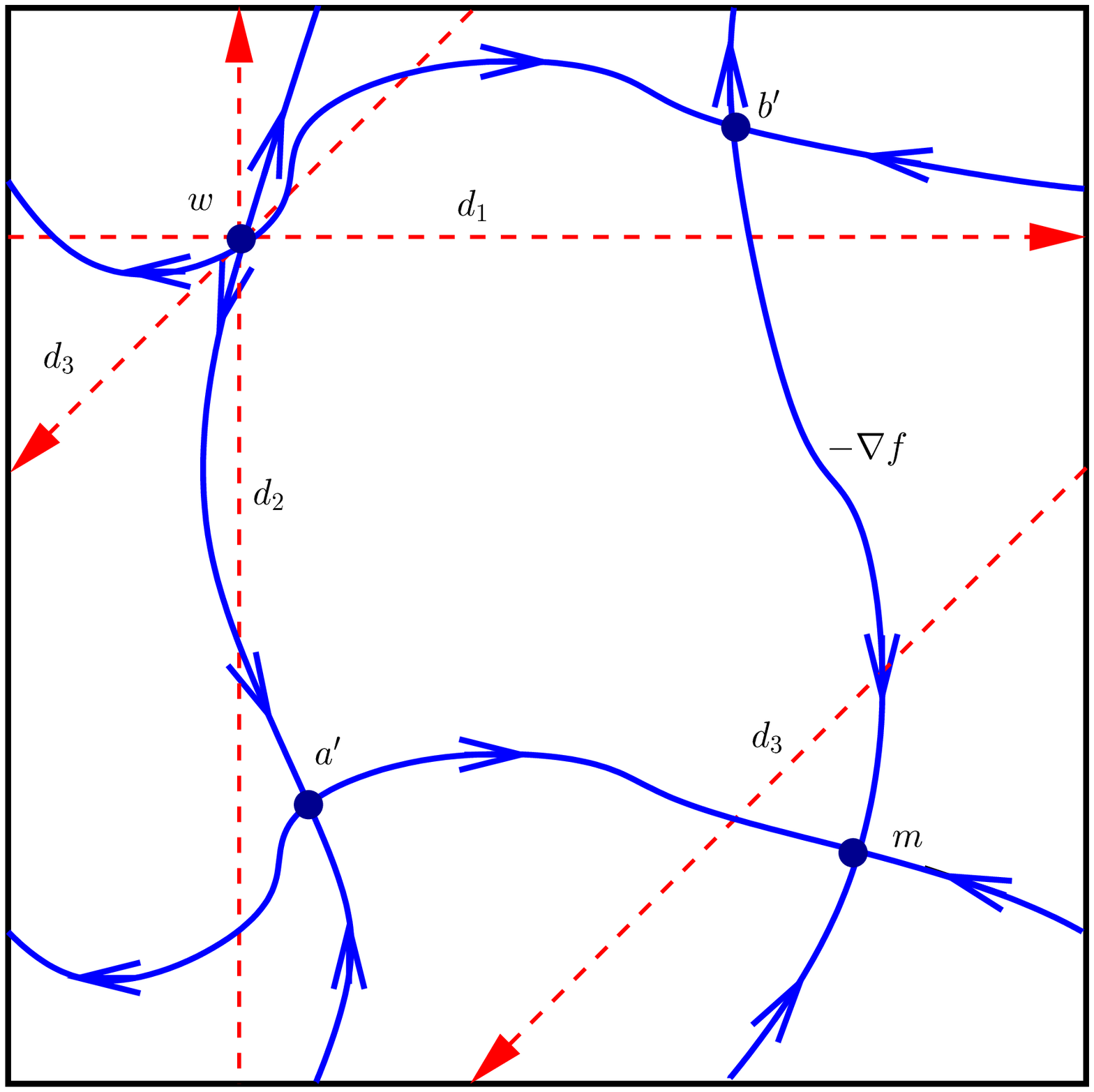, width=0.6\linewidth}
      \end{center}
      \caption{Trajectories of $-\nabla f$ and holomorphic disks on
        $\mathbb{T}^2_{\textnormal{clif}}$.}
      \label{f:clif-traj-1}
   \end{figure}
   By possibly perturbing the function $f$ slightly we may assume that
   the unstable manifold of $a'$ intersects $d_{2}$ and $d_{3}$ in a
   single point and is disjoint from $d_{1}$.  Similarly, we may
   assume that the unstable manifold of $b'$ intersect $d_{1}$ and
   $d_{3}$ in a single point and is disjoint from $d_{2}$. With these
   choices the pearl complex $(\mathcal{C}(f, J, \rho), d)$ is well
   defined. Here we take $J$ to be the standard complex structure of
   ${\mathbb{C}}P^2$ or a small perturbation of it and $\rho$ a
   generic Riemannian metric on $\mathbb{T}^2_{\textnormal{clif}}$.
   We claim that the differential $d$ of the pearl complex vanishes.
   Indeed for dimension reasons we can write $d = \partial_0 +
   \partial_1 t$, where $\partial_0$ is the Morse differential and
   $\partial_1: \mathcal{C}_* \to \mathcal{C}_{*+1}$ is an operator
   that counts the contribution of the pearly trajectories involving
   $J$-holomorphic disks with Maslov index $2$. Since $f$ is perfect
   $\partial_0=0$ hence $d=\partial_1 t$. As we have already seen
   before $QH_*(\mathbb{T}^2_{\textnormal{clif}}) \cong
   (H(\mathbb{T}^2_{\textnormal{clif}};\mathbb{Z}_2) \otimes
   \Lambda)_*$. It follows that $\partial_1 = 0$ too since otherwise
   we would have $\dim_{\mathbb{Z}_2}
   QH_i(\mathbb{T}^2_{\textnormal{clif}}) < \dim_{\mathbb{Z}_2}
   (H(\mathbb{T}^2_{\textnormal{clif}};\mathbb{Z}_2) \otimes
   \Lambda)_i$ for some $i$, a contradiction. This proves that $d=0$.

   It is instructive to give a more direct proof of the fact that
   $d=0$ based on the specific knowledge of the $\mu=2$ -- holomorphic
   disks.  For this purpose we first note that $da'=0=db'$. This is
   because the only two possibilities for $da'$ are $da'=0$ and
   $da'=wt$ and, as there are precisely two disks that go through $w$
   and intersect the unstable manifold of $a'$ and each of them
   intersects it in exactly one point, we see that we are in the first
   case. The same argument applies to $b'$. A similar computation
   shows that $dm=0$.  Finally, $dw = 0$ for degree reasons.

   Summarizing the above, $d=0$ hence $QH_*(f,J,\rho) =
   H_*(\mathcal{C}(f,J,\rho), d) = \mathcal{C}_*(f,J,\rho)$. From now
   on we will view $m, a', b', w$ as generators (over $\Lambda$) of
   $QH_*(\mathbb{T}^2_{\textnormal{clif}})$. Note that $m$ depends on
   the choice of $f$ in the sense that if we take another perfect
   Morse function $g$ with minimum $\tilde{m}$ then $\tilde{m}$ might
   give an element of $QH_0(\mathbb{T}^2_{\textnormal{clif}})$ which
   is different than $m$. On the other hand $a', b', w \in QH$ are
   canonical.

   We now discuss the product. For degree reasons we have $a'\ast b'=
   m+ \epsilon wt$, $b'\ast a'=m+\epsilon' wt$ with $\epsilon,
   \epsilon' \in \mathbb{Z}_2$. Of course, $\epsilon$ is the number
   modulo $2$ of disks going - in order ! - through the following
   points: one point in the unstable manifold of $a'$ then $w$ and,
   finally one point in the unstable manifold of $b'$.  Similarly,
   $\epsilon'$ is the number modulo $2$ of disks going in order
   through a point in the unstable manifold of $b'$, $w$ and then a
   point in the unstable manifold of $a'$. There is a single disk
   through $w$ which also intersects both the unstable manifolds of
   $a'$ and $b'$ - the disk $d_{3}$. However, the order in which the
   three types of points appear on the boundary of this disk implies
   that precisely one of $\epsilon$ and $\epsilon'$ is non-zero. Which
   one of the two is non-zero is, obviously, a matter of convention
   and we will take here $\epsilon\not=0$. Notice that we also have
   $a\ast a=\delta wt$ with $\delta\in\{0,1\}$.  To estimate this
   product we need to use a second Morse function on
   $\mathbb{T}^2_{\textnormal{clif}}$,
   $g:\mathbb{T}^2_{\textnormal{clif}}to \R$. We will take this
   function to be perfect also and in such a way that the critical
   points of index one - denoted by $a''$ and $b''$ - have unstable
   and stable manifolds that are ``parallel'' copies of the respective
   stable and unstable manifolds of $f$. See
   figure~\ref{f:clif-traj-2}.
   \begin{figure}[htbp]
      \begin{center}
         \epsfig{file=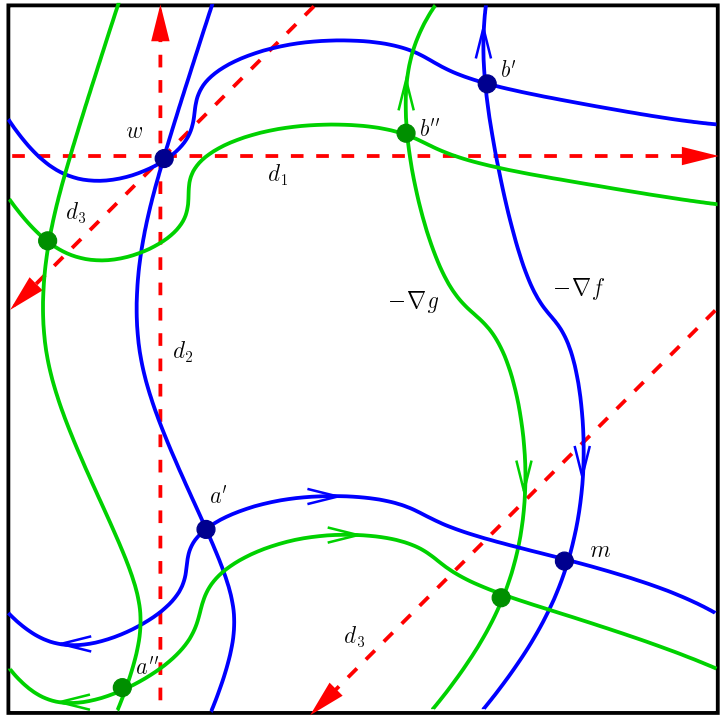, width=0.6\linewidth}
      \end{center}
      \caption{Trajectories of $-\nabla f$, $-\nabla g$ and
        holomorphic disks on $\mathbb{T}^2_{\textnormal{clif}}$.}
      \label{f:clif-traj-2}
   \end{figure}
   Now, there are precisely two pseudo-holomorphic disks that go
   through $w$ as well as both the unstable manifolds of $a'$ and
   $a''$: the disks $d_{2}$ and $d_{3}$. It is at this point that we
   use the fact that $[d_{2}]=b$, $[d_{3}]=-a-b$. Indeed, this means
   that the order in which these two disks pass through these three
   points is opposite.  Thus, exactly one of these disks will
   contribute to $\delta$ and so $\delta=1$. A similar argument shows
   $b\ast b= wt$. The formula for $m*m$ follows now from the
   associativity of the product. Indeed:
   $$m*m = (a*b+wt)*(b*a) = a*(b*b)*a + b*a t = mt+wt^2.$$
   (Recall
   that we are working over $\mathbb{Z}_2$.)

   We now determine what is the map
   $$\tilde{\phi}:QH_{\ast}(\mathbb{T}^2_{\textnormal{clif}})\to
   QH_{\ast}(\mathbb{T}^2_{\textnormal{clif}})$$
   which is induced by
   $\bar{\phi}$. For degree reasons we have $\tilde{\phi}(w)=w$,
   $\tilde{\phi}(a)=b$, $\tilde{\phi}(b)=a$ and, by Corollary
   \ref{cor:symm}, we know that $\tilde{\phi}$ is a morphism of
   algebras over $QH(M)$ (from this also follows immediately that
   $\tilde{\phi}(m)=m+wt$).

   Finally, we compute $h\ast a$ and $h\ast b$. We have, $h\ast
   a=h\ast \tilde{\phi}(b)=\tilde{\phi}(h\ast b)$. Now $h\ast a=
   (u_{1}a + u_{2}b)t$ with $u_{1},u_{2}\in \mathbb{Z}_2$ which
   implies that $h\ast b= (u_{1}b + u_{2} a)t$. As in Corollary
   \ref{T:Lag-CPn} we also have that $h \ast (-):
   H_{1}(\mathbb{T}^2_{\textnormal{clif}};\mathbb{Z}_2)\to
   H_{1}(\mathbb{T}^2_{\textnormal{clif}};\mathbb{Z}_2)t$ is an
   isomorphism.  This implies that precisely one of $u_{1},u_{2}$ is
   non zero.  Assume that $u_{1}=0$ and $u_{2}=1$.  Then $h\ast a=
   bt$, $h\ast(h\ast a)=at^{2}$ and $h\ast(h\ast(h\ast a))=bt^{3}$
   which is not possible because $h^{\ast 3}=[\C P^{2}] t^{3}$ (where,
   $[\C P^{2}]$ denotes the fundamental class of $\C P^{2}$) and $[\C
   P^{2}]\ast a=a$. Thus we are left with $u_{1}=1$, $u_{0}=0$ as
   claimed. It is now easy to estimate $h\ast w$. Indeed, $h\ast w
   t=h\ast( a\ast a)= (h\ast a)\ast a=(a\ast a)t=wt^{2}$.  Similarly
   $h\ast m=h\ast (b\ast a)=(h\ast b)\ast a= mt$. Finally, point iii
   is an immediate consequence of the first two points and of
   formula~\eqref{eq:inclusion_mod} from Theorem~\ref{thm:alg_main} in
   \S \ref{sec:algstr}.

   We turn to the proof for the positive quantum homology
   $Q^{+}H(\mathbb{T}^2_{\textnormal{clif}})$. Recall that there is a
   canonical homomorphism $Q^{+}H(L) \to QH(L)$ induced from the
   inclusion $\mathcal{C}^{+}(f,J) \subset \mathcal{C}(f,J)$.  But in
   our case ($L = \mathbb{T}^2_{\textnormal{clif}}$) this homomorphism
   is actually an injection.  This is because for the choices of $f$
   and $J$ made above the differential $d$ vanishes.  As all the
   quantum operations (the product as well as the external product)
   involve only positive powers of $t$ the statement on
   $Q^{+}H(\mathbb{T}^2_{\textnormal{clif}})$ follows.
\end{proof}

\begin{rem}
   \begin{enumerate}
     \item It is easy to see that the augmentation in the case of the
      Clifford torus in $\C P^{2}$ is the obvious one.  However, the
      duality map $\tilde\eta: Q^{+}H_{\ast}(L)\to Q^{+}H^{\ast}(L)$
      (as defined in \ref{susubsubsec:duality}) verifies
      $\tilde\eta(m)=w^{\ast}+m^{\ast}t$, $\tilde\eta(a)=b^{\ast}$,
      $\tilde\eta(b)=a^{\ast}$ and $\tilde\eta(w)=m^{\ast}$.  Here, by
      a slight abuse in notation, $Q^{+}H^{\ast}(L)=
      H_{\ast}((\mathcal{C}^{+}(L))^{\ast})$.
     \item A careful inspection of the proof of
      Proposition~\ref{prop:clifford} above shows that essentially the
      same argument can be applied for any monotone $2$-dimensional
      Lagrangian torus $L$, as long as we have enough information on
      the holomorphic disks of Maslov index $2$ passing through a
      generic point in $L$. This will be further explored in
      \S\ref{Sb:explicit-quant-T2} below.
   \end{enumerate}
\end{rem}

\subsubsection{Lagrangians that look like
  $\mathbb{R}P^n \subset {\mathbb{C}}P^n$} \label{Sb:rpn-cpn}

The real projective space $\mathbb{R}P^n$ viewed as a submanifold of
${\mathbb{C}}P^n$ is a monotone Lagrangian minimal Maslov number
$=n+1$. Note that when $n \geq 2$,
$H_1(\mathbb{R}P^n;\mathbb{Z})=\mathbb{Z}_2$ hence
$2H_1(\mathbb{R}P^n;\mathbb{Z})=0$. It turns out that the latter
condition on its own imposes very strong restrictions on Lagrangians
in ${\mathbb{C}}P^n$. As we will see below Lagrangians $L$ with
$2H_1(L;\mathbb{Z})=0$ have strong topological similarities to
$\mathbb{R}P^n$. Moreover, as we will see in \S\ref{Sb:quant-rpn-cpn}
their quantum structures are are almost determined by that condition.
We start with topological restrictions.
\begin{prop} \label{T:L-rpn}
   Let $L \subset {\mathbb{C}}P^n$ be a Lagrangian submanifold with
   $2H_1(L;\mathbb{Z})=0$. Then:
   \begin{enumerate}
     \item $N_L = n+1$. \label{I:L-rpn-NL}
     \item $H_i(L;\mathbb{Z}_2) \cong \mathbb{Z}_2$ for every $0 \leq
      i \leq n$, and $H_1(L;\mathbb{Z}) \cong \mathbb{Z}_2$.
      \label{I:L-rpn-Hi}
     \item There exists a canonical isomorphisms of graded vector
      spaces $QH_*(L) \cong (H(L;\mathbb{Z}_2) \otimes \Lambda)_*$.
      Hence by~\ref{I:L-rpn-NL},~\ref{I:L-rpn-Hi} we have
      $QH_j(L)\cong \mathbb{Z}_2$ for every $j \in \mathbb{Z}$.
      \label{I:L-rpn-QH}
     \item Let $\alpha_i \in H_i(L;\mathbb{Z}_2)$ be the generator.
      Then $\alpha_{n-2} \cap (-):H_i(L:\mathbb{Z}_2) \to
      H_{i-2}(L;\mathbb{Z}_2)$ is an isomorphism for every $2 \leq i
      \leq n$. Thus we have $\alpha_i \cap \alpha_{n-2} =
      \alpha_{i-2}$ for every $2 \leq i \leq n$.  Moreover
      $\alpha_{n-2}=h \cap_L [L]$, where $[L]=\alpha_n \in
      H_n(L;\mathbb{Z}_2)$ is the fundamental class and $\cap_L$
      stands for the exterior cap product between elements of
      $H_*({\mathbb{C}}P^n;\mathbb{Z}_2)$ and $H_*(L;\mathbb{Z}_2)$.
      \label{I:L-rpn-cap}
     \item When $n=$\,even, $\alpha_{n-1} \cap (-):
      H_i(L;\mathbb{Z}_2) \to H_{i-1}(L;\mathbb{Z}_2)$ is an
      isomorphism for every $1 \leq i \leq n$. In particular
      $H_*(L;\mathbb{Z}_2)$ is generated by $\alpha_{n-1} \in
      H_{n-1}(L;\mathbb{Z}_2)$. \label{I:L-rpn-even}
     \item Let $\textnormal{inc}_*:H_i(L;\mathbb{Z}_2) \to
      H_i({\mathbb{C}}P^n;\mathbb{Z}_2)$ be the homomorphism induced
      by the inclusion $L \subset {\mathbb{C}}P^n$. Then
      $\textnormal{inc}_*$ is an isomorphism for every $0 \leq
      i=$\,even $\leq n$.
      \label{I:L-rpn-incl}
   \end{enumerate}
\end{prop}

\begin{rem}
   \begin{enumerate}
     \item Proposition ~\ref{T:L-rpn} has already been established in
      the past. Statements~\ref{I:L-rpn-Hi},~\ref{I:L-rpn-QH} have
      been proved by Seidel~\cite{Se:graded} using the theory of
      graded Lagrangian submanifolds. An alternative approach which
      also proves statements~\ref{I:L-rpn-cap},~\ref{I:L-rpn-even} has
      been given by Biran~\cite{Bi:Nonintersections}. Below we give a
      different proof based on our theory.
     \item Other than $\mathbb{R}P^n$ we are not aware of other
      Lagrangian submanifolds $L\subset {\mathbb{C}}P^n$ with
      $2H_1(L;\mathbb{Z})=0$. Note however that in $\mathbb{C}P^3$
      there exists a Lagrangian submanifold $L^3$, not diffeomorphic
      to $\mathbb{R}P^3$, with $H_i(L;\mathbb{Z}_2)=\mathbb{Z}_2$ for
      every $i$. This Lagrangian is the quotient of $\mathbb{R}P^3$ by
      the dihedral group $D_3$. It has $H_1(L;\mathbb{Z})\cong
      \mathbb{Z}_4$. This example is due to Chiang~\cite{Chiang:RP3}.
   \end{enumerate}
\end{rem}

\begin{proof}[Proof of Proposition ~\ref{T:L-rpn}]
   Since $2H_1(L;\mathbb{Z})=0$ it is easy to see that $L \subset
   {\mathbb{C}}P^n$ is monotone. Moreover, a simple computation shows
   that the minimal Maslov number of $L$ is $N_L = k(n+1)$ for some $k
   \geq 1$.

   Let $f:L \to \mathbb{R}$ be a Morse function with exactly one local
   minimum $x_0$ and one local maximum $x_n$. Let $CM_*(f) =
   \mathbb{Z}_2 \langle \textnormal{Crit}(f) \rangle$ be graded by
   Morse indices. Denote by $\mathcal{C}_* = (CM(f) \otimes
   \Lambda)_*$ the string of pearls complex. The differential $d:
   \mathcal{C}_* \to \mathcal{C}_{*-1}$ can be written as $d = \sum_{j
     \geq 0} \partial_j \otimes t^j$, where
   \begin{equation} \label{Eq:del-j}
      \partial_j: CM_*(f) \longrightarrow CM_{*-1+jN_L}(f)
   \end{equation}
   counts trajectories of pearls with total Maslov number $= jN_L$.
   Note that since $L$ is monotone $\partial_0$ is just the
   Morse-homology differential.

   We now prove statement~\ref{I:L-rpn-NL}. If $k \geq 2$ then $N_L =
   k(n+1) > n+1$ hence by formula~\eqref{Eq:del-j} we have
   $\partial_j=0$ for every $j \geq 1$ and we obtain:
   \begin{equation*}
      QH_i(L) =
      \begin{cases}
         H_i(L;\mathbb{Z}_2), & 0 \leq i \leq n \\
         0, & n+1 \leq i \leq N_L-1
      \end{cases}
   \end{equation*}
   But this contradicts the $2$-periodicity asserted by
   Corollary~\ref{T:Lag-CPn}. Thus $k=1$ and $N_L=n+1$. This proves
   statement~\ref{Eq:del-j}. Note that this also implies that
   $H_1(L;\mathbb{Z}) \neq 0$ (for otherwise $N_L=2(n+1)$).  Since
   $2H_1(L;\mathbb{Z}) = 0$ we have $H_1(L;\mathbb{Z}_2) =
   H_1(L;\mathbb{Z}) \otimes \mathbb{Z}_2 \neq 0$. We will use this
   below.

   We prove statements~\ref{I:L-rpn-Hi},~\ref{I:L-rpn-QH}.  Consider
   the operator $\partial_1: C_*(f) \to C_{*+n}(f)$. Clearly
   $\partial_1(x)=0$ for every $x \in \textnormal{Crit}(f)$ with $|x|
   \neq 0$. Thus $d=\partial_0$ on $\mathcal{C}_j$ for every $j$ that
   satisfies $j \not \equiv 0(\bmod{n+1})$ and $j \not \equiv n
   (\bmod{n+1})$.  It follows that
   \begin{equation} \label{Eq:QHL-il}
      QH_{i+l(n+1)}(L) \cong H_i(L) t^{-l}, \quad
      \forall \, 0< i< n, l \in \mathbb{Z}.
   \end{equation}
   Next, consider the value of $\partial_1(x_0)$. There are two
   possibilities:
   \begin{enumerate}[(i)]
     \item $\partial_1(x_0)=x_n$.
     \item $\partial_1(x_0)=0$.
   \end{enumerate}
   We claim that possibility~i is impossible. Indeed by standard Morse
   theory $\partial_0(CM_1(f))=0$ and $\partial_0(CM_n(f))=0$,
   therefore if $\partial_1(x_0)=x_n$ then
   $$d:\mathcal{C}_0=\mathbb{Z}_2 x_0 \longrightarrow \mathcal{C}_{-1}
   = \mathbb{Z}_2 x_n t$$
   is an isomorphism hence $QH_0(L)=0$,
   $QH_{-1}(L)=0$. By Corollary~\ref{T:Lag-CPn} we obtain $QH_j(L)=0$
   for every $j \in \mathbb{Z}$. On the other hand
   by~\eqref{Eq:QHL-il} $QH_1(L)\cong H_1(L;\mathbb{Z}_2)$ and we have
   just seen that $H_1(L;\mathbb{Z}_2) \neq 0$. A contradiction. This
   proves that $\partial_1(x_0)=0$. It follows that $d=\partial_0$
   hence $QH_*(L) \cong (H(L;\mathbb{Z}_2) \otimes \Lambda)_*$. In
   particular $QH_0(L) \cong \mathbb{Z}_2$ and $QH_{-1}(L) \cong
   QH_n(L) \cong \mathbb{Z}_2$. By Corollary~\ref{T:Lag-CPn} $QH_j(L)
   \cong \mathbb{Z}_2$ for every $j \in \mathbb{Z}$. We also conclude
   that $H_i(L;\mathbb{Z}_2) \cong QH_i(L) \cong \mathbb{Z}_2$ for
   every $0 \leq i \leq n$. Finally note that since
   $2H_1(L;\mathbb{Z}) = 0$ we have $H_1(L;\mathbb{Z}) \cong
   \mathbb{Z}_2^{\oplus r}$ for some $r \geq 0$, hence
   $H_1(L;\mathbb{Z}_2) = \mathbb{Z}_2^{\oplus r}$.  But we have seen
   that $H_1(L;\mathbb{Z}_2) \cong \mathbb{Z}_2$ hence $r=1$ and
   $H_1(L;\mathbb{Z}) \cong \mathbb{Z}_2$. This completes the proof of
   statements~\ref{I:L-rpn-Hi},~\ref{I:L-rpn-QH}.

   To prove statement~\ref{I:L-rpn-cap} recall formula~\eqref{Eq:qm-*}
   of \S\ref{S:qm} by which the quantum module operation $QH_l(M)
   \otimes QH_j(L) \to QH_{l+j-2n}(L)$ is defined. As $N_L = n+1$ and
   the degree of the hyperplane class $h$ is $2n-2$, it follows from
   that formula that $h * \alpha = h \cap_L \alpha$ for every $\alpha
   \in QH_i(L) \cong H_i(L;\mathbb{Z}_2)$ for $2 \leq i \leq n$.
   Denote by $\alpha_i \in H_i(L;\mathbb{Z}_2)$ the generator. By
   Corollary~\ref{T:Lag-CPn} it follows that $h*\alpha_n = h \cap_L
   \alpha_n \neq 0$ hence $h * \alpha_n = \alpha_{n-2}$. Next, recall
   that $\alpha_n \in QH_n(L)$ is the unity hence
   $$\alpha_{n-2}* \alpha_i= (h*\alpha_n)*\alpha_i = h*(\alpha_n *
   \alpha_i) = h * \alpha_i=\alpha_{i-2}, \quad \forall\, 2 \leq i
   \leq n.$$
   But by formula~\eqref{eq:prod} of \S\ref{subsec:product}
   $\alpha_{n-2} * \alpha_i = \alpha_{n-2} \cap \alpha_i$ for every $2
   \leq i \leq n$. This proves statement~\ref{I:L-rpn-cap}.

   To prove statement~\ref{I:L-rpn-even} note that by
   Corollary~\ref{T:Lag-CPn}
   $$h^{*(\frac{n}{2}+1)} * \alpha_n = \alpha_{n-1} t.$$
   Therefore
   $$\alpha_{n-1}*\alpha_i = (h^{*(\frac{n}{2}+1)} * \alpha_n) *
   \alpha_i t^{-1} = (h^{*(\frac{n}{2}+1)}*\alpha_i) t^{-1} =
   \alpha_{i-1}, \quad \forall\, 1\leq i \leq n,$$
   where the last
   equality also follows from Corollary~\ref{T:Lag-CPn}. But by
   looking at the Morse indices in formula~\eqref{eq:prod} of
   \S\ref{subsec:product} we conclude again that $\alpha_{n-1} *
   \alpha_i = \alpha_{n-1} \cap \alpha_i$ for every $1\leq i\leq n$.
   Thus $\alpha_{n-1} \cap \alpha_i = \alpha_{i-1}$.

   Finally, statement~\ref{I:L-rpn-incl} will follow immediately from
   Proposition~\ref{T:quant-incl-rpn-cpn} below. We therefore postpone
   the proof.
\end{proof}

\subsubsection{Quantum structures} \label{Sb:quant-rpn-cpn}
Let $L \subset {\mathbb{C}}P^n$ be a Lagrangian with
$2H_1(L;\mathbb{Z})=0$. In view of Proposition~\ref{T:L-rpn} denote by
$\alpha_i \in QH_i(L)$ the generator for every $i \in \mathbb{Z}$.
According to this notation we have $\alpha_{i+l(n+1)} = \alpha_i
t^{-l}$ for every $i, l \in \mathbb{Z}$ and by
Proposition~\ref{T:L-rpn}, $h*\alpha_i = \alpha_{i-2}$ for every $i
\in \mathbb{Z}$.

\begin{prop} \label{T:quant-rpn-cpn}
   Let $L \subset {\mathbb{C}}P^n$ be as above. Let $k,j \in
   \mathbb{Z}$.  If one of $k,j$ is odd then:
   \begin{equation} \label{Eq:alpha_k-alpha-j}
      \alpha_k * \alpha_j = \alpha_{j+k-n}.
   \end{equation}
   The same formula holds for every $k,j \in \mathbb{Z}$ (regardless
   of their parity) in each of the following cases:
   \begin{enumerate}
     \item When $n=$\,even.
     \item When $L$ is diffeomorphic to $\mathbb{R}P^n$.
     \item More generally, when $\alpha_{n-1} \cap \alpha_{n-1} \neq
      0$. (c.f. statement~\ref{I:L-rpn-even} of
      Theorem~\ref{T:L-rpn}.)
   \end{enumerate}
\end{prop}

\begin{proof}
   \textbf{Assume $k=$\,odd.}  By Corollary~\ref{T:Lag-CPn}:
   \begin{align*}
      \alpha_k * \alpha_j & = (h^{*(-\frac{k+1}{2})} *
      \alpha_{-1})*\alpha_j = (h^{*(-\frac{k+1}{2})} *
      \alpha_{n}T)*\alpha_j = (h^{*(-\frac{k+1}{2})} * \alpha_j) t \\
      & = \alpha_{j+k+1}t = \alpha_{j+k-n}.
   \end{align*}


   \smallskip \noindent \textbf{Assume $j=$\,odd.} The proof is
   similar to the case $k=$\,odd since for every $a \in
   QH({\mathbb{C}}P^n)$ we have $a*(\alpha_k*\alpha_j) =
   \alpha_k*(a*\alpha_j)$.

   \smallskip \noindent \textbf{Assume $n=$\,even.} We may assume that
   $k=$\,even. Then we have $\alpha_k = h^{*\frac{n-k}{2}}* \alpha_n$
   hence:
   $$\alpha_k * \alpha_j = (h^{*\frac{n-k}{2}}* \alpha_n)*\alpha_j =
   h^{*\frac{n-k}{2}}* \alpha_j = \alpha_{j+k-n}.$$

   \smallskip \noindent \textbf{Assume $n=$\,odd and $\alpha_{n-1}
     \cap \alpha_{n-1} \neq 0$.} In view of the above we may assume
   that $k,j$ are both even.  By the definition of the quantum product
   we have $\alpha_{n-1} * \alpha_{n-1} = \alpha_{n-1} \cap
   \alpha_{n-1} = \alpha_{n-2}$, where the last equality follows from
   the fact that $H_{n-2}(L;\mathbb{Z}_2) = \mathbb{Z}_2\alpha_{n-2}$.
   Since $k+1$ and $j+1$ are both odd then by what we have proved
   above:
   $$\alpha_k * \alpha_j =
   (\alpha_{k+1}*\alpha_{n-1})*(\alpha_{n-1}*\alpha_{j+1}) =
   \alpha_{k+1}* \alpha_{n-2} * \alpha_{j+1} =
   \alpha_{k-1}*\alpha_{j+1} = \alpha_{k+j-n}.$$
\end{proof}

The next result describes the quantum inclusion map $i_L:QH_*(L) \to
QH_*({\mathbb{C}}P^n)$. Denote by $a_j \in
H_j({\mathbb{C}}P^n;\mathbb{Z}_2)$ the generator, $0 \leq j \leq 2n$.
Thus
\begin{equation*}
   a_j =
   \begin{cases}
      0, & j=\textnormal{odd} \\
      h^{\cap(n-\frac{j}{2})}, & j=\textnormal{even}
   \end{cases}
\end{equation*}
\begin{prop} \label{T:quant-incl-rpn-cpn}
   Let $L \subset {\mathbb{C}}P^n$ be as above.
   \begin{enumerate}
     \item If $n=$\,even then:
      \begin{align*}
         & i_L(\alpha_{2k})=a_{2k}, \quad  \forall\, 0 \leq 2k \leq n, \\
         & i_L(\alpha_{2k+1}) = a_{2k+n+2} t,\quad \forall\, 1 \leq
         2k+1 \leq n-1.
      \end{align*}
     \item If $n=$\,odd then:
      \begin{align*}
         & i_L(\alpha_{2k}) = a_{2k} +
         a_{2k+n+1} t, \quad \forall\, 0 \leq 2k \leq n, \\
         & i_L(\alpha_{2k+1})=0, \quad \forall\, k.
      \end{align*}
   \end{enumerate}
\end{prop}
\begin{proof}
   By our notation $\alpha_0 = [\textnormal{point}] \in
   H_0(L;\mathbb{Z}_2) \cong QH_0(L)$, $a_0 = [\textnormal{point}] \in
   H_0({\mathbb{C}}P^n;\mathbb{Z}_2)$. Recall also that $h^{*n}=a_0$
   and $h^{*(n+1)} = u t^2 $.

   From the definition of $i_L$ (see \S\ref{Sb:q-inclusion}) it
   follows by a simple computation that:
   \begin{equation} \label{Eq:iL-alpha-0}
      i_L(\alpha_0) = a_0 + b t, \quad
      \textnormal{for some } b \in H_{n+1}({\mathbb{C}}P^n;\mathbb{Z}_2).
   \end{equation}
   We claim that:
   \begin{equation} \label{Eq:iL-alpha-2k}
      i_L(\alpha_{2k}) = a_{2k} + h^{*(-k)}*b t, \quad \forall\,
      0 \leq 2k \leq n.
   \end{equation}
   Indeed, by~\eqref{Eq:iL-alpha-0} we have:
   \begin{align*}
      i_L(\alpha_{2k}) &= i_L(h^{*(-k)} * \alpha_0) =
      h^{*(-k)}*i_L(\alpha_0) = h^{*(-k)}*(a_0 + b t) \\
      &= h^{*(-k)}*(h^{*n} + b t) = h^{*(n-k)} + h^{*(-k)}*b t =
      a_{2k} + h^{*(-k)}*b t.
   \end{align*}
   Next we claim that :
   \begin{equation} \label{Eq:iL-alpha-n-1}
      i_L(\alpha_{n-1-2r}) = a_{2n-2r} t + h^{*(r+1)}*b, \quad
      \forall\, 0 \leq 2r \leq n-1.
   \end{equation}
   Indeed by~\eqref{Eq:iL-alpha-0} we have:
   \begin{align*}
      i_L(\alpha_{n-1-2r}) & = i_L(h^{*(r+1)}*\alpha_{n+1}) =
      i_L(h^{*(r+1)}*\alpha_0 t^{-1}) = h^{*(r+1)}*(a_0 + b t)t^{-1} \\
      & = (h^{*(r+1)}*h^{*n} + h^{*(r+1)}*b t)t^{-1} = (h^{*r}*
      u t^2 + h^{*(r+1)}*b t)t^{-1} \\
      & = a_{2n-2r} t + h^{*(r+1)}*b.
   \end{align*}

   Suppose that $n=$\,even. Then $b=0$ since $\deg b = n+1=$\,odd.
   Statement~1 of the theorem follows immediately
   from~\eqref{Eq:iL-alpha-2k} and~\eqref{Eq:iL-alpha-n-1}.

   Assume now that $n=$\,odd. Comparing~\eqref{Eq:iL-alpha-2k}
   to~\eqref{Eq:iL-alpha-n-1} with $2k=n-1$, $r=0$, we obtain that
   $b=a_{n+1}$. The statement about $i_L(\alpha_{2k})$ follows now
   from~\eqref{Eq:iL-alpha-2k}. As for $i_L(\alpha_{2k+1})$, it is $0$
   since for $n=$\,odd, $QH_{2k+1}({\mathbb{C}}P^n)=0$ for every $k$.

   Finally note that statement~\ref{I:L-rpn-incl} of
   Proposition~\ref{T:L-rpn} follows from the above since in our case
   we have $\textnormal{inc}_* = {i_L}_{|_{t=0}}$.
\end{proof}

\subsubsection{Existence of holomorphic disks satisfying constrains}
\label{Sb:enum-rpn}

Denote by $\mathcal{J}=\mathcal{J}({\mathbb{C}}P^n,
\omega_{\textnormal{FS}})$ the space of
$\omega_{\textnormal{FS}}$-compatible almost complex structures on
${\mathbb{C}}P^n$.
\begin{prop} \label{P:enum-rpn-2}
   Let $L \subset {\mathbb{C}}P^n$ be a Lagrangian submanifold with
   $2H_1(L;\mathbb{Z})=0$. Assume that one of the following conditions
   is satisfied:
   \begin{enumerate}
     \item $n=$\,even.
     \item $L$ is diffeomorphic to $\mathbb{R}P^n$.
     \item More generally, $\alpha_{n-1} \cap \alpha_{n-1} \neq 0$.
   \end{enumerate}
   Let $x',x'' \in L$ two distinct points. Then for every $J \in
   \mathcal{J}$ there exists a $J$-holomorphic disk $u:(D, \partial D)
   \to ({\mathbb{C}}P^n, L)$ with $\mu([u])=n+1$ and $u(\partial D)
   \ni x', x''$. For a generic choice of $J \in \mathcal{J}$ the
   number of such disks with $u(-1)=x'$, $u(1)=x''$, up to
   parametrizations fixing $-1,1 \in D$, is $\geq 2$ and even.
\end{prop}
\begin{rem} \label{R:enum-rpn-2}
   We will see in the proof that if $C^{n-1} \subset L$ is an
   $(n-1)$-dimensional $\mathbb{Z}_2$-cycle with
   $[C^{n-1}]=\alpha_{n-1}$ and such that $x', x'' \notin C^{n-1}$,
   then for generic $J \in \mathcal{J}$ there exist two
   $J$-holomorphic disks $u_1, u_2:(D, \partial D) \to
   ({\mathbb{C}}P^n, L)$ with the following properties:
   \begin{enumerate}[(i)]
     \item $\mu([u_1])=\mu([u_2])=n+1$.
     \item $u_1(-1)=x', \quad u_1(1)=x'', \quad u_1(i) \in C^{n-1}$.
     \item $u_2(-1)=x', \quad u_2(1)=x'', \quad u_2(-i) \in C^{n-1}$.
     \item $u_2$ is not obtained from $u_1$ by a reparametrization
      fixing $-1,1 \in D$.
   \end{enumerate}
   Moreover, the number of disks $u_1$ (resp. $u_2$) as above is odd.
\end{rem}

\begin{proof}[Proof of Proposition~\ref{P:enum-rpn-2}]
   By Proposition~\ref{T:quant-rpn-cpn} we have:
   \begin{equation} \label{Eq:prod-0-n-1}
      \alpha_0 * \alpha_{n-1} = \alpha_{-1} = \alpha_n t.
   \end{equation}
   Let $f_1, f_2, f_3:L \to \mathbb{R}$ be a generic triple of Morse
   functions with the same critical points (and the same indices at
   each critical point) and such that the $f_i$'s have exactly one
   local minimum $x'$ and one local maximum $x''$. Let $\rho_1,
   \rho_2, \rho_3$ be a generic triple of Riemannian metrics on $L$.
   Choose a generic $J \in \mathcal{J}$.

   Denote by $CM_*(f_i)$ the Morse complex of $f_i$ (with respect to
   $\rho_i$).  Choose a cycle $y \in CM_{n-1}(f_3)$ that represents
   $\alpha_{n-1} \in H_{n-1}(L;\mathbb{Z}_2)$. Recall from the proof
   of Theorem~\ref{T:L-rpn} that the Floer differential coincides with
   the Morse differential thus $\alpha_0=[x']$, $\alpha_n=[x'']$,
   $\alpha_{n-1}=[y]$ in Floer homology. It follows
   from~\eqref{Eq:prod-0-n-1} that $x'*y = x'' t$. (Note that due to
   dimension we cannot have additional boundary terms on the
   right-hand side.) Therefore, in the notation of~\eqref{eq:prod}
   (see \S\ref{subsec:product}) there exists a class $A \in
   H_2({\mathbb{C}}P^n,L;\mathbb{Z})$ with $\mu(A)=n+1$ and a critical
   point $y_0$ participating in $y \in CM_{n-1}(f_3)$ such that
   $\mathcal{P}(x',y_0,x'';A,J) \neq \emptyset$. Let $(l_1, l_2, l_3,
   u) \in \mathcal{P}(x',y_0,x'';A,J)$. As $x'$ is a minimum and $x''$
   a maximum their unstable and stable manifolds are $W_{x'}^u =
   \{x'\}$, $W_{x''}^s=\{x''\}$ respectively. Moreover since
   $\mu(A)=n+1$ which is the minimal Maslov number it follows that the
   only possible configuration that $(l_1, l_2, l_3, u)$ can take is
   the following (see figure~\ref{f:disk-min-max}):
   \begin{itemize}
     \item $l_1$ is the constant trajectory at $x'$.
     \item $l_2$ is a (negative) gradient trajectory (without
      $J$-holomorphic disks) emanating from $y_0$.
     \item $l_3$ is the constant trajectory at $x''$.
     \item The $J$-holomorphic disk $u$ is not constant, hence it has
      $\mu([u])=n+1$ and $u(e^{2\pi i /3})=x'$, $u(1) \in l_2$,
      $u(e^{4 \pi i /3})=x''$.
   \end{itemize}
   We have proved that for generic $J \in \mathcal{J}$ there exists a
   $J$-holomorphic disk $u$ with $\mu([u])=n+1$ and $u(\partial D)
   \ni x'', x'$.  As $N_L=n+1$ it follows from Gromov compactness
   theorem that there exists such a disk for every $J \in
   \mathcal{J}$.
   \begin{figure}[htbp]
      \begin{center}
         \epsfig{file=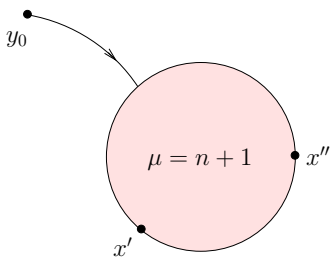, width=0.3\linewidth}
      \end{center}
      \caption{A $J$-holomorphic disk from $\mathcal{P}(x',y_0,x'';A,J)$.}
      \label{f:disk-min-max}
   \end{figure}

   Note that after a suitable reparametrization we may assume that
   $u(-1)=x'$, $u(1)=x''$. It remains to prove that for generic $J$
   the number of such disks is even. For this purpose we use the
   notation from the proof of Proposition~\ref{T:L-rpn}. Recall the
   operator $\partial_1:CM_*(f) \to CM_{*+n}(f)$. Since $x'$ is a
   minimum and $x''$ a maximum, $\partial_1(x')$ counts the number
   (mod $2$) of $J$-holomorphic disks $u$, up to parametrization, with
   $\mu([u])=n+1$ and $u(\partial D) \ni x'', x'$. As we have seen
   in the proof of Proposition~\ref{T:L-rpn}, $\partial_1(x')=0$,
   hence this number is even.

   Finally, we prove the statement of Remark~\ref{R:enum-rpn-2}.
   First we claim that for every $J$-holomorphic disk $u$ with
   $\mu[u]=n+1$ we have $[u(\partial D)] = \alpha_1 \in
   H_1(L;\mathbb{Z}_2)$.  To see this note that $[u(\partial D)] \neq
   0 \in H_1(L;\mathbb{Z})$ for otherwise $\mu([u])=2(n+1)$. But by
   Proposition~\ref{T:L-rpn} $H_1(L;\mathbb{Z}) \cong \mathbb{Z}_2$.
   Thus $[u(\partial D)] \neq 0$ also in $H_1(L;\mathbb{Z}_2)$.

   Consider the space $\mathcal{M}(x',x'')$ of $J$-holomorphic disks
   $u: (D, \partial D) \to ({\mathbb{C}}P^n, L)$ with $\mu([u])=n+1$
   and $u(-1)=x'$, $u(1) = x''$. The group $G_{-1,1} \subset
   \textnormal{Aut}(D)$ acts on this space. Denote by
   $\mathcal{P}(x',x'')$ its quotient. Note that for generic $J$,
   $\mathcal{P}(x',x'')$ is a finite set.

   Denote by $\gamma_1, \gamma_2 \subset \partial D$ the arcs
   $\gamma_1 = \{ e^{\pi i t} \}_{0 \leq t \leq \pi}$ and $\gamma_2 =
   \{ e^{\pi i t} \}_{\pi \leq t \leq 2\pi}$.  Let $y \in
   CM_{n-1}(f_3)$ be a cycle representing $\alpha_{n-1}$.  The union
   of the unstable submanifolds corresponding to the critical points
   in $y$ is a pseudo cycle homologous to $C^{n-1}$.  By choosing $J
   \in \mathcal{J}$ generic we may assume that $C^{n-1}$ is in general
   position with respect to the evaluation map $\mathcal{M}(x',x'')
   \to L$ evaluating $\tilde{v} \in \mathcal{M}(x',x'')$ at a marked
   point $p_0 \in \partial D \setminus \{ \pm 1\}$. For every $v \in
   \mathcal{P}(x',x'')$ put
   $$n_1(v) = \#_{\mathbb{Z}_2} \bigl( \tilde{v}(\gamma_1) \cap
   C^{n-1}\bigr), \quad n_2(v) = \#_{\mathbb{Z}_2} \bigl(
   \tilde{v}(\gamma_2) \cap C^{n-1}\bigr),$$
   where $\tilde{v} \in
   \mathcal{M}(x',x'')$ parametrizes $v$.

   With this notation the coefficient of $x'' t$ in $x'*y$ is $\sum_{v
     \in \mathcal{P}(x',x'')} n_1(v)$.  Similarly, the coefficient of
   $x'' t$ in $y*x'$ is $\sum_{v \in \mathcal{P}(x',x'')} n_2(v)$. But
   $x'*y = y*x' = x'' t$ hence: $$\sum_{v \in \mathcal{P}(x',x'')}
   n_1(v) = 1, \quad \sum_{v \in \mathcal{P}(x',x'')} n_2(v)=1.$$
   Next
   note that for every $v \in \mathcal{P}(x',x'')$
   $$n_1(v)+n_2(v) = \#_{\mathbb{Z}_2} \bigl( v(\partial D) \cap
   C^{n-1} \bigr) = \alpha_1 \cap \alpha_{n-1}=1.$$
   Let $u_1 \in
   \mathcal{P}(x',x'')$ with $n_1(u_1)=1$. Then $n_2(u_1)=0$ hence
   there exists $u_2 \neq u_1$ with $n_2(u_2)=1$. After suitable
   reparametrizations of $u_1, u_2$ we obtain two disks with the
   properties claimed in Remark~\ref{R:enum-rpn-2}.
\end{proof}

In case $L$ does not satisfy one of the 3 conditions from
Proposition~\ref{P:enum-rpn-2} we still have the following weaker
statement.  (Note that this is theoretical, since we do not know
Lagrangians in ${\mathbb{C}}P^n$ with $2H_1(L;\mathbb{Z})=0$ other
than $\mathbb{R}P^n$.)
\begin{prop} \label{P:enum-rpn-1}
   Let $L \subset {\mathbb{C}}P^n$ be a Lagrangian submanifold with
   $2H_1(L;\mathbb{Z})=0$. Then for every $x \in L$ and every $J \in
   \mathcal{J}$ there exists a $J$-holomorphic disk $u:(D,\partial D)
   \to ({\mathbb{C}}P^n, L)$ with $\mu([u])=n+1$ and $u(\partial D)
   \ni x$.
\end{prop}
\begin{proof}
   The proof is similar to that of Proposition~\ref{P:enum-rpn-1},
   only that now we use the identity $\alpha_0 * \alpha_1 =
   \alpha_{1-n} = \alpha_2 t$.
\end{proof}

The next result deals with existence of holomorphic disks satisfying
mixed constrains, i.e. some of the marked points are on the boundary
and some in the interior.
\begin{prop} \label{P:enum-mixed-rpn}
   Let $L \subset {\mathbb{C}}P^n$ be a Lagrangian with
   $2H_1(L;\mathbb{Z})=0$.
   \begin{enumerate}
     \item \label{I:enum-rpn-1} For every $p \in {\mathbb{C}}P^n
      \setminus L$ and every $J \in \mathcal{J}$ there exists a simple
      $J$-holomorphic disk $u:(D,\partial D) \to ({\mathbb{C}}P^n, L)$
      with $\mu([u])=n+1$ and $u(0)=p$.
     \item \label{I:enum-rpn-2} Let $x \in L$, $p \in {\mathbb{C}}P^n
      \setminus L$. Then for generic $J \in \mathcal{J}$ there exists
      a simple $J$-holomorphic disk $u:(D,\partial D) \to
      ({\mathbb{C}}P^n, L)$ with $\mu([u])=2n+2$ and $u(0)=p$,
      $u(-1)=x$.
     \item \label{I:enum-rp2} Suppose $n=2$ and $L$ is diffeomorphic
      to $\mathbb{R}P^2$.  Let $x',x'' \in L$ be two distinct points
      and $p \in {\mathbb{C}}P^2 \setminus L$. Then for generic $J$
      there exists a simple $J$-holomorphic disk $u:(D, \partial D)
      \to ({\mathbb{C}}P^2, L)$ with $\mu([u])=6$ and $u(-1)=x'$,
      $u(1)=x''$ and $u(0)=p$. Moreover the number of such disks is
      odd.
   \end{enumerate}
\end{prop}
\begin{proof}
   Statement~\ref{I:enum-rpn-1} follows from
   Proposition~\ref{T:quant-incl-rpn-cpn} since
   \begin{equation*}
      i_L(\alpha_{n-1}) =
      \begin{cases}
         u t, & n=\,\textnormal{even}\\
         a_{n-1}+u t, & n=\,\textnormal{odd},
      \end{cases}
   \end{equation*}
   where $u = a_{2n} \in QH_{2n}({\mathbb{C}}P^n)$ is the fundamental
   class.

   We turn to statement~\ref{I:enum-rpn-2}. We start with the
   following identity:
   \begin{equation} \label{Eq:a0-alpha0}
      a_0*\alpha_0 = h^{*n}*\alpha_0 = \alpha_{-2n}=\alpha_2 t^2.
   \end{equation}
   Let $f: L \to \mathbb{R}$ be a Morse function with exactly one
   local minimum at the point $x$. Let $g: {\mathbb{C}}P^n \to
   \mathbb{R}$ be a Morse function with exactly one local minimum at
   the point $p$. Choose a cycle $y \in C_2(f)$ representing
   $\alpha_2$. Choose a generic $J \in \mathcal{J}$.

   From formula~\eqref{Eq:a0-alpha0} it follows that $p*x =y t^2 +
   \textnormal{boundary terms}$. By the definition of the $*$
   operation (see formula~\eqref{Eq:qm} in \S\ref{Sb:external-op}) it
   follows that there exists a vector of non-zero classes
   $\mathbf{A}=(A_1, \ldots, A_l)$ with $\mu(\mathbf{A})=2(n+1)$ such
   that one of the spaces $\mathcal{P}_I(p,x,y';\mathbf{A},J)$,
   $\mathcal{P}_{I'}(p,x,y';\mathbf{A},J)$ is not empty, where $y'$ is
   a critical point participating in $y$. Note that since $p$ is a
   minimum the only trajectory of $-\textnormal{grad} g$ emanating
   from $p$ is constant. Since $p \notin L$ we must have
   $\mathcal{P}_{I'}(p,x,y';\mathbf{A},J) = \emptyset$, thus
   $\mathcal{P}_I(p,x,y';\mathbf{A},J) \neq \emptyset$.  As
   $N_L=2(n+1)$ we have $l \leq 2$.

   We claim that $l=1$. To see this first note that since $x$ is a
   minimum the trajectories of $-\textnormal{grad} f$ emanating from
   $x$ are constant. Therefore an element
   $\mathcal{P}_I(p,x,y';\mathbf{A},J)$ looks like one of the three
   types in figure~\ref{f:a0-alpha0}. In the first case we have a
   (simple) disk $u_1:(D, \partial D) \to ({\mathbb{C}}P^n, L)$ with
   $\mu([u_1])=n+1$ and such that $u(-1)=x$, $u(0)=p$. A simple
   computation shows that the dimension of this configuration is
   negative hence cannot occur for generic $J$. In the second case we
   have a $J$-holomorphic disk $u_2$ with $u_2(0)=p$, $u_2(1) \in
   W_{y'}^s$ which again is a configuration of dimension $-1$ hence
   cannot occur for generic $J$.  We are thus left with the last
   possibility ($l=1$) in which we have a $J$-holomorphic disk $u$
   with $u(0)=p$ and $u(-1)=x$. Note that the disk $u$ is simple since
   by Proposition~\ref{P:qm-axy} all the $J$-holomorphic disks coming
   from $\mathcal{P}_I(p,x,y';\mathbf{A},J)$ are simple.
   \begin{figure}[htbp]
      \begin{center}
         \epsfig{file=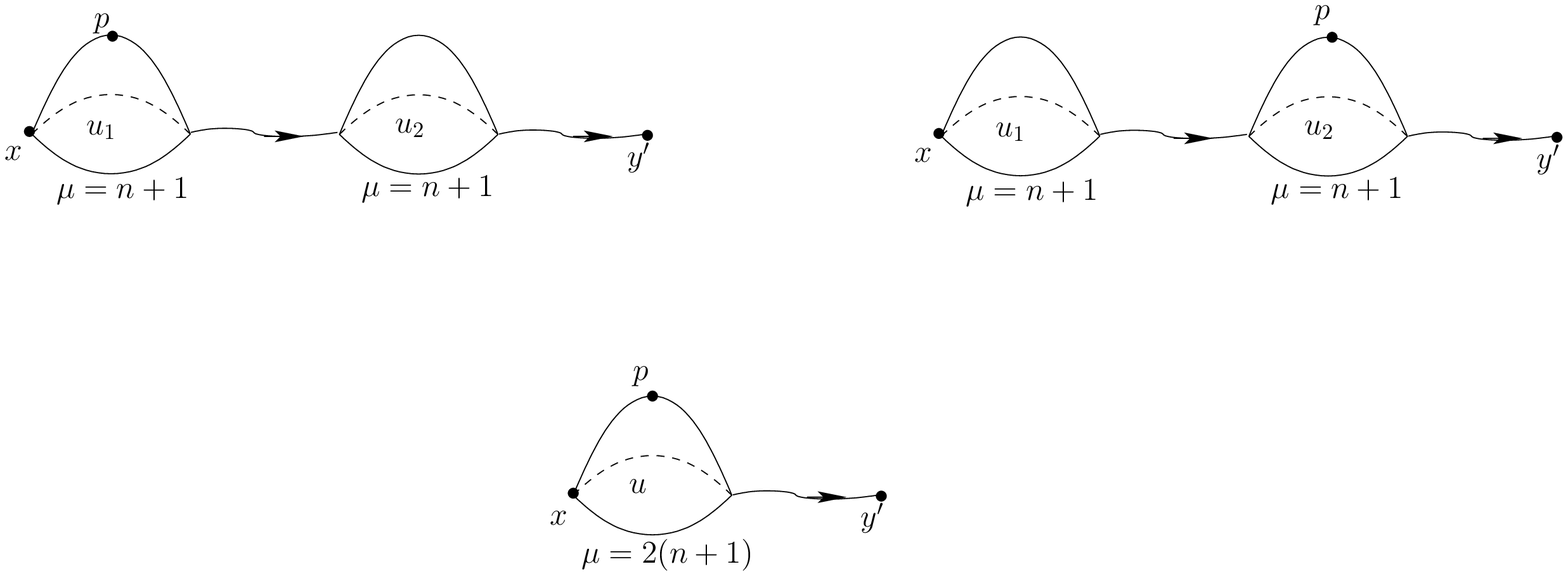, width=0.8\linewidth}
      \end{center}
      \caption{Possible elements of $\mathcal{P}_I(p,x,y';\mathbf{A},J)$.}
      \label{f:a0-alpha0}
   \end{figure}

   To prove statement~\ref{I:enum-rp2} we argue as above but we
   require in addition the function $f$ to have a single local
   maximum. To keep compatibility with the notation we have already
   used denote the point $x'$ by $x$ and $x''$ by $y'$.  Then the
   proof of statement~\ref{I:enum-rpn-2} gives us a simple
   $J$-holomorphic disk with $\mu([u])=2(2+1)=6$ and $u(-1)=x$, $u(1)
   \in W_{y'}^s$ and $u(0)=p$. Note that now $y'$ is the maximum
   (since $n=2$), hence $W_{y'}^s=\{y'\}$ and we obtain $u(1)=y'$.
\end{proof}

\begin{rem} \label{R:rpn-cpn} In case $L$ is the standard embedding
   $\mathbb{R}P^n \subset {\mathbb{C}}P^n$ some of the statements of
   Proposition~\ref{P:enum-rpn-2}-\ref{P:enum-mixed-rpn} can be proved
   by more direct methods.  This is typically done by degenerating to
   the standard almost complex structure $J_{\textnormal{std}}$ (for
   which the statement is obviously true) and passing to generic $J$
   using Gromov compactness theorem.
\end{rem}

\subsubsection{Further questions and remarks} \label{Sb:further-cpn}

In view of the results above one is led to ask the following
questions.  {\sl Do there exist Lagrangian submanifolds $L \subset
  {\mathbb{C}}P^n$ with $2H_1(L;\mathbb{Z})=0$ different than
  $\mathbb{R}P^n$? Are all Lagrangian embeddings of $\mathbb{R}P^n
  \hookrightarrow {\mathbb{C}}P^n$ Hamiltonianly isotopic to the
  standard embedding?}

It is obvious that using the Fukaya $A_{\infty}$-category
theory~\cite{FO3} or alternatively the theory of cluster homology
~\cite{Cor-La:Cluster-1, Cor-La:Cluster-2} one can get more
enumerative results on $J$-holomorphic disks with boundary on
Lagrangians as above.  It would be interesting to test this even in
dimension $n=2$. For example, do these techniques imply that given $5$
points $x_1, \ldots, x_5 \in L$ then for generic $J$ there exists a
$J$-holomorphic disk $u:(D, \partial D) \to ({\mathbb{C}}P^2, L)$ with
$\mu([u])=6$ and $u(\partial D) \ni x_1, \ldots, x_5$?  It would of
course be also interesting to understand the relation of the above to
the works of Welschinger~\cite{Welsch:Enumerative,
  Welsch:Invariants-1, Welsch:Invariants-2} and of
Solomon~\cite{Solomon:intersection}.

\

\subsection{Lagrangian submanifolds of the quadric} \label{Sb:quadric}
Let $Q \subset {\mathbb{C}}P^{n+1}$ be a smooth complex
$n$-dimensional quadric, where $n \geq 2$. More specifically we can
write $Q$ as the zero locus $Q= \{ z \in {\mathbb{C}}P^{n+1} \mid
q(z)=0\}$ of a homogeneous quadratic polynomial $q$ in the variables
$[z_0: \cdots: z_{n+1}] \in {\mathbb{C}}P^{n+1}$, where $q$ defines a
quadratic form of maximal rank. We endow $Q$ with the symplectic
structure induced from ${\mathbb{C}}P^{n+1}$. (Recall that we use the
normalization that the symplectic structure $\omega_{\textnormal{FS}}$
of ${\mathbb{C}}P^{n+1}$ satisfies $\int_{\mathbb{C}P^1}
\omega_{\textnormal{FS}} = \pi$.) When $n \geq 3$ we have by Lefschetz
theorem $H^2(Q;\mathbb{R}) \cong \mathbb{R}$, therefore by Moser
argument all K\"{a}hler forms on $Q$ are symplectically equivalent
up to a constant factor. When $n=2$, $Q \subset {\mathbb{C}}P^3$ is
symplectomorphic to $(\mathbb{C}P^1 \times \mathbb{C}P^1, \omega_{FS}
\oplus \omega_{FS})$. Also note that the symplectic structure on $Q$
(in any dimension) does not depend (up to symplectomorphism) on the
specific choice of the defining polynomial $q$ (this follows from
Moser argument too since the space of smooth quadrics is connected).

\subsubsection{Topology of the quadric} \label{Sb:top-quadric}
The quadric has the following homology:
\begin{equation*}
   H_i(Q;\mathbb{Z}) \cong
   \begin{cases}
      0 & \textnormal{if }i=\textnormal{odd} \\
      \mathbb{Z} & \textnormal{if } i=\textnormal{even} \neq n
   \end{cases}
\end{equation*}
Moreover, when $n=$\,even, $H_n(Q;\mathbb{Z})\cong \mathbb{Z} \oplus
\mathbb{Z}$. To see the generators of $H_n(Q;\mathbb{Z})$, write
$n=2k$. There exist two families $\mathcal{F}, \mathcal{F}'$ of
complex $k$-dimensional planes lying on $Q$
(see~\cite{Gr-Ha:alg-geom}). Let $P \in \mathcal{F}$, $P' \in
\mathcal{F}'$ be two such planes belonging to different families.  Put
$a=[P]$, $b=[P']$. Then $H_n(Q;\mathbb{Z}) = \mathbb{Z}a \oplus
\mathbb{Z} b$ and $h^{\cap k}=a+b$. Moreover, we have:
\begin{equation} \label{Eq:cap-quad}
   \begin{aligned}
      & \textnormal{for } k=\textnormal{odd}: & a\cap b = 1,\; a\cap a
      =
      b\cap b =0, \\
      & \textnormal{for } k=\textnormal{even}: & a\cap b = 0, \; a\cap
      a = b\cap b =1.
   \end{aligned}
\end{equation}

\subsubsection{Quantum homology of the quadric} \label{Sb:qhom-quadric}
Let $h \in H_{2n-2}(Q;\mathbb{Z})$ be the class of a hyperplane
section, $p \in H_0(Q;\mathbb{Z})$ the class of a point and $u \in
H_{2n}(Q;\mathbb{Z})$ the fundamental class. We will first describe
the quantum cohomology over $\mathbb{Z}$. Define
$\Lambda^{\mathbb{Z}}_* = \mathbb{Z}[t, t^{-1}]$ where $\deg t =
-N_L$.  Here $N_L$ is the minimal Maslov number of a Lagrangian
submanifold that will appear later on. Note that that $c_1(Q) =
nPD(h)$, hence $N_L | 2n$. Let
$QH_*(Q;\Lambda^{\mathbb{Z}})=(H(Q;\mathbb{Z}) \otimes
\Lambda^{\mathbb{Z}})_*$ be the quantum homology endowed with the
quantum cap product $*$.

\begin{prop}[See~\cite{Beau:quant}] \label{P:qhom-quadric}
   The quantum product satisfies the following identities:
   \begin{align*}
      & h^{*j}=h^{\cap j} \; \forall \, 0 \leq j \leq n-1, \quad
      h^{*n}=2p+2u  t^{2n/N_L}, \\
      & h^{*(n+1)}=4h t^{2n/N_L}, \quad p*p=u t^{4n/N_L}.
   \end{align*}
   When $n=$\,even we have the following additional identities:
   \begin{enumerate}
     \item $h*a=h*b$.
     \item If $n/2=$\,odd then $a*b=p$, $a*a=b*b=u t^{2n/N_L}$.
     \item If $n/2=$\,even then $a*a=b*b=p$, $a*b=u t^{2n/N_L}$.
   \end{enumerate}
\end{prop}
\begin{proof}
   The first three identities and the fact that $h*a=h*b$ are proved
   in~\cite{Beau:quant}. To prove the remaining two identities write
   $n=2k$. Recall from~\cite{Beau:quant} that
   $$(a-b)*(a-b) = (a-b) \cap (a-b) \frac{1}{2}(h^{*n}-4u t^{2n/NL}) =
   (a-b)\cap (a-b)(p-u t^{2n/N_L}).$$
   Substituting~\eqref{Eq:cap-quad}
   in this we obtain:
   \begin{equation} \label{Eq:qprod-quad-1}
      (a-b)*(a-b) = (-1)^k 2 (p-u t^{2n/N_L}).
   \end{equation}
   On the other hand we have $h^{*k} = h^{\cap k} = a+b$, hence
   \begin{equation} \label{Eq:qprod-quad-2}
      (a+b)*(a+b) = h^{*n} = 2p + 2u  t^{2n/N_L}.
   \end{equation}
   Next we claim that $a*a = b*b$. Indeed $a*a-b*b = (a+b)*(a-b) =
   h^{*k}*(a-b)=0$. The desired identities follow from this together
   with~\eqref{Eq:qprod-quad-1},~\eqref{Eq:qprod-quad-2}.
\end{proof}

\subsubsection{Topological restrictions on Lagrangian submanifolds}
\label{Sb:quad-lag}
The quadric $Q$ has Lagrangian spheres. To see this write $Q$ as $Q =
\{z_0^2 + \cdots + z_n^2 = z_{n+1}^2\} \subset {\mathbb{C}}P^{n+1}$.
Then $L = \{ [z_0: \cdots : z_{n+1}] \in Q \mid z_i \in \mathbb{R},
\forall\, i \}$ is a Lagrangian sphere. The following theorem shows
that for $n=$\,even, at least homologically, this is the only type of
Lagrangian with $H_1(L;\mathbb{Z})=0$.
\begin{prop} \label{T:quad-L}
   Assume $\dim_{\mathbb{C}}Q=$\,even. Let $L \subset Q$ be a
   Lagrangian submanifold with $H_1(L;\mathbb{Z})=0$. Then
   $H_*(L;\mathbb{Z}_2) \cong H_*(S^n;\mathbb{Z}_2)$.
\end{prop}
\begin{proof}
   Since $H_1(L;\mathbb{Z})=0$ we have $N_L=2n$. Let
   $\Lambda=\mathbb{Z}_2[t,t^{-1}]$, where $\deg t=-2n$.  Since
   $N_L=2n>n+1$ we have $QH_*(L)\cong (H(L;\mathbb{Z}_2)\otimes
   \Lambda)_*$, hence for every $q \in \mathbb{Z}$, $0 \leq r <2n$ we
   have:
   \begin{equation} \label{Eq:QH-quad}
      QH_{2nq+r}(L) \cong
      \begin{cases}
         H_r(L;\mathbb{Z}_2) & \textnormal{if } 0 \leq r \leq n\\
         0 & \textnormal{if } n+1 \leq r \leq 2n-1
      \end{cases}
   \end{equation}
   Denote by $QH_*(Q)=(H(Q;\mathbb{Z}_2) \otimes \Lambda)_*$ the
   quantum homology of the quadric (over $\mathbb{Z}_2$). Reducing
   modulo $2$ the identities from Proposition~\ref{P:qhom-quadric} it
   follows that $a\in QH_n(Q)$ is an invertible element. Therefore
   $a*:QH_i(L) \to QH_{i-n}$ is an isomorphism for every $i \in
   \mathbb{Z}$.  It now easily follows from~\eqref{Eq:QH-quad} that
   $H_i(L;\mathbb{Z}_2)=0$ for every $0<i<n$.
\end{proof}
We are not aware of existence of any Lagrangian submanifolds in $Q$
with $H_1(L;\mathbb{Z})=0$ which are not diffeomorphic to a sphere,
and it is tempting to conjecture that spheres are indeed the only
example.

Theorem~\ref{T:quad-L} can be also proved by Seidel's method of graded
Lagrangian submanifolds~\cite{Se:graded}. Indeed for $n=$\,even the
quadric has a Hamiltonian $S^1$-action which induces a shift by $n$ on
$QH_*(L)$. To see this write $n=2k$ and write $Q$ as $Q=\{
\sum_{j=0}^k z_jz_{j+1+k}=0\}$. Then $S^1$ acts by $t\cdot [z_0:
\cdots : z_{2k+1}] = [tz_0: \cdots: tz_k:z_{k+1}: \cdots: z_{2k+1}]$.
A simple computation of the weights of the action at a fixed point
gives a shift of $n$ on graded Lagrangian submanifolds in the sense
of~\cite{Se:graded}.

When $n=$\,odd our methods (as well as those of~\cite{Se:graded}) do
not seem to yield a similar result to Theorem~\ref{T:quad-L}. However
the works of Buhovsky~\cite{Bu:TSn} and of Seidel~\cite{Se:TSn} may be
an evidence that such a result should hold.

\subsubsection{Quantum structures for Lagrangians of the quadric}
\label{Sb:quad-quant-L}

Let $L \subset Q$ be a Lagrangian submanifold with
$H_1(L;\mathbb{Z})=0$. Assume that $n=\dim_{\mathbb{C}} Q \geq 2$. As
$N_L = 2n > n+1$ we have a canonical isomorphism $QH_*(L) \cong
(H(L;\mathbb{Z}_2) \otimes \Lambda)_*$, where $\Lambda_* =
\mathbb{Z}_2[t,t^{-1}]$ with $\deg t = -2n$. Denote by $\alpha_0 \in
QH_0(L)$ the class of a point and by $\alpha_n \in QH_n(L)$ the
fundamental class. Denote by $p \in QH_0(Q)$ the class of a point and
by $u \in QH_{2n}(Q)$ the fundamental class. Denote by $i_L:QH_*(L)
\to QH_*(L)$ the quantum inclusion map. With this notation we have:
\begin{prop} \label{T:quad-quant-L}
   Let $L \subset Q$ be as above. Then:
   \begin{enumerate}
     \item $p*\alpha_0 = \alpha_0 t$, $p*\alpha_n = \alpha_n t$.
      \label{quad-p-alpha}
     \item $i_L(\alpha_0) = p - u t$. \label{quad-iL}
     \item If $n$ is even then $\alpha_0 * \alpha_0 = \alpha_n t$.
      \label{quad-prod}
   \end{enumerate}
\end{prop}
\begin{proof}
   By Proposition~\ref{P:qhom-quadric} $p \in QH_0(Q)$ is an
   invertible element, hence $p*:QH_i(L) \to QH_{i-2n}(L)$ is an
   isomorphism for every $i$. But $QH_0(L) \cong \mathbb{Z}_2
   \alpha_0$ and $QH_{-2n}(L) \cong \mathbb{Z}_2 \alpha_0 t$.
   Therefore $p*\alpha_0 = \alpha_0 t$. The statement on $p* \alpha_n$
   is similar. This proves~\ref{quad-p-alpha}.

   We prove~\ref{quad-iL}. It easily follows from the definition of
   the quantum inclusion map that $$i_L(\alpha_0) = p + e u t$$
   for
   some $e \in \mathbb{Z}_2$. Clearly $h * \alpha_0 = 0$ since $h *
   \alpha_0$ belongs to $QH_{-2}(L) \cong QH_{2n-2}(L)=0$.  Therefore
   we have
   $$0=i_L(h*\alpha_0) = h*(p+e u t) = h*p + eh t.$$
   On the other hand a simple computation based on the identities of
   Proposition~\ref{P:qhom-quadric} gives $h*p=h t$. It follows
   that $e=-1$. This proves point~\ref{quad-iL} of the theorem.

   We prove~\ref{quad-prod}. By Proposition~\ref{P:qhom-quadric} when
   $n=$\,even the element $a \in QH_n(Q)$ is invertible (even if we
   work with coefficients in $\mathbb{Z}_2$). Therefore $a*\alpha_n
   =\alpha_0$ and $a*\alpha_0 = \alpha_n t$. It follows that
   $$\alpha_0* \alpha_0 = (a*\alpha_n)*\alpha_0 = a*
   (\alpha_n*\alpha_0) = a*\alpha_0 = \alpha_n t.$$
\end{proof}

Denote $\mathcal{J}$ the space of almost complex structures compatible
with the symplectic structure of $Q$. The following corollary is a
straightforward consequence of Theorem~\ref{T:quad-quant-L}.
\begin{cor} \label{C:quad-disks}
   Let $L \subset Q$ be a Lagrangian submanifold with
   $H_1(L;\mathbb{Z})=0$. Assume $n=\dim_{\mathbb{C}}Q \geq 2$. Then
   the following holds:
   \begin{enumerate}
     \item Let $x \in L$ and $z \in Q \setminus L$. Then for every $J
      \in \mathcal{J}$ there exists a $J$-holomorphic disk $u:(D,
      \partial D) \to (Q,L)$ with $u(-1)=x$, $u(0)=z$ and
      $\mu([u])=2n$.
     \item Assume that $n=$\,even. Let $x', x'', x'' \in L$. Then for
      every $J \in \mathcal{J}$ there exists a $J$-holomorphic disk
      $u:(D, \partial D) \to (Q, L)$ with $u(e^{2\pi i/3})=x'$,
      $u(1)=x''$, $u(e^{4 \pi i/3})=x'''$ and $\mu([u])=2n$.
   \end{enumerate}
\end{cor}

\subsubsection{Lagrangians with $2H_1(L;\mathbb{Z})=0$}
\label{Sb:quad-2-tors}

Let $L \subset Q$ be a Lagrangian submanifold with $H_1(L;\mathbb{Z})
\neq 0$ but $2H_1(L;\mathbb{Z})=0$. Assume that $n > 2$. The quadric
contains several types of such Lagrangians. To see this write the
quadric as $Q=\{ z_0^2 + \cdots + z_{n+1}^2 = 0\}$. For every $0 \leq
r \leq n$ put $L_r = \{[z_0: \ldots: z_{n+1}] \in Q \mid z_0, \ldots,
z_r \in \mathbb{R}, \quad z_{r+1}, \ldots, z_{n+1} \in i
\mathbb{R}\}$.  An easy computation shows that the $L_r$'s are
Lagrangian submanifolds of $Q$ (so called ``real quadrics'') and that
$L_r$ is diffeomorphic to $(S^r \times S^{n-r})/\mathbb{Z}_2$, where
$\mathbb{Z}_2$ acts on both factors by the antipode map. Note that
when $1<r<n-1$ we have $H_1(L_r;\mathbb{Z}) \cong \mathbb{Z}_2$ and
$H_2(L_r;\mathbb{Z}_2) \neq 0$.

A Lagrangian $L$ as above must be monotone. Moreover, when $n=$\,even
its minimal Maslov number is $N_L=n$. To see this note that since the
minimal Chern number of $Q$ is $n$ and $2H_1(L;\mathbb{Z})=0$ we have
$2N_L = 2kn$ for some $k \in \mathbb{N}$. We claim that $k=1$. Indeed
if $k > 1$ then $N_L \geq 2n > n+1$ hence $QH_*(L) \cong
H_*(L;\mathbb{Z}_2) \otimes \mathbb{Z}_2[t, t^{-1}]$ where $\deg t =
-N_L \leq -2n$. In particular $QH_1(L) \cong H_1(L;\mathbb{Z}_2) \neq
0$ and $QH_{n+1}(L)=0$. On the other hand this is impossible since
$a*:QH_{n+1}(L) \to QH_1(L)$ is an isomorphism. A contradiction. This
proves that $k=1$, hence $N_L=n$.

As $N_L=n$ we will work now with the ring $\Lambda_* =
\mathbb{Z}_2[t,t^{-1}]$ where $\deg t =-n$. Put
$QH_*(Q)=(H(Q;\mathbb{Z}_2) \otimes \Lambda)_*$. Note that with this
grading the variable $t$ will appear in the identities of
Proposition~\ref{P:qhom-quadric} in power $2$.

\begin{prop} \label{T:quad-L-2H_1=0}
   Let $L^n \subset Q$ be a Lagrangian submanifold with
   $H_1(L;\mathbb{Z}) \neq 0$ but $2H_1(L;\mathbb{Z})=0$. Assume that
   $n=\dim_{\mathbb{C}}Q \geq 4$ and if $n=$\,odd assume also that
   $N_L=n$. Suppose there exists $2 \leq j \leq n-2$ such that
   $H_j(L;\mathbb{Z}_2) \neq 0$. Then there exists an isomorphism of
   graded vector spaces $QH_*(L) \cong (H(L;\mathbb{Z}_2) \otimes
   \Lambda)_*$. This isomorphism is canonical for general degree $*$
   except perhaps for $* \equiv 0 (\bmod{n})$, but there is a
   canonical injection $H_n(L;\mathbb{Z}_2) \hookrightarrow QH_n(L)$
   identifying the fundamental class $\alpha_n \in H_n(L)$ with the
   unit of $QH(L)$. Moreover, let $\alpha_0 \in QH_0(L)$ be an element
   such that $\{ \alpha_0, \alpha_n t\}$ form a basis for $QH_0(L)
   \cong H_0(L;\mathbb{Z}_2) \oplus H_n(L;\mathbb{Z}_2)t$.  Then one
   of the following two possibilities occurs:
   \begin{enumerate}
     \item $p*\alpha_0 = \alpha_n t^3, \;\; p*\alpha_n =
      \alpha_0 t$. \label{I:p-alpha-1}
     \item $p*\alpha_0 = \alpha_0 t^2 + s \alpha_n
      t^3, \;\; p*\alpha_n = r \alpha_0 t + \alpha_n t^2$, for
      some $s,r \in \mathbb{Z}_2$ with $sr=0$. \label{I:p-alpha-2}
   \end{enumerate}
   In case $n=$\,odd we also have:
   \begin{enumerate}[(1')]
     \item Either $i_L(\alpha_0)=p$, \;\; $i_L(\alpha_n)=u t$;
      \label{I:iLalpha-1} or
     \item $i_L(\alpha_0)=p+u t^2$, \;\; $i_L(\alpha_n) = r u t$,
      \label{I:iLalpha-2}
   \end{enumerate}
   Where~\ref{I:iLalpha-1}' corresponds to
   possibility~\ref{I:p-alpha-1} above and~\ref{I:iLalpha-2}' to
   possibility~\ref{I:p-alpha-2}.
\end{prop}

Before proving Proposition~\ref{T:quad-L-2H_1=0}, here is an immediate
corollary of it:
\begin{cor} \label{C:quad-disks-2H_1=0}
   Let $L \subset Q$ be a Lagrangian submanifold as in
   Theorem~\ref{T:quad-L-2H_1=0}. Let $x \in L$ and $z \in Q \setminus
   L$. Then for every $J \in \mathcal{J}$ there exists a
   $J$-holomorphic disk $u:(D, \partial D) \to (Q,L)$ with $u(-1)=x$,
   $u(0)=z$ and $\mu([u])=n$.
\end{cor}

\begin{proof}[Proof of Proposition ~\ref{T:quad-L-2H_1=0}]
   We first show that $QH_*(L) \cong (H(L;\mathbb{Z}_2) \otimes
   \Lambda)_*$. Let $f:L \to \mathbb{R}$ be a Morse function and
   denote by $CM_*(f) = \mathbb{Z}_2 \langle \textnormal{Crit}(f)
   \rangle$ the Morse complex of $f$. Put $CF_* = (C(f) \otimes
   \Lambda)_*$. Since $N_L = n$ the Floer differential $d$ can be
   written as $d = \partial_0 + \partial_1 t$, where
   $\partial_0:CM_*(f) \to CM_{*-1}(f)$ is the Morse differential and
   $\partial_1$ is an operator acting as $\partial_1: CM_*(f) \to
   CM_{*-1+n}(f)$.  Put $E_*^{1} = (H(L;\mathbb{Z}_2) \otimes
   \Lambda)_*$.  By the results of~\cite{Bi:Nonintersections},
   $\partial_1 t$ descends to a differential $d_1$ defined on the
   homology $E^1_* = (H(C(f), \partial_0) \otimes \Lambda)_* \cong
   (H(L;\mathbb{Z}_2) \otimes \Lambda)_*$, which has the form $d_1 =
   \delta_1 t$ where $\delta_1: H_*(L;\mathbb{Z}_2) \to
   H_{*-1+n}(L;\mathbb{Z}_2)$.  Moreover the homology $H_*(d_1)$ is
   isomorphic to $QH_*(L)$.  By grading reasons $\delta_1$ is zero on
   $H_k(L;\mathbb{Z}_2)$ for every $k \geq 2$.

   We claim that $\delta_1=0$ (also on $H_0(L;\mathbb{Z}_2)$ and
   $H_1(L;\mathbb{Z}_2)$). To prove this we use the work of
   Buhovsky~\cite{Bu:products}, by which $\delta_1$ satisfies Leibniz
   rule with respect to the cap product on $H_*(L;\mathbb{Z}_2)$, i.e.
   $\delta_1 (\alpha \cap \beta) = \delta_1(\alpha) \cap \beta +
   \alpha \cap \delta_1(\beta)$ for every $\alpha, \beta \in
   H_*(L;\mathbb{Z}_2)$. By assumption there exists a non-trivial
   element $\alpha_j \in H_j(L;\mathbb{Z}_2)$ for some $2 \leq j \leq
   n-2$. By Poincar\'{e} duality there exists $\alpha_{n-j} \in
   H_{n-j}(L;\mathbb{Z}_2)$ such that $\alpha_j \cap \alpha_{n-j} =
   \gamma$ where $\gamma \in H_0(L;\mathbb{Z}_2)$ is the class of a
   point. As $\delta_1(\alpha_j)$, $\delta_1(\alpha_{n-j})=0$ it
   follows from Leibniz rule that $\delta_1(\gamma)=0$. For degree
   reasons $\gamma$ cannot be a $d_1$-boundary. Therefore $QH_0(L)
   \neq 0$. Next, note that $\delta_1$ maps $H_1(L;\mathbb{Z}_2)$ to
   $H_n(L;\mathbb{Z}_2)=\mathbb{Z}_2 \alpha_n$, where $\alpha_n$ is
   the fundamental class. If $\delta_1$ were not $0$ on
   $H_1(L;\mathbb{Z}_2)$ we would get that $\alpha_n =
   \delta_1(\alpha_1)$ for some $\alpha_1 \in H_1(L;\mathbb{Z}_2)$,
   hence $[\alpha_n]$ would be $0$ in $QH(L)$. Since $[\alpha_n]$ is
   the unity in $QH_*(L)$ it follows that $QH_*(L)=0$, a
   contradiction.  This proves that $\delta_1$ vanishes also on
   $H_1(L;\mathbb{Z}_2)$. Summarizing, we have $\delta_1=0$ and
   therefore $QH_*(L) \cong (H(L) \otimes \Lambda)_*$ as claimed.
   Note that $QH_0(L) \cong H_0(L;\mathbb{Z}_2) \oplus
   H_n(L;\mathbb{Z}_2)t$. The statement on the canonicity of these
   isomorphisms for various values of $*$ follows easily from degree
   reasons.

   We turn to the proof of the other statements of the theorem.  Write
   \begin{equation} \label{Eq:p-alpha}
      p*\alpha_0 = s_0  \alpha_0 t^2 + s_n \alpha_n t^3,
      \quad
      p*\alpha_n = t_0 \alpha_0 t + t_n \alpha_n t^2,
   \end{equation}
   for some $s_0, s_n, t_0, t_n \in \mathbb{Z}_2$. By
   Proposition~\ref{P:qhom-quadric} we have $p*p = u t^4$.
   Multiplying both sides of equations~\eqref{Eq:p-alpha} by $p$ and
   comparing coefficients on both sides we obtain three possibilities
   for the values of $s_0, s_n, t_0, t_n$:
   \begin{enumerate}
     \item $s_0=0$, $s_n=1$, $t_0=1$, $t_n=0$.
     \item $s_0=1$, $t_0=0$, $t_n=1$.
     \item $s_0=1$, $t_0=1$, $s_n=0$, $t_n=1$.
   \end{enumerate}
   The first case leads to possibility~\ref{I:p-alpha-1} of the
   theorem. The two other two cases lead to
   possibility~\ref{I:p-alpha-2}.

   We now prove identities~\ref{I:iLalpha-1}',~\ref{I:iLalpha-2}'
   assuming $n=$\,odd. We can write $i_L(\alpha_n)=g u t$ for some $g
   \in \mathbb{Z}_2$. Since $\{ \alpha_0, \alpha_n t\}$ are linearly
   independent it is easy to see that $\epsilon_L(\alpha_0)=1$. It
   follows that $i_L(\alpha_0)$ has the form $i_L(\alpha_0) = p + d u
   t^2$ for some $d \in \mathbb{Z}_2$.  Using the fact that
   $i_L(p*\alpha_j) = p*i_L(\alpha_j)$ and the
   identities~\ref{I:p-alpha-1},~\ref{I:p-alpha-2} just proved, gives
   the desired result on $d,g$.
\end{proof}

\begin{rem} \label{R:quad-2H1=0}
   \begin{enumerate}
     \item Using similar methods, when $n=$\,even one can obtain
      information on $a*\alpha_0$, $b*\alpha_0$, $a* \alpha_n$, $a*
      \alpha_n$ and also on $i_L(\alpha_0)$, $i_L(\alpha_n)$. We leave
      these computations to the reader.
     \item The identities for $p*\alpha_j$ in
      Theorem~\ref{T:quad-L-2H_1=0} suggest that for some Lagrangians
      in $Q$ a stronger version of Corollary~\ref{C:quad-disks-2H_1=0}
      should hold. An interesting case seems to be for example when
      identity~\ref{I:p-alpha-1} of Theorem~\ref{T:quad-L-2H_1=0}
      holds.
     \item As was shown above, when $n=$\,even, $H_1(L;\mathbb{Z})=0$,
      $2H_1(L;\mathbb{Z}) \neq 0$ we have $N_L=n$. A priori this need
      not be the case when $n=$\,odd (i.e. it may happen that
      $N_L=2n$). However we are not aware of any examples of this
      sort.
   \end{enumerate}
\end{rem}

\subsection{Complete intersections} \label{Sb:cinter}
Some of the result in this section and in \S\ref{Sb:alg-qh} are
somewhat non-rigorous. The reason for this is that we need to use here
Floer homology with coefficients in $\mathbb{Q}$ rather than
$\mathbb{Z}_2$. For this end one has to orient the moduli spaces of
pearls arising in our constructions.  Orientations of moduli spaces of
pseudo-holomorphic disks in the context of Floer homology have been
worked out by Fukaya, Oh, Ohta and Ono~\cite{FO3}. Still, it remains
to check whether some ingredients of the theory of the present paper,
e.g. the quantum module structure and the quantum inclusion maps, are
indeed compatible with the coherent orientations arising in Floer
homology over $\mathbb{Q}$. While it seems very likely that the two
theories are compatible we have not checked all the details. On the
other hand we found it worth presenting here possible applications of
our theory over $\mathbb{Q}$.  Below we state such applications with
the convention that theorems and corollaries are marked with $*$
whenever they depend on the validity of our theory over $\mathbb{Q}$.

\

Let $X \subset {\mathbb{C}}P^{n+r}$ be a smooth complete intersection
of degree $(d_1, \ldots, d_r)$, i.e. a transverse intersection of $r$
complex hypersurfaces of degrees $d_1, \ldots, d_r$ in
${\mathbb{C}}P^{n+r}$. Note that $\dim_{\mathbb{C}}X = n$. We will
assume that $n\geq 3$ hence by the Lefschetz hyperplane section
theorem we have $\dim H^2(X;\mathbb{R}) \cong \mathbb{R}$. It follows
by Moser argument that all K\"{a}hler forms on $X$ are symplectically
equivalent up to a constant factor. We endow $X$ with the symplectic
structure induced from ${\mathbb{C}}P^{n+r}$. We will also assume that
$X$ is non-linear, i.e. that at least one of the $d_i$'s is $>1$.  Put
$N=n+r+1-\sum_{i=1}^r d_i$. When $N>0$, $X$ is monotone with minimal
Chern number $N$.

Note that such a complete intersection $X$ must have Lagrangian
submanifolds $L$ with $H_1(L;\mathbb{Z})=0$. Indeed, $X$ can be
degenerated to a variety with isolated singularities, hence by
symplectic Picard-Lefschetz theory (see e.g.~\cite{Ar:monodromy,
  Do:polynom, Se:vcycles-mut, Se:thesis, Bi:algf}) $X$ must have
Lagrangian spheres.

Let $\mathbb{K}$ be one of the fields $\mathbb{Z}_2$ or $\mathbb{Q}$.
Let $\Lambda^{\mathbb{K}}_*=\mathbb{K}[t,t^{-1}]$ where $\deg t =
-2N$. Let $QH_*(X;\Lambda^{\mathbb{K}}) = (H(X;\mathbb{K}) \otimes
\Lambda^{\mathbb{K}})_*$ be the quantum homology of $X$ with
coefficients in $\Lambda^{\mathbb{K}}$. Let $L \subset X$ be a
Lagrangian submanifold with $H_1(L;\mathbb{Z})=0$. Clearly $L$ is
monotone with minimal Maslov number $N_L=2N$. Denote by
$QH_*(L;\Lambda^{\mathbb{K}})$ the Floer homology of $L$.  When
$\mathbb{K}=\mathbb{Q}$ we assume that $L$ is orientable and
relatively spin (see~\cite{FO3}). Note that when $n \geq 2\sum_{i=1}^r
d_i -2r$ we have $2N \geq n+2$ hence there exists a canonical
isomorphism $QH_*(L;\Lambda^{\mathbb{K}}) \cong (H(L;\mathbb{K})
\otimes \Lambda^{\mathbb{K}})_*$. In this case denote by $\alpha_0 \in
H_0(L;\mathbb{K})$ the class of a point in $L$, by $p \in
H_0(X;\mathbb{K})$ the class of a point in $X$ and by $h \in
H_{2n-2}(X;\mathbb{K})$ the class of a hyperplane section of $X$.
\begin{propspec}
   Let $L \subset X$ be a Lagrangian submanifold as above. Assume that
   $n \geq 3$ and $n \geq 2\sum_{i=1}^r d_i -2r+1$. Then
   $$i_L(\alpha_0) = p - \bigl(\prod_{i=1}^r(d_i-1)!\bigr)
   h^{\cap(n-N)} t.$$
   In particular for every $x \in L$ and
   every almost complex structure $J$ compatible with the symplectic
   structure of $X$ there exists a $J$-holomorphic disk $u:(D,
   \partial D) \to (X, L)$ with $\mu([u])=2N$ and $x \in u(\partial
   D)$. In fact, the number of such disks intersecting a complex
   $N$-dimensional hyperplane in $X$ is $\prod_{i=1}^r(d_i-1)!$, when
   counted appropriately.
\end{propspec}
The identity on $i_L(\alpha_0)$ is completely rigorous when
$\mathbb{K}=\mathbb{Z}_2$ but the coefficient of $h^{\cap(n-N)}t$ is
$0 \in \mathbb{Z}_2$ unless $d_i \leq 2$ for every $i$.

\begin{proof}
   Put $\mathbb{K} = \mathbb{Q}$.  Since $4N > 2n$ we have
   $QH_0(X;\Lambda^{\mathbb{K}}) = H_0(X;\mathbb{K}) \oplus
   H_{2N}(X;\mathbb{K}) t$.  Also note that since $2N \geq
   n+3$ we have $H_{2N}(X;\mathbb{K}) \cong \mathbb{K} h^{\cap(n-N)}$.
   Therefore we can write
   \begin{equation} \label{Eq:cinter-iL}
      i_L(\alpha_0) = p + \tau h^{\cap (n-N)} t
   \end{equation}
   for some $\tau \in \mathbb{K}$.

   We first claim that $h*\alpha_0 = 0$. Indeed $h*\alpha_0 \in
   QH_{-2}(L;\Lambda^{\mathbb{K}}) \cong H_{2N-2}(L;\mathbb{K}) t$.
   But $2N-2 \geq n+1$ hence $H_{2N-2}(L;\mathbb{K})=0$.

   A straightforward computation based on the results of
   Beauville~\cite{Beau:quant} gives:
   \begin{equation*}
      h*p = \bigl( \prod_{i=1}^r (d_i-1)! \bigr) h^{\cap (n-N+1)} t,
      \quad
      h^{*j} = h^{\cap j} \; \forall\,\, 0 \leq j \leq N.
   \end{equation*}
   We now obtain from~\eqref{Eq:cinter-iL} that $$0 = i_L(h*\alpha_0)
   = h*p + \tau(h*h^{(n-N)})T = \prod_{i=1}^r (d_i-1)!h^{\cap (n-N+1)}
   t + \tau h^{\cap (n-N+1)} t.$$
   Here we have used the fact that
   $n-N+1 \leq N$ hence $h*h^{(n-N)} = h^{\cap (n-N+1)} \neq 0$. It
   immediately follows that $\tau = - \prod_{i=1}^r (d_i-1)!$.
\end{proof}

\subsection{Algebraic properties of quantum homology and
  Lagrangian submanifolds} \label{Sb:alg-qh}

Here we present relations between algebraic properties of the quantum
homology of an (ambient) symplectic manifold and the kind of
Lagrangians it contains.  Interestingly enough the existence of
certain types of Lagrangian submanifolds (e.g. spheres) imposes strong
restrictions on the quantum homology of the ambient manifold.

Let $(M^{2n},\omega)$ be a spherically monotone symplectic manifold
with minimal Chern number $N$. Let $\mathbb{K}$ be either $\mathbb{Q}$
or $\mathbb{Z}_2$. Let $\Lambda^{\mathbb{K}}_*$ be either
$\mathbb{K}[t,t^{-1}]$ or the field $\mathbb{K}[t]]$ of formal Laurent
series with finitely many terms having negative powers of $t$. Define
$\deg t = -2N$. Let $QH_*(M;\Lambda^{\mathbb{K}}) = (H(M;\mathbb{K})
\otimes \Lambda^{\mathbb{K}})_*$ be the quantum homology over
$\mathbb{K}$ endowed with the quantum cap product. Let
$$\textnormal{pr}_l:QH_*(M;\Lambda^{\mathbb{K}}) \longrightarrow
QH_l(M;\Lambda^{\mathbb{K}}) = \bigoplus_{k \in \mathbb{Z}}
H_{l+2Nk}(M;\mathbb{K}) t^k$$
be the projection on the
degree-$l$ component of $QH$.

\begin{propspec} \label{T:Lag-sphere-inv}
   Suppose that $(M^{2n},\omega)$ has a Lagrangian sphere. Assume that
   $n=\dim_{\mathbb{C}}M \geq 2$ and $N \nmid (n+1)$.  Let $a \in
   QH_*(M;\Lambda^{\mathbb{K}})$ be an invertible element.  Then
   either $\textnormal{pr}_l(a) \neq 0$ for some $l$ with $l \equiv
   2n(\bmod{2N})$ or there exist two indices $l_1,l_2$ with $l_1
   \equiv n (\bmod{2N})$, $l_2 \equiv 3n(\bmod{2N})$ such that
   $\textnormal{pr}_{l_1}(a)$, $\textnormal{pr}_{l_2}(a) \neq 0$.  If
   we assume in addition that $N \nmid n$ then the indices $l_1$,
   $l_2$ must be distinct hence the only invertible elements $a \in
   QH(M;\Lambda^{\mathbb{K}})$ of pure degree satisfy $\deg{a} \equiv
   2n (\bmod{2N})$.

   As before, the statement of the Theorem is completely rigorous when
   $\mathbb{K}=\mathbb{Z}_2$.
\end{propspec}
\begin{proof}
   Let $L \subset M$ be a Lagrangian sphere.  Clearly $L$ is monotone
   with minimal Maslov number $2N$. Since $N \nmid n+1$, a standard
   argument as in~\cite{Bi:Nonintersections} shows that Oh's spectral
   sequence collapses at stage $r=1$, hence
   $QH_*(L;\Lambda^{\mathbb{K}}) \cong (H(S^n;\mathbb{K}) \otimes
   \Lambda^{\mathbb{K}})_*$. Choose two non-zero elements $\alpha_0
   \in QH_0(L)$, $\alpha_n \in QH_n(L)$ (e.g. under the previous
   isomorphism we can take $\alpha_0$, $\alpha_n$ to correspond to the
   class of a point in $H_0(S^n;\mathbb{K})$ and to the fundamental
   class in $H_n(S^n; \mathbb{K})$). Note that when $n \geq 3$, $L
   \approx S^n$ is automatically relatively spin hence $QH(L)$ is well
   defined in case $\mathbb{K}=\mathbb{Q}$. The same holds for $n=2$
   since the 2'nd Stiefel-Whitney class of $S^2$ vanishes.

   Write $a = \sum_{j=k_1}^{k_2} a_j$, for some $k_1\leq k_2$, where
   $a_j=\textnormal{pr}_j(a)$.  Since $a$ is invertible there exists
   $l_1$ such that $a_{l_1}*\alpha_n \neq 0$. As $a_{l_1}*\alpha_n \in
   QH_{l_1-n}(L)$ it follows that either $l_1 \equiv 2n(\bmod{2N})$ or
   $l_1 \equiv n(\bmod{2N})$. Similarly there exists $l_2$ such that
   $\alpha_{l_2} * \alpha_0 \neq 0$. This leads to the following two
   possibilities: either $l_2 \equiv 2n(\bmod{2N})$ or $l_2 \equiv
   3n(\bmod{2N})$.
\end{proof}

\begin{propspec} \label{T:0-div}
   Let $(M,\omega)$ be a spherically monotone (resp. monotone)
   symplectic manifold. Suppose that $(M,\omega)$ contains a
   Lagrangian submanifold $L$ which is simply connected (resp. has
   $H_1(L;\mathbb{Z})=0$) and is relatively spin. Assume that the
   minimal Chern number $N$ of $(M,\omega)$ satisfies $N \geq
   \frac{n}{2}+1$ where $n=\dim_{\mathbb{C}}M$. Then for every
   $3n+1-2N \leq l \leq 2n-1$ all elements of
   $QH_l(M;\Lambda^{\mathbb{K}})$ are divisors of $0$.

   As before, the statement of the Theorem is completely rigorous when
   $\mathbb{K}=\mathbb{Z}_2$.
\end{propspec}
\begin{proof}
   Clearly $L$ is monotone with minimal Maslov number $N_L = 2N \geq
   n+2$. Therefore there exists a canonical isomorphism
   $QH_*(L;\Lambda^{\mathbb{K}}) \cong (H(L;\mathbb{K}) \otimes
   \Lambda^{\mathbb{K}})_*$. In particular
   $QH_j(L;\Lambda^{\mathbb{K}}) = 0 \;\; \forall\, n+1-2N \leq j \leq
   -1$.

   Let $i_L:QH_*(L;\Lambda^{\mathbb{K}}) \longrightarrow
   QH_*(M;\Lambda^{\mathbb{K}})$ be the quantum inclusion homomorphism
   (see \S\ref{Sb:q-inclusion}). Let $x \in H_0(L;\mathbb{K})$ be the
   class of a point. Put $a = i_L(x) \in
   QH_0(M;\Lambda^{\mathbb{K}})$. Note that $a \neq 0$ since $i_L(x)$
   has the form $p + \sum_{i \geq 1}a_i t^i$, where $p \in
   H_0(M;\mathbb{K})$ is the class of a point and $a_i \in
   H_{2iN}(M;\mathbb{K})$.

   We claim that for every $b \in QH_l(M;\Lambda^{\mathbb{K}})$ with
   $3n+1-2N \leq l \leq 2n-1$ we have $b*a=0$. Indeed let $b$ be such
   an element. Then $b*a = i_L(b*x)$. But $b*x \in
   QH_{l-2n}(L;\Lambda^{\mathbb{K}})=0$ since $n+1-2N \leq l-2n \leq
   -1$.
\end{proof}

\subsubsection*{Semi-simplicity of quantum homology}
Recall that a commutative algebra $A$ over a field $\mathbb{F}$ is
called semi-simple if it splits into a direct sum of finite
dimensional vector spaces over $\mathbb{F}$, $A = A_1 \oplus \cdots
\oplus A_r$ such that the splitting is compatible with the
multiplication of $A$ (i.e.  $(a_1, \ldots, a_r)\cdot(a'_1, \ldots,
a'_r) = (a_1 a'_1, \ldots, a_r a'_r)$) and such that each $A_i$ is a
field with respect to the ring structure induced from $A$.
\begin{rem} \label{R:semi-simple}
   There exist several different notions of semi-simplicity in the
   context of quantum homology, or more generally in the context of
   Frobenius algebras and Frobenius manifolds (see
   e.g.~\cite{Dubrovin:geom-2d-tft, Ko-Ma:GW, Tian-Xu:semisimple}).
   The semi-simplicity we use here was first considered in the context
   of quantum homology by Abrams~\cite{Abrams:semi-simple}. It is in
   general not equivalent to semi-simplicity in the sense of
   Dubrovin~\cite{Dubrovin:geom-2d-tft} since we work with a different
   coefficient ring.
\end{rem}

Let $(M^{2n},\omega)$ be a spherically monotone symplectic manifold
with minimal Chern number $N$. For simplicity we will work here with
the even quantum homology
$$QH_*^{\textnormal{ev}}(M; \mathbb{F}) = \bigoplus_{i \in \mathbb{Z}}
QH_{2i}(M;\mathbb{F}) = \bigoplus_{j=0}^n H_{2j}(M;\mathbb{Q}) \otimes
\mathbb{F},$$
where $\mathbb{F}$ is the field $\mathbb{F} =
\mathbb{Q}[t]]$ with $\deg t = -2N$. We could work here also with the
full quantum homology but then semi-simplicity should be considered in
the framework of skew-commutative algebras.

\begin{propspec} \label{T:semi-simple-1}
   Let $(M^{2n}, \omega)$ be a monotone (resp. spherically monotone)
   symplectic manifold with minimal Chern number $N \geq
   \frac{n}{2}+1$. Assume that $QH^{ev}_*(M;\mathbb{F})$ is
   semi-simple. Let $L \subset (M, \omega)$ be a Lagrangian
   submanifold which is relatively spin and has $H_1(L;\mathbb{Z})=0$
   (resp. $\pi_1(L)=1$). Put
   \begin{align*}
      & I = \{ 0 \leq i \leq n \mid i \equiv ln (\bmod{2N})
      \textnormal{ for some } l \in \mathbb{Z} \}, \\
      & J = \{ 0 \leq j \leq n \mid j \equiv q + 2nk (\bmod{2N})
      \textnormal{ for some } n+1 \leq q \leq 2N-1, k \in \mathbb{Z}\}, \\
      & J' = \{ 0 \leq j \leq n \mid j \equiv 1 + nk' (\bmod{2N})
      \textnormal{ for some } k' \in \mathbb{Z} \}, \\
      & J'' = \{ 0 \leq j \leq n \mid j \equiv -1 + nk'' (\bmod{2N})
      \textnormal{ for some } k'' \in \mathbb{Z} \}.
   \end{align*}
   Then $H_i(L;\mathbb{Q}) \cong \mathbb{Q}$ for every $i \in I$ and
   $H_i(L;\mathbb{Q})=0$ for every $j \in J \cup J' \cup J''$.
\end{propspec}
We will give the proof of Proposition~\ref{T:semi-simple-1} later in this
section. In the meanwhile here is an immediate corollary.
\begin{corspec} \label{C:semi-simple-1}
   Let $(M^{2n}, \omega)$ be a monotone (resp. spherically monotone)
   symplectic manifold with minimal Chern number $N$. Assume that
   $\frac{3n+1}{4} \leq N \leq n-1$ and that $(M,\omega)$ has a
   Lagrangian submanifold $L$ which is relatively spin and has
   $H_1(L;\mathbb{Z})=0$ (resp. $\pi_1(L)=1$). Then
   $QH^{ev}_*(M;\mathbb{F})$ is not semi-simple.
\end{corspec}
\begin{proof}[Proof of Corollary~\ref{C:semi-simple-1}]
   Suppose that $QH^{ev}_*(M;\mathbb{F})$ is semi-simple. Let $J$ be
   the set of indices defined in Theorem~\ref{T:semi-simple-1}.  As
   $\frac{3n+1}{4} \leq N \leq n-1$, we have $n+2 \leq 3n-2N \leq
   2N-1$ hence $n \in J$ (take $q = 3n-2N$, $k=-1$).  By
   Theorem~\ref{T:semi-simple-1}, $H_n(L;\mathbb{Q})=0$ which is
   impossible since $L$ is orientable. A contradiction.
\end{proof}

Here is another restriction on semi-simplicity.
\begin{propspec} \label{T:semi-simple-shpere}
   Let $(M^{2n},\omega)$ be a spherically monotone symplectic manifold
   of $\dim_{\mathbb{C}}M \geq 2$. Suppose that $(M,\omega)$ contains
   a Lagrangian sphere and that its minimal Chern number $N$ satisfies
   $N \nmid n$ and $N \nmid (n+1)$. Then $QH_*^{ev}(M;\mathbb{F})$ is
   not semi-simple.
\end{propspec}
In order to prove
Propositions~\ref{T:semi-simple-shpere},~\ref{T:semi-simple-1} we will
need some preparations regarding semi-simplicity.  Let
$$\eta:QH_*(M;\mathbb{F}) \longrightarrow H_0(M;\mathbb{Q}) \otimes
\mathbb{F} \cong \mathbb{F}$$
be the projection. The identification of
the last isomorphism here is made via $H_0(M;\mathbb{Q}) =
\mathbb{Q}p$, where $p$ is the class of a point.  The projection
$\eta$ assigns to an element $a \in
QH^{\textnormal{ev}}_*(M;\mathbb{F})$ a power series in $t$ which is
the coefficient of $a$ at $p$. Define a pairing
\begin{equation*}
   \Delta:QH^{\textnormal{ev}}_*(M;\mathbb{F}) \otimes
   QH^{\textnormal{ev}}_*(M;\mathbb{F}) \longrightarrow \mathbb{F},
   \quad \Delta(a,b) =  \eta (a*b).
\end{equation*}
It is well known that $\Delta$ is a non-degenerate pairing
(see~\cite{McD-Sa:Jhol-2}). Let $e_1, \ldots, e_{\nu}$ be a basis over
$\mathbb{F}$ of $QH^{\textnormal{ev}}_*(M;\mathbb{F})$ and denote by
$e_1^{\#}, \ldots, e_{\nu}^{\#}$ the dual basis with respect to
$\Delta$.  Define an element $\mathcal{E}_Q^{ev}(M) \in
QH^{\textnormal{ev}}_*(M;\mathbb{F})$ by
\begin{equation} \label{Eq:q-euler}
   \mathcal{E}_Q^{ev}(M) = \sum_{i=1}^{\nu} e_i*e_i^{\#}.
\end{equation}
This element is called the even quantum Euler class. It does not
depend on the choice of the basis $e_1, \ldots, e_{\nu}$. A Theorem
due to Abrams~\cite{Abrams:semi-simple} asserts that {\sl
  $QH^{\textnormal{ev}}_*(M;\mathbb{F})$ is semi-simple iff
  $\mathcal{E}_Q^{ev}(M)$ is invertible}. Of course, it is possible to
define a more complete quantum Euler class by taking an alternate sum
as in~\eqref{Eq:q-euler} over a basis of the whole of
$QH_*(M;\mathbb{F})$ (not just the even part). In this setting the
quantum Euler is indeed a deformation of the classical Euler class,
hence its name.

We will now need the following proposition.  Denote by $\cdot:
H_*(M;\mathbb{Q}) \otimes H_*(M;\mathbb{Q}) \to \mathbb{Q}$ the
classical intersection pairing, with the convention that $a \cdot b =
0$ whenever $a,b$ are elements of pure degree with $\deg(a) + \deg(b)
\neq 2n$.
\begin{prop} \label{P:QH-eta}
   Let $e',e'' \in H_*(M;\mathbb{Q})$ and view $e',e''$ as elements of
   $QH_*(M;\mathbb{F})$. Then $\eta(e'*e'') = e' \cdot e''$.
\end{prop}
\begin{proof}
   Without loss of generality we may assume that $e'$, $e''$ have pure
   degrees. Then $\eta(e'*e'') = st^j$ for some $j \geq 0$ and $s \in
   \mathbb{Q}$. We may also assume that $s \neq 0$ since otherwise the
   statement is obvious.

   We claim that $j=0$ . To prove this, choose generic
   $\mathbb{Q}$-cycles $C', C''$ representing $e', e''$.  Suppose that
   $j>0$ and $s \neq 0$. Then for a generic $\omega$-compatible almost
   complex structure $J$ there exists a simple $J$-holomorphic
   rational curve $u$ passing through $C'$ and $C''$ with
   $c_1([u])=jN$. Denote by $\mathcal{M}(J)$ the space of simple
   $J$-holomorphic curves in the class $[u]$. Consider the evaluation
   map:
   \begin{equation*}
      ev: \bigl(\mathcal{M}(J)\times
      \mathbb{C}P^1 \times \mathbb{C}P^1 \bigr)/G
      \longrightarrow M \times M, \quad ev(w,p_1, p_2) = (w(p_1), w(p_2)),
   \end{equation*}
   where $G = \textnormal{Aut}(\mathbb{C}P^1)$.  As $ev^{-1}(C' \times
   C'') \neq \emptyset$ we obtain that $\dim(ev^{-1}(C' \times C'') =
   \deg(e') + \deg(e'')-2n + 2Nj-2 \geq 0$. On the other hand since
   $\eta(e'*e'') = st^j$, we have $\deg(e')+\deg(e'') - 2n = -2Nj$. A
   contradiction. This proves that $j=0$, hence $\eta(e'*e'') = s = e'
   \cdot e''$.
\end{proof}

\begin{lem} \label{L:qeuler-deg-0}
   The even quantum Euler class is an element of pure degree $0$.
\end{lem}
\begin{proof}
   Let $e_1, \ldots, e_{\nu}$ be elements of pure degree which form a
   basis of $H^{ev}_*(M;\mathbb{Q})$ over $\mathbb{Q}$. Let $e_1^{\#},
   \ldots, e_{\nu}^{\#} \in H^{ev}_*(M;\mathbb{Q})$ be a dual basis to
   $e_1, \ldots, e_{\nu}$ with respect to the classical intersection
   pairing $(-) \cdot (-)$. Thus $e_i \cdot e_j^{\#} = \delta_{i,j}$.
   Note that both $\{e_i\}$ and $\{e_i^{\#}\}$ are bases of
   $QH^{ev}_*(M,\mathbb{F})$ over $\mathbb{F}$.  We claim that the
   basis $\{e_i^{\#}\}$ is dual to $\{e_i\}$ also with respect to the
   pairing $\Delta$.  (This is, by the way, contrary to what is
   written in~\cite{Abrams:semi-simple}.)  Indeed, by
   Proposition~\ref{P:QH-eta} we have $\Delta(e_i,e_j^{\#}) = \eta(e_i
   * e_j^{\#}) = e_i\cdot e_j^{\#} = \delta_{i,j}$.

   Finally, recall that the even quantum Euler class does not depend
   on the choice of the basis, hence $\mathcal{E}_Q^{ev}(M) =
   \sum_{i=1}^{\nu} e_i*e_i^{\#}$. Since $\deg({e_i^{\#}}) =
   2n-\deg(e_i)$ it follows that $\mathcal{E}_Q^{ev}(M)$ has degree
   $0$.
\end{proof}

We are now in position to prove
Propositions~\ref{T:semi-simple-1},~\ref{T:semi-simple-shpere}.
\begin{proof}[Proof of Proposition~\ref{T:semi-simple-shpere}]
   By Proposition~\ref{T:Lag-sphere-inv} $QH_*(M;\mathbb{F})$ has no
   invertible elements of pure degree $\not\equiv 2n (\bmod{2N})$.  As
   $0 \not\equiv 2n (\bmod{2N})$, the even quantum Euler class
   $\mathcal{E}_Q^{ev}(M)$ is not invertible, hence
   $QH^{ev}_*(M;\mathbb{F})$ is not semi-simple.
\end{proof}

\begin{proof}[Proof of Proposition~\ref{T:semi-simple-1}]
   The Lagrangian $L^n \subset M^{2n}$ is monotone and has minimal
   Maslov number $2N$. By assumption $2N \geq n+2$, hence there exists
   a canonical isomorphism $QH_*(L;\mathbb{F}) \cong (H(L;\mathbb{Q})
   \otimes \mathbb{F})_*$, or more specifically:
   \begin{equation} \label{Eq:semi-simple-QH}
      \begin{aligned}
         & QH_{i+2kN}(L;\mathbb{F}) \cong H_i(L;\mathbb{Q}), \quad
         \forall
         \, 0 \leq i \leq n, k \in \mathbb{Z}, \\
         & QH_{j+2kN}(L;\mathbb{F})=0, \quad \forall \, n+1 \leq j
         \leq 2N-1, k\in \mathbb{Z}.
      \end{aligned}
   \end{equation}
   Since $QH^{ev}_*(M,\mathbb{F})$ is semi-simple the even quantum
   Euler class $\mathcal{E}_Q^{ev}(M) \in QH_0(M;\mathbb{F})$ is
   invertible.  Therefore exterior multiplication by
   $\mathcal{E}_Q^{ev}(M)$ gives isomorphisms $QH_l(L;\mathbb{F})
   \cong QH_{l-2n}(L;\mathbb{F})$ for every $l \in \mathbb{Z}$. The
   rest of the proof follows from these isomorphisms together
   with~\eqref{Eq:semi-simple-QH} and the fact that for every $k \in
   \mathbb{Z}$ we have:
   \begin{align*}
      & QH_{2kN}(L;\mathbb{F}) \cong H_0(L;\mathbb{Q}) \cong
      \mathbb{Q}, \quad QH_{n+2kN}(L;\mathbb{F}) \cong
      H_n(L;\mathbb{Q}) \cong
      \mathbb{Q}, \\
      & QH_{1+2kN}(L;\mathbb{F}) \cong H_1(L;\mathbb{Q})=0, \quad
      QH_{n-1+2kN}(L;\mathbb{F}) \cong H_{n-1}(L;\mathbb{Z})=0.
   \end{align*}
\end{proof}

\subsubsection*{Examples}
Here are a few examples of symplectic manifolds with semi-simple
quantum homology.
\begin{enumerate}
  \item ${\mathbb{C}}P^n$.
  \item The smooth complex quadric $Q$ (defined in \S\ref{Sb:quadric}
   above). See~\cite{Abrams:semi-simple} for the computation of the
   quantum Euler class.
  \item Complex Grassmannians. See~\cite{Abrams:semi-simple} for the
   proof.
  \item Semi-simplicity is preserved when taking products.  More
   precisely, let $(M_1, \omega_1)$, $(M_2, \omega_2)$ be two
   spherically monotone symplectic manifolds with the same
   proportionality factor between the $\omega_i|_{\pi_2}$'s and the
   $c_1(M_i)|_{\pi_2}$'s. Then $(M_1 \times M_2, \omega_1 \oplus
   \omega_2)$ is also spherically monotone. If $(M_1, \omega_1)$ and
   $(M_2, \omega_2)$ both have semi-simple quantum homology then so
   does $(M_1 \times M_2, \omega_1 \oplus \omega_2)$. This follows
   from the quantum K\"{u}nneth formula (see e.g.~\cite{McD-Sa:Jhol-2}).
   Of course, here one should consider semi-simplicity in the
   skew-commutative framework (see~\cite{Abrams:semi-simple}).
\end{enumerate}

Let $X \subset {\mathbb{C}}P^{n+r}$ be a smooth complete intersection
of degree $(d_1, \ldots, d_r)$ as in \S\ref{Sb:cinter}.  Let us now
examine the semi-simplicity of complete intersections in relation to
Lagrangian submanifolds. Assume that $2(\sum_{i=1}^r d_i -1) \leq n$
and that $n \geq 3$. Put $d = \prod_{i=1}^r d_i$. Consider the
following cases:

\paragraph{1. $d=1$} In this case $X \cong {\mathbb{C}}P^n$.
The quantum homology $QH_*(X;\mathbb{F})$ in this case is a field (see
\S\ref{S:cpn}), hence semi-simple.

\paragraph{2. $d=2$} In this case $X$ is isomorphic to a smooth quadric.
A direct computation of the quantum Euler class shows that it is an
invertible element (see~\cite{Abrams:semi-simple}) hence
$QH_*(X;\mathbb{F})$ is semi-simple.

\paragraph{3. $d=3$} In this case either $d_i \geq 3$ for some $i$
or there exist at least two $d_i$'s that are $\geq 2$. Therefore the
minimal Chern number $N$ of $X$ satisfies $N = n+1- \sum_{i=1}^r
(d_i-1) \leq n-1$. The assumption that $2(\sum_{i=1}^r d_i -1) \leq n$
implies that we also have $\frac{n}{2}+1 \leq N$. It follows that $N
\nmid n$ and $N \nmid n+1$. Next, note that $X$ has a Lagrangian
sphere.  This follows from symplectic Picard-Lefschetz theory since
for $d \geq 2$, $X$ can be degenerated to a variety with isolated
singularities. By Theorem~\ref{T:semi-simple-shpere} we conclude that
$QH_*(X;\mathbb{F})$ is not semi-simple. Of course, this can be also
verified by computing the quantum Euler class
(see~\cite{Abrams:semi-simple}).

It is not clear to us how large is the class of symplectic manifolds
with semi-simple quantum homology. It seems that it is in fact a
rather restricted class of manifolds.

Finally, let us mention that Entov and
Polterovich~\cite{En-Po:rigid-subsets} have also found restrictions on
semi-simplicity of quantum homology related to Lagrangian
submanifolds. The methods they use are based on the theory of spectral
numbers and are quite different than ours.

\subsection{Gromov radius and Relative symplectic packing} \label{Sb:pack}

Let $(M^{2n}, \omega)$ be a $2n$-dimensional symplectic manifold and
$L \subset M$ a Lagrangian submanifold. Denote by $B(r) \subset
\mathbb{R}^{2n}$ the closed $2n$-dimensional Euclidean ball of radius
$r$ endowed with the standard symplectic structure
$\omega_{\textnormal{std}}$ of $\mathbb{R}^{2n}$. Denote by
$B_{\mathbb{R}}(r) \subset B(r)$ the ``real'' part of $B(r)$, i.e.
$B_{\mathbb{R}}(r) = B(r) \cap (\mathbb{R}^n \times 0)$. Note that
$B_{\mathbb{R}}(r)$ is Lagrangian in $B(r)$.

By a {\em relative symplectic embedding} $\varphi: (B(r),
B_{\mathbb{R}}(r)) \to (M,L)$ of a ball in $(M,L)$ we mean a
symplectic embedding $\varphi: B(r) \to (M, \omega)$ which satisfies
the following conditions:
\begin{enumerate}
  \item $\varphi(B_{\mathbb{R}}(r)) \subset L$.
  \item $\varphi(x) \notin L$ for every $x \in B(r) \setminus
   B_{\mathbb{R}}(r)$.
   \label{I:rel-pack-2}
\end{enumerate}
These two conditions can be rewritten as ``$\varphi^{-1}(L) =
B_{\mathbb{R}}(r)$''.  Condition~(\ref{I:rel-pack-2}) may look strange
at first sight. We will explain its role soon (see
Remark~\ref{R:rel-pack-2nd-condition}).

In analogy with the (absolute) Gromov radius, we define here the
Gromov radius of $L \subset M$ to be
$$Gr(L) = \sup \{ r \mid \exists \textnormal{ a relative symplectic
  embedding } (B(r), B_{\mathbb{R}}(r)) \to (M,L)\}.$$
We will
consider also symplectic embeddings of balls in the complement of $L$,
i.e. symplectic embeddings $\psi:B(r) \to (M\setminus L, \omega)$.  We
denote the Gromov radius of $(M \setminus L, \omega)$ by
$Gr(M\setminus L)$ i.e.
$$Gr(M\setminus L) = \sup \{ r \mid \exists \textnormal{ a symplectic
  embedding } B(r) \to (M \setminus L) \}.$$

Denote by $\mathcal{J}(M,\omega)$ the space of almost complex
structure on $M$ which are compatible with $\omega$.
\begin{prop} \label{P:disk-pack-1} Let $L \subset (M, \omega)$ be a
   Lagrangian submanifold. Let $E\,', E'' > 0$.
   \begin{enumerate}
     \item Suppose that there exists a dense subset $\mathcal{J}'
      \subset \mathcal{J}(M,\omega)$ and a dense subset $\mathcal{U}
      \subset M$ such that for every $J \in \mathcal{J}'$ and every $p
      \in \mathcal{U}'$ there exists a non-constant $J$-holomorphic
      disk $u:(D, \partial D) \to (M,L)$ with $u(\textnormal{Int\,}D)
      \ni p$ and $\textnormal{Area}_{\omega}(u) \leq E\,'$. Then $$\pi
      Gr(M \setminus L)^2 \leq E\,'.$$ \label{I:pack-compl}
     \item Suppose that there exists a dense subset $\mathcal{J}''
      \subset \mathcal{J}(M,\omega)$ and a dense subset $\mathcal{U}''
      \subset L$ such that for every $J \in \mathcal{J}''$ and every
      $q \in \mathcal{U}''$ there exists a non-constant
      $J$-holomorphic disk $u:(D, \partial D) \to (M,L)$ with
      $u(\partial D) \ni q$ and $\textnormal{Area}_{\omega}(u) \leq
      E''$.  Then
      $$\frac{\pi Gr(L)^2}{2} \leq E''.$$ \label{I:pack-rel}
   \end{enumerate}
\end{prop}

\begin{proof}[Proof of Proposition~\ref{P:disk-pack-1}]
   The proof is based on an argument of Gromov from~\cite{Gr:phol}.
   Variants of the proof below can be found in~\cite{Bar-Cor:Serre,
     Cor-La:Cluster-1}.

   We begin with the proof of statement~\ref{I:pack-rel}. Let
   $\varphi:(B(r), B_{\mathbb{R}}(r)) \to (M,L)$ be a relative
   symplectic embedding. Let $J_0$ be the standard complex structure
   of $B(r)$. Let $J \in \mathcal{J}(M,\omega)$ be an almost complex
   structure which extends the complex structure $\varphi_*(J_0)$
   defined on the image of the ball $\varphi(B(r))$.

   Put $q=\varphi(0) \in L$. By Gromov compactness theorem there
   exists a non-constant $J$-holomorphic disk $u:(D, \partial D) \to
   (M, L)$ with $u(\partial D) \ni q$ and
   $\textnormal{Area}_{\omega}(u) \leq E''$. Put
   $S'=\varphi^{-1}(u(D))$ and let $S''=\overline{S'}$ be the complex
   conjugate copy of $S'$. Put $S=S' \cup S''$, $S^{\circ} = S \cap
   \textnormal{Int\,} B(r)$. Clearly $S^{\circ}$ is an analytic
   subvariety of $(\textnormal{Int\,}B(r), J_0)$. Note also that due
   to condition~(\ref{I:rel-pack-2}) (in the definition of relative
   symplectic embedding) the subvariety $S^{\circ} \subset
   \textnormal{Int\,}B(r)$ is properly embedded (with respect to the
   induced topology from $\textnormal{Int\,}B(r)$). By the Lelong
   inequality~\cite{Gr-Ha:alg-geom} we have $$1 \leq
   \textnormal{mult}_0 S^{\circ} \leq
   \frac{\textnormal{Area}_{\omega_{\textnormal{std}}}(S^{\circ})}{\pi
     r^2} =
   \frac{2\textnormal{Area}_{\omega_{\textnormal{std}}}(S')}{\pi r^2}
   \leq \frac{2E\,'}{\pi r^2}.$$
   This proves
   statement~\ref{I:pack-rel}.

   The proof of Statement~\ref{I:pack-compl} is very similar (but now
   no reflection argument is needed).
\end{proof}

\begin{rem} \label{R:rel-pack-2nd-condition} Let us explain the role
   of the condition~(\ref{I:rel-pack-2}) in the definition of relative
   symplectic embeddings. If this condition is dropped and we only
   require that $\varphi^{-1}(L) \supset B_{\mathbb{R}}(r)$ then the
   problem becomes equivalent to the problem of absolute symplectic
   embeddings. The point is that given any symplectic embedding
   $\varphi:B(r) \to M$, there exists another symplectic embedding
   $\varphi': B(r) \to M$ with $\varphi'(B_{\mathbb{R}}(r)) \subset
   L$. Indeed by a straightforward argument there exists a Hamiltonian
   diffeomorphism $h:(M,\omega) \to (M, \omega)$ such that $h \circ
   \varphi (B_{\mathbb{R}}(r)) \subset L$. Put $\varphi' = h \circ
   \varphi$.

   Thus if we drop condition~(\ref{I:rel-pack-2}) in the definition of
   relative symplectic embedding we lose the effect of the presence of
   the Lagrangian submanifold $L$. This can be easily illustrated
   already in dimension $2n=2$. Let $L \subset \mathbb{R}^2$ be a
   circle of radius $r$. Clearly $Gr(L)=\pi r^2$.  However, it is easy
   to see that for every $R>0$ there exists a symplectic embedding
   $\varphi:B^2(R) \to \mathbb{R}^2$ with $\varphi(B_\mathbb{R}(r))
   \subset L$.
\end{rem}

Before we continue would like to make a general remark concerning the
Gromov radius of Lagrangians. Given a monotone Lagrangian $L \subset
(M,\omega)$ we denote by $\tau$ the monotonicity constant $\tau =
\tfrac{\omega}{\mu}\big|_{\pi_2(M,L)}$ (see
formula~\eqref{eq:monotonicity}).

\begin{rem} \label{R:pack-general} Assume that the monotone Lagrangian
   $L^n \subset (M^{2n},\omega)$ has the property that
   $QH_{\ast}(L)=0$. We claim that this implies that for any almost
   complex structure $J$, through each point in $L$ passes a
   non-trivial $J$-holomorphic of area at most $(n+1)\tau$. This
   argument appears in \cite{Cor-La:Cluster-1} (where it is described
   mainly in the more delicate cluster setup) and goes as follows. Let
   $f$ be a Morse function on $L$ with a single maximum, $x_n$. We
   know from Remark \ref{rem:diff_m} that $dx_n=0$ in
   $\mathcal{C}(L;f,J)$. Thus $x_n$ is a boundary in this complex. But
   for this to happen some moduli space
   $\mathcal{P}(y,x_n;\mathbf{A},J)$ has to be non-trivial.  To have
   $x_n t^{\mu(\mathbf{A})/N_L}=dx+\ldots$ we need that
   $|x|-1=n-\mu(\mathbf{A})$ which means $\mu(\mathbf{A}) \leq n+1$.
   Thus, there is a $J$-holomorphic disk passing through $x_n$ of
   Maslov class at most $n+1$ and of area at most $(n+1)\tau$. As
   $x_n$ can be chosen as generic point of $L$ this shows the claim.

   In case $L$ is displaceable a variant of this argument that also
   takes into account the action filtration on the Floer complex shows
   that the area of these disks can be bounded by the displacement
   energy, $E(L)$.

   The existence of these disks implies by
   Proposition~\ref{P:disk-pack-1} that, in general we have
   $$\pi Gr(L)^{2}/2\leq (n+1)\tau$$
   and, in the displaceable case, we
   also have
   $$\pi Gr(L)^{2}/2\leq E(L)~.~$$
   The last inequality (in the
   displaceable case) can be proven using the cluster machinery even
   if $L$ is not monotone but is relative spin, orientable and
   $H_{\ast}(L;\Q)=0$ for all even $\ast$ different from $0$ and
   $dim(L)$ - see again \cite{Cor-La:Cluster-1} for this and more
   details on this argument.
\end{rem}

The following Corollary bounds the Gromov radius for Lagrangian tori.
Recall from Proposition~\ref{P:criterion-QH=0-1} that for a monotone
Lagrangian torus $T \subset (M, \omega)$ we either have $QH_*(T) =
(H(T) \otimes \Lambda)_*$ or $QH_*(T)=0$. In the latter case
Proposition~\ref{P:criterion-QH=0-1} implies that for generic $J$ and
a generic point $x \in T$ there exists a $J$-holomorphic disk $u:(D,
\partial D) \to (M,T)$ with $u(\partial D) \ni x$ and $\mu([u])=2$.
Combining this with Proposition~\ref{P:disk-pack-1} we the following.
\begin{cor}\label{C:disks_mon_torus}
   If $T \subset (M,\omega)$ is a monotone Lagrangian torus with
   monotonicity constant $\tau$ and $QH_{\ast}(T)\not=
   H_{\ast}(T)\otimes \La$, then
   $$\frac{\pi Gr(T)^{2}}{2}\leq 2\tau.$$
\end{cor}

Endow ${\mathbb{C}}P^n$ with the standard K\"{a}hler symplectic
structure $\omega_{\textnormal{FS}}$, normalized so that
$\int_{\mathbb{C}P^1} \omega_{\textnormal{FS}}=\pi$. Note that with
this normalization $({\mathbb{C}}P^n \setminus {\mathbb{C}}P^{n-1},
\omega_{\textnormal{FS}}) \cong (\textnormal{Int\,} B(1),
\omega_{\textnormal{std}})$, hence the (absolute) Gromov radius of
${\mathbb{C}}P^n$ is $Gr({\mathbb{C}}P^n) = 1$.
\begin{cor}\label{C:disks_compl}
   Let $L \subset {\mathbb{C}}P^n$ be a monotone Lagrangian with
   minimal Maslov number $N_L$.
   \begin{enumerate}
     \item If $QH_{\ast}(L)\not=0$, then we have:
      $$Gr({\mathbb{C}}P^n \setminus L)^2 \leq
      \frac{[\frac{2n}{N_L}]N_L}{2(n+1)}.$$
      In particular,
      $$Gr({\mathbb{C}}P^n \setminus L)^2 \leq \frac{n}{n+1}.$$
      \label{I:disks_compl-1}
     \item Suppose that $QH_{\ast}(L)\cong H_{\ast}(L)\otimes \La$,
      then:
      $$\frac{Gr(L)^{2}}{2}+ Gr({\mathbb{C}}P^n \setminus L)^2 \leq
      1.$$
      \label{I:disks_compl-2}
   \end{enumerate}
\end{cor}
\begin{rem} \label{R:QH=H}
   \begin{enumerate}
     \item In most of the examples below we will
      be in a special situation in which the pearl and Morse differentials agree, $d=\partial_{0}$,
      which greatly simplifies the proof of the Corollary ( $\partial_{0}$ is the Morse differential).
      More precisely,
      suppose that $QH_*(L) \cong (H(L;\mathbb{Z}_2) \otimes
      \Lambda)_*$ and that $L$ admits a perfect Morse function $f: L
      \to \mathbb{R}$.  In this case the pearl complex differential
      $d$ satisfies $d=\partial_0=0$.
      Indeed, since $\partial_0=0$ we
      have $\mathcal{C}_*(L;f,J) \cong (H(L;\mathbb{Z}_2) \otimes
      \Lambda)_*$.  Therefore if $d \neq 0$ we would have $\dim
      QH_i(L) < \dim (H(L;\mathbb{Z}_2) \otimes \Lambda)_i$ for some
      $i$.  A contradiction.
     \item Note that the term $[\frac{2n}{N_L}]$ in the first
      inequality of Corollary~\ref{C:disks_compl} cannot be $0$. This
      is because for any monotone Lagrangian submanifold $L \subset
      {\mathbb{C}}P^n$ we have $N_L \leq n+1$. This can be easily
      proved by the techniques from~\cite{Bi:Nonintersections}.
   \end{enumerate}
\end{rem}

\begin{proof}[Proof of Corollary~\ref{C:disks_compl}]
   We will use here the module structure described in \S\ref{S:qm}. As
   in that section we fix a Morse function $f:L\to \R$ and a Morse
   function $g:\C P^{n}\to \R$.  We will fix the following notation
   for the coefficients of the module operation. For $a\in \Crit(g),
   x\in \Crit(f)$ we write
   \begin{equation} \label{Eq:a*x}
      a\ast x= \sum_{y}n(a,x,y)y +
      \sum_{y,\mathbf{A}} n(a,x,y;\mathbf{A},J)
      yt^{\mubar(\mathbf{A})}
   \end{equation}
   as in formula~\eqref{Eq:qm-*}.

   We assume that $f$ has a single maximum denoted by $x_n$. We also
   assume that $g$ is a perfect Morse function so that we may identify
   its critical points with the generators of $H_{\ast}(\C P^{n})$.
   In particular, the minimum of $g$, will be denoted by $p$ and is
   identified with $h^{\cap n}$ where $h=[{\mathbb{C}}P^{n-1}] \in
   H_{2n-2}({\mathbb{C}}P^n;\mathbb{Z}_2)$ is the homology class of
   the hyperplane. We recall from the quantum homology of $\C P^{n}$
   that $h * p =h^{*(n+1)}=ut^{2(n+1)/N_L}$ where $u \in
   H_{2n}({\mathbb{C}}P^n;\mathbb{Z}_2)$ is the fundamental class of
   $\C P^{n}$. By Remark \ref{rem:prod_topclass} we know that as
   $QH_{\ast}(L)\not=0$, $x_n$ can not be a boundary and hence its
   homology class $w=[x_n]$ does not vanish in $QH_{\ast}(L)$.  Put
   $\alpha=p \ast w \in QH_{-n}(L)$. Notice that
   \begin{equation}\label{eq:cap_proj}
      h\ast \alpha = h \ast (p \ast w) =
      h^{\ast (n+1)} \ast w = (u \ast w)t^{2(n+1)/N_L}
      = wt^{2(n+1)/N_L}.
   \end{equation}
   We deduce that $\alpha \not = 0$. For dimension reasons the
   classical terms $n(p,x_n,-)$ in formula~\eqref{Eq:a*x} for $p*x_n$
   vanish. Therefore we can write
   $$\alpha=p \ast w = \Bigl[ \sum_{z,\mathbf{A}} n(p,x_n,z;
   \mathbf{A},J)zt^{\mubar(\mathbf{A})} \Bigr]$$
   where $z\in
   \textnormal{Crit}(f)$ and $[-]$ indicates taking homology classes.
   We have $-n=|z t^{\mubar(\mathbf{A})}| = |z|-\mu(\mathbf{A})$,
   $0\leq |z|\leq n$ so that if $n(p,x_n,z;\mathbf{A},J)\not=0$, as
   $p$ is the minimum of $g$, we deduce that there exists a
   $J$-holomorphic disk of Maslov index at most $2n$ which passes
   through $p$ and with boundary on $L$. As the Maslov index comes in
   multiples of $N_L$ the Maslov index of that disk is in fact $\leq
   [\frac{2n}{N_L}]N_L$. As $p$ may be chosen generically and the
   monotonicity constant $\tau$ is here $\pi/(2n+2)$ we deduce the
   statement at~(\ref{I:disks_compl-1}).

   For statement~(\ref{I:disks_compl-2}) we first pursue with the
   analysis of equation~\eqref{eq:cap_proj} under the additional asumption that the
   function $f$ is so that the differential of the pearl complex coincides with the Morse differential,
   $d=\partial_{0}$.  As mentioned in Remark \ref{R:QH=H} this condition is often satisfied
   and the proof in this case is relatively straightforward.

   We write:
   \begin{align*}
      w \ t^{2(n+1)/N_L} = h \ast \alpha = & \Bigl[
      \sum_{z,\mathbf{A}} n(p,x_n,z;\mathbf{A},J) (h \ast z)
      t^{\mubar(\mathbf{A})} \Bigr] \\
      = & \Bigl[ \sum_{z,\mathbf{A},y} n(p,x_n,z;\mathbf{A},J)
      n(h,z,y)y t^{\mubar(\mathbf{A})}  \\
      & + \sum_{z,\mathbf{A},z',\mathbf{A'}} n(p,x_n,z;\mathbf{A},J)
      n(h,z,z';\mathbf{A'},J)z'
      t^{\mubar(\mathbf{A})+\mubar(\mathbf{A'})} \Bigr],
   \end{align*}
   This means that in the chain complex $\mathcal{C}(L;f,J)$ we have:
   \begin{equation} \label{Eq:xn-sum}
      \begin{aligned}
         x_n & - \sum_{z,\mathbf{A},y} n(p,x_n,z;\mathbf{A},J)
         n(h,z,y)y t^{(\mu(\mathbf{A})-2(n+1))/N_L} \\
         & - \sum_{z,\mathbf{A},z',\mathbf{A'}}
         n(p,x_n,z;\mathbf{A},J) n(h,z,z';\mathbf{A},J)z'
         t^{(\mu(\mathbf{A})+\mu(\mathbf{A'})-2(n+1))/N_L} \in
         \textnormal{image\,}(d).
      \end{aligned}
   \end{equation}
   Since $d=\partial_0$
   and since $x_n$ is not homologous to any element in
   $$\bigoplus_{j\neq 0} \mathbb{Z}_2 x_n t^j \bigoplus \,
   (\mathbb{Z}_2 \langle \textnormal{Crit}(f) \setminus \{x_n\}
   \rangle \otimes \Lambda)$$
   it follows that $x_n$ must appear in one
   of the two sums in formula~\eqref{Eq:xn-sum}. But the $y$'s in the
   first sum of~\eqref{Eq:xn-sum} satisfy $|y| = |z| - 2 \leq n-2$,
   hence $y \neq x_n$. Therefore $x_n$ must appear as one of the
   summands in the second sum of~\eqref{Eq:xn-sum}. In other words
   there exist $z$, and vectors $\mathbf{A}, \mathbf{A'}$ of non-zero
   classes in $H_2(M,L)$ such that:
   \begin{enumerate}
     \item $n(p,x_n,z;\mathbf{A},J) \neq 0$, $n(h,z,x_n;\mathbf{A'},J)
      \neq 0$.
     \item $\mu(\mathbf{A})+ \mu(\mathbf{A'})=2(n+1)$.
   \end{enumerate}
   This means that through $p$ (the minimum of $h$) passes a
   non-constant $J$-holomorphic disk of Maslov index at most $\mu(\mathbf{A})$
   and through $x_n$ passes the boundary of non-constant
   $J$-holomorphic disk of Maslov index at most $\mu(\mathbf{A'}) = 2(n+1) -
   \mu(\mathbf{A})$. By Proposition~\ref{P:disk-pack-1} we have:
   $$\frac{\pi Gr(L)^{2}}{2}+ \pi Gr(\C P^{n},L)^{2} \leq \tau
   \mu(\mathbf{A}) + \tau \mu(\mathbf{A'}) = \tau(2n+2) = \pi.$$

   To conclude the proof of the point (2) of the proposition we now need to
   describe an argument in the absence of the condition $d=\partial_{0}$. To do
   so we will make use of the minimal model technology as described in
   \S\ref{subsec:minimal}. Thus, recall that there is a complex $(\mathcal{C}_{min},\delta)$, $\mathcal{C}_{min}=H(L;\Z_{2})\otimes \La$, and
   chain morphisms $\phi:\mathcal{C}(L;f,J)\to \mathcal{C}_{min}$, $\psi :\mathcal{C}_{min}\to \mathcal{C}(L;f,J)$ so that $\delta_{0}=0$ and $\phi\circ\psi=id$, $\phi$, $\psi$ ( respectively, $\phi_{0}$, $\psi_{0}$)
   induce isomorphisms in pearl (respectively, Morse) homology. As in Remark \ref{rem:product_min}
   we use the applications $\phi$ and $\psi$ to transport the module structure on $\mathcal{C}_{min}$.
   For $u\in \Crit(h)$ and $x\in H_{\ast}(L;\Z_{2})$ this module structure
   has the form:
   \begin{equation}\label{eq:min_module_pro}
   u\ast x=\sum_{k,y\in H_{\ast}(L;\Z_{2})} n'(u,x,y;k)yt^{k}\end{equation}
   so that
   $u\ast x=\phi (u\ast \psi (x))$. It is important to note that, as the minimal pearl model
   is constructed (in \S\ref{subsec:proof_minimal}) with coefficients in $\La^{+}$ we have that,
   for degree reasons, $\psi([L])=x_{n}$, $\phi(x_{n})=[L]$ (where $[L]\in H_{n}(L;\Z_{2})$ is the fundamental class).
   It follows that all the argument above can now be applied to $\mathcal{C}_{min}$ instead of
   $\mathcal{C}(L;f,J)$.  It leads to the fact that there is $z\in H_{\ast}(L;\Z_{2})$ for which there
   are coefficients $n'(-,-,-)$ which verify:
   \begin{enumerate}
     \item[(1')]$n'(p,[L],z;k) \neq 0$, $n'(h,z,[L];k')
      \neq 0$.
     \item[(2')] $(k+k')N_{L}=2(n+1)$.
   \end{enumerate}
  The condition $QH(L)\cong H(L;\Z_{2})\otimes \La$ means that $\delta=0$ in
  $\mathcal{C}_{min}$ (see Remark \ref{rem:product_min}). With this in mind
we take another look at equation (\ref{eq:min_module_pro}) and taking into account that
 both $\phi$ and $\psi$ are defined over $\La^{+}$ we deduce from $n'(p,[L],z;k)\not=0$
 that there exists $u\in \Crit(f)$ and $\mathbf{A}$ so that $\mu(\mathbf{A})\leq kN_{L}$ and $n(p,x_{n},u;\mathbf{A},J)\not=0$ and so through the minimum of $h$ passes a disk of Maslov class at most $kN_{L}$.
 We now want to interpret the equation $n'(h,z,[L];k')\not=0$. This obviously means that in $H_{\ast}(\mathcal{C}_{min})$ we have $h\ast z=\ldots+[L]t^{k'}+\ldots$.
  In other words, if we write $H_{\ast}(\mathcal{C}_{min})=
 H(L;\Z_{2})\otimes \La= E\oplus \Z_{2}<[L]t^{k'}>$ where $E$ is some complement
 of the $\Z_{2}$ vector space $\Z_{2}<[L]t^{k'}>$ we have $h\ast z=\xi+[L]t^{k'}$ with $\xi\in E$.
 This implies, in homology,
 $h\ast\psi(z)=\psi(\xi)+ [x_{n}]t^{k'}$ so that at the chain level
 we have $h\ast \psi(z)+dw=x_{n}t^{k'}+\psi(\xi)$.
 We also split $\mathcal{C}(L;f,J)=E'\oplus \Z_{2}<x_{n}t^{k'}>$
 and we consider two cases. First suppose that there is some cycle $\xi'\in E'$ so that $[\xi']=[x_{n}t^{k'}]$.
 The only possibility for that to happen is that there is some $u''\in\Crit(f)$, $\mathbf{A}''$
 with $\mu(\mathbf{A}'')\leq k' N_{L}$ so that $n(u'',x_{n};\mathbf{A}'',J)\not=0$  where $n(-,-)$ are the coefficients of the differential in the pearl complex $\mathcal{C}(L;f,J)$. In this case, it follows that through the maximum of $f$ passes a disk of Maslov class at most $k'N_{L}$ and the proof ends as in the special
 case treated before.  Assume therefore from now on that for all the cycles $\xi'$ in $E'$ we have
 $[\xi']\not=[x_{n}t^{k'}]$. In particular, this means that $x_{n}t^{s}$ does not appear
 in the boundary expression of any critical point of $f$ for $s\leq k'$.  At this point we need
  to recall the explicit construction of the map $\phi$: it is defined on a basis
  of $\Z_{2}<\Crit(f)>$ which is chosen so that it is formed by three types of elements $x$,'s with
  $\partial_{0} x=0$, $y$'s with $\partial_{0} y=0$ and $y'$'s with $\partial_{0} y'=y$  and is given by
  $\phi(x)=x$, $\phi(y')=0$, $\phi{y}=\phi^{n-s+1}(dy'-y)$ (the last equality appears
  in the inductive step - see again \S\ref{subsec:proof_minimal}). The critical point $x_{n}$ is
  itself a generator of type $x$ and it is the only one in dimension $n$.
  Under our assumption, this means that we have $\phi(E')=E$.  Indeed, the only thing to check is that
  each generator of type $y$ is sent to $E$ but this follows inductively because  $dy'$ does
  not contain $x_{n}t^{s}$ for $s\leq k'$.

  We now write $\psi(\xi)=\xi''+\epsilon x_{n}t^{k'}$
 where $\epsilon$ is $0$ or $1$ and $\xi''\in E'$.
  In case $\epsilon =1$ we have, in homology, $h\ast \psi(z)=[\xi'']$ and, so
  $\phi(h\ast \psi (z))=h\ast z=\phi [\xi'']=\xi+[L]t^{k}$ which contradicts $\phi(\xi'')\in E$
  (we use here again $\delta=0$).

  Thus we are left to discuss the case when $\epsilon=0$. Therefore we have
  $h\ast \psi(z)+dw=x_{n}t^{k'}+\xi''$ and given our assumption this means that
 there exists $u'\in\Crit(f)$, $\mathbf{A}'$ with $\mu(\mathbf{A}')\leq k'N_{L}$ so that
 $n(h,u',x_{n};\mathbf{A}',J)\not=0$. It again follows that through the maximum of $f$ passes a disk
 of Maslov class at most $k'N_{L}$ and concludes the proof.
\end{proof}

\subsection{Examples} \label{Sb:ex-gr-radius}

Endow ${\mathbb{C}}P^n$ with the standard K\"{a}hler symplectic
structure $\omega_{\textnormal{FS}}$, normalized so that
$\int_{\mathbb{C}P^1} \omega_{\textnormal{FS}}=\pi$.
\begin{cor} \label{C:pack-2H_1=0} Let $L \subset {\mathbb{C}}P^n$ be a
   Lagrangian submanifold with $2H_1(L;\mathbb{Z})=0$. Then
   $Gr({\mathbb{C}}P^n \setminus L)^2 \leq \frac{1}{2}$.
\end{cor}
\begin{proof}
   This follows immediately from Proposition~\ref{P:enum-mixed-rpn} in
   conjunction with Proposition~\ref{P:disk-pack-1}.

   Alternatively, one can use Proposition~\ref{T:L-rpn} by which $N_L =
   n+1$ and $QH_*(L) \neq 0$. Then Corollary~\ref{C:disks_compl}
   implies that: $$Gr({\mathbb{C}}P^n \setminus L)^2 \leq
   \frac{[\frac{2n}{n+1}](n+1)}{2(n+1)} = \frac{1}{2}.$$
\end{proof}
When $L$ is the standard real projective space $\mathbb{R}P^n \subset
{\mathbb{C}}P^n$ the inequality in Corollary~\ref{C:pack-2H_1=0} has
been proved before by Biran~\cite{Bi:Barriers} by different methods.
In this case the inequality turns out to be sharp. Note also that (by
an explicit construction) we have $Gr(\mathbb{R}P^n)=1$, thus the
inequality in Corollary~\ref{C:disks_compl}-(\ref{I:disks_compl-2}) is
sharp in this case.

Consider the $n$-dimensional Clifford torus
$\mathbb{T}^n_{\textnormal{clif}} \subset {\mathbb{C}}P^n$. (See
Example~\ref{Sb:tclif-QH-1}.)
\begin{cor} \label{C:clif-gr}
   $Gr(\mathbb{T}^n_{\textnormal{clif}})^2 \leq \frac{2}{n+1}, \quad
   Gr({\mathbb{C}}P^n \setminus \mathbb{T}^n_{\textnormal{clif}})^2 =
   \frac{n}{n+1}.$
\end{cor}
\begin{remnonum}
   By an explicit packing construction, communicated to us by
   Buhovsky, it seems that the first inequality in
   Corollary~\ref{C:clif-gr} is in fact sharp, i.e.
   $Gr(\mathbb{T}^n_{\textnormal{clif}})^2 = \frac{2}{n+1}$.
\end{remnonum}
\begin{proof}[Proof of Corollary~\ref{C:clif-gr}]
   A straightforward construction using moment maps (in the spirit
   of~\cite{Karshon:pack, Tr:pack-constr}) shows that
   $Gr({\mathbb{C}}P^n \setminus \mathbb{T}_{\textnormal{clif}}) \geq
   \frac{n}{n+1}$. On the other hand, by Example~\ref{Sb:tclif-QH-1}
   we have $QH_*(\mathbb{T}^n_{\textnormal{clif}}) \cong
   (H(\mathbb{T}^n_{\textnormal{clif}};\mathbb{Z}_2) \otimes
   \Lambda)_*$ hence by
   Corollary~\ref{C:disks_compl}-(\ref{I:disks_compl-1}) we have
   $Gr({\mathbb{C}}P^n \setminus \mathbb{T}^n_{\textnormal{clif}})
   \leq \frac{n}{n+1}$. This proves that $Gr({\mathbb{C}}P^n \setminus
   \mathbb{T}^n_{\textnormal{clif}}) = \frac{n}{n+1}$.

   The first inequality follows now from the fact that
   $QH_*(\mathbb{T}^n_{\textnormal{clif}}) \cong
   (H(\mathbb{T}^n_{\textnormal{clif}};\mathbb{Z}_2) \otimes
   \Lambda)_*$ in conjunction with
   Corollary~\ref{C:disks_compl}-(\ref{I:disks_compl-2}) and
   Remark~\ref{R:QH=H}. An alternative, more direct, proof goes as
   follows. By Example~\ref{Sb:tclif-QH-1} it follows that for generic
   $J \in \mathcal{J}({\mathbb{C}}P^n)$ and $x \in
   \mathbb{T}^n_{\textnormal{clif}}$ there exists a $J$-holomorphic
   disk $u:(D, \partial D) \to ({\mathbb{C}}P^n, L)$ with
   $\mu([u])=2$. This disk has area $\frac{\pi}{n+1}$. By
   Proposition~\ref{P:disk-pack-1}-(\ref{I:pack-rel}) we obtain
   $Gr(\mathbb{T}^n_{\textnormal{clif}}) \leq \frac{1}{2(n+1)}$.
\end{proof}

\subsection{Mixed symplectic packing} \label{Sb:mixed-pack}
Let $l, m \geq 0$ and $r_1, \ldots, r_l>0$, $\rho_1, \ldots, \rho_m
>0$.  A mixed symplectic packing of $(M,L)$ by balls of radii $(r_1,
\ldots, r_l; \rho_1, \ldots, \rho_m)$ is given by $l$ relative
symplectic embeddings $\varphi_i: (B(r_i), B_{\mathbb{R}}) \to (M,L)$,
$i=1, \ldots, l$, and $m$ symplectic embeddings $\psi_j:B(r_j) \to (M
\setminus L, \omega)$, $j=1, \ldots, m$, such that the images of all
the $\varphi_i$ and $\psi_j$ are mutually disjoint, i.e.
\begin{enumerate}
  \item $\varphi_{i'} (B(r_{i'})) \cap \varphi_{i''} (B(r_{i''})) =
   \emptyset$ for every $i' \neq i''$.
  \item $\psi_{j'} (B(r_{j'})) \cap \psi_{j''} (B(r_{j''})) =
   \emptyset$ for every $j' \neq j''$.
  \item $\varphi_i(B(r_i)) \cap \psi_j(B(r_j)) = \emptyset$ for every
   $i,j$.
\end{enumerate}

The following can be proved in a similar way to
Proposition~\ref{P:disk-pack-1}.
\begin{prop} \label{P:disk-pack-2}
   Let $L \subset (M, \omega)$ be a Lagrangian submanifold and $E>0$.
   Suppose that there exists a dense subset $\mathcal{J}_* \subset
   \mathcal{J}(M, \omega)$, a dense subset of $m$-tuples $\mathcal{U}'
   \subset (M \setminus L)^{\times m}$, and a dense subset of
   $l$-tuples $\mathcal{U}'' \subset L^{\times l}$ such that for every
   $J \in \mathcal{J}_*$, $(q_1, \ldots, q_l) \in \mathcal{U}'$,
   $(p_1, \ldots, p_m) \in \mathcal{U}''$ there exists a non-constant
   $J$-holomorphic disk $u:(D, \partial D) \to (M, L)$ with
   $u(\textnormal{Int\,}D) \ni p_1, \ldots, p_m$, $u(\partial D) \ni
   q_1, \ldots, q_l$ and $\textnormal{Area}_{\omega}(u) \leq E$. Then
   for every mixed symplectic packing of $(M,L)$ by balls of radii
   $(r_1, \ldots, r_l; \rho_1, \ldots, \rho_m)$ we have:
   $$\sum_{i=1}^l \frac{\pi r_i^2}{2} + \sum_{j=1}^m \pi \rho_j^2 \leq
   E.$$
\end{prop}

\subsection{Examples} \label{Sb:ex-mixed-pack}

\begin{cor} \label{C:mixed-pack-tclif}
   For every mixed symplectic packing of $({\mathbb{C}}P^2,
   \mathbb{T}^2_{\textnormal{clif}})$ by two balls of radii $(r;\rho)$
   we have $\frac{1}{2}r^2 + \rho^2 \leq \frac{2}{3}$. In particular
   if the two balls are assumed to have the same radius $r=\rho$ then
   $r^2 \leq \frac{4}{9}$.
\end{cor}
Note that the inequality in Corollary~\ref{C:mixed-pack-tclif} is
stricter than the inequalities for (absolute) packing of
${\mathbb{C}}P^2$ by $2$ balls. Indeed, for every symplectic packing
of ${\mathbb{C}}P^2$, $B(r) \coprod B(\rho) \hookrightarrow
{\mathbb{C}}P^2$ the (optimal) packing inequality reads $r^2 + \rho^2
\leq 1$ (See~\cite{Gr:phol, MP:Pack}). In particular, if the balls are
of equal radii $r=\rho$ the latter inequality reads $r^2\leq
\frac{1}{2}$ (and this is sharp), whereas the mixed packing inequality
gives $r^2 \leq \frac{4}{9}$.

We do not know whether the inequality in
Corollary~\ref{C:mixed-pack-tclif} is sharp.

\begin{proof}[Proof of Corollary~\ref{C:mixed-pack-tclif}]
   We will use here the notation of \S\ref{Ss:clifford} and
   Proposition~\ref{prop:clifford}. Fix two Morse functions
   $g:{\mathbb{C}}P^2 \to \mathbb{R}$,
   $f:\mathbb{T}^2_{\textnormal{clif}} \to \mathbb{R}$ and two generic
   Riemannian metrics $\nu$, $\nu_{_{\mathbb{T}}}$ on
   ${\mathbb{C}}P^2$ and $\mathbb{T}^2_{\textnormal{clif}}$. We assume
   that both $g$ and $f$ are perfect Morse functions so that we can
   identify their critical points with homology classes in
   $H_*({\mathbb{C}}P^2;\mathbb{Z}_2)$ and
   $H_*(\mathbb{T}^2_{\textnormal{clif}};\mathbb{Z}_2)$. Denote by $h
   \in H_2({\mathbb{C}}P^2;\mathbb{Z}_2)$ the class of the hyperplane
   and by $p \in H_0({\mathbb{C}}P^2; \mathbb{Z}_2)$ the class of a
   point. Denote by $m$ the minimum of $f$. For every generic $J \in
   \mathcal{J}({\mathbb{C}}P^2)$, $m$ defines an element, still
   denoted $m$, in $QH_0(\mathbb{T}^2_{\textnormal{clif}};
   f,\nu_{_{\mathbb{T}}},J)$. (See the discussion in the beginning of
   \S\ref{Sb:clif-enum}.) According to Proposition~\ref{prop:clifford}
   we have
   $$p*m = h*h*m = h*mt= m t^2.$$
   The pearly trajectories that
   potentially contribute to this computation appear in
   figure~\ref{f:mixed-pack-1}.
   \begin{figure}[htbp]
      \begin{center}
         \epsfig{file=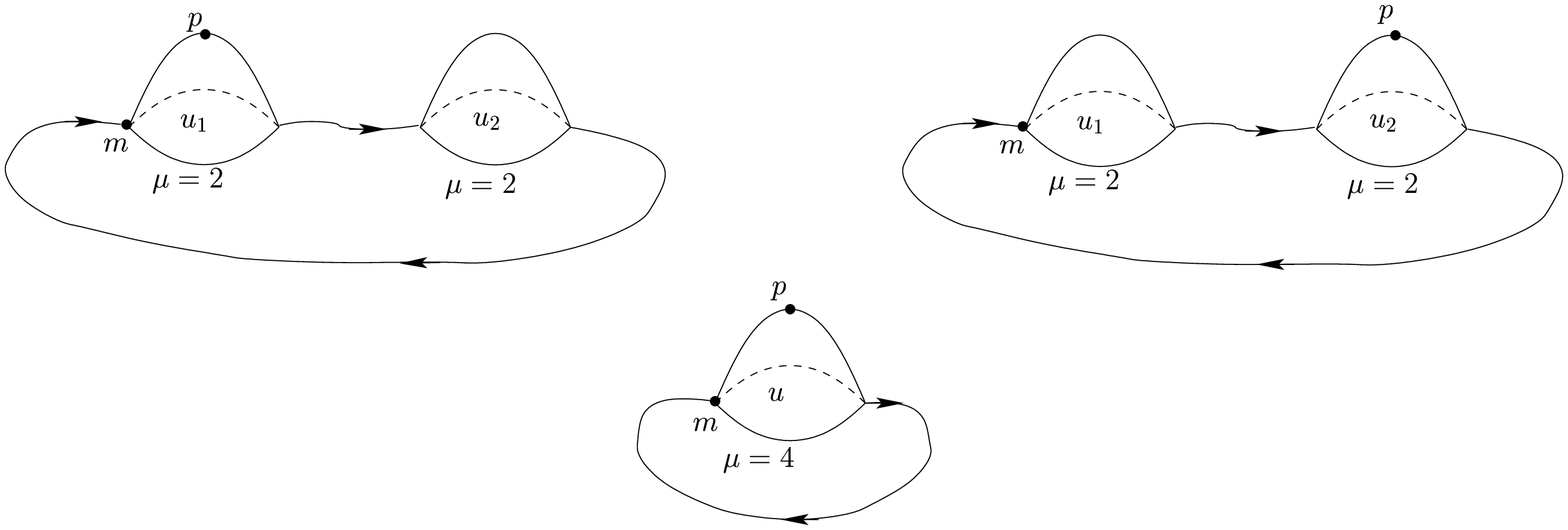, width=1.0\linewidth}
      \end{center}
      \caption{}
      \label{f:mixed-pack-1}
   \end{figure}

   It follows that for generic $J$ one of the following three
   possibilities occur:
   \begin{enumerate}
     \item There exists a $J$-holomorphic disk $u_1:(D, \partial D)
      \to ({\mathbb{C}}P^2, \mathbb{T}^2_{\textnormal{clif}})$ with
      $\mu([u_1])=2$ and $u_1(-1)=m$, $u_1(0)=p$.
     \item There exist two $J$-holomorphic disks $u_1, u_2: (D,
      \partial D) \to ({\mathbb{C}}P^2,
      \mathbb{T}^2_{\textnormal{clif}})$ with
      $\mu([u_1])=\mu([u_2])=2$ and $u_1(-1)=m$, $u_2(0)=p$.
     \item There exists a $J$-holomorphic disk $u:(D, \partial D) \to
      ({\mathbb{C}}P^2, \mathbb{T}^2_{\textnormal{clif}})$ with
      $\mu([u])=4$ and $u(-1)=m$, $u(0)=p$.
   \end{enumerate}
   Note that the in cases~(1),~(2) we have
   $\textnormal{Area}_{\omega_{\textnormal{FS}}}(u_1) =
   \textnormal{Area}_{\omega_{\textnormal{FS}}}(u_2) = 2\tau = 1/3$
   and in case~(3)
   $\textnormal{Area}_{\omega_{\textnormal{FS}}}(u)=4\tau = 2/3$.

   Let $\varphi: B(r) \to ({\mathbb{C}}P^2,
   \mathbb{T}^2_{\textnormal{clif}})$, $\psi: B(\rho) \to
   {\mathbb{C}}P^2 \setminus \mathbb{T}^2_{\textnormal{clif}}$ be a
   mixed symplectic packing of $({\mathbb{C}}P^2,
   \mathbb{T}^2_{\textnormal{clif}})$ by two balls of radii
   $(r;\rho)$. Take $m$ to be $\varphi(0)$ and $p=\psi(0)$. Arguing as
   in the proof of Propositions~\ref{P:disk-pack-1}
   and~\ref{P:disk-pack-2} we obtain the following inequality:
   $$\frac{1}{2}\pi r^2 + \pi \rho^2 \leq \max \{2\tau, 2\tau+2\tau,
   4\tau\} = 4 \tau = \frac{2\pi}{3}.$$
\end{proof}
The proof of Corollary~\ref{C:mixed-pack-tclif} suggests that for
mixed packing of $({\mathbb{C}}P^n, \mathbb{T}^n_{\textnormal{clif}})$
by two balls of radii $(r; \rho)$ the following packing inequality
should hold: $\frac{1}{2}r^2 + \rho^2 \leq \frac{n}{n+1}$.

\

Consider the smooth complex quadric $Q \subset {\mathbb{C}}P^{n+1}$
endowed with the symplectic structure $\omega$ induced from
${\mathbb{C}}P^{n+1}$ (See \S\ref{Sb:quadric}). With our normalization
of the symplectic structure on ${\mathbb{C}}P^{n+1}$ the symplectic
structure on $Q$ and the first Chern class $c_1$ (of the tangent
bundle of $Q$) have the following relation $c_1 = n[\omega]/\pi$.

\begin{cor} \label{C:mixed-pack-quad-1} Let $L \subset Q$ be a
   Lagrangian submanifold with $H_1(L;\mathbb{Z})=0$ (e.g. a
   Lagrangian sphere). Then for every mixed symplectic packing of
   $(Q,L)$ with $2$ balls of radii $(r;\rho)$ we have: $\frac{1}{2}r^2
   + \rho^2 \leq 1$.
\end{cor}
Note that inequality in Corollary~\ref{C:mixed-pack-quad-1} is
stricter than the absolute packing inequalities, at least in dimension
$2n=4$. Indeed in that case $(Q, \omega) \cong (\mathbb{C}P^1 \times
\mathbb{C}P^1, \omega_{\textnormal{FS}} \oplus
\omega_{\textnormal{FS}})$ and the (optimal) absolute packing
inequalities for two balls $B(r) \coprod B(\rho) \hookrightarrow
(\mathbb{C}P^1 \times \mathbb{C}P^1, \omega_{\textnormal{FS}} \oplus
\omega_{\textnormal{FS}})$ read $r^2 \leq 1$, $\rho^2 \leq 1$.

We do not know whether the inequality in
Corollary~\ref{C:mixed-pack-quad-1} is sharp.

\begin{proof}[Proof of Corollary~\ref{C:mixed-pack-quad-1}]
   This follows from Corollary~\ref{C:quad-disks} and
   Proposition~\ref{P:disk-pack-2}.
\end{proof}

\subsection{Further questions}
\label{Sb:further-questions-pack}

The results in \S\ref{Sb:pack}-~\ref{Sb:ex-mixed-pack} give rise to
several questions. The first one is whether the packing inequalities
above are sharp. This is especially relevant for the results in
Corollaries~\ref{C:disks_compl},~\ref{C:pack-2H_1=0},
~\ref{C:clif-gr},~\ref{C:mixed-pack-tclif},~\ref{C:mixed-pack-quad-1}.
This would require to obtain lower bounds on the radii of balls in
relative/mixed symplectic packing. Such bounds can sometimes be
obtained by an explicit packing construction, but one would like a
more systematic method. In the theory of absolute symplectic packing
an essential ingredient is the symplectic blowing up and down
constructions (see~\cite{MP:Pack}). It seems relevant to establish a
relative (with respect to a Lagrangian $L$) version of the symplectic
blowing up and down constructions. Another important ingredient in the
theory of absolute symplectic packing consists of criteria for
realizing $2$-dimensional cohomology classes by symplectic/K\"{a}hler
forms such as the Nakai-Moishezon criterion (see~\cite{Bi:Pack,
  Bi:ECM2000, Bi:Stbl-Pack} for symplectic analogues of this). In the
relative version one would expect to have analogous criteria with the
additional requirement that the resulting symplectic form makes a
given submanifold Lagrangian.

Another interesting question is what happens for larger number of
balls. This would require to establish existence of holomorphic disks
(or surfaces with boundary) passing through many points on $L$. It
seems likely that our methods combined with the $A_{\infty}$ approach
of Fukaya-Oh-Ohta-Ono \cite{FO3} or the cluster homology approach of
Cornea-Lalonde \cite{Cor-La:Cluster-1} would be relevant for this
purpose. In the same direction, it would be interesting to find out
whether the packing obstructions disappear in the relative case once
the number of balls becomes large enough, as happens in the absolute
case in dimension $4$ (see~\cite{Bi:Pack, Bi:Stbl-Pack}).

Some of the results on existence of holomorphic disks in
\S\ref{Sb:enum-rpn} and~\ref{Sb:quad-quant-L} give rise to redundant
packing inequalities in the sense that they coincide (or can be
derived from) the absolute packing inequalities. For example,
Proposition~\ref{P:enum-mixed-rpn} implies that for every relative
symplectic packing $$(B(r_1), B_{\mathbb{R}}(r_1)) \coprod (B(r_2),
B_{\mathbb{R}}(r_2)) \hookrightarrow ({\mathbb{C}}P^n,
\mathbb{R}P^n)$$
we have $\frac{1}{2}r_1^2 + \frac{1}{2}r_2^2 \leq
\frac{1}{2}$. However this is precisely the (optimal) absolute packing
inequality for two balls in ${\mathbb{C}}P^n$ (see~\cite{Gr:phol}).
In fact, an explicit construction as in~\cite{Karshon:pack} shows that
for $l\leq 3$ balls there is no difference between relative and
absolute symplectic packing (i.e. the same packing inequalities hold
for both cases). We do not know whether this continues to hold for
general $l$ (even in dimension $2n=4$). It would also be interesting
to find a geometric explanation to why the relative symplectic packing
problem coincides with the absolute one for some Lagrangians while for
others it gives stricter restrictions.

\subsection{Quantum product on tori and enumerative geometry}
\label{Sb:clif-enum}
The goal of this section is to give a geometric interpretation of the
quantum cap product for Lagrangian tori in terms of enumeration of
holomorphic disks. Related results in this direction have recently
been obtained for torus fibres of Fano toric manifolds by
Cho~\cite{Cho:counting-disks} by a different approach. For simplicity
we consider here only the case of $2$-dimensional tori.

Let $L^2 \subset (M^4, \omega)$ be a monotone Lagrangian torus with
minimal Maslov number $N_L=2$. Assume that $QH_*(L) \neq 0$. By
Proposition~\ref{P:criterion-QH=0-1} we have $QH_*(L) \cong
(H(L;\mathbb{Z}_2) \otimes \Lambda)_*$. It follows that
\begin{align}
   & QH_0(L) \cong H_0(L;\mathbb{Z}_2) \oplus H_2(L;\mathbb{Z}_2) t,
   \label{Eq:QH_0} \\
   & QH_1(L) \cong H_1(L;\mathbb{Z}_2).
   \label{Eq:QH_1}
\end{align}
We will see in a moment that the isomorphism~\eqref{Eq:QH_1} is
canonical. However, the splitting in~\eqref{Eq:QH_0} is {\em not
  canonical} in the sense that it is not compatible with the canonical
identifications in Morse homology. As we will see below, the
``$H_2(L;\mathbb{Z}_2) t$'' part is canonical but the inclusion in
$QH_0(L)$ of the summand $H_0(L;\mathbb{Z}_2)$ actually depends on the
Morse function $f:L \to \mathbb{R}$ and the almost complex structure
$J$ used to compute $QH(L)$. Let us explain this point in more detail.

\begin{figure}[htbp]
   \begin{center}
      \epsfig{file=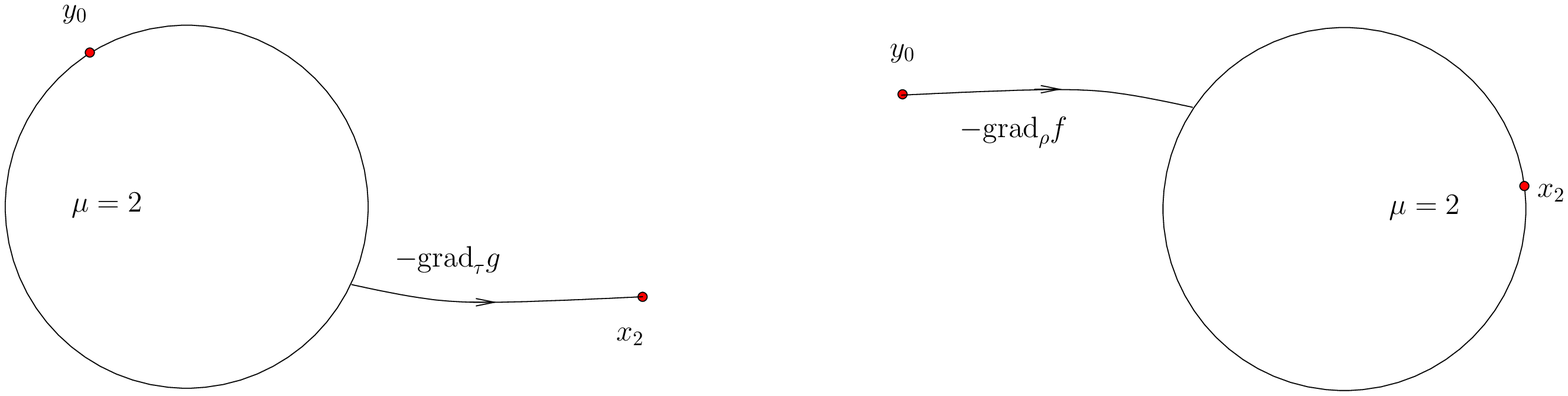, width=1.0\linewidth}
   \end{center}
   \caption{Pearly trajectories contributing to the coefficient $\epsilon$.}
   \label{f:QH-iso-1}
\end{figure}

Let $f,g:L \to \mathbb{R}$ be two {\em perfect} Morse functions. Fix
two generic Riemannian metrics $\rho, \tau$ on $L$ and a generic
almost complex structure $J \in \mathcal{J}(M,\omega)$. Denote:
\begin{align*}
   & \textnormal{Crit}_0(f)=\{x_0\}, \quad
   \textnormal{Crit}_1(f)=\{x'_1, x''_1\}, \quad
   \textnormal{Crit}_2(f) = x_2, \\
   & \textnormal{Crit}_0(g)=\{y_0\}, \quad
   \textnormal{Crit}_1(g)=\{y'_1, y''_1\}, \quad
   \textnormal{Crit}_2(g) = y_2.
\end{align*}
Since $QH_*(L) \cong (H(L;\mathbb{Z}_2) \otimes \Lambda)_*$ the Floer
differential vanishes and we have:
\begin{align*}
   & QH_0(f,\rho,J) = \mathbb{Z}_2 x_0 \oplus \mathbb{Z}_2 x_2 t,
   \quad
   QH_1(f,\rho, J) = \mathbb{Z}_2 x'_1 \oplus \mathbb{Z}_2 x''_1, \\
   & QH_0(g,\tau,J) = \mathbb{Z}_2 y_0 \oplus \mathbb{Z}_2 y_2 t,
   \quad QH_1(g,\tau, J) = \mathbb{Z}_2 y'_1 \oplus \mathbb{Z}_2
   y''_1.
\end{align*}
Denote by $\phi^M : H_*(g,\tau) \to H_*(f, \rho)$ the canonical
isomorphisms of Morse homologies. Denote by $\phi^F : QH_*(g,\tau,J)
\to QH_*(f, \rho, J)$ the canonical isomorphisms of Floer homologies
as described in \S\ref{subsubsec:inv}. Then for degree reasons
$\phi^F$ coincides with $\phi^M$ on $QH_1(g, \tau, J)$.  Clearly we
also have $\phi^M(y_0)=x_0$ and $\phi^M(y_2) = x_2$.  Moreover
$\phi^F(y_2) = x_2$ since $y_2, x_2$ are the unities of the respective
Floer homologies. However it might happen that $\phi^F(y_0) \neq
\phi^M(y_0)$. More precisely, following the description in
\S\ref{subsubsec:inv} we have $\phi^F(y_0) = x_0 + \epsilon x_2 t$
where the coefficient $\epsilon \in \mathbb{Z}_2$ is determined by
counting the number of pearly trajectories appearing in
figure~\ref{f:QH-iso-1}.

\begin{figure}[htbp]
   \begin{center}
      \epsfig{file=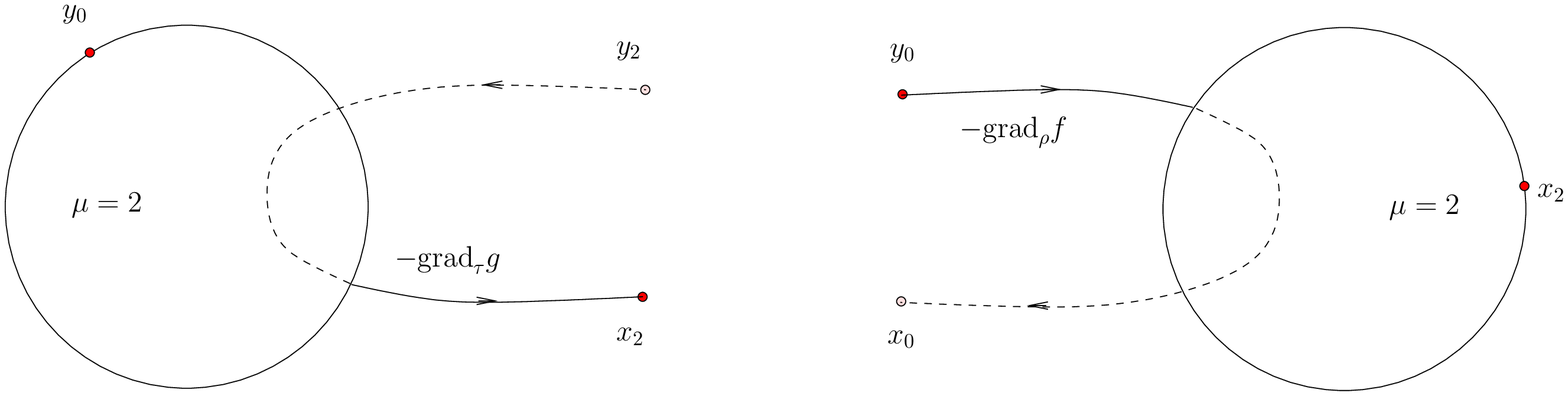, width=1.0\linewidth}
   \end{center}
   \caption{}
   \label{f:QH-iso-2}
\end{figure}

As $f,g$, $\rho, \tau$, $J$ are taken to be generic we may assume that
$y_0 \in W^s_{x_0}(-\textnormal{grad}_{\rho} f)$, $x_2 \in
W^u_{y_2}(-\textnormal{grad}_{\tau} g)$ and that no two points from
$x_0, y_0, x_2, y_2$ lie on the boundary of the same $J$-holomorphic
disk with Maslov number $2$. Denote by $\gamma_{y_2, x_2}(g)$ the
$-\textnormal{grad}_{\tau}(g)$ trajectory connecting $y_2$ to $x_2$
and by $\gamma_{y_0, x_0}(f)$ the $-\textnormal{grad}_{\rho}(f)$
trajectory connecting $y_0$ to $x_0$.  Using this and the notation
from~\eqref{Eq:delta-x} we have (see figure~\ref{f:QH-iso-2}):
\begin{equation} \label{Eq:coef-eps-1}
   \epsilon = \#_{\mathbb{Z}_2}
   \bigl(\delta_{y_0}(J) \cap \gamma_{y_2, x_2}(g)\bigr)
   +  \#_{\mathbb{Z}_2}
   \bigl(\delta_{x_2}(J) \cap \gamma_{y_0, x_0}(f)\bigr).
\end{equation}

By Proposition~\ref{P:criterion-QH=0-1} the cycles $\delta_{y_0}(J)$
and $\delta_{x_2}(J)$ are $0$ in the homology $H_1(L;\mathbb{Z}_2)$.
Therefore, the path $\gamma_{y_2, x_2}(g)$ in
formula~\eqref{Eq:coef-eps-1} can be replaced by any path $\ell(y_2,
x_2) \subset L$ connecting $y_2$ to $x_2$ (as long as the intersection
of that path with $\delta_{y_0}(J)$ is transverse).  Similarly
$\gamma_{y_0, x_0}(f)$ can be replaced by any path $\ell(y_0,x_0)$.

Summarizing the above we have:
\begin{prop} \label{P:QH-fg}
   The isomorphism $\phi^F: QH_*(g,\tau,J) \longrightarrow
   QH_*(f,\rho,J)$ is given by:
   $$\phi^F(y) = \phi^M(y) \; \forall y \in QH_1(g,\tau,J), \quad
   \phi^F(y_2) = x_2, \quad \phi^F(y_0) = x_0 + \epsilon x_2 t,$$
   where
   \begin{equation} \label{Eq:coef-eps-2}
   \epsilon = \#_{\mathbb{Z}_2} \bigl(\delta_{y_0}(J)
   \cap \ell(y_2, x_2)\bigr) + \#_{\mathbb{Z}_2} \bigl(\delta_{x_2}(J)
   \cap \ell(y_0, x_0)\bigr)
   \end{equation}
   and $\ell(y_2, x_2), \ell(y_0, x_0) \subset L$ are any two paths in
   $L$ joining $y_2$ to $x_2$, resp.  $y_0$ to $x_0$, and meeting the
   cycles $\delta_{y_0}(J)$, resp.  $\delta_{x_2}(J)$, transversely.
   In particular, the isomorphism $QH_1(L) \cong H_1(L;\mathbb{Z}_2)$
   as well as the inclusion $H_2(L;\mathbb{Z}_2) \to
   QH_2(L;\mathbb{Z}_2)$ are canonical.
\end{prop}
\begin{rem} \label{R:QH-fg}
   \begin{enumerate}
     \item From formula~\eqref{Eq:coef-eps-2} it follows that once we
      fix a generic $J$ and take $f,g$ with the critical points $x_2$
      close enough to $y_2$ and $x_0$ close enough to $y_0$ then the
      coefficient $\epsilon$ is $0$. Thus in this case
      $\phi^F=\phi^M$.
     \item One can derive a similar formula to~\eqref{Eq:coef-eps-2}
      when we have two different almost complex structures $J_0, J_1$.
      In this case one takes a generic path $\{J_t\}$ connecting $J_0$
      to $J_1$.  Then there is an additional contribution to the
      coefficient $\epsilon$ coming from $J_t$-holomorphic disks (for
      some $t$) with Maslov index $2$ whose boundaries pass through
      the points $y_0, x_2$.
   \end{enumerate}
\end{rem}

Let $L$ be as above. Denote by $w \in QH_2(L)=H_2(L;\mathbb{Z}_2)$ the
fundamental class. Choose $m \in QH_0(L;\mathbb{Z}_2)$ so that $\{m, w
t\}$ consists of a basis for $QH_0(L)$.  As $m*m \in QH_{-2}(L) =
QH_0(L) \otimes t$ we can write
\begin{equation} \label{Eq:m*m-1}
   m*m = s_1 m t + s_2 w t^2,
\end{equation}
for some $s_1, s_2 \in \mathbb{Z}_2$. We first claim that the
coefficient $s_1$ does not depend on the choice of the element $m$.
Indeed, if we replace $m$ by $m'=m+w t$ then
$$m'*m' = m*m + w t^2 = s_1m t + (s_2+1)w t^2 = s_1 m' t +
(s_1+1+s_2)w t^2.$$
This also shows that when $s_1=1$ the coefficient
$s_2$ is also independent of the choice of $m$.

\begin{prop} \label{T:m*m-1}
   The coefficient $s_1$ can be computed as follows.  Let $J \in
   \mathcal{J}$ be a generic almost complex structure. Let $p_1, p_2,
   p_3 \in L$ be a generic triple of points. Let $\ell(p_1,p_2)$,
   $\ell(p_2, p_3)$, $\ell(p_3, p_1) \subset L$ be three paths such
   that $\ell(p_i,p_j)$ connects $p_i$ to $p_j$ and intersects
   $\delta_{p_k}(J)$ transversely, where $p_k \in \{p_1, p_2, p_3\}$
   is the third point (i.e. $k \neq i,j$). Then:
   \begin{equation} \label{Eq:s1}
      s_1 = \#_{\mathbb{Z}_2} \bigl(\delta_{p_1}(J) \cap \ell(p_2,
      p_3)\bigr) + \#_{\mathbb{Z}_2} \bigl(\delta_{p_2}(J) \cap
      \ell(p_3, p_1)\bigr) + \#_{\mathbb{Z}_2} \bigl(\delta_{p_3}(J)
      \cap \ell(p_1,
      p_2)\bigr).
   \end{equation}
   In case $s_1=1$ the coefficient $s_2$ can be computed as follows.
   Denote by $n_4(p_1, p_2, p_3)$ the number modulo $2$ of simple
   $J$-holomorphic disks $u:(D, \partial D) \to (M, L)$ with
   $\mu([u])=4$ and such that $u(e^{-2\pi i l/3})=p_l$ for every
   $1\leq l \leq 3$. Then:
   \begin{equation} \label{Eq:s2}
      s_2 = \#_{\mathbb{Z}_2} \bigl(\delta_{p_3}(J) \cap \ell(p_1,
         p_2)\bigr) \cdot \#_{\mathbb{Z}_2}
         \bigl(\delta_{p_1}(J) \cap \ell(p_2, p_3)\bigr) +
         n_4(p_1, p_2, p_3).
   \end{equation}
   Moreover, when $s_1=1$ we have the following identities:
   \begin{equation} \label{Eq:delta-p123}
      \begin{aligned}
         & \#_{\mathbb{Z}_2} \bigl(\delta_{p_2}(J) \cap \ell(p_3,
         p_1)\bigr) \cdot \#_{\mathbb{Z}_2}
         \bigl(\delta_{p_3}(J) \cap \ell(p_1, p_2)\bigr) \\
         = \, & \#_{\mathbb{Z}_2} \bigl(\delta_{p_3}(J) \cap \ell(p_1,
         p_2)\bigr) \cdot \#_{\mathbb{Z}_2}
         \bigl(\delta_{p_1}(J) \cap \ell(p_2, p_3)\bigr) \\
         = \, & \#_{\mathbb{Z}_2} \bigl(\delta_{p_1}(J) \cap \ell(p_2,
         p_3)\bigr) \cdot \#_{\mathbb{Z}_2} \bigl(\delta_{p_2}(J) \cap
         \ell(p_3, p_1)\bigr).
      \end{aligned}
   \end{equation}
   \begin{equation} \label{Eq:n4}
       n_4(p_1, p_2, p_3) = n_4(p_1, p_3, p_2).
   \end{equation}
\end{prop}
\begin{rem} \label{R:m*m}
   \begin{enumerate}
     \item It is interesting to note that the numbers $n_4(p_1, p_2,
      p_3)$ depend on the position of the three points $p_1, p_2, p_3$
      (and of course on $J$). Thus the number of $\mu=4$ simple
      $J$-holomorphic disks with boundary passing through $3$ points
      on $L$ is {\em not a symplectic invariant}.
     \item The $J$-holomorphic disks counted by $n_4(p_1, p_2, p_3)$
      are different than those counted by $n_4(p_1, p_3, p_2)$.
      Nevertheless these two numbers (mod $2$) are equal. Thus the
      number of $J$-holomorphic disks of Maslov index $4$ whose
      boundary passes through $p_1, p_2, p_3$ in any possible order is
      even.
     \item Essentially the same result as Proposition~\ref{T:m*m-1} holds
      for any orientable monotone Lagrangian $L$ with $N_L=2$ and
      $\textnormal{genus}(L)>0$.  However, other than tori we are not
      aware of any examples of such Lagrangians. Note that if such a
      Lagrangian $L^2 \subset M^4$ exists then $M$ cannot have $b_2^+
      = 1$. By classification results this means that $M$ cannot be a
      monotone symplectic manifold (i.e. $c_1(M) = \lambda [\omega]
      \in H^2(M;\mathbb{R})$ for some $\lambda>0$) or even
      birationally equivalent to it. Similarly $M$ cannot be a (blow
      up of a) ruled surface.
   \end{enumerate}
\end{rem}

\begin{proof}[Proof of Proposition~\ref{T:m*m-1}]
   Fix a generic $J \in \mathcal{J}(M, \omega)$. Let $f,g:L \to
   \mathbb{R}$ be two perfect Morse functions. Fix a generic
   Riemannian metric on $L$. Denote by $x_0$ (resp. $y_0$) the minimum
   of $f$ (resp. $g$) and by $x_2$ (resp. $y_2$) the maximum of $f$
   (resp. $g$). We can choose the functions $f,g$ so that $x_0=p_1$,
   $y_0=p_2$, $x_2=p_3$ and so that $y_2$ is very close to $x_2=p_3$.

   With this data fixed we have two versions of the quantum cap
   product:
   \begin{align*}
      *_{f,g} & : QH(f,J) \otimes QH(g,J) \to QH(f,J), \\
      * & : QH(f,J) \otimes QH(f,J) \to QH(f,J).
   \end{align*}
   The relation between these products is that for $x \in QH(f,J)$, $y
   \in QH(g,J)$, we have $x*_{f,g}y = x*\phi^F(y)$, where $\phi^F$ is
   the isomorphism from Proposition~\ref{P:QH-fg}.

   Viewing $x_0, x_2, y_0, y_2$ as elements of $QH(f,J)$ and $QH(g,J)$
   we can write:
   \begin{equation} \label{Eq:x0*y0} x_0 *_{f,g} y_0 = r_1 x_0 t + r_2
      x_2 t, \quad \textnormal{for some } r_1, r_2 \in \mathbb{Z}_2.
   \end{equation}
   From the definition of the quantum cap product it follows that the
   coefficient $r_1$ is given by counting trajectories as in
   figure~\ref{f:coef-s1-1}.
   \begin{figure}[htbp]
      \begin{center}
         \epsfig{file=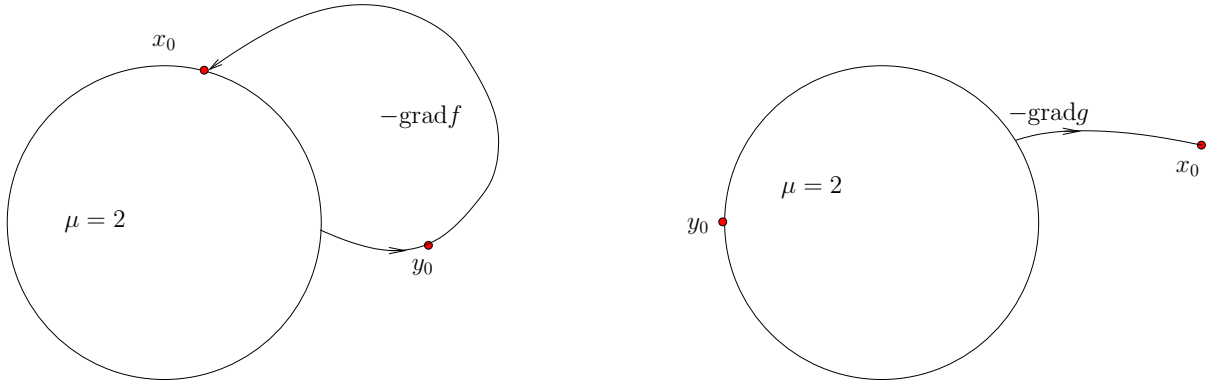, width=1.0\linewidth}
      \end{center}
      \caption{Pearly trajectories contributing to the coefficient
        $r_1$.}
      \label{f:coef-s1-1}
   \end{figure}
   Denote by $\gamma_{x_2,y_0}(f)$ the $-\textnormal{grad}f$
   trajectory connecting $x_2$ to $y_0$ and by $\gamma_{y_2,x_0}(g)$
   the $-\textnormal{grad} g$ trajectory connecting $y_2$ to $x_0$.
   (by taking the Riemannian metric generic we may assume that $y_0$
   and $x_0$ lie on trajectories as above). It follows that
   \begin{equation} \label{Eq:r_1-1} r_1 = \#_{\mathbb{Z}_2} \bigl(
      \delta_{x_0}(J) \cap \gamma_{x_2,y_0}(f) \bigr) +
      \#_{\mathbb{Z}_2} \bigl( \delta_{y_0}(J) \cap
      \gamma_{y_2,x_0}(g) \bigr).
   \end{equation}
   See figure~\ref{f:coef-s1-2}.
   \begin{figure}[htbp]
      \begin{center}
         \epsfig{file=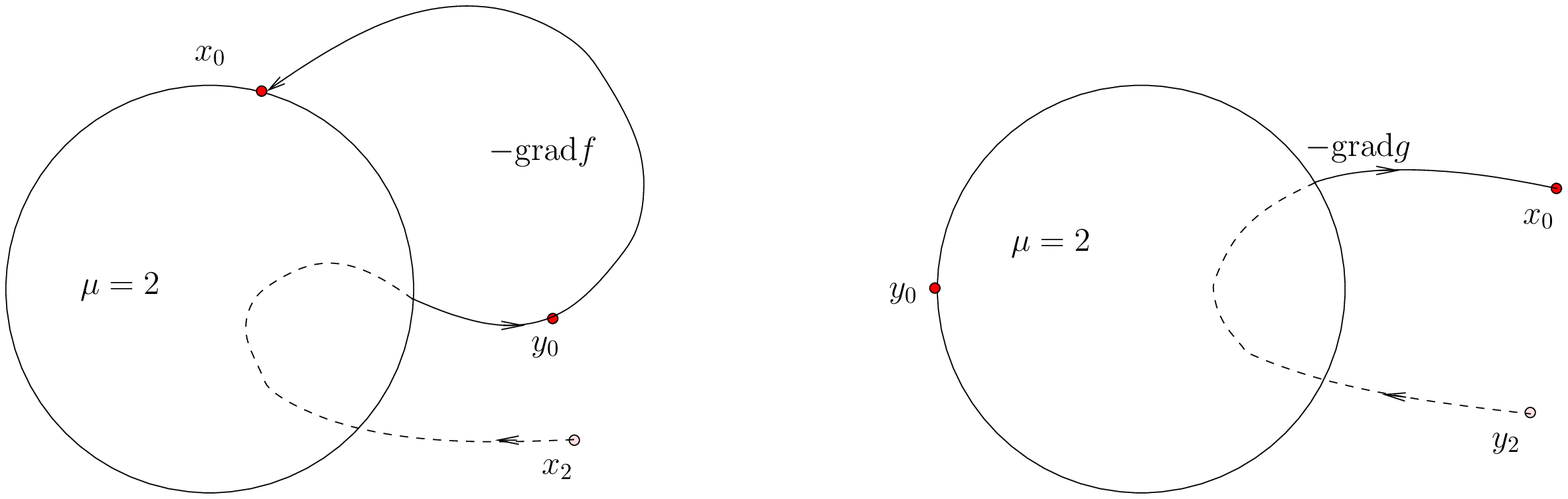, width=1.0\linewidth}
      \end{center}
      \caption{}
      \label{f:coef-s1-2}
   \end{figure}
   Since $\delta_{x_0}(J)$, $\delta_{y_0}(J)$ are null-homologous,
   $\gamma_{x_2,y_0}(f)$, $\gamma_{y_2,x_0}(g)$ from
   formula~\eqref{Eq:r_1-1} can be replaced by any paths $\ell(x_2,
   y_0)$, $\ell(y_2, x_0)$ in $L$ connecting $x_2$ to $y_0$ and $y_2$
   to $x_0$. Moreover, since $y_2$ can be chosen arbitrary close to
   $x_2$ we can replace $\ell(y_2, x_0)$ by any path $\ell(x_2, x_0)$
   connecting $x_2$ to $x_0$. Thus
   \begin{equation} \label{Eq:r_1-2} r_1 = \#_{\mathbb{Z}_2} \bigl(
      \delta_{x_0}(J) \cap \ell(x_2,y_0) \bigr) + \#_{\mathbb{Z}_2}
      \bigl(\delta_{y_0}(J) \cap \ell(x_2,x_0) \bigr).
   \end{equation}
   By Proposition~\ref{P:QH-fg} $$\phi^F(y_0) = x_0 +
   \#_{\mathbb{Z}_2} \bigl( \delta_{x_2}(J) \cap \ell(y_0,x_0)
   \bigr)x_2t.$$
   Here we have used the fact that $y_2$ is very close
   to $x_2$. Therefore:
   \begin{align*}
      x_0 * x_0 & = x_0 *_{f,g} (\phi^F)^{-1}(x_0) = x_0 *_{f,g}
      \Bigl( y_0 + \#_{\mathbb{Z}_2} \bigl( \delta_{x_2}(J) \cap
      \ell(y_0,x_0)
      \bigr)y_2t \Bigr) \\
      & = \Bigl( r_1+ \#_{\mathbb{Z}_2} \bigl( \delta_{x_2}(J) \cap
      \ell(y_0,x_0) \bigr) \Bigr) x_0t + r_2x_2 t^2 \\
      & = \Bigl( \#_{\mathbb{Z}_2} \bigl( \delta_{x_0}(J) \cap
      \ell(x_2,y_0) \bigr) + \#_{\mathbb{Z}_2} \bigl( \delta_{y_0}(J)
      \cap \ell(x_2,x_0) \bigr) + \#_{\mathbb{Z}_2} \bigl(
      \delta_{x_2}(J) \cap \ell(y_0,x_0) \bigr) \Bigr)x_0t \\
      & \;\;\; + r_2 x_2 t^2.
   \end{align*}
   Finally, take the element $m$ to be $x_0$ and recall that
   $x_0=p_1$, $y_0=p_2$, $x_2=p_3$. Formula~\eqref{Eq:s1} follows.

   We turn to the proof of formula~\eqref{Eq:s2}. From the
   computations above we see that $s_2=r_2$, i.e. the coefficient of
   $x_2t^2$ in $x_0 *_{f,g} y_0$. By the definition of the quantum cap
   product this coefficient counts the number of trajectories that
   appear in figure~\ref{f:coef-s2-1}.
   \begin{figure}[htbp]
      \begin{center}
         \epsfig{file=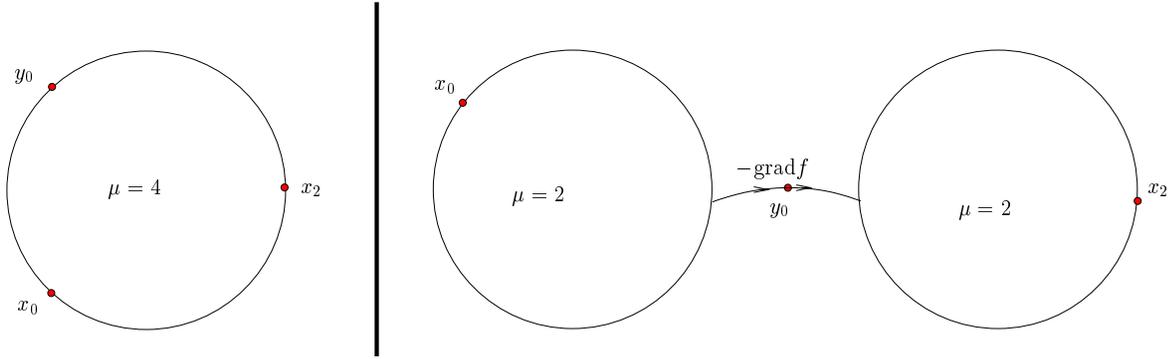, width=1.0\linewidth}
      \end{center}
      \caption{Pearly trajectories contributing to the coefficient
        $s_2=r_2$.}
      \label{f:coef-s2-1}
   \end{figure}
   The trajectories in the left-hand side of figure~\ref{f:coef-s2-1}
   contribute $n_4(x_0, y_0, x_2)=n_4(p_1, p_2, p_3)$ to the
   coefficient $r_2$.  As for the trajectories in the right-hand side,
   a similar argument to what we had before shows that their number
   is:
   \begin{align*}
      \#_{\mathbb{Z}_2}\bigl(\delta_{x_0}(J) \cap \gamma_{x_2,
        y_0}(f)\bigr) \cdot \#_{\mathbb{Z}_2}\bigl(\delta_{x_2}(J)
      \cap \gamma_{y_0,
        x_0}(f)\bigr) \\
      = \#_{\mathbb{Z}_2}\bigl(\delta_{x_0}(J) \cap \ell(x_2,
      y_0)\bigr) \cdot \#_{\mathbb{Z}_2}\bigl(\delta_{x_2}(J) \cap
      \ell(y_0, x_0)\bigr) \\
      = \#_{\mathbb{Z}_2}\bigl(\delta_{p_1}(J) \cap \ell(p_2, p_3)
      \bigr) \cdot \#_{\mathbb{Z}_2}\bigl(\delta_{p_3}(J) \cap
      \ell(p_1, p_2)\bigr).
   \end{align*}
   See figure~\ref{f:coef-s2-2}.  This completes the proof of
   formula~\eqref{Eq:s2}.
   \begin{figure}[htbp]
      \begin{center}
         \epsfig{file=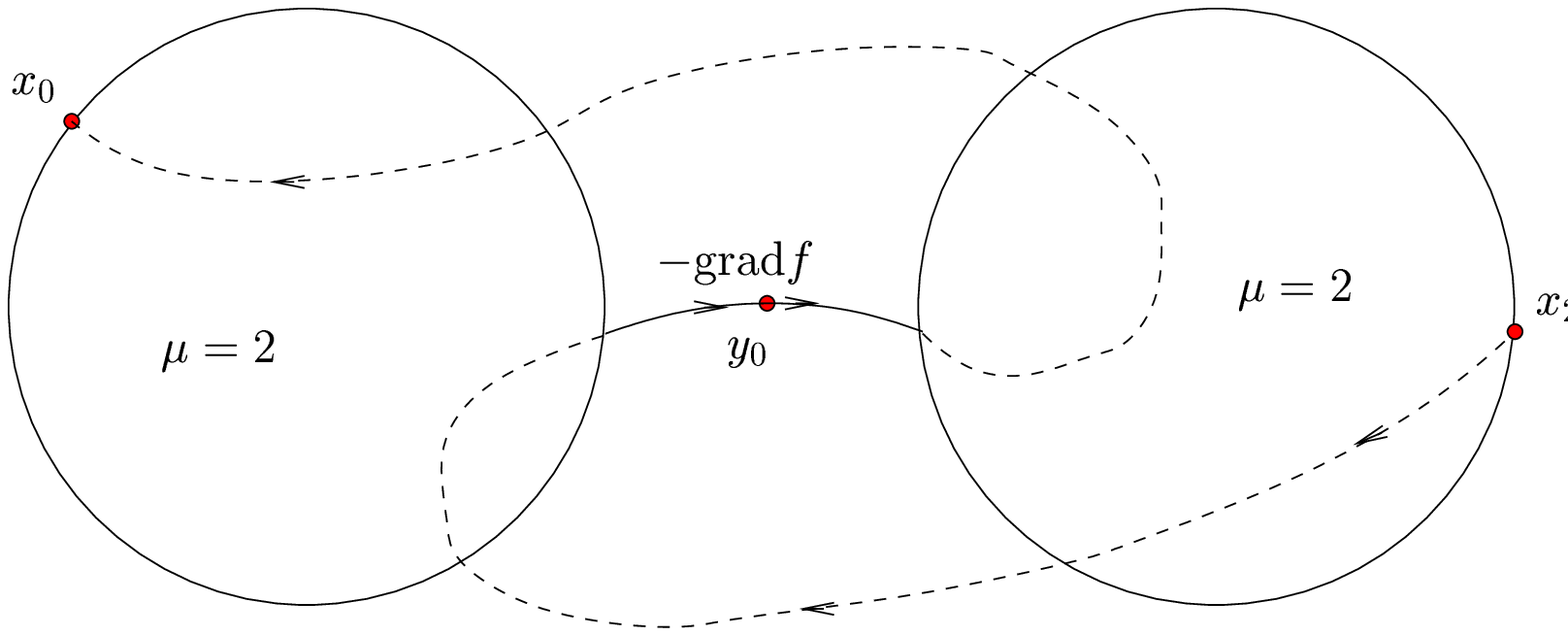, width=0.6\linewidth}
      \end{center}
      \caption{}
      \label{f:coef-s2-2}
   \end{figure}

   Formulae~\eqref{Eq:delta-p123} and~\eqref{Eq:n4} follow from the
   fact that when $s_1=1$ the coefficient $s_2$ does not depend on the
   choice of the element $m$. Thus we can take the points $x_0, x_2,
   y_2$ to be any permutation of the points $p_1, p_2, p_3$ and
   formula~\eqref{Eq:s2} will still give us the same number $s_2$.
   Note also that when $s_1=1$ it follows from formula~\eqref{Eq:s1}
   that either all three numbers $$\#_{\mathbb{Z}_2}
   \bigl(\delta_{p_1}(J) \cap \ell(p_2, p_3)\bigr), \quad
   \#_{\mathbb{Z}_2} \bigl(\delta_{p_2}(J) \cap \ell(p_3, p_1)\bigr),
   \quad \#_{\mathbb{Z}_2} \bigl(\delta_{p_3}(J) \cap \ell(p_1,
   p_2)\bigr)$$
   are $1$, or precisely two of them are $0$ and one of
   them is $1$.
\end{proof}

Let us examine in view of Proposition~\ref{T:m*m-1} the case of the
$2$-dimensional Clifford torus $L = \mathbb{T}^2_{\textnormal{clif}}
\subset {\mathbb{C}}P^2$. By Proposition~\ref{prop:clifford} there
exist generators $a,b$ of $H_1(L;\mathbb{Z}_2)$ and a basis for
$QH_0(L)$ of the form $\{ m, w t \}$ such that $$a*b = m+w t, \quad
b*a = m, \quad a*a=b*b=w t.$$
Using the associativity of the quantum
product we obtain: $m*a = b*a*a = b*w t = b t$.  Therefore:
$$m*m = b*a*b*a = b*(m+w t)*a = b*(b t + a t) = m t + w t^2.$$
Thus
$s_1 = s_2 = 1$. Since $s_1=1$ the coefficient $s_2=1$ is independent
of the choice of $m$.  It follows from Proposition~\ref{T:m*m-1} that the
following three numbers
$$\#_{\mathbb{Z}_2} \bigl(\delta_{p_1}(J) \cap \ell(p_2, p_3)\bigr),
\quad \#_{\mathbb{Z}_2} \bigl(\delta_{p_2}(J) \cap \ell(p_3,
p_1)\bigr), \quad \#_{\mathbb{Z}_2} \bigl(\delta_{p_3}(J) \cap
\ell(p_1, p_2)\bigr)$$
can be either all $1$ or exactly one of them is
$1$ and the other two are $0$. Moreover, in the second case there
exists a simple $J$-holomorphic disk $u$ with $\mu([u])=4$ and such
that $u(\partial D) \ni p_1, p_2, p_3$.

\begin{figure}[htbp]
   \begin{center}
      \epsfig{file=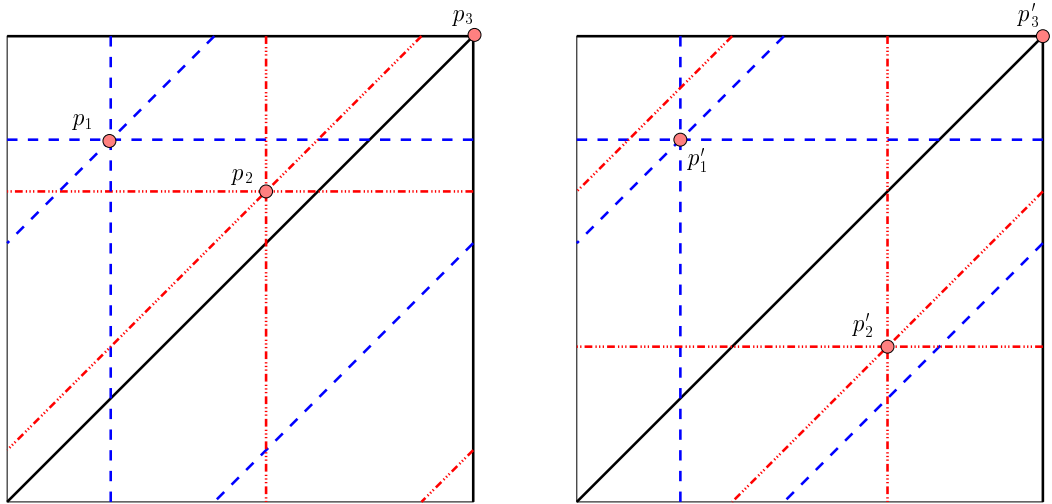, width=0.9\linewidth}
   \end{center}
   \caption{On the left $n_4(p_1, p_2, p_3)=1$. On the right
     $n_4(p'_1, p'_2, p'_3)=0$.}
   \label{f:T2-3pts-1}
\end{figure}

It is instructive to consider the case when $J=J_0$ is the standard
complex structure of ${\mathbb{C}}P^2$ (or a small perturbation of
it). In this case, for every $p \in L$ the cycle $\delta_{p}(J_0)$
consist of three embedded circles passing through $p$ (see
\S\ref{Sb:tclif-QH-1}). It is easy to see that $L \setminus
\delta_{p}(J_0)$ has {\em two} connected components.  Thus, if $p_1,
p_2$ lie in the same connected component of $L \setminus
\delta_{p_3}(J_0)$ then for a small enough perturbation $J$ of $J_0$
we will have $n_4(p_1, p_2, p_3)=1$. Of course, it is possible to find
a different configurations of points $p'_1, p'_2, p'_3$ for which
$n_4(p'_1, p'_2, p'_3)=0$. Figure~\ref{f:T2-3pts-1} shows these two
possibilities. In this figure the torus
$\mathbb{T}^2_{\textnormal{clif}}$ is represented as a square with
opposite sides identified. The three lines through each of the points
$p_i$ represent the boundaries of the three $J_0$-holomorphic disks
passing through $p_i$. Using formula~\eqref{Eq:s2} and the fact that
$s_2=1$ it is easy to see that $n_4(p_1, p_2, p_3)=1$ while $n_4(p'_1,
p'_2, p'_3)=0$.

\bigskip We continue with analyzing the coefficient $s_2$ from
formula~\eqref{Eq:m*m-1}. As discussed above the coefficient $s_2$ may
depend on the choice of $m$ (this happens only when $s_1=0$). A
natural class of choices for $m$ can be made as follows. Choose a
generic almost complex structure $J$ and a generic pair of points $p,
q \in L$. Choose a Riemannian metric $\rho$ on $L$ and a perfect Morse
function $f: L \to \mathbb{R}$ so that $p$ is its single (local)
minimum and $q$ is its single (local) maximum. As noted before the
differential of the pearl complex vanishes hence we can view $p$ as an
element of $QH_0(f,\rho,J)$. It follows from
formula~\eqref{Eq:coef-eps-2} of Proposition~\ref{P:QH-fg} that this
element is independent of the function $f$ and metric $\rho$ as long
as $f$ has $p$ and $q$ as its minimum and maximum. (Of course, once
$J$ is fixed, a slight perturbation of $p$ and $q$ will define the
same element.) Thus a choice of a generic pair of points $p,q$ and $J$
determines in a canonical way an element $m(p,q,J) = [p] \in QH_0(L)$.
Put $w = [q]$. Clearly $\{ m(p,q,J), w t \}$ form a basis for
$QH_0(L)$. Note also that different functions $f$ as above may give
rise to different generators for $QH_1(L) \cong H_1(L;\mathbb{Z}_2)$
however, the element $m(p,q,J)$ is not affected by this.

Fix a generic $J$ and a generic pair of points $p,q$. Denote by
$s_2(p,q,J)$ the coefficient $s_2$ corresponding to the choice of
$m=m(p,q,J)$. The following Proposition shows how the coefficient
$s_2$ changes when we change $p,q$ to a another pair $p',q'$, hence
$m(p,q,J)$ to $m(p',q',J)$.
\begin{prop} \label{P:s2-pq}
   $s_2(p',q',J) = s_2(p,q,J) + \eta(s_1+1)$, where $\eta$ is given
   by:
   $$\eta = \#_{\mathbb{Z}_2} \bigl( \delta_{p'}(J) \cap
   \ell(q',q)\bigr) + \#_{\mathbb{Z}_2} \bigl( \delta_{q}(J) \cap
   \ell(p',p) \bigr).$$
   In particular, when $q=q'$ we have $\eta =
   \#_{\mathbb{Z}_2} \bigl(\delta_{q}(J) \cap \ell (p',p) \bigr)$.
\end{prop}
\begin{proof}
   Put $m=m(p,q,J)$, $m'=m(p',q',J)$. By Proposition~\ref{P:QH-fg}, we
   have $m'=m + \eta wt$. A direct computation (over $\mathbb{Z}_2$ !)
   gives:
   \begin{align*}
      m'*m' & = m*m + \eta wt^2 = s_1 mt + s_2(p,q,J)wt^2 + \eta wt^2 \\
      & = s_1m't + s_1 \eta wt^2 + s_2(p,q,J)wt^2 + \eta wt^2 \\
      & = s_1 m' t + \bigl( s_2(p,q,J) + \eta(s_1+1) \bigr) wt^2.
   \end{align*}
\end{proof}

\begin{prop} \label{T:m*m-2}
   With the notation of Proposition~\ref{T:m*m-1} we have:
   \begin{equation} \label{Eq:s2-pq-n4}
      s_2(p_1,p_3,J) = \#_{\mathbb{Z}_2} \bigl(\delta_{p_3}(J) \cap \ell(p_1,
         p_2)\bigr) \cdot \#_{\mathbb{Z}_2}
         \bigl(\delta_{p_1}(J) \cap \ell(p_2, p_3)\bigr) +
         n_4(p_1, p_2, p_3).
   \end{equation}
\end{prop}
The proof is essentially the same as the proof of
Proposition~\ref{T:m*m-1}.

\subsubsection{Explicit formulae for the quantum product of
  $2$-dimensional Lagrangian tori} \label{Sb:explicit-quant-T2}

Here we develop formulae which allow us to reproduce the quantum cap
product for every (monotone) $2$-dimensional Lagrangian torus $L$ from
minimal information on holomorphic disks with boundaries on $L$.  As
it turns out, it is enough to know the number of Maslov index $2$ -
holomorphic disks passing through a generic point on $L$ in every
homology class.

Let $L^2 \subset (M^4,\omega)$ be a $2$-dimensional monotone
Lagrangian torus with minimal Maslov number $N_L=2$. Assume that
$QH_*(L) \neq 0$ (which is equivalent to $QH_*(L) \cong
(H(L;\mathbb{Z}_2) \otimes \Lambda)_*$. See
Proposition~\ref{P:criterion-QH=0-1}).

Fix generators $a, b$ of the {\em integral} homology
$H_1(L;\mathbb{Z})$ so that $H_1(L;\mathbb{Z}) = \mathbb{Z} a \oplus
\mathbb{Z} b$. Define a function $\nu:\mathbb{Z} \oplus \mathbb{Z} \to
\mathbb{Z}_2$ as follows:
\begin{equation} \label{Eq:nu-function}
   \nu(k,l) = \sum_{\substack{A \in \mathcal{E}_2, \\
       \partial A = ka+lb}} \deg_{\mathbb{Z}_2} ev_{A,J}.
\end{equation}
where $J \in \mathcal{J}_{\textnormal{reg}}$ is a generic almost
complex structure. In other words, $\nu(k,l)$ counts, mod $2$, the
number of $J$-holomorphic disks of Maslov index $\mu=2$ passing
through a generic point in $L$ and whose boundaries realize the
homology class $ka+lb \in H_1(L;\mathbb{Z})$. By the discussion at the
beginning of \S\ref{Sb:criteria-QH}, $\nu(k,l)$ does not depend on the
choice of $J \in \mathcal{J}_{\textnormal{reg}}$. Moreover,
$\nu(k,l)=0$ for all but a finite number of pairs $(k,l)$.

Denote by $w \in H_2(L;\mathbb{Z}_2)$ the fundamental class. Let $m
\in QH_0(L;\mathbb{Z}_2)$ be an element so that $\{m, w t\}$ forms a
basis for $QH_0(L)$. Then we can write
\begin{equation} \label{Eq:quant-prod-formulae-1}
   \begin{aligned}
      & a*a = \alpha w t, \quad b*b = \beta w t, \\
      & a*b = m + \gamma' w t, \quad b*a = m + \gamma'' w t, \\
      & a*b + b*a = (\gamma'+\gamma'') w t,
   \end{aligned}
\end{equation}
for some $\alpha, \beta, \gamma', \gamma'' \in \mathbb{Z}_2$.

\begin{prop} \label{T:nu-function-product}
   The coefficients $\alpha, \beta$ are given by:
   $$\alpha = \sum_{k,l} \nu(k,l) \frac{l(l+1)}{2} \,(\bmod{2}), \quad
   \beta = \sum_{k,l} \nu(k,l) \frac{k(k+1)}{2} \,(\bmod{2}).$$
   The
   sum $\gamma'+\gamma''$ is independent of the choice of the element
   $m$ and we have:
   $$\gamma'+\gamma'' = \sum_{k,l} \nu(k,l)kl \,(\bmod{2}).$$
\end{prop}
The proof is given in \S\ref{Sb:prf-nu-function} below. Note that
using the coefficients $\alpha, \beta, \gamma', \gamma''$ we can
recover the quantum product. Indeed, a simple computation based
on~\eqref{Eq:quant-prod-formulae-1} gives:
\begin{equation} \label{Eq:quant-prod-formulae-2}
   \begin{aligned}
      & m*a = \alpha b t+\gamma''a t, \quad a*m = \alpha
      b t + \gamma' a t \\
      & m*b = \beta a t + \gamma' b t, \quad b*m = \beta
      a t + \gamma'' b t \\
      & m*m = (\gamma'+\gamma'') m t + (\alpha\beta +
      \gamma'\gamma'')w t^2.
   \end{aligned}
\end{equation}
Notice that the quantum product is commutative iff
$\gamma'+\gamma''=0$ (recall that we work here over $\mathbb{Z}_2$).
Thus the function $\nu(k,l)$ determines whether or not the quantum
product is commutative. It is also worth noting that when
$\gamma'+\gamma''=1$ we must have $\gamma'\gamma''=0$, hence in this
case $m*m = m t + \alpha\beta w t^2$.

Let us apply the above to two examples.  We start with our favorite
example, the $2$-dimensional Clifford torus $L
=\mathbb{T}^2_{\textnormal{clif}} \subset {\mathbb{C}}P^2$.  We will
use here the notation of \S\ref{Sb:tclif-QH-1}. We take the basis
$a,b$ for $H_1(L;\mathbb{Z})$ to be $a=\partial A_1$, $b=\partial
A_2$. With this choice the function $\nu:\mathbb{Z} \oplus \mathbb{Z}
\to \mathbb{Z}_2$ is:
\begin{align*}
   & \nu(1,0)=1, \quad \nu(0,1)=1, \quad \nu(-1,-1)=1, \\
   & \nu(k,l)=0 \textnormal{ for all other } k,l.
\end{align*}
It easily follows from Proposition~\ref{T:nu-function-product} that
$\alpha=\beta=1$ and that $\gamma'+\gamma''=1$. It follows that
precisely one of $\gamma'$, $\gamma''$ equals $0$ and the other one
equals $1$. We thus recover again the quantum product for
$\mathbb{T}^2_{\textnormal{clif}}$.

Our second example is the split Lagrangian torus in $S^2 \times S^2$.
Endow $S^2 \times S^2$ with the split symplectic form $\omega =
\omega_{S^2} \oplus \omega_{S^2}$, where $\omega_{S^2}$ is the
standard symplectic form of $S^2$. Let $Eq \subset S^2$ be the
equator.  Then $L = Eq \times Eq \subset S^2 \times S^2$ is a monotone
Lagrangian with $N_L=2$. Denote by $D_0, D_1$ the two oriented disk
obtained from $S^2 \setminus Eq$. Put
$$A_0=[D_0 \times pt], A_1 = [D_1 \times pt], B_0=[pt \times D_0], B_1
= [pt \times D_1] \in H_2(S^2 \times S^2, L;\mathbb{Z}).$$
Clearly
$H^D_2 = H_2(S^2\times S^2, L; \mathbb{Z})\cong \mathbb{Z} A_0 \oplus
\mathbb{Z}A_1 \oplus \mathbb{Z} B_0 \oplus \mathbb{Z} B_1$.

Let $J_0 = j_0 \oplus j_0$ be the standard split complex structure
where $j_0$ is the complex structure of $S^2 \cong \mathbb{C}P^1$. Let
$u_0, u_1:(D,\partial) \to (S^2, Eq)$ be the obvious $j_0$-holomorphic
disks of Maslov index $2$ parametrizing $D_0$ and $D_1$. It is easy to
see that the only $J_0$-holomorphic disks with Maslov index $2$ are in
one of the classes $A_0, A_1, B_0, B_1$ and in fact they are all given
(up to reparametrization by an element of $\textnormal{Aut}(D)$) by
$$u_0 \times pt, \quad u_1 \times pt, \quad pt \times u_0, \quad pt
\times u_1.$$
Moreover $J_0$ is regular for all the classes $A \in
H_2^D$ with Maslov index $2$. It follows that for every $p=(x,y) \in
L$ we have $\delta_p(J_0) = Eq \times x - Eq \times x + y \times Eq -
y \times Eq$, hence $D_1=0$ (See~\eqref{Eq:delta-x}
and~\eqref{Eq:D1}). By Proposition~\ref{P:criterion-QH=0-1} we have
$QH_*(L) \cong (H(L) \otimes \Lambda)_*$. (This can be easily verified
also by the Floer K\"{u}nneth formula).

We now fix the following generators $a=[Eq \times pt], b=[pt \times
Eq] \in H_1(L;\mathbb{Z})$. Let $m \in QH_0(L)$ be an element such
that $\{m, w t\}$ forms a basis for $QH_0(L)$. The function
$\nu:\mathbb{Z} \oplus \mathbb{Z} \to \mathbb{Z}_2$ is:
\begin{align*}
   & \nu(1,0)=1, \quad \nu(-1,0)=1, \quad, \nu(0,1)=1, \quad \nu(0,-1)=1, \\
   & \nu(k,l)=0 \textnormal{ for all other } k,l.
\end{align*}
By Proposition~\ref{T:nu-function-product} we have $\alpha=1$, $\beta=1$
and $\gamma'+\gamma''=0$ (i.e. $\gamma'=\gamma''$). It follows that
the coefficient $s_1$ in $m*m$ is $0$. This means that the coefficient
$s_2$ may depend on the choice of $m$. To choose $m$ let $p,q \in L$
be two generic points. Let $J_{\epsilon}$ be a small enough
perturbation of $J_0$. This defines the element $m =
m(p,q;J_{\epsilon})$. A straightforward computation shows that $a*b =
m$ and $b*a=m$, hence $\gamma'=\gamma''=0$.
Using~\eqref{Eq:quant-prod-formulae-1}
and~\eqref{Eq:quant-prod-formulae-2} we obtain (for the above choice
of $m$):
\begin{align*}
   & a*a=b*b=w t, \quad a*b=b*a=m, \\
   & a*m=m*a=b t, \quad b*m=m*b=a t, \\
   & m*m = w t^2.
\end{align*}
Thus $s_2(p,q,J_{\epsilon})=1$.

Let us interpret the value of $s_2$ using Proposition~\ref{T:m*m-2}.  Put
$p_1=p$, $p_3=q$ and let $p_2 \in L \setminus \delta_{p_3}(J_0)$ be a
generic point. Note that $L \setminus \delta_{p_3}(J_0)$ is connected
(it is an open rectangle). If the perturbation $J_{\epsilon}$ of $J_0$
is small enough then the cycle $\delta_{p_3}(J_{\epsilon})$ will be
close to $\delta_{p_3}(J_0)$ in the $C^0$-topology. Therefore, for
small enough perturbation the point $p_1$ and $p_2$ will still remain
in the same connected component of $L \setminus
\delta_{p_3}(J_{\epsilon})$. In particular
$\#_{\mathbb{Z}_2}\bigl(\delta_{p_3}(J) \cap \ell(p_1, p_2)\bigr)=0$.
It follows from Proposition~\ref{T:m*m-2} that $n_4(p_1, p_2, p_3)=1$.  In
other words, we obtain the following.
\begin{cor} \label{C:disks-on-L-in-S2xS2}
   Let $p_1,p_2, p_3 \in L$ be generic points. Then there exists a
   neighbourhood $\mathcal{U} \subset \mathcal{J}$ of $J_0$ (which
   depends on $p_1, p_2, p_3$) such that for generic $J \in
   \mathcal{U}$ there exists a simple $J$-holomorphic disk $u:(D,
   \partial D) \to (S^2\times S^2, L)$ with $\mu(L)=4$ and such that
   $u(e^{-2\pi i l/3})=p_l$ for every $1\leq l \leq 3$.
\end{cor}

\subsubsection{Proof of Proposition~\ref{T:nu-function-product}}
\label{Sb:prf-nu-function}
Let $f:L \to \mathbb{R}$ be a perfect Morse function and $\rho$ a
Riemannian metric on $L$. Denote by $x_2$ the unique maximum of $f$,
by $x_0$ the unique minimum and by $x'_1, x''_1$ the critical points
of index $1$.  Since $a, b$ generate $H_1(L;\mathbb{Z})$ the pair
$(f,\rho)$ can be chosen so that the closures of the unstable
submanifolds $\overline{W}_{x'_1}^u$, $\overline{W}_{x''_1}^u$
represent the homology classes $\pm a, \pm b$. Fix an orientation on
$L$ so that $a\cdot b = 1$ and pick orientations $\mathbf{o}'$,
$\mathbf{o}''$ on $\overline{W}_{x'_1}^u$, $\overline{W}_{x''_1}^u$ so
that $[\overline{W}_{x'_1}^u]=a$, $[\overline{W}_{x''_1}^u]=b$.

Having fixed the above let $g:L \to \mathbb{R}$ be a small
perturbation of $f$ chosen in the following way. Denote by $y_0, y'_1,
y''_1, y_2$ the critical points of $g$ (with the convention that
subscripts denote the Morse index). We choose $g$ so that $W_{y'_1}^u$
is a ``parallel'' translate of $W_{x'_1}^u$ pushed in a direction
$\vec{\mathbf{n}}'$ normal to $W_{x'_1}^u$ such that the pair
$(\mathbf{o}', \vec{\mathbf{n}}' )$ give the positive orientation on
$L$.  Similarly we require that $W_{y''_1}^u$ is obtained by pushing
$W_{x''_1}^u$ in the direction $\vec{\mathbf{n}}''$ so that the pair
$(\mathbf{o}'', \vec{\mathbf{n}}'')$ give the negative orientation on
$L$.

By choosing $g$ to be close enough to $f$ we may assume that the
canonical quasi-isomorphism between $\mathcal{C}_*(f,J)$ and
$\mathcal{C}_*(g,J)$ is in fact a base preserving isomorphism (i.e.
$x_i$ is mapped to $y_i$, $x'_1$ to $y'_1$ and $x''_1$ to $y''_1$).
Thus in order to compute the coefficient $\alpha$ in $a*a = \alpha w
t$ we will compute $[x'_1] * [y'_1]$.

Choose a generic $J \in \mathcal{J}(M,\omega)$. Let $u:(D, \partial)
\to (M, L)$ be a $J$-holomorphic disk with Maslov index $2$ that
contributes to $\nu(k,l)$, i.e. $u(\partial D) \ni x_2$ and
$[u(\partial D)] = ka + lb$. After reparametrization we may assume
that $u(1)=x_2$. The contribution of this disk to the coefficient
$\alpha$ in $[x'_1] * [y'_1]$ is given by the number, mod $2$, of
pairs $(\theta_g, \theta_f)$ such that:
$$0 < \theta_g < \theta_f < 2\pi, \quad u(e^{i\theta_g}) \in
W^u_{y'_1}, \quad u(e^{i\theta_f}) \in W^u_{x'_1}.$$
In other words we
have to mark on the boundary $u(\partial D)$ of the disk $u$ two types
of points: the points hit by a $-\textnormal{grad} f$ trajectory
coming from $x'_1$ and the points hit by a $-\textnormal{grad} g$
trajectory coming from $y'_1$. Let us call the first set ``points of
type $f$'' and the second ``points of type $g$''. Then we have to
count how many pairs of points $(q, p)$, $p$ of type $f$ and $q$ of
type $g$, are there such that $q$ appears ``before'' $p$ along
$u(\partial D)$. We denote this number $\alpha_u$.  The wording
``before'' means that we take $x_2 \in u(\partial D)$ as the origin
and use the standard orientation of $\partial D$ for ordering.

\begin{figure}[htbp]
   \begin{center}
      \epsfig{file=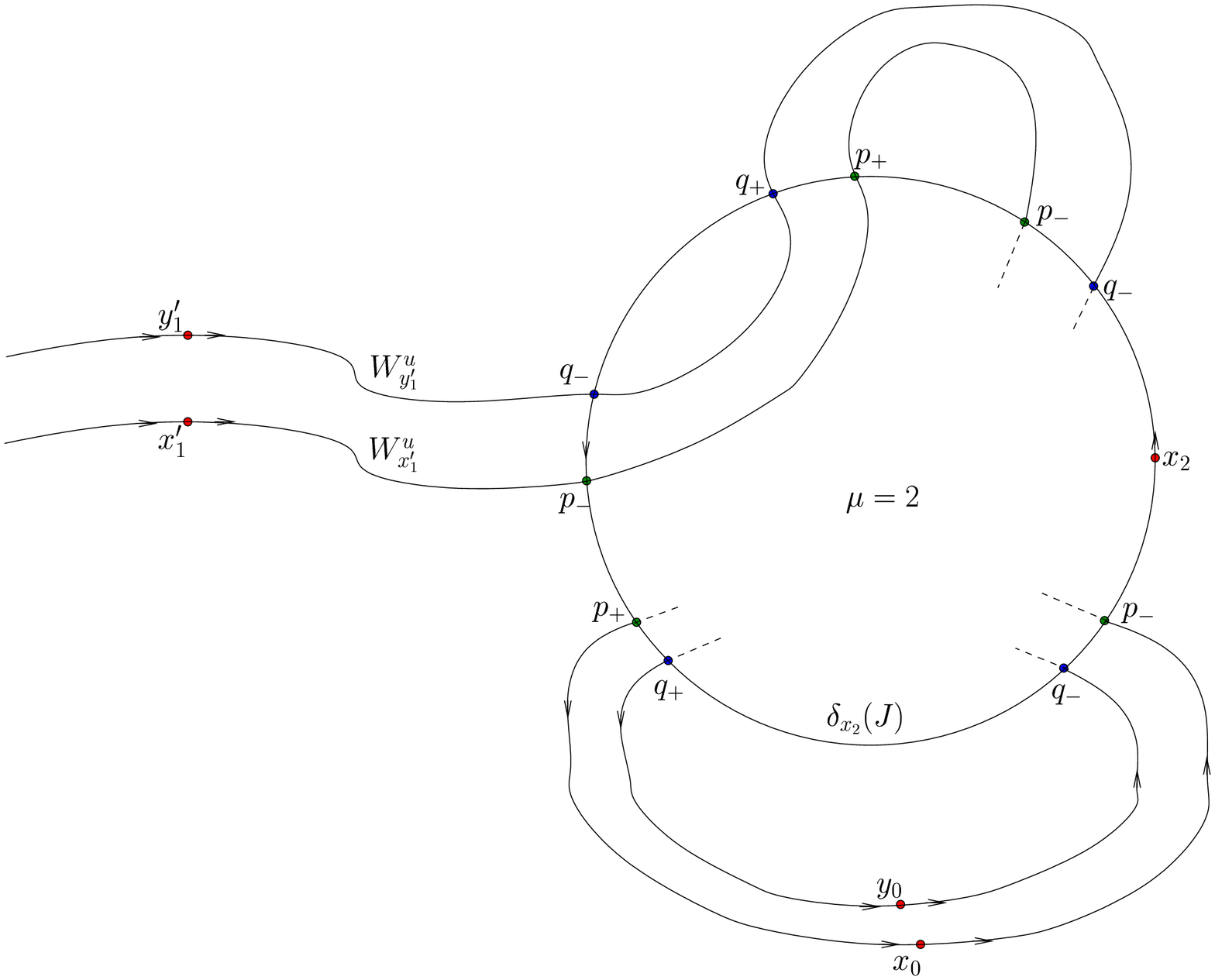, width=0.7\linewidth}
   \end{center}
   \caption{}
   \label{f:traj-aa}
\end{figure}

To compute $\alpha_u$ first note that points of types $f$ and $g$
appear in pairs (provided that the function $g$ is close enough to
$f$). In fact, each time $W_{x'_1}^u$ intersects positively the
boundary of the disk $u(\partial D)$ then $W_{y'_1}^u$ intersects it
positively too (in a nearby point) and we we obtain two points
$p_{+}$, $q_{+}$, of types $f$ and $g$, where $p_{+}$ comes before
$q_{+}$. When $W_{x'_1}^u$ intersects $u(\partial D)$ negatively we
obtain two points $q_{-}$, $p_{-}$, of types $g$ and $f$, where
$q_{-}$ comes before $p_{-}$. This is described in
figure~\ref{f:traj-aa}. Note that the arrows along $W^u_{x'_1},
W^u_{y'_1}$ in this figure might be a bit confusing. These arrows
represent the orientation $\mathbf{o}'$ discussed above, {\em not} the
direction of the gradient (or minus gradient) flow.

Denote by $i_{+}$, $i_{-}$ the number of positive, resp. negative,
intersection points of $W^u_{x'_1}$ with $u(\partial D)$. As for the
contribution to $\alpha_u$, each point $q_{\pm}$ of type $g$ (positive
or negative) can be paired to any point $p_{\pm}$ coming ``after'' it.
In addition, each point $q_{-}$ appearing in a negative intersection
as above can also be paired with its adjacent point $p_{-}$. Therefore
\begin{equation} \label{Eq:alpha_u}
   \begin{aligned}
      \alpha_u & = (\sum_{j=i_{+} + i_{-}} j ) + i_{-} \\
      & = \frac{(i_{+}+i_{-})(i_{+}+i_{-}+1)}{2} + i_{-} =
      \frac{i_{+}^2 + 2i_{+}i_{-}+i_{-}^2 + i_{+}+i_{-}}{2}+i_{-}.
   \end{aligned}
\end{equation}
On the other hand since $[u(\partial D)]=ka+lb$ we have $i_{+} - i_{-}
= [W_{x'_1}] \cdot [u(\partial D)] = l$. Therefore we obtain from
~\eqref{Eq:alpha_u}:
\begin{align*}
   \alpha_u & \equiv \frac{i_{+}^2-2i_{+}i_{-}+i_{-}^2+i_{+}-i_{-}}{2}
   (\bmod{2}) \equiv \frac{(i_{+}-i_{-})(i_{+}-i_{-}+1)}{2} (\bmod{2}) \\
   & \equiv \frac{l(l+1)}{2} (\bmod{2}).
\end{align*}
It follows that the contribution to coefficient $\alpha$ of the disks
$u$ with $[u(\partial D)] = ka + lb$ is
$\nu(k,l)\frac{l(l+1)}{2}(\bmod{2})$, hence
$$\alpha = \sum_{k,l} \nu(k,l) \frac{l(l+1)}{2} (\bmod{2}).$$
The
formula for the coefficient $\beta$ is proved in a similar way.

The formula for $\gamma'+\gamma''$ is more straightforward.  Let
$u:(D, \partial D) \to (M,L)$ be a $J$-holomorphic disk with Maslov
index $2$, with $u(1)=x_2$ and such that $[u(\partial D)]=ka+lb$.
Denote by $i'$, $i''$ the number of intersection points (counted
without signs) of $W_{x'_1}^u$, resp.  $W^u_{x''_1}$ with $u(\partial
D)$. Note that if $g$ is chosen to be close enough to $f$ we have:
$$i' = \# (W_{y'_1}^u \cap u(\partial D)), \quad i'' = \# (W_{y''_1}^u
\cap u(\partial D)).$$
Since $a*b + b*a = [x'_1]*[y''_1] +
[x''_1]*[y'_1] = (\gamma'+\gamma'')wt$ the contribution of the disk
$u$ to the coefficient $\gamma'+\gamma''$ is $i'i''$.  (This is
because we are computing the symmetric expression $a*b+b*a$, hence the
order of the $3$ points on $u(\partial D)$ involved in the quantum
product does not matter.) Since $[\overline{W}_{x'_1}^u] \cdot [u] =
l$, $[\overline{W}_{x''_1}^u] \cdot [u] = -k$ we have $i' \equiv l
(\bmod{2})$, $i'' \equiv k(\bmod{2})$. It follows that
$$\gamma'+\gamma'' \equiv \sum_{k,l} \nu(k,l)lk \,(\bmod{2}).$$
\Qed



%


\end{document}